\documentclass{book}

\input{diagrams}
\diagramstyle[repositionpullbacks,midshaft,PostScript=dvips,nohug,scriptlabels]
\usepackage{epsfig}

\usepackage{makeidx}
\makeindex

\newcommand{\mcm}[3]{\newcommand{#1}[#2]{{\ensuremath{#3}}}}

\usepackage{latexsym}
\usepackage{amssymb}

\mcm{\emptybk}{0}{\:\:}
\mcm{\blank}{0}{(\emptybk)}
\mcm{\dashbk}{0}{-}
\mcm{\hyph}{0}{\mbox{-}}

\mcm{\diagspace}{0}{\mbox{\hspace{2em}}}

\newcommand{\bref}[1]{(\ref{#1})}
\newcommand{\ucontents}[2]{\addcontentsline{toc}{#1}{\numberline{}{#2}}}

\mcm{\mb}{1}{\mathbf{#1}}
\mcm{\mc}{1}{\mathcal{#1}}
\mcm{\mi}{1}{\mathit{#1}}
\mcm{\mr}{1}{\mathrm{#1}}
\mcm{\cat}{1}{\mc{#1}}
\mcm{\fcat}{1}{\mb{#1}}
\mcm{\ovln}{1}{\overline{#1}}
\mcm{\twid}{1}{\widetilde{#1}}

\newcommand{\slsh}{/\linebreak[0]}

\newcommand{\dt}{.\linebreak[0]}

\mcm{\sub}{0}{\,\subseteq\,}
\mcm{\such}{0}{\:|\:}
\mcm{\without}{0}{\setminus}

\mcm{\ladj}{0}{\,\dashv\,}
\mcm{\of}{0}{\raisebox{0.08ex}{\ensuremath{\scriptstyle\circ}}}
\mcm{\sof}{0}{\raisebox{0.08ex}{\ensuremath{\scriptscriptstyle\circ}}}

\mcm{\bdry}{0}{\partial}
\mcm{\blob}{0}{\raisebox{.3ex}{\ensuremath{\scriptscriptstyle{\bullet}}}}
\newcommand{\epsln}{\varepsilon}
\mcm{\implies}{0}{\,\Rightarrow\,}

\mcm{\Hom}{0}{\mr{Hom}}
\mcm{\ob}{0}{\mr{ob}\,}		% As binary operation
\mcm{\op}{0}{\mr{op}}

\mcm{\comp}{0}{\mi{comp}}
\mcm{\id}{0}{\mi{id}}
\mcm{\ids}{0}{\mi{ids}}

\mcm{\Ab}{0}{\fcat{Ab}}
\mcm{\Alg}{0}{\fcat{Alg}}
\mcm{\Bicat}{0}{\fcat{Bicat}}
\mcm{\Bim}{1}{\fcat{Bim}(#1)}
\mcm{\Cat}{0}{\fcat{Cat}}
\mcm{\fc}{0}{\fcat{fc}}
\mcm{\Gph}{0}{\fcat{Gph}}
\mcm{\Graph}{0}{\fcat{Graph}}
\mcm{\Multicat}{0}{\fcat{Multicat}}
\mcm{\One}{0}{\fcat{1}}
\mcm{\Set}{0}{\fcat{Set}}
\mcm{\Span}{0}{\fcat{Span}}
\mcm{\Struc}{0}{\fcat{Struc}}
\mcm{\Sym}{0}{\fcat{Sym}}
\mcm{\Top}{0}{\fcat{Top}}
\mcm{\UBicat}{0}{\fcat{UBicat}}

\mcm{\integers}{0}{\mathbb{Z}}

\mcm{\range}{2}{#1,\,\ldots\,,#2}
\mcm{\tuplebts}{1}{(#1)}
\mcm{\bftuple}{2}{\tuplebts{\range{#1}{#2}}}
\mcm{\tuple}{3}{\tuplebts{\range{#1,#2}{#3}}}

\mcm{\eend}{2}{#1[#2]}
\mcm{\ehom}{3}{#1[#2,#3]}
\mcm{\ftrcat}{2}{[#1,#2]}

\mcm{\goesto}{0}{\,\longmapsto\,}
\mcm{\goiso}{0}{\goby{\diso}}
\mcm{\monic}{0}{\rMonic}
\mcm{\og}{0}{\lTo}
\mcm{\ogby}{1}{\lTo^{#1}}

\mcm{\oppair}{2}{\pile{\rTo^{\scriptstyle #1}\\ \lTo_{\scriptstyle #2}}}
\mcm{\parpair}{2}{\pile{\rTo^{\scriptstyle #1}\\ \rTo_{\scriptstyle #2}}}
\mcm{\parpairu}{0}{\pile{\rTo\\ \rTo}}

\mcm{\vslob}{3}
	{\left.
	\begin{diagram}[height=1.5em]
	#1		\\
	\dTo>{\,#2}	\\
	#3		\\
	\end{diagram}
	\right.}

\newarrow{Equals}=====
\newarrow{Get}....>
\newarrow{Goesto}|---{->}
\newarrow{Incl}C--->
\newarrow{Mod}--+->
\newarrow{Monic}{vee}---{vee}
\newarrow{NT}===={=>}

\newenvironment{tree}
	{\begin{diagram}[height=1em,width=.75em,abut,noPS,tight]}	
	{\end{diagram}}

\newcommand{\dn}{\dLine}
\newcommand{\lt}[1]{\ldLine(#1,2)}
\newcommand{\rt}[1]{\rdLine(#1,2)}
\mcm{\node}{0}{\bullet}
\mcm{\enode}{0}{\circ}
\mcm{\nl}{1}{\stackrel{\textstyle #1}{\node}}

\mcm{\diso}{0}{\sim}
\mcm{\vdiso}{0}{\wr}
\newcommand{\pullshape}
	{\setlength{\unitlength}{1em}
	\begin{picture}(2,5)(-1,-5)
	\put(0,-5){\line(1,1){1}}
	\put(0,-5){\line(-1,1){1}}
	\end{picture}}
\newcommand{\Spbk}{\overprint{\raisebox{-2.5em}{\pullshape}}}

\newcommand{\lbl}[1]{\label{#1}} 

\newcommand{\chapterquote}[2]{%
\parbox[b]{0.7\textwidth}{\textit{#1}}%
\parbox[b]{0.3\textwidth}{\raggedleft #2}%
\vspace*{1ex}\\}
\usepackage{graphics}
\newcommand{\da}[3]{\cell{#1}{#2}{c}{\rotatebox{#3}{$\Downarrow$}}}

\mcm{\Eee}{0}{\cat{E}}
\mcm{\Cartpr}{0}{\pr{\Eee}{T}}
\mcm{\mult}{0}{\mi{mult}}
\mcm{\unit}{0}{\mi{unit}}
\mcm{\PD}{1}{\fcat{PD}_{#1}}
\mcm{\Tr}{0}{\fcat{Tr}}
\mcm{\pr}{2}{\tuplebts{#1,#2}}
\mcm{\graph}{4}{\spaan{#1}{T #2}{#2}{#3}{#4}}
\mcm{\spaan}{5}{#2 \ogby{#4} #1 \goby{#5} #3}
\mcm{\bktdvslob}{3}
	{\left(
	\begin{diagram}[height=1.5em,scriptlabels]
	#1		\\
	\dTo>{\,#2}	\\
	#3		\\
	\end{diagram}
	\right)}
\mcm{\slob}{3}{(#1 \goby{#2} #3)}
\newarrow{Edge}---->		% for opetopes
\mcm{\UBilax}{0}{\fcat{UBicat}_\mr{lax}}
\mcm{\UBiwk}{0}{\fcat{UBicat}_\mr{wk}}
\mcm{\UBistr}{0}{\fcat{UBicat}_\mr{str}}
\mcm{\Bilax}{0}{\fcat{Bicat}_\mr{lax}}
\mcm{\Biwk}{0}{\fcat{Bicat}_\mr{wk}}
\mcm{\Bistr}{0}{\fcat{Bicat}_\mr{str}}
\mcm{\rotsub}{0}{\cup \raisebox{0.1em}{$\scriptstyle{|}$}}
\mcm{\pd}{0}{\fcat{pd}}
\mcm{\rep}{1}{\widehat{#1}}
\mcm{\tr}{0}{\fcat{tr}}
\mcm{\END}{0}{\fcat{End}}
\mcm{\HOM}{0}{\fcat{Hom}}
\newcommand{\ditto}{,,}
\newcommand{\place}[3]{\put(#1,#2){\makebox(0,0)[c]{#3}}}
\mcm{\act}{1}{\mi{act}_{#1}}
\mcm{\ofdim}{1}{\,\of_{#1}\,}

\mcm{\strcat}{1}{\fcat{Str\hyph}#1\fcat{\hyph Cat}}
\mcm{\wkcat}{1}{\fcat{Wk\hyph}#1\fcat{\hyph Cat}}
\mcm{\strc}{1}{#1\hyph\Cat}		% for use in Appendix free-strict
\mcm{\strtuplecat}{1}{\fcat{Str\hyph}#1\fcat{\hyph tuple\hyph Cat}}

\mcm{\UMClax}{0}{\fcat{UMonCat}_\mr{lax}}
\mcm{\UMCwk}{0}{\fcat{UMonCat}_\mr{wk}}
\mcm{\UMCstr}{0}{\fcat{UMonCat}_\mr{str}}
\mcm{\MClax}{0}{\fcat{MonCat}_\mr{lax}}
\mcm{\MCwk}{0}{\fcat{MonCat}_\mr{wk}}
\mcm{\MCstr}{0}{\fcat{MonCat}_\mr{str}}
\mcm{\st}{0}{\fcat{st}}
\mcm{\Operad}{0}{\fcat{Operad}}
\mcm{\CAT}{0}{\fcat{CAT}}
\newcommand{\thmname}[1]{\textbf{\textup{(#1)}}} 	% The trouble with
\mcm{\oppairu}{0}{\pile{\rTo\\ \lTo}}
\mcm{\ctr}{0}{\fcat{ctr}}
\mcm{\upr}{1}{[#1]}		% Object of Delta
\mcm{\lwr}{1}{\mathbf{#1}}	% Object of \scat{D}
\mcm{\Sp}{2}{#1_{(#2)}}
\newenvironment{fcdiagram}
	{\begin{diagram}[height=2em]}
	{\end{diagram}}
\newcommand{\disjt}{\amalg}	% {+}
\mcm{\Pd}{1}{\fcat{Pd}_{#1}}
\mcm{\TR}{0}{\fcat{TR}}
\newenvironment{quotedthm}[1]
{\begin{trivlist} \item \textbf{#1}\ \itshape}
{\end{trivlist}}
\mcm{\slind}{1}{\twid{#1}}	% For CJ-type induced functor to slice
\mcm{\Fam}{0}{\fcat{Fam}}
\mcm{\gluing}{0}{\downarrow} 	% This one just included to get spacing
\newarrow{Globopd}----{triangle}
\mcm{\newmnd}{1}{#1^\sharp}	% For free enriched categories in appendix
\mcm{\reals}{0}{\mathbb{R}}
\mcm{\complexes}{0}{\mathbb{C}}
\mcm{\cheesydisk}{1}{B^{[#1]}}	% Fraction of a disk, for Swiss cheese
\mcm{\wej}{0}{\vee}
\mcm{\SymOpd}{0}{\fcat{S}}
\newcommand{\minihead}[1]{\subsubsection*{\textmd{\it #1}}}	% For
\mcm{\ldisks}{0}{\fcat{D}}
\newcommand{\astyle}[1]{#1}	% `Axiom style' for
\newenvironment{centredpic}%
{\begin{array}{c}
\setlength{\unitlength}{1em}}%
{\end{array}}
\newcommand{\vx}{\bullet}	% Vertex of tree
\mcm{\gm}{1}{T_{(#1)}}
\mcm{\trunc}{1}{(#1)}
\newcommand{\latin}[1]{\textit{#1}}	% Prima facie, sui generis
\newcommand{\vs}{vs}			% Concise OED has it unitalicized
\newcommand{\tcs}[1]{#1}	% Two-cell label style for ch:fcm
\newenvironment{notes}{\subsection*{Notes}\small}{}
\newcommand{\bkthack}{\mbox{\rule{0mm}{2.5ex}}}	
\newenvironment{scriptarray}%
{\begin{array}{l}}%
{\end{array}}
\newenvironment{scriptarrayc}%
{\begin{array}{c}}%
{\end{array}}
\newcommand{\nent}{\rotatebox{45}{$\Rightarrow$}}
\newcommand{\nwnt}{\rotatebox{-45}{$\Leftarrow$}}
\newcommand{\swnt}{\rotatebox{45}{$\Leftarrow$}}
\newcommand{\sent}{\rotatebox{-45}{$\Rightarrow$}}
\newcommand{\neeq}{\rotatebox{45}{$=$}}
\newcommand{\nuzeropic}{{\!\!\!\begin{array}{c}\raisebox{-0.5ex}{\epsfig{file=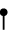}}\end{array}\!\!\!}}
\newcommand{\nutwopic}{{\!\!\!\begin{array}{c}\epsfig{file=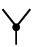}\end{array}\!\!\!}}
\newcommand{\lambdapic}{{\!\!\!\begin{array}{c}\epsfig{file=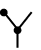}\end{array}\!\!\!}}
\newcommand{\rhopic}{{\!\!\!\begin{array}{c}\epsfig{file=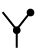}\end{array}\!\!\!}}
\newcommand{\assleftpic}{{\!\!\!\begin{array}{c}\epsfig{file=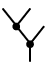}\end{array}\!\!\!}}
\newcommand{\assrightpic}{{\!\!\!\begin{array}{c}\epsfig{file=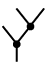}\end{array}\!\!\!}}
\mcm{\utree}{0}{\, | \,}
\newcommand{\glo}[1]{\label{glo:#1}}
\newcommand{\gr}[1]{\pageref{glo:#1}}
\newlength{\templength}    % all-purpose temporary length
\newlength{\templengthtwo} % all-purpose temporary length
\newenvironment{glosslist}%
{\small\begin{tabbing}
\hspace*{1\templength}\=\kill}%
{\end{tabbing}\normalsize}
\newcommand{\guse}[3]%
{\ensuremath{#2}\>\parbox{\templengthtwo}{#3\hfill #1}\\}

\mcm{\cod}{0}{\mr{cod}} 
\mcm{\dom}{0}{\mr{dom}}
\mcm{\elt}{0}{\in}	
\mcm{\Mon}{0}{\fcat{Mon}}
\newcommand{\demph}[1]{\textbf{\textup{#1}}}
\newcommand{\done}{\hfill\ensuremath{\Box}}
\newenvironment{prooflike}[1]{\begin{trivlist}\item\textbf{#1}\ }
{\end{trivlist}}
\newenvironment{proof}{\begin{prooflike}{Proof}}{\end{prooflike}}
\newenvironment{slopeydiag}
	{\begin{diagram}[width=1.7em,height=1.7em,scriptlabels]}
	{\end{diagram}}
\mcm{\End}{0}{\fcat{End}}
\mcm{\scat}{1}{\mathbb{#1}}
\newcommand{\go}{\rTo\linebreak[0]}
\newcommand{\goby}[1]{\rTo^{#1}\linebreak[0]}
\newcommand{\iso}{\cong}
\newcommand{\nat}{\mathbb{N}}	
\newcommand{\eqv}{\simeq}
\newcommand{\url}[1]{\texttt{#1}}
\newlength{\hdwidth}	% the width of a component
\newlength{\hdvert}	% the overall vertical measurement
\newlength{\hddrop}	% the distance a labelled glob protrudes below the
\newlength{\hdbaredrop}	% the distance from the textline to the bottom of the
\newlength{\hdoffset}	% the distance from the centre of the glob to the
\newlength{\hdtemp}	% temporary register

\newcommand{\present}[1]{%
\makebox[1\hdwidth]{%
\rule[-1\hddrop]{0ex}{1\hdvert}%
\raisebox{-1\hdbaredrop}{#1}}}

\newcommand{\presentl}[1]{%
\makebox[1\hdwidth][l]{%
\rule[-1\hddrop]{0ex}{1\hdvert}%
\raisebox{-1\hdbaredrop}{#1}}}

\newcommand{\presentr}[1]{%
\makebox[1\hdwidth][r]{%
\rule[-1\hddrop]{0ex}{1\hdvert}%
\raisebox{-1\hdbaredrop}{#1}}}

\newcommand{\sidespic}[1]{%
\settowidth{\hdtemp}{\ensuremath{#1}}%
\addtolength{\hdwidth}{1\hdtemp}}

\newcommand{\abovepic}[1]{%
\settoheight{\hdtemp}{\ensuremath{#1}}%
\addtolength{\hdvert}{1\hdtemp}%
\settodepth{\hdtemp}{\ensuremath{#1}}%
\addtolength{\hdvert}{1\hdtemp}}

\newcommand{\belowpic}[1]{%
\settoheight{\hdtemp}{\ensuremath{#1}}%
\addtolength{\hdvert}{1\hdtemp}%
\addtolength{\hddrop}{1\hdtemp}%
\settodepth{\hdtemp}{\ensuremath{#1}}%
\addtolength{\hdvert}{1\hdtemp}%
\addtolength{\hddrop}{1\hdtemp}}

\newcommand{\cell}[4]{\put(#1,#2){\makebox(0,0)[#3]{\ensuremath{#4}}}}
\mcm{\zmark}{0}{\scriptstyle{\bullet}}

\newcommand{\ginitdims}[2]{%		% LARGE GLOBULAR VERSION
\setlength{\unitlength}{1em}%		% unitlength = 1em
\setlength{\hdoffset}{.25\unitlength}%	% globular offset = .25em
\setlength{\hdwidth}{#1\unitlength}%	% width as specified
\setlength{\hdvert}{#2\unitlength}%	% vert = #2
\setlength{\hddrop}{.5\hdvert}%		% 
\addtolength{\hddrop}{-1\hdoffset}%	% 
\setlength{\hdbaredrop}{1\hddrop}%	% hddrop = drop = half(vert) - offset
\addtolength{\hdvert}{.6\unitlength}%	% total extra clearance of .6em...
\addtolength{\hddrop}{.3\unitlength}}	% ...half of which is at bottom

\newcommand{\gsinitdims}[2]{%		% SMALL GLOBULAR VERSION
\setlength{\unitlength}{0.5em}%		% unitlength = 0.5em
\setlength{\hdoffset}{.25\unitlength}%	% globular offset = .25em
\setlength{\hdwidth}{#1\unitlength}%	% width as specified
\setlength{\hdvert}{#2\unitlength}%	% vert = #2
\setlength{\hddrop}{.5\hdvert}%		% 
\addtolength{\hddrop}{-1\hdoffset}%	% 
\setlength{\hdbaredrop}{1\hddrop}%	% hddrop = drop = half(vert) - offset
\addtolength{\hdvert}{.6\unitlength}%	% total extra clearance of .3em...
\addtolength{\hddrop}{.3\unitlength}}	% ...half of which is at bottom

\newcommand{\pregzero}[1]{%
\begin{picture}(0.8,0.4)(-0.4,-0.4)
\cell{0}{-0.4}{t}{#1}%
\cell{0}{0}{c}{\zmark}%
\end{picture}}

\mcm{\gzero}{1}{%
\ginitdims{0.8}{.6}%
\belowpic{#1}%
\sidespic{#1}%	
\present{\pregzero{#1}}}

\newcommand{\pregfst}[1]{%
\begin{picture}(0.5,0.2)(-0.5,-0.2)%
\cell{-0.1}{-0.2}{tr}{#1}%
\cell{0}{0}{c}{\zmark}%
\end{picture}}

\mcm{\gfst}{1}{%
\ginitdims{0.5}{0.4}%
\sidespic{#1}%
\belowpic{#1}%
\presentr{\pregfst{#1}}}

\newcommand{\preglst}[1]{%
\begin{picture}(0.5,0.2)(0,-0.2)%
\cell{0.1}{-0.2}{tl}{#1}%
\cell{0.05}{0}{c}{\zmark}%
\end{picture}}

\mcm{\glst}{1}{%
\ginitdims{.5}{.4}%
\sidespic{#1}%
\belowpic{#1}%
\presentl{\preglst{#1}}}

\newcommand{\preglft}[1]{%
\begin{picture}(0,0.2)(0,-0.2)%
\cell{-0.1}{-0.2}{tr}{#1}%
\cell{0.05}{0}{c}{\zmark}%
\end{picture}}

\mcm{\glft}{1}{%
\ginitdims{0}{.4}%
\belowpic{#1}%
\present{\preglft{#1}}}

\newcommand{\pregrgt}[1]{%
\begin{picture}(0,0.2)(0,-0.2)%
\cell{0.1}{-0.2}{tl}{#1}%
\cell{0.05}{0}{c}{\zmark}%
\end{picture}}

\mcm{\grgt}{1}{%
\ginitdims{0}{.4}%
\belowpic{#1}%
\present{\pregrgt{#1}}}

\newcommand{\pregblw}[1]{%
\begin{picture}(0,0.3)(0,-0.3)
\cell{0}{-0.3}{t}{#1}%
\cell{0.05}{0}{c}{\zmark}%
\end{picture}}

\mcm{\gblw}{1}{%
\ginitdims{0}{.6}%
\belowpic{#1}%
\present{\pregblw{#1}}}

\newcommand{\pregfbw}[1]{%
\begin{picture}(0,0.65)(0,-0.65)
\cell{0}{-0.65}{t}{#1}%
\cell{0.05}{0}{c}{\zmark}%
\end{picture}}

\mcm{\gfbw}{1}{%
\ginitdims{0}{1.3}%
\belowpic{#1}%
\present{\pregfbw{#1}}}

\newcommand{\pregone}[1]{%
\begin{picture}(5,0.4)(0,-0.2)%
\cell{2.5}{0.2}{b}{#1}%
\put(0,0){\vector(1,0){4.9}}% Shortened length from 5 doing survey
\end{picture}}

\mcm{\gone}{1}{%
\ginitdims{5}{0.4}%
\abovepic{#1}%
\present{\pregone{#1}}}

\newcommand{\pregtwo}[3]{%
\begin{picture}(5,3.4)(0,-0.2)%
\cell{2.5}{3.2}{b}{#1}%
\cell{2.5}{-.2}{t}{#2}%
\cell{2.7}{1.5}{l}{#3}%
\qbezier(0,1.5)(2.5,4.5)(5,1.5)%
\qbezier(0,1.5)(2.5,-1.5)(5,1.5)%
\put(5,1.5){\vector(1,-1){0}}%
\put(5,1.5){\vector(1,1){0}}%
\put(2.5,2.5){\vector(0,-1){2}}%
\end{picture}}

\mcm{\gtwo}{3}{%
\ginitdims{5}{3.4}%
\abovepic{#1}%
\belowpic{#2}%
\present{\pregtwo{#1}{#2}{#3}}}

\newcommand{\pregtwoop}[3]{%
\begin{picture}(5,3.4)(0,-0.2)%
\cell{2.5}{3.2}{b}{#1}%
\cell{2.5}{-.2}{t}{#2}%
\cell{2.7}{1.5}{l}{#3}%
\qbezier(0,1.5)(2.5,4.5)(5,1.5)%
\qbezier(0,1.5)(2.5,-1.5)(5,1.5)%
\put(0,1.5){\vector(-1,-1){0}}%
\put(0,1.5){\vector(-1,1){0}}%
\put(2.5,2.5){\vector(0,-1){2}}%
\end{picture}}

\mcm{\gtwoop}{3}{%
\ginitdims{5}{3.4}%
\abovepic{#1}%
\belowpic{#2}%
\present{\pregtwoop{#1}{#2}{#3}}}

\newcommand{\pregtwoco}[3]{%
\begin{picture}(5,3.4)(0,-0.2)%
\cell{2.5}{3.2}{b}{#1}%
\cell{2.5}{-.2}{t}{#2}%
\cell{2.7}{1.5}{l}{#3}%
\qbezier(0,1.5)(2.5,4.5)(5,1.5)%
\qbezier(0,1.5)(2.5,-1.5)(5,1.5)%
\put(5,1.5){\vector(1,-1){0}}%
\put(5,1.5){\vector(1,1){0}}%
\put(2.5,0.5){\vector(0,1){2}}%
\end{picture}}

\mcm{\gtwoco}{3}{%
\ginitdims{5}{3.4}%
\abovepic{#1}%
\belowpic{#2}%
\present{\pregtwoco{#1}{#2}{#3}}}

\newcommand{\pregtwocoop}[3]{%
\begin{picture}(5,3.4)(0,-0.2)%
\cell{2.5}{3.2}{b}{#1}%
\cell{2.5}{-.2}{t}{#2}%
\cell{2.7}{1.5}{l}{#3}%
\qbezier(0,1.5)(2.5,4.5)(5,1.5)%
\qbezier(0,1.5)(2.5,-1.5)(5,1.5)%
\put(0,1.5){\vector(-1,-1){0}}%
\put(0,1.5){\vector(-1,1){0}}%
\put(2.5,0.5){\vector(0,1){2}}%
\end{picture}}

\mcm{\gtwocoop}{3}{%
\ginitdims{5}{3.4}%
\abovepic{#1}%
\belowpic{#2}%
\present{\pregtwocoop{#1}{#2}{#3}}}

\newcommand{\pregthree}[5]{%
\begin{picture}(5,5.4)(0,-1.2)%
\cell{2.5}{4.2}{b}{#1}%
\cell{1.5}{1.7}{b}{#2}%
\cell{2.5}{-1.2}{t}{#3}%
\cell{2.7}{2.75}{l}{#4}%
\cell{2.7}{0.25}{l}{#5}%
\qbezier(0,1.5)(2.5,6.5)(5,1.5)%
\qbezier(0,1.5)(2.5,-3.5)(5,1.5)%
\put(0,1.5){\vector(1,0){5}}%
\put(2.5,3.5){\vector(0,-1){1.5}}%
\put(2.5,1){\vector(0,-1){1.5}}%
\put(5,1.5){\vector(1,-3){0}}%
\put(5,1.5){\vector(1,3){0}}%
\end{picture}}

\mcm{\gthree}{5}{%
\ginitdims{5}{5.4}%
\abovepic{#1}%
\belowpic{#3}%
\present{\pregthree{#1}{#2}{#3}{#4}{#5}}}

\newcommand{\pregfour}[7]{%
\begin{picture}(5,8.4)(0,-2.7)%
\cell{2.5}{5.7}{b}{#1}%
\cell{1.5}{2.8}{b}{#2}%
\cell{1.5}{0.2}{t}{#3}%
\cell{2.5}{-2.7}{t}{#4}%
\cell{2.7}{4.25}{l}{#5}%
\cell{2.7}{1.5}{l}{#6}%
\cell{2.7}{-1.25}{l}{#7}%
\qbezier(0,1.5)(2.5,9.5)(5,1.5)%
\qbezier(0,1.5)(2.5,4)(5,1.5)%
\qbezier(0,1.5)(2.5,-1)(5,1.5)%
\qbezier(0,1.5)(2.5,-6.5)(5,1.5)%
\put(2.5,5.25){\vector(0,-1){2}}%
\put(2.5,2.5){\vector(0,-1){2}}%
\put(2.5,-0.25){\vector(0,-1){2}}%
\put(5,1.5){\vector(1,-4){0}}%
\put(5,1.5){\vector(4,-3){0}}%
\put(5,1.5){\vector(4,3){0}}%
\put(5,1.5){\vector(1,4){0}}%
\end{picture}}

\mcm{\gfour}{7}{%
\ginitdims{5}{8.4}%
\abovepic{#1}%
\belowpic{#4}%
\present{\pregfour{#1}{#2}{#3}{#4}{#5}{#6}{#7}}}

\newcommand{\pregdots}[6]{%
\begin{picture}(5,8.4)(0,-2.7)%
\cell{2.5}{5.7}{b}{#1}%
\cell{1.5}{2.8}{b}{#2}%
\cell{1.5}{0.2}{t}{#3}%
\cell{2.5}{-2.7}{t}{#4}%
\cell{2.7}{4.25}{l}{#5}%
\cell{2.7}{-1.25}{l}{#6}%
\qbezier(0,1.5)(2.5,9.5)(5,1.5)%
\qbezier(0,1.5)(2.5,4)(5,1.5)%
\qbezier(0,1.5)(2.5,-1)(5,1.5)%
\qbezier(0,1.5)(2.5,-6.5)(5,1.5)%
\put(2.5,5.25){\vector(0,-1){2}}%
\put(2.5,-0.25){\vector(0,-1){2}}%
\cell{2.5}{1.9}{c}{\vdots}%
\put(5,1.5){\vector(1,-4){0}}%
\put(5,1.5){\vector(4,-3){0}}%
\put(5,1.5){\vector(4,3){0}}%
\put(5,1.5){\vector(1,4){0}}%
\end{picture}}

\mcm{\gdots}{6}{%
\ginitdims{5}{8.4}%
\abovepic{#1}%
\belowpic{#4}%
\present{\pregdots{#1}{#2}{#3}{#4}{#5}{#6}}}

\newcommand{\pregthreecell}[5]{%
\begin{picture}(8,5)(-4,-2.5)%
\cell{0}{2.5}{b}{#1}%
\cell{0}{-2.5}{t}{#2}%
\cell{-1.7}{0}{r}{#3}%
\cell{1.7}{0}{l}{#4}%
\cell{0}{0.2}{b}{#5}%
\qbezier(-4,0)(0,4.2)(4,0)%
\qbezier(-4,0)(0,-4.2)(4,0)%
\qbezier(-0.5,1.8)(-2.5,0)(-0.5,-1.8)%
\qbezier(0.5,1.8)(2.5,0)(0.5,-1.8)%
\put(-1,0){\vector(1,0){2}}%
\put(4,0){\vector(1,-1){0}}%
\put(4,0){\vector(1,1){0}}%
\put(-0.5,-1.8){\vector(1,-1){0}}%
\put(0.5,-1.8){\vector(-1,-1){0}}%
\end{picture}}

\mcm{\gthreecell}{5}{%
\ginitdims{8}{5}%
\abovepic{#1}%
\belowpic{#2}%
\present{\pregthreecell{#1}{#2}{#3}{#4}{#5}}}

\mcm{\gzeros}{1}{%
\gsinitdims{0.8}{.6}%
\belowpic{\scriptstyle #1}%
\sidespic{\scriptstyle #1}%	
\present{\pregzero{\scriptstyle #1}}}

\mcm{\gfsts}{1}{%
\gsinitdims{0.5}{0.4}%
\sidespic{\scriptstyle #1}%
\belowpic{\scriptstyle #1}%
\presentr{\pregfst{\scriptstyle #1}}}

\mcm{\glsts}{1}{%
\gsinitdims{.5}{.4}%
\sidespic{\scriptstyle #1}%
\belowpic{\scriptstyle #1}%
\presentl{\preglst{\scriptstyle #1}}}

\mcm{\glfts}{1}{%
\gsinitdims{0}{.4}%
\belowpic{\scriptstyle #1}%
\present{\preglft{\scriptstyle #1}}}

\mcm{\grgts}{1}{%
\gsinitdims{0}{.4}%
\belowpic{\scriptstyle #1}%
\present{\pregrgt{\scriptstyle #1}}}

\mcm{\gblws}{1}{%
\gsinitdims{0}{.6}%
\belowpic{#1}%
\present{\pregblw{\scriptstyle #1}}}

\mcm{\gfbws}{1}{%
\gsinitdims{0}{1.3}%
\belowpic{\scriptstyle #1}%
\present{\pregfbw{\scriptstyle #1}}}

\mcm{\gones}{1}{%
\gsinitdims{5}{0.4}%
\abovepic{\scriptstyle #1}%
\present{\pregone{\scriptstyle #1}}}

\mcm{\gtwos}{3}{%
\gsinitdims{5}{3.4}%
\abovepic{\scriptstyle #1}%
\belowpic{\scriptstyle #2}%
\present{\pregtwo{\scriptstyle #1}{\scriptstyle #2}{\scriptstyle #3}}}

\mcm{\gtwoops}{3}{%
\gsinitdims{5}{3.4}%
\abovepic{\scriptstyle #1}%
\belowpic{\scriptstyle #2}%
\present{\pregtwoop{\scriptstyle #1}{\scriptstyle #2}{\scriptstyle #3}}}

\mcm{\gtwocos}{3}{%
\gsinitdims{5}{3.4}%
\abovepic{\scriptstyle #1}%
\belowpic{\scriptstyle #2}%
\present{\pregtwoco{\scriptstyle #1}{\scriptstyle #2}{\scriptstyle #3}}}

\mcm{\gtwocoops}{3}{%
\gsinitdims{5}{3.4}%
\abovepic{\scriptstyle #1}%
\belowpic{\scriptstyle #2}%
\present{\pregtwocoop{\scriptstyle #1}{\scriptstyle #2}{\scriptstyle #3}}}

\mcm{\gthrees}{5}{%
\gsinitdims{5}{5.4}%
\abovepic{\scriptstyle #1}%
\belowpic{\scriptstyle #3}%
\present{\pregthree{\scriptstyle #1}{\scriptstyle #2}{\scriptstyle #3}{\scriptstyle #4}{\scriptstyle #5}}}

\mcm{\gfours}{7}{%
\gsinitdims{5}{8.4}%
\abovepic{\scriptstyle #1}%
\belowpic{\scriptstyle #4}%
\present{\pregfour{\scriptstyle #1}{\scriptstyle #2}{\scriptstyle
#3}{\scriptstyle #4}{\scriptstyle #5}{\scriptstyle #6}{\scriptstyle #7}}} 

\mcm{\gdotss}{7}{%
\gsinitdims{5}{8.4}%
\abovepic{\scriptstyle #1}%
\belowpic{\scriptstyle #4}%
\present{\pregdots{\scriptstyle #1}{\scriptstyle #2}{\scriptstyle
#3}{\scriptstyle #4}{\scriptstyle #5}{\scriptstyle #6}}}

\mcm{\gzersu}{0}{%
\gsinitdims{0}{.6}%
\present{\pregblw{}}}

\mcm{\gfstsu}{0}{%
\gsinitdims{0.5}{0.4}%
\presentr{\pregfst{}}}

\mcm{\glstsu}{0}{%
\gsinitdims{0.5}{0.4}%
\presentl{\preglst{}}}

\mcm{\gonesu}{0}{%
\gsinitdims{5}{0.4}%
\present{\pregone{}}}

\mcm{\gtwosu}{0}{%
\gsinitdims{5}{3.4}%
\present{\pregtwo{}{}{}}}

\mcm{\gthreesu}{0}{%
\gsinitdims{5}{5.4}%
\present{\pregthree{}{}{}{}{}}}

\mcm{\gfoursu}{0}{%
\gsinitdims{5}{8.4}%
\present{\pregfour{}{}{}{}{}{}{}}}

\mcm{\gdotssu}{0}{%
\gsinitdims{5}{8.4}%
\present{\pregdots{}{}{}{}{}{}}}

\newcommand{\pregtwomult}[3]{%
\begin{picture}(5,3.4)(0,-0.2)%
\cell{2.5}{3.2}{b}{#1}%
\cell{2.5}{-.2}{t}{#2}%
\cell{2.8}{1.5}{l}{#3}%
\qbezier(0,1.5)(2.5,4.5)(5,1.5)%
\qbezier(0,1.5)(2.5,-1.5)(5,1.5)%
\put(5,1.5){\vector(1,-1){0}}%
\put(5,1.5){\vector(1,1){0}}%
\put(2.4,2.5){\line(0,-1){1.8}}%
\put(2.6,2.5){\line(0,-1){1.8}}%
\cell{2.51}{0.4}{b}{\vee}%
\end{picture}}

\mcm{\gtwomult}{3}{%
\ginitdims{5}{3.4}%
\abovepic{#1}%
\belowpic{#2}%
\present{\pregtwomult{#1}{#2}{#3}}}

\newcommand{\pregthreemult}[5]{%
\begin{picture}(5,5.4)(0,-1.2)%
\cell{2.5}{4.2}{b}{#1}%
\cell{1.5}{1.7}{b}{#2}%
\cell{2.5}{-1.2}{t}{#3}%
\cell{2.8}{2.75}{l}{#4}%
\cell{2.8}{0.25}{l}{#5}%
\qbezier(0,1.5)(2.5,6.5)(5,1.5)%
\qbezier(0,1.5)(2.5,-3.5)(5,1.5)%
\put(0,1.5){\vector(1,0){5}}%
\put(2.4,3.5){\line(0,-1){1.3}}%
\put(2.6,3.5){\line(0,-1){1.3}}%
\cell{2.51}{1.9}{b}{\vee}%
\put(2.4,1){\line(0,-1){1.3}}%
\put(2.6,1){\line(0,-1){1.3}}%
\cell{2.51}{-0.6}{b}{\vee}%
\put(5,1.5){\vector(1,-3){0}}%
\put(5,1.5){\vector(1,3){0}}%
\end{picture}}

\mcm{\gthreemult}{5}{%
\ginitdims{5}{5.4}%
\abovepic{#1}%
\belowpic{#3}%
\present{\pregthreemult{#1}{#2}{#3}{#4}{#5}}}

\newcommand{\pregthreecellmult}[5]{%
\begin{picture}(8,5)(-4,-2.5)%
\cell{0}{2.5}{b}{#1}%
\cell{0}{-2.5}{t}{#2}%
\cell{-1.8}{0}{r}{#3}%
\cell{1.8}{0}{l}{#4}%
\cell{0}{0.3}{b}{#5}%
\qbezier(-4,0)(0,4.2)(4,0)%
\qbezier(-4,0)(0,-4.2)(4,0)%
\qbezier(-0.6,1.8)(-2.6,0)(-0.6,-1.8)%
\qbezier(-0.4,1.8)(-2.4,0)(-0.5,-1.7)%
\cell{-0.6}{-1.8}{b}{\lrcorner}%
\qbezier(0.4,1.8)(2.4,0)(0.5,-1.7)%
\qbezier(0.6,1.8)(2.6,0)(0.6,-1.8)%
\cell{0.65}{-1.8}{b}{\llcorner}%
\put(-1,0.15){\line(1,0){1.7}}%
\put(-1,0){\line(1,0){2}}%
\put(-1,-0.15){\line(1,0){1.7}}%
\cell{1.15}{0}{r}{>}%
\put(4,0){\vector(1,-1){0}}%
\put(4,0){\vector(1,1){0}}%
\end{picture}}

\mcm{\gthreecellmult}{5}{%
\ginitdims{8}{5}%
\abovepic{#1}%
\belowpic{#2}%
\present{\pregthreecellmult{#1}{#2}{#3}{#4}{#5}}}

\newcommand{\pregtwomultsu}{%
\begin{picture}(5,3.4)(0,-0.2)%
\qbezier(0,1.5)(2.5,4.5)(5,1.5)%
\qbezier(0,1.5)(2.5,-1.5)(5,1.5)%
\put(5,1.5){\vector(1,-1){0}}%
\put(5,1.5){\vector(1,1){0}}%
\cell{2.5}{1.5}{c}{\Downarrow}%
\end{picture}}

\mcm{\gtwomultsu}{0}{%
\gsinitdims{5}{3.4}%
\present{\pregtwomultsu}}

\newcommand{\pregthreemultsu}{%
\begin{picture}(5,5.4)(0,-1.2)%
\qbezier(0,1.5)(2.5,6.5)(5,1.5)%
\qbezier(0,1.5)(2.5,-3.5)(5,1.5)%
\put(0,1.5){\vector(1,0){5}}%
\cell{2.5}{2.75}{c}{\Downarrow}%
\cell{2.5}{0.25}{c}{\Downarrow}%
\put(5,1.5){\vector(1,-3){0}}%
\put(5,1.5){\vector(1,3){0}}%
\end{picture}}

\mcm{\gthreemultsu}{0}{%
\gsinitdims{5}{5.4}%
\present{\pregthreemultsu}}

\newcommand{\pregfourmultsu}{%
\begin{picture}(5,8.4)(0,-2.7)%
\qbezier(0,1.5)(2.5,9.5)(5,1.5)%
\qbezier(0,1.5)(2.5,4)(5,1.5)%
\qbezier(0,1.5)(2.5,-1)(5,1.5)%
\qbezier(0,1.5)(2.5,-6.5)(5,1.5)%
\cell{2.5}{4.25}{c}{\Downarrow}%
\cell{2.5}{1.5}{c}{\Downarrow}%
\cell{2.5}{-1.25}{c}{\Downarrow}%
\put(5,1.5){\vector(1,-4){0}}%
\put(5,1.5){\vector(4,-3){0}}%
\put(5,1.5){\vector(4,3){0}}%
\put(5,1.5){\vector(1,4){0}}%
\end{picture}}

\mcm{\gfourmultsu}{0}{%
\gsinitdims{5}{8.4}%
\present{\pregfourmultsu}}

\newcommand{\pregthreecellu}{%
\begin{picture}(5,3.4)(-0.5,-0.2)%
\qbezier(-.5,1.5)(2,4.5)(4.5,1.5)%
\qbezier(-.5,1.5)(2,-1.5)(4.5,1.5)%
\qbezier(1.5,2.7)(0.5,1.5)(1.5,0.3)%
\qbezier(2.5,2.7)(3.5,1.5)(2.5,0.3)%
\put(1.3,1.5){\vector(1,0){1.4}}%
\put(4.5,1.5){\vector(1,-1){0}}%
\put(4.5,1.5){\vector(1,1){0}}%
\put(1.5,0.3){\vector(2,-3){0}}%
\put(2.5,0.3){\vector(-2,-3){0}}%
\end{picture}}

\mcm{\gthreecellu}{0}{%
\ginitdims{5}{3.4}%
\present{\pregthreecellu}}

\newcommand{\pregtwocentre}[3]{%
\begin{picture}(5,3.4)(0,-0.2)%
\cell{2.5}{3.2}{b}{#1}%
\cell{2.5}{-.2}{t}{#2}%
\cell{2.5}{1.5}{c}{#3}%
\qbezier(0,1.5)(2.5,4.5)(5,1.5)%
\qbezier(0,1.5)(2.5,-1.5)(5,1.5)%
\put(5,1.5){\vector(1,-1){0}}%
\put(5,1.5){\vector(1,1){0}}%
\put(2.5,2.5){\vector(0,-1){2}}%
\end{picture}}

\mcm{\gtwocentre}{3}{%
\ginitdims{5}{3.4}%
\abovepic{#1}%
\belowpic{#2}%
\present{\pregtwocentre{#1}{#2}{#3}}}

\newcommand{\pregspecialone}[9]{%
\begin{picture}(8,8)(-4,-4)%
\cell{0}{3.9}{b}{#1}%
\cell{-2}{-0.2}{t}{#2}%
\cell{0}{-3.9}{t}{#3}%
\cell{-1.5}{1.1}{r}{#4}%
\cell{0.2}{1.5}{l}{#5}%
\cell{1.5}{1.1}{l}{#6}%
\cell{0.2}{-2}{l}{#7}%
\cell{-0.9}{2.3}{b}{#8}%
\cell{0.9}{2.3}{b}{#9}%
\qbezier(-4,0)(0,8)(4,0)%
\qbezier(-4,0)(0,-8)(4,0)%
\qbezier(-0.5,3.4)(-3.5,2)(-0.5,0.6)%
\qbezier(0.5,3.4)(3.5,2)(0.5,0.6)%
\put(-4,0){\vector(1,0){8}}%
\put(0,3.4){\vector(0,-1){2.8}}%
\put(0,-0.8){\vector(0,-1){2.4}}%
\put(-1.5,2.2){\vector(1,0){1.2}}%
\put(0.3,2.2){\vector(1,0){1.2}}%
\put(4,0){\vector(1,-2){0}}%
\put(4,0){\vector(1,2){0}}%
\put(-0.5,0.6){\vector(2,-1){0}}%
\put(0.5,0.6){\vector(-2,-1){0}}%
\end{picture}}

\mcm{\gspecialone}{9}{%
\ginitdims{8}{8}%
\abovepic{#1}%
\belowpic{#3}%
\present{\pregspecialone{#1}{#2}{#3}{#4}{#5}{#6}{#7}{#8}{#9}}}

\newcommand{\pregspecialtwo}{%
\begin{picture}(5,3.4)(0,-0.2)%
\qbezier(0,1.5)(2.5,4.5)(5,1.5)%
\qbezier(0,1.5)(2.5,-1.5)(5,1.5)%
\qbezier(1.7,2.5)(0,1.5)(1.7,0.5)%
\qbezier(3.3,2.5)(5,1.5)(3.3,0.5)%
\put(5,1.5){\vector(1,-1){0}}%
\put(5,1.5){\vector(1,1){0}}%
\put(1.7,0.5){\vector(3,-2){0}}%
\put(3.3,0.5){\vector(-3,-2){0}}%
\put(2.5,2.5){\vector(0,-1){2}}%
\put(1.2,1.5){\vector(1,0){1}}%
\put(2.8,1.5){\vector(1,0){1}}%
\end{picture}}

\mcm{\gspecialtwo}{0}{%
\ginitdims{5}{3.4}%
\present{\pregspecialtwo}}

\newcommand{\pregspecialthree}{%
\begin{picture}(5,5.4)(0,-1.2)%
\qbezier(0,1.5)(2.5,6.5)(5,1.5)%
\qbezier(0,1.5)(2.5,-3.5)(5,1.5)%
\qbezier(2,3.5)(1,2.75)(2,2)%
\qbezier(3,3.5)(4,2.75)(3,2)%
\qbezier(2,1)(1,0.25)(2,-0.5)%
\qbezier(3,1)(4,0.25)(3,-0.5)%
\put(0,1.5){\vector(1,0){5}}%
\put(2,2.75){\vector(1,0){1}}%
\put(2,0.25){\vector(1,0){1}}%
\put(5,1.5){\vector(1,-3){0}}%
\put(5,1.5){\vector(1,3){0}}%
\put(2,2){\vector(1,-1){0}}%
\put(3,2){\vector(-1,-1){0}}%
\put(2,-0.5){\vector(1,-1){0}}%
\put(3,-0.5){\vector(-1,-1){0}}%
\end{picture}}

\mcm{\gspecialthree}{0}{%
\ginitdims{5}{5.4}%
\present{\pregspecialthree}}

\newcommand{\pregonew}[1]{%
\begin{picture}(8,0.4)(0,-0.2)%
\cell{4}{0.2}{b}{#1}%
\put(0,0){\vector(1,0){8}}%
\end{picture}}

\mcm{\gonew}{1}{%
\ginitdims{8}{0.4}%
\abovepic{#1}%
\present{\pregonew{#1}}}

\newcommand{\preghole}{%
\begin{picture}(5,0)(0,0)%
\end{picture}}

\mcm{\ghole}{0}{%
\gsinitdims{5}{0}%
\present{\preghole}}

\newcommand{\pregtwodotty}[3]{%
\begin{picture}(5,3.4)(0,-0.2)%
\cell{2.5}{3.2}{b}{#1}%
\cell{2.5}{-.2}{t}{#2}%
\cell{2.7}{1.5}{l}{#3}%
\qbezier(0,1.5)(2.5,4.5)(5,1.5)%
\qbezier(0,1.5)(2.5,-1.5)(5,1.5)%
\put(5,1.5){\vector(1,-1){0}}%
\put(5,1.5){\vector(1,1){0}}%
\multiput(2.5,2.5)(0,-0.25){7}{\makebox(0,0)[c]{$\cdot$}}
\put(2.5,0.5){\vector(0,-1){0}}%
\end{picture}}

\mcm{\gtwodotty}{3}{%
\ginitdims{5}{3.4}%
\abovepic{#1}%
\belowpic{#2}%
\present{\pregtwodotty{#1}{#2}{#3}}}

\newcommand{\preghappy}[1]{%
\begin{picture}(5,3.4)(0,-0.2)%
\cell{2.5}{-.2}{t}{#1}%
\qbezier(0,1.5)(2.5,-1.5)(5,1.5)%
\put(5,1.5){\vector(1,1){0}}%
\end{picture}}

\mcm{\ghappysu}{0}{%
\gsinitdims{5}{3.4}%
\present{\preghappy{}}}

\newcommand{\pregunhappy}[1]{%
\begin{picture}(5,3.4)(0,-0.2)%
\cell{2.5}{3.2}{b}{#1}%
\qbezier(0,1.5)(2.5,4.5)(5,1.5)%
\put(5,1.5){\vector(1,-1){0}}%
\end{picture}}

\mcm{\gunhappysu}{0}{%
\gsinitdims{5}{3.4}%
\present{\pregunhappy{}{}{}}}

\newcommand{\pregtwowide}[3]{%
\begin{picture}(10,3.4)(0,-0.2)%
\cell{5}{3.2}{b}{#1}%
\cell{5}{-.2}{t}{#2}%
\cell{5.2}{1.5}{l}{#3}%
\qbezier(0,1.5)(5,4.5)(10,1.5)%
\qbezier(0,1.5)(5,-1.5)(10,1.5)%
\put(10,1.5){\vector(2,-1){0}}%
\put(10,1.5){\vector(2,1){0}}%
\put(5,2.5){\vector(0,-1){2}}%
\end{picture}}

\mcm{\gtwowidesu}{0}{%
\gsinitdims{10}{3.4}%
\present{\pregtwowide{}{}{}}}

\newcommand{\cinitdims}[2]{%		% CELLULAR VERSION
\setlength{\unitlength}{1em}%		% unitlength = 1em
\setlength{\hdoffset}{.35\unitlength}%	% cellular offset = .35em
\setlength{\hdwidth}{#1\unitlength}%	% width as specified
\setlength{\hdvert}{#2\unitlength}%	% vert = #2
\setlength{\hddrop}{.5\hdvert}%		% 
\addtolength{\hddrop}{-1\hdoffset}%	% 
\setlength{\hdbaredrop}{1\hddrop}%	% hddrop = drop = half(vert) - offset
\addtolength{\hdvert}{.6\unitlength}%	% total extra clearance of .6em...
\addtolength{\hddrop}{.3\unitlength}}	% ...half of which is at bottom

\newcommand{\hdc}[1]{\scriptstyle #1}	% style for 1-cell labels in
\newcommand{\precone}[1]{%
\begin{picture}(4.2,0.4)(-0.3,-0.2)%
\cell{1.8}{0.2}{b}{\hdc #1}%
\put(0,0){\vector(1,0){3.6}}%
\end{picture}}

\mcm{\cone}{1}{%
\cinitdims{4.2}{0.4}%
\abovepic{\hdc #1}%
\present{\precone{#1}}}

\newcommand{\prectwo}[3]%
{\begin{picture}(4.2,3.4)(-0.1,-0.2)%
\cell{2}{3.2}{b}{\hdc #1}%
\cell{2}{-0.2}{t}{\hdc #2}%
\cell{2.2}{1.5}{l}{#3}%
\qbezier(0,2)(2,4)(4,2)%
\qbezier(0,1)(2,-1)(4,1)%
\put(4,2){\vector(1,-1){0}}%
\put(4,1){\vector(1,1){0}}%
\put(2,2.5){\vector(0,-1){2}}%
\end{picture}}

\mcm{\ctwo}{3}{%
\cinitdims{4.2}{3.4}%
\abovepic{\hdc #1}%
\belowpic{\hdc #2}%
\present{\prectwo{#1}{#2}{#3}}}

\newcommand{\precthree}[5]{%
\begin{picture}(4.2,5.4)(-0.1,-0.2)%
\cell{2}{5.2}{b}{\hdc #1}%
\cell{1}{2.7}{b}{\hdc #2}%
\cell{2}{-.2}{t}{\hdc #3}%
\cell{2.2}{3.75}{l}{#4}%
\cell{2.2}{1.25}{l}{#5}%
\qbezier(0,3)(2,7)(4,3)%
\qbezier(0,2)(2,-2)(4,2)%
\put(0,2.5){\vector(1,0){4}}%
\put(2,4.5){\vector(0,-1){1.5}}%
\put(2,2){\vector(0,-1){1.5}}%
\put(4,3){\vector(1,-3){0}}%
\put(4,2){\vector(1,3){0}}%
\end{picture}}

\mcm{\cthree}{5}{%
\cinitdims{4.2}{5.4}%
\abovepic{\hdc #1}%
\belowpic{\hdc #3}%
\present{\precthree{#1}{#2}{#3}{#4}{#5}}}

\newcommand{\precthreecell}[5]{%
\begin{picture}(8.2,5)(-4.1,-2.5)%
\cell{0}{2.5}{b}{\hdc #1}%
\cell{0}{-2.5}{t}{\hdc #2}%
\cell{-1.7}{0}{r}{\hdc #3}%
\cell{1.7}{0}{l}{\hdc #4}%
\cell{0}{0.2}{b}{#5}%
\qbezier(-4,0.5)(0,4)(4,0.5)%
\qbezier(-4,-0.5)(0,-4)(4,-0.5)%
\qbezier(-0.5,2)(-2.5,0)(-0.5,-2)%
\qbezier(0.5,2)(2.5,0)(0.5,-2)%
\put(-1,0){\vector(1,0){2}}%
\put(4,0.5){\vector(1,-1){0}}%
\put(4,-0.5){\vector(1,1){0}}%
\put(-0.5,-2){\vector(1,-1){0}}%
\put(0.5,-2){\vector(-1,-1){0}}%
\end{picture}}

\mcm{\cthreecell}{5}{%
\cinitdims{8.2}{5}%
\abovepic{\hdc #1}%
\belowpic{\hdc #2}%
\present{\precthreecell{#1}{#2}{#3}{#4}{#5}}}

\newcommand{\prectwomult}[3]%
{\begin{picture}(4.2,3.4)(-0.1,-0.2)%
\cell{2}{3.2}{b}{\hdc #1}%
\cell{2}{-0.2}{t}{\hdc #2}%
\cell{2.3}{1.5}{l}{#3}%
\qbezier(0,2)(2,4)(4,2)%
\qbezier(0,1)(2,-1)(4,1)%
\put(4,2){\vector(1,-1){0}}%
\put(4,1){\vector(1,1){0}}%
\put(1.9,2.5){\line(0,-1){1.8}}%
\put(2.1,2.5){\line(0,-1){1.8}}%
\cell{2.01}{0.4}{b}{\vee}%
\end{picture}}

\mcm{\ctwomult}{3}{%
\cinitdims{4.2}{3.4}%
\abovepic{\hdc #1}%
\belowpic{\hdc #2}%
\present{\prectwomult{#1}{#2}{#3}}}

\newcommand{\prectwomultcoop}[3]%
{\begin{picture}(4.2,3.4)(-0.1,-0.2)%
\cell{2}{3.2}{b}{\hdc #1}%
\cell{2}{-0.2}{t}{\hdc #2}%
\cell{2.3}{1.5}{l}{#3}%
\qbezier(0,2)(2,4)(4,2)%
\qbezier(0,1)(2,-1)(4,1)%
\put(0,2){\vector(-1,-1){0}}%
\put(0,1){\vector(-1,1){0}}%
\put(1.9,2.3){\line(0,-1){1.8}}%
\put(2.1,2.3){\line(0,-1){1.8}}%
\cell{2.01}{2.0}{b}{\wedge}%
\end{picture}}

\mcm{\ctwomultcoop}{3}{%
\cinitdims{4.2}{3.4}%
\abovepic{\hdc #1}%
\belowpic{\hdc #2}%
\present{\prectwomultcoop{#1}{#2}{#3}}}

\newcommand{\precthreemult}[5]{%
\begin{picture}(4.2,5.4)(-0.1,-0.2)%
\cell{2}{5.2}{b}{\hdc #1}%
\cell{1}{2.7}{b}{\hdc #2}%
\cell{2}{-.2}{t}{\hdc #3}%
\cell{2.3}{3.75}{l}{#4}%
\cell{2.3}{1.25}{l}{#5}%
\qbezier(0,3)(2,7)(4,3)%
\qbezier(0,2)(2,-2)(4,2)%
\put(0,2.5){\vector(1,0){4}}%
\put(1.9,4.5){\line(0,-1){1.3}}%
\put(2.1,4.5){\line(0,-1){1.3}}%
\cell{2.01}{2.9}{b}{\vee}%
\put(1.9,2){\line(0,-1){1.3}}%
\put(2.1,2){\line(0,-1){1.3}}%
\cell{2.01}{0.4}{b}{\vee}%
\put(4,3){\vector(1,-3){0}}%
\put(4,2){\vector(1,3){0}}%
\end{picture}}

\mcm{\cthreemult}{5}{%
\cinitdims{4.2}{5.4}%
\abovepic{\hdc #1}%
\belowpic{\hdc #3}%
\present{\precthreemult{#1}{#2}{#3}{#4}{#5}}}

\newcommand{\precthreecellmult}[5]{%
\begin{picture}(8.2,5)(-4.1,-2.5)%
\cell{0}{2.5}{b}{\hdc #1}%
\cell{0}{-2.5}{t}{\hdc #2}%
\cell{-1.8}{0}{r}{\hdc #3}%
\cell{1.8}{0}{l}{\hdc #4}%
\cell{0}{0.3}{b}{#5}%
\qbezier(-4,0.5)(0,4)(4,0.5)%
\qbezier(-4,-0.5)(0,-4)(4,-0.5)%
\qbezier(-0.6,2)(-2.6,0)(-0.6,-2)%
\qbezier(-0.4,2)(-2.4,0)(-0.5,-1.9)%
\cell{-0.6}{-2}{b}{\lrcorner}%
\qbezier(0.4,2)(2.4,0)(0.5,-1.9)%
\qbezier(0.6,2)(2.6,0)(0.6,-2)%
\cell{0.65}{-2}{b}{\llcorner}%
\put(-1,0.15){\line(1,0){1.7}}%
\put(-1,0){\line(1,0){2}}%
\put(-1,-0.15){\line(1,0){1.7}}%
\cell{1.15}{0}{r}{>}%
\put(4,0.5){\vector(1,-1){0}}%
\put(4,-0.5){\vector(1,1){0}}%
\end{picture}}

\mcm{\cthreecellmult}{5}{%
\cinitdims{8.2}{5}%
\abovepic{\hdc #1}%
\belowpic{\hdc #2}%
\present{\precthreecellmult{#1}{#2}{#3}{#4}{#5}}}

\newcommand{\presplitcoeqlhs}[3]%
{\begin{picture}(4.4,1.8)(-0.1,-0.2)%
\cell{2}{1.8}{b}{\scriptstyle{#1}}%
\cell{2}{1.2}{t}{\scriptstyle{#2}}%
\cell{2}{-0.2}{t}{\scriptstyle{#3}}%
\qbezier(0,1)(2,-1)(4,1)%
\put(0,1.7){\vector(1,0){4}}%
\put(0,1.3){\vector(1,0){4}}%
\put(0,1){\vector(-1,1){0}}%
\end{picture}}

\mcm{\splitcoeqlhs}{3}{%
\cinitdims{4.4}{3.4}%
\abovepic{#1}%
\belowpic{#3}%
\present{\presplitcoeqlhs{#1}{#2}{#3}}}

\newcommand{\presplitcoeqrhs}[2]%
{\begin{picture}(4.4,1.8)(-0.1,-0.2)%
\cell{2}{1.7}{b}{\scriptstyle{#1}}%
\cell{2}{-0.2}{t}{\scriptstyle{#2}}%
\qbezier(0,1)(2,-1)(4,1)%
\put(0,1.5){\vector(1,0){4}}%
\put(0,1){\vector(-1,1){0}}%
\end{picture}}

\mcm{\splitcoeqrhs}{2}{%
\cinitdims{4.4}{3.4}%
\abovepic{#1}%
\belowpic{#2}%
\present{\presplitcoeqrhs{#1}{#2}}}

\newcommand{\prebundleint}[3]%
{\begin{picture}(4.4,1.8)(-0.1,-0.2)%
\cell{2}{1.6}{b}{\scriptstyle{#1}}%
\cell{2}{0.4}{b}{\scriptstyle{#2}}
\cell{2}{-0.35}{t}{\scriptstyle{#3}}%
\put(0,1.5){\vector(1,0){4}}%
\qbezier(0,1.25)(2,-0.75)(4,1.25)%
\qbezier(0,0.85)(2,-1.15)(4,0.85)%
\put(0,1.25){\vector(-1,1){0}}%
\put(0,0.85){\vector(-1,1){0}}%
\end{picture}}

\mcm{\bundleint}{3}{%
\cinitdims{4.4}{3.4}%
\abovepic{#1}%
\belowpic{#3}%
\present{\prebundleint{#1}{#2}{#3}}}

\newcommand{\prectwocentre}[3]%
{\begin{picture}(4.2,3.4)(-0.1,-0.2)%
\cell{2}{3.2}{b}{\hdc #1}%
\cell{2}{-0.2}{t}{\hdc #2}%
\cell{2}{1.5}{c}{#3}%
\qbezier(0,2)(2,4)(4,2)%
\qbezier(0,1)(2,-1)(4,1)%
\put(4,2){\vector(1,-1){0}}%
\put(4,1){\vector(1,1){0}}%
\put(2,2.5){\vector(0,-1){2}}%
\end{picture}}

\mcm{\ctwocentre}{3}{%
\cinitdims{4.2}{3.4}%
\abovepic{\hdc #1}%
\belowpic{\hdc #2}%
\present{\prectwocentre{#1}{#2}{#3}}}

\newcommand{\prectwodotty}[3]%
{\begin{picture}(4.2,3.4)(-0.1,-0.2)%
\cell{2}{3.2}{b}{\hdc #1}%
\cell{2}{-0.2}{t}{\hdc #2}%
\cell{2.2}{1.5}{l}{#3}%
\qbezier(0,2)(2,4)(4,2)%
\qbezier(0,1)(2,-1)(4,1)%
\put(4,2){\vector(1,-1){0}}%
\put(4,1){\vector(1,1){0}}%
\multiput(2,2.5)(0,-0.25){7}{\makebox(0,0)[c]{$\cdot$}}
\put(2,0.5){\vector(0,-1){0}}%
\end{picture}}

\mcm{\ctwodotty}{3}{%
\cinitdims{4.2}{3.4}%
\abovepic{\hdc #1}%
\belowpic{\hdc #2}%
\present{\prectwodotty{#1}{#2}{#3}}}

\newcommand{\prectwoop}[3]%
{\begin{picture}(4.2,3.4)(-0.1,-0.2)%
\cell{2}{3.2}{b}{\hdc #1}%
\cell{2}{-0.2}{t}{\hdc #2}%
\cell{2.2}{1.5}{l}{#3}%
\qbezier(0,2)(2,4)(4,2)%
\qbezier(0,1)(2,-1)(4,1)%
\put(0,2){\vector(-1,-1){0}}%
\put(0,1){\vector(-1,1){0}}%
\put(2,2.5){\vector(0,-1){2}}%
\end{picture}}

\mcm{\ctwoop}{3}{%
\cinitdims{4.2}{3.4}%
\abovepic{\hdc #1}%
\belowpic{\hdc #2}%
\present{\prectwoop{#1}{#2}{#3}}}

\newcommand{\prectwopar}[4]{%
\begin{picture}(4.2,3.4)(-0.1,-0.2)%
\cell{2}{3.2}{b}{\hdc #1}%
\cell{2}{-0.2}{t}{\hdc #2}%
\cell{1.6}{1.5}{r}{#3}%
\cell{2.4}{1.5}{l}{#4}%
\qbezier(0,2)(2,4)(4,2)%
\qbezier(0,1)(2,-1)(4,1)%
\put(4,2){\vector(1,-1){0}}%
\put(4,1){\vector(1,1){0}}%
\put(1.8,2.5){\vector(0,-1){2}}%
\put(2.2,2.5){\vector(0,-1){2}}%
\end{picture}}

\mcm{\ctwopar}{4}{%
\cinitdims{4.2}{3.4}%
\abovepic{\hdc #1}%
\belowpic{\hdc #2}%
\present{\prectwopar{#1}{#2}{#3}{#4}}}

\newcommand{\precthreein}[5]{%
\begin{picture}(4.2,5.4)(-0.1,-0.2)%
\cell{2}{5.2}{b}{\hdc #1}%
\cell{1}{2.7}{b}{\hdc #2}%
\cell{2}{-.2}{t}{\hdc #3}%
\cell{2.2}{3.75}{l}{#4}%
\cell{2.2}{1.25}{l}{#5}%
\qbezier(0,3)(2,7)(4,3)%
\qbezier(0,2)(2,-2)(4,2)%
\put(0,2.5){\vector(1,0){4}}%
\put(2,4.5){\vector(0,-1){1.5}}%
\put(2,0.5){\vector(0,1){1.5}}%
\put(4,3){\vector(1,-3){0}}%
\put(4,2){\vector(1,3){0}}%
\end{picture}}

\mcm{\cthreein}{5}{%
\cinitdims{4.2}{5.4}%
\abovepic{\hdc #1}%
\belowpic{\hdc #3}%
\present{\precthreein{#1}{#2}{#3}{#4}{#5}}}

\newcommand{\precthreecellpar}[6]{%
\begin{picture}(8.2,5)(-4.1,-2.5)%
\cell{0}{2.5}{b}{\hdc #1}%
\cell{0}{-2.5}{t}{\hdc #2}%
\cell{-1.7}{0}{r}{\hdc #3}%
\cell{1.7}{0}{l}{\hdc #4}%
\cell{0}{0.4}{b}{#5}%
\cell{0}{-0.4}{t}{#6}%
\qbezier(-4,0.5)(0,4)(4,0.5)%
\qbezier(-4,-0.5)(0,-4)(4,-0.5)%
\qbezier(-0.5,2)(-2.5,0)(-0.5,-2)%
\qbezier(0.5,2)(2.5,0)(0.5,-2)%
\put(-1,0.2){\vector(1,0){2}}%
\put(-1,-0.2){\vector(1,0){2}}%
\put(4,0.5){\vector(1,-1){0}}%
\put(4,-0.5){\vector(1,1){0}}%
\put(-0.5,-2){\vector(1,-1){0}}%
\put(0.5,-2){\vector(-1,-1){0}}%
\end{picture}}

\mcm{\cthreecellpar}{6}{%
\cinitdims{8.2}{5}%
\abovepic{\hdc #1}%
\belowpic{\hdc #2}%
\present{\precthreecellpar{#1}{#2}{#3}{#4}{#5}{#6}}}

\newcommand{\prectwov}[5]{%
\begin{picture}(3.4,4.2)(0.8,0.9)%
\cell{2.5}{5.1}{b}{#1}%
\cell{2.5}{0.9}{t}{#2}%
\cell{0.8}{3}{r}{\hdc #3}%
\cell{4.2}{3}{l}{\hdc #4}%
\cell{2.5}{3.2}{b}{#5}%
\qbezier(2,5)(0,3)(2,1)%
\qbezier(3,5)(5,3)(3,1)%
\put(2,1){\vector(1,-1){0}}%
\put(3,1){\vector(-1,-1){0}}%
\put(1.5,3){\vector(1,0){2}}%
\end{picture}}

\mcm{\ctwov}{5}{%
\cinitdims{3.4}{4.2}%
\abovepic{#1}%
\belowpic{#2}%
\sidespic{\hdc #3}%
\sidespic{\hdc #4}%
\present{\prectwov{#1}{#2}{#3}{#4}{#5}}}

\newcommand{\precthreecellv}[7]{%
\begin{picture}(5,8.2)(0.5,-1.6)%
\cell{3}{6.6}{b}{#1}%
\cell{3}{-1.6}{t}{#2}%
\cell{0.5}{2.5}{r}{\hdc #3}%
\cell{5.5}{2.5}{l}{\hdc #4}%
\cell{3}{4.2}{b}{\hdc #5}%
\cell{3}{0.8}{t}{\hdc #6}%
\cell{3.2}{2.5}{l}{#7}%
\qbezier(3.5,6.5)(7,2.5)(3.5,-1.5)%
\qbezier(2.5,6.5)(-1,2.5)(2.5,-1.5)%
\put(2.5,-1.5){\vector(1,-1){0}}%
\put(3.5,-1.5){\vector(-1,-1){0}}%
\qbezier(1,3)(3,5)(5,3)%
\qbezier(1,2)(3,0)(5,2)%
\put(5,3){\vector(1,-1){0}}%
\put(5,2){\vector(1,1){0}}%
\put(3,3.5){\vector(0,-1){2}}%
\end{picture}}

\mcm{\cthreecellv}{7}{%
\cinitdims{5}{8.2}%
\abovepic{#1}%
\belowpic{#2}%
\sidespic{\hdc #3}%
\sidespic{\hdc #4}%
\present{\precthreecellv{#1}{#2}{#3}{#4}{#5}{#6}{#7}}}

\newcommand{\oinitdims}[2]{%		% OPETOPIC VERSION
\setlength{\unitlength}{1em}%		% unitlength = 1em
\setlength{\hdoffset}{.25\unitlength}%	% globular offset = .25em
\setlength{\hdwidth}{#1\unitlength}%	% width as specified
\setlength{\hdvert}{#2\unitlength}%	% vert = #2
\setlength{\hddrop}{.5\hdvert}%		% 
\addtolength{\hddrop}{-1\hdoffset}%	% 
\setlength{\hdbaredrop}{1\hddrop}%	% hddrop = drop = half(vert) - offset
\addtolength{\hdvert}{.6\unitlength}%	% total extra clearance of .6em...
\addtolength{\hddrop}{.3\unitlength}}	% ...half of which is at bottom

\newcommand{\pretopez}[2]{%
\begin{picture}(2.6,2.3)(-1.3,-2.2)% 
\cell{0}{-2.2}{t}{#1}%
\cell{0}{-1.2}{c}{#2}%
\qbezier(0,0)(-2,-2)(0,-2)%
\qbezier(0,0)(2,-2)(0,-2)%
\put(0,0){\vector(-1,1){0}}%
\end{picture}}

\mcm{\topez}{2}{%
\oinitdims{2.6}{2.3}% 
\belowpic{#1}%
\present{\pretopez{#1}{#2}}}

\newcommand{\pretopea}[3]{%
\begin{picture}(4,1.9)(-2,-0,2)%
\cell{0}{1.7}{b}{#1}%
\cell{0}{-0.2}{t}{#2}%
\cell{0}{0.7}{c}{#3}%
\qbezier(-2,0)(0,3)(2,0)%
\put(-2,0){\vector(1,0){4}}%
\put(2,0){\vector(2,-3){0}}%
\end{picture}}

\mcm{\topea}{3}{%
\oinitdims{4}{1.9}%
\abovepic{#1}%
\belowpic{#2}%
\present{\pretopea{#1}{#2}{#3}}}

\newcommand{\pretopeb}[4]{%
\begin{picture}(4,2.2)(-2,-0.2)%
\cell{-1.1}{1}{br}{#1}%
\cell{1.1}{1}{bl}{#2}%
\cell{0}{-0.2}{t}{#3}%
\cell{0}{0.8}{c}{#4}%
\put(-2,0){\vector(1,1){2}}%
\put(0,2){\vector(1,-1){2}}%
\put(-2,0){\vector(1,0){4}}%
\end{picture}}

\mcm{\topeb}{4}{%
\oinitdims{4}{2.2}%
\belowpic{#3}%
\present{\pretopeb{#1}{#2}{#3}{#4}}}

\newcommand{\pretopec}[5]{%
\begin{picture}(4,2.2)(-2,-0.2)%
\cell{-1.8}{1}{br}{#1}%
\cell{0}{2.2}{b}{#2}%
\cell{1.8}{1}{bl}{#3}%
\cell{0}{-0.2}{t}{#4}%
\cell{0}{0.8}{c}{#5}%
\put(-2,0){\vector(1,2){1}}%
\put(-1,2){\vector(1,0){2}}%
\put(1,2){\vector(1,-2){1}}%
\put(-2,0){\vector(1,0){4}}%
\end{picture}}

\mcm{\topec}{5}{%
\oinitdims{4}{2.2}%
\sidespic{#1}%
\abovepic{#2}%
\sidespic{#3}%
\belowpic{#4}%
\present{\pretopec{#1}{#2}{#3}{#4}{#5}}}

\newcommand{\pretoped}[6]{%
\begin{picture}(4,2.5)(-2,-0.2)%
\cell{-2}{0.6}{br}{#1}%
\cell{-0.7}{2.2}{br}{#2}%
\cell{0.7}{2.2}{bl}{#3}%
\cell{2}{0.6}{bl}{#4}%
\cell{0}{-0.2}{t}{#5}%
\cell{0}{0.8}{c}{#6}%
\put(-2,0){\vector(1,3){0.5}}%
\put(-1.5,1.5){\vector(3,2){1.5}}%
\put(0,2.5){\vector(3,-2){1.5}}%
\put(1.5,1.5){\vector(1,-3){0.5}}%
\put(-2,0){\vector(1,0){4}}%
\end{picture}}

\mcm{\toped}{6}{%
\oinitdims{4}{2.5}%
\sidespic{#1}%
\abovepic{#2}%
\abovepic{#3}%
\sidespic{#4}%
\belowpic{#5}%
\present{\pretoped{#1}{#2}{#3}{#4}{#5}{#6}}}

\newcommand{\pretopeeu}[1]{%		% This one's unlabelled except for
\begin{picture}(4,2.8)(-2,-0.2)%
\cell{0}{1}{c}{#1}%
\put(-2,0){\vector(1,4){0.4}}%
\put(-1.6,1.6){\vector(1,1){1}}%
\put(-0.6,2.6){\vector(1,0){1.2}}%
\put(0.6,2.6){\vector(1,-1){1}}%
\put(1.6,1.6){\vector(1,-4){0.4}}%
\put(-2,0){\vector(1,0){4}}%
\end{picture}}

\mcm{\topeeu}{1}{%
\oinitdims{4}{2.5}%
\abovepic{}
\belowpic{}
\present{\pretopeeu{#1}}}

\newcommand{\pretopeq}[5]{%
\begin{picture}(4,2.5)(-2,-0.2)%
\cell{-2}{0.6}{br}{#1}%
\cell{-1}{2.2}{br}{#2}%
\cell{2}{0.6}{bl}{#3}%
\cell{0}{-0.2}{t}{#4}%
\cell{0}{0.8}{c}{#5}%
\put(-2,0){\vector(1,3){0.5}}%
\put(-1.5,1.5){\vector(1,1){1}}%
\cell{0.9}{2.3}{c}{\ddots}
\put(1.5,1.5){\vector(1,-3){0.5}}%
\put(-2,0){\vector(1,0){4}}%
\end{picture}}

\mcm{\topeq}{5}{%
\oinitdims{4}{2.5}%
\sidespic{#1}%
\abovepic{#2}%
\sidespic{#3}%
\belowpic{#4}%
\present{\pretopeq{#1}{#2}{#3}{#4}{#5}}}

\newcommand{\pretopebase}[1]{%
\begin{picture}(4,0.4)(0,-0.2)%
\cell{2}{0.2}{b}{#1}%
\put(0,0){\vector(1,0){4}}%
\end{picture}}

\mcm{\topebase}{1}{%
\oinitdims{4}{0.4}%
\abovepic{#1}%
\present{\pretopebase{#1}}}

\newcommand{\pretopeavar}[5]{%
\begin{picture}(4.4,1.9)(-2.2,-0.2)%
\cell{-2.2}{0}{br}{#1}
\cell{2.2}{0}{bl}{#2}
\cell{0}{1.7}{b}{#3}%
\cell{0}{-0.2}{t}{#4}%
\cell{0.2}{0.7}{l}{#5}%
\qbezier(-2,0)(0,3)(2,0)%
\put(-2,0){\vector(1,0){4}}%
\cell{-2}{0}{c}{\zmark}%
\cell{2}{0}{c}{\zmark}%
\put(2,0){\vector(2,-3){0}}%
\put(0,1.2){\vector(0,-1){1}}
\end{picture}}

\mcm{\topeavar}{5}{%
\oinitdims{4}{1.9}%
\sidespic{#1}%
\sidespic{#2}%
\abovepic{#3}%
\belowpic{#4}%
\sidespic{\zmark}%
\present{\pretopeavar{#1}{#2}{#3}{#4}{#5}}}

\newcommand{\pretopezn}[2]{%
\begin{picture}(2.6,2.3)(-1.3,-2.2)% 
\cell{0}{-2.2}{t}{#1}%
\cell{0}{-1.2}{c}{#2}%
\qbezier(0,0)(-2,-2)(0,-2)%
\qbezier(0,0)(2,-2)(0,-2)%
\end{picture}}

\mcm{\topezn}{2}{%
\oinitdims{2.6}{2.3}% 
\belowpic{#1}%
\present{\pretopezn{#1}{#2}}}

\newcommand{\pretopean}[3]{%
\begin{picture}(4,1.9)(-2,-0,2)%
\cell{0}{1.7}{b}{#1}%
\cell{0}{-0.2}{t}{#2}%
\cell{0}{0.7}{c}{#3}%
\qbezier(-2,0)(0,3)(2,0)%
\put(-2,0){\line(1,0){4}}%
\end{picture}}

\mcm{\topean}{3}{%
\oinitdims{4}{1.9}%
\abovepic{#1}%
\belowpic{#2}%
\present{\pretopean{#1}{#2}{#3}}}

\newcommand{\pretopebn}[4]{%
\begin{picture}(4,2.2)(-2,-0.2)%
\cell{-1.1}{1}{br}{#1}%
\cell{1.1}{1}{bl}{#2}%
\cell{0}{-0.2}{t}{#3}%
\cell{0}{0.8}{c}{#4}%
\put(-2,0){\line(1,1){2}}%
\put(0,2){\line(1,-1){2}}%
\put(-2,0){\line(1,0){4}}%
\end{picture}}

\mcm{\topebn}{4}{%
\oinitdims{4}{2.2}%
\belowpic{#3}%
\present{\pretopebn{#1}{#2}{#3}{#4}}}

\newcommand{\pretopecn}[5]{%
\begin{picture}(4,2.2)(-2,-0.2)%
\cell{-1.8}{1}{br}{#1}%
\cell{0}{2.2}{b}{#2}%
\cell{1.8}{1}{bl}{#3}%
\cell{0}{-0.2}{t}{#4}%
\cell{0}{0.8}{c}{#5}%
\put(-2,0){\line(1,2){1}}%
\put(-1,2){\line(1,0){2}}%
\put(1,2){\line(1,-2){1}}%
\put(-2,0){\line(1,0){4}}%
\end{picture}}

\mcm{\topecn}{5}{%
\oinitdims{4}{2.2}%
\sidespic{#1}%
\abovepic{#2}%
\sidespic{#3}%
\belowpic{#4}%
\present{\pretopecn{#1}{#2}{#3}{#4}{#5}}}

\newcommand{\pretopedn}[6]{%
\begin{picture}(4,2.5)(-2,-0.2)%
\cell{-2}{0.6}{br}{#1}%
\cell{-0.7}{2.2}{br}{#2}%
\cell{0.7}{2.2}{bl}{#3}%
\cell{2}{0.6}{bl}{#4}%
\cell{0}{-0.2}{t}{#5}%
\cell{0}{0.8}{c}{#6}%
\put(-2,0){\line(1,3){0.5}}%
\put(-1.5,1.5){\line(3,2){1.5}}%
\put(0,2.5){\line(3,-2){1.5}}%
\put(1.5,1.5){\line(1,-3){0.5}}%
\put(-2,0){\line(1,0){4}}%
\end{picture}}

\mcm{\topedn}{6}{%
\oinitdims{4}{2.5}%
\sidespic{#1}%
\abovepic{#2}%
\abovepic{#3}%
\sidespic{#4}%
\belowpic{#5}%
\present{\pretopedn{#1}{#2}{#3}{#4}{#5}{#6}}}

\newcommand{\pretopeqn}[5]{%
\begin{picture}(4,2.5)(-2,-0.2)%
\cell{-2}{0.6}{br}{#1}%
\cell{-1}{2.2}{br}{#2}%
\cell{2}{0.6}{bl}{#3}%
\cell{0}{-0.2}{t}{#4}%
\cell{0}{0.8}{c}{#5}%
\put(-2,0){\line(1,3){0.5}}%
\put(-1.5,1.5){\line(1,1){1}}%
\cell{0.9}{2.3}{c}{\ddots}
\put(1.5,1.5){\line(1,-3){0.5}}%
\put(-2,0){\line(1,0){4}}%
\end{picture}}

\mcm{\topeqn}{5}{%
\oinitdims{4}{2.5}%
\sidespic{#1}%
\abovepic{#2}%
\sidespic{#3}%
\belowpic{#4}%
\present{\pretopeqn{#1}{#2}{#3}{#4}{#5}}}

\newcommand{\pretopebasen}[1]{%
\begin{picture}(4,0.4)(0,-0.2)%
\cell{2}{0.2}{b}{#1}%
\put(0,0){\line(1,0){4}}%
\end{picture}}

\mcm{\topebasen}{1}{%
\oinitdims{4}{0.4}%
\abovepic{#1}%
\present{\pretopebasen{#1}}}

\newcommand{\tinputabv}[1]{%
\begin{picture}(1,0.4)(0,0)
\cell{0.6}{0.4}{br}{#1}%
\put(0,0){\line(1,0){1}}
\end{picture}}

\newcommand{\tinputlft}[1]{%
\begin{picture}(1,0)(0,0)
\cell{-0.4}{0}{r}{#1}%
\put(0,0){\line(1,0){1}}
\end{picture}}

\newcommand{\tinputslft}[2]{%
\begin{picture}(1.4,3.8)(-0.4,-1.9)
\cell{0}{1.5}{l}{\tinputlft{#1}}%
\cell{0}{-1.5}{l}{\tinputlft{#2}}%
\cell{0.2}{0.3}{c}{\vdots}%
\end{picture}}

\newcommand{\tinputsslft}[3]{%
\begin{picture}(1.4,3.8)(-0.4,-1.9)
\cell{0}{1.5}{l}{\tinputlft{#1}}%
\cell{0}{0.7}{l}{\tinputlft{#2}}%
\cell{0}{-1.5}{l}{\tinputlft{#3}}%
\cell{0.2}{-0.2}{c}{\vdots}%
\end{picture}}

\newcommand{\toutputrgt}[1]{%
\begin{picture}(1,0)(-1,0)
\cell{0.4}{0}{l}{#1}%
\put(0,0){\line(-1,0){1}}%
\end{picture}}

\newcommand{\tinputslftsmall}[2]{%
\begin{picture}(1.4,3.2)(-0.4,-1.6)
\cell{0}{1.1}{l}{\tinputlft{#1}}%
\cell{0}{-1.1}{l}{\tinputlft{#2}}%
\cell{0.2}{0.3}{c}{\vdots}%
\end{picture}}

\newcommand{\tinputvert}[1]{%
\begin{picture}(0,1)(0,-1)
\cell{0}{0.4}{b}{#1}%
\put(0,0){\line(0,-1){1}}
\end{picture}}

\newcommand{\tinputsvert}[2]{%
\begin{picture}(3.8,1)(-1.9,0)
\cell{-1.5}{0}{b}{\tinputvert{#1}}
\cell{1.5}{0}{b}{\tinputvert{#2}}
\cell{0}{0.8}{c}{\cdots}
\end{picture}}

\newcommand{\tinputssmallvert}[2]{%
\begin{picture}(3.2,1)(-1.6,0)
\cell{-1.1}{0}{b}{\tinputvert{#1}}
\cell{1.1}{0}{b}{\tinputvert{#2}}
\cell{0}{0.8}{c}{\cdots}
\end{picture}}

\newcommand{\tinputssvert}[3]{%
\begin{picture}(3.8,1)(-1.9,0)
\cell{-1.5}{0}{b}{\tinputvert{#1}}
\cell{-0.6}{0}{b}{\tinputvert{#2}}
\cell{1.5}{0}{b}{\tinputvert{#3}}
\cell{0.45}{0.6}{c}{\cdots}
\end{picture}}

\newcommand{\toutputvert}[1]{%
\begin{picture}(0,1)(0,0)
\cell{0}{-0.4}{t}{#1}%
\put(0,0){\line(0,1){1}}
\end{picture}}

\newcommand{\tusual}[1]{%
\begin{picture}(4,4)(-2,-2)
\cell{-0.2}{0}{c}{#1}%
\put(-2,-2){\line(0,1){4}}%
\put(-2,2){\line(2,-1){4}}%
\put(2,0){\line(-2,-1){4}}%
\end{picture}}

\newcommand{\tsmall}[1]{%
\begin{picture}(3,3)(-1.5,-1.5)
\cell{-0.2}{0}{c}{#1}%
\put(-1.5,-1.5){\line(0,1){3}}%
\put(-1.5,1.5){\line(2,-1){3}}%
\put(1.5,0){\line(-2,-1){3}}%
\end{picture}}

\newcommand{\tmid}[1]{%
\begin{picture}(6,6)(-3,-3)
\cell{-0.2}{0}{c}{#1}%
\put(-3,-3){\line(0,1){6}}%
\put(-3,3){\line(2,-1){6}}%
\put(3,0){\line(-2,-1){6}}%
\end{picture}}

\newcommand{\tusualvert}[1]{%
\begin{picture}(4,4)(-2,-2)
\cell{0}{0.75}{c}{#1}%
\put(-2,2){\line(1,0){4}}%
\put(-2,2){\line(1,-2){2}}%
\put(2,2){\line(-1,-2){2}}%
\end{picture}}

\newcommand{\tsmallvert}[1]{%
\begin{picture}(3,3)(-1.5,-1.5)
\cell{0}{0.5}{c}{#1}%
\put(-1.5,1.5){\line(1,0){3}}%
\put(-1.5,1.5){\line(1,-2){1.5}}%
\put(1.5,1.5){\line(-1,-2){1.5}}%
\end{picture}}

\newcommand{\twiggleleft}{%
\begin{picture}(0.5,0.8)(-0.5,0)
\qbezier(0,0)(-0.5,0.4)(0,0.8)
\end{picture}}

\newcommand{\twiggleright}{%
\begin{picture}(0.5,0.8)(0,0)
\qbezier(0,0)(0.5,0.4)(0,0.8)
\end{picture}}

\newcommand{\twiggly}[1]{%
\begin{picture}(4.5,4)(0,-2)
\put(4.5,0){\line(-2,-1){4}}
\put(4.5,0){\line(-2,1){4}}
\cell{0.5}{2}{tr}{\twiggleleft}
\cell{0.5}{1.2}{tl}{\twiggleright}
\cell{0.5}{0.4}{tr}{\twiggleleft}
\cell{0.5}{-0.4}{tl}{\twiggleright}
\cell{0.5}{-1.2}{tr}{\twiggleleft}
\cell{2.5}{0}{c}{#1}
\end{picture}}

\newcommand{\tusualdotty}[1]{%
\begin{picture}(4,4)(-2,-2)
\cell{-0.2}{0}{c}{#1}%
\qbezier[20](-2,-2)(-2,0)(-2,2)
\qbezier[24](-2,2)(0,1)(2,0)
\qbezier[24](2,0)(0,-1)(-2,-2)
\end{picture}}

\newcommand{\tnelabel}[1]{%
\begin{picture}(1.2,1.5)(0,0)%
\cell{0}{0}{bl}{\qbezier(0,0)(0.1,1.3)(1.0,1.4)}
\cell{0}{0}{c}{\scriptstyle\bullet}
\cell{1.2}{1.5}{l}{#1}
\end{picture}}

\newcommand{\tselabel}[1]{%
\begin{picture}(1.2,1.5)(0,-1.5)%
\cell{0}{0}{tl}{\qbezier(0,0)(0.1,-1.3)(1.0,-1.4)}
\cell{0}{0}{c}{\scriptstyle\bullet}
\cell{1.2}{-1.5}{l}{#1}
\end{picture}}

\newtheorem{thm}{Theorem}[section]
\newtheorem{propn}[thm]{Proposition}
\newtheorem{lemma}[thm]{Lemma}
\newtheorem{cor}[thm]{Corollary}
\newtheorem{lotsofremarks}[thm]{Remarks}
\newenvironment{remarks}[1]
	{\begin{lotsofremarks}\label{#1}\end{lotsofremarks}\begin{enumerate}}
	{\end{enumerate}}
\newtheorem{claim}[thm]{Claim}

\newtheorem{predefn}[thm]{Definition}
\newenvironment{defn}{\begin{predefn}\upshape}{\end{predefn}}
\newtheorem{preexample}[thm]{Example}
\newenvironment{example}{\begin{preexample}\upshape}{\end{preexample}}

\renewcommand{\theenumi}{\alph{enumi}}

\renewcommand{\qbeziermax}{150}

\usepackage{fancyhdgs}
\begin{document}

\begin{titlepage}
\thispagestyle{empty}

{\centering

\vspace*{20mm}
{\makebox[0em]{\Huge\bf Higher Operads, Higher Categories}}\\
\vspace*{20mm}
{\huge Tom Leinster}\\
\vspace*{12mm}
{\large\it William Hodge Fellow}\\
{\large\it Institut des Hautes \'Etudes Scientifiques}\\
}

\clearpage
\thispagestyle{empty}

\noindent
\copyright\ 2003 Tom Leinster
\medskip\noindent\\
The author hereby asserts his moral right always to be identified as the
author of this work in accordance with the provisions of the UK Copyright,
Designs and Patent Acts~(1988).\\

\bigskip\noindent
This is a preprint version of a book, to appear as
\begin{quote}
Tom Leinster, \textit{Higher Operads, Higher Categories}, London
Mathematical Society Lecture Notes Series, Cambridge University Press, ISBN
0-521-53215-9. 
\end{quote}
Publication is expected around October 2003.  Details will appear at the
Cambridge University Press website, \url{www.cambridge.org}.  You can
pre-order now at online book stores.
\medskip\noindent\\
By special agreement with Cambridge University Press, the final version of
this text will appear on the ArXiv one year after its publication in
traditional form.\\

\bigskip\noindent
I would be grateful to hear of any errors: \url{leinster@ihes.fr}.  A list
will be maintained at \url{www\dt ihes\dt fr\slsh $\sim$leinster}.

\clearpage
\thispagestyle{empty}

\vspace*{30mm}
\centering\large
\itshape
To\\
\ \\
Martin Hyland\\
Wilson Sutherland\\
and\\
Mr Bull\\

\clearpage
\thispagestyle{empty}
\ 

\end{titlepage}

\tableofcontents

\frontmatter

\chapter{Diagram of Interdependence}

\begin{center}
\epsfig{file=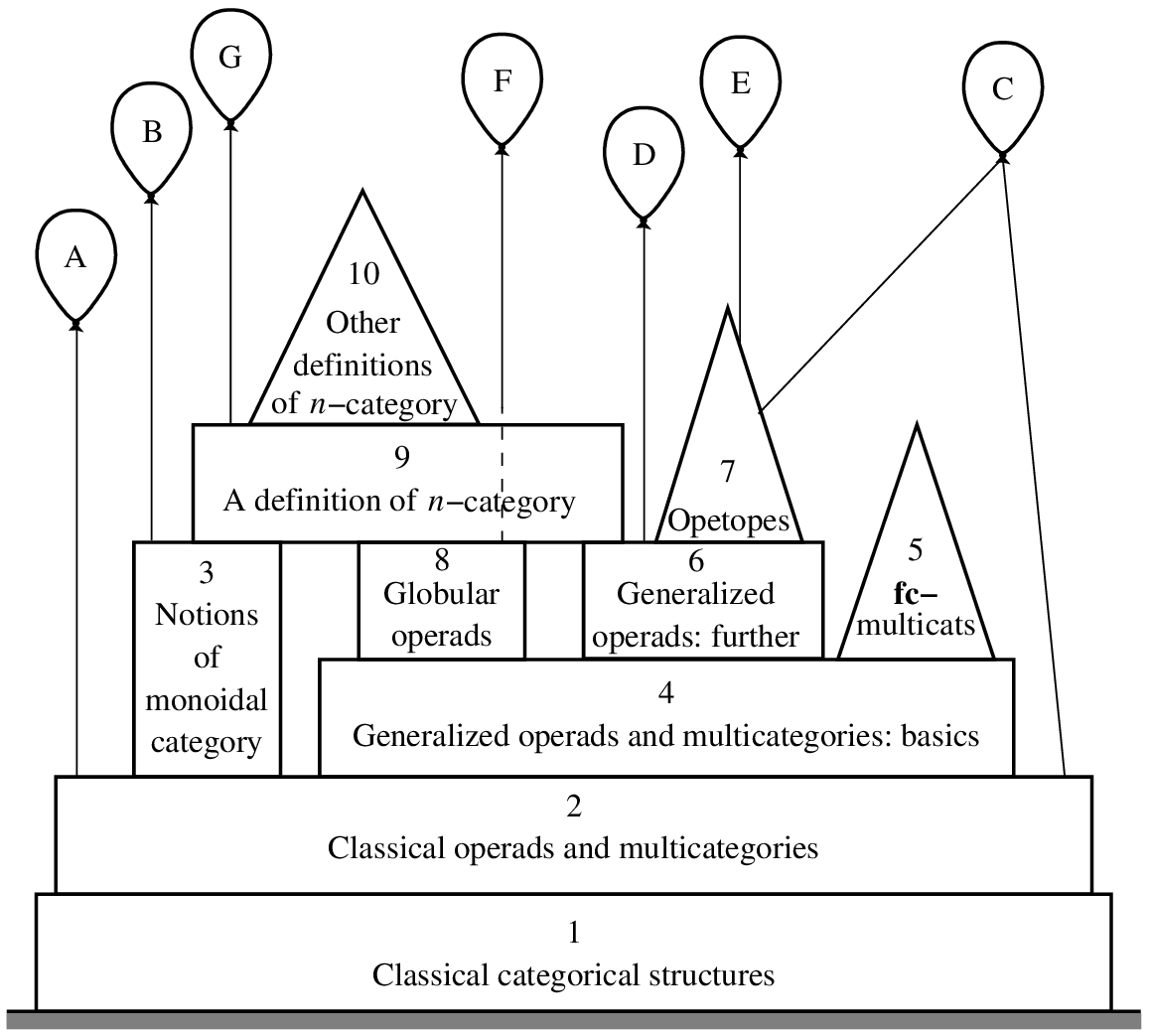}
\end{center}

\vspace*{3ex}
\noindent
Balloons are appendices, giving support in the form of proofs.  They can be
omitted if some results are taken on trust.

\chapter{Acknowledgements}

The ideas presented here have been at the centre of my research for most of
the last six years, during which time I have been influenced---for the
better, I think---by more people than I can thank formally here.  Above all
I feel moulded by the people who were around me during my time in
Cambridge.  I want to thank Martin Hyland in particular: his ways of
thought have been a continual inspiration, and I feel very lucky to have
been his student and to have shared a department with him for so long.  I
have also gained greatly from conversations with Eugenia Cheng, Peter
Johnstone, and Craig Snydal.  Ivan Smith and Dick Thomas have cheerfully
acted as consultants for daft geometry questions.

Ian Grojnowski suggested this Lecture Notes series to me, and I am very
glad he did: I could not have wished for a more open-minded, patient and
helpful editor than Jonathan Walthoe at CUP.  I would also like to thank
Nigel Hitchin, the series editor.

I am very grateful to Andrea Hesketh for sound strategic advice.

Some of the quotations starting the chapters were supplied, directly or
indirectly, by Sean Carmody, David Corfield, and Colin Davey.  Paul-Andr\'e
Melli\`es gave me important information on Swiss cheese.

Almost all of the software used in the preparation of this book was free,
not just financially, but also in the libertarian sense: freedom to take it
apart, alter it, and propagate it, like a piece of mathematics.  I am
grateful to the thousands of developers who brought about the truly
remarkable situation that made this possible.  I also acknowledge the use
of Paul Taylor's excellent commutative diagrams package.

I started writing this when I was Laurence Goddard Research Fellow at St
John's College, Cambridge, and finished when I was William Hodge Fellow at
the Institut des Hautes \'Etudes Scientifiques.  The index was compiled
while I was visiting the University of Chicago at the invitation of Peter
May.  I am immensely grateful to St John's and the IH\'ES for the
opportunities they have given me, and for their unwavering dedication to
research.  In particular, I thank the IH\'ES for giving me the chance to
live in Paris and absorb some Russian culture.

\chapter{Introduction}

\chapterquote{%
It must be admitted that the use of geometric intuition has no logical
necessity in mathematics, and is often left out of the formal presentation
of results.  If one had to construct a mathematical brain, one would
probably use resources more efficiently than creating a visual system.  But
the system is there already, it is used to great advantage by human
mathematicians, and it gives a special flavor to human mathematics.}{%
Ruelle~\cite{Rue}}

\noindent
Higher-dimensional category theory is the study of a zoo of exotic
structures: operads, $n$-categories, multicategories, monoidal categories,
braided monoidal categories, and more.  It is intertwined with the study of
structures such as homotopy algebras ($A_\infty$-categories,
$L_\infty$-algebras, $\Gamma$-spaces, \ldots), $n$-stacks, and $n$-vector
spaces, and draws it inspiration from areas as diverse as topology, quantum
algebra, mathematical physics, logic, and theoretical computer science.

No surprise, then, that the subject has developed chaotically.  The rush
towards formalizing certain commonly-imagined concepts has resulted in an
extraordinary mass of ideas, employing diverse techniques from most of the
subject areas mentioned.  What is needed is a transparent, natural, and
practical language in which to express these ideas.

The main aim of this book is to present one.  It is the language of
generalized operads.  It is introduced carefully, then used to give simple
descriptions of a variety of higher categorical structures.

I hope that by the end, the reader will be convinced that generalized
operads provide as appropriate a language for higher-dimensional category
theory as vector spaces do for linear algebra, or sheaves for algebraic
geometry.  Indeed, the reader may also come to share the feeling that
generalized operads are as applicable and pervasive in mathematics at large
as are $n$-categories, the usual focus of higher-dimensional category
theorists.

Here are some of the structures that we will study, presented informally.

Let $n\in\nat$.  An \demph{$n$-category}%
\index{n-category@$n$-category!definitions of!informal}
consists of \demph{$0$-cells}%
\index{cell!n-category@of $n$-category}
(objects) $a, b, \ldots$, \demph{$1$-cells}
(arrows) $f, g, \ldots$, \demph{$2$-cells} (arrows between arrows) $\alpha,
\beta, \ldots$, \demph{$3$-cells} (arrows between arrows between arrows)
$\Gamma, \Delta, \ldots$, and so on, all the way up to \demph{$n$-cells},
together with various composition operations.  The cells are usually drawn
like this:
\[
\gzero{a},
\ \ 
\gfst{a}\gone{f}\glst{b},
\ \ 
\gfst{a}\gtwomult{f}{g}{\alpha}\glst{b},
\ \ 
\gfst{a}\gthreecellmult{f}{g}{\alpha}{\beta}{\Gamma}\glst{b},
\ \ 
\ldots.
\]
Typical example: for any topological space $X$ there is an $n$-category whose
$k$-cells are maps from the closed $k$-dimensional ball into $X$.  A
$0$-category%
\index{zero-category@0-category}
is just a set, and a $1$-category%
\index{one-category@1-category}
just an ordinary category.

A \demph{multicategory}%
\index{multicategory!informal definition of}
consists of objects $a, b, \ldots$, arrows $\theta,
\phi, \ldots$, a composition operation, and identities, just like an
ordinary category, the difference being that the domain of an arrow is not
just a single object but a finite sequence of them.  An arrow is therefore
drawn as
\[
\begin{centredpic}
\begin{picture}(8,4)(-2,-2)
\cell{0}{0}{l}{\tusual{\theta}}
\cell{0}{0}{r}{\tinputsslft{a_1}{a_2}{a_k}}
\cell{4}{0}{l}{\toutputrgt{a}}
\end{picture}
\end{centredpic}
\]
(where $k\in\nat$), and composition turns a tree of arrows into a single
arrow.  Vector spaces and linear maps form a category; vector spaces and
multilinear maps form a multicategory.

An \demph{operad}%
\index{operad!informal definition of}
is a multicategory with only one object.
Explicitly, an operad consists of a set $P(k)$ for each $k\in\nat$, whose
elements are thought of as `$k$-ary operations' and drawn as
\[
\begin{centredpic}
\begin{picture}(8,4)(-2,-2)
\cell{0}{0}{l}{\tusual{\theta}}
\cell{0}{0}{r}{\tinputsslft{}{}{}}
\cell{4}{0}{l}{\toutputrgt{}}
\end{picture}
\end{centredpic}
\]
with $k$ input wires on the left, together with a rule for composing the
operations and an identity operation.  Example: for any vector space $V$
there is an operad whose $k$-ary operations are the linear maps $V^{\otimes
k} \go V$.

\begin{figure}\centering
\[
\begin{array}{ccc}
\begin{array}{c}
\setlength{\unitlength}{1em}
\begin{picture}(8,4)(-2,-2)
\cell{0}{0}{l}{\tusual{\theta}}
\cell{0}{1.5}{r}{\tinputlft{}}
\cell{0}{0.75}{r}{\tinputlft{}}
\cell{0}{0}{r}{\tinputlft{}}
\cell{0}{-0.75}{r}{\tinputlft{}}
\cell{0}{-1.5}{r}{\tinputlft{}}
\cell{4}{0}{l}{\toutputrgt{}}
\end{picture}
\end{array}
&&
\begin{array}{c}
\epsfig{file=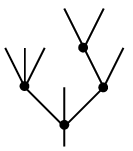}
\end{array}
\goby{\textstyle\theta}
\begin{array}{c}
\epsfig{file=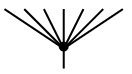}
\end{array}
\\
\\
\begin{array}{c}
\setlength{\unitlength}{1em}
\begin{picture}(5.4,3.6)
\put(0,0){\line(1,0){5.4}}
\put(0,1.8){\line(1,0){5.4}}
\put(0,3.6){\line(1,0){5.4}}
\put(0,0){\line(0,1){3.6}}
\put(1.8,0){\line(0,1){3.6}}
\put(3.6,0){\line(0,1){3.6}}
\put(5.4,0){\line(0,1){3.6}}
\cell{0.9}{0.9}{c}{\Downarrow}
\cell{2.7}{0.9}{c}{\Downarrow}
\cell{4.5}{0.9}{c}{\Downarrow}
\cell{0.9}{2.7}{c}{\Downarrow}
\cell{2.7}{2.7}{c}{\Downarrow}
\cell{4.5}{2.7}{c}{\Downarrow}
\end{picture}
\end{array}
\goby{\textstyle\theta}
\begin{array}{c}
\setlength{\unitlength}{1em}
\begin{picture}(1.8,1.8)
\put(0,0){\line(1,0){1.8}}
\put(0,1.8){\line(1,0){1.8}}
\put(0,0){\line(0,1){1.8}}
\put(1.8,0){\line(0,1){1.8}}
\cell{0.9}{0.9}{c}{\Downarrow}
\end{picture}
\end{array}
&&
\gfstsu\gthreemultsu\gzersu\gfourmultsu\glstsu
\goby{\textstyle\theta}
\gfstsu\gtwomultsu\glstsu
\end{array}
\]
\caption{Operations $\theta$ in four different types of generalized operad}
\label{fig:T-examples}
\end{figure}

Operads describe operations%
\index{operation}%
\index{arrow}
that take a finite sequence of things as input
and produce a single thing as output.  A finite sequence is a 1-dimensional
entity, so operads can be used, for example, to describe the operation of
composing a (1-dimensional) string of arrows in a (1-)category.  But if we
are interested in higher-dimensional structures such as $n$-categories
then we need a more general notion of operad, one where the inputs of an
operation can form a higher-dimensional shape---a grid, perhaps, or a tree,
or a so-called pasting diagram.  For each choice of `input type'%
\index{input type}
$T$, there
is a class of \demph{$T$-operads}.%
\index{generalized operad!informal definition of}
 A $T$-operad consists of a family of
operations whose inputs are of the specified type, together with a rule for
composition; for instance, if the input type $T$ is `finite sequences' then
a $T$-operad is an ordinary operad.  Fig.~\ref{fig:T-examples} shows
typical operations $\theta$ in a $T$-operad, for four different choices of
$T$.
Similarly, there are \demph{$T$-multicategories},%
\index{generalized multicategory!informal definition of}
where the shapes at the
domain and codomain of arrows are labelled with the names of objects.
These are the `generalized operads' and `generalized multicategories'
at the heart of this book.

The uniting feature of all these structures is that they are purely
algebraic in definition, yet near-impossible to understand without drawing
or visualizing pictures.  They are inherently geometrical.  

A notorious problem in this subject is the multiplicity of definitions of
$n$-category.%
\index{n-category@$n$-category!definitions of!comparison}
 Something like a dozen different definitions have been
proposed, and there are still very few precise results stating equivalence
between any of them.  This is not quite the scandal it may seem: it is hard to
say what `equivalence'%
\index{equivalence!definitions of n-category@of definitions of $n$-category}
should even mean.  Suppose that Professors X and Y
each propose a definition of $n$-category.  To compare their definitions,
you find a way of taking one of X's $n$-categories and deriving from it
one of Y's $n$-categories, and vice-versa, then you try to show that
doing one process then the other gets you back to where you started.  It
is, however, highly unrealistic to expect that you will get back to
\emph{exactly} where you started.  For most types of mathematical
structure, getting back to somewhere isomorphic to your starting point
would be a reasonable expectation.  But for $n$-categories, as we shall
see, this is still unrealistic: the canonical notion of equivalence%
\index{equivalence!n-categories@of $n$-categories}
of
$n$-categories is much weaker than isomorphism.  Finding a precise
definition of equivalence for a given definition of $n$-category can be
difficult.  Indeed, many of the proposed definitions of $n$-category did
not come with accompanying proposed definitions of equivalence, and this
gap must be almost certainly be filled before any comparison results can be
proved.

Is this all `just language'?  There would be no shame if it were: language
can have the most profound effect.  New language can make new concepts
thinkable, and make old, apparently obscure, concepts suddenly seem natural
and obvious.  But there is no clear line between mathematical language and
`real' mathematics.  For example, we will see that a 3-category%
\index{three-category@3-category!degenerate}
with only
one 0-cell and one 1-cell is precisely a braided monoidal category,%
\index{monoidal category!braided!degenerate}
and
that the free braided monoidal category on one object is the sequence
$(B_n)_{n\in\nat}$ of braid%
\index{braid}
groups.  So if $n$-categories are just
language, not `real' mathematical objects, then the same is true of the
braid groups, which describe configurations of knotted string.  The
distinction begins to look meaningless.

Here is a summary of the contents.

\subsection*{Motivation for topologists}

Topology and higher-dimensional category theory are intimately related.
The diagrams that one cannot help drawing when thinking about higher
categorical structures can very often be taken literally as pieces of
topology.  We start with an informal discussion of the connections between
the two subjects.  This includes various topological examples of
$n$-categories, and an account of how the world of $n$-categories is a
mirror of the world of homotopy groups of spheres.

\subsection*{Part~\ref{part:background}: Background}

We will build on various `classical' notions.  Those traditionally
considered the domain of category theorists are in
Chapter~\ref{ch:classical}: ordinary categories, bicategories, strict
$n$-categories, and enrichment.  Classical operads and multicategories have
Chapter~\ref{ch:om} to themselves.  They should be viewed as categorical
structures too, although, anomalously, operads are best known to homotopy
theorists and multicategories to categorical logicians.

The familiar concept of monoidal (tensor) category can be formulated in a
remarkable number of different ways.  We look at several in
Chapter~\ref{ch:monoidal}, and prove them equivalent.  Monoidal categories
can be identified with one-object 2-categories, so this is a microcosm of
the comparison of different definitions of $n$-category.

\subsection*{Part~\ref{part:operads}: Operads} 

This introduces the central idea of the text: that
of generalized (`higher') operad and multicategory.  The definitions---of
generalized operad and multicategory, and of algebra for a generalized operad
or multicategory---are stated and explained in
Chapter~\ref{ch:gom-basics}, and some further theory is developed in
Chapter~\ref{ch:gom-further}.  

There is a truly surprising theory of enrichment%
\index{enrichment!generalized multicategory@of generalized multicategory}%
\index{generalized multicategory!enriched}
for generalized
multicategories---it is not at all the routine extension of traditional
enriched category theory that one might expect.  This was to have formed
Part~IV of the book, but for reasons of space it was (reluctantly) dropped.  A
summary of the theory, with pointers to the original papers, is
in Section~\ref{sec:enr-mtis}.  

The rest of Part~\ref{part:operads} is made up of examples and
applications.  Chapter~\ref{ch:fcm} is devoted to so-called
\fc-multicategories,%
\index{fc-multicategory@$\fc$-multicategory}
which are generalized multicategories for a certain choice of input shape.
They turn out to provide a clean setting for some familiar categorical
constructions that have previously been encumbered by technical
restrictions.  In Chapter~\ref{ch:opetopic} we look at opetopic sets,
structures analogous to simplicial sets and used in the definitions of
$n$-category proposed by Baez, Dolan, and others.  Again, the language of
higher operads provides a very clean approach; we also find ourselves drawn
inexorably into higher-dimensional topology.

\subsection*{Part~\ref{part:n-categories}: $n$-Categories}

Using the language of generalized operads, some of the proposed definitions
of $n$-category are very simple to state.  We start by concentrating on one
in particular, in which an $n$-category is defined as an algebra for a
certain globular operad.  A globular operad%
\index{globular operad!informal definition of}
is a $T$-operad for a certain
choice of input type $T$; the associated diagrams are complexes of disks,
as in the last arrow $\theta$ of Fig.~\ref{fig:T-examples}.
Chapter~\ref{ch:globular} explains what globular operads are in pictorial
terms.  In Chapter~\ref{ch:a-defn} we choose a particular globular operad,
define an $n$-category as an algebra for it, and explore the implications
in some depth.  

The many proposed definitions of $n$-category are not as dissimilar as they
might at first appear.  We go through most of them in
Chapter~\ref{ch:other-defns}, drawing together the common threads.

\subsection*{Appendices}

This book is mostly about description: we develop language in which
structures can be described simply and naturally, accurately reflecting
their geometric reality.  In other words, we mostly avoid the convolutions
and combinatorial complexity often associated with higher-dimensional
category theory.  Where things run less smoothly, and in other situations
where a lengthy digression threatens to disrupt the flow of the main text,
the offending material is confined to an appendix.  As long as a few
plausible results are taken on trust, the entire main text can be read and
understood without looking at any of the appendices.

\paragraph*{}

A few words on terminology are needed.  There is a distinction between
`weak' and `strict' $n$-categories,%
\index{n-category@$n$-category!weak vs. strict@weak \vs.\ strict}
as will soon be explained.  For many
years only the strict ones were considered, and they were known simply as
`$n$-categories'.  More recently it came to be appreciated that weak
$n$-categories are much more abundant in nature, and many authors now use
`$n$-category' to mean the weak version.  I would happily join in, but for
the following obstacle: in most parts of this book that concern
$n$-categories, both the weak and the strict versions are involved and
discussed in close proximity.  It therefore seemed preferable to be
absolutely clear and say either `weak $n$-category' or `strict
$n$-category' on every occasion.  The only exceptions are in this
Introduction and the Motivation for Topologists, where the modern
convention is used.

The word `operad'%
\index{operad!usage of word}
will be used in various senses.  The most primitive kind
of operad is an operad of sets without symmetric group action, and this is
our starting point (Chapter~\ref{ch:om}).  We hardly ever consider operads
equipped with symmetric group actions, and when we do we call them
`symmetric operads'; see p.~\pageref{p:sym-warning} for a more
comprehensive warning.

Any finite sequence $x_1, \ldots, x_n$ of elements of a monoid has a
product $x_1 \cdots x_n$.  When $n=0$, this is the the unit element.
Similarly, an identity arrow in a category can be regarded as the composite
of a zero-length string of arrows placed end to end.  I have taken the view
throughout that there is nothing special about units or identities; they
are merely nullary%
\index{nullary!composite}
products or composites.  Related to this is a small but
important convention: the natural%
\index{natural number}
numbers, $\nat$,%
\glo{nat}
start at zero.

\chapter{Motivation for Topologists}

\chapterquote{%
I'm a goddess \emph{and} a nerd!}{%
Bright~\cite{Bri}}

\noindent
Higher-dimensional category theory \emph{can} be treated as a purely
algebraic subject, but that would be missing the point.  It is inherently
topological in nature: the diagrams that one naturally draws to illustrate
higher-dimensional structures can be taken quite literally as pieces of
topology.  Examples of this are the braidings in a braided%
\index{monoidal category!braided}
monoidal
category and the pentagon%
\index{pentagon}
appearing in the definitions of both monoidal
category and $A_\infty$-space.%
\index{A-@$A_\infty$-!space}

This section is an informal description of what higher-dimensional category
theory is and might be, and how it is relevant to topology.  Grothendieck,
for instance, suggested that tame topology should be the study of
$n$-groupoids;%
\index{n-groupoid@$n$-groupoid}
others have hoped that an $n$-category of cobordisms between cobordisms
between \ldots\ will provide a clean setting for topological quantum field
theory; and there is convincing evidence that the whole world of
$n$-categories is a mirror of the world of homotopy groups of spheres.

There are no real theorems, proofs, or definitions here.  But to whet your
appetite, here is a question to which we will reach an answer by the end:

\paragraph*{Question} What is the close connection between the following two
facts? 
\begin{description}
\item[A] No-one ever got into trouble for leaving out the brackets in a
tensor product of several objects (abelian groups, etc.).  For instance, it
is safe to write $A\otimes B\otimes C$ instead of $(A\otimes B)\otimes C$
or $A\otimes (B\otimes C)$.%
\index{associativity}

\item[B] There exist non-trivial knots%
\index{knot}
in $\reals^3$.
\end{description}

\section*{The very rough idea}
\ucontents{section}{The very rough idea}
\index{arrow|(}

In ordinary category theory we have diagrams of objects and arrows such as
\[
\gfstsu\ \gonesu\ \gzersu\ \gonesu\ \gzersu\ \gonesu\ \glstsu\ .
\]
We can imagine more complex category-like structures in which there
are diagrams such as
\[
\epsfig{file=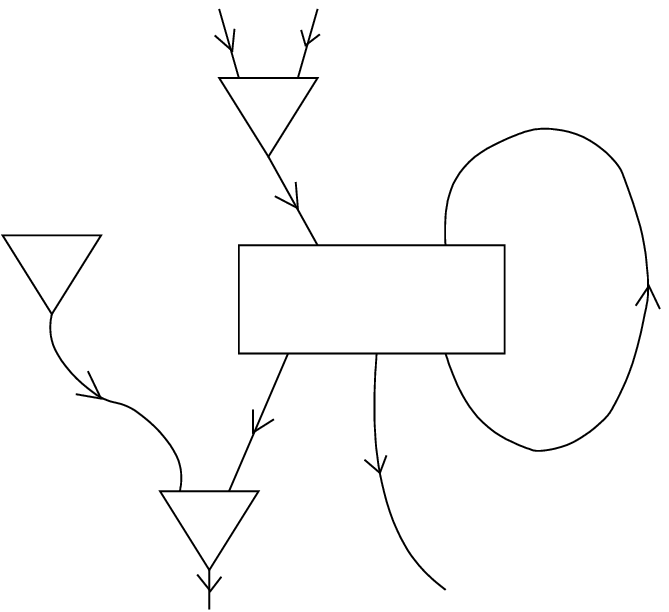}.
\]
This looks like an electronic circuit diagram or a flow chart; the unifying
idea is that of `information flow'. It can be redrawn as
\[
\epsfig{file=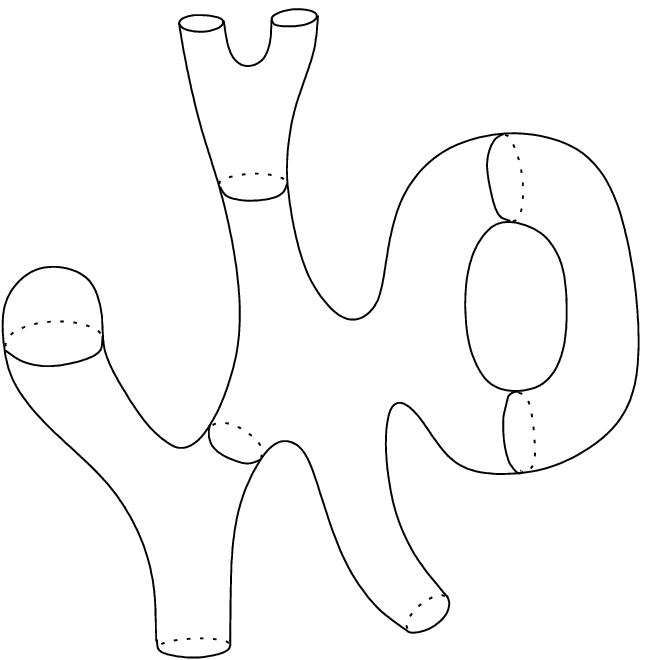,height=62mm},
\]
which looks like a surface or a diagram from topological quantum field
theory.%
\index{topological quantum field theory}

We can also use diagrams like this to express algebraic laws such as
commutativity:
\begin{center}
\setlength{\unitlength}{1mm}
\begin{picture}(55,27)(0,3)
% main article
\cell{5}{6}{bl}{\epsfig{file=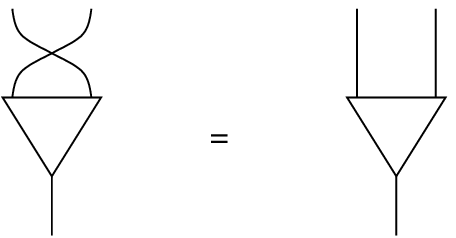}}
% LHS
\cell{5.5}{29.2}{r}{x}
\cell{15}{29}{l}{y}
\cell{5.5}{22}{r}{y}
\cell{15}{22.2}{l}{x}
\cell{10}{5}{t}{y\cdot x}
% RHS
\cell{40.5}{29.2}{r}{x}
\cell{50}{29}{l}{y}
\cell{45}{5}{t}{x\cdot y}
\end{picture}
.
\end{center}
The fact that two-dimensional TQFTs%
\index{topological quantum field theory}
are the same as commutative Frobenius%
\index{Frobenius algebra}%
\index{algebra!Frobenius}
algebras is an example of an explicit link between the spatial and
algebraic aspects of diagrams like these.

Moreover, if we allow crossings, as in the commutativity diagram or as in
\[
\epsfig{file=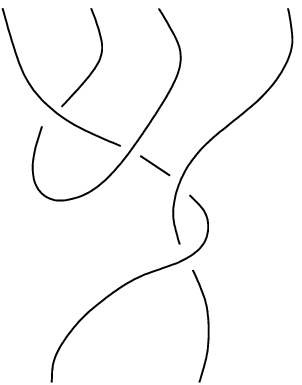,height=25mm,width=50mm},
\]
then we obtain pictures looking like knots; and as we shall see, there are
indeed relations between knot theory and higher categorical structures.

So the idea is:
\begin{trivlist} \item
\fbox{\parbox{0.98\textwidth}{\centering%
Ordinary category theory uses $1$-dimensional arrows $\go$\\
Higher-dimensional category theory uses higher-dimensional arrows}}
\end{trivlist}
The natural topology of these higher-dimensional arrows is what makes
higher-dimensional category theory an inherently topological subject.%
\index{arrow|)}

We will be concerned with structures such as operads, generalized operads
(of which the variety familiar to homotopy theorists is a basic special
case), multicategories, various flavours of monoidal categories, and
$n$-categories; in this introduction I have chosen to concentrate on
$n$-categories.  Terminology: an $n$-category (or `higher-dimensional
category')%
\index{higher-dimensional category@`higher-dimensional category'}
is not a special kind of category, but a generalization of the
notion of category; compare the usage of `quantum group'.  A $1$-category%
\index{one-category@1-category}
is the same thing as an ordinary category, and a $0$-category%
\index{zero-category@0-category}
is just a
set.

\section*{$n$-categories}
\ucontents{section}{$n$-categories}

Here is a very informal

\paragraph*{`Definition'} Let $n\geq 0$.  An
\demph{$n$-category}%
\index{n-category@$n$-category!definitions of!informal}
consists of 
\begin{itemize}
\item \demph{0-cells}%
\index{cell!n-category@of $n$-category}
or \demph{objects}, $A, B, \ldots$
\item \demph{1-cells} or \demph{morphisms}, drawn as
$A \cone{f} B$ 
\item \demph{2-cells} $A \ctwomult{f}{g}{\alpha} B$ (`morphisms between
morphisms') 
\item \demph{3-cells}
$A \cthreecellmult{f}{g}{\alpha}{\beta}{\Gamma} B$ 
(where the arrow
labelled $\Gamma$ is meant to be going in a direction perpendicular to the
plane of the paper)
\item \ldots
\item all the way up to \demph{$n$-cells}
\item various kinds of \demph{composition}, e.g.
\[
\begin{array}{ccc}
A \cone{f} B \cone{g} C					&\textrm{gives}
&A \cone{g\of f} C 			\\
% &&\textrm{(as usual)}			\\
A \cthreemult{f}{g}{h}{\alpha}{\beta} B			&\textrm{gives}
&A \ctwomult{f}{h}{\beta\of\alpha} B	
\label{p:2-cell-comps}
\\
A \ctwomult{f}{g}{\alpha} A' \ctwomult{f'}{g'}{\alpha'}
A'' 							&\textrm{gives}
&A \ctwomult{f'\of f}{g'\of g}{\!\!\!\!\!\!\!\!\alpha' *\alpha} A'',
\end{array}
\]
and so on in higher dimensions; and similarly \demph{identities}.
\end{itemize}
These compositions are required to `all fit together nicely'---a phrase
hiding many subtleties. 
\demph{$\omega$-categories}%
\index{omega-category@$\omega$-category!informal definition of}
(also known as \demph{$\infty$-categories}) are
defined similarly, by going on up the dimensions forever instead of stopping
at $n$.

\paragraph*{} There is nothing forcing us to make the cells spherical
here.  We could, for instance, consider cubical structures, in which
2-cells look like
\begin{diagram}[size=2em,abut]
\bullet	&\rTo		&\bullet\\
\dTo	&\Downarrow	&\dTo	\\
\bullet	&\rTo		&\bullet\makebox[0em]{.}	
\end{diagram}
This is an inhabitant of the `zoo of structures' mentioned earlier, but is
not an $n$-category as such.  (See~\ref{sec:cl-strict} and~\ref{sec:wk-dbl}
for more on this particular structure.)

\paragraph*{Critical Example} Any topological space $X$ gives rise to an
$\omega$-category $\Pi_\omega X$%
\glo{fundomega}
(its \demph{fundamental $\omega$-groupoid}),%
\index{fundamental!omega-groupoid@$\omega$-groupoid}
 in which
\begin{itemize}
\item 0-cells are points of $X$, drawn as $\ \blob$
\item 1-cells are paths%
\index{path}
in $X$ (maps $[0,1] \go X$), drawn as
$\gfstsu\gonesu\glstsu$ ---though whether that is meant to be a picture in
the space $X$ or the $\omega$-category $\Pi_\omega X$ is deliberately
ambiguous; the idea is to blur the distinction between geometry and
algebra
\item 2-cells are homotopies of paths (relative to endpoints), drawn as
$\gfstsu\gtwomultsu\glstsu$
\item 3-cells are homotopies%
\index{homotopy!higher}
of homotopies of paths (that is, suitable maps
$[0,1]^3 \go X$)
\item \ldots
\item composition is by pasting paths and homotopies.
\end{itemize}
(The word `groupoid'%
\index{groupoid}
means that all cells of dimension higher than zero are
invertible.)

$\Pi_\omega X$ should contain all the information you want about $X$ if
your context is `tame topology'.  In particular, you should be able to
compute from it the homotopy, homology and cohomology of $X$.  You can also
truncate after $n$ steps in order to obtain $\Pi_n X$,%
\glo{fundn}
the
\demph{fundamental $n$-groupoid} of $X$; for instance, $\Pi_1 X$ is the
familiar fundamental groupoid.

\paragraph{Alert}%
\index{n-category@$n$-category!weak vs. strict@weak \vs.\ strict|(}
As you may have noticed, composition in $\Pi_\omega X$
is not genuinely associative; nor is it unital, and nor are the cells
genuinely invertible%
\index{invertibility}
(only up to homotopy).  We are therefore interested in
\demph{weak} $n$-categories, where the `fitting together nicely' only
happens up to some kind of equivalence, rather than \demph{strict}
$n$-categories, where associativity and so on hold in the strict sense.

To define strict $n$-categories precisely turns out to be easy.  To define
weak $n$-categories,%
\index{n-category@$n$-category!definitions of}
we face the same kind of challenge as algebraic
topologists did in the 1960s, when they were trying to state the exact
sense in which a loop%
\index{loop space}
space is a topological group.%
\index{group!topological}
 It is clearly not a group in the literal sense, as composition of paths is
not associative; but it is associative up to homotopy, and if you pick
specific homotopies to do this job then these homotopies obey laws of their
own---or at least, obey them up to homotopy; and so on.  At least two
precise formulations of `group up to (higher) homotopy' became popular:
Stasheff's%
\index{Stasheff, Jim}
$A_\infty$-spaces%
\index{A-@$A_\infty$-!space}
and Segal's%
\index{Segal, Graeme}
special $\Delta$-spaces.%
\index{Delta-space@$\Delta$-space}%
\index{special}%
\index{simplicial space}
 (More exactly,
these are notions of monoid or semigroup up to homotopy; the inverses%
\index{invertibility}
are
dealt with separately.)

The situation for weak $n$-categories is similar but more extreme: there
are something like a dozen proposed definitions and, as mentioned in the
Introduction, not much has been proved about how they relate to one
another.  Happily, we can ignore all this here and work informally.  This
means that nothing in the rest of this section is true with any degree of
certainty or accuracy.

At this point you might be thinking: can't we do away with this difficult
theory of weak $n$-categories and just stick to the strict ones?  The
answer is: if you're interested in topology, no.  The difference between
the weak and strict theories is genuine and nontrivial: for while it is
true that every weak $2$-category is equivalent%
\index{coherence!n-categories@for $n$-categories}
to some strict one, and so
it is also true that homotopy $2$-types%
\index{homotopy!type}
can be modelled by strict
$2$-groupoids, neither of these things is true in dimensions $\geq 3$.  For
instance, there exist spaces $X$ (such as the 2-sphere%
\index{sphere}
$S^2$) for which the
weak $3$-category%
\index{fundamental!3-groupoid}
$\Pi_3 X$ is not equivalent%
\index{coherence!tricategories@for tricategories}
to any strict $3$-category.

For the rest of this section, `$n$-category' will mean `weak $n$-category'.
The strict ones are very much the lesser-spotted species.%
\index{n-category@$n$-category!weak vs. strict@weak \vs.\ strict|)}

\paragraph*{Some More Examples} of $\omega$-categories:
\begin{description}
\item[\Top]%
\glo{Top}
This is very similar to the $\Pi_\omega$ example above.  \Top\
has:
\begin{itemize}
\item 0-cells: topological spaces
\item 1-cells: continuous maps
\item 2-cells $X \ctwomult{f}{g}{} Y$: homotopies
between $f$ and $g$
\item 3-cells: homotopies between homotopies (that is, suitable maps
$[0,1]^2 \times X \go Y$)
\item \ldots
\item composition as expected.
\end{itemize}

\item[\fcat{ChCx}]
This $\omega$-category has:
\begin{itemize}
\item 0-cells: chain complexes%
\index{chain complex!omega-category of complexes@$\omega$-category of complexes}
(of abelian groups, say)
\item 1-cells: chain maps
\item 2-cells: chain homotopies%
\index{chain homotopy}
\item 3-cells $A\cthreecellmult{f}{g}{\alpha}{\beta}{\Gamma}B$: homotopies
between homotopies,% 
\lbl{p:ch-hty-hty}\index{homotopy!higher}
that is, maps $\Gamma: A \go B$ of degree $2$ such that
$d\Gamma - \Gamma d = \beta - \alpha$
\item \ldots
\item composition: more or less as expected, but some choices are involved.
For instance, if you try to write down the composite of two chain
homotopies $\gfstsu\gtwomultsu\gzersu\gtwomultsu\glstsu$ then you will find
that there are two equally reasonable ways of doing it: one `left-handed',%
\index{handedness}
one `right-handed'.  This is something like choosing the parametrization
when deciding how to compose two loops in a space (usual choice: do
everything at double speed).  Somehow the fact that there is no canonical
choice means that the resulting $\omega$-category is bound to be weak.
\end{itemize}
In a reasonable world there ought to be a weak $\omega$-functor
$\fcat{Chains}: \fcat{Top} \go \fcat{ChCx}$.

\item[Cobord]%
\index{cobordism|(}\index{manifold!corners@with corners|(}%
\index{topological quantum field theory|(}
This is an $\omega$-category of cobordisms.  
\begin{itemize}
\item 0-cells: 0-manifolds, where `manifold' means `compact, smooth, oriented
manifold'.  A typical 0-cell is 
$\ \stackrel{\uparrow}{\blob}\ \ 
\stackrel{\downarrow}{\blob}\ \ 
\stackrel{\uparrow}{\blob}\ \ 
\stackrel{\uparrow}{\blob}\ $.%
\item 1-cells: 1-manifolds with corners, that is, cobordisms between
0-manifolds, such as
\begin{center}
\setlength{\unitlength}{1mm}
\begin{picture}(30,30)(0,-5)
\cell{0}{0}{bl}{\epsfig{file=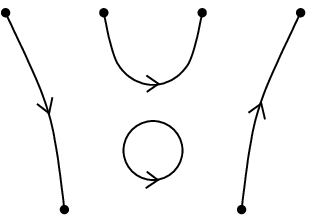}}
\cell{0.5}{21.5}{b}{\scriptstyle{\downarrow}}
\cell{10.5}{21.5}{b}{\scriptstyle{\downarrow}}
\cell{20.5}{21.5}{b}{\scriptstyle{\uparrow}}
\cell{30.5}{21.5}{b}{\scriptstyle{\uparrow}}
\cell{6.5}{-0.5}{t}{\scriptstyle{\downarrow}}
\cell{24.5}{-0.5}{t}{\scriptstyle{\uparrow}}
\end{picture}
\end{center}
(this being a 1-cell from a 4-point 0-manifold to a 2-point 0-manifold).
Atiyah-Segal-style TQFT stops here and takes \emph{isomorphism classes} of
the 1-cells just described, to make a category.  We avoid this (unnatural?)
quotienting out and carry on up the dimensions.
\item 2-cells: 2-manifolds with corners, such as
\[
\begin{array}[c]{c}
\epsfig{file=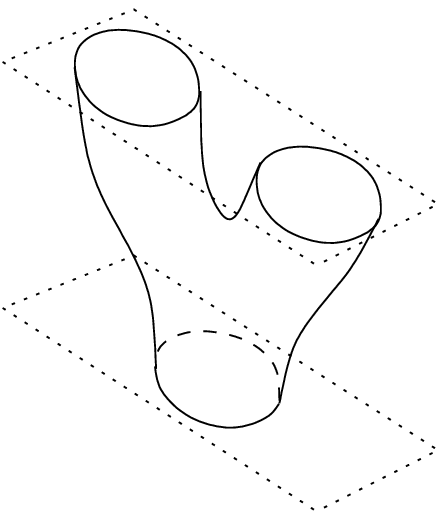}	
\end{array}
\ \ \textrm{or}\ \ 	
\begin{array}[c]{c}
\epsfig{file=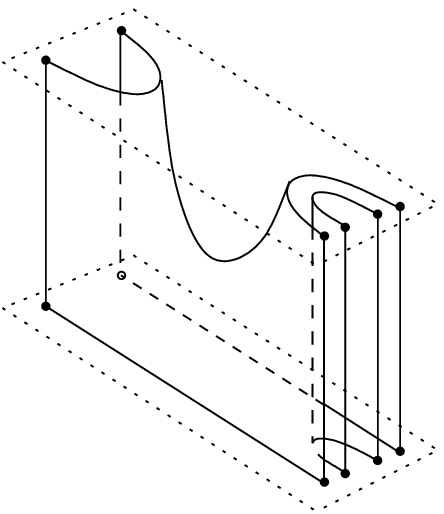}	
\end{array}
\]
(leaving all the orientations off).  The right-hand diagram shows a 2-cell
$L \ctwomult{M}{M'}{N} L'$, where 
\[
L = \ \blob\ \ \blob\ ,
\diagspace 
L'=\ \blob\ \ \blob\ \ \blob\ \ \blob\ , 
\diagspace
M=\begin{array}[c]{c}
\epsfig{file=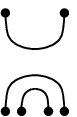}
\end{array},
\diagspace
M'=\begin{array}[c]{c}
\epsfig{file=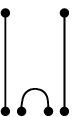}
\end{array},  
\]
and $N$ is the disjoint union of the two 2-dimensional sheets.
Khovanov~\cite{KhoFVI}%
\index{Khovanov, Mikhail}
discusses TQFTs with corners in the language of
2-categories; his approach stops here and take isomorphism classes of the
2-cells just described, to make a 2-category.  Again, we do not quotient
out but keep going up the dimensions.
\item 3-cells, 4-cells, \ldots\ are defined similarly
\item composition is gluing of manifolds.
\end{itemize}
Some authors discuss `extended TQFTs' using the notion of $n$-vector
space.%
\index{n-vector space@$n$-vector space}
A $0$-vector space is a complex number, a $1$-vector space is an ordinary
complex vector space, and $n$-vector spaces for higher $n$ are something
more sophisticated.  See the Notes below for references.%
\index{cobordism|)}\index{manifold!corners@with corners|)}
\index{topological quantum field theory|)}
\end{description}

\section*{Stabilization}
\ucontents{section}{Stabilization}
\index{stabilization|(}\index{n-category@$n$-category!degenerate|(}

So far we have seen that topological structures provide various good examples
of $n$-categories, and that alone might be enough to convince you that
$n$-categories are interesting from a topological point of view.  But the
relationship between topology and higher-dimensional category theory is
actually much more intimate than that.  To see how, we analyse certain
types of degenerate $n$-categories.  It will seem at first as if this is a
purely formal exercise, but before long the intrinsic topology will begin
to shine through.

\paragraph*{Some Degeneracies}

\begin{itemize}

\item A category%
\index{category!degenerate}
\cat{C} with only one object is the same thing as a
monoid%
\lbl{p:degen-cat-monoid}%
\index{monoid!degenerate category@as degenerate category}
(semigroup with unit) $M$.  For if the single object of \cat{C} is
called $\star$, say, then \cat{C} just consists of the set
$\Hom(\star,\star)$ together with a binary operation of
composition and a unit element $1$, obeying the usual axioms.  So we have:
\[
\begin{array}{rcl}
\textrm{morphism in }\cat{C}	&=	&\textrm{element of }M		\\
\of \textrm{ in } \cat{C}	&=	&\cdot \textrm{ in }M.
\end{array}
\ \ \ 
\begin{array}{c}
\setlength{\unitlength}{1em}
\begin{picture}(7.4,2.6)(-3.7,-0.2)
% Object
\cell{0}{0}{c}{\star}
% LHS
\qbezier(-0.5,0)(-2.5,0)(-2.5,1.3)
\qbezier(-2.5,1.3)(-2.5,2.4)(-1.2,2.4)
\qbezier(-1.2,2.4)(-0.4,2.2)(-0.2,0.5)
\put(-0.2,0.5){\vector(0,-1){0}}
\cell{-2.7}{1.4}{r}{x}
% RHS
\qbezier(0.5,0)(2.5,0)(2.5,1.3)
\qbezier(2.5,1.3)(2.5,2.4)(1.2,2.4)
\qbezier(1.2,2.4)(0.4,2.2)(0.2,0.5)
\put(0.5,0.0){\vector(-1,0){0}}
\cell{2.7}{1.4}{l}{y}
\end{picture}
\end{array}
\]

\item A 2-category%
\index{two-category@2-category!degenerate}
\cat{C} with only one $0$-cell is the same thing as a
monoidal category \cat{M}.  (Private thought: if \cat{C} has only one
$0$-cell then there are only interesting things happening in the top two
dimensions, so it must be \emph{some} kind of one-dimensional structure.)
This works as follows:
\begin{eqnarray*}
1\textrm{-cell in }\cat{C}	&=	&\textrm{object of }\cat{M}	\\
2\textrm{-cell in }\cat{C}	&=	&\textrm{morphism of }\cat{M}	\\
\textrm{composition } \ \gzersu\gonesu\gzersu\gonesu\gzersu\  \textrm{ in }
\cat{C} 
&=	
&\otimes \textrm{ of objects in }\cat{M}				\\
\textrm{composition } \ \gzersu\gthreemultsu\gzersu\  \textrm{ in } \cat{C}
&=	
&\of \textrm{ of morphisms in }\cat{M}.				\\
\end{eqnarray*}

\item A monoidal category%
\index{monoidal category!degenerate}
\cat{C} with only one object is\ldots\ well, if
we forget the monoidal structure for a moment then, as we have just seen,
it is a monoid whose elements are the morphisms of \cat{C} and whose
multiplication is the composition in \cat{C}.  Now, the monoidal structure
on \cat{C} provides not only a tensor product for objects, but also a
tensor product for morphisms: so the set of morphisms of \cat{C} has a
second multiplication on it, $\otimes$.  So a one-object monoidal category
is a set $M$ equipped with two monoid structures that are in some sense
compatible (because of the axioms on a monoidal category).  A well-known
result~(\ref{lemma:EH})%
\index{Eckmann--Hilton argument}
says that in this situation, the two
multiplications are in fact equal and commutative.  So, a one-object
monoidal category is a commutative monoid.

This is, essentially, the argument often used to prove that the higher
homotopy groups are abelian, or that the fundamental group of a topological
group is abelian.  In fact, we can deduce that $\pi_2$%
\index{homotopy!group}
is abelian from our
`results' so far:

\textbf{Corollary:} $\pi_2(X,x_0)$ is abelian, for any space $X$ with
basepoint $x_0$. 

\textbf{Proof:} The 2-category $\Pi_2 X$%
\index{fundamental!2-groupoid}
has a sub-2-category whose only
0-cell is $x_0$, whose only 1-cell is the constant path at $x_0$, and whose
2-cells are all the possible ones from $\Pi_2 X$---that is, are the
homotopies from the constant path to itself, that is, are the elements of
$\pi_2(X,x_0)$.  This sub-2-category is a 2-category with only one 0-cell
and one 1-cell, that is, a monoidal category with only one object, that is,
a commutative monoid.

\item Next consider a 3-category%
\index{three-category@3-category!degenerate}
with only one 0-cell and one 1-cell.  We
have not looked at (weak) 3-categories in enough detail to work this out
properly, but it turns out that such a 3-category is the same thing as a
braided monoidal category.  By definition, a \demph{braided%
\index{monoidal category!braided}
monoidal
category} is a monoidal category equipped with a map (a \demph{braiding})
\[
A\otimes B \goby{\beta_{A,B}} B\otimes A
\]
for each pair $(A,B)$ of objects, satisfying axioms \emph{not} including
that
\[
(A\otimes B \goby{\beta_{A,B}} B\otimes A \goby{\beta_{B,A}} A\otimes B)
= 1.
\]
The canonical example of a braided monoidal category (in fact, the braided
monoidal category freely generated by a single object) is \fcat{Braid}.%
\glo{Braid}
 This
has:
\begin{itemize}
\item objects: natural numbers $0, 1, \ldots$
\item morphisms: braids,%
\index{braid}
for instance
\[
\begin{array}[c]{c}
\epsfig{file=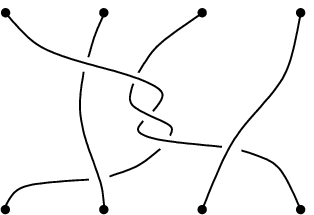}
\end{array}
\diagspace\diagspace\diagspace
\begin{diagram}[height=10mm]
4 \\ \dTo \\ 4 \\
\end{diagram}
\]
(taken up to homotopy); there are no morphisms $m \go n$ when $m \neq n$
\item tensor: placing side-by-side (which on objects means addition)
\item braiding: left over right, for instance
\[
\begin{array}[c]{c}
\epsfig{file=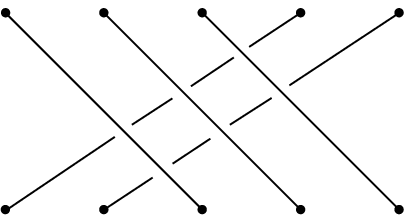}
\end{array}
\diagspace\diagspace
\begin{diagram}[height=10mm]
3+2 \\ \dTo<{\beta_{3,2}} \\ 2+3 \\
\end{diagram}
\]
(Note that $\beta_{n,m} \of \beta_{m,n}$ is not the identity braid.)
\end{itemize}

\item We are rapidly getting out of our depth, but nevertheless: we have
already considered $n$-categories that are only interesting in the top two
dimensions for $n= 1$, $2$, and $3$.  These are categories, monoidal
categories, and braided monoidal categories respectively.  What next?  For
$r\geq 4$, an $r$-category with only one $i$-cell for each $i<r-1$ is,
people believe, the same as a symmetric%
\index{monoidal category!symmetric}
monoidal category (that is, a
braided monoidal category in which $\beta_{B,A} \of \beta_{A,B} = 1$ for
all $A,B$).  So the situation has stabilized\ldots\ and this is meant to
make you start thinking of stabilization phenomena in homotopy.
\end{itemize}

\paragraph*{The Big Picture}  Let us assemble this information on
degeneracies systematically.  Define an \demph{$m$-monoidal $n$-category}%
\index{monoidal n-category@monoidal $n$-category}%
\index{m-monoidal n-category@$m$-monoidal $n$-category}
to be an $(m+n)$-category with only one $i$-cell for each $i<m$.  (This is
certainly some kind of $n$-dimensional structure, as there are only
interesting cells in the top $(n+1)$ dimensions.)  Here is what $m$-monoidal
$n$-categories are for some low values of $m$ and $n$, laid out in the
so-called `periodic%
\index{periodic table}
table'.  Explanation follows.

\begin{center}
% \begin{trivlist} \item
% \hspace*{-10mm}
\setlength{\unitlength}{0.92mm}
\begin{picture}(131,88)(-4,-11)
\place{69}{76}{$n$}
\place{-1}{32}{$m$}
\put(2,66){\line(1,0){125}}
\put(10,0){\line(0,1){72}}
\place{6}{62}{$0$}
\place{6}{54}{$1$}
\place{6}{44}{$2$}
\place{6}{34}{$3$}
\place{6}{24}{$4$}
\place{6}{14}{$5$}
\place{6}{4}{$6$}
\place{24}{70}{$0$}
\place{24}{62}{set}
\place{24}{54}{monoid}
\place{24}{46}{commutative}
\place{24}{42}{monoid}
\place{24}{34}{\ditto}
\place{24}{24}{\ditto}
\place{24}{14}{\ditto}
\place{24}{4}{\ditto}
\place{54}{70}{$1$}
\place{54}{62}{category}
\place{54}{56}{monoidal}
\place{54}{52}{category}
\place{54}{46}{braided}
\place{54}{42}{mon cat}
\place{54}{36}{symmetric}
\place{54}{32}{mon cat}
\place{54}{24}{\ditto}
\place{54}{14}{\ditto}
\place{54}{4}{\ditto}
\place{84}{70}{$2$}
\place{84}{62}{2-category}
\place{84}{56}{monoidal}
\place{84}{52}{2-category}
\place{84}{46}{braided}
\place{84}{42}{mon 2-cat}
\place{84}{34}{rhubarb}
\place{84}{26}{symmetric}
\place{84}{22}{mon 2-cat}
\place{84}{14}{\ditto}
\place{84}{4}{\ditto}
\place{114}{70}{$3$}
\place{114}{62}{3-category}
\place{114}{56}{monoidal}
\place{114}{52}{3-category}
\place{114}{46}{braided}
\place{114}{42}{mon 3-cat}
\place{114}{34}{rhubarb}
\place{114}{24}{rhubarb}
\place{114}{16}{symmetric}
\place{114}{12}{mon 3-cat}
\place{114}{4}{\ditto}
\qbezier[10](24,54)(39,58)(54,62)
\qbezier[20](24,44)(54,53)(84,62)
\qbezier[30](24,34)(69,48)(114,62)
\qbezier[30](24,24)(69,39)(114,54)
\qbezier[30](24,14)(69,29)(114,44)
\qbezier[30](24,4)(69,19)(114,34)
\qbezier[20](54,4)(84,14)(114,24)
\qbezier[10](84,4)(99,9)(114,14)
\qbezier(28,-11)(35.5,-8.5)(43,-6)
\put(28,-11){\vector(-4,-1){0}}
\cell{45}{-9}{l}{\textrm{:\ \ \ take just the one-object structures}}
\label{p:periodic-table}
\end{picture}
% \end{trivlist}
\end{center}

In the first row ($m=0$), a $0$-monoidal $n$-category is simply an
$n$-category: it is not monoidal at all.

In the next row ($m=1$), a $1$-monoidal $n$-category is a monoidal%
\index{monoidal n-category@monoidal $n$-category}
$n$-category, in other words, an $n$-category equipped with a tensor
product that is associative and unital up to equivalence of a suitable
kind.  For instance, a $1$-monoidal $0$-category is a one-object category
(a monoid), and a $1$-monoidal $1$-category is a one-object $2$-category (a
monoidal category).  A monoidal $2$-category can be \emph{defined} as a
one-object 3-category, or can be defined directly as a 2-category with
tensor.

We see from these examples, or from the general definition of $m$-monoidal
$n$-category, that going in the direction $\swarrow$ means restricting to
the one-object structures.

Now look at the third row ($m=2$).  We have already seen that a degenerate
monoidal category is a commutative monoid and a doubly-degenerate
3-category is a braided monoidal category.  It is customary to keep writing
`braided monoidal $n$-category' all along the row, but you can regard this
as nothing more than name-calling.

Next consider the first \emph{column} ($n=0$).  A one-object braided%
\index{monoidal category!braided!degenerate}
monoidal category is going to be a commutative monoid together with a
little extra data (for the braiding) satisfying some axioms, but in fact
this is trivial and we do not get anything new: in some sense, `you can't
get better than a commutative monoid'.  This gives the entry for $m=3,
n=0$, and the same applies all the way down the rest of the column.

A similar story can be told for the second column ($n=1$).  We saw---or
rather, I claimed---that for $m\geq 3$, an $m$-monoidal $1$-category is
just a symmetric monoidal category.  So again the column stabilizes, and
again the point of stabilization is `the most symmetric thing possible'.

The same goes in subsequent columns.  The `rhubarbs' could be replaced by
more terminology---for instance, the first would become `sylleptic%
\index{monoidal category!sylleptic}
monoidal
2-category'---but the details are not important here.

The main point is that the table stabilizes for $m\geq n+2$---just like
$\pi_{m+n}(S^m)$.  So if you overlaid a table of the homotopy groups%
\index{homotopy!group!sphere@of sphere}
of
spheres onto the table above then they would stabilize at the same points.
There are arguments to see why this should be so (and I remind you that
this is all very informal and by no means completely understood).  Roughly,
the fact that the archetypal braided monoidal category \fcat{Braid} is not
symmetric comes down to the fact that you cannot usually translate two
1-dimensional affine subspaces of 3-dimensional space past each other, and
this is the same kind of dimensional calculation as you make when proving
that the homotopy groups of spheres stabilize.%
\index{stabilization|)}

\paragraph{Answer} to the initial question:

\begin{description}
\item[A] Every weak 2-category is equivalent%
\index{coherence!bicategories@for bicategories}
to a strict one.
In particular, every (weak) monoidal category is equivalent to a strict one.
So, for instance, we can pretend that the monoidal category of abelian groups
is strict, making $\otimes$ strictly associative.%
\index{associativity}

\item[B] \emph{Not} every weak 3-category is equivalent%
\index{coherence!tricategories@for tricategories}
to a strict one.
We can construct a counterexample from what we have just done (details
aside).  Facts:
\begin{itemize}
\item a weak 3-category with one $0$-cell and one $1$-cell is a braided%
\index{monoidal category!braided}
monoidal category
\item a strict 3-category with one $0$-cell and one $1$-cell is a strict
symmetric monoidal category
\item any braided monoidal category equivalent to a symmetric monoidal
category is itself symmetric.
\end{itemize}
It follows that any non-symmetric braided monoidal category is a weak
3-category not equivalent to a strict one.  The canonical example is
\fcat{Braid} itself, which is non-symmetric precisely because the overpass
$\begin{array}{c} 
\epsfig{file=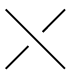}
\end{array}$
cannot be deformed to the underpass 
$\begin{array}{c}
\epsfig{file=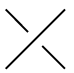}
\end{array}$
in $\reals^3$.
\end{description}%
\index{n-category@$n$-category!degenerate|)}

\begin{notes}

Much has been written on the various interfaces between topology and higher
category theory.  I will just mention a few texts that I happen to have
come across.

Grothendieck%
\index{Grothendieck, Alexander}
puts the case that tame topology is really the study of
$\omega$-groupoids%
\index{omega-groupoid@$\omega$-groupoid}
in his epic~\cite{GroPS} letter to Quillen.  A seminal paper of
 Shum~\cite{Shum}%
\index{Shum, Mei Chee}
establishes connections between higher
categorical structures and knot%
\index{knot}
theory; see also Yetter~\cite{Yet}.

Another introduction to higher categories from a topological viewpoint,
with many similar themes to this one, is the first half of
Baez~\cite{BzIN}.%
\index{Baez, John}

Specifically 2-categorical approaches to topological quantum field theory%
\index{topological quantum field theory}
can be found in Tillmann~\cite{Till}%
\index{Tillmann, Ulrike}
and Khovanov~\cite{KhoFVI}.%
\index{Khovanov, Mikhail}
$n$-vector spaces%
\index{n-vector space@$n$-vector space}
are explained in Kapranov%
\index{Kapranov, Mikhail}
and Voevodsky~\cite{KV},%
\index{Voevodsky, Vladimir}
and
their possible role in topological field theory is discussed in
Lawrence~\cite{Law}.%
\index{Lawrence, Ruth}

The periodic%
\index{periodic table}
table is an absolutely fundamental object of mathematics, only
discovered quite recently (Baez%
\index{Baez, John}
and Dolan~\cite{BDHDA0}),%
\index{Dolan, James}
although
foreshadowed in the work of Breen%
\index{Breen, Lawrence}
and of Street%
\index{Street, Ross}
and the Australian school.  That the table stabilizes for $m\geq n+2$ is
the `Stabilization%
\index{Stabilization Hypothesis}
Hypothesis'.  To state it precisely one needs to set up all the appropriate
definitions first.  A form of it has been proved by Simpson~\cite{SimOBB},%
\index{Simpson, Carlos}
and a heuristic argument in the semistrict case has been given by
Crans~\cite[3.8]{CraBSS}.%
\index{Crans, Sjoerd}

I have made one economy with the truth: a (weak) monoidal%
\index{monoidal category!degenerate}
category with
only one object is not exactly a commutative monoid,%
\index{Eckmann--Hilton argument}
but rather a
commutative monoid equipped with a distinguished invertible element; see
Leinster~\cite[1.6(vii)]{GECM} for the reason.

One interesting idea not mentioned above is a higher categorical approach
to non-abelian cohomology:%
\index{cohomology}%
\index{homology}
specifically, $n$th cohomology should have
coefficients in an $n$-category.  This is explained in
Street~\cite[Introduction]{StrAOS}.%
\index{Street, Ross}

A serious and, of course, highly recommended survey of the proposed
definitions of weak $n$-category,%
\index{n-category@$n$-category!definitions of}
including ten such definitions, is my
own~\cite{SDN}.

\end{notes}

\mainmatter

\part{Background}	
\label{part:background}

\chapter{Classical Categorical Structures}
\lbl{ch:classical}

\chapterquote{%
We will need to use some very simple notions of category theory, an
esoteric subject noted for its difficulty and irrelevance}{%
Moore and Seiberg~\cite{MoSe}}

\noindent
You might imagine that you would need to be on top of the whole of ordinary
category theory before beginning to attempt the higher-dimensional version.
Happily, this is not the case.  The main prerequisite for this book is
basic categorical language, such as may be found in most introductory texts
on the subject.  Except in the appendices, we will need few actual
theorems.

The purpose of this chapter is to recall some familiar categorical ideas
and to explain some less familiar ones.  Where the boundary lies depends,
of course, on the reader, but very little here is genuinely new.
Section~\ref{sec:cats} is on ordinary, `1-dimensional', category theory,
and is a digest of the concepts that will be used later on.  Impatient
readers will want to skip immediately to~\ref{sec:mon-cats}, monoidal
categories.  This covers the basic concepts and two kinds of coherence
theorem.  \ref{sec:cl-enr} is a short section on categories enriched in
monoidal categories.  We need enrichment in the next
section,~\ref{sec:cl-strict}, on strict $n$-categories and strict
$\omega$-categories.  This sets the scene for later chapters, where we
consider the much more profound and interesting \emph{weak} $n$-categories.
Finally, in~\ref{sec:bicats}, we discuss bicategories, the best-known
notion of weak 2-category, including coherence and their (not completely
straightforward) relation to monoidal categories.

Examples of all these structures are given.  Topological spaces and chain
complexes are, as foreshadowed in the Motivation for Topologists, a
recurring theme.

\section{Categories}
\lbl{sec:cats}

This section is a sketch of the category theory on which the rest of the
text is built.  I have also taken the opportunity to state some notation
and some small results that will eventually be needed.

In later chapters I will assume familiarity with the language of
categories, functors, and natural transformations, and the basics of limits
and adjunctions.  I will also use the basic language of monads (sometimes
still called `triples').  As monads are less well known than the other
concepts, and as they will be central to this text, I have included a short
introduction to them below.  

Given objects $A$ and $B$ of a category $\cat{A}$, I write $\cat{A}(A,B)$%
\glo{homset}\index{hom-set}
for the set of maps (or morphisms, or arrows) from $A$ to $B$, in
preference to the less informative $\Hom(A,B)$.  The opposite%
\index{opposite category}%
\index{category!opposite}
or dual%
\index{dual!category@of category}
of
$\cat{A}$ is written $\cat{A}^\op$.%
\glo{catop}
 Isomorphism%
\index{isomorphism}
between objects in a
category is written $\iso$.%
\glo{iso}
 For any two categories $\cat{A}$ and
$\cat{B}$, there is a category $\ftrcat{\cat{A}}{\cat{B}}$%
\glo{ftrcat}\index{functor!category}
whose objects
are functors from $\cat{A}$ to $\cat{B}$ and whose maps are natural
transformations.

My set-theoretic%
\index{foundations}
morals are lax; I have avoided questions of `size'
whenever possible.  Where the issue is unavoidable, I have used the words
\demph{small}%
\index{small}
and \demph{large}%
\index{large}
to mean `forming a set' and `forming a
proper class' respectively.  Some readers may prefer to re-interpret this
using universes.  A category is \demph{small} if its collection of arrows
is small.

The category of sets is written as $\Set$%
\glo{Set}
and the category of small
categories as $\Cat$;%
\glo{Cat}
occasionally I refer to $\CAT$, the (huge)
category of all categories.  There are functors
\[
\begin{diagram}[width=3em,scriptlabels]
\Cat	&\pile{\lTo^D\\ \rTo~{\mr{ob}}\\ \lTo_I}	&\Set,	\\
\end{diagram}
\glo{disccat}\glo{obcat}\glo{indisccat}
\]
where
\begin{itemize}
\item $D$ sends a set $A$ to the \demph{discrete}%
\index{category!discrete}
category on $A$,
whose object-set is $A$ and all of whose maps are identities
\item $\mr{ob}$ sends a category to its set of objects%
\index{objects functor}
\item $I$ sends a set $A$ to the \demph{indiscrete}%
\lbl{p:indiscrete}\index{category!indiscrete}
category on $A$, whose object-set is $A$ and which has precisely one map $a
\go b$ for each $a, b \in A$; all maps are necessarily isomorphisms.
\end{itemize}

We use \demph{limits}%
\index{limit}
(`inverse limits', `projective limits') and
\demph{colimits}%
\index{colimit}
(`direct limits', `inductive limits').  Binary product%
\index{product}
is
written as $\times$%
\glo{binprod}
and arbitrary product as $\prod$;%
\glo{arbprod}
dually, binary
coproduct%
\index{coproduct}
(sum)%
\index{sum}
is written as $+$%
\glo{bincoprod}
and arbitrary coproduct as $\coprod$.%
\glo{arbcoprod}
Nullary products are terminal%
\index{terminal object}
(final)%
\index{final object}
objects, written as $1$;%
\glo{terminal}
in
particular, $1$ often denotes a one-element set.  The unique map from an
object $A$ of a category to a terminal object $1$ of that category is
written $!: A \go 1$.%
\glo{bang}

We will make particular use of \demph{pullbacks}%
\index{pullback}
(fibred%
\index{fibred product}
products).
Pullback squares are indicated by right-angle marks:%
\glo{rightangle}
\[
\begin{diagram}[size=2em]
P\SEpbk		&\rTo	&B	\\
\dTo		&	&\dTo	\\
A		&\rTo	&C.	
\end{diagram}
\]
We also write $P = A \times_C B$.%
\glo{pbtimes}
In later chapters dozens of elementary
manipulations of diagrams involving pullback squares are left to the
hypothetical conscientious reader; almost all are made easy by the
following invaluable lemma.
\begin{lemma}[Pasting Lemma] \lbl{lemma:pasting}\index{Pasting Lemma}
Take a commutative diagram of shape
\[
\begin{diagram}[size=2em]
\cdot	&\rTo	&\cdot	&\rTo	&\cdot	\\
\dTo	&	&\dTo	&	&\dTo	\\
\cdot	&\rTo	&\cdot	&\rTo	&\cdot	\\
\end{diagram}
\]
in some category, and suppose that the right-hand square is a pullback.
Then the left-hand square is a pullback if and only if the outer rectangle
is a pullback.
\done
\end{lemma}

An \demph{adjunction}%
\index{adjunction}
is a pair of functors $\cat{A} \oppair{F}{G} \cat{B}$
together with an isomorphism
\begin{equation}	\label{eq:adjn}
\cat{B}(FA, B) \goiso \cat{A}(A, GB)
\end{equation}
natural in $A \in \cat{A}$ and $B\in \cat{B}$.  Then $F$ is \demph{left
adjoint} to $G$, $G$ is \demph{right adjoint} to $F$, and I write $F\ladj
G$.%
\glo{ladj}
In most of the examples that we meet, $G$ is a forgetful functor and
$F$ the corresponding free functor.  A typical example is that $\cat{A}$ is
the category of sets, $\cat{B}$ is the category of \demph{monoids}%
\index{monoid}
(sets
equipped with an associative binary multiplication and a two-sided unit),
$G$ forgets the monoid structure, and $F$ sends a set $A$ to the monoid%
\index{monoid!free}
\[
FA = \coprod_{n\in\nat} A^n
\]
of finite sequences of elements of $A$, whose multiplication is
concatenation of sequences. 

Take an adjunction as above, and let $A \in \cat{A}$: then applying
isomorphism~\bref{eq:adjn} to $1_{FA} \in \cat{B}(FA, FA)$ yields a map
$\eta_A: A \go GFA$.  The resulting natural transformation $\eta: 1_\cat{A}
\go GF$ is the \demph{unit}%
\index{unit!adjunction@of adjunction}
of the adjunction.  Dually, there is a
\demph{counit}%
\index{counit}
$\epsln: FG \go 1_\cat{B}$, and the unit and counit satisfy
the so-called triangle%
\index{triangle!identities}
identities (Mac Lane~\cite[IV.1(9)]{MacCWM}).  In fact, an adjunction can
 equivalently be defined as a quadruple $(F, G, \eta, \epsln)$ where
\begin{equation}	\label{eq:adjn-data}
\cat{A} \goby{F} \cat{B},
\diagspace
\cat{B} \goby{G} \cat{A},
\diagspace
1_\cat{A} \goby{\eta} GF, 
\diagspace
FG \goby{\epsln} 1_\cat{B}
\end{equation}
and $\eta$ and $\epsln$ satisfy the triangle identities.

Equivalence of categories can be formulated in several ways.  By
definition, an \demph{equivalence}%
\index{equivalence!categories@of categories}
of categories $\cat{A}$ and $\cat{B}$
consists of functors and natural transformations~\bref{eq:adjn-data} such
that $\eta$ and $\epsln$ are isomorphisms.  An \demph{adjoint%
\index{adjoint equivalence}
equivalence}
is an adjunction $(F, G, \eta, \epsln)$ that is also an equivalence.  A
functor $F: \cat{A} \go \cat{B}$ is \demph{essentially%
\index{essentially surjective on objects}
surjective on
objects} if for all $B \in \cat{B}$ there exists $A \in \cat{A}$ such that
$FA \iso B$.
\begin{propn}	\lbl{propn:eqv-eqv}
The following conditions on a functor $F: \cat{A} \go \cat{B}$ are
equivalent: 
\begin{enumerate}
\item \lbl{item:eqv-eqv-adjt-eqv}
there exist $G$, $\eta$ and $\epsln$ such that $(F, G, \eta, \epsln)$ is an
adjoint equivalence 
\item \lbl{item:eqv-eqv-eqv}
there exist $G$, $\eta$ and $\epsln$ such that $(F, G, \eta, \epsln)$ is an
equivalence 
\item \lbl{item:eqv-eqv-ffe}
$F$ is full, faithful, and essentially surjective on objects.
\end{enumerate}
\end{propn}
\begin{proof}
See Mac Lane~\cite[IV.4.1]{MacCWM}.  \done
\end{proof}
If the conditions of the proposition are satisfied then the functor $F$ is
called an \demph{equivalence}.%
\index{equivalence!categories@of categories}
 If the categories $\cat{A}$ and $\cat{B}$
are equivalent then we write $\cat{A} \eqv \cat{B}$%
\glo{eqvcat}
(in contrast to
$\cat{A} \iso \cat{B}$, which denotes isomorphism).

\index{monad|(}
Monads are a remarkably economical formalization of the notion of
`algebraic%
\index{algebraic theory}
theory', traditionally formalized by universal algebraists in
various rather concrete and inflexible ways.  For example, there is a monad
corresponding to the theory of rings, another monad for the theory of
complex Lie algebras, another for the theory of topological groups, another
for the theory of strict $10$-categories, and so on, as we shall see.

A monad on a category $\cat{A}$ can be defined as a monoid in the monoidal
category $(\ftrcat{\cat{A}}{\cat{A}}, \of, 1_\cat{A})$ of endofunctors on
$\cat{A}$.
% category $(\End(\cat{A}), \of, 1_\cat{A})$ of endofunctors on $\cat{A}$.
Explicitly:
\begin{defn}	\lbl{defn:monad}
A \demph{monad} on a category $\cat{A}$ consists of a functor $T: \cat{A}
\go \cat{A}$ together with natural transformations
\[
\mu: T\of T \go T,
\diagspace
\eta: 1_\cat{A} \go T,
\]
called the \demph{multiplication}%
\index{multiplication of monad}
and \demph{unit}%
\index{unit!monad@of monad}
respectively, such that
the diagrams
\[
\begin{slopeydiag}
	&	&T\of T\of T	&	&	\\
	&\ldTo<{T \mu}&		&\rdTo>{\mu T}&	\\
T\of T	&	&		&	&T\of T	\\
	&\rdTo<\mu&		&\ldTo>\mu&	\\
	&	&T		&	&	\\
\end{slopeydiag}
\diagspace
\begin{slopeydiag}
	&		&T\of 1		\\
	&\ldTo<{T \eta}	&		\\
T\of T	&		&\dTo>\id	\\
	&\rdTo<\mu	&		\\
	&		&T		\\
\end{slopeydiag}
\diagspace
\begin{slopeydiag}
1\of T		&		&	\\
		&\rdTo>{\eta T}	&	\\
\dTo<\id	&		&T\of T	\\
		&\ldTo>\mu	&	\\
T		&		&	\\
\end{slopeydiag}
\]
commute (the \demph{associativity} and \demph{unit} laws).
\end{defn}

Any adjunction
\[
\begin{diagram}[height=2em]
\cat{B}	\\
\uTo<F \ladj \dTo>G	\\
\cat{A}	\\
\end{diagram}
\]
induces a monad $(T, \mu, \eta)$ on $\cat{A}$: take $T= G\of F$, $\eta$
to be the unit of the adjunction, and 
\[
\mu = G\epsln F: GFGF \go GF
\]
where $\epsln$ is the counit.  Often $\cat{A}$ is the category of sets and
$\cat{B}$ is a category of `algebras' of some kind, as in the following
examples. 

\begin{example}	\lbl{eg:monad-monoid}
Take the free-forgetful%
\index{monoid!free}
adjunction between the category of monoids and the
category of sets, as above.  Let $(T, \mu, \eta)$ be the induced monad.
Then $TA$ is $\coprod_{n\in\nat} A^n$, the set of finite sequences of
elements of $A$, for any set $A$.  The multiplication $\mu$ strips inner
brackets from double sequences:
\[
\begin{array}{rcl}
T(TA)		&\goby{\mu_A}		&TA,	\\

((a_1^1, \ldots, a_1^{k_1}), \ldots, (a_n^1, \ldots, a_n^{k_n}))	&
\goesto	&
(a_1^1, \ldots, a_1^{k_1}, \ldots, a_n^1, \ldots, a_n^{k_n})	
\end{array}
\]
($n, k_i \in \nat, a_i^j \in A$).  The unit $\eta$ forms sequences of
length 1:
\[
\begin{array}{rcl}
A	&\goby{\eta_A}	&TA,	\\
a	&\goesto	&(a).
\end{array}
\]
\end{example}

\begin{example}	\lbl{eg:monad-R-mod}
Fix a ring $R$.  One can form the free $R$-module%
\index{module!ring@over ring}
on any given set, and
conversely one can take the underlying set of any $R$-module, giving an
adjunction 
\[
\begin{diagram}[height=2em]
R\hyph\fcat{Mod}	\\
\uTo<F \ladj \dTo>G	\\
\Set	\\
\end{diagram}%
\glo{RMod}
\]
hence a monad $(T, \mu, \eta)$ on $\Set$.  Explicitly, if $A$ is a set then
$TA$ is the set of formal $R$-linear combinations of elements of $A$.  The
multiplication $\mu$ realizes a formal linear combination of formal linear
combinations as a single formal linear combination, and the unit $\eta$
realizes an element of a set $A$ as a trivial linear combination of
elements of $A$.
\end{example}

\begin{example}	\lbl{eg:monad-alg-thy}
The same goes for all other `algebraic%
\index{algebraic theory}
theories': groups, Lie algebras,
Boolean algebras, \ldots.  The functor $T$ sends a set $A$ to the set of
formal words in the set $A$ (which in some cases, such as that of groups,
is cumbersome to describe).  There is no need for the ambient category
$\cat{A}$ to be $\Set$: the theory of topological groups,%
\index{group!topological}
for instance,
gives a monad on the category $\Top$%
\glo{Topcat}
of topological spaces.
\end{example}

A monad is meant to be an algebraic theory, so if we are handed a monad
then we ought to be able to say what its `models' are.  For instance, if we
are handed the monad of~\ref{eg:monad-R-mod} then its `models' should be
exactly $R$-modules.  Formally, if $T = (T, \mu, \eta)$ is a monad on a
category $\cat{A}$ then a \demph{$T$-algebra}%
\index{algebra!monad@for monad}
is an object $A \in \cat{A}$
together with a map $h: TA \go A$ compatible with the multiplication and
unit of the monad: see Mac Lane~\cite[VI.2]{MacCWM} for the axioms.  In the
case of~\ref{eg:monad-R-mod}, a $T$-algebra is a set $A$ equipped with a
function
\[
h:
\{ \textrm{formal }R \textrm{-linear combinations of elements of }A \}
\go 
A
\]
satisfying some axioms, and this does indeed amount exactly to an
$R$-module.%
\index{module!ring@over ring}

The category of algebras for a monad $T = (T, \mu, \eta)$ on a category
$\cat{A}$ is written $\cat{A}^T$.  There is an evident forgetful functor
$\cat{A}^T \go \cat{A}$, this has a left adjoint (forming `free%
\index{algebra!monad@for monad!free}
$T$-algebras'), and the monad on $\cat{A}$ induced by this adjunction is
just the original $T$.  So every monad arises from an adjunction, and
informally we have
\[
\{\textrm{monads on } \cat{A} \}
\subset
\{\textrm{adjunctions based on } \cat{A} \}.
\]
The inclusion is proper: not every adjunction is of the form $\cat{A}^T
\pile{\rTo_\top \\ \lTo} \cat{A}$ just described.  For instance, the
forgetful functor $\Top \go \Set$ has a left adjoint (forming discrete
spaces); the induced monad on $\Set$ is the identity, whose category of
algebras is merely $\Set$, and $\Set \not\eqv \Top$.  The adjunctions that
do arise from monads are called \demph{monadic}.%
\index{monadic adjunction}
 All of the adjunctions in
Examples \ref{eg:monad-monoid}--\ref{eg:monad-alg-thy} are monadic, and the
non-monadicity of the adjunction $\Top \pile{\rTo_\top \\ \lTo} \Set$
expresses the thought that topology%
\index{topology vs. algebra@topology \vs.\ algebra}
is not algebra.%
\index{monad|)}

\index{presheaf|(}
Presheaves will be important.  A \demph{presheaf} on a category $\cat{A}$
is a functor $\cat{A}^\op \go \Set$.  Any object $A \in \cat{A}$ gives rise
to a presheaf $\cat{A}(\dashbk, A)$ on $\cat{A}$, and this defines a
functor
\[
\begin{array}{rcl}
\cat{A}		&\go		&\ftrcat{\cat{A}^\op}{\Set}	\\
A		&\goesto	&\cat{A}(\dashbk, A),
\end{array}
\]
the \demph{Yoneda embedding}.%
\index{Yoneda!embedding}
 It is full and faithful.  This follows from
the Yoneda Lemma,%
\index{Yoneda!Lemma}
which states that if $A \in \cat{A}$ and $X$ is a
presheaf on $\cat{A}$ then natural transformations $\cat{A}(\dashbk, A) \go
X$ correspond one-to-one with elements of $XA$.

If $\Eee$ is a category and $S$ a set then there is a category $\Eee^S$,%
\glo{catpower}
a
power%
\index{power of category}
of $\Eee$, whose objects are $S$-indexed families of objects of
$\Eee$.  On the other hand, if $\Eee$ is a category and $E$ an object of
$\Eee$ then there is a \demph{slice%
\index{slice!category}
category} $\cat{E}/E$,%
\glo{slicecat}
whose objects
are maps $D \goby{p} E$ into $E$ and whose maps are commutative triangles.
If $S$ is a set then there is an equivalence of categories
\begin{equation}	\label{eq:Set-slice-power}
\Set^S \eqv \Set/S,
\end{equation}
given in one direction by taking the disjoint union of an $S$-indexed
family of sets, and in the other by taking fibres of a set over $S$.

There is an analogue of~\bref{eq:Set-slice-power} in which the set $S$ is
replaced by a category.  Fix a small category $\scat{A}$.  The replacement
for $\Set^S$ is $\ftrcat{\scat{A}^\op}{\Set}$, but what should replace the
slice category $\Set/S$?  First note that any presheaf $X$ on $\scat{A}$
gives rise to a category $\scat{A}/X$,%
\glo{catelts}
the \demph{category of elements}%
\lbl{p:defn-caty-elts}\index{category!elements@of elements}
of $X$, whose objects are pairs $(A, x)$ with $A \in \scat{A}$
and $x \in XA$ and whose maps
\[
(A', x') \go (A, x)
\]
are maps $f: A' \go A$ in $\scat{A}$ such that $x' = (Xf)(x)$.  There is an
evident forgetful functor $\scat{A}/X \go \scat{A}$,
the \demph{Grothendieck fibration}%
\index{fibration!Grothendieck}
of $X$.  This is an example of a
\demph{discrete fibration},%
\index{fibration!discrete}
that is, a functor $G: \cat{D} \go \cat{C}$
such that
\begin{quote}
  for any object $D \in \cat{D}$ and map $C' \goby{p} GD$ in $\cat{C}$,
  there is a unique map $D' \goby{q} D$ in $\cat{D}$ such that $Gq = p$. 
\end{quote}
Discrete fibrations over $\cat{C}$ (that is, with codomain $\cat{C}$) can
be made into a category $\fcat{DFib}(\cat{C})$%
\glo{DFibcat}
in a natural way, and this
is the desired generalization of slice category.  We then have an
equivalence
\[
\ftrcat{\scat{A}^\op}{\Set}
\eqv
\fcat{DFib}(\scat{A}),
\]
given in one direction by taking categories of elements, and in the other
by taking fibres in a suitable sense.  There is also a dual notion of
\demph{discrete opfibration},%
\lbl{p:defn-cl-d-opfib}\index{fibration!discrete opfibration} 
and an equivalence
\[
\ftrcat{\scat{A}}{\Set}
\eqv
\fcat{DOpfib}(\scat{A}).%
\glo{DOpfibcat}
\]

A \demph{presheaf%
\index{presheaf!category}
category} is a category equivalent to
$\ftrcat{\scat{A}^\op}{\Set}$ for some small $\scat{A}$.  The class of
presheaf categories is closed under slicing:
\begin{propn}	\lbl{propn:pshf-slice}
Let $\scat{A}$ be a small category and $X$ a presheaf on $\scat{A}$.  Then
there is an equivalence of categories
\[
\ftrcat{\scat{A}^\op}{\Set} / X
\eqv
\ftrcat{(\scat{A}/X)^\op}{\Set}.
\]
\ \done
\end{propn}%
\index{presheaf|)}

\index{internal!algebraic structure|(}%
Finally, we will need just a whisper of internal
category theory.  If
$\cat{A}$ is any category with finite products then an \demph{(internal)%
\index{group!internal}
group}
in $\cat{A}$ consists of an object $A \in \cat{A}$ together with
maps $m: A \times A \go A$ (multiplication), $e: 1 \go A$ (unit), and $i: A
\go A$ (inverses), such that certain diagrams expressing the group axioms
commute.  Thus, a group in $\Set$ is an ordinary group, a group in the
category of smooth manifolds is a Lie group, and so on.  A similar
definition pertains for algebraic structures other than groups.  Categories
themselves can be defined in this way: if $\cat{A}$ is any category with
pullbacks then an \demph{(internal)%
\index{category!internal}
category} $C$ in $\cat{A}$ is a diagram
\[
\begin{slopeydiag}
	&	&C_1	&	&	\\
	&\ldTo<\dom&	&\rdTo>\cod&	\\
C_0	&	&	&	&C_0	\\
\end{slopeydiag}%
\glo{C0cat}\glo{C1cat}\glo{domcat}\glo{codcat}
\]
in $\cat{A}$ together with maps
\[
C_1 \times_{C_0} C_1 \goby{\comp} C_1,
\diagspace
C_0 \goby{\ids} C_1%
\glo{compcat}\glo{idscat}
\]
in $\cat{A}$, satisfying certain axioms.  Here $C_1 \times_{C_0} C_1$
is the pullback%
\lbl{p:defn-caty-pb}
\[
\begin{slopeydiag}
	&	&C_1 \times_{C_0} C_1 \Spbk	&	&	\\
	&\ldTo	&				&\rdTo	&	\\
C_1	&	&				&	&C_1	\\
	&\rdTo<\cod&				&\ldTo>\dom&	\\
	&	&C_0.				&	&	\\
\end{slopeydiag}
\]
When $\cat{A} = \Set$ we recover the usual notion of small category: $C_0$
and $C_1$ are the sets of objects and of arrows, $\dom$ and $\cod$ are the
domain and codomain functions, $C_1 \times_{C_0} C_1$ is the set of
composable pairs of arrows, $\comp$ and $\ids$ are the functions
determining binary composition and identity maps, and the axioms specify
the domain and codomain of composites and identities and express
associativity and identity laws.  When $\cat{A} = \Top$ we obtain a notion
of `topological%
\index{topological category}\index{category!topological}
category',
in which both the set of objects and the set of
arrows carry a topology.  For instance, given any space $X$, there is a
topological category $C = \Pi_1 X$%
\glo{Pi1top}\index{fundamental!1-groupoid}
in which $C_0 = X$ and $C_1$ is $X^{[0,1]}/\sim$, the space of all paths in
$X$ factored out by path homotopy relative to endpoints.%
\index{internal!algebraic structure|)}

\section{Monoidal categories}
\lbl{sec:mon-cats}

Monoidal categories come in a variety of flavours: strict, weak, plain,
braided, symmetric.  We look briefly at strict monoidal categories but
spend most time on the more important weak case and on the coherence
theorem: every weak monoidal category is equivalent to a strict one.

In the terminology of the previous section, a \demph{strict%
\index{monoidal category!strict}
monoidal
category} is an internal monoid in $\Cat$, that is, a category $\cat{A}$
equipped with a functor
\[
\begin{array}{rrcl}
\otimes:	&\cat{A} \times \cat{A}	&\go		&\cat{A},	\\
		&(A, B)			&\goesto	&A\otimes B
\end{array}%
\glo{otimes}
\]
and an object $I\in \cat{A}$,%
\glo{unitobj}
obeying strict associativity and unit laws:
\[
(A\otimes B)\otimes C 
=
A\otimes (B\otimes C),
\diagspace
I\otimes A
=
A,
\diagspace
A\otimes I
= 
A
\]
for all objects $A, B, C \in \cat{A}$, and similarly for morphisms.
Functoriality of $\otimes$ encodes the `interchange%
\index{interchange}
laws':
\begin{equation}	\label{eq:mon-interchange}
(g' \of f') \otimes (g \of f)
=
(g' \otimes g) \of (f' \otimes f)
\end{equation}
for all maps $g', f', f, f$ for which these composites make sense, and $1_A
\otimes 1_B = 1_{A \otimes B}$ for all objects $A$ and $B$.

Since one is usually not interested in equality of objects in a category,
only in isomorphism, strict monoidal categories are quite rare.

\begin{example}	\lbl{eg:str-mon-endo}
The category $\ftrcat{\cat{C}}{\cat{C}}$ of endofunctors on a given
category $\cat{C}$ has a strict monoidal structure given by composition (as
$\otimes$) and 
% the identity functor 
$1_{\cat{C}}$ (as $I$).
\end{example}

\begin{example}	\lbl{eg:str-mon-D}
Given a natural number $n$ (possibly $0$), let $\lwr{n}$%
\glo{lwrn}
denote the
$n$-element set $\{1, \ldots, n\}$ with its usual total order.  Let
$\scat{D}$%
\glo{augsimplexcat}\index{augmented simplex category $\scat{D}$}
be the category whose objects are the natural numbers and whose maps $m\go
n$ are the order-preserving functions $\lwr{m} \go \lwr{n}$.  This is the
`augmented simplex category', one object bigger than the standard
topologists' $\Delta$, and is equivalent to the category of (possibly
empty) finite totally ordered sets.  It has a strict monoidal structure
given by addition and $\lwr{0}$.
\end{example}

\begin{example}	\lbl{eg:str-mon-comm}
A category with only one object is just a monoid
(p.~\pageref{p:degen-cat-monoid}): if the category is called $\cat{A}$ and
its single object is called $\star$ then the monoid is the set
$M = \cat{A}(\star, \star)$ with composition $\of$ as multiplication and the
identity $1 = 1_\star$ as unit.  A one-object strict monoidal%
\index{monoidal category!degenerate}
category therefore
consists of a set $M$ with monoid structures $(\of, 1)$ and $(\otimes, 1)$
(the latter being tensor of arrows in the monoidal category), such
that~\bref{eq:mon-interchange} holds for all $g', f', g, f \in M$.
Lemma~\ref{lemma:EH} below tells us that this forces the binary operations
$\of$ and $\otimes$ to be equal and commutative.  So a one-object strict
monoidal category is just a commutative%
\index{monoid!commutative}
monoid.
\end{example}

\begin{lemma}[Eckmann--Hilton~\cite{EH}]	\lbl{lemma:EH}%
\index{Eckmann--Hilton argument}
Suppose that $\of$ and $\otimes$ are binary operations on a set $M$,
satisfying~\bref{eq:mon-interchange} for all $g', f', g, f \in M$, and
suppose that $\of$ and $\otimes$ share a two-sided unit.  Then $\of =
\otimes$ and $\of$ is commutative.
% and associative.  
\end{lemma}
\begin{proof}
Write $1$ for the unit.  Then for $g, f \in M$,
\[
g \of f
=
(g \otimes 1) \of (1 \otimes f)
=
(g \of 1) \otimes (1 \of f)
=
g \otimes f,
\]
so $\of = \otimes$, and 
\[
g \of f
=
(1 \otimes g) \of (f \otimes 1)
=
(1 \of f) \otimes (g \of 1)
=
f \otimes g,
\]
so $\of$ is commutative.  
\done
\end{proof}

Much more common are weak monoidal categories, usually just called
`monoidal categories'.
\begin{defn}	\lbl{defn:mon-cat}
A \demph{(weak) monoidal category}%
\index{monoidal category!classical}
is a category $\cat{A}$ together with a
functor 
% \[
$
\otimes: \cat{A} \times \cat{A} \go \cat{A},
$
% \]
% 
an object $I\in\cat{A}$ (the \demph{unit}), and isomorphisms
\[
(A\otimes B) \otimes C
\rTo^{\alpha_{A, B, C}}_{\diso}
A\otimes (B\otimes C),
\diagspace
I \otimes A
\rTo^{\lambda_A}_{\diso}
A,
\diagspace
A\otimes I
\rTo^{\rho_A}_{\diso}
A%
\glo{monalpha}\glo{monlambda}\glo{monrho}
\]%
\index{associativity!isomorphism}%
\index{unit!isomorphism}%
natural in $A, B, C \in \cat{A}$ (\demph{coherence%
\index{coherence!isomorphism}
isomorphisms}), such
that the following diagrams commute for all $A, B, C, D \in \cat{A}$:
\[
\begin{array}{c}
\begin{diagram}[width=4em,height=2em,scriptlabels,tight,noPS]
	&	&	&
(A\otimes B) \otimes (C\otimes D)
				&	&	&	\\
	&	&
\ruTo(3,2)<{\alpha_{A\otimes B, C, D}}
			&	&	
\rdTo(3,2)>{\alpha_{A, B, C\otimes D}}
					&	&	\\
((A\otimes B) \otimes C)\otimes D
	&	&	&	&	&	&
A \otimes (B\otimes (C \otimes D))	\\
	&\rdTo(1,2)<{\alpha_{A, B, C} \otimes 1_D}
		&	&	&	&
\ruTo(1,2)>{1_A \otimes \alpha_{B, C, D}}
						&	\\
	&	
(A\otimes (B\otimes C))\otimes D
		&	&
\rTo_{\alpha_{A, B\otimes C, D}}
			&	&	
A\otimes ((B\otimes C)\otimes D)
					&	&	\\
\end{diagram}
\\
\\
\index{pentagon}\index{triangle!coherence axiom}
\begin{diagram}[size=2em,scriptlabels,noPS]
(A\otimes I) \otimes B
	&	&
\rTo^{\alpha_{A, I, B}}
			&	&
A\otimes (I\otimes B)	
					\\
	&
\rdTo<{\rho_A \otimes 1_B}	
		&	&
\ldTo>{1_A \otimes \lambda_B}	
				&	\\
	&	&
A\otimes B.	
			&	&	\\
\end{diagram}
\end{array}
\]
\end{defn}
The pentagon and triangle axioms ensure that `all diagrams' constructed out
of coherence isomorphisms commute.  This is one form of the coherence
theorem, discussed below.

\begin{example}
Strict monoidal categories can be identified with monoidal categories in
which all the components of $\alpha$, $\lambda$ and $\rho$ are identities.
\end{example}

\begin{example}
Let $\cat{A}$ be a category in which all finite%
\index{category!finite product}\index{monoidal category!cartesian}
products exist.
Choose a particular terminal object $1$, and for each $A, B \in
\cat{A}$ a particular product diagram $A \ogby{\mr{pr}_1} A\times B
\goby{\mr{pr}_2} B$.%
\glo{pri}
 Then $\cat{A}$ acquires a monoidal structure
with $\otimes = \times$ and $I = 1$; the maps $\alpha$, $\lambda$ and
$\rho$ are the canonical ones.
\end{example}

\begin{example}
For any commutative ring $R$, the category of $R$-modules%
\index{module!ring@over ring}
is monoidal with respect to the usual tensor $\otimes_R$ and unit object
$R$.
\end{example}

\begin{example}	\lbl{eg:mon-cat-loops}
Take a topological space with basepoint.  There is a monoidal category
whose objects are loops%
\index{loop space}
on the basepoint and whose maps are homotopy classes of loop homotopies
(relative to the basepoint).  We have to take homotopy classes so that the
ordinary categorical composition obeys associativity and identity laws.
Tensor is concatenation of loops (on objects) and gluing of homotopies (on
maps).  The coherence isomorphisms are the evident reparametrizations.
\end{example}

Earlier we met the notion of (internal) algebraic structures, such as
groups, in a category with finite products.  There is no clear way to
extend this to arbitrary monoidal categories, since to express an axiom
such as $x \cdot x^{-1} = 1$ diagrammatically requires the
product-projections.  We can, however, define a \demph{monoid}%
\index{monoid!monoidal category@in monoidal category}
in a
monoidal category $\cat{A}$ as an object $A$ together with maps
\[
m: A \otimes A \go A,
\diagspace
e: I \go A
\]
such that associativity and unit diagrams similar to those in
Definition~\ref{defn:monad} commute.  With the obvious notion of map, this
gives a category $\Mon(\cat{A})$%
\glo{Monofcat}
of monoids in $\cat{A}$.  When $\cat{A}$
is the category of sets, with product as monoidal structure, this is the
usual category of monoids.

There are various notions of map between monoidal categories.  In what
follows we use $\otimes$, $I$, $\alpha$, $\lambda$, and $\rho$ to denote
the monoidal structure of both the categories concerned.   
\begin{defn}	\lbl{defn:mon-ftr}
Let $\cat{A}$ and $\cat{A'}$ be monoidal categories.  A \demph{lax%
\index{monoidal functor!classical|(}
monoidal
functor} $F = (F, \phi): \cat{A} \go \cat{A'}$ is a functor $F: \cat{A} \go
\cat{A'}$ together with \demph{coherence%
\index{coherence!map}
maps}
\[
\phi_{A, B}: FA \otimes FB \go F(A\otimes B),
\diagspace
\phi_\cdot: I \go FI
\]
in $\cat{A'}$, the former natural in $A, B \in \cat{A}$, such that the
following diagrams commute for all $A, B, C \in \cat{A}$:
\[
\begin{array}{c}
\begin{diagram}[size=2em,scriptlabels]
(FA \otimes FB) \otimes FC	&\rTo^{\phi_{A,B} \otimes 1_{FC}}	&
F(A\otimes B) \otimes FC	&\rTo^{\phi_{A\otimes B, C}}		&
F((A\otimes B)\otimes C)	\\
\dTo<{\alpha_{FA, FB, FC}}	&					&
				&					&
\dTo>{F\alpha_{A, B, C}}	\\
FA \otimes (FB \otimes FC)	&\rTo_{1_{FA} \otimes \phi_{B,C}}	&
FA \otimes F(B \otimes C)	&\rTo_{\phi_{A, B\otimes C}}		&
F(A\otimes (B\otimes C))	\\
\end{diagram}
\\
\\
\begin{diagram}[size=2em,scriptlabels]
FA \otimes I	&\rTo^{1_{FA} \otimes \phi_\cdot}	&
% FA \otimes I	&\rTo^{1 \otimes \phi_\cdot}	&
FA \otimes FI	&\rTo^{\phi_{A, I}}			&
F(A\otimes I)	\\
\dTo<{\rho_{FA}}&					&
		&					&
\dTo>{F\rho_A}	\\
FA		&					&
\rEquals	&					&
FA		\\
\end{diagram}
% 
% \diagspace
\\
\\
\begin{diagram}[size=2em,scriptlabels]
I\otimes FA	&\rTo^{\phi_\cdot \otimes 1_{FA}}	&
% I\otimes FA	&\rTo^{\phi_\cdot \otimes 1}	&
FI \otimes FA	&\rTo^{\phi_{I, A}}			&
F(I\otimes A)	\\
\dTo<{\lambda_{FA}}&					&
		&					&
\dTo>{F\lambda_A}\\
FA		&					&
\rEquals	&					&
FA.		\\
\end{diagram}
\end{array}
\]
A \demph{colax monoidal functor} $F = (F, \phi): \cat{A} \go \cat{A'}$ is a
functor $F: \cat{A} \go \cat{A'}$ together with maps
\[
\phi_{A, B}: F(A \otimes B) \go FA \otimes FB,
\diagspace
\phi_\cdot: FI \go I
\]
satisfying axioms dual to those above.  A \demph{weak} (respectively,
\demph{strict}) \demph{monoidal functor} is a lax monoidal functor $(F,
\phi)$ in which all the maps $\phi_{A, B}$ and $\phi_\cdot$ are
isomorphisms (respectively, identities).
\end{defn}
We write $\MClax$%
\glo{MClax}
for the category of monoidal categories and lax maps, and similarly
$\fcat{MonCat}_\mr{colax}$, $\MCwk$, and $\MCstr$.  There are various
alternative systems of terminology; in particular, what we call weak
monoidal functors are sometimes called `strong%
\index{monoidal functor!classical|)}
monoidal functors' or just
`monoidal functors'.
\begin{example}	%\lbl{eg:lax-mon-Ab}
The forgetful functor $U: \Ab \go \Set$%
\glo{Ab}
from abelian groups to sets has a
lax monoidal structure with respect to the usual monoidal structures on
$\Ab$ and $\Set$, given by the canonical maps
\[
UA \times UB \go U(A \otimes B)
\]
($A, B \in \Ab$) and the map $1 \go U\integers$ picking out $0\in\integers$.
\end{example}

\begin{example}	\lbl{eg:mon-cat-action}
Let $\cat{C}$ be a category and $\cat{A}$ a monoidal category.  A weak
monoidal functor from $\cat{A}$ to the monoidal category
$\ftrcat{\cat{C}}{\cat{C}}$ of endofunctors on
$\cat{C}$~(\ref{eg:str-mon-endo}) is called an \demph{action}%
\index{action!monoidal category@of monoidal category}
of $\cat{A}$
on $\cat{C}$, and amounts to a functor $\cat{A} \times \cat{C} \go \cat{C}$
together with coherence isomorphisms satisfying axioms.
\end{example}

To state one of the forms of the coherence theorem we will need a notion of
equivalence of monoidal categories, and for this we need in turn a notion
of transformation.
\begin{defn}
Let $(F, \phi), (G, \psi): \cat{A} \go \cat{A'}$ be lax monoidal functors.
A \demph{monoidal%
\index{monoidal transformation!classical}
transformation} $(F, \phi) \go (G, \psi)$ is a natural
transformation $\sigma: F \go G$ such that the following diagrams commute
for all $A, B \in \cat{A}$:
\[
\begin{diagram}[size=2em,scriptlabels]
FA \otimes FB	&\rTo^{\sigma_A \otimes \sigma_B}	&GA \otimes GB	\\
\dTo<{\phi_{A,B}}&					&
\dTo>{\psi_{A,B}}\\
F(A\otimes B)	&\rTo_{\sigma_{A\otimes B}}		&G(A\otimes B)	\\
\end{diagram}
\diagspace
\begin{diagram}[size=2em,scriptlabels]
I			&\rEquals	&I			\\
\dTo<{\phi_\cdot}	&		&\dTo>{\psi_\cdot}	\\
FI			&\rTo_{\sigma_I}&GI.			\\
\end{diagram}
\]
\end{defn}

A weak monoidal functor $(F,\phi)$ is called an \demph{equivalence}%
\index{equivalence!classical monoidal categories@of classical monoidal categories} 
of
monoidal categories if it satisfies the conditions of the following
proposition.  
% Condition~\bref{item:mon-eqv-long} is the more morally
% correct, and condition~\bref{item:mon-eqv-short} the more convenient.
% 
\begin{propn}	\lbl{propn:mon-eqv-eqv}
The following conditions on a weak monoidal functor $(F, \phi): \cat{A} \go
\cat{A'}$ are equivalent:
\begin{enumerate}
\item	\lbl{item:mon-eqv-long}
there exist a weak monoidal functor $(G, \psi): \cat{A'} \go \cat{A}$ and
invertible monoidal transformations
\[
\eta: 1_{\cat{A}} \go (G, \psi)\of (F, \phi),
\diagspace
\epsln: (F, \phi) \of (G, \psi) \go 1_{\cat{A'}}
\]
\item	\lbl{item:mon-eqv-short}
the functor $F$ is an equivalence of categories.
\end{enumerate}
\end{propn}
\begin{proof}
If $F$ is an equivalence of categories then by
Proposition~\ref{propn:eqv-eqv}, there exist a functor $G$ and
transformations $\eta$ and $\epsln$ such that $(F, G, \eta, \epsln)$ is an
adjoint equivalence.  It is easy to verify that $G$ acquires a weak
monoidal structure and that $\eta$ and $\epsln$ are then invertible
monoidal transformations.
% The opposite implication is trivial.
\done
\end{proof}

\index{coherence!monoidal categories@for monoidal categories!classical|(} 
A coherence%
\index{coherence}
theorem is, roughly, a description of a structure that makes it
more manageable.%
\lbl{p:coherence-discussion}
For example, one coherence theorem for monoidal categories is that all
diagrams built out of the coherence isomorphisms commute.  Another is that
any weak monoidal category is equivalent to some strict monoidal category.
All non-trivial applications of monoidal categories rely on a coherence
theorem in some form; the axioms as they stand are just too unwieldy.
Indeed, one might argue that the `all diagrams commute' principle should be
an explicit part of the definition of monoidal category, and we will take
this approach when we come to define higher-dimensional categorical
structures.

`All diagrams commute' can be made precise in various ways.  A very direct
statement is in Mac%
\index{Mac Lane, Saunders}
Lane~\cite[VII.2]{MacCWM}, and a less direct (but
sharper) statement is~\ref{thm:diag-coh-mc} below.  A typical instance is
that the diagram
\begin{equation}	\label{diag:typical-coh}
\begin{diagram}[height=2em,width=2em,scriptlabels]
		&		&
(A\otimes (I\otimes B)) \otimes C
					&		&	\\
		&
\ldTo<{\alpha_{A, I, B}^{-1} \otimes 1}
				&	&
\rdTo(2,3)>{\alpha_{A, I\otimes B, C}}	
							&	\\
((A\otimes I) \otimes B)\otimes C	
		&		&	&		&	\\
\dTo<{(\rho_A \otimes 1)\otimes 1}	
		&		&	&		&
A \otimes ((I\otimes B) \otimes C)				\\
(A\otimes B)\otimes C		
		&		&	&
\ldTo(2,3)>{1\otimes (\lambda_B \otimes 1)}		
							&	\\
		&
\rdTo<{\alpha_{A, B, C}}	
				&	&		&	\\
		&		&
A\otimes (B\otimes C)	
					&		&	\\
\end{diagram}
\end{equation}
commutes for all objects $A, B, C$ of a monoidal category.  We will soon
see how this follows from the alternative form of the coherence theorem:
\begin{thm}[Coherence for monoidal categories]
\lbl{thm:coh-mon-wk-str}
Every monoidal category is equivalent to some strict monoidal category.
\end{thm}
Here is how \emph{not}%
\index{coherence!how not to prove}
to prove this: take a monoidal category $\cat{A}$, form a quotient strict
monoidal category $\cat{A'}$ by turning isomorphism into equality, and show
that the natural map $\cat{A} \go \cat{A'}$ is an equivalence.  To see why
this fails, first note that a monoidal category may have the property that
any two isomorphic objects are equal (and so, in particular, the tensor
product of objects is strictly associative and unital), but even so need
not be strict---the coherence isomorphisms need not be identities.  An
example can be found in Mac Lane~\cite[VII.1]{MacCWM}.  If $\cat{A}$ has
this property then identifying isomorphic objects of $\cat{A}$ has no
effect at all.  One might attempt to go further by identifying the
coherence isomorphisms with identities---for instance, identifying the two
maps
\[
(A\otimes B) \otimes C 
\parpair{\alpha_{A, B, C}}{1}
A \otimes (B\otimes C)
\]
---to make a strict monoidal category $\cat{A'}$; but then the quotient map
$\cat{A} \go \cat{A'}$ is not faithful, so not an equivalence.

\begin{prooflike}{Sketch proof of~\ref{thm:coh-mon-wk-str}}
This is a modification of Joyal%
\index{Joyal, Andr\'e|(}%
\index{Street, Ross!coherence for monoidal categories@on coherence for monoidal categories|(}
and Street's proof~\cite[1.4]{JSBTC}.  Let
$\cat{A}$ be a monoidal category.  We define a strict monoidal category
$\cat{A'}$ and a monoidal equivalence $\mathbf{y}: \cat{A} \go \cat{A'}$.
An object of $\cat{A'}$ is a pair $(E, \delta)$ where $E$ is an endofunctor
of the (unadorned) category $\cat{A}$ and $\delta$ is a family of
isomorphisms 
\[
\left(
\delta_{A, B}: (EA) \otimes B \goiso E(A\otimes B)
\right)_{A, B \in \cat{A}}
\]
natural in $A$ and $B$ and satisfying the evident coherence axioms.  Tensor
in $\cat{A'}$ is
\[
(E', \delta') \otimes (E, \delta) 
=
(E'\of E, \delta'')
\]
where $\delta''$ is defined in the only sensible way; then $\cat{A'}$ is a
strict monoidal category.  The functor $\mathbf{y}$ is given by
\[
\mathbf{y}(Z) = (Z\otimes\dashbk, \alpha_{Z, \dashbk, \dashbk}).
\]
It is weak monoidal and full, faithful and essentially surjective on
objects, so by~\ref{propn:mon-eqv-eqv} an equivalence of monoidal
categories.  \done
\end{prooflike}
Joyal and Street motivate their proof as a generalization of the Cayley%
\index{Cayley representation}%
\index{representation theorem}
Theorem representing any group as a group of permutations.  We find another
way of looking at it when we come to bicategories~(\ref{sec:bicats}).

Now let us deduce that `all diagrams commute', or at least, by way of
example, that diagram~\bref{diag:typical-coh} commutes.  For any objects
$A$, $B$ and $C$ of any monoidal category, write 
\[
(A \otimes (I\otimes B)) \otimes C
\parpair{\chi_{A, B, C}}{\omega_{A, B, C}}
A \otimes (B\otimes C)
\]
for the two composite coherence maps shown in~\bref{diag:typical-coh}.  Now
take a particular monoidal category $\cat{A}$ and objects $A, B, C \in
\cat{A}$; we want to conclude that $\chi_{A, B, C} = \omega_{A, B, C}$.  We
have a monoidal equivalence $(F, \phi)$ from $\cat{A}$ to a strict monoidal
category $\cat{A'}$, and a serially commutative diagram
\[
\begin{diagram}[size=2em,scriptlabels]
(FA \otimes (I\otimes FB)) \otimes FC	&
\pile{\rTo^{\chi_{FA, FB, FC}}\\ \rTo_{\omega_{FA, FB, FC}}}
					&
FA \otimes (FB \otimes FC)		\\
\dTo<{\phi}				&
					&
\dTo>{\phi}				\\
F( (A\otimes (I\otimes B)) \otimes C )	&
\pile{\rTo^{F \chi_{A, B, C}}\\ \rTo_{F \omega_{A, B, C}}}
		 			&
F( A \otimes (B\otimes C) )		\\
\end{diagram}
\]
where the maps labelled $\phi$ are built out of $\phi_\cdot$ and various
$\phi_{D, E}$'s.  (`Serially commutative' means that both the top and the
bottom square commute.)  Since $\cat{A'}$ is strict, $\chi_{FA, FB, FC} =
\omega_{FA, FB, FC}$; then since the $\phi$'s are isomorphisms, $F\chi_{A,
B, C} = F\omega_{A, B, C}$; then since $F$ is faithful, $\chi_{A, B, C} =
\omega_{A, B, C}$, as required.

There are similar diagrammatic coherence theorems for monoidal functors,
saying, for instance, that if $(F, \phi)$ is a lax monoidal functor then
any two maps
\[
((FA \otimes FB) \otimes I) \otimes (FC \otimes FD)
\parpairu
F(A \otimes ((B\otimes C) \otimes D))
\]
built out of copies of the coherence maps are equal.  The form of the
codomain is important, being $F$ applied to a product of objects; in
contrast, the coherence maps can be assembled to give two maps
\[
FI 
\parpairu 
FI \otimes FI
\]
that are in general not equal.  See Lewis~\cite{Lew}%
\index{Lewis, Geoffrey}
for both this
counterexample and a precise statement of coherence for monoidal functors.

Many everyday monoidal categories have a natural symmetric structure.
Formally, a symmetric monoidal category is a monoidal category $\cat{A}$
together with a specified isomorphism $\gamma_{A, B}: A \otimes B \go B
\otimes A$ for each pair $(A, B)$ of objects, satisfying coherence axioms.
There are various coherence theorems for symmetric monoidal categories
(Mac%
\index{Mac Lane, Saunders}
Lane~\cite{MacNAC}, Joyal and Street~\cite{JSBTC}).  Beware, however, that
the symmetry isomorphisms $\gamma_{A, B}$ cannot be turned into identities:
every symmetric monoidal category is equivalent to some symmetric strict
monoidal category, but not usually to any strict symmetric monoidal
category.  The latter structures---commutative monoids in $\Cat$---are
rare.%
\index{coherence!monoidal categories@for monoidal categories!classical|)} 

More general than symmetric monoidal categories are braided%
\index{monoidal category!braided}
monoidal
categories, mentioned in the Motivation for Topologists.  See Joyal%
\index{Joyal, Andr\'e|)}
and
Street~\cite{JSBTC}%
\index{Street, Ross!coherence for monoidal categories@on coherence for monoidal categories|)}
for the definitions and coherence theorems, and
Gordon,%
\index{Gordon, Robert}%
\index{tricategory}
Power and Street~\cite{GPS} for a 3-dimensional perspective.

\section{Enrichment}
\lbl{sec:cl-enr}

In many basic examples of categories $\cat{C}$, the hom-sets $\cat{C}(A,
B)$ are richer than mere sets.  For instance, if $\cat{C}$ is a category of
chain complexes then $\cat{C}(A, B)$ is an abelian group, and if $\cat{C}$
is a suitable category of topological spaces then $\cat{C}(A, B)$ is itself
a space.

This idea is called `enrichment' and can be formalized in various ways.
The best-known, enrichment in a monoidal category, is presented here.  We
will see later~(\ref{sec:enr-mtis}) that it is not the most natural or
general formalization, but it serves a purpose before we reach that point.

\begin{defn}	\lbl{defn:V-gph}
Let $\cat{V}$ be a category.  A \demph{$\cat{V}$-graph}%
\index{graph!enriched}
$X$ is a set $X_0$
together with a family $(X(x,x'))_{x,x' \in X_0}$ of objects of $\cat{V}$.
A \demph{map of $\cat{V}$-graphs} $f: X \go Y$ is a function $f_0: X_0 \go
Y_0$ together with a family of maps
\[
\left( X(x, x') \goby{f_{x, x'}} Y(f_0 x, f_0 x') \right)_{x, x' \in X_0}.
\]
We usually write both $f_0$ and $f_{x, x'}$ as just $f$.  The category of
$\cat{V}$-graphs is written $\cat{V}\hyph\Gph$.%
\glo{VGph}
\end{defn}
A $\Set$-graph is, then, an ordinary directed graph.  A category is a
 directed graph equipped with composition and identities, suggesting the
 following definition.
\begin{defn}	\lbl{defn:cl-enr-cat}
Let $(\cat{V}, \otimes, I)$ be a monoidal category.  A \demph{category
enriched in $\cat{V}$},%
\index{enrichment!category@of category!monoidal category@in monoidal category} 
or \demph{\cat{V}-enriched
category}, is a
$\cat{V}$-graph $A$ together with families of maps
\[
\left(
A(b, c) \otimes A(a, b) \goby{\comp_{a,b,c}} A(a,c)
\right)_{a, b, c \in A_0},
\diagspace
\left(
I \goby{\ids_a} A(a,a)
\right)_{a \in A_0}%
\glo{compenr}\glo{idsenr}
\]
in $\cat{V}$ satisfying associativity and identity axioms (expressed as
commutative diagrams).  A \demph{\cat{V}-enriched%
\index{enrichment!functor@of functor}%
\index{functor!enriched}
functor} $F: A \go B$ is
a map of the underlying $\cat{V}$-graphs commuting with the composition
maps $\comp_{a,b,c}$ and identity maps $\ids_a$.  This defines a category
$\cat{V}\hyph\Cat$.%
\glo{VCat}
\end{defn}
A $(\Set, \times, 1)$-enriched category is, of course, an ordinary (small)
category, and an $(\Ab, \otimes, \integers)$-enriched category is an
$\Ab$-category%
\index{Ab-category@$\Ab$-category}
in the sense of homological algebra.  A one-object
$\cat{V}$-enriched category is a monoid%
\index{monoid!monoidal category@in monoidal category}
in the monoidal category $\cat{V}$.
Compare and contrast internal%
\index{internal!enriched@\vs.\ enriched}%
\index{enrichment!internal@\vs.\ internal}
and enriched categories: in the case of
topological spaces, for instance, an internal category in $\Top$ is an
ordinary category equipped with a topology on the set of objects and a
topology on the set of all arrows, whereas a category enriched in $\Top$ is
an ordinary category equipped with a topology on each hom-set.

Any lax monoidal functor $Q = (Q, \phi): \cat{V} \go \cat{W}$ induces a
functor
\[
Q_*: \cat{V}\hyph\Cat \go \cat{W}\hyph\Cat.%
\glo{starenr}
\]
In particular, if $\cat{V}$ is any monoidal category then the functor
$\cat{V}(I, \dashbk): \cat{V} \go \Set$ has a natural lax monoidal
structure, and the induced functor defines the \demph{underlying%
\index{enrichment!category@of category!underlying category}
category}
of a $\cat{V}$-enriched category.  This does exactly what we would expect
in the familiar cases of $\cat{V}$.

In the next section we will enrich in categories $\cat{V}$ whose monoidal
structure is ordinary (cartesian) product.  The following result will be
useful; its proof is straightforward.
\begin{propn}	\lbl{propn:fin-prod-enr}
\begin{enumerate}
\item	\lbl{item:fin-prod-enr-cat}
  If $\cat{V}$ is a category with finite products then the category
  $\cat{V}\hyph\Cat$ also has finite products.
\item	\lbl{item:fin-prod-enr-ftr}
  If $Q: \cat{V} \go \cat{W}$ is a finite-product-preserving functor
  between categories with finite products then the induced functor $Q_*:
  \cat{V}\hyph\Cat \go \cat{W}\hyph\Cat$ also preserves finite products.
\done
\end{enumerate}
\end{propn}

The theory of categories enriched in a monoidal category can be taken much
further: see Kelly~\cite{KelBCE}, for instance.  Under the assumption that
$\cat{V}$ is symmetric monoidal closed and has all limits and colimits,
very large parts of ordinary category theory can be extended to the
$\cat{V}$-enriched context.

\section{Strict $n$-categories}
\lbl{sec:cl-strict}

Strict $n$-categories are not encountered nearly as often as their weak
cousins, but there are nevertheless some significant examples.

We start with a very short definition of strict $n$-category, and some
examples.  This definition, being iterative, can seem opaque, so we then
formulate a much longer but equivalent definition providing a complementary
viewpoint.  We then look briefly at the infinite-dimensional case, strict
$\omega$-categories, and at the cubical analogue of strict $n$-categories.

The definition of strict $n$-category uses enrichment%
\index{enrichment!category@of category!finite product category@in finite product category}%
\index{enrichment!define n-category@to define $n$-category} 
in a category with
finite products.  
\begin{defn}	\lbl{defn:strict-n-cat-enr}
Let $\left( \strcat{n} \right)_{n\in\nat}$%
\glo{strcatn}
be the sequence of categories
given inductively by
\[
\strcat{0} = \Set,
\diagspace
\strcat{(n+1)} = (\strcat{n})\hyph\Cat.
\]
A \demph{strict $n$-category}%
\index{n-category@$n$-category!strict}
is an object of $\strcat{n}$, and a
\demph{strict $n$-functor}%
\index{n-functor, strict@$n$-functor, strict}
is a map in $\strcat{n}$.
\end{defn}
This makes sense by
Proposition~\ref{propn:fin-prod-enr}\bref{item:fin-prod-enr-cat}.  

Strict $0$-categories are sets and strict $1$-categories are categories.  A
strict $2$-category%
\index{two-category@2-category!strict}
$A$ consists of a set $A_0$, a category $A(a,b)$ for
each $a, b \in A_0$, composition functors as in~\ref{defn:cl-enr-cat}, and
an identity object of $A(a,a)$ for each $a\in A_0$, all obeying
associativity and identity laws.  

\begin{example}
There is a (large) strict 2-category $\cat{A}$ in which $\cat{A}_0$ is the
class of topological%
\index{topological space!two-category of spaces@2-category of spaces}
spaces and, for spaces $X$ and $Y$, $\cat{A}(X, Y)$ is
the category whose objects are continuous maps $X \go Y$ and whose arrows
are homotopy
classes%
\index{homotopy!classes}
of homotopies.  (We need to take homotopy classes so
that composition in $\cat{A}(X, Y)$ is associative and unital.)  The
composition functors
\begin{equation}	\label{eq:2-cat-comp-ftr}
\cat{A}(Y, Z) \times \cat{A}(X, Y) \go \cat{A}(X, Z)
\end{equation}
and the identity objects $1_X \in \cat{A}(X, X)$ are the obvious ones.
\end{example}

\begin{example}
Similarly, there is a strict 2-category $\cat{A}$ in which $\cat{A}_0$ is
the class of chain%
\index{chain complex!two-category of complexes@2-category of complexes}
complexes and $\cat{A}(X, Y)$ is the category whose objects are chain maps
$X \go Y$ and whose arrows are homotopy classes%
\index{homotopy!classes}
of chain homotopies, in the
sense of p.~\pageref{p:ch-hty-hty}.  (This time we need to take homotopy
classes in order that the composition functors~\bref{eq:2-cat-comp-ftr}
really are functorial.)
\end{example}

\begin{example}	\lbl{eg:str-2-cat-Cat}
There is, self-referentially, a strict 2-category $\Cat$%
\glo{Cat2cat}%
\index{two-category@2-category!categories@of categories}%
\index{category!two-category of@2-category of}
of categories.
Here $\Cat_0$ is the class of small categories and $\Cat(C, D)$ is the
functor category $\ftrcat{C}{D}$.  In fact, there is for each $n\in\nat$ a
strict $(n+1)$-category%
\index{n-category@$n$-category!n-ZZZcategory of@$(n+1)$-category of}
of strict $n$-categories: it can be proved by
induction that for each $n$ the category $\strcat{n}$ is cartesian closed,
which implies that it is naturally enriched in itself, in other words,
forms a strict $(n+1)$-category.  Sensitive readers may find this shocking:
the entire $(n+1)$-category of strict $n$-categories can be extracted from
the mere $1$-category.
\end{example}

We now build up to an alternative and more explicit definition of strict
$n$-category.  

A category can be regarded as a directed graph%
\index{graph!directed}
with structure.  The
most obvious $n$-dimensional analogue of a directed graph uses spherical or
`globular' shapes, as in the following definition.
\begin{defn}
Let $n\in\nat$.  An \demph{$n$-globular%
\index{n-globular set@$n$-globular set}
set} $X$ is a diagram 
\[
X(n)
\parpair{s}{t}
X(n-1)
\parpair{s}{t}
\diagspace 
\cdots
\diagspace 
\parpair{s}{t}
X(0)%
\glo{globsce}\glo{globtgt}
\]
of sets and functions, such that
\begin{equation}	\label{eq:glob}
s(s(x)) = s(t(x)),
\diagspace
t(s(x)) = t(t(x))
\end{equation}
for all $m\in \{2, \ldots, n \}$ and $x\in X(m)$.
\end{defn}

Alternatively, an $n$-globular set is a presheaf on the category
$\scat{G}_n$%
\glo{Gn}
generated by objects and arrows
\[
n
\pile{\lTo^{\scriptstyle \sigma_n}\\ \lTo_{\scriptstyle \tau_n}}
n-1
\pile{\lTo^{\scriptstyle \sigma_{n-1}}\\ \lTo_{\scriptstyle \tau_{n-1}}}
\diagspace
\cdots
\diagspace
\pile{\lTo^{\scriptstyle \sigma_1}\\ \lTo_{\scriptstyle \tau_1}}
0
\]
subject to equations
\[
\sigma_m \of \sigma_{m-1} = \tau_m \of \sigma_{m-1},
\diagspace
\sigma_m \of \tau_{m-1} = \tau_m \of \tau_{m-1}
\]
($m \in \{ 2, \ldots, n \}$).  The category of $n$-globular sets can then
be defined as the presheaf category $\ftrcat{\scat{G}_n^\op}{\Set}$.

Let $X$ be an $n$-globular set.  Elements of $X(m)$ are called
\demph{$m$-cells}%
\index{cell!globular set@of globular set}
of $X$ and drawn as labels on an $m$-dimensional disk.
Thus, $a \in X(0)$ is drawn as
\[
\gzero{a}
\]
(and sometimes called an \demph{object}%
\index{object}
rather than a 0-cell), and $f \in
X(1)$ is drawn as
\[
\gfst{a} \gone{f} \glst{b}
\]
where $a = s(f)$ and $b = t(f)$.  We call $s(x)$ the \demph{source}%
\index{source}
of $x$,
and $t(x)$ the \demph{target};%
\index{target}
these are alternative names for `domain' and
`codomain'.  A 2-cell $\alpha \in X(2)$ is drawn as
\[
\gfst{a} \gtwo{f}{g}{\alpha} \glst{b}
\]
where
\[
f = s(\alpha),
\diagspace
g = t(\alpha),
\diagspace
a = s(f) = s(g),
\diagspace
b = t(f) = t(g).
\]
That $s(f) = s(g)$ and $t(f) = t(g)$ follows from the globularity
equations~\bref{eq:glob}.  A 3-cell $x \in X(3)$ is drawn as
\[
\gfst{a} \gthreecell{f}{g}{\alpha}{\beta}{x} \glst{b}
\]
where $\alpha = s(x)$, $\beta = t(x)$, and so on.  Sometimes, as in the
Motivation for Topologists, I have used double-shafted arrows%
\index{arrow!shafts of}\index{cell!depiction of}
for 2-cells,
triple-shafted arrows for 3-cells, and so on, and sometimes, as here, I
have stuck to single-shafted arrows for cells of all dimensions; this is
a purely visual choice.

\begin{example}	\lbl{eg:n-glob-set-Pi}
Let $n\in\nat$ and let $S$ be a topological space: then there is an
$n$-globular%
\index{n-globular set@$n$-globular set!space@from space}%
\index{topological space!n-globular set from@$n$-globular set from}
set in which an $m$-cell is a labelled $m$-dimensional disk in
$S$.  Formally, let $D^m$%
\glo{disk}
be the closed $m$-dimensional Euclidean disk%
\index{disk}
(ball), and consider the diagram
\[
D^n
\pile{\lTo\\ \lTo}
D^{n-1}
\pile{\lTo\\ \lTo}
\diagspace
\cdots 
\diagspace 
\pile{\lTo\\ \lTo}
D^0 = 1
\]
formed by embedding $D^{m-1}$ as the upper or lower cap of $D^m$.  This is
a functor $\scat{G}_n \go \Top$ (an `$n$-coglobular space'), and so induces
a functor $\scat{G}_n^\op \go \Set$, the $n$-globular set
\[
\Top(D^n, S)
\pile{\rTo\\ \rTo}
\Top(D^{n-1}, S)
\pile{\rTo\\ \rTo}
\diagspace
\cdots 
\diagspace 
\pile{\rTo\\ \rTo}
\Top(D^0, S).
\]
\end{example}

\begin{example}	\lbl{eg:n-glob-set-ch-cx}
Analogously, any non-negatively graded chain%
\index{chain complex!n-globular set from@$n$-globular set from}
complex $C$ of abelian groups
give rise to an $n$-globular%
\index{n-globular set@$n$-globular set!chain complex@from chain complex}
set $X$ for each $n\in\nat$.  An $m$-cell of
$X$ is not quite just an element of $C_m$; for instance, we regard an
element $f$ of $C_1$ as a 1-cell $a \go b$ for any $a, b \in C_0$ such that
$d(f) = b - a$.  In general, an element of $X(m)$ is a $(2m+1)$-tuple
\[
\mathbf{c} =
(c_m, c_{m-1}^-, c_{m-1}^+, c_{m-2}^-, c_{m-2}^+, \ldots, c_0^-, c_0^+)
\]
where $c_m \in C_m$, $c_p^-, c_p^+ \in C_p$, and 
\[
d(c_m) = c_{m-1}^+ - c_{m-1}^-,
\diagspace
d(c_p^-) = d(c_p^+) = c_{p-1}^+ - c_{p-1}^-
\]
for all $p \in \{1, \ldots, m-1\}$.  We then put
\[
s(\mathbf{c}) = 
(c_{m-1}^-, c_{m-2}^-, c_{m-2}^+, \ldots, c_0^-, c_0^+)
\]
and dually the target.
\end{example}

In the following alternative definition, a strict $n$-category is an
$n$-globular set equipped with identities and various binary composition
operations, satisfying various axioms.  Identities are simple: every
$p$-cell $x$ has an identity $(p+1)$-cell $1_x$ on it ($0\leq p < n$), as
in
\[
\gzero{a} \ \goesto\  \gfst{a}\gone{1_a}\glst{a},
\diagspace
\gfst{a}\gone{f}\glst{b} \ \goesto\  \gfst{a}\gtwo{f}{f}{1_f}\glst{b}
\]
($p = 0, 1$).  There are $m$ different binary composition operations for
$m$-cells.  When $m=1$ this is ordinary categorical composition.  We saw
the 2 possibilities for composing 2-cells on p.~\pageref{p:2-cell-comps}.
The 3 ways of composing 3-cells are drawn as
\[
\gfstsu\gthreecellu\gzersu\gthreecellu\glstsu,
\diagspace
\gfstsu\gspecialthree\glstsu,
\diagspace
\gfstsu\gspecialtwo\glstsu.
\]

To express this formally we write, for any $n$-globular set $A$ and $0\leq
p\leq m\leq n$, 
\[
A(m) \times_{A(p)} A(m) 
=
\{
(x', x) \in A(m) \times A(m)
\such
t^{m-p}(x) = s^{m-p}(x') 
\},
\]
the set of pairs of $m$-cells with the potential to be joined along
$p$-cells.

\begin{defn}	\lbl{defn:strict-n-cat-glob}
Let $n\in\nat$.  A \demph{strict%
\index{n-category@$n$-category!strict}
$n$-category} is an $n$-globular set $A$
equipped with
\begin{itemize}
\item a function $\ofdim{p}: A(m) \times_{A(p)} A(m) \go A(m)$%
\glo{ofdim}
for each $0\leq
  p < m \leq n$; we write $\ofdim{p}(x', x)$ as $x' \ofdim{p} x$ and call
  it a \demph{composite} of $x$ and $x'$
\item a function $i: A(p) \go A(p+1)$%
\glo{istrncat}
for each $0 \leq p < n$; we write $i(x)$ as $1_x$%
\glo{onestrncat}
and call it the
\demph{identity} on $x$,
\end{itemize}
satisfying the following axioms:
\begin{enumerate}
\item	\lbl{item:s-t-comp}
(sources and targets of composites) if $0\leq p < m \leq n$ and $(x',
x) \in A(m) \times_{A(p)} A(m)$ then
\[
\begin{array}{llll}
s(x' \ofdim{p} x) = 
s(x) 			&	
\textrm{and}			&
t(x' \ofdim{p} x) = 
t(x') 			&
\textrm{if }
p = m-1	\\
s(x' \ofdim{p} x) = 
s(x') \ofdim{p} s(x)	&
\textrm{and}			&
t(x' \ofdim{p} x) = 
t(x') \ofdim{p} t(x)	&
\textrm{if }
p \leq m-2
\end{array}
\]
\item	%\lbl{item:s-t-ids}
(sources and targets of identities) if $0 \leq p < n$ and $x \in
  A(p)$ then $s(1_x) = x = t(1_x)$
\item (associativity) if $0\leq p < m \leq n$ and $x, x', x''
  \in A(m)$ with $(x'', x'), (x', x) \in A(m) \times_{A(p)} A(m)$ then
  \[
  (x'' \ofdim{p} x') \ofdim{p} x 
  = 
  x'' \ofdim{p} (x' \ofdim{p} x)
  \]
\item (identities) if $0\leq p < m \leq n$ and $x\in A(m)$ then 
  \[
  i^{m-p}(t^{m-p}(x)) \ofdim{p} x
  = 
  x
  =
  x \ofdim{p} i^{m-p}(s^{m-p}(x))
  \]
\item (binary interchange)%
\index{interchange}
if $0\leq q < p < m \leq n$ and $x, x', y, y'
  \in A(m)$ with
  \[
  (y', y), (x', x) \in A(m) \times_{A(p)} A(m),
  \ \ 
  (y', x'), (y, x) \in A(m) \times_{A(q)} A(m)
  \]
  then 
  \[
  (y' \ofdim{p} y) \ofdim{q} (x' \ofdim{p} x) 
  = 
  (y' \ofdim{q} x') \ofdim{p} (y \ofdim{q} x)
  \]
\item (nullary interchange)%
\index{interchange}
if $0\leq q < p < n$ and $(x', x) \in A(p)
  \times_{A(q)} A(p)$ then $1_{x'} \ofdim{q} 1_x = 1_{x' \sof_q x}$.
\end{enumerate}
If $A$ and $B$ are strict $n$-categories then a \demph{strict $n$-functor}%
\index{n-functor, strict@$n$-functor, strict}
is a map $f: A \go B$ of the underlying $n$-globular sets commuting with
composition and identities.  This defines a category $\strcat{n}$ of strict
$n$-categories. 
\end{defn}

\begin{propn}	\lbl{propn:str-n-cats-comparison}
The categories $\strcat{n}$ defined in~\ref{defn:strict-n-cat-enr}
and~\ref{defn:strict-n-cat-glob} are equivalent.
\end{propn}
\begin{prooflike}{Sketch proof}
We first compare the underlying graph structures.  Define for each
$n\in\nat$ the category $n\hyph\Gph$%
\glo{nGph}
of \demph{$n$-graphs}%
\index{n-graph@$n$-graph}
by
\[
0\hyph\Gph = \Set,
\diagspace
(n+1)\hyph\Gph = (n\hyph\Gph)\hyph\Gph.
\]
An $(n+1)$-globular set amounts to a graph of $n$-globular sets: precisely,
an $(n+1)$-globular set $X$ corresponds to the graph $(X(a,b))_{a, b \in
X(0)}$ where $X(a,b)$ is the $n$-globular set defined by
\[
(X(a,b))(m) 
=
\{
x \in X(m+1)
\such
s^{m+1}(x) = a,\ 
t^{m+1}(x) = b
\}.
\]
So by induction, $n\hyph\Gph \eqv \ftrcat{\scat{G}_n^\op}{\Set}$.

Now we bring in the algebra.  Given a strict $(n+1)$-category $A$ in the
sense of~\ref{defn:strict-n-cat-glob}, the functions $\ofdim{p}$ and $i:
A(p) \go A(p+1)$ taken over $1\leq p < n$ give a strict $n$-category
structure on $A(a,b)$ for each $a, b \in A(0)$.  This determines a graph
$(A(a,b))_{a, b \in A(0)}$ of strict $n$-categories.  Moreover, the
functions $\ofdim{0}$ and $i: A(0) \go A(1)$ give this graph the structure
of a category enriched in $\strcat{n}$.  With a little work we find that a
strict $(n+1)$-category in the sense of~\ref{defn:strict-n-cat-glob} is, in
fact, \emph{exactly} a category enriched in $\strcat{n}$, and by induction
we are done.  \done
\end{prooflike}

Some examples of $n$-categories are more easily described using one
definition than the other.  It also helps to have the terminology of both
at hand.  We can now say that the strict 2-category $\Cat$
of~\ref{eg:str-2-cat-Cat} has categories as 0-cells, functors as 1-cells,
and natural transformations as 2-cells.  There are two kinds of composition
of natural transformations.  We usually write the (`vertical')%
\index{composition!vertical}
composite
of natural transformations
\[
C \cthree{F}{G}{H}{\alpha}{\beta} C'
\]
as $\beta\of\alpha$%
\glo{vertcompnat}
(rather than the standard $\beta\ofdim{1}\alpha$) and
the (`horizontal')%
\index{composition!horizontal}
composite of natural transformations
\[
C \ctwo{F}{G}{\alpha} C' \ctwo{F'}{G'}{\alpha'} C''
\]
as $\alpha' * \alpha$%
\glo{horizcompnat}
(rather than $\alpha' \ofdim{0} \alpha$).

\begin{example}	\lbl{eg:str-n-cat-ch-cx}
Let $n\in\nat$ and let $C$ be a chain%
\index{chain complex!n-category from@$n$-category from}
complex, as
in~\ref{eg:n-glob-set-ch-cx}.  Then the $n$-globular set $X$ arising from
$C$ has the structure of a strict $n$-category:%
\index{n-category@$n$-category!chain complex@from chain complex}
if $(\mathbf{d},
\mathbf{c}) \in X(m) \times_{X(p)} X(m)$ then $\mathbf{d} \ofdim{p}
\mathbf{c} = \mathbf{e}$ where, for $0\leq r \leq n$ and $\sigma \in \{ -,
+ \}$,
\[
e_r^\sigma
=
\left\{
\begin{array}{ll}
d_r^\sigma + c_r^\sigma	&\textrm{if } r > p	\\
c_p^-			&\textrm{if } r = p \textrm{ and } \sigma = -	\\
d_p^+			&\textrm{if } r = p \textrm{ and } \sigma = +	\\
c_r^\sigma = d_r^\sigma	&\textrm{if } r < p.
\end{array}
\right.
\]
\end{example}

Strict $\omega$-categories can also be defined in both styles.  The
globular or `global' definition is obvious.  Let $\scat{G}$%
\glo{G}
be the category
generated by objects and arrows
\[
\cdots
\diagspace
\pile{\lTo^{\scriptstyle \sigma_{m+1}}\\ \lTo_{\scriptstyle \tau_{m+1}}}
m
\pile{\lTo^{\scriptstyle \sigma_m}\\ \lTo_{\scriptstyle \tau_m}}
m-1
\pile{\lTo^{\scriptstyle \sigma_{m-1}}\\ \lTo_{\scriptstyle \tau_{m-1}}}
\diagspace
\cdots
\diagspace
\pile{\lTo^{\scriptstyle \sigma_1}\\ \lTo_{\scriptstyle \tau_1}}
0
\]
subject to the usual globularity equations.  Then
$\ftrcat{\scat{G}^\op}{\Set}$ is the category of \demph{globular%
\index{globular set}
sets}, and
the category $\strcat{\omega}$%
\glo{strcatomega}
of \demph{strict $\omega$-categories}%
\index{omega-category@$\omega$-category!strict}
is
defined just as in~\ref{defn:strict-n-cat-glob} but without the upper limit
of $n$.  

\begin{example}	\lbl{eg:str-omega-ch-cx}
Any chain%
\index{chain complex!omega-category from@$\omega$-category from}
complex $C$ gives rise to a strict $\omega$-category%
\index{omega-category@$\omega$-category!chain complex@from chain complex}
$X$, just as
in the previous example.  In fact, $X$ is an abelian group in
$\strcat{\omega}$, and in this way the category of abelian groups in
$\strcat{\omega}$ is equivalent to the category of non-negatively graded
chain complexes of abelian groups.
\end{example}

For the enriched or `local' definition, we first define a sequence
\[
\label{p:forgetful-strict-n}
\cdots 
\goby{S_{n+1}} 
\strcat{(n+1)}
\goby{S_n}
\strcat{n}
\goby{S_{n-1}}
\diagspace
\cdots
\diagspace
% \goby{S_1}
\Cat
\goby{S_0}
\Set%
\glo{Sstrncat}
\]
of finite-product-preserving functors by $S_0 = \mr{ob}$ and $S_{n+1} =
(S_n)_*$, which is possible by
Proposition~\ref{propn:fin-prod-enr}\bref{item:fin-prod-enr-ftr}.  We then
take $\strcat{\omega}$ to be the limit of this diagram in $\CAT$.  It is
easy to prove
\begin{propn}
The two categories $\strcat{\omega}$ just defined are equivalent.
\done
\end{propn}

Strict $n$-categories use globular shapes; there is also a cubical
analogue.  Again, there is a short inductive definition and a longer
explicit version.  The short form uses internal rather than enriched
categories.  If $\cat{V}$ is a category with pullbacks then we write
$\Cat(\cat{V})$%
\glo{Catinternal}
for the category of internal%
\index{category!internal}
categories in $\cat{V}$,
which, it can be shown, also has pullbacks.  
\begin{defn}
Let $(\strtuplecat{n})_{n\in\nat}$ be the sequence of categories given
inductively by
\[
\strtuplecat{0} = \Set,
\diagspace
\strtuplecat{(n+1)} = \Cat(\strtuplecat{n}).
\]
A \demph{strict $n$-tuple category}%
\index{n-tuple category@$n$-tuple category!strict}
(or `strict cubical
$n$-category') is an object of $\strtuplecat{n}$. 
\end{defn}

A strict single ($=$ 1-tuple) category is just a category.  A strict
double%
\index{double category!strict}
($=$ 2-tuple) category $D$ is a diagram  
\[
\begin{slopeydiag}
	&	&D_1	&	&	\\
	&\ldTo<\dom&	&\rdTo>\cod&	\\
D_0	&	&	&	&D_0	\\
\end{slopeydiag}
\]
of categories and functors, with extra structure.  The objects of $D_0$ are
called the \demph{0-cells}%
\index{cell!strict double category@of strict double category}
or \demph{objects} of $D$, the maps in $D_0$ are
the \demph{vertical 1-cells} of $D$, the objects of $D_1$ are the
\demph{horizontal 1-cells} of $D$, and the maps in $D_1$ are the
\demph{2-cells} of $D$, as in the picture
\begin{equation}	\label{diag:str-dbl-two-cell}
\begin{fcdiagram}
a	&\rTo^m			&a'	\\
\dTo<f	&\Downarrow\tcs{\theta}	&\dTo>{f'}\\
b	&\rTo_p			&b'	\\
\end{fcdiagram}
\end{equation}
where $a \goby{f} b$, $a' \goby{f'} b'$ are maps in $D_0$ and $m
\goby{\theta} p$ is a map in $D_1$, with $\dom(m) = a$, $\dom(\theta)=f$,
and so on.  The `extra structure' consists of various kinds of composition
and identities, so that vertical 1-cells can be composed vertically,
horizontal 1-cells can be composed horizontally, and 2-cells can be
composed both vertically and horizontally.  Thus, any $p \times q$ grid of
2-cells has a unique 2-cell composite, for any $p, q \in \nat$.

\begin{example}	\lbl{eg:2-cub-glob}
A strict double category in which all vertical 1-cells (or all horizontal
1-cells) are identities is just a strict 2-category.  
\end{example}

\begin{example}
A strict 2-category $A$ gives rise to a strict double category $D$ in two
other ways: take the 0-cells of $D$ to be those of $A$, both the
vertical and the horizontal 1-cells of $D$ to be the 1-cells of $A$, and
the 2-cells~\bref{diag:str-dbl-two-cell} of $D$ to be the 2-cells
\[
\begin{fcdiagram}
a	&\rTo^m			&a'	\\
\dTo<f	&\swnt\tcs{\theta}	&\dTo>{f'}\\
b	&\rTo_p			&b'	\\
\end{fcdiagram}
\]
in $A$, or the same with the arrow for $\theta$ reversed.
\end{example}

The longer definition of strict $n$-tuple category is omitted since we will
not need it, but goes roughly as follows.  As for the globular case, an
$n$-tuple category is a presheaf with extra algebraic structure.  

Let $\scat{H} = \scat{G}_1 = (1 \pile{\lTo^\sigma\\ \lTo_\tau} 0)$, so that
a presheaf on $\scat{H}$ is a directed graph.%
\index{graph!directed}
 If $n\in\nat$ then an
\demph{$n$-cubical%
\index{cubical!set}\index{n-cubical set@$n$-cubical set}
set} is a presheaf on $\scat{H}^n$.  An $n$-cubical set
consists of a set $X(M)$ for each $M \sub \{1, \ldots, n\}$, and a function
$X(\xi): X(M) \go X(P)$ for each $P \sub M$ and function $\xi: M\without P
\go \{-, +\}$, satisfying functoriality axioms.  For instance, a 2-cubical
set $X$ consists of sets
\begin{eqnarray*}
X(\emptyset)	&=	&\{ 0 \textrm{-cells} \}	\\
X(\{1\})	&=	&\{ \textrm{vertical } 1 \textrm{-cells} \}	\\
X(\{2\})	&=	&\{ \textrm{horizontal } 1 \textrm{-cells} \}	\\
X(\{1, 2\})	&=	&\{ 2 \textrm{-cells} \}
\end{eqnarray*}
with various source and target functions between them. 

A strict $n$-tuple category can then be defined as an $n$-cubical set $A$
together with a composition function
\[
A(P \cup \{ m \}) \times_{A(P)} A(P \cup \{ m \})
\go
A(P \cup \{ m \}) 
\]
and an identity function
\[
A(P)
\go
A(P \cup \{ m \}) 
\]
for each $P \sub \{1, \ldots, n\}$ and $m \in \{1, \ldots, n\} \without P$,
satisfying axioms.  Strict $n$-categories can be identified with strict
$n$-tuple categories whose underlying $n$-cubical set is degenerate%
\index{n-tuple category@$n$-tuple category!degenerate}
in a
certain way, generalizing Example~\ref{eg:2-cub-glob}.

\section{Bicategories}
\lbl{sec:bicats}

Bicategories are to strict 2-categories as weak monoidal categories are to
strict monoidal categories.  There are many other formalizations of the
idea of weak 2-category (see~\ref{sec:notions-bicat}), but bicategories are
the oldest and best-known.

We define a bicategory as a category `weakly enriched' in $\Cat$.  An
alternative definition as a 2-globular set with structure can be found in
B\'enabou~\cite{Ben},%
\index{Benabou, Jean@B\'enabou, Jean}
where bicategories were introduced.

\begin{defn}	\lbl{defn:cl-bicat}
A \demph{bicategory}%
\index{bicategory!classical}
$\cat{B}$ consists of
\begin{itemize}
\item a class $\cat{B}_0$, whose elements are called the \demph{objects} or
\demph{$0$-cells}%
\index{cell!classical bicategory@of classical bicategory}
of $\cat{B}$
\item for each $A, B \in \cat{B}_0$, a category $\cat{B}(A, B)$, whose
  objects $f$ are called the \demph{1-cells} of $\cat{B}$ and written 
  $
  \gfsts{A} \gones{f} \glsts{B}
  $
  and whose arrows $\gamma$ are called the \demph{2-cells} of $\cat{B}$ and
  written 
  \[
  \gfst{A} \gtwo{f}{g}{\gamma} \glst{B}
  \]
\item for each $A, B, C \in \cat{B}_0$, a functor
\[
\cat{B}(B, C) \times \cat{B}(A, B) \go \cat{B}(A, C)
\]
(\demph{composition}), written 
\begin{eqnarray*}
(g, f)			&\goesto	&g \of f,	\\
(\delta, \gamma)	&\goesto	&\delta * \gamma	
\end{eqnarray*}%
\glo{bicatstar}%
on 1-cells $f, g$ and 2-cells $\gamma, \delta$
\item for each $A \in \cat{B}_0$, an object $1_A \in \cat{B}(A, A)$ (the
  \demph{identity} on $A$)
\item for each triple 
  $
  \gfsts{A}\gones{f}\gblws{B}\gones{g}\gblws{C}\gones{h}\glsts{D}
  $
  of 1-cells, an isomorphism
  \[
  \gfst{A}
  \gtwo{(h\of g) \of f}{h \of (g\of f)}{\!\!\!\!\!\alpha_{h, g, f}}
  \glst{D}
  \]%
\glo{bicatass}%
in $\cat{B}(A, D)$ (the \demph{associativity
coherence isomorphism})%
\index{coherence!isomorphism}%
\index{associativity!isomorphism}
\item for each 1-cell $\gfsts{A}\gones{f}\glsts{B}$, isomorphisms
  \[
  \gfst{A}
  \gtwo{1_B \of f}{f}{\lambda_f}
  \glst{B},
  \diagspace
  \gfst{A}
  \gtwo{f \of 1_A}{f}{\rho_f}
  \glst{B}
  \]%
\glo{bicatlambda}\glo{bicatrho}%
in $\cat{B}(A, B)$ (the \demph{unit coherence isomorphisms})%
\index{unit!isomorphism}
\end{itemize}
such that the coherence isomorphisms $\alpha_{f, g, h}$, $\lambda_f$, and
$\rho_f$ are natural in $f$, $g$, and $h$ and satisfy pentagon and triangle
axioms like those in~\ref{defn:mon-cat} (replacing $((A\otimes B) \otimes
C) \otimes D$ with $((k\of h) \of g) \of f$, etc.).
\end{defn}

Functoriality of composition encodes `interchange%
\index{interchange}
laws', just as for
monoidal categories (\bref{eq:mon-interchange},
p.~\pageref{eq:mon-interchange}): the two evident derived composites of a
diagram of shape
\[
\gfstsu
\gthreesu 
\gzersu
\gthreesu 
\glstsu
\]
are equal, and similarly for 2-cells formed from diagrams
$
\gfstsu
\gonesu %{f}
\gzersu
\gonesu %{g}
\glstsu.
$

\begin{example}	\lbl{eg:bicat-mon-cat}
A bicategory%
\index{bicategory!degenerate}
with only one object is just a monoidal category.
\end{example}

\begin{example}
Strict 2-categories can be identified with bicategories in which all the
components of $\alpha$, $\lambda$ and $\rho$ are identities.
\end{example}

\begin{example}	\lbl{eg:bicat-Pi}
Any topological space $S$ gives rise to a bicategory $\Pi_2 S$,%
\glo{Pitwobicat}
the
\demph{fundamental%
\index{fundamental!2-groupoid}
2-groupoid} of $S$.  (`2-groupoid' means that all the
2-cells are isomorphisms and all the 1-cells are equivalences, in the sense
defined below.)  The objects of $\Pi_2 S$ are the points of $S$.  The
1-cells $a \go b$ are the paths%
\index{path}
from $a$ to $b$, that is, the maps $f: [0,1] \go S$ satisfying $f(0) = a$
and $f(1) = b$.  The 2-cells are the homotopy classes of path homotopies,
relative to endpoints.  Any particular point $s \in S$ determines a
one-object full sub-bicategory of $\Pi_2 S$ (that is, the sub-bicategory
whose only object is $s$ and with all possible 1- and 2-cells), and this is
the monoidal category of~\ref{eg:mon-cat-loops}.

Notice, incidentally, that a chain%
\index{chain complex!topological space@\vs.\ topological space}%
\index{topological space!chain complex@\vs.\ chain complex}
complex gives rise to a strict
$n$-category for each $n$
(\ref{eg:str-n-cat-ch-cx},~\ref{eg:str-omega-ch-cx}) but a space gives rise
to only \emph{weak}%
\index{n-category@$n$-category!weak vs. strict@weak \vs.\ strict}
structures.  This reflects the difference in difficulty
between homology%
\index{homology}
and homotopy.%
\index{homotopy!homology@\vs.\ homology}

The fundamental $\omega$-groupoid of a space is constructed
in~\ref{eg:wk-omega-cat-Pi}.  
\end{example}

\begin{example}	\lbl{eg:bicat-Ring}%
\index{bicategory!rings@of rings}
There is a bicategory $\cat{B}$ in which objects are rings, 1-cells $A \go
B$ are $(B, A)$-bimodules, and 2-cells are maps of bimodules.  (A
\demph{$(B, A)$-bimodule}%
\index{module!bimodule over rings}
is an abelian group $M$ equipped with a left
$B$-module structure and a right $A$-module structure satisfying 
$
(b \cdot m) \cdot a = b \cdot (m \cdot a)
$
for all $b\in B$, $m\in M$, $a\in A$.)  Composition is tensor:%
\index{tensor!ring@over ring}
if $M$ is a
$(B, A)$-bimodule and $N$ a $(C, B)$-bimodule then $N \otimes_B M$ is a
$(C, A)$-bimodule. 
\end{example}

\begin{example}
An $n$-category has $2^n$ duals,%
\index{dual!n-category@of $n$-category}%
\index{dual!bicategory@of bicategory}%
\index{bicategory!dual}
including the original article.  For a
bicategory $\cat{B}$, the dual obtained by reversing the 1-cells is
traditionally called $\cat{B}^\op$%
\glo{bicatop}
and that obtained by reversing the
2-cells is called $\cat{B}^{\mr{co}}$.%
\glo{bicatco}
 So the following four pictures show
corresponding 2-cells in $\cat{B}$, $\cat{B}^\op$, $\cat{B}^{\mr{co}}$ and
$\cat{B}^{\mr{co}\,\op} = \cat{B}^{\op\,\mr{co}}$:
\[
\gfsts{A}\gtwos{f}{g}{\gamma}\glsts{B},
\diagspace
\gfsts{A}\gtwoops{f}{g}{\gamma}\glsts{B},
\diagspace
\gfsts{A}\gtwocos{f}{g}{\gamma}\glsts{B},
\diagspace
\gfsts{A}\gtwocoops{f}{g}{\gamma}\glsts{B}.
\]
\end{example}

Since $\Cat$ is a bicategory~(\ref{eg:str-2-cat-Cat}), we may take
definitions from category theory and try to imitate them in an arbitrary
bicategory $\cat{B}$.  For example, a \demph{monad}%
\lbl{p:defn-monad-in-bicaty}\index{monad!bicategory@in bicategory}
in $\cat{B}$ is an object $A$ of $\cat{B}$ together with a 1-cell $A
\goby{t} A$ and 2-cells
\[
\mu: t \of t \go t, 
\diagspace
\eta: 1_A \go t,
\]
rendering commutative the diagrams of~\ref{defn:monad} (with $T$'s changed
to $t$'s, $T\mu$ changed to $1_t * \mu$, and so on).  An
\demph{adjunction}%
\index{adjunction!bicategory@in bicategory}
in $\cat{B}$ is a pair $(A, B)$ of objects together with 1- and 2-cells
\begin{equation}	\label{eq:bicat-adjn-data}
A \goby{f} B,
\diagspace
B \goby{g} A,
\diagspace
1_A \goby{\eta} g \of f,
\diagspace
f\of g \goby{\epsln} 1_B
\end{equation}
satisfying the triangle identities,%
\index{triangle!identities}
\[
(\epsln * 1_f) \of (1_f * \eta) = 1_f,
\diagspace
(1_g * \epsln) \of (1_g * \eta) = 1_g.
\]%
\index{equivalence!bicategory@in bicategory|(}%
An \demph{equivalence} between objects $A$ and $B$ of $\cat{B}$ is a
quadruple $(f, g, \eta, \epsln)$ of cells as in~\bref{eq:bicat-adjn-data}
such that $\eta$ and $\epsln$ are isomorphisms (in their respective
hom-categories).  An \demph{adjoint equivalence}%
\index{adjoint equivalence!bicategory@in bicategory}
is a quadruple that is
both an adjunction and an equivalence.  A 1-cell $f$ in a bicategory is
called an \demph{equivalence} if it satisfies the equivalent conditions of
the following result, which generalizes Proposition~\ref{propn:eqv-eqv}%
(\ref{item:eqv-eqv-adjt-eqv}--\ref{item:eqv-eqv-eqv}):
\begin{propn}	\lbl{propn:bicat-eqv-eqv}
Let $\cat{B}$ be a bicategory.  The following conditions on a 1-cell $f: A
\go B$ in $\cat{B}$ are equivalent:
\begin{enumerate}
\item %\lbl{item:bicat-eqv-eqv-adjt-eqv}
there exist $g$, $\eta$ and $\epsln$ such that $(f, g, \eta, \epsln)$ is an
adjoint equivalence 
\item %\lbl{item:bicat-eqv-eqv-eqv}
there exist $g$, $\eta$ and $\epsln$ such that $(f, g, \eta, \epsln)$ is an
equivalence. 
\end{enumerate}
\end{propn}
\begin{proof}
More is true: given an equivalence $(f, g, \eta, \epsln)$, there is a
unique $\epsln'$ such that $(f, g, \eta, \epsln')$ is an adjoint
equivalence, namely
\[
f \of g
\goby{\epsln^{-1} * 1_{f\sof g}}
f \of g \of f \of g
\goby{1_f * \eta^{-1} * 1_g}
f \of g
\goby{\epsln}
1_B.
\]
Brackets have been omitted, as if $\cat{B}$ were a strict 2-category; the
conscientious reader can both fill in the missing coherence isomorphisms
and verify that $\eta$ and $\epsln'$ satisfy the triangle identities (a
long but elementary exercise).
\done
\end{proof}
We write $A \eqv B$%
\glo{eqvinbicat}
if there exists an equivalence $A \go B$.%
\index{equivalence!bicategory@in bicategory|)}

Strict 2-categories form a strict 3-category~(\ref{eg:str-2-cat-Cat}),%
\index{n-category@$n$-category!n-ZZZcategory of@$(n+1)$-category of}
so
we would expect there to be notions of functor between bicategories and
transformation between functors, and then a further notion of map between
transformations.  On the other hand, one-object%
\index{bicategory!degenerate}
bicategories---monoidal
categories---only form a 2-category: there are various notions of functor
between monoidal categories, and a notion of transformation between
monoidal functors, but nothing more after that.  We will soon see how this
apparent paradox is resolved.

For functors everything goes smoothly. 
\begin{defn}
Let $\cat{B}$ and $\cat{B'}$ be bicategories.  A \demph{lax functor}%
\index{functor!classical bicategories@of classical bicategories}%
\index{bicategory!classical!functor of}
$F = (F, \phi): \cat{B} \go \cat{B'}$ consists of
\begin{itemize}
  \item a function $F_0: \cat{B}_0 \go \cat{B}'_0$, usually just written
  $F$
  \item for each $A, B \in \cat{B}_0$, a functor $F_{A, B}: \cat{B}(A, B)
  \go \cat{B'}(FA, FB)$, usually also written $F$
  \item for each composable pair $(f, g)$ of 1-cells in $\cat{B}$, a 2-cell
  $\phi_{g, f}: Fg \of Ff \go F(g\of f)$ 
  \item for each $A \in \cat{B}_0$, a 2-cell $\phi_A: 1_{FA} \go
  F1_A$ 
\end{itemize}
satisfying naturality and coherence axioms analogous to those of
Definition~\ref{defn:mon-ftr}.  \demph{Colax}, \demph{weak} and
\demph{strict functors} are also defined analogously. 
\end{defn}

\begin{example}	\lbl{eg:mon-cat-bicat-ftr}
Let $\cat{A}$ and $\cat{A'}$ be monoidal categories, and $\Sigma\cat{A}$%
\glo{SigmaMCbi}%
\index{suspension!monoidal category@of monoidal category}
and $\Sigma\cat{A'}$ the corresponding one-object bicategories.  Then lax
monoidal functors $\cat{A} \go \cat{A'}$ are exactly%
\index{bicategory!degenerate}
lax functors
$\Sigma\cat{A} \go \Sigma\cat{A'}$, and similarly for colax, weak and
strict functors. 
\end{example}

The evident composition of lax functors between bicategories is, perhaps
surprisingly, strictly associative and unital.  This means that we have a
category $\Bilax$%
\glo{Bilax}
of bicategories and lax functors, and subcategories $\Biwk$ and $\Bistr$.
The same goes for colax functors, but we concentrate on the lax (and above
all, the weak) case.

\begin{defn}
Let $(F, \phi), (G, \psi): \cat{B} \go \cat{B'}$ be lax functors between
bicategories.  A \demph{lax transformation}%
\index{transformation!classical bicategories@of classical bicategories}%
\index{bicategory!classical!transformation of}
$\sigma: F \go G$ consists of
\begin{itemize}
  \item for each $A \in \cat{B}_0$, a 1-cell $\sigma_A: FA \go GA$
  \item for each 1-cell $f: A \go B$ in $\cat{B}$, a 2-cell
  \[
  \begin{diagram}
  FA			&\rTo^{Ff}		&FB		\\
  \dTo<{\sigma_A}	&\nent \tcs{\sigma_f}	&\dTo>{\sigma_B}\\
  GA			&\rTo_{Gf}		&GB		\\
  \end{diagram}
  \]	
\end{itemize}
such that $\sigma_f$ is natural in $f$ and satisfies coherence axioms as in
Street~\cite[p.~568]{StrCS} or Leinster~\cite[1.2]{BB}.  A \demph{colax
transformation} is the same but with the direction of $\sigma_f$ reversed.
\demph{Weak} and \demph{strict transformations} are lax transformations in
which all the $\sigma_f$'s are, respectively, isomorphisms or identities.
\end{defn}

\begin{example}	\lbl{eg:mon-cat-bicat-transf}
Let $F, G: \cat{A} \go \cat{A'}$ be lax monoidal functors between monoidal
categories.  As in~\ref{eg:mon-cat-bicat-ftr}, these correspond%
\index{bicategory!degenerate}
to lax
functors $\Sigma F, \Sigma G: \Sigma\cat{A} \go \Sigma\cat{A'}$ between
bicategories.  A lax transformation $\Sigma F \go \Sigma G$ is an object
$S$ of $\cat{A'}$ together with a map
\[
\sigma_X: GX \otimes S \go S \otimes FX
\]
for each object $X$ of $\cat{A}$, satisfying coherence axioms.  (We are
writing $S = \sigma_\star$, where $\star$ is the unique object of
$\Sigma\cat{A}$.)

Transformations of one-object bicategories are, therefore, more general
than transformations of monoidal categories.  Monoidal transformations can
be identified with colax bicategorical transformations $\sigma$ for which
$\sigma_\star = I$.  Even putting aside the reversal of direction, we see
that what seems appropriate for bicategories in general seems inappropriate
in the one-object case.  We find out more in~\ref{cor:mon-eqv-bieqv}
and~\ref{eg:coh-mon-bi} below.
\end{example}

\begin{defn}
Let 
\[
\cat{B}
\ctwopar{F}{G}{\sigma}{\twid{\sigma}}
\cat{B'}
\]
be lax transformations between lax functors between bicategories.  A
\demph{modification}%
\index{modification!classical}\index{bicategory!classical!modification of}
$\Gamma: \sigma \go \twid{\sigma}$ consists of a 2-cell $\Gamma_A: \sigma_A
\go \twid{\sigma}_A$ for each $A \in \cat{B}_0$, such that $\Gamma_A$ is
natural in $A$ and satisfies a coherence axiom (Street~\cite[p.~569]{StrCS}
or Leinster~\cite[1.3]{BB}).
\end{defn}

Functors, transformations, and modifications can be composed in various
ways.  Let us consider just weak functors and transformations from now on.
It is straightforward to show that for any two bicategories $\cat{B}$ and
$\cat{B'}$, there is a functor bicategory $\ftrcat{\cat{B}}{\cat{B'}}$%
\glo{ftrbicat}\index{functor!bicategory}
whose objects are the weak functors from $\cat{B}$ to $\cat{B'}$, whose
1-cells are weak transformations, and whose 2-cells are modifications.
This is not usually a strict 2-category, because composing 1-cells in
$\ftrcat{\cat{B}}{\cat{B'}}$ involves composing 1-cells in $\cat{B'}$, but
it is a strict 2-category if $\cat{B'}$ is.  In particular, there is a
strict 2-category $\ftrcat{\cat{B}^\op}{\Cat}$ of `presheaves' on any
bicategory $\cat{B}$.

We might expect bicategories to form a weak 3-category, if we knew what one
was.  Gordon,%
\index{Gordon, Robert}
Power and Street gave an explicit definition of tricategory%
\index{tricategory}
(weak 3-category) in~\cite{GPS}, and showed that bicategories form a
tricategory.  It is worth noting, however, that there are two equally
sensible and symmetrically%
\index{handedness}
opposite ways of making a bicategory into a
tricategory, because there are two ways of horizontally composing weak
transformations, as is easily verified.  

\index{coherence!bicategories@for bicategories!classical|(}%
There are coherence theorems for bicategories analogous to those for
monoidal categories.  `All diagrams commute' is handled in exactly the same
way.  For the other statement of coherence we need the right notion of
equivalence of bicategories.  

Let $(F, \phi): \cat{B} \go \cat{B'}$ be a weak functor between
bicategories.  We call $F$ a \demph{local%
\index{local equivalence}\index{equivalence!local}
equivalence} if for each $A, B
\in \cat{B}_0$, the functor $F_{A, B}: \cat{B}(A, B) \go \cat{B}'(FA, FB)$
is an equivalence of categories, and \demph{essentially%
\index{essentially surjective on objects}
surjective on
objects} if for each $A' \in \cat{B}'_0$, there exists $A \in \cat{B}_0$
such that $FA \eqv A'$.  We call $F$ a \demph{biequivalence}%
\index{biequivalence!classical}
if it
satisfies the equivalent conditions of the following result, analogous to
Proposition~\ref{propn:eqv-eqv}%
\bref{item:eqv-eqv-eqv}--\bref{item:eqv-eqv-ffe}:
\begin{propn}	\lbl{propn:bieqv-eqv}
The following conditions on a weak functor $F: \cat{B} \go \cat{B'}$
between bicategories are equivalent:
\begin{enumerate}
\item	\lbl{item:bieqv-eqv-bieqv}
  there exists a weak functor $G: \cat{B'} \go \cat{B}$ such that
  $1_\cat{B} \eqv G \of F$ in $\ftrcat{\cat{B}}{\cat{B}}$ and 
  $F \of G \eqv 1_\cat{B'}$ in $\ftrcat{\cat{B'}}{\cat{B'}}$
\item	\lbl{item:bieqv-eqv-lee}
  $F$ is a local equivalence and essentially surjective on objects.
\end{enumerate}
\end{propn}
\begin{prooflike}{Sketch proof}
\bref{item:bieqv-eqv-bieqv} $\Rightarrow$ \bref{item:bieqv-eqv-lee} is
straightforward.  For the converse, choose for each $A' \in \cat{B'}$ an
object $GA' \in \cat{B}$ together with an \emph{adjoint} equivalence
between $FGA'$ and $A'$, which is possible by~\ref{propn:bicat-eqv-eqv}.
Then the remaining constructions and checks are straightforward, if
tedious.  \done
\end{prooflike}

Condition~\bref{item:bieqv-eqv-bieqv} of the Proposition says that there is
a system of functors, \emph{transformations and modifications} relating
$\cat{B}$ and $\cat{B'}$.  It therefore seems quite plausible that in the
one-object case, biequivalence is a looser relation than monoidal
equivalence.  Nevertheless,
\begin{cor}	\lbl{cor:mon-eqv-bieqv}%
\index{bicategory!degenerate}
Two monoidal categories are monoidally equivalent if and only if the
corresponding one-object bicategories are biequivalent. 
\end{cor}
\begin{proof}
The monoidal and bicategorical notions of weak functor are the
same~(\ref{eg:mon-cat-bicat-ftr}), so this follows from
conditions~\bref{item:bieqv-eqv-lee} of
Propositions~\ref{propn:mon-eqv-eqv} and~\ref{propn:bieqv-eqv}. 
\done
\end{proof}

\begin{thm}[Coherence for bicategories]	\lbl{thm:coh-bi-wk-str}
Every bicategory is biequivalent to some strict 2-category.
\end{thm}
\begin{prooflike}{Sketch proof}
Let $\cat{B}$ be a bicategory.  There is a weak functor
\[
\mathbf{y}: \cat{B} \go \ftrcat{\cat{B}^\op}{\Cat},
\]
the analogue of the Yoneda%
\index{Yoneda!embedding}
embedding for categories, sending an object $A$
of $\cat{B}$ to the weak functor
\[
\cat{B}(\dashbk, A): \cat{B}^\op \go \Cat
\]
and so on.  Just as the ordinary Yoneda embedding is full and faithful,
$\mathbf{y}$ is a local equivalence.  So if $\cat{B'}$ is the
sub-strict-2-category of $\ftrcat{\cat{B}^\op}{\Cat}$ consisting of all the
objects in the image of $\mathbf{y}$ and all the 1- and 2-cells between
them then $\mathbf{y}$ defines a biequivalence from $\cat{B}$ to
$\cat{B'}$.  
\done
\end{prooflike}

\begin{example}	\lbl{eg:coh-mon-bi}%
\index{bicategory!degenerate}
Corollary~\ref{cor:mon-eqv-bieqv} enables us to deduce coherence for
monoidal categories~(\ref{thm:coh-mon-wk-str}) from coherence for
bicategories.  In fact, the proof of~\ref{thm:coh-mon-wk-str} is exactly
the proof of~\ref{thm:coh-bi-wk-str} in the one-object case.  Let $\cat{B}$
be a bicategory with single object $\star$ and let $\cat{A}$ be the
corresponding monoidal category.  Then $\cat{B'}$ is the sub-2-category of
$\ftrcat{\cat{B}^\op}{\Cat}$ whose single object is the representable
functor $\cat{B}(\dashbk, \star)$, and a weak transformation from
$\cat{B}(\dashbk, \star)$ to itself is an object $(E, \delta)$ of
$\cat{A'}$ as in the proof of~\ref{thm:coh-mon-wk-str}, and so on.
\end{example}%
\index{coherence!bicategories@for bicategories!classical|)}%

\begin{notes}

Almost everything in this chapter is very well-known (in the sense of the
phrase to which mathematicians are accustomed).  I have used Mac
Lane~\cite{MacCWM} as my main reference for ordinary category theory, with
Borceux~\cite{Borx1, Borx2} as backup.  The standard text on enriched
category theory is Kelly~\cite{KelBCE}.  For higher category theory,
Street~\cite{StrCS} provides a useful survey and reference list.  Strict
$n$-categories were introduced by Ehresmann~\cite{Ehr}.%
\index{Ehresmann, Charles}

I thank Bill Fulton and Ross Street for an enlightening exchange on
equivalence of bicategories, and Nathalie Wahl for a useful conversation on
how not%
\index{coherence!how not to prove}
to prove coherence.

\end{notes}

\chapter{Classical Operads and Multicategories}
\lbl{ch:om}

\chapterquote{%
Some `pictures' are not really pictures, but rather are windows to Plato's%
\index{Plato}
heaven}{%
Brown~\cite{BrownJ}}

\noindent
Where category theory has arrows, higher-dimensional category theory has
higher-dimensional arrows.  One of the simplest examples of a
`higher-dimensional arrow' is one like
\[
\begin{centredpic}
\begin{picture}(6,4)(-1,-2)
\cell{0}{0}{l}{\tusual{}}
\cell{0}{0}{r}{\tinputsslft{}{}{}}
\cell{4}{0}{l}{\toutputrgt{}}
\end{picture}
\end{centredpic}
.
\]
Think of this as a box with $n$ input wires coming in on the left, where
$n$ is any natural number, and one output wire emerging on the right.  (For
instance, when $n=1$ this is just an arrow as in an ordinary category.)
With this in mind, it is easy to imagine what composition of such arrows
might look like: outputs of one arrow attach to inputs of another.

A categorical structure with arrows like this is called a multicategory.
(Multicategories and $n$-categories are not the same!)  A very familiar
example: the objects (drawn as labels on wires) are vector spaces, and the
arrows are multilinear maps.  The special case of a multicategory where
there is only one object---that is, the wires are unlabelled---is
particularly interesting.  Such a structure is called an operad; for a
basic example, fix a topological space $X$ and define an arrow with $n$
inputs to be a continuous map $X^n \go X$.

There is a curiously widespread impression that operads are frighteningly
complicated structures.  Among users of operads, there is a curiously
widespread impression that multicategories---usually known to them as
`coloured operads'---are some obscure and esoteric elaboration of the basic
notion of operad.  I hope this chapter will correct both impressions.  Both
structures are as natural as can be, and, if one draws some pictures, very
simple to understand.

We start with multicategories~(\ref{sec:cl-mtis}) then specialize to
operads~(\ref{sec:cl-opds}).  (This order of presentation may convince
sceptical operad-theorists that multicategories are natural structures in
their own right.)  The basic definitions are given, with a broad range of
examples.  An assortment of further topics on operads and multicategories
is covered in~\ref{sec:om-further}.

\paragraph*{Warning}%
\lbl{p:sym-warning}\index{operad!usage of word}\index{operad!symmetric}  
Readers already familiar with operads may be used to them coming equipped
with symmetric group actions.  In this text operads \emph{without}
symmetries are the default.  This is partly to fit with the convention that
monoidal categories are by default non-symmetric, and rings, groups and
monoids non-commutative, but is mostly for reasons that will emerge later.
So:
\begin{quote}\centering\bf
`Operad' means what is sometimes called `non-$\Sigma$ operad' or
`non-symmetric operad'.
\end{quote}
Operads equipped with symmetries will be called `symmetric operads'.

\section{Classical multicategories}
\lbl{sec:cl-mtis}%

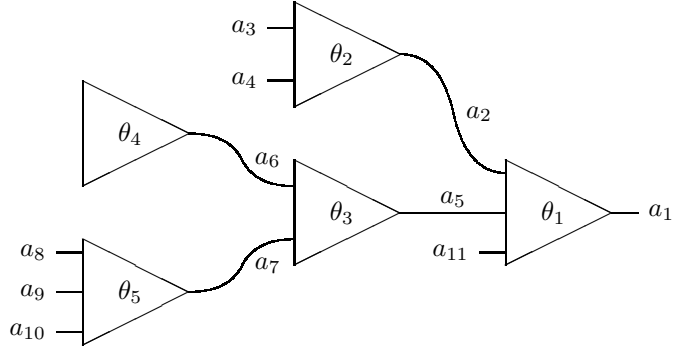
\begin{figure}	
% \[
\begin{center}
\setlength{\unitlength}{1em}
\begin{picture}(24,13)(0,-5)
% transistors
\cell{18}{0}{l}{\tusual{\theta_1}}
\cell{10}{6}{l}{\tusual{\theta_2}}
\cell{10}{0}{l}{\tusual{\theta_3}}
\cell{2}{3}{l}{\tusual{\theta_4}}
\cell{2}{-3}{l}{\tusual{\theta_5}}
% short wires
\cell{22}{0}{l}{\toutputrgt{a_1}}
\cell{10}{7}{r}{\tinputlft{a_3}}
\cell{10}{5}{r}{\tinputlft{a_4}}
\cell{2}{-1.5}{r}{\tinputlft{a_8}}
\cell{2}{-3}{r}{\tinputlft{a_9}}
\cell{2}{-4.5}{r}{\tinputlft{a_{10}}}
\cell{18}{-1.5}{r}{\tinputlft{a_{11}}}
% long wires
\qbezier(18,1.5)(16.5,1.5)(16,3.75)
\qbezier(14,6)(15.5,6)(16,3.75)
\cell{16.5}{3.75}{l}{a_2}
\put(18,0){\line(-1,0){4}}
\cell{16}{0.2}{b}{a_5}
\qbezier(10,1)(8.5,1)(8,2)
\qbezier(6,3)(7.5,3)(8,2)
\cell{8.5}{2}{l}{a_6}
\qbezier(10,-1)(8.5,-1)(8,-2)
\qbezier(6,-3)(7.5,-3)(8,-2)
\cell{8.5}{-2}{l}{a_7}
\end{picture}
\end{center}
% \hand{45}{50}
% \]
\caption{Composable diagram of arrows in a multicategory}
\label{fig:random-multi-diagram}
\end{figure}

Let us make the description above more precise.  
A category consists of objects, arrows between objects, and a way of
composing arrows.  Precisely the same description applies to
multicategories.  The only difference 
% between the two structures
lies in the shape of the arrows: in a category, an arrow looks like
\[
a \goby{\theta} b,
\]
with one object as its domain and one object as its codomain, whereas in a
multicategory, an arrow looks like
\begin{equation}	\label{diag:multi-trans-arrow}
\begin{centredpic}
\begin{picture}(8,4)(-2,-2)
\cell{0}{0}{l}{\tusual{\theta}}
\cell{0}{0}{r}{\tinputsslft{a_1}{a_2}{a_n}}
\cell{4}{0}{l}{\toutputrgt{a}}
\end{picture}
\end{centredpic}
\end{equation}
($n\in\nat$), with a finite sequence of `input' objects as its domain and one
`output' object as its codomain.  
% Put another way, an arrow in a multicategory has
% several inputs and a single output.  
Arrows can be composed when outputs
are joined to inputs, which for categories means that any string of arrows
\[
a_0 \goby{\theta_1} a_1 \goby{\theta_2} 
\diagspace \cdots \diagspace
\goby{\theta_n} a_n
\]
has a well-defined composite, and for multicategories means (more
interestingly in geometrical terms) that any tree of arrows such as that in
Fig.~\ref{fig:random-multi-diagram}
has a well-defined composite, in this case of the form
\[
\begin{centredpic}
\begin{picture}(10,6)(-2,-3)
% transistor
\cell{0}{0}{l}{\tmid{}}
% tags
\cell{0}{2.5}{r}{\tinputlft{a_3}}
\cell{0}{1.5}{r}{\tinputlft{a_4}}
\cell{0}{0.5}{r}{\tinputlft{a_8}}
\cell{0}{-0.5}{r}{\tinputlft{a_9}}
\cell{0}{-1.5}{r}{\tinputlft{a_{10}}}
\cell{0}{-2.5}{r}{\tinputlft{a_{11}}}
\cell{6}{0}{l}{\toutputrgt{a_1.}}
\end{picture}
\end{centredpic}
\]
Perhaps the most familiar example is where the objects are vector spaces
and the arrows are multilinear maps.  Commonly the multicategory structure
is obscured by the device of considering multilinear maps as linear maps
out of a tensor product; but in many situations it is the multicategory,
not the monoidal category, that is fundamental.

\begin{defn}	\lbl{defn:cl-multi}
A \demph{multicategory}%
\index{multicategory}
$C$ consists of
\begin{itemize}
\item a class $C_0$,%
\glo{C0multicat}
whose elements are called the \demph{objects} of $C$
\item for each $n\in\nat$ and $a_1, \ldots, a_n, a \in C_0$, a class
$C(a_1, \ldots, a_n; a)$,%
\glo{multicathomset}
whose elements $\theta$ are called \demph{arrows}
or \demph{maps} and depicted as in~\bref{diag:multi-trans-arrow} or as
\begin{equation}	\label{eq:multi-arrow}
a_1, \ldots, a_n \goby{\theta} a
\end{equation}
\item for each $n, k_1, \ldots, k_n \in \nat$ and $a, a_i, a_i^j \in C_0$,
a function (Fig.~\ref{fig:multi-comp})
\[
\begin{array}{r}
C(a_1, \ldots, a_n; a) \times
C(a_1^1, \ldots, a_1^{k_1}; a_1) \times \cdots \times
C(a_n^1, \ldots, a_n^{k_n}; a_n)\\
\go
C(a_1^1, \ldots, a_1^{k_1}, \ldots, a_n^1, \ldots, a_n^{k_n}; a),
\end{array}
\]
called \demph{composition} and written
\[
(\theta, \theta_1, \ldots, \theta_n) \goesto 
\theta \of (\theta_1, \ldots, \theta_n)
\glo{multicomposite}%
\]%
\item for each $a\in C_0$, an element $1_a \in C(a;a)$,%
\glo{multiidentity}
called the
\demph{identity} on $a$
\end{itemize}
satisfying
\begin{itemize}
\item associativity:
\begin{eqnarray*}
&&\theta \of 
\left(
\theta_1 \of (\theta_1^1, \ldots, \theta_1^{k_1}),
\ldots,
\theta_n \of (\theta_n^1, \ldots, \theta_n^{k_n})
\right)\\
&=&
\left(\theta \of (\theta_1, \ldots, \theta_n)\right) \of
(\theta_1^1, \ldots, \theta_1^{k_1}, 
\ldots, 
\theta_n^1, \ldots, \theta_n^{k_n})
\end{eqnarray*}
whenever $\theta, \theta_i, \theta_i^j$ are arrows for which these
composites make sense
\item identity:
\[
\theta \of (1_{a_1}, \ldots, 1_{a_n}) 
= 
\theta 
= 
1_a \of (\theta)
\]
whenever $\theta: a_1, \ldots, a_n \go a$ is an arrow.
\end{itemize}
\end{defn}
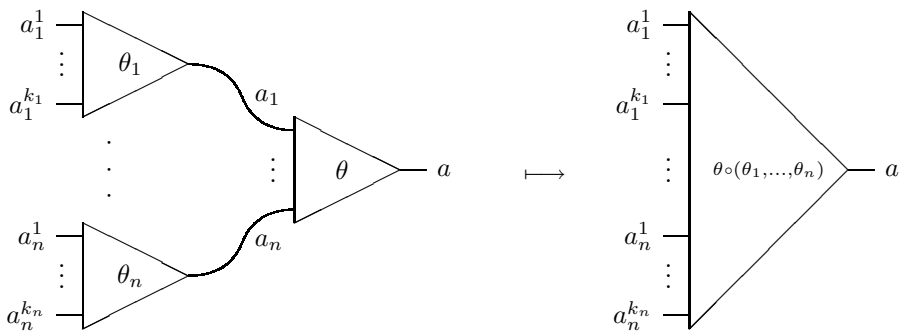
\begin{figure}
\[
% DOMAIN
%
\begin{centredpic}
\begin{picture}(16.5,12)(-0.5,-6)
% transistors
\cell{10}{0}{l}{\tusual{\theta}}
\cell{2}{4}{l}{\tusual{\theta_1}}
\cell{2}{-4}{l}{\tusual{\theta_n}}
% short wires
\cell{14}{0}{l}{\toutputrgt{a}}
\cell{2}{4}{r}{\tinputslft{a_1^1}{a_1^{k_1}}}
\cell{2}{-4}{r}{\tinputslft{a_n^1}{a_n^{k_n}}}
% long wires
\qbezier(10,1.5)(8.5,1.5)(8,2.75)
\qbezier(6,4)(7.5,4)(8,2.75)
\cell{8.5}{2.75}{l}{a_1}
\qbezier(10,-1.5)(8.5,-1.5)(8,-2.75)
\qbezier(6,-4)(7.5,-4)(8,-2.75)
\cell{8.5}{-2.75}{l}{a_n}
% ellipses
\cell{9.2}{0.3}{c}{\vdots}
\cell{3}{0}{c}{\cdot}
\cell{3}{1}{c}{\cdot}
\cell{3}{-1}{c}{\cdot}
\end{picture}
\end{centredpic}
\mbox{\hspace{2em}}
\goesto
\mbox{\hspace{2em}}
%
% CODOMAIN
%
\begin{centredpic}
\begin{picture}(10,12)(0,-6)
% transistor
\put(2,-6){\line(0,1){12}}
\put(8,0){\line(-1,-1){6}}
\put(8,0){\line(-1,1){6}}
\cell{5}{0}{c}{\scriptstyle\theta \sof (\theta_1, \ldots, \theta_n)}
% \cell{5}{0}{c}{\theta \circ (\theta_1, \ldots, \theta_n)}
% tags
\cell{2}{4}{r}{\tinputslft{a_1^1}{a_1^{k_1}}}
\cell{2}{-4}{r}{\tinputslft{a_n^1}{a_n^{k_n}}}
\cell{1.2}{0.3}{c}{\vdots}
\cell{8}{0}{l}{\toutputrgt{a}}
\end{picture}
\end{centredpic}
\]
% \hand{50}{51}
\caption{Composition in a multicategory}
\label{fig:multi-comp}
\end{figure}

Operads are precisely multicategories with only one object, and are the
subject of the next section.  For now we look at multicategories with many
objects---often a proper class of them.

\begin{example}	\lbl{eg:multi-unary}
A multicategory in which every arrow is \demph{unary}%
\index{unary}
(that is, of the
form~\bref{eq:multi-arrow} with $n=1$) is the same thing as a category. 
\end{example}

\begin{example}	\lbl{eg:multi-mon}%
Any monoidal category $(A, \otimes)$ has an \demph{underlying
multicategory}%
\index{multicategory!underlying}
 $C$.  It has the same objects as $A$, and a map
\[
a_1, \ldots, a_n \go a
\]
in $C$ is a map
\[
a_1 \otimes \cdots \otimes a_n \go a
\]
in $A$.  Composition in $C$ is derived from composition and tensor in $A$.

If $A$ is a non-strict monoidal category then there is some ambiguity in
the meaning of `$a_1 \otimes\cdots\otimes a_n$'.  We will address this
properly in~\ref{sec:non-alg-notions}, but for now let us just choose
a particular bracketing, e.g.\ $((a_1 \otimes a_2) \otimes a_3) \otimes
a_4$ for $n=4$.

Given a commutative ring $R$, the monoidal category of $R$-modules%
\index{module!ring@over ring}
with their usual tensor gives rise in this way to a multicategory of
$R$-modules and $R$-multilinear maps.  Similarly, given a category $A$ with
finite products there is an underlying multicategory $C$ with the same
objects as $A$ and with
\[
C(a_1, \ldots, a_n; a) = A(a_1 \times\cdots\times a_n, a).
\]
We will make particular use of the multicategory $C$ coming from the
category $A = \Set$;%
\index{multicategory!sets@of sets}
we write $C = \Set$ too.
\end{example}

\begin{example}	\lbl{eg:multi-some-of-mon}
More generally, if $A$ is a monoidal category and $C_0$ a collection of
objects of $A$ (not necessarily closed under tensor) then there is a
multicategory $C$ whose class of objects is $C_0$ and whose arrows are
defined as in the previous example. 
\end{example}

It is important to realize that \emph{not} every multicategory%
\lbl{p:not-ESO}
is the underlying%
\index{multicategory!underlying}
multicategory of a monoidal
category~(\ref{eg:multi-mon}).  Example~\ref{eg:multi-some-of-mon} makes
this reasonably clear.  Another way of seeing it is to consider one-object
multicategories (operads, next section): if $A$ is a monoidal category with
single object $a$ then in the underlying multicategory $C$, the cardinality
of the hom-set $C(a, \ldots, a; a)$ is independent of the number of $a$'s
in the domain.  But this independence property certainly does not hold for
all one-object multicategories, as numerous examples in
Section~\ref{sec:cl-opds} show.

\begin{example}%
\index{multicategory!underlying}
Sometimes tensor is irrelevant; sometimes it does not even exist.  For
instance, let $V$ be a symmetric monoidal category and let $\Ab(V)$ be the
category of abelian groups in $V$.  In general $V$ will not have enough
colimits to define a tensor product on $\Ab(V)$ and so make it into a
monoidal category.  But it is always possible to define multilinear maps
between abelian groups in $V$, giving $\Ab(V)$ the structure of a
multicategory.
\end{example}

\begin{example}	\lbl{eg:multi-pretend-coproducts}
If $A$ is a category with finite coproducts%
\index{operad!coproducts@from coproducts}
then Example~\ref{eg:multi-mon}
reveals that there is a multicategory $C$ with the same objects as $A$,
with
\[
C(a_1, \ldots, a_n; a)
=
% A(a_1 + \cdots + a_n, a)
% \iso
A(a_1, a) \times\cdots\times A(a_n, a),
\]
and with composition given by
\begin{eqnarray*}
&&
(f_1, \ldots, f_n) \of 
((f_1^1, \ldots, f_1^{k_1}), \ldots, (f_n^1, \ldots,  f_n^{k_n}))\\
&=&
(f_1 \of f_1^1, \ldots, f_1 \of f_1^{k_1}, \ldots,
f_n \of f_n^1, \ldots, f_n \of f_n^{k_n})
\end{eqnarray*}
($f_i \in A(a_i, a)$, $f_i^j \in A(a_i^j, a_i)$).  Indeed, these formulas
make sense and define a multicategory $C$ for any category $A$ whatsoever.
An arrow of $C$ can be drawn as
\[
\begin{centredpic}
\begin{picture}(10,6)(0,-3)
% transistor
\cell{2}{0}{l}{\tmid{}}
% tags
\cell{2}{2}{r}{\tinputlft{a_1}}
\cell{2}{1}{r}{\tinputlft{a_2}}
\cell{2}{-2}{r}{\tinputlft{a_n}}
\cell{8}{0}{l}{\toutputrgt{a.}}
% lines inside transistor...
\qbezier(2,2)(5,1)(8,0)
\qbezier(2,1)(5,0.5)(8,0)
\qbezier(2,-2)(5,-1)(8,0)
% ...and their arrowheads
\cell{3.5}{1.4}{br}{\rotatebox{35}{$\lrcorner$}}
\cell{3.5}{0.6}{br}{\rotatebox{30}{$\lrcorner$}}
\cell{3.5}{-1.3}{tr}{\rotatebox{70}{$\lrcorner$}}
% ...and their labels
\cell{3}{2.1}{c}{\scriptstyle f_1}
\cell{3}{0.4}{c}{\scriptstyle f_2}
% \cell{3}{-1.3}{c}{\scriptstyle f_n}
\cell{3}{-2}{c}{\scriptstyle f_n}
% ellipses
\cell{1.2}{-.3}{c}{\vdots}
\cell{4}{-0.1}{c}{\vdots}
\end{picture}
\end{centredpic}
% \hand{30}{52}.
\]
\end{example}

\begin{example}
A category in which each hom-set has at most one element is the same thing
as a preordered%
\index{order!generalized}%
\index{multicategory!poset-like}
set.  Similarly, a multicategory $C$ in which $C(a_1,
\ldots, a_n; a)$ has at most one element for each $n, a_1, \ldots, a_n, a$
is some kind of generalized poset.  It amounts to a class $X = C_0$ of
elements together with an $(n+1)$-ary relation $\leq_n$ on $X$
%  \sub X^n \times X$
for each $n\in\nat$, satisfying generalized reflexivity and transitivity
axioms.  For instance, let $d\in\nat$, let $X = \reals^d$, and define
$\leq_n$ by $(a_1, \ldots, a_n) \leq_n a$ if and only if $a$ is in the
convex%
\index{convex hull}
hull of $\{ a_1, \ldots, a_n \}$.  This gives a `generalized poset';
the axioms express basic facts about convex hulls.
\end{example}

The final example is more substantial.

\begin{example}	\lbl{eg:multi-cheese}
The \demph{Swiss%
\index{multicategory!Swiss cheese}
cheese multicategory} $\fcat{SC}$ has the natural numbers
as its objects, with $d\in\nat$ thought of as the $d$-dimensional disk.
Maps are configurations of disks, half-disks, quarter-disks, and so on, as
now described.

Let $\reals^\infty$ be the set of sequences $(x_n)_{n=1}^{\infty}$ of real
numbers.  For each $d\in\nat$, let
\[
\cheesydisk{d} = 
\{ \mathbf{x} \in \reals^\infty
\such
\sum_{n=1}^{\infty} x_n^2 \leq 1 
\textrm{, and } x_n \geq 0 \textrm{ for all } n>d \},
\]
and define a sequence of subsets $F^d_1, F^d_2, \ldots$ of $\cheesydisk{d}$
by
\[
F^d_n =
\left\{
\begin{array}{ll}
\emptyset	&\textrm{ if } n\leq d	\\
\{ \mathbf{x} \in \cheesydisk{d}
\such
x_n = 0 \}	&\textrm{ if } n>d.
\end{array}
\right.
\]
(The idea is that if the hyperplane $x_n = 0$ contains one of the flat
faces of $\cheesydisk{d}$ then $F^d_n$ is that flat face, and otherwise
$F^d_n$ is empty.)  Let $G$ be the group of bijections on $\reals^\infty$
of the form $\mathbf{x} \goesto \mathbf{a} + \lambda\mathbf{x}$ for some
$\mathbf{a} \in \reals^\infty$ and $\lambda > 0$.  Then there is a category
$A$ in which the objects are the natural numbers and a map $d' \go d$ is an
element $\alpha\in G$ such that
\[
\alpha\cheesydisk{d'} \sub \cheesydisk{d},
\diagspace
\alpha F^{d'}_n \sub F^d_n \textrm{ for all } n\geq 1.
\]
(So $\alpha$ sends $\cheesydisk{d'}$ into $\cheesydisk{d}$, with flat faces
mapping into flat faces.)  By Example~\ref{eg:multi-pretend-coproducts},
there arises a multicategory $C$ with object-set $\nat$ and with
\[
C(d_1, \ldots, d_k; d) = A(d_1, d) \times\cdots\times A(d_k, d).
\]
The Swiss cheese multicategory $\fcat{SC}$ is the sub-multicategory of $C$
consisting of the same objects but only those arrows
\[
(\alpha_1, \ldots, \alpha_k) \in C(d_1, \ldots, d_k; d)
\]
for which the images $\alpha_i \cheesydisk{d_i}$ and $\alpha_j
\cheesydisk{d_j}$ are disjoint whenever $i\neq j$.  

A few calculations reveal that $\fcat{SC}(d_1, \ldots, d_k; d)$ is only
non-empty when $d_i \geq d$ for each $i$, and that the description of this
hom-set is unchanged if we replace $\reals^\infty$ by $\reals^{\max\{d_1,
\ldots, d_k\} }$ throughout.  It follows that, for instance, a map
\[
1, 1, 2, 2, 2, 3 \go 1
\]
(Fig.~\ref{fig:cheese}) 
\begin{figure}
\begin{center}
\setlength{\unitlength}{1mm}
\begin{picture}(104,70)
% cheese.eps is 72mm wide and 70mm high
\cell{0}{0}{bl}{\epsfig{file=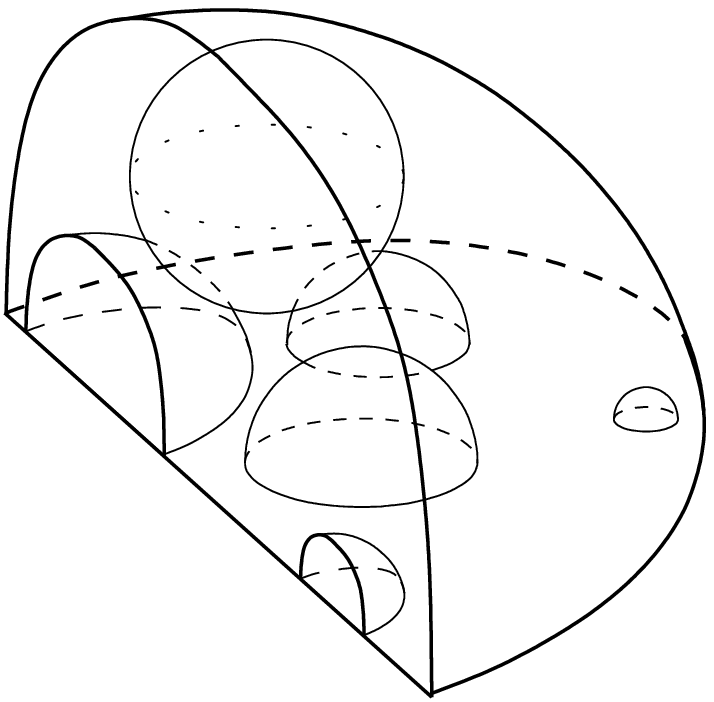}}
\cell{92}{35}{l}{\epsfig{file=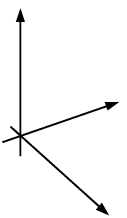}}
\cell{105}{23}{c}{x_1}
\cell{107}{36.5}{c}{x_2}
\cell{94}{47}{c}{x_3}
\end{picture}
\end{center}
% \hand{70}{53}
\caption{Swiss cheese map $1, 1, 2, 2, 2, 3 \protect\go 1$}
\label{fig:cheese}
\end{figure}
is a configuration of two quarter-balls, three half-balls, and one whole
ball inside the unit quarter-ball in $\reals^3$, with the (straight)
$1$-dimensional face of each quarter-ball lying on the $1$-dimensional face
of the unit quarter-ball and the (flat) $2$-dimensional face of each
half-ball lying on the 2-dimensional face $x_3=0$ of the unit quarter-ball.
The six little fractions of balls must be disjoint.

For each $d\in\nat$ there is a sub-multicategory $\fcat{SC}_d$ consisting
of only the two objects $d$ and $d+1$ (and all the maps between them).
This was Voronov's%
\index{Voronov, Alexander}
original `Swiss-cheese%
\index{operad!Swiss cheese}
operad' (a `2-coloured operad');
see Voronov~\cite{VorSCO} and Kontsevich~\cite[2.5]{KonOMD}.%
\index{Kontsevich, Maxim}
 The
one-object sub-multicategories of $\fcat{SC}$ are more famous: for each
$d\in\nat$, the sub-multicategory consisting of $d$ (and all maps) is
exactly the little%
\index{operad!little disks}
$d$-dimensional disks
operad~(\ref{eg:opd-little-disks}).
\end{example}

Many of the multicategories mentioned have a natural symmetric structure:
there is a bijection
\[
\dashbk\cdot\sigma:
C(a_1, \ldots, a_n; a)
\goiso
C(a_{\sigma(1)}, \ldots, a_{\sigma(n)}; a)
\]%
\glo{dotsigma}%
for each $a_1, \ldots, a_n, a \in C_0$ and permutation $\sigma\in S_n$.%
\glo{symgp}
 A
`symmetric%
\index{multicategory!symmetric}
multicategory'%
\lbl{p:sym-mti-informal}
is defined as a multicategory equipped with such a family of bijections,
satisfying axioms.  We give the precise definition in~\ref{defn:sym-mti};
see also Appendix~\ref{app:sym}.  For example, if $A$ is a symmetric
monoidal category then its underlying multicategory $C$ is also symmetric.

Similarly, many of these examples are also naturally `enriched':%
\lbl{p:sym-enr-mti}%
\index{enrichment!plain multicategory@of plain multicategory!symmetric monoidal category@in symmetric monoidal category}
$C(a_1, \ldots, a_n; a)$ is more than a mere set.  For instance, if $C$ is
vector spaces and multilinear maps then $C(a_1, \ldots, a_n; a)$ naturally
has the structure of a vector space itself, and in some of the oldest
examples of operads the `hom-sets' are topological spaces (as we shall
see).  One way to formalize this is to allow $C(a_1, \ldots, a_n; a)$ to be
an object of some chosen symmetric%
\index{monoidal category!symmetric!enrichment in}
monoidal category $\cat{V}$ (generalizing from $\cat{V} = \Set$).  This is
actually not the most natural generalization, essentially because the
tensor product on $\cat{V}$ is redundant; we come back to this
in~\ref{sec:enr-mtis}.

\begin{defn}
Let $C$ and $C'$ be multicategories.  A \demph{map%
\index{multicategory!map of}
of multicategories} $f:
C \go C'$ consists of a function $f_0: C_0 \go C'_0$ (usually just written
$f$) together with a function
\[
C(a_1, \ldots, a_n; a) \go C'(f(a_1), \ldots, f(a_n); f(a))
\]
(written $\theta \goesto f(\theta)$) for each $a_1, \ldots, a_n, a \in
C_0$, such that composition and identities are preserved.  The category
$\Multicat$%
\glo{Multicat}
consists of small multicategories and maps between them.
\end{defn}

\begin{example}	\lbl{eg:map-mti-mon}%
\index{multicategory!underlying}
Let $A$ and $A'$ be monoidal categories, with respective underlying
multicategories $C$ and $C'$.  Then a map $C \go C'$ of multicategories is
precisely a lax monoidal functor $A \go A'$.
\end{example}

\begin{example}
Monoids can be described as maps.  A \demph{monoid}%
\index{monoid!multicategory@in multicategory}
in a multicategory $C$
consists of an object $a$ of $C$ together with arrows
\[
a, a \goby{\mu} a,
\diagspace
\cdot \goby{\eta} a
\]
(where the domain of $\eta$ is the empty sequence) obeying associativity
and identity laws.  The terminal%
\index{multicategory!terminal}
multicategory $1$ consists of one object,
$\star$, say, and one arrow
\[
\underbrace{\star, \ldots, \star}_n \go \star
\]
for each $n\in\nat$.  A map from $1$ into a multicategory $C$ therefore
consists of an object $a$ of $C$ together with an arrow
\[
\underbrace{a, \ldots, a}_n \goby{\mu_n} a
\]
for each $n\in\nat$, obeying `all possible laws', and this is exactly a
monoid in $C$. 
\end{example}

In ordinary category theory, functors into $\Set$ play an important role.
The same goes for multicategories.
\begin{defn}	\lbl{defn:alg-multi}
Let $C$ be a multicategory.  An \demph{algebra for $C$},%
\index{algebra!plain multicategory@for plain multicategory}%
\index{multicategory!algebra for}
or
\demph{$C$-algebra}, is a map from $C$ into the multicategory $\Set$
of Example~\ref{eg:multi-mon}. 
\end{defn}
The name makes sense if a multicategory is regarded as an algebraic%
\index{algebraic theory!plain multicategory as}
theory
with as many sorts as there are objects; it also generalizes the standard
terminology for operads.  Explicitly, a $C$-algebra $X$ consists of
\begin{itemize}
\item for each object $a$ of $C$, a set $X(a)$, whose elements $x$ can be
drawn as
\[
\setlength{\unitlength}{1em}
\begin{picture}(6.5,4)(0,-2)
% transistor
\cell{0}{0}{l}{\twiggly{x}}
% tag
% tag
\cell{4.5}{0}{l}{\toutputrgt{a}}
\end{picture}
% \hand{20}{54}
\]
\item for each map $\theta: a_1, \ldots, a_n \go a$ in $C$, a function
\[
\ovln{\theta} = X(\theta): 
X(a_1) \times\cdots\times X(a_n)
\go
X(a)
\]%
\glo{thetabar}%
(Fig.~\ref{fig:alg-action}),
\end{itemize}
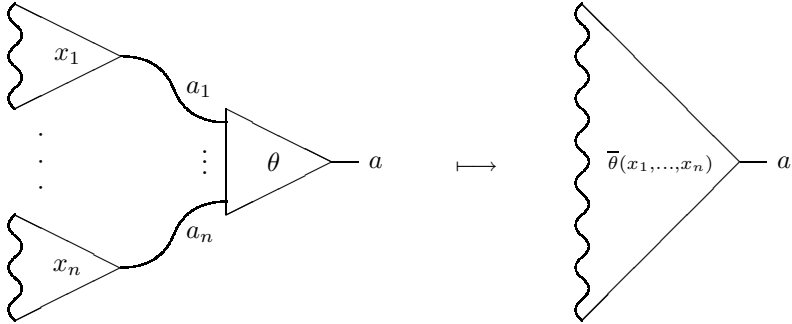
\begin{figure}
\[
% DOMAIN
%
\begin{centredpic}
\begin{picture}(14.5,12)(0.5,-6)
% transistors
\cell{9}{0}{l}{\tusual{\theta}}
\cell{0.5}{4}{l}{\twiggly{x_1}}
\cell{0.5}{-4}{l}{\twiggly{x_n}}
% short wires
\cell{13}{0}{l}{\toutputrgt{a}}
% long wires
\qbezier(9,1.5)(7.5,1.5)(7,2.75)
\qbezier(5,4)(6.5,4)(7,2.75)
\cell{7.5}{2.75}{l}{a_1}
\qbezier(9,-1.5)(7.5,-1.5)(7,-2.75)
\qbezier(5,-4)(6.5,-4)(7,-2.75)
\cell{7.5}{-2.75}{l}{a_n}
% ellipses
\cell{8.2}{0.3}{c}{\vdots}
\cell{2}{0}{c}{\cdot}
\cell{2}{1}{c}{\cdot}
\cell{2}{-1}{c}{\cdot}
\end{picture}
\end{centredpic}
\mbox{\hspace{2em}}
\goesto
\mbox{\hspace{2em}}
%
% CODOMAIN
%
\begin{centredpic}
\begin{picture}(8.5,12)(1.5,-6)
% transistor
\cell{2}{6}{tr}{\twiggleleft}
\cell{2}{5.2}{tl}{\twiggleright}
\cell{2}{4.4}{tr}{\twiggleleft}
\cell{2}{3.6}{tl}{\twiggleright}
\cell{2}{2.8}{tr}{\twiggleleft}
\cell{2}{2}{tl}{\twiggleright}
\cell{2}{1.2}{tr}{\twiggleleft}
\cell{2}{0.4}{tl}{\twiggleright}
\cell{2}{-0.4}{tr}{\twiggleleft}
\cell{2}{-1.2}{tl}{\twiggleright}
\cell{2}{-2}{tr}{\twiggleleft}
\cell{2}{-2.8}{tl}{\twiggleright}
\cell{2}{-3.6}{tr}{\twiggleleft}
\cell{2}{-4.4}{tl}{\twiggleright}
\cell{2}{-5.2}{tr}{\twiggleleft}
\put(8,0){\line(-1,-1){6}}
\put(8,0){\line(-1,1){6}}
\cell{5}{0}{c}{\scriptstyle \ovln{\theta}(x_1, \ldots, x_n)}
% tags
\cell{8}{0}{l}{\toutputrgt{a}}
\end{picture}
\end{centredpic}
\]
% \hand{40}{55}
\caption{Action of a multicategory on an algebra}
\label{fig:alg-action}
\end{figure}
satisfying axioms of the same shape as the associativity and second
identity axioms in Definition~\ref{defn:cl-multi}.  This explicit form
makes it clear that a \demph{map of $C$-algebras},%
\lbl{p:map-of-algs}
$\alpha: X \go Y$, should be defined as a family of functions
\[
\left(X(a) \goby{\alpha_a} Y(a)\right)_{a\in C_0} 
\]
satisfying the evident compatibility condition; so we have a category
$\Alg(C)$%
\glo{Algplainmulti}
of $C$-algebras.  An alternative definition uses the notion of
transformation between maps between multicategories, as described
in~\ref{sec:om-further}.  

% All of this can also be repeated in the more
% elaborate settings of symmetric and/or enriched multicategories.  

\begin{example}	\lbl{eg:alg-mon}
If $A$ is a strict monoidal category then an algebra for its underlying
multicategory is, by~\ref{eg:map-mti-mon}, just a lax monoidal functor from
$A$ to $(\Set, \times, 1)$.
\end{example}

\begin{example}	\lbl{eg:alg-multi-unary}
If $C$ is a multicategory in which all arrows are
unary~(\ref{eg:multi-unary}), and so essentially just a category, then
$\Alg(C)$ is the ordinary functor category $\ftrcat{C}{\Set}$.
\end{example}

\begin{example}	\lbl{eg:alg-multi-Cayley}
As the pictures suggest, there is for each multicategory $C$ an algebra $X$
defined by taking $X(a)$ to be the set $C(;a)$ of arrows in $C$ from the
empty sequence into $a$.  When $C$ is the multicategory of modules%
\index{module!ring@over ring}
over
some commutative ring $R$ (Example~\ref{eg:multi-mon}), this is the evident
forgetful map $C \go \Set$.
\end{example}

\begin{example}	\lbl{eg:alg-multi-End}
Any family $(X(a))_{a\in S}$ of sets, indexed over any set $S$, gives rise
to a multicategory $\End(X)$%
\glo{Endplainmulti}
with object-set $S$, the \demph{endomorphism%
\index{endomorphism!plain multicategory}%
\index{multicategory!endomorphism}
multicategory} of $X$.  This
is constructed by transporting back from the multicategory $\Set$: an arrow
$a_1, \ldots, a_n \go a$ in $\End(X)$ is a function
\[
X(a_1) \times\cdots\times X(a_n) \go X(a),
\]
and composition and identities are as in $\Set$.  An algebra%
\index{algebra!plain multicategory@for plain multicategory}%
\index{multicategory!algebra for}
for a
multicategory $C$ can then be defined, equivalently, as a family
$(X(a))_{a\in C_0}$ of sets together with a multicategory map $f: C \go
\End(X)$ fixing the objects.
\end{example}

Set-valued functors on ordinary categories are also the same thing as
discrete opfibrations.%
\index{fibration!discrete opfibration}
 There is a parallel theory of opfibrations for
multicategories (\ref{sec:non-alg-notions},~\ref{sec:alg-fibs},
Leinster~\cite{FM}), but this will not be of central concern.

\section{Classical operads}
\lbl{sec:cl-opds}%
\index{multicategory!operad@\vs.\ operad|(}

An operad is a multicategory with only one object.  In a sense there is no
more to be said: the definitions of map between operads, algebra
for an operad, and so on, are just special cases of the definitions for
multicategories.  

On the other hand, operads have a distinctive feel to them and are worth
considering in their own right.  The analogy to keep in mind is monoids
versus categories.  A monoid is nothing but a one-object category, and many
basic monoid-theoretic concepts are specializations of category-theoretic
concepts; nevertheless, monoids still form a natural and interesting class
of structures.  The same goes for operads and multicategories.  Conversely,
a category can be regarded as a many-object monoid, and a multicategory as
a many-object operad; this leads to the alternative name `coloured operad'
for multicategory.

Most of this section is examples.  But first we give the definitions,
describe equivalent `explicit' forms, and discuss the ever-important
question of terminology.

\begin{defn}	\lbl{defn:plain-opd}
An \demph{operad} is a multicategory $C$ with exactly one object.  The
category $\Operad$%
\glo{Operad}
is the full subcategory of $\Multicat$ whose objects are
the operads.
% If $P$ is an operad then $\Alg(P)$, the
% \demph{category of algebras for $P$}, is the category of algebras for the
% corresponding multicategory.  
\end{defn}

Let $C$ be an operad, with single object called $\star$.  Then $C$ consists
of one set
\[
P(n) = C(\underbrace{\star, \ldots, \star}_n; \star)
\]
for each $n\in\nat$, together with composition and identities satisfying
axioms.  Of course, it makes no difference (up to isomorphism) what the
single object is called, so an operad can equivalently be described as
consisting of
\begin{itemize}
\item a sequence $(P(n))_{n\in\nat}$ of sets, whose elements $\theta$ will
be called the \demph{$n$-ary%
\index{n-ary@$n$-ary}
operations}%
\index{operation}%
\index{arrow}
of $P$ and drawn as
\[
n
\left\{
\mbox{\rule[-4ex]{0em}{8ex}}
\right. 
% \!\!\!\!\!\!
\begin{centredpic}
\begin{picture}(6,4)(-1,-2)
\cell{0}{0}{l}{\tusual{\theta}}
\cell{0}{0}{r}{\tinputsslft{}{}{}}
\cell{4}{0}{l}{\toutputrgt{}}
\end{picture}
\end{centredpic}
\]
\item for each $n, k_1, \ldots, k_n \in \nat$, a function
\[
\begin{array}{rcl}
P(n) \times P(k_1) \times\cdots\times P(k_n)	&
\go	&
P(k_1 + \cdots + k_n)	\\
(\theta, \theta_1, \ldots, \theta_n)	&
\goesto	&
\theta \of (\theta_1, \ldots, \theta_n),
\end{array}
\]%
\glo{opdcomposite}%
called \demph{composition} (as in Fig.~\ref{fig:multi-comp} but without the
labels $a$, $a_i$, $a_i^j$) 
\item an element $1 = 1_P \in P(1)$,%
\glo{opdidentity}
called the \demph{identity}, 
\end{itemize}
satisfying associativity and identity axioms as in the definition of
multicategory~(\ref{defn:cl-multi}).  This direct description, or something
like it, is what one usually sees as the definition of operad.  Similarly,
a map $f: P \go Q$ of operads consists of a family
\[
(f_n: P(n) \go Q(n))_{n\in\nat}
\]
of functions, preserving composition and identities.  

All of this works just as well if the $P(n)$'s are allowed to be objects of
some chosen symmetric monoidal category $\cat{V}$, rather than necessarily
being sets; these are \demph{operads%
\index{operad!symmetric monoidal category@in symmetric monoidal category}
in $\cat{V}$}, or
\demph{$\cat{V}$-operads},%
\lbl{p:defn-V-Operad}
and form a category $\cat{V}\hyph\Operad$.%
\glo{VOperad}
 We noted a more general version
of this for multicategories, where the appropriate terminology was
`multicategory enriched%
\index{enrichment!plain multicategory@of plain multicategory!symmetric monoidal category@in symmetric monoidal category}
in $\cat{V}$'; we also said there that symmetric
monoidal categories are actually not the most natural
setting~(\ref{sec:enr-mtis}).

Symmetries%
\index{operad!symmetric}
are often present in operads.  Again, this was noted in the more
general context of multicategories.  A symmetric structure on an operad $P$
consists, then, of an action of the symmetric group $S_n$ on $P(n)$ for
each $n\in\nat$, subject to axioms; the precise definition is
in~\ref{defn:sym-mti}.  Most authors use the term `operad' to mean what I
call a symmetric operad, and `non-symmetric' or `non-$\Sigma$' operad for
what I call a (plain) operad.  For us the non-symmetric case will be by far
the more important.  The generalized (multicategories and) operads that we
consider for much of this book are in some sense a more advanced
replacement for symmetric operads; see p.~\pageref{p:sym-discussion} for
further explanation of this point of view.

Many authors use the term `coloured%
\lbl{p:col-opd}\index{operad!coloured}
operad' to mean multicategory.  The `colours' are the objects, and the idea
is that a multicategory is a more complicated version of an operad.  The
analogous usage would make a category a `coloured monoid' and a groupoid a
`coloured group'.  I prefer the term `multicategory' because it emphasizes
that we are dealing with a categorical structure, with objects and arrows.
% (Thus linear algebra concerns the category of vector spaces,
% and multilinear algebra the multicategory of vector spaces.)
Moreover, although thinking of the objects as colours is practical when
there is only a small, finite, number of them, it becomes somewhat baroque
otherwise: in the multicategory of abelian groups, for instance, one has to
paint each group a different colour (the real numbers are green, the cyclic
group of order 10 is pink, and so on).  This is not entirely a frivolous
issue: the coloured viewpoint has given rise to eyebrow-raising statements
such as
\begin{quote}	\lbl{p:MSS-quote}
  [\ldots] it serves us well to have a subtle generalization of operad
  known as a bicolored operad.  Still more colorful operads can be defined,
  but they are currently not of great importance
\end{quote}
(Markl,%
\index{Markl, Martin}
Shnider%
\index{Shnider, Steve}
and Stasheff~\cite[p.~115]{MSS});%
\index{Stasheff, Jim}
of course, most of the
coloured operads that one meets every day 
($\Set$, $\Top$, $R\hyph\fcat{Mod}$, \ldots)
% (sets, spaces, modules, \ldots)
have not just more than two colours, but a proper class of them.  The
quotation above is analogous to, and indeed
includes~(\ref{eg:multi-unary}), the statement that no important category
has more than two objects.%
\index{multicategory!operad@\vs.\ operad|)}

Algebras%
\index{operad!algebra for}\index{algebra!plain operad@for plain operad}
for operads are very important.  Translating the multicategorical
definition~(\ref{defn:alg-multi}) into explicit language, an algebra for an
operad $P$ is a set $X$ together with a function
\[
\ovln{\theta}: X^n \go X
\]%
\glo{thetabaropd}%
for each $n\in\nat$ and $\theta\in P(n)$, satisfying the evident axioms.
If $P$ is a symmetric operad then a further axiom must be satisfied.  If
$P$ is an operad in a symmetric monoidal category $\cat{V}$ then $X$ is not
a set but an object of $\cat{V}$, and an algebra structure on $X$ is a
family $(P(n) \otimes X^{\otimes n} \go X)_{n\in\nat}$ of maps, satisfying
axioms.  (More generally, $X$ could be an object of a monoidal category
either acted on by $\cat{V}$ or enriched in $\cat{V}$.)  It is clear in all
situations what a map of algebras is.

The first group of examples follows the slogan `an operad is an algebraic%
\index{algebraic theory!plain operad as|(}%
\index{operad!algebraic theory@as algebraic theory|(}
theory'. 

\begin{example}	\lbl{eg:opd-terminal}
The terminal%
\index{operad!terminal}
operad $1$ has exactly one $n$-ary operation for each
$n\in\nat$.  An algebra for $1$ is a set $X$ together with a function
\[
\begin{array}{rcl}
X^n	&\go	&X	\\
(x_1, \ldots, x_n)	&\goesto	&(x_1 \cdot\ldots\cdot x_n)
\end{array}
\]
for each $n\in\nat$, satisfying the axioms
\begin{eqnarray*}
((x_1^1 \cdot \ldots \cdot x_1^{k_1}) 
\cdot \ldots \cdot 
(x_n^1 \cdot \ldots \cdot x_n^{k_n}))	&
=	&
(x_1^1 \cdot \ldots \cdot x_n^{k_n})			\\
x	&=	&(x).			
\end{eqnarray*}
Hence $\Alg(1)$ is the category of monoids.%
\index{monoid!operad-algebra@as operad-algebra}

If $1$ is considered as a \emph{symmetric} operad then there is a further
axiom on its algebras $X$:
\[
(x_{\sigma(1)} \cdot \ldots \cdot x_{\sigma(n)})	
=	
(x_1 \cdot \ldots \cdot x_n) 
\]
for all $\sigma\in S_n$ (as will follow from the definition,
p.~\pageref{p:defn-sym-alg}).  This means that algebras for $1$ are now
\emph{commutative} monoids.
\end{example}

\begin{example}
Various sub-operads of $1$ are commonly encountered.  The smallest operad
$P$ is given by $P(1) = 1$ and $P(n) = \emptyset$ for $n\neq 1$; its
algebras are merely sets.  The unique operad $P$ satisfying $P(0) =
\emptyset$ and $P(n) = 1$ for $n\geq 1$ has semigroups as its algebras.  (A
\demph{semigroup}%
\index{semigroup}
is a set equipped with an associative binary operation; a
monoid is a semigroup with identity.)  As a kind of dual, the unique operad
$P$ satisfying $P(n) = 1$ for $n\leq 1$ and $P(n) = \emptyset$ for $n > 1$
has as its algebras \demph{pointed%
\index{pointed set}
sets} (sets equipped with a basepoint).
\end{example}

\begin{example}
Let $M$ be a monoid.  Then there is an operad $P$ in which $P(1) = M$,
$P(n)= \emptyset$ for $n\neq 1$, and the composition and identity are the
multiplication and unit of $M$.  An algebra for $P$ consists of a set $X$
together with a function $\ovln{\theta}: X \go X$ for each $\theta\in M$,
satisfying axioms; in other words, it is a set with a left $M$-action.

This example is the one-object case of a multicategorical example that we
have already seen: a multicategory all of whose arrows are unary is just a
category~(\ref{eg:multi-unary}), and its algebras are just functors
from that category into $\Set$~(\ref{eg:alg-multi-unary}).
\end{example}

\begin{example} \lbl{eg:opd-sr}
We have seen so far that the theories of monoids, of semigroups, of pointed
sets, and of $M$-sets (for a fixed monoid $M$) can all be described by
operads.  The natural question is: which algebraic
theories can be
described by operads?  The answer is: the strongly%
\index{strongly regular theory}
regular finitary
theories.  Here is the definition; the proof is
in~\ref{sec:opds-alg-thys}. 

A finitary algebraic theory is \demph{strongly regular} if it can be
presented by operations and strongly regular equations.  In turn, an
equation (made up of variables and finitary operation symbols) is
\demph{strongly regular} if the same variables appear in the same order,
without repetition, on each side. So all of
\[
(x\cdot y) \cdot z = x\cdot (y\cdot z),
\diagspace
x\cdot 1 = x,
\diagspace
(x^y)^z = x^{y\cdot z},
\]
but none of 
\[
x\cdot 0 = 0,
\diagspace
x\cdot y = y\cdot x,
\diagspace
x\cdot (y+z) = x\cdot y + x\cdot z,
\]
are strongly regular.  For instance, the theory of monoids is strongly
regular.  One would guess that the theories of commutative monoids and of
groups are not strongly regular, because their usual presentations involve
the equations $x\cdot y = y\cdot x$ and $x^{-1} \cdot x = 1$ respectively,
neither of which is strongly regular.  For the moment the possibility
remains that these theories can be presented by some devious selection of
operations and strongly regular equations, but in~\ref{eg:mon-CJ} we will
see a method for proving that a given theory is not strongly regular, and
in particular it can be applied to confirm that no such devious selections
exist.

The name `strongly regular' is due to Carboni and Johnstone~\cite{CJ}, and,
as they explain, is something of an accident.
\end{example}

\begin{example}	\lbl{eg:opd-sr-enr}
A much wider range of algebraic
theories is covered if symmetries are
allowed and if the $P(n)$'s are allowed to be objects of a symmetric
monoidal category $\cat{V}$ instead of just sets.  For instance, if
$\cat{V}$ is the category of vector spaces over some field then there is a
symmetric operad $P$ in $\cat{V}$ generated by one element $\theta\in P(2)$
subject to the equations
\begin{eqnarray*}
\theta + \theta\cdot\tau	&=	&0	\\
\theta \of (1, \theta)
+ (\theta \of (1, \theta))\cdot \sigma
+ (\theta \of (1, \theta))\cdot \sigma^2
				&=	&0
\end{eqnarray*}
where $\tau\in S_2$ is a 2-cycle, $\sigma\in S_3$ is a 3-cycle, and the
action of $S_n$ on $P(n)$ is denoted by a dot; $P$-algebras are exactly
Lie%
\index{Lie algebra}
algebras.  Operads of this kind will not be of direct concern here,
our generalization of the notion of operad being in a different direction,
but they are very much in use.
\end{example}%
\index{algebraic theory!plain operad as|)}
\index{operad!algebraic theory@as algebraic theory|)}

The most familiar examples of multicategories are those underlying%
\index{multicategory!underlying}
monoidal
categories~(\ref{eg:multi-mon}) or, more generally, their full
sub-multicategories~(\ref{eg:multi-some-of-mon}).  Here is the one-object
case.

\begin{example}	\lbl{eg:opd-from-comm-mon}%
\index{monoidal category!degenerate}%
\index{operad!commutative monoid@from commutative monoid}
A one-object strict monoidal category is a commutative
monoid~(\ref{eg:str-mon-comm}), so any commutative%
\index{monoid!commutative!operad from}
monoid $(A,+,0)$ has an
`underlying' operad $P$.  Concretely, $P(n)=A$ for all $n$, composition is
\[
\begin{array}{rcl}
P(n) \times P(k_1) \times\cdots\times P(k_n)	&\go	&
P(k_1 +\cdots + k_n)	\\
(a, a_1, \ldots, a_n)	&\goesto	&a + a_1 + \cdots + a_n,
\end{array}
\]
and the identity is $0\in P(1)$.
\end{example}

\begin{example}	\lbl{eg:opd-End}%
\index{endomorphism!plain operad|(}%
\index{operad!endomorphism|(}
For any object $b$ of a monoidal category $B$, there is an operad structure
on the sequence of sets $(B(b^{\otimes n}, b))_{n\in\nat}$, given by
substitution.  This is a special case of
Example~\ref{eg:multi-some-of-mon}.  We call it $\End(b)$,%
\glo{Endplainoperad}
the
\demph{endomorphism operad} of $b$.  In the case $B = \Set$ this is
compatible with the $\End$ notation of~\ref{eg:alg-multi-End}, and
specializing the observations there, an algebra%
\index{operad!algebra for}\index{algebra!plain operad@for plain operad}
for an operad $P$ amounts
to a set $X$ together with a map $P \go \End(X)$ of operads.
\end{example}

\begin{example}
The \demph{operad of curves}%
\index{operad!curves@of curves}
$P$ is defined by
\[
P(n) =
\{ \textrm{smooth maps } \reals \go \reals^n \}
\]
and substitution.  If $B$ denotes the monoidal category of smooth
manifolds%
\index{manifold!monoidal category of}
and smooth maps, with product as monoidal structure, then $P$ is the
endomorphism operad of $\reals$ in $B^\op$.  
\end{example}

\begin{example}
Fix a commutative ring $k$.  Substituting polynomials into variables gives
an operad structure on the sequence of sets $(k[X_1, \ldots,
X_n])_{n\in\nat}$; this is the endomorphism operad of the object $k[X]$
of the monoidal category $(\textrm{commutative }k\textrm{-algebras})^\op$.%
\index{algebra!ring@over ring}
\end{example}

\begin{example}
The same construction works for any algebraic%
\index{algebraic theory!operad of free algebras}%
\index{operad!free algebras@of free algebras}
theory: if $T$ is a monad on
$\Set$ then there is a natural operad structure on $(T(n))_{n\in\nat}$.
The first $n$ here denotes an $n$-element set, so $T(n)$ is the set of
words in $n$ variables.  Informally, composition is substitution of words;
formally, the composition functions can be written down in terms of the
monad structure on $T$ (exercise).  This is the endomorphism operad of the
free algebra on one generator in the opposite of the category of algebras,
where tensor is coproduct of algebras.  

In particular, this gives operad structures on each of $(n)_{n\in\nat}$
(from the theory of sets, that is, the identity monad), $(k[X_1, \ldots,
X_n])_{n\in\nat}$ (from the theory of $k$-algebras, the previous example),
and $(k^n)_{n\in\nat}$ (from the theory of $k$-modules).
% , and $(2^n)_{n\in\nat}$
% (from the theory of Boolean algebras).
\end{example}

\begin{example}	\lbl{eg:opd-Riemann}
There is an operad $P$ in which $P(n)$ is the set of isomorphism classes of
Riemann%
\index{Riemann surface}\index{manifold!operad of}
surfaces whose boundaries are identified with the disjoint union of
$(n+1)$ copies of the circle $S^1$ (in order, with the first $n$ thought of
as inputs and the last as an output); composition is gluing.  Many
geometric and topological variants of this example exist.  See the end
of~\ref{sec:struc} for remarks on a more sophisticated version in which
there is no need to quotient out by isomorphism.
\end{example}

\begin{example}	\lbl{eg:opd-monoid-powers}
We saw in~\ref{eg:multi-pretend-coproducts} that given any category $A$, we
can pretend that $A$ has coproducts%
\index{operad!coproducts@from coproducts}
to obtain a multicategory $C$ with the
same objects as $A$.  The one-object version is: given any monoid $M$,
there is an operad $P$ with $P(n) = M^n$ and composition
\begin{eqnarray*}
&&
(\alpha_1, \ldots, \alpha_n) \of 
((\alpha_1^1, \ldots, \alpha_1^{k_1}), \ldots, 
(\alpha_n^1, \ldots,  \alpha_n^{k_n}))\\
&=&
(\alpha_1 \alpha_1^1, \ldots, \alpha_1 \alpha_1^{k_1}, \ldots,
\alpha_n \alpha_n^1, \ldots, \alpha_n \alpha_n^{k_n})
\end{eqnarray*}
($\alpha_i, \alpha_i^j \in M$).  Some specific instances are
in~\ref{eg:opd-fractals} and~\ref{eg:opd-little-disks} below.
\end{example}

\begin{figure}
\begin{center}
\setlength{\unitlength}{1mm}
\begin{picture}(101,120)(-21,0)
% Peano label
\cell{-17}{4}{br}{\textrm{(d)}}
% Peano generator
\cell{0}{0}{bl}{%
\begin{picture}(30,35)(-15,-20)
\cell{-15}{0}{c}{\zmark}
\cell{0}{0}{c}{\zmark}
\cell{0}{15}{c}{\zmark}
\cell{0}{-15}{c}{\zmark}
\cell{15}{0}{c}{\zmark}
\cell{-15}{-2}{t}{z_0}
\cell{0}{-2}{t}{z_1, z_3, z_5}
\cell{0}{13}{t}{z_2}
\cell{0}{-17}{t}{z_4}
\cell{15}{-2}{t}{z_6}
\end{picture}}
% Peano curve
\cell{50}{5}{bl}{\epsfig{file=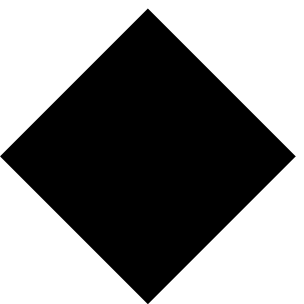}}
% Sierpinski label
\cell{-17}{54}{br}{\textrm{(c)}}
% Sierpinski generator
\cell{0}{50}{bl}{%
\begin{picture}(30,31)(-15,-5)
\cell{-15}{0}{c}{\zmark}
\cell{0}{0}{c}{\zmark}
\cell{-7.5}{13}{c}{\zmark}
\cell{7.5}{13}{c}{\zmark}
\cell{15}{0}{c}{\zmark}
\cell{-15}{-2}{t}{z_0}
\cell{0}{-2}{t}{z_1, z_4}
\cell{-7.5}{11}{t}{z_2}
\cell{7.5}{11}{t}{z_3}
\cell{15}{-2}{t}{z_5}
\end{picture}}
% Sierpinski curve
\cell{50}{55}{bl}{\epsfig{file=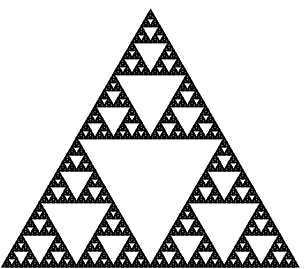}}
% Koch label
\cell{-17}{95}{br}{\textrm{(b)}}
% Koch generator
\cell{0}{91}{bl}{%
\begin{picture}(30,14)(-15,-5)
\cell{-15}{0}{c}{\zmark}
\cell{-5}{0}{c}{\zmark}
\cell{0}{8.6}{c}{\zmark}
\cell{5}{0}{c}{\zmark}
\cell{15}{0}{c}{\zmark}
\cell{-15}{-2}{t}{z_0}
\cell{-5}{-2}{t}{z_1}
\cell{0}{6.6}{t}{z_2}
\cell{5}{-2}{t}{z_3}
\cell{15}{-2}{t}{z_4}
\end{picture}}
% Koch curve
\cell{50}{96}{bl}{\epsfig{file=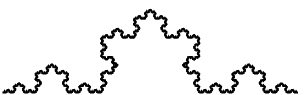}}
% Dedekind label
\cell{-17}{119}{br}{\textrm{(a)}}
% Dedekind generator
\cell{0}{115}{bl}{%
\begin{picture}(30,5)(-15,-5)
\cell{-15}{0}{c}{\zmark}
\cell{0}{0}{c}{\zmark}
\cell{15}{0}{c}{\zmark}
\cell{-15}{-2}{t}{z_0}
\cell{0}{-2}{t}{z_1}
\cell{15}{-2}{t}{z_2}
\end{picture}}
% Dedekind curve
\put(50,120){\line(1,0){30}}
\end{picture}%
\index{Koch curve}%
\index{Sierpinski gasket@Sierpi\'nski gasket}%
\index{Peano curve@P\'eano curve}%
\end{center}
% \hand{110}{56}
\caption{On the left, $(z_0, \ldots, z_n) \in Q(n)$, and on the right, its
  unique fixed point in $\mc{K}$: (a)~interval, (b)~Koch curve,
  (c)~Sierpi\'nski gasket, (d)~P\'eano curve} 
\label{fig:fractals}
\end{figure}

\begin{example}	\lbl{eg:opd-fractals}
Let $G$ be the group of affine automorphisms of the complex plane---maps of
the form $z \goesto a + bz$ with $a, b \in \complexes$ and $b\neq 0$.  (So
$G$ is generated by translations, rotations and dilatations $z \goesto
\lambda z$ with $\lambda> 0$.)  Let $P$ be the operad $(G^n)_{n\in\nat}$
defined in the previous example.  Then $P$ has a sub-operad $Q$ given by
\begin{eqnarray*}
Q(n)	&
=	& 
\{
(\alpha_1, \ldots, \alpha_n) \in G^n \such	
0 = \alpha_1(0),\, \alpha_1(1) = \alpha_2(0),\, \ldots,\,
\\&&
\alpha_{n-1}(1) = \alpha_n(0),\, \alpha_n(1) = 1
\}
\end{eqnarray*}
($n\geq 1$) and $Q(0) = \emptyset$.  Since $G$ acts freely and transitively
on the set of ordered pairs of distinct points in the plane, we have
\begin{equation}	\label{eq:planar-seqs}
Q(n) \iso 
\{
(z_0, \ldots, z_n) \in \complexes^{n+1}
\such
0 = z_0 \neq z_1 \neq \cdots \neq z_n = 1
\},
\end{equation}
and we call $Q$ the \demph{operad of finite%
\index{operad!finite planar sequences@of finite planar sequences}
planar sequences}.  The set
$\mc{K}$ of nonempty compact subsets of $\complexes$ is a $Q$-algebra, with
action
\[
\begin{array}{rrcl}
\ovln{(\alpha_1, \ldots, \alpha_n)}:	&
\mc{K}^n	&
\go	&
\mc{K}	\\
	&
(S_1, \ldots, S_n)	&
\goesto	&
\alpha_1 S_1 \cup \cdots \cup \alpha_n S_n.	\\
\end{array}
\]

Given any operad $R$, algebra $X$ for $R$, and operation $\theta \in R(n)$,
we can study the \demph{fixed%
\index{fixed point}
points} of $\theta$, that is, the elements
$x\in X$ satisfying $\ovln{\theta}(x, \ldots, x) = x$.  Here this is the
study of affinely self-similar%
\index{self-similarity}\index{fractal}
planar sets, and a theorem of
Hutchinson~\cite{Hut}%
\index{Hutchinson, John}
implies that any $(\alpha_1, \ldots, \alpha_n) \in
Q(n)$ for which each $\alpha_i$ is a contraction (or in terms
of~\bref{eq:planar-seqs}, any $(z_0, \ldots, z_n)$ for which $|z_{i+1} -
z_i| < 1$ for each $i$) has a unique fixed point in $\mc{K}$.  Some of
these fixed points are shown in Fig.~\ref{fig:fractals}.%
\index{endomorphism!plain operad|)}%
\index{operad!endomorphism|)}
\end{example}

\index{loop space|(}
Operads rose to fame for their role in loop space theory.  Let $\Top_*$%
\glo{Topstar}
be
the category of topological spaces with a distinguished basepoint.  For any
$Y \in \Top_*$, the set $\Top_*(S^1, Y)$ of basepoint-preserving maps from
the circle $S^1$ into $Y$, endowed with the canonical ($=$ compact-open)
topology, is called the \demph{loop space}%
\lbl{p:defn-loop-space}
on $Y$ and written $\Omega(Y)$.%
\glo{Omegaloop}
 This defines an endofunctor $\Omega$ of
$\Top_*$.  If $d\in\nat$ then $\Omega^d(Y) \iso \Top_*(S^d, Y)$;%
\glo{dsphere}
a space
homeomorphic to $\Omega^d(Y)$ for some $Y$ is called a \demph{$d$-fold loop
space}.

Loop spaces are nearly monoids.%
\index{monoid!topological}
 A binary multiplication on $\Omega(Y)$ is
a rule for composing two loops in $Y$; the standard choice is to travel the
first then the second, each at double speed.  The unit is the constant
loop.  The associativity and unit laws are obeyed not quite exactly, but up
to homotopy.  Moreover, if one chooses particular homotopies to witness
this, then these homotopies obey laws of their own---not quite exactly, but
up to homotopy; and so \latin{ad infinitum}.  A loop space therefore admits
an algebraic structure of a rather complex kind, and it is to describe this
complex structure that operads are so useful.  In the following group of
examples we meet various operads for which loop spaces, and more generally
$d$-fold loop spaces, are naturally algebras.

\begin{example}	\lbl{eg:opd-univ-loop}
The coproduct in $\Top_*$ is the wedge product $\wej$%
\glo{wedge}
(disjoint union with basepoints identified).  This makes $\Top_*$ into a
monoidal category, so for each $d\in\nat$ there is an endomorphism operad
$\fcat{U}_d$ given by $\fcat{U}_d(k) = \Top_*(S^d, (S^d)^{\wej k})$.  Any
$d$-fold loop space is naturally a $\fcat{U}_d$-algebra via the evident
maps
\begin{eqnarray*}
\fcat{U}_d(k) \times (\Omega^d(Y))^k	&
\iso	&
\Top_*(S^d, (S^d)^{\wej k}) \times \Top_*((S^d)^{\wej k}, Y)	\\
&\go	& 
\Top_*(S^d, Y)	\\
&\iso	& 
\Omega^d(Y).
\end{eqnarray*}
Borrowing the terminology of Salvatore~\cite{Sal},%
\index{Salvatore, Paolo}
$\fcat{U}_d$ is the
\demph{universal%
\index{operad!universal for loop spaces}
operad for $d$-fold loop spaces}.  
\end{example}

\begin{example}	\lbl{eg:opd-little-disks}
Let $d\in\nat$ and let $G$ be the group of transformations $\alpha$ of
$\reals^d$ of the form $\alpha(\mb{x}) = \mb{a} + \lambda \mb{x}$, with
$\mb{a} \in \reals^d$ and $\lambda > 0$.  Denote by $D(\mb{a}, \lambda)$
the closed disk
(ball) in $\reals^d$ with centre $\mb{a}$ and radius
$\lambda$.  Then $(G^k)_{k\in\nat}$ is naturally an operad
(Example~\ref{eg:opd-monoid-powers}), and the \demph{little%
\index{operad!little disks}
$d$-disks
operad} $\ldisks_d$%
\glo{littleddisks}
is the sub-operad defined by
\begin{eqnarray*}
\ldisks_d(k)	&
=	&
\{ (\alpha_1, \ldots, \alpha_k) \in G^k \such	
\textrm{the images of } D(\mb{0}, 1) \textrm{ under }
\alpha_1, \ldots, \alpha_k \textrm{ are }
\\&&
\textrm{disjoint subsets of }
D(\mb{0}, 1) \}.
\end{eqnarray*}
Since $G$ acts freely and transitively on the set $\{D(\mb{a}, \lambda)
\such \mb{a} \in \reals^d, \lambda>0 \}$ of disks, $\ldisks_d(k)$ may be
identified with the set of configurations of $d$ ordered disjoint `little'
disks inside the unit disk: the $i$th little disk is $\alpha_i D(\mb{0},
1)$.  Fig.~\ref{fig:little-disks}
\begin{figure}
\centering
$\begin{array}{c}
\setlength{\unitlength}{1mm}
\begin{picture}(55,61)(-1,0)
\cell{0}{0}{bl}{\epsfig{file=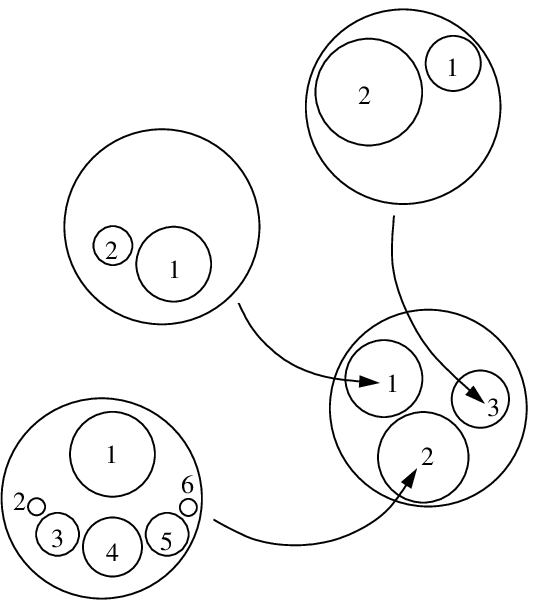}}
\cell{8}{47}{c}{\theta_1}
\cell{1}{19}{c}{\theta_2}
\cell{31}{58}{c}{\theta_3}
\cell{52}{11}{c}{\theta}
\end{picture}
\end{array}$
\diagspace
composes to
\diagspace
$\begin{array}{c}
\setlength{\unitlength}{1mm}
\begin{picture}(20,20)
\cell{0}{0}{bl}{\epsfig{file=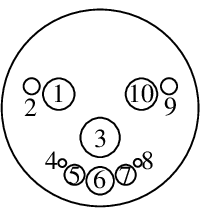}}
\cell{10}{-2}{t}{\theta \of (\theta_1, \theta_2, \theta_3)}
\end{picture}
\end{array}$
% \hand{70}{57}
\caption{Window to Plato's heaven}%
\index{Plato}
\label{fig:little-disks}
\end{figure}
shows an example of composition  
\[
\ldisks_2(3) \times
(\ldisks_2(2) \times \ldisks_2(6) \times \ldisks_2(2)) 
\go
\ldisks_2(10)
\]
in the little 2-disks operad $\ldisks_2$.

The same can be done with cubes instead of disks; this leaves the homotopy
type of $\ldisks_d(k)$ unchanged.  $\ldisks_d(k)$ is also homotopy
equivalent to the space of configurations of $k$ ordered, distinct points
in $\reals^d$, but there is no obvious way to put an operadic composition
on this sequence of spaces.

Any $d$-fold loop space $\Omega^d(Y)$ is naturally a $\ldisks_d$-algebra.
Concretely, a configuration of $k$ little $d$-disks shows how to glue
together $k$ based maps $S^d \go Y$ to make a single based map $S^d \go Y$
(noting that $S^d$ is $D(\mb{0}, 1)$ with its boundary collapsed to a
point).  Abstractly, any element $(\alpha_1, \ldots, \alpha_k)$ of
$\ldisks_d(k)$ defines a continuous injection from the disjoint union of
$k$ copies of $D(\mb{0}, 1)$ into a single copy of $D(\mb{0}, 1)$.
Collapsing the boundaries and taking the inverse gives a based map $S^d \go
(S^d)^{\wej k}$.  This process embeds $\ldisks_d$ as a sub-operad of
$\fcat{U}_d$~(\ref{eg:opd-univ-loop}),%
\index{operad!universal for loop spaces}
and the inclusion $\ldisks_d \rIncl \fcat{U}_d$ induces a functor
$\Alg(\fcat{U}_d) \go \Alg(\ldisks_d)$ in the opposite direction.  So any
$d$-fold loop space is a $\ldisks_d$-algebra.  In fact, the converse is
almost true---roughly, any algebra for $\ldisks_d$ (as a symmetric operad
in $\Top$) is a $d$-fold loop space---see Adams~\cite[Ch.~2]{Ad}, for
instance.
\end{example}

\begin{example}	\lbl{eg:opd-associahedra}
The first really important operad to be considered in topology was
Stasheff's~\cite{StaHAHI}%
\index{Stasheff, Jim}
operad $K$%
\glo{assK}
of associahedra.%
\index{associahedron}
 (The terms `operad'
and `associahedra' came later.)  This is a non-symmetric, topological
operad; $K(n)$ is an $(n-2)$-dimensional solid polyhedron whose vertices
are indexed by the $n$-leafed, planar, binary, rooted trees.  Any loop
space is naturally a $K$-algebra.  See~\ref{sec:trees} below and Markl,
Shnider and Stasheff~\cite{MSS} for more.
\end{example}

\begin{example}	\lbl{eg:opd-Trimble}%
\index{categorical algebra for operad}%
\index{operad!category over}%
\index{operad!path reparametrizations@of path reparametrizations}
Suppose we are interested in paths%
\index{path!loop@\vs.\ loop}
rather than based loops: then the
appropriate replacement for the space $\fcat{U}_1(k) = \Top_*(S^1,
(S^1)^{\wej k})$ of~\ref{eg:opd-univ-loop} is the space
\[
E(k) = 
\{ \gamma \in \Top([0,1], [0,k]) \such
\gamma(0) = 0 \textrm{ and } \gamma(1) = k \}.
\]
If $Y$ is any space, and $Y(y, y')$ denotes the space of
$[0,1]$-parametrized paths from $y$ to $y'$ in $Y$, then there is a natural
map 
\[
E(k) \times Y(y_0, y_1) \times\cdots\times Y(y_{k-1}, y_k)
\go
Y(y_0, y_k)
\]
for each $y_0, \ldots, y_k \in Y$.  Here we have something like a
$E$-algebra, in that these maps satisfy axioms resembling very closely
those for an algebra for an operad, but it is not an algebra as such; we
will meet the appropriate language when we come to generalized operads in
Part~\ref{part:operads}.  The operad $E$ and its `nearly-algebras' are used
in Trimble's%
\index{Trimble, Todd}
proposed definition of weak $n$-category
(Leinster~\cite{SDN} and~\ref{sec:alg-defns-n-cat} below).
\end{example}%
\index{loop space|)}

A miscellaneous example:
\begin{example}%
\index{multicategory!maps of operads@for maps of operads}
Any operad $P$ gives rise to a `bicoloured operad' (2-object
multicategory), $\fcat{Map}_P$, that has sometimes been found useful (e.g.,
Markl, Shnider and Stasheff~\cite[\S 2.9]{MSS}).  Call the colours
(objects) $0$ and $1$, and for $\epsln_1, \ldots, \epsln_n, \epsln \in
\{0,1\}$, put
\[
\fcat{Map}_P (\epsln_1, \ldots, \epsln_n; \epsln)
=
\left\{
\begin{array}{ll}
P(n)		&\textrm{if } \epsln_i \leq \epsln \textrm{ for each }i	\\
\emptyset	&\textrm{otherwise}.
\end{array}
\right.
\]
Then a $\fcat{Map}_P$-algebra is a map $f: X \go Y$ of $P$-algebras. 

To place this in context, let $A$ be any category and $C$ any
multicategory.  Write $\ovln{A}$ for the multicategory obtained by
pretending that $A$ has coproducts~(\ref{eg:multi-pretend-coproducts}):
then there is an isomorphism of categories
\begin{equation}	\label{eq:exp-transpose-Alg}
\Alg(\ovln{A} \times C) \iso \ftrcat{A}{\Alg(C)}.
\end{equation}
This is easy to show directly; alternatively, once we have defined
transformations of multicategories and hence a category $\ftrcat{C}{D}$
for any multicategory $D$ (see~\ref{sec:om-further}), it follows
by taking $D=\Set$ in the more general isomorphism
\[
\ftrcat{\ovln{A} \times C}{D} \iso \ftrcat{A}{\ftrcat{C}{D}}
\]
Let $\fcat{2}$ be the arrow category $(\blob \go \blob)$:
then~\bref{eq:exp-transpose-Alg} with $A = \fcat{2}$ says that an algebra
for $\ovln{\fcat{2}} \times C$ is a map of $C$-algebras.  Indeed, the
2-object multicategory $\fcat{Map}_P$ defined originally is simply
$\ovln{\fcat{2}} \times P$.
\end{example}

The next example provides the language for defining symmetric%
\index{operad!symmetric|(}\index{multicategory!symmetric|(}
operads and multicategories.

\begin{example}	\lbl{eg:opd-Sym}
The sequence $(S_n)_{n\in\nat}$ consisting of the underlying sets of the
symmetric%
\index{symmetric group}
groups is naturally an operad.  We call it the \demph{operad%
\index{operad!symmetries@of symmetries}
of
symmetries}, $\SymOpd$.%
\glo{SymOpd}
 Fig.~\ref{fig:sym-comp}
\begin{figure}
\[
\begin{centredpic}
\begin{picture}(8,14.9)(0,0)
% thick lines
\thicklines
\put(0,11.6){\framebox(4,1.8){}}
\put(0,6.2){\framebox(4,3.4){}}
\put(0,0){\framebox(4,4.2){}}
\qbezier(4,12.5)(6,10.2)(8,7.9)
\qbezier(4,7.9)(6,5.1)(8,2.3)
\qbezier(4,2.3)(6,7.4)(8,12.5)
% thin lines
\thinlines
%   top box
\put(0,12.9){\line(5,-1){4}}
\put(0,12.1){\line(5,1){4}}
%   middle box
\put(0,9.1){\line(1,0){4}}
\put(0,8.3){\line(5,-2){4}}
\put(0,7.5){\line(5,1){4}}
\put(0,6.7){\line(5,1){4}}
%   bottom box
\put(0,3.7){\line(5,-4){4}}
\put(0,2.9){\line(5,1){4}}
\put(0,2.1){\line(5,1){4}}
\put(0,1.3){\line(5,1){4}}
\put(0,0.5){\line(5,1){4}}
% labels
\cell{2}{14.4}{t}{\rho_1}
\cell{2}{10.6}{t}{\rho_2}
\cell{2}{5.2}{t}{\rho_3}
\cell{6}{14.4}{t}{\sigma}
\end{picture}
\end{centredpic}
\diagspace
\parbox{5em}{\centering composes to give}
\diagspace
\begin{centredpic}
\begin{picture}(8,14.9)(0,0)
% lines for top box
\put(0,12.9){\line(5,-1){4}}
\put(0,12.1){\line(5,1){4}}
\multiput(4,12.9)(0,-0.8){2}{\qbezier(0,0)(2,-2.3)(4,-4.6)}
% lines for middle box
\put(0,9.1){\line(1,0){4}}
\put(0,8.3){\line(5,-2){4}}
\put(0,7.5){\line(5,1){4}}
\put(0,6.7){\line(5,1){4}}
\multiput(4,9.1)(0,-0.8){4}{\qbezier(0,0)(2,-2.8)(4,-5.6)}
% lines for bottom box
\put(0,3.7){\line(5,-4){4}}
\put(0,2.9){\line(5,1){4}}
\put(0,2.1){\line(5,1){4}}
\put(0,1.3){\line(5,1){4}}
\put(0,0.5){\line(5,1){4}}
\multiput(4,3.7)(0,-0.8){5}{\qbezier(0,0)(2,5.1)(4,10.2)}
% label
\cell{4}{14.9}{t}{\sigma\of (\rho_1, \rho_2, \rho_3)}
\end{picture}
\end{centredpic}
\]
% \hand{85}{58}
\caption{Composition in the operad of symmetries}
\label{fig:sym-comp}
\end{figure}
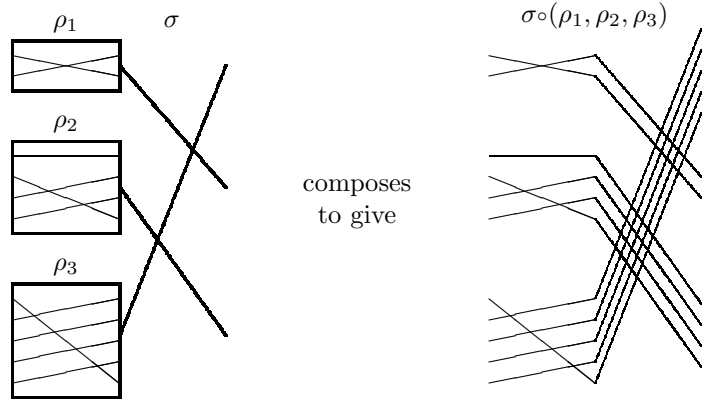
shows an example of composition
\[
\begin{array}{rcl}
S_3 \times (S_2 \times S_4 \times S_5)	&\go	&S_{11},	\\
(\sigma, \rho_1, \rho_2, \rho_3)	&\goesto&
\sigma \of (\rho_1, \rho_2, \rho_3)
\end{array}
\]
with
\[
\begin{array}{cc}
\sigma = 
\left(
\begin{array}{ccc}
1	&2	&3	\\
2	&3	&1	
\end{array}
\right),
\\ 
\ \\
\rho_1 = 
\left(
\begin{array}{cc}
1	&2	\\
2	&1	
\end{array}
\right),
\ 
\rho_2 = 
\left(
\begin{array}{cccc}
1	&2	&3	&4	\\
1	&4	&2	&3		
\end{array}
\right),
\ 
\rho_3 = 
\left(
\begin{array}{ccccc}
1	&2	&3	&4	&5	\\
5	&1	&2	&3	&4		
\end{array}
\right)
\end{array}
\]
(that is, $\sigma(1) = 2$, $\sigma(2) = 3$, etc), and
\[
\sigma\of(\rho_1, \rho_2, \rho_3)
=
\left(
\begin{array}{ccccccccccc}
1 &2 &3 &4 &5 &6 &7 &8 &9 &10&11\\
7 &6 &8 &11&9 &10&5 &1 &2 &3 &4
\end{array}
\right).
\]
Formally, let $\sigma\in S_n, \rho_1 \in S_{k_1}, \ldots, \rho_n \in
S_{k_n}$: then for $1\leq i\leq n$ and $1\leq j\leq k_i$, 
\[
\sigma\of(\rho_1, \ldots, \rho_n) (k_1 +\cdots + k_{i-1} + j)
=
k_{\sigma^{-1}(1)} + \cdots + k_{\sigma^{-1}(\sigma(i)-1)} + \rho_i(j).
\]
This gives $\SymOpd$ the structure of an operad.

A different construction takes the $n$-ary operations to be the total
orders%
\index{order!operad of orders}\index{operad!orders@of orders}
on the set $\{1, \ldots, n\}$ and composition to be lexicographic
combination.  In this formulation it is clear that composition is
associative and unital.  This operad of total orders is isomorphic to the
operad of symmetries---but for the proof, beware that you need to use the
right one out of the two obvious bijections between $\{\textrm{total orders
on } \{1, \ldots, n\} \}$ and $S_n$.  It is also homotopy equivalent, in a
suitable sense, to the little intervals%
\index{operad!little intervals}
operad $\ldisks_1$.
\end{example}

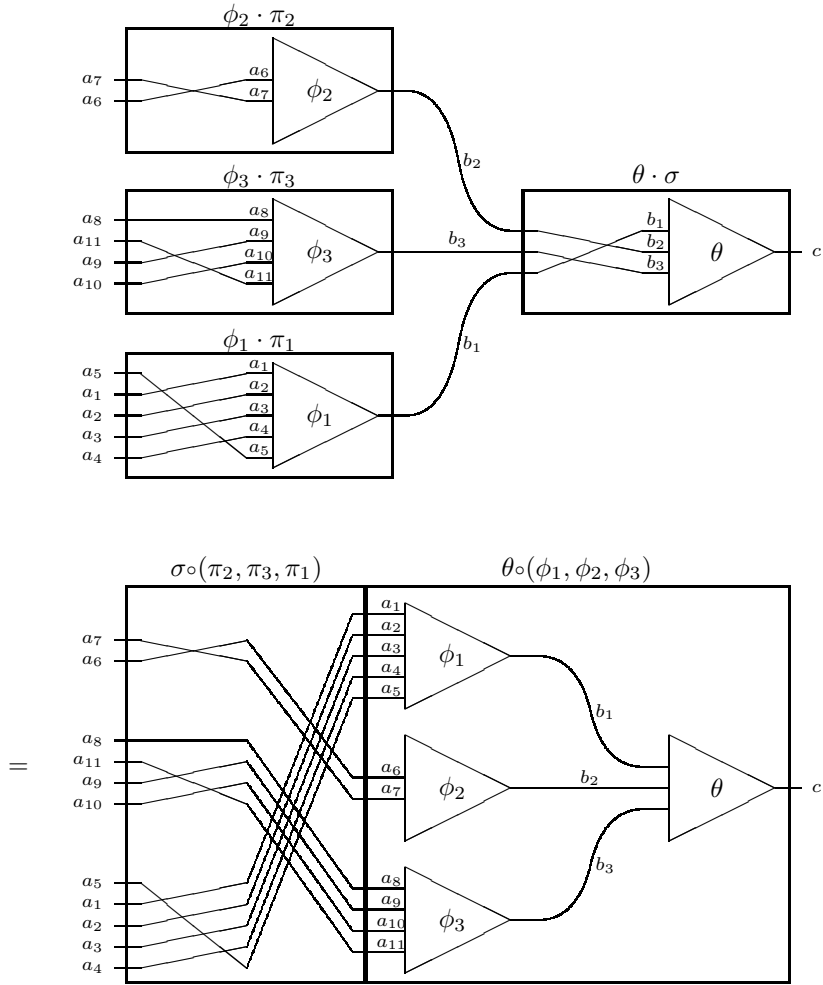
\begin{figure}%
\centering
\setlength{\unitlength}{1em}%
\begin{picture}(31,37)(-1,0)
\cell{2}{37}{tl}{%
\begin{picture}(28,18.1)(0,-0.5)
% Leftmost input tags
%   top box
\cell{2}{14.7}{r}{\tinputlft{{\scriptstyle a_7}}}
\cell{2}{13.9}{r}{\tinputlft{{\scriptstyle a_6}}}
%   middle box
\cell{2}{9.4}{r}{\tinputlft{{\scriptstyle a_8}}}
\cell{2}{8.6}{r}{\tinputlft{{\scriptstyle a_{11}}}}
\cell{2}{7.8}{r}{\tinputlft{{\scriptstyle a_9}}}
\cell{2}{7.0}{r}{\tinputlft{{\scriptstyle a_{10}}}}
%   bottom box
\cell{2}{3.6}{r}{\tinputlft{{\scriptstyle a_5}}}
\cell{2}{2.8}{r}{\tinputlft{{\scriptstyle a_1}}}
\cell{2}{2.0}{r}{\tinputlft{{\scriptstyle a_2}}}
\cell{2}{1.2}{r}{\tinputlft{{\scriptstyle a_3}}}
\cell{2}{0.4}{r}{\tinputlft{{\scriptstyle a_4}}}
% Lefthand permutations
%   top box
\put(2,14.7){\line(5,-1){4}}
\put(2,13.9){\line(5,1){4}}
%   middle box
\put(2,9.4){\line(1,0){4}}
\put(2,8.6){\line(5,-2){4}}
\put(2,7.8){\line(5,1){4}}
\put(2,7.0){\line(5,1){4}}
%   bottom box
\put(2,3.6){\line(5,-4){4}}
\put(2,2.8){\line(5,1){4}}
\put(2,2.0){\line(5,1){4}}
\put(2,1.2){\line(5,1){4}}
\put(2,0.4){\line(5,1){4}}
% Input tags to lefthand transistors
%   top box
\cell{7}{14.7}{br}{\tinputabv{}}
\cell{7}{13.9}{br}{\tinputabv{}}
%   middle box
\cell{7}{9.4}{br}{\tinputabv{}}
\cell{7}{8.6}{br}{\tinputabv{}}
\cell{7}{7.8}{br}{\tinputabv{}}
\cell{7}{7.0}{br}{\tinputabv{}}
%   bottom box
\cell{7}{3.6}{br}{\tinputabv{}}
\cell{7}{2.8}{br}{\tinputabv{}}
\cell{7}{2.0}{br}{\tinputabv{}}
\cell{7}{1.2}{br}{\tinputabv{}}
\cell{7}{0.4}{br}{\tinputabv{}}
% Labels on those tags
%   top box
\cell{6.5}{14.8}{b}{{\scriptstyle a_6}}
\cell{6.5}{14.0}{b}{{\scriptstyle a_7}}
%   middle box
\cell{6.5}{9.5}{b}{{\scriptstyle a_8}}
\cell{6.5}{8.7}{b}{{\scriptstyle a_9}}
\cell{6.5}{7.9}{b}{{\scriptstyle a_{10}}}
\cell{6.5}{7.1}{b}{{\scriptstyle a_{11}}}
%   bottom box
\cell{6.5}{3.7}{b}{{\scriptstyle a_1}}
\cell{6.5}{2.9}{b}{{\scriptstyle a_2}}
\cell{6.5}{2.1}{b}{{\scriptstyle a_3}}
\cell{6.5}{1.3}{b}{{\scriptstyle a_4}}
\cell{6.5}{0.5}{b}{{\scriptstyle a_5}}
% Lefthand transistors
\cell{7}{14.3}{l}{\tusual{\phi_2}}
\cell{7}{8.2}{l}{\tusual{\phi_3}}
\cell{7}{2.0}{l}{\tusual{\phi_1}}
% Output tags on lefthand transistors
\cell{11}{14.3}{l}{\toutputrgt{}}
\cell{11}{8.2}{l}{\toutputrgt{}}
\cell{11}{2.0}{l}{\toutputrgt{}}
% Joining curves in central column of picture
\qbezier(12,14.3)(13.5,14.3)(14,11.65)
\qbezier(16,9.0)(14.5,9.0)(14,11.65)
\put(12,8.2){\line(1,0){4}}
\qbezier(12,2.0)(13.5,2.0)(14,4.7)
\qbezier(16,7.4)(14.5,7.4)(14,4.7)
% Labels on those curves
\cell{14.2}{11.65}{l}{{\scriptstyle b_2}}
\cell{14}{8.4}{b}{{\scriptstyle b_3}}
\cell{14.2}{4.7}{l}{{\scriptstyle b_1}}
% Tags just left of righthand permutation
\cell{17}{9.0}{r}{\tinputlft{}}
\cell{17}{8.2}{r}{\tinputlft{}}
\cell{17}{7.4}{r}{\tinputlft{}}
% Righthand permutation
\put(17,9.0){\line(5,-1){4}}
\put(17,8.2){\line(5,-1){4}}
\put(17,7.4){\line(5,2){4}}
% Input tags on righthand transistor
\cell{22}{9.0}{r}{\tinputlft{}}
\cell{22}{8.2}{r}{\tinputlft{}}
\cell{22}{7.4}{r}{\tinputlft{}}
% Labels on those tags
\cell{21.5}{9.1}{b}{{\scriptstyle b_1}}
\cell{21.5}{8.3}{b}{{\scriptstyle b_2}}
\cell{21.5}{7.5}{b}{{\scriptstyle b_3}}
% Righthand transistor
\cell{22}{8.2}{l}{\tusual{\theta}}
% Output tag on righthand transistor
\cell{26}{8.2}{l}{\toutputrgt{{\scriptstyle c}}}
% Boxes
\thicklines
\put(1.5,12.0){\framebox(10,4.6){}}
\put(1.5,5.9){\framebox(10,4.6){}}
\put(1.5,-0.3){\framebox(10,4.6){}}
\put(16.5,5.9){\framebox(10,4.6){}}
% Labels on boxes
\cell{6.5}{17.6}{t}{\phi_2 \cdot \pi_2}
\cell{6.5}{11.5}{t}{\phi_3 \cdot \pi_3}
\cell{6.5}{5.3}{t}{\phi_1 \cdot \pi_1}
\cell{21.5}{11.5}{t}{\theta \cdot \sigma}
\end{picture}}
\cell{-1}{8.05}{l}{=}
\cell{2}{0}{bl}{%
\begin{picture}(28,16.1)
% Leftmost input tags
%   top group
\cell{2}{12.9}{r}{\tinputlft{{\scriptstyle a_7}}}
\cell{2}{12.1}{r}{\tinputlft{{\scriptstyle a_6}}}
%   middle group
\cell{2}{9.1}{r}{\tinputlft{{\scriptstyle a_8}}}
\cell{2}{8.3}{r}{\tinputlft{{\scriptstyle a_{11}}}}
\cell{2}{7.5}{r}{\tinputlft{{\scriptstyle a_9}}}
\cell{2}{6.7}{r}{\tinputlft{{\scriptstyle a_{10}}}}
%   bottom group
\cell{2}{3.7}{r}{\tinputlft{{\scriptstyle a_5}}}
\cell{2}{2.9}{r}{\tinputlft{{\scriptstyle a_1}}}
\cell{2}{2.1}{r}{\tinputlft{{\scriptstyle a_2}}}
\cell{2}{1.3}{r}{\tinputlft{{\scriptstyle a_3}}}
\cell{2}{0.5}{r}{\tinputlft{{\scriptstyle a_4}}}
% Permutation
%   lines for top group
\put(2,12.9){\line(5,-1){4}}
\put(2,12.1){\line(5,1){4}}
\multiput(6,12.9)(0,-0.8){2}{\qbezier(0,0)(2,-2.6)(4,-5.2)}
%   lines for middle group
\put(2,9.1){\line(1,0){4}}
\put(2,8.3){\line(5,-2){4}}
\put(2,7.5){\line(5,1){4}}
\put(2,6.7){\line(5,1){4}}
\multiput(6,9.1)(0,-0.8){4}{\qbezier(0,0)(2,-2.8)(4,-5.6)}
%   lines for bottom group
\put(2,3.7){\line(5,-4){4}}
\put(2,2.9){\line(5,1){4}}
\put(2,2.1){\line(5,1){4}}
\put(2,1.3){\line(5,1){4}}
\put(2,0.5){\line(5,1){4}}
\multiput(6,3.7)(0,-0.8){5}{\qbezier(0,0)(2,5.1)(4,10.2)}
% Tags to right of permutation
\multiput(10,13.9)(0,-0.8){5}{\line(1,0){2}}
\multiput(10,7.7)(0,-0.8){2}{\line(1,0){2}}
\multiput(10,3.5)(0,-0.8){4}{\line(1,0){2}}
% Labels on those tags
%   top transistor
\cell{11.5}{14.0}{b}{{\scriptstyle a_1}}
\cell{11.5}{13.2}{b}{{\scriptstyle a_2}}
\cell{11.5}{12.4}{b}{{\scriptstyle a_3}}
\cell{11.5}{11.6}{b}{{\scriptstyle a_4}}
\cell{11.5}{10.8}{b}{{\scriptstyle a_5}}
%   middle transistor
\cell{11.5}{7.8}{b}{{\scriptstyle a_6}}
\cell{11.5}{7.0}{b}{{\scriptstyle a_7}}
%   bottom transistor
\cell{11.5}{3.6}{b}{{\scriptstyle a_8}}
\cell{11.5}{2.8}{b}{{\scriptstyle a_9}}
\cell{11.5}{2.0}{b}{{\scriptstyle a_{10}}}
\cell{11.5}{1.2}{b}{{\scriptstyle a_{11}}}
% Lefthand transistors
\cell{12}{12.3}{l}{\tusual{\phi_1}}
\cell{12}{7.3}{l}{\tusual{\phi_2}}
\cell{12}{2.3}{l}{\tusual{\phi_3}}
% Output tags of lefthand transistors
\cell{16}{12.3}{l}{\toutputrgt{}}
\cell{16}{7.3}{l}{\toutputrgt{}}
\cell{16}{2.3}{l}{\toutputrgt{}}
% Joining curves between transistors
\qbezier(17,12.3)(18.5,12.3)(19,10.2)
\qbezier(21,8.1)(19.5,8.1)(19,10.2)
\put(17,7.3){\line(1,0){4}}
\qbezier(17,2.3)(18.5,2.3)(19,4.4)
\qbezier(21,6.5)(19.5,6.5)(19,4.4)
% Labels on those curves
\cell{19.2}{10.2}{l}{{\scriptstyle b_1}}
\cell{19}{7.4}{b}{{\scriptstyle b_2}}
\cell{19.2}{4.4}{l}{{\scriptstyle b_3}}
% Input tags to righthand transistor
\cell{22}{8.1}{r}{\tinputlft{}}
\cell{22}{7.3}{r}{\tinputlft{}}
\cell{22}{6.5}{r}{\tinputlft{}}
% Righthand transistor
\cell{22}{7.3}{l}{\tusual{\theta}}
% Output tag of righthand transistor
\cell{26}{7.3}{l}{\toutputrgt{{\scriptstyle c}}}
% Boxes
\thicklines
\put(1.5,0){\framebox(9,14.9){}}
\put(10.5,0){\framebox(16,14.9){}}
% Labels on boxes
\cell{6}{16.1}{t}{\sigma \of (\pi_2,\pi_3,\pi_1)}
\cell{18.5}{16.1}{t}{\theta \of (\phi_1, \phi_2, \phi_3)}
\end{picture}}
\end{picture}
% \hand{150}{59}
\caption{Symmetric multicategory axiom}
\label{fig:sym-mti-axiom}
\end{figure}

As discussed on p.~\pageref{p:sym-mti-informal}, a symmetric structure on a
multicategory $C$ should consist of a map
\begin{equation}	\label{eq:sym-action}
\dashbk\cdot\sigma:
C(a_1, \ldots, a_n; b)
\go
C(a_{\sigma(1)}, \ldots, a_{\sigma(n)}; b)
\end{equation}
for each $a_1, \ldots, a_n, b \in C_0$ and $\sigma\in S_n$.  These maps
should satisfy the obvious axioms
\begin{equation}	\label{eq:sym-axioms-obvious}
(\theta \cdot \sigma) \cdot \rho = \theta \cdot (\sigma\rho),
\diagspace
\theta = \theta \cdot 1_{S_n}
\end{equation}
($\theta\in C(a_1, \ldots, a_n; b)$, $\sigma, \rho \in S_n$), which
guarantee that $\dashbk\cdot\sigma$ is a bijection.  The symmetric action
should also be compatible with composition in $C$, as for instance in
Fig.~\ref{fig:sym-mti-axiom}.
In general, we want
\begin{eqnarray}
&&(\theta\cdot \sigma) \of 
(\phi_{\sigma(1)} \cdot \pi_{\sigma(1)}, \ldots, 
\phi_{\sigma(n)} \cdot \pi_{\sigma(n)})
\nonumber\\
&=&
(\theta \of (\phi_1, \ldots, \phi_n))
\cdot
(\sigma \of (\pi_{\sigma(1)}, \ldots, \pi_{\sigma(n)}))
\label{eq:sym-mti-axiom}
\end{eqnarray}
whenever $\theta, \phi_1, \ldots, \phi_n$ are maps in $C$ and $\sigma,
\pi_1, \ldots, \pi_n$ are permutations for which these expressions make
sense.  The permutation $\sigma \of (\pi_{\sigma(1)}, \ldots,
\pi_{\sigma(n)})$ on the right-hand side is the composite in $\SymOpd$.
(It is easier to get axiom~\bref{eq:sym-mti-axiom} right for
multicategories than for the special case of operads---the different
objects should stop us from writing down nonsense.)

\begin{defn}	\lbl{defn:sym-mti}
A \demph{symmetric multicategory} is a multicategory $C$ together with a
map~\bref{eq:sym-action} for each $a_1, \ldots, a_n, b \in C_0$ and
$\sigma\in S_n$, satisfying the axioms in~\bref{eq:sym-axioms-obvious}
and~\bref{eq:sym-mti-axiom}.  A \demph{map of symmetric multicategories} is
a map $f$ of multicategories such that $f(\theta\cdot\sigma) = f(\theta)
\cdot\sigma$ whenever $\theta$ is an arrow of $C$ and $\sigma$ a
permutation for which this makes sense.  The category of symmetric
multicategories is written $\fcat{SymMulticat}$.%
\glo{SymMulticat}
 A \demph{symmetric
operad} is a one-object symmetric multicategory.
\end{defn}

Any symmetric monoidal category is naturally a symmetric multicategory,%
\index{multicategory!underlying}
via the symmetry maps
\[
\sigma\cdot\dashbk:
a_{\sigma(1)} \otimes\cdots\otimes a_{\sigma(n)}
\goiso
a_1 \otimes\cdots\otimes a_n.
\]
This is true in particular of the category of sets.  An \demph{algebra}%
\lbl{p:defn-sym-alg}%
\index{algebra!symmetric multicategory@for symmetric multicategory}%
\index{multicategory!symmetric!algebra for}
for a symmetric multicategory $C$ is a map $C \go \Set$ of symmetric
multicategories.  In general, $C$ has more algebras when regarded as a
non-symmetric multicategory than when the symmetries are taken into
account.

An equivalent definition of symmetric multicategory is given in
Appendix~\ref{app:sym}: `fat symmetric multicategories', in many ways more
graceful.

\begin{example}
The operad $\SymOpd$%
\glo{SymOpdsym}
of symmetries%
\index{operad!symmetries@of symmetries}
becomes a symmetric operad by
multiplication in the symmetric groups. 
\end{example}

\begin{example}	\lbl{eg:sym-multi-for-opds}%
\index{multicategory!symmetric!operads@for operads}
There is a symmetric multicategory $\cat{O}$ whose algebras (as a
\emph{symmetric} multicategory) are exactly operads.  This example will be
done informally; we replace it with a precise construction later.

The objects of $\cat{O}$ are the natural numbers.  To define the arrows we
use finite, rooted, planar trees%
\index{tree!vertices ordered@with vertices ordered}
in which each vertex may have any natural
number of branches (including $0$) coming up out of it.  An element of
$\cat{O}(m_1, \ldots, m_k; n)$ is an $n$-leafed tree with $k$ vertices that
are totally ordered in such a way that the $i$th vertex has $m_i$ branches
coming up out of it.  For example,
\begin{equation}	\label{eq:223}
\cat{O}(2,2;3) 
= 
\left\{
\setlength{\unitlength}{.5em}
% 
% TREE 1
% 
\begin{array}{c}
\begin{picture}(6,6)(0,0)
% bottom layer
\put(4,0){\line(0,1){2}}
\cell{4}{2}{c}{\vx}
% middle layer
\put(4,2){\line(-1,1){2}}
\put(4,2){\line(1,1){2}}
\cell{2}{4}{c}{\vx}
% top layer
\put(2,4){\line(-1,1){2}}
\put(2,4){\line(1,1){2}}
% labels
\cell{3.5}{2}{r}{\scriptstyle 1}
\cell{1.5}{4}{r}{\scriptstyle 2}
\end{picture}
\end{array},
% 
% TREE 2
% 
\begin{array}{c}
\begin{picture}(6,6)(0,0)
% bottom layer
\put(4,0){\line(0,1){2}}
\cell{4}{2}{c}{\vx}
% middle layer
\put(4,2){\line(-1,1){2}}
\put(4,2){\line(1,1){2}}
\cell{2}{4}{c}{\vx}
% top layer
\put(2,4){\line(-1,1){2}}
\put(2,4){\line(1,1){2}}
% labels
\cell{3.5}{2}{r}{\scriptstyle 2}
\cell{1.5}{4}{r}{\scriptstyle 1}
\end{picture}
\end{array},
% 
% TREE 3
% 
\begin{array}{c}
\begin{picture}(6,6)(0,0)
% bottom layer
\put(2,0){\line(0,1){2}}
\cell{2}{2}{c}{\vx}
% middle layer
\put(2,2){\line(-1,1){2}}
\put(2,2){\line(1,1){2}}
\cell{4}{4}{c}{\vx}
% top layer
\put(4,4){\line(-1,1){2}}
\put(4,4){\line(1,1){2}}
% labels
\cell{2.5}{2}{l}{\scriptstyle 1}
\cell{4.5}{4}{l}{\scriptstyle 2}
\end{picture}
\end{array},
% 
% TREE 4
% 
\begin{array}{c}
\begin{picture}(6,6)(0,0)
% bottom layer
\put(2,0){\line(0,1){2}}
\cell{2}{2}{c}{\vx}
% middle layer
\put(2,2){\line(-1,1){2}}
\put(2,2){\line(1,1){2}}
\cell{4}{4}{c}{\vx}
% top layer
\put(4,4){\line(-1,1){2}}
\put(4,4){\line(1,1){2}}
% labels
\cell{2.5}{2}{l}{\scriptstyle 2}
\cell{4.5}{4}{l}{\scriptstyle 1}
\end{picture}
\end{array}
\right\}.
% \hand{12}{60}.
\end{equation}
Composition is substitution of trees into vertices (much as in the little
disks operad), the identity on $n$ is the $n$-leafed tree with only one
vertex, and the symmetric group action is by permutation of the order of the
vertices.  An $\cat{O}$-algebra consists of a set $P(n)$ for each
$n\in\nat$ together with a map
\[
P(m_1) \times\cdots\times P(m_k) \go P(n)
\]
for each element of $\cat{O}(m_1, \ldots, m_k; n)$, satisfying axioms, and
this is exactly an operad.  For example, the first element $\alpha$ of
$\cat{O}(2,2;3)$ listed in~\bref{eq:223} induces the function
\[
\begin{array}{rrcl}
\ovln{\alpha}:	&P(2) \times P(2)	&\go	&P(3)	\\
		&(\theta, \theta')	&\goesto&
\theta \of (\theta', 1),	
\end{array}
\]
part of the operadic structure of $P$.  If $\sigma$ is the nontrivial
element of $S_2$ then the second element of $\cat{O}(2,2;3)$ listed
in~\bref{eq:223} is $\alpha\cdot\sigma$, and a consequence of $P$ being an
algebra for $\cat{O}$ as a \emph{symmetric} multicategory is that
$\ovln{\alpha\cdot\sigma} (\theta', \theta) = \ovln{\alpha} (\theta,
\theta')$.

Similarly, there is a symmetric multicategory $\cat{O}'$ whose algebras are
symmetric operads; it is the same as $\cat{O}$ except that the trees are
equipped with an ordering of the leaves as well as the vertices.  And more
generally, for any set $S$ there are symmetric multicategories $\cat{O}_S$
and $\cat{O}'_S$ whose algebras are, respectively, non-symmetric and
symmetric multicategories with object-set $S$; the object-sets of both
$\cat{O}_S$ and $\cat{O}'_S$ are $(\coprod_{n\in\nat} S^n) \times S$.

There is no symmetric multicategory whose algebras are \emph{all}
multicategories.  There is, however,%
\index{multicategory!symmetric vs. generalized@symmetric \vs.\ generalized!multicategory for multicategories}
a \emph{generalized} multicategory
with this property, as we shall see.
\end{example}%
\index{operad!symmetric|)}\index{multicategory!symmetric|)}

\section{Further theory}
\lbl{sec:om-further}

So far we have seen the basic definitions in the theory of multicategories
and operads, and some examples.  Here we consider a few further topics.
Most are special cases of constructions for generalized multicategories
that we meet later; some are generalizations of concepts familiar for
ordinary categories.

We start with two alternative ways of defining (operad and) multicategory,
both in use; they go by the names of `circle-$i$' ($\of_i$) and `PROs'
respectively.  Then we extend three concepts of category theory to
multicategories: the free category on a graph, transformations between maps
between categories, and modules over categories.

\minihead{Circle-$i$}%
\index{circle-i@$\of_i$ (`circle-$i$')}%
\index{composition!circle-i@$\of_i$ (`circle-$i$')}

The `circle-$i$' method takes composition of diagrams of shape
\[
\begin{centredpic}
\begin{picture}(15.5,12)(-1,-6)
% leftmost tags
\cell{2}{4}{r}{\tinputslft{b_1}{b_{i-1}}}
\cell{2}{0}{r}{\tinputslft{a_1}{a_n}}
\cell{2}{-4}{r}{\tinputslft{b_{i+1}}{b_m}}
% lefthand transistor
\cell{2}{0}{l}{\tusual{\phi}}
% joining wires
\qbezier(2,5.5)(3.875,5.5)(4.5,4.75)
\qbezier(7,4)(5.125,4)(4.5,4.75)
\qbezier(2,2.5)(3.875,2.5)(4.5,1.75)
\qbezier(7,1)(5.125,1)(4.5,1.75)
\put(6,0){\line(1,0){2}}
\cell{6.8}{0.1}{b}{b_i}
\qbezier(2,-5.5)(3.875,-5.5)(4.5,-4.75)
\qbezier(7,-4)(5.125,-4)(4.5,-4.75)
\qbezier(2,-2.5)(3.875,-2.5)(4.5,-1.75)
\qbezier(7,-1)(5.125,-1)(4.5,-1.75)
% inputs to righthand transistor
\cell{8}{2.5}{r}{\tinputslft{}{}}
\cell{8}{-2.5}{r}{\tinputslft{}{}}
% righthand transistor
\put(8,4.5){\line(0,-1){9}}
\put(8,4.5){\line(1,-1){4.5}}
\put(8,-4.5){\line(1,1){4.5}}
\cell{10}{0}{c}{\theta}
% output tag
\cell{12.5}{0}{l}{\toutputrgt{c}}
\end{picture}
\end{centredpic},
% \hand{50}{61},
\]
along with identities, to be the basic operations in a multicategory.  Thus,
a multicategory $C$ can be defined as a set $C_0$ of objects together
with hom-sets $C(a_1, \ldots, a_n; a)$, a function
\[
\begin{array}{rrcl}
\begin{array}[b]{r}\of_i:\\ \,\end{array}	&
\begin{array}[b]{l}
C(b_1, \ldots, b_m; c) \\
\times C(a_1, \ldots, a_n; b_i)	
\end{array}
&
\go	&
\begin{array}[t]{l}
C(b_1, \ldots, b_{i-1}, a_1, \ldots, a_n, \\
b_{i+1}, \ldots, b_m; c)	
\end{array}
\\
	&
(\theta, \phi)	&
\goesto	&
\begin{array}{l}\theta \ofdim{i} \phi\end{array} 
\end{array}
\]%
\glo{circlei}%
for each $1 \leq i \leq m$, $n\in\nat$, $a_1, \ldots, a_n, b_1, \ldots, b_m
\in C_0$, and an element $1_a \in C(a;a)$ for each $a\in C_0$, satisfying
certain axioms.  This is an equivalent definition: given a multicategory in
the usual sense, we put
\[
\theta \of_i \phi =
(1_{b_1}, \ldots, 1_{b_{i-1}}, \phi, 1_{b_{i+1}}, \ldots, 1_{b_m}),
\]
and given a multicategory in the new sense, the composite maps $\theta\of
(\phi_1, \ldots, \phi_n)$ can be built using $n$ operations of the form
$\of_i$.  

This was, in fact, Lambek's%
\index{Lambek, Joachim}
original definition of
multicategory~\cite[p.~103]{LamDSCII}.  His motivating example was that of
a deductive%
\index{deductive system}
system: objects are statements, maps $a_1, \ldots, a_n \go b$
are deductions of $b$ from $a_1, \ldots, a_n$, and the $\of_i$ operation is
Gentzen cut.%
\index{cut}
 The style of definition is also useful if for some reason one
does not want one's multicategories to have identities%
\lbl{p:pseudo-operads}\index{pseudo-operad}%
\index{operad!pseudo-}
(as in Markl, Shnider and Stasheff~\cite[p.~45]{MSS}): for with all the
$\of_i$'s (but not identities) one can build an operation for composing
diagrams in the shape of any non-trivial tree, whereas with the usual
$\of$'s (but not identities) one only obtains the non-trivial trees whose
leaves are at uniform height.  We stick firmly to the original definition,
as that is what is generalized to give the all-important definition of
generalized multicategory.  The $\of_i$ definition does not generalize in
the same way: consider, for instance, the $\fc$-multicategories of
Chapter~\ref{ch:fcm}.

\minihead{PROs and PROPs}%
\index{monoidal category!multicategory@\vs.\ multicategory|(}

To reach the second alternative definition of multicategory we consider how
multicategories are related to strict monoidal categories.  (The weak case
is left until~\ref{sec:non-alg-notions}.)  As we saw in~\ref{eg:multi-mon},
there is a forgetful functor
\[
\fcat{StrMonCat} \goby{U} \Multicat
\]%
\glo{StrMonCat}%
where the domain is the category of strict monoidal categories and strict
monoidal functors.  This has a left adjoint 
\[
\Multicat \goby{F} \fcat{StrMonCat}.
\]
Given a multicategory $C$, the objects (respectively, arrows) of $F(C)$ are
finite ordered sequences of objects (respectively, arrows) of $C$, and the
tensor product in $F(C)$ is concatenation of sequences.  So a typical arrow
\[
(a_1,a_2,a_3,a_4,a_5) \go (a'_1,a'_2,a'_3)
\]
in $F(C)$ looks like
\begin{equation}	\label{diag:arrows-in-mon-cat} 
\begin{centredpic}
\begin{picture}(7.2,12.4)(0,-0.2)
% input tags
\cell{2}{11.2}{r}{\tinputlft{a_1}}
\cell{2}{10}{r}{\tinputlft{a_2}}
\cell{2}{8.8}{r}{\tinputlft{a_3}}
\cell{2}{3.0}{r}{\tinputlft{a_4}}
\cell{2}{1.0}{r}{\tinputlft{a_5}}
% transistors
\multiput(2,8.4)(0,-4){3}{\line(0,1){3.2}}
\multiput(5.2,10)(0,-4){3}{\line(-2,-1){3.2}}
\multiput(5.2,10)(0,-4){3}{\line(-2,1){3.2}}
% labels
\cell{3.4}{10}{c}{\theta_1}
\cell{3.4}{6}{c}{\theta_2}
\cell{3.4}{2}{c}{\theta_3}
% output tags
\cell{5.2}{10}{l}{\toutputrgt{a'_1}}
\cell{5.2}{6}{l}{\toutputrgt{a'_2}}
\cell{5.2}{2}{l}{\toutputrgt{a'_3}}
% box
\thicklines
\put(1.5,0){\framebox(4.2,12){}}
\end{picture}
\end{centredpic}
\end{equation}
where $\theta_1: a_1, a_2, a_3 \go a'_1$ in $C$, etc.  An arrow $(a_1,
\ldots, a_n) \go (a)$ in $F(C)$ is simply an arrow $a_1, \ldots, a_n \go a$
in $C$.

The monoidal categories that arise freely from multicategories can be
characterized intrinsically, and this makes it possible to redefine a
multicategory as a monoidal category with certain properties.  Let $C$ be a
multicategory.  First, the strict monoidal category $F(C)$ has the property
that its underlying monoid of objects is the free monoid on the set $C_0$.
Second, let $S$ be a set and let $A$ be a strict monoidal category whose
monoid of objects is the free monoid on $S$; then for any elements $b_1,
\ldots, b_m, a_1, \ldots, a_n$ of $S$, tensor product in $A$ defines a map
\begin{equation}	\label{eq:PRO-map}
\begin{array}{r}\displaystyle
\displaystyle\coprod_{a_1^1, \ldots, a_n^{k_n}}
\left(
A((a_1^1, \ldots, a_1^{k_1}), (a_1))
\times\cdots\times
A((a_n^1, \ldots, a_n^{k_n}), (a_n))
\right)
\\
\go
A((b_1, \ldots, b_m), (a_1, \ldots, a_n))
\end{array}
\end{equation}
where the union is over all $n, k_1, \ldots, k_n \in\nat$ and $a_i^j \in S$
such that there is an equality of formal sequences
\[
(a_1^1, \ldots, a_1^{k_1}, \ldots, a_n^1, \ldots, a_n^{k_n})
=
(b_1, \ldots, b_m).
\]
The crucial point is that when $A=F(C)$, the map~\bref{eq:PRO-map} is
always a bijection.

A \demph{PRO}%
\index{PRO}
is a pair $(S,A)$ where $S$ is a set and $A$ is a strict
monoidal category, such that
\begin{itemize}
\item the monoid of objects of $A$ is equal to the free monoid on $S$
\item for all $b_1, \ldots, b_m, a_1, \ldots, a_n \in S$, the canonical
  map~\bref{eq:PRO-map} is a bijection.
\end{itemize}
A \demph{map of PROs} $(u,f): (S, A) \go (S', A')$ is a function $u: S \go
S'$ together with a strict monoidal functor $f: A \go A'$ such that the
objects-function $f_0: A_0 \go A'_0$ is the result of applying the free
monoid functor to $u$.  This gives a category $\fcat{PRO}$ of PROs.  There
is a forgetful functor $\fcat{PRO} \go \fcat{StrMonCat}$, and $F$ lifts in
the obvious way to give a functor $\twid{F}: \Multicat \go \fcat{PRO}$.
\begin{propn}
The functor $\twid{F}: \Multicat \go \fcat{PRO}$ is an equivalence.
\end{propn}
\begin{proof}
Given a PRO $(S,A)$, there is a multicategory $C$ with $C_0 = S$ and
\[
C(a_1, \ldots, a_n; a) = 
A((a_1, \ldots, a_n), (a)),
\]
and the condition that the maps~\bref{eq:PRO-map} are bijections implies
that $\twid{F}(C) \iso (S,A)$.  The rest of the proof is straightforward.
\done
\end{proof}

The same kind of equivalence can be established for symmetric
multicategories and monoidal categories.  The symmetric analogue of a PRO
is a PROP.  These structures were introduced by Adams%
\index{Adams, Frank}
and Mac%
\index{Mac Lane, Saunders}
Lane (Mac
Lane~\cite{MacNAC}) and developed by Boardman%
\index{Boardman, Michael}
and Vogt~\cite{BV};%
\index{Vogt, Rainer}
the names
stand for `PROduct (and Permutation) category'.  Boardman and Vogt called a
pair $(S,A)$ an `$S$-coloured PRO(P)', and paid particular attention to the
single-coloured case.  Spelling it out: the category $\fcat{sPRO}$ of
single-coloured PROs has as objects those strict monoidal categories $A$
whose underlying monoid of objects is equal to $(\nat,+,0)$ and for which
the canonical map
\[
\coprod_{k_1 + \cdots + k_n = m}
A(k_1, 1) \times\cdots\times A(k_n, 1)
\go 
A(m,n)
\]
is a bijection for all $m,n\in\nat$, and as arrows those strict monoidal
functors that are the identity on objects.  We have immediately:
\begin{cor}
The functor $\twid{F}$ restricts to an equivalence of categories %\linebreak
$\Operad \go \fcat{sPRO}$.  
\done
\end{cor}%
\index{monoidal category!multicategory@\vs.\ multicategory|)}

\minihead{Free multicategories}%
\index{multicategory!free}

Free structures are the formal origin of much of the geometry in this
subject.  A basic case is that free multicategories are made out of trees,
as now explained.

A \demph{multigraph}%
\index{multigraph}
is a set $X_0$ together with a set $X(a_1, \ldots,
a_n; a)$ for each $n\in\nat$ and $a_1, \ldots, a_n, a \in X_0$.  Forgetting
composition and identities gives a functor $U: \Multicat \go
\fcat{Multigraph}$.  This has a left adjoint $F$, the free%
\lbl{p:free-mti-ftr}
multicategory functor, which can be described as follows.  Let $X$ be a
multigraph.  The free multicategory $FX$ on $X$ has the same objects:
$(FX)_0 = X_0$.  Its arrows are formal gluings of arrows of $X$, that is, the
hom-sets of $FX$ are generated recursively by the clauses
\begin{itemize}
\item\label{p:free-plain-clauses}
if $a \in X_0$ then $1_a \in (FX)(a;a)$
\item if $\xi\in X(a_1, \ldots, a_n; a)$ and 
\[
\theta_1 \in (FX)(a_1^1, \ldots, a_1^{k_1}; a_1), 
\ \ldots,\  
\theta_n \in (FX)(a_n^1, \ldots, a_n^{k_n}; a_n)
\]
then $\xi \of (\theta_1, \ldots, \theta_n) \in
(FX)(a_1^1, \ldots, a_n^{k_n}; a)$. 
\end{itemize}
Here $1_a$ and $\xi \of (\theta_1, \ldots, \theta_n)$ are just formal
expressions, but also make it clear how identities and composition in $FX$
are to be defined.  A typical arrow in $FX$ is
\[
\xi_1 \, \of \,
(\xi_2 \of (1_{a_3}, 1_{a_4}), \,
\xi_3 \of (\xi_4 \of (), 
	\xi_5 \of (1_{a_8}, 1_{a_9}, 1_{a_{10}})), \,
1_{a_{11}})
\]
where
\[
\xi_1 \in X(a_2, a_5, a_{11}; a_1),
\diagspace
\xi_2 \in X(a_3, a_4; a_2),
\]
and so on, naturally drawn as in Fig.~\ref{fig:random-multi-diagram}
(p.~\pageref{fig:random-multi-diagram}) with $\xi_i$'s in place of
$\theta_i$'s.  The multigraph $X$ is embedded in $FX$ by sending $\xi \in
X(a_1, \ldots, a_n; a)$ to
\[
\xi\of (1_{a_1}, \ldots, 1_{a_n}) \in (FX)(a_1, \ldots, a_n; a).
\]

\begin{example}	\lbl{eg:opd-of-trees}
The \demph{operad $\tr$%
\glo{tr}
of trees}%
\index{tree!operad of}\index{operad!trees@of trees}
is defined as $F1$, the free
multicategory on the terminal multigraph.  Explicitly, the sets $\tr(n)$
($n\in\nat$) are generated recursively by
\begin{itemize}
\item $\tr(1)$ has an element $\utree$%
\glo{utree}
\item if $n, k_1, \ldots, k_n \in\nat$ and $\tau_1 \in \tr(k_1), \ldots,
\tau_n\in\tr(k_n)$, then $\tr(k_1 + \cdots + k_n)$ has an element $(\tau_1,
\ldots, \tau_n)$.%
\glo{tupleoftrees}
\end{itemize}
Here $\utree$ is a formal symbol and $(\tau_1, \ldots, \tau_n)$ a formal
$n$-tuple.  The elements of $\tr(n)$ are called \demph{$n$-leafed trees},
and drawn as diagrams with $n$ edges coming into the top and one edge (the
\demph{root})%
\index{root of tree}
emerging from the bottom, in the following way:
\begin{itemize}
\item $\utree\in\tr(1)$ is drawn as $\utree$
\item if $\tau_1 \in \tr(k_1), \ldots, \tau_n \in \tr(k_n)$, and if
%   \drk{vertical transistor labelled with $\tau_i$} 
\[
\setlength{\unitlength}{1em}
\begin{picture}(3,5)(-1.5,0)
\cell{0}{4}{b}{\tinputssmallvert{}{}}
\cell{0}{4}{t}{\tsmallvert{\tau_i}}
\cell{0}{1}{t}{\toutputvert{}}
\end{picture}
\]
represents the diagram of $\tau_i$, then the tree $(\tau_1, \ldots,
\tau_n)$ is drawn as
\begin{equation}	\label{diag:joined-trees}
\begin{centredpic}
\begin{picture}(15,6)(0,0)
% bottom layer
\put(7.5,0){\line(0,1){1}}
\cell{7.5}{1}{c}{\vx}
% middle layer
\put(7.5,1){\line(-6,1){6}}
\put(7.5,1){\line(-2,1){2}}
\put(7.5,1){\line(6,1){6}}
\cell{9.5}{1.8}{c}{\cdots}
% transistor layer
\cell{1.5}{5}{t}{\tsmallvert{\tau_1}}
\cell{5.5}{5}{t}{\tsmallvert{\tau_2}}
\cell{9.5}{4}{c}{\cdots}
\cell{13.5}{5}{t}{\tsmallvert{\tau_3}}
% inputs to transistors
\cell{1.5}{5}{b}{\tinputssmallvert{}{}}
\cell{5.5}{5}{b}{\tinputssmallvert{}{}}
\cell{13.5}{5}{b}{\tinputssmallvert{}{}}
\end{picture}
\end{centredpic}.
% \hand{25}{63}.
\end{equation}
\end{itemize}
For example, $\tr(3)$ has an element $(\utree, \utree, \utree)$, drawn as
\[
% \drk{3-leafed corolla},
\begin{centredpic}
\begin{picture}(2,2)(0,0)
% lower layer
\put(1,0){\line(0,1){1}}
\cell{1}{1}{c}{\vx}
% upper layer
\put(1,1){\line(-1,1){1}}
\put(1,1){\line(0,1){1}}
\put(1,1){\line(1,1){1}}
\end{picture}
\end{centredpic},
\]
and $\tr(4)$ has an element $((\utree,\utree,\utree),\utree)$, drawn as
\[
\begin{centredpic}
\begin{picture}(3,3)(0,0)
% bottom layer
\put(2,0){\line(0,1){1}}
\cell{2}{1}{c}{\vx}
% middle layer
\put(2,1){\line(-1,1){1}}
\put(2,1){\line(1,1){1}}
\cell{1}{2}{c}{\vx}
% top layer
\put(1,2){\line(-1,1){1}}
\put(1,2){\line(0,1){1}}
\put(1,2){\line(1,1){1}}
\end{picture}
\end{centredpic}.
\]
These diagrams are just like Fig.~\ref{fig:random-multi-diagram}
(p.~\pageref{fig:random-multi-diagram}), but rotated and unlabelled.

Some special cases can trap the unwary.  For $n=1$ and $n=0$,
diagram~\bref{diag:joined-trees} looks like
\[
% \drk{pic of }(\tau_1)
\setlength{\unitlength}{1em}
\begin{array}[b]{c}
\begin{picture}(3,6)(0,0)
% bottom layer
\put(1.5,0){\line(0,1){1}}
\cell{1.5}{1}{c}{\vx}
% output of transistor
\put(1.5,1){\line(0,1){1}}
% transistor
\cell{1.5}{5}{t}{\tsmallvert{\tau_1}}
% inputs to transistor
\cell{1.5}{5}{b}{\tinputssmallvert{}{}}
\end{picture}
\end{array}
\diagspace
\textrm{ and } 
\diagspace
\begin{array}[b]{c}
\begin{picture}(0,1)(0,0)
\put(0,0){\line(0,1){1}}
\cell{0}{1}{c}{\vx}
\end{picture}
\end{array}
\]
respectively.  The tree on the right is an element of $\tr(0)$, that is,
has $0$ leaves; a leaf is an edge \emph{without} a vertex at its upper
end.  Note in particular that the trees
\[
% \drk{pics of } (), \utree, (\utree), ((\utree)), \ldots
\setlength{\unitlength}{1em}
\begin{array}[b]{c}
\begin{picture}(0,1)(0,0)
\put(0,0){\line(0,1){1}}
\cell{0}{1}{c}{\vx}
\end{picture}
\end{array},
\diagspace
\begin{array}[b]{c}
\begin{picture}(0,1)(0,0)
\put(0,0){\line(0,1){1}}
\end{picture}
\end{array},
\diagspace
\begin{array}[b]{c}
\begin{picture}(0,2)(0,0)
\put(0,0){\line(0,1){1}}
\cell{0}{1}{c}{\vx}
\put(0,1){\line(0,1){1}}
\end{picture}
\end{array},
\diagspace
\begin{array}[b]{c}
\begin{picture}(0,3)(0,0)
\put(0,0){\line(0,1){1}}
\cell{0}{1}{c}{\vx}
\put(0,1){\line(0,1){1}}
\cell{0}{2}{c}{\vx}
\put(0,2){\line(0,1){1}}
\end{picture}
\end{array},
\diagspace
\ldots
\]
are all different.  Formally, the first is the element $()$ of $\tr(0)$,
and the rest are elements of $\tr(1)$, namely, $\utree$,
$(\utree)$, $((\utree))$, \ldots.  Moral: the vertices matter.  

The embedding $1 \go F1$ picks out 
\[
\nu_n 
= 
(\utree, \ldots, \utree) 
=
% \drk{pic of corolla}
\begin{centredpic}
\begin{picture}(3,2)(-1.5,0)
% lower layer
\put(0,0){\line(0,1){1}}
\cell{0}{1}{c}{\vx}
% upper layer
\put(0,1){\line(-3,2){1.5}}
\cell{0}{1.8}{c}{\cdots}
\put(0,1){\line(3,2){1.5}}
\end{picture}
\end{centredpic}
\in 
\tr(n),
\]%
\glo{corollan}%
the \demph{$n$-leafed corolla},%
\index{corolla}
for each $n\in\nat$.  The operadic
composition in $\tr$ is `grafting' (gluing roots to leaves), with unit
$\utree\in\tr(1)$.  

As should be apparent, `tree' is used to mean finite, rooted, planar tree.
Non-planar trees arise similarly from \emph{symmetric} operads.  We will
examine planar trees in detail, including how they form a category,
in~\ref{sec:trees}.
\end{example}

\begin{example}	\lbl{eg:opd-of-cl-trees}
Let $X$ be the multigraph with a single object $\star$ and in which
$X(\star, \ldots, \star; \star)$ has one element if there are 0 or 2 copies
of $\star$ to the left of the semi-colon, and no elements otherwise.  For
reasons that will emerge in the next chapter, $FX$ is called the
\demph{operad of classical%
\index{tree!classical}\index{operad!classical trees@of classical trees}
trees}, $\ctr$.%
\glo{ctr}
 It is the sub-operad of $\tr$
containing just those trees in which each vertex has either 0 or 2 vertices
coming up out of it.  Explicitly, the sets $\ctr(n)$ ($n\in\nat$) are
generated recursively by
\begin{itemize}
\item $\ctr(1)$ has an element $\utree$
\item $\ctr(0)$ has an element 
% $\drk{pic of 0-leafed corolla}$
$\begin{centredpic}
\begin{picture}(0,1)(0,0)
\put(0,0){\line(0,1){1}}
\cell{0}{1}{c}{\vx}
\end{picture}
\end{centredpic}$
\item if $k_1, k_2 \in \nat$, $\tau_1 \in \ctr(k_1)$, and $\tau_2 \in
  \ctr(k_2)$, then $\ctr(k_1 + k_2)$ has an element $(\tau_1, \tau_2)$.
\end{itemize}
\end{example}

\minihead{Transformations}

Mac%
\index{Mac Lane, Saunders}
Lane recounts that he and Eilenberg%
\index{Eilenberg, Samuel}
started category theory in order to
enable them to talk about natural transformations.  Transformations for
multicategories will not be so important here but are still worth a look.
The definition is suggested by the definition of a map between algebras
(p.~\pageref{p:map-of-algs}).

\begin{defn}
Let $C \parpair{f}{f'} D$ be a pair of maps between multicategories.  A
\demph{transformation}%
\index{transformation!plain multicategories@of plain multicategories}%
\index{multicategory!transformation of}
$\alpha: f \go f'$ is a family $\left( f(a)
\goby{\alpha_a} f'(a) \right)_{a\in C_0}$ of unary maps in $D$, such that
(Fig.~\ref{fig:cl-transf-axiom}) 
\[
\alpha_a \of (f(\theta))
=
f'(\theta) \of (\alpha_{a_1}, \ldots, \alpha_{a_n})
\]
for every map $a_1, \ldots, a_n \goby{\theta} a$ in $C$.  
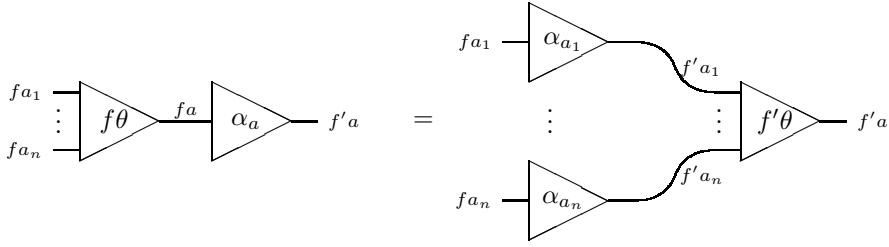
\begin{figure}%
\setlength{\unitlength}{1em}%
\centering
\begin{picture}(33.5,9)(0,-4.5)
% \cell{0}{4}{l}{\cdot}
% \cell{33.5}{4}{l}{\cdot}
\cell{0}{0}{l}{%
\begin{picture}(13,3)(-1,-1.5)
\cell{2}{0}{r}{\tinputslftsmall{\scriptstyle fa_1}{\scriptstyle fa_n}}
\cell{2}{0}{l}{\tsmall{f\theta}}
\put(5,0){\line(1,0){2}}
\cell{6}{0.1}{b}{\scriptstyle fa}
\cell{7}{0}{l}{\tsmall{\alpha_a}}
\cell{10}{0}{l}{\toutputrgt{\scriptstyle f'a}}
\end{picture}}
\cell{16}{0}{c}{=}
\cell{18}{0}{l}{%
\begin{picture}(15.5,9)(0,-4.5)
% leftmost tags
\cell{2}{3}{r}{\tinputlft{\scriptstyle fa_1}}
\cell{2}{-3}{r}{\tinputlft{\scriptstyle fa_n}}
% lefthand transistors
\cell{2}{3}{l}{\tsmall{\alpha_{a_1}}}
\cell{2.8}{0.3}{c}{\vdots}
\cell{2}{-3}{l}{\tsmall{\alpha_{a_n}}}
% output tags of lefthand transistors
\cell{5}{3}{l}{\toutputrgt{}}
\cell{5}{-3}{l}{\toutputrgt{}}
% joining curves
\qbezier(6,3)(7.125,3)(7.5,2.05)
\qbezier(9,1.1)(7.875,1.1)(7.5,2.05)
\qbezier(6,-3)(7.125,-3)(7.5,-2.05)
\qbezier(9,-1.1)(7.875,-1.1)(7.5,-2.05)
% labels on joining curves
\cell{7.7}{2.05}{l}{\scriptstyle f'a_1}
\cell{7.7}{-2.05}{l}{\scriptstyle f'a_n}
% input tags to righthand transistor
\cell{10}{0}{r}{\tinputslftsmall{}{}}
% righthand transistor
\cell{10}{0}{l}{\tsmall{f'\theta}}
% output tag
\cell{13}{0}{l}{\toutputrgt{\scriptstyle f'a}}
\end{picture}}
\end{picture}
% 
% \hand{45}{64}
\caption{Axiom for a transformation}
\label{fig:cl-transf-axiom}
\end{figure}
\end{defn}
Transformations compose in the evident ways, making $\Multicat$%
\glo{Multicat2cat}
into a
strict 2-category.%
\index{multicategory!two-category of@2-category of}
 In particular, there is a category $\ftrcat{C}{D}$%
\index{multicategory!functor categories for}%
\index{functor!category!of multicategories}
for
any multicategories $C$ and $D$, consisting of maps $C\go D$ and
transformations, and when $D=\Set$ this is $\Alg(C)$.  

Actually, $\Multicat$ is not just a 2-category: in the language of
Chapter~\ref{ch:fcm}, it is an $\fc$-multicategory.  One of the ingredients
missing in the 2-category but present in the $\fc$-multicategory is
modules, which we consider next.

\minihead{Modules}

Recall that given categories $C$ and $D$, a \demph{$(D,C)$-module}%
\lbl{p:defn-cat-module}%
\index{module!categories@over categories}
(also called a bimodule, profunctor or distributor) is a functor $X: C^\op
\times D \go \Set$.  We write $X: C \rMod D$.  When $C$ and $D$ are monoids
(one-object categories), $X$ is a set with a left $D$-action and a
compatible right $C$-action; in general, $X$ is a family $(X(c,d))_{c\in C,
d\in D}$ of sets `acted on' by the arrows of $C$ and $D$.

\begin{defn}	%\lbl{defn:cl-mti-module}
Let $C$ and $D$ be multicategories.  A \demph{$(D,C)$-module}%
\index{module!multicategories@over multicategories}
$X$, written
$X: C \rMod D$, consists of
\begin{itemize}
\item for each $a_1, \ldots, a_n \in C$ and $b\in D$, a set $X(a_1, \ldots,
  a_n; b)$ (Fig.~\ref{fig:module}(a))
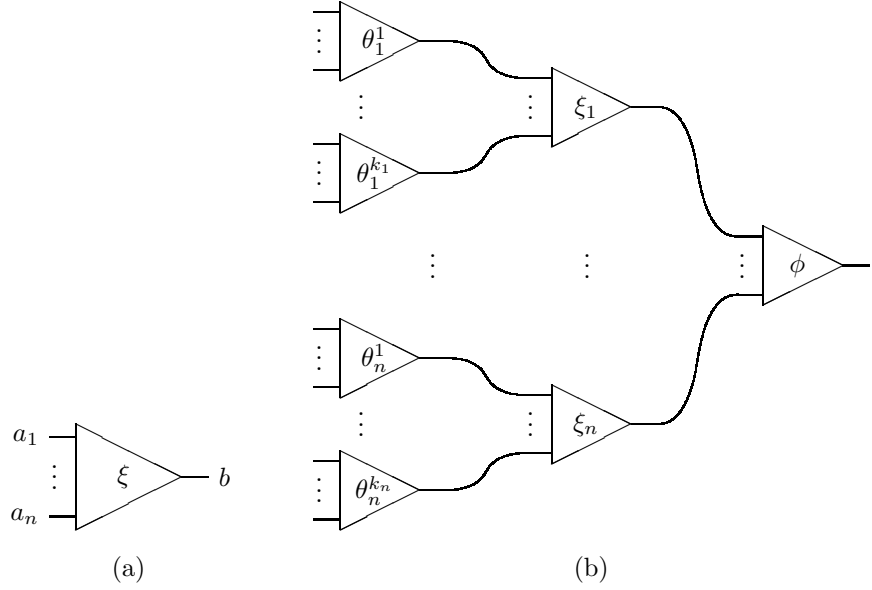
\begin{figure}%
\setlength{\unitlength}{1em}%
\centering
\begin{picture}(32,22)(0,-12)
\cell{0}{-10}{bl}{%
\begin{picture}(8,4)(0,-2)
\cell{2}{0}{r}{\tinputslft{a_1}{a_n}}
\cell{2}{0}{l}{\tusual{\xi}}
\cell{6}{0}{l}{\toutputrgt{b}}
\end{picture}}
\cell{4}{-12}{b}{\textrm{(a)}}
\cell{11}{-10}{bl}{%
\begin{picture}(21,20)(0,-10)
% leftmost input tags
\cell{1}{8.5}{r}{\tinputslftsmall{}{}}
\cell{1}{3.5}{r}{\tinputslftsmall{}{}}
\cell{1}{-3.5}{r}{\tinputslftsmall{}{}}
\cell{1}{-8.5}{r}{\tinputslftsmall{}{}}
% first column of transistors
\cell{1}{8.5}{l}{\tsmall{\theta_1^1}}
\cell{1}{3.5}{l}{\tsmall{\theta_1^{k_1}}}
\cell{1}{-3.5}{l}{\tsmall{\theta_n^1}}
\cell{1}{-8.5}{l}{\tsmall{\theta_n^{k_n}}}
% first column of ellipses
\cell{1.8}{6.3}{c}{\vdots}
\cell{4.5}{0.3}{c}{\vdots}
\cell{1.8}{-5.7}{c}{\vdots}
% output tags of first column of transistors
\cell{4}{8.5}{l}{\toutputrgt{}}
\cell{4}{3.5}{l}{\toutputrgt{}}
\cell{4}{-3.5}{l}{\toutputrgt{}}
\cell{4}{-8.5}{l}{\toutputrgt{}}
% lefthand joining curves
\qbezier(5,8.5)(6.125,8.5)(6.5,7.8)
\qbezier(8,7.1)(6.875,7.1)(6.5,7.8)
\qbezier(5,3.5)(6.125,3.5)(6.5,4.2)
\qbezier(8,4.9)(6.875,4.9)(6.5,4.2)
\qbezier(5,-8.5)(6.125,-8.5)(6.5,-7.8)
\qbezier(8,-7.1)(6.875,-7.1)(6.5,-7.8)
\qbezier(5,-3.5)(6.125,-3.5)(6.5,-4.2)
\qbezier(8,-4.9)(6.875,-4.9)(6.5,-4.2)
% input tags of second column of transistors
\cell{9}{6}{r}{\tinputslftsmall{}{}}
\cell{9}{-6}{r}{\tinputslftsmall{}{}}
% second column of transistors
\cell{9}{6}{l}{\tsmall{\xi_1}}
\cell{9}{-6}{l}{\tsmall{\xi_n}}
% ellipsis in second column of transistors
\cell{10.35}{0.3}{c}{\vdots}
% output tags of second column of transistors
\cell{12}{6}{l}{\toutputrgt{}}
\cell{12}{-6}{l}{\toutputrgt{}}
% righthand joining curves
\qbezier(13,6)(14.125,6)(14.5,3.55)
\qbezier(16,1.1)(14.875,1.1)(14.5,3.55)
\qbezier(13,-6)(14.125,-6)(14.5,-3.55)
\qbezier(16,-1.1)(14.875,-1.1)(14.5,-3.55)
% input tags of rightmost transistor
\cell{17}{0}{r}{\tinputslftsmall{}{}}
% rightmost transistor
\cell{17}{0}{l}{\tsmall{\phi}}
% rightmost output tag
\cell{20}{0}{l}{\toutputrgt{}}
\end{picture}}
\cell{21.5}{-12}{b}{\textrm{(b)}}
\end{picture}
% \hand{50}{65}
\caption{(a) `Element' of a module, (b)~compatibility of left and right
actions} 
\label{fig:module}  
\end{figure}
\item for each $a_i^j \in C$ and $b_i, b \in D$, a function
\begin{eqnarray*}
\begin{array}[b]{r}
D(b_1, \ldots, b_n; b) \times
X(a_1^1, \ldots, a_1^{k_1}; b_1) \times\cdots \\
\times 
X(a_n^1, \ldots, a_n^{k_n}; b_n)	
\end{array}
&
\go	&
X(a_1^1, \ldots, a_n^{k_n}; b),	\\
(\phi, \xi_1, \ldots, \xi_n)	&
\goesto	&
\phi \cdot (\xi_1, \ldots, \xi_n)
\end{eqnarray*}
\item for each $a_i^j, a_i \in C$ and $b\in D$, a function
\begin{eqnarray*}
\begin{array}[b]{r}
X(a_1, \ldots, a_n; b) \times
C(a_1^1, \ldots, a_1^{k_1}; a_1) \times\cdots \\
\times C(a_n^1, \ldots, a_n^{k_n}; a_n)	
\end{array}
&
\go	&
X(a_1^1, \ldots, a_n^{k_n}; b),	\\
(\xi, \theta_1, \ldots, \theta_n)	&
\goesto	&
\xi\cdot (\theta_1, \ldots, \theta_n),
\end{eqnarray*}
\end{itemize}
satisfying the evident axioms for compatibility of the two actions with
composition and identities in $D$ and $C$, together with a further axiom
stating compatibility with each other (Fig.~\ref{fig:module}(b)):
\[
(\phi \cdot (\xi_1, \ldots, \xi_n))
\cdot 
(\theta_1^1, \ldots, \theta_n^{k_n})
=
\phi \cdot
(\xi_1 \cdot (\theta_1^1, \ldots, \theta_1^{k_1}),
\ldots,
\xi_n \cdot (\theta_n^1, \ldots, \theta_n^{k_n}))
\]
whenever these expressions make sense.  
\end{defn}

When $C$ and $D$ have only unary arrows, this is the usual definition of
module between categories.  When $C$ and $D$ are operads, $X$ is a sequence
$(X(n))_{n\in\nat}$ of sets with left $D$-action and right $C$-action.
When $C=D$, taking $X(a_1, \ldots, a_n; b) = C(a_1, \ldots, a_n; b)$ gives
a canonical module $C \rMod C$.  For any $C$ and $D$, there is an obvious
notion of \demph{map between $(D,C)$-modules}, making $(D,C)$-modules into
a category.

Just as for rings, it is fruitful to consider one-sided modules.  Thus,
when $D$ is a multicategory, a \demph{left $D$-module} is a family
$(X(b))_{b\in D}$ of sets together with a left $D$-action---nothing other
than a $D$-algebra.  When $C$ is a multicategory, a \demph{right%
\lbl{p:cl-rt-modules}
$C$-module} is a family $(X(a_1, \ldots, a_n))_{a_1, \ldots, a_n \in C}$ of
sets together with a right $C$-action; these structures have sometimes been
considered in the special case of operads (see Voronov~\cite[\S 1]{VorSCO},
for instance).

\begin{notes}

The story of operads%
\index{operad!history}
and multicategories%
\index{multicategory!history}
is a typical one in mathematics,
strewn with failures of communication between specialists in different
areas.  Lazard%
\index{Lazard, Michel}
seems to have been the first person to have published the
basic idea, in work on formal group laws~\cite{Laz}; his `analyseurs'%
\index{analyseur}
are close to what were later dubbed operads.  Lambek~\cite{LamDSCII}
introduced multicategories, in the context of logic and linguistics.  He
says that B\'enabou%
\index{Benabou, Jean@B\'enabou, Jean}
and Cartier%
\index{Cartier, Pierre}
had both considered multicategories previously;
indeed, the idea might have occurred to anyone who knew what both a
category and a multilinear map were.  Perhaps a year or two later (but
effectively in a parallel universe), the special case of operads became
very important in homotopy theory, with the work of Boardman%
\index{Boardman, Michael}
and
Vogt~\cite{BV}%
\index{Vogt, Rainer}
(who approached them via PROPs) and May~\cite{MayGIL}%
\index{May, Peter}
(who
gave operads their name).  They in turn were building on the work of
Stasheff~\cite{StaHAHI},%
\index{Stasheff, Jim}
who defined the operad of associahedra without
isolating the operad concept explicitly.

As far as I can tell, all three parties (Lazard, Lambek, and the homotopy
theorists) were unaware of the work of the others until much later.  Users
of multicategories did not realize how interesting the one-object case was;
users of operads were slow to see how simple and natural was the
many-object case.  It seems to have been more than twenty years before the
appearance of the first paper containing both Lambek and Boardman--Vogt or
May in its bibliography: Beilinson%
\index{Beilinson, Alexander}
and Drinfeld~\cite{BeDr}.%
\index{Drinfeld, Vladimir}

Activity in operads and multicategories surged in the mid-1990s with, among
other things, the emergence of mirror symmetry and higher category theory.
See, for example, Loday, Stasheff and Voronov~\cite{Ren} and Baez and
Dolan~\cite{BDHDA3}.  Some important structures that I have not included
here are cyclic%
\index{operad!cyclic}
and modular%
\index{operad!modular}
operads (Getzler%
\index{Getzler, Ezra}
and
Kapranov~\cite{GeKaCOC, GeKaMO});%
\index{Kapranov, Mikhail}
cyclic and modular multicategories can be
defined similarly.

Some authors on operads make a special case of nullary%
\index{nullary!operation}
operations.
May~\cite[1.1]{MayGIL} insisted in the case of topological operads $P$ that
$P(0)$ should have only one element.  Markl, Shnider and
Stasheff~\cite{MSS} do not have a $P(0)$ at all; they start at $P(1)$.
They also consider pseudo-operads,%
\index{pseudo-operad}%
\index{operad!pseudo-}
which do not have an identity $1_P \in
P(1)$: see p.~\pageref{p:pseudo-operads} above.  

I thank Paolo Salvatore for illumination on the various operads for which
loop spaces are algebras.

Example~\ref{eg:multi-cheese} is probably a misnomer: the Swiss%
\index{multicategory!Swiss cheese}%
\index{operad!Swiss cheese}
Cheese
Board has recently been running an advertising campaign under the banner
`Si c'est trou\'e, c'est pas Suisse'.

\end{notes}

\chapter{Notions of Monoidal Category}
\lbl{ch:monoidal}

\chapterquote{%
In one sense all he ever wanted to be was someone with many nicknames}{%
Marcus~\cite{Marcus}}

\noindent
The concept of monoidal category is in such widespread use that one might
expect---or hope, at least---that its formalization would be thoroughly
understood.  Nevertheless, it is not.  Here we look at five different
possible definitions, plus one infinite family of definitions, of monoidal
category.  We prove equivalence results between almost all of them.

Apart from the wish to understand a common mathematical structure, there is
a reason for doing this motivated by higher-dimensional category theory.  A
monoidal category in the traditional sense is, as observed
in~\ref{sec:bicats}, 
% in~\ref{eg:bicat-mon-cat}, 
the same thing as a bicategory with only one object.  Similarly, any
proposed definition of weak $n$-category gives rise to a notion of monoidal
category, defined as a one-object weak 2-category.  So if we want to be
able to compare the (many) proposed definitions%
\index{n-category@$n$-category!definitions of!comparison}
of weak $n$-category then
we will certainly need a firm grip on the various notions of monoidal
category and how they are related.  (This is something like a physicist's
toy model: a manageably low-dimensional version of a higher-dimensional
system.)  Of course, if two definitions of weak $n$-category happen to
induce equivalent definitions of monoidal category then this does not imply
their equivalence in the general case, but the surprising variety of
different notions of monoidal category means that it is a surprisingly good
test.

In the classical definition of monoidal category, any pair $(X_1, X_2)$ of
objects has a specified tensor product $X_1 \otimes X_2$, and there is also
a specified unit object.  In~\ref{sec:unbiased} we consider a notion of
monoidal category in which any sequence $X_1, \ldots, X_n$ of objects has a
specified product; these are called `unbiased monoidal categories'.  More
generally, one might have a category equipped with various different tensor
products of various (finite) arities, but as long there are enough
isomorphisms between derived products, this should make essentially no
difference.  We formalize and prove this in~\ref{sec:alg-notions}.  

The definitions described so far are `algebraic' in that tensor is an
operation.  For instance, in the classical definition, two objects $X_1$
and $X_2$ give rise to an actual, specified, object, $X_1 \otimes X_2$, not
just an isomorphism class of objects.  This goes somewhat against
intuition: when using products of sets, for instance, it is not of the
slightest importance to remember the standard set-theoretic definition of
ordered pair ($(x_1, x_2) = \{ \{x_1\}, \{x_1, x_2\} \}$, as it happens);
all that matters is the universal property of the product.  So perhaps it
is better to use a notion of monoidal category in which the tensor product
of objects is only defined up to canonical isomorphism.  Three such
`non-algebraic' notions are considered in~\ref{sec:non-alg-notions}, one
using multicategories and the other two (closely related) using simplicial
objects.

Almost everything we do for monoidal categories could equally be done for
bicategories, as discussed in~\ref{sec:notions-bicat}.  The extension is
mostly routine---there are few new ideas---but the explanations are a
little easier in the special case of monoidal categories.  

Some generality is sacrificed in the theorems asserting equivalence of
different notions of monoidal category: these are equivalences of the
categories of monoidal categories, but we could have made them equivalences
of 2-categories (in other words, included monoidal transformations).  I
have not aimed to exhaust the subject, or the reader.

\section{Unbiased monoidal categories}
\lbl{sec:unbiased}

A classical monoidal category is a category $A$ equipped with a functor
$A^n \go A$ for each of $n=2$ and $n=0$ (the latter being the unit object
$I$), together with coherence data (Definition~\ref{defn:mon-cat}).  An
unbiased%
\index{monoidal category!unbiased}
monoidal category is the same, but with $n$ allowed to be any
natural number.  So in an unbiased monoidal category any four objects $a_1,
a_2, a_3, a_4$ have a specified tensor product
\[
(a_1 \otimes a_2 \otimes a_3 \otimes a_4),
\]
but in a classical monoidal category there are only derived products such
as
\[
(a_1 \otimes a_2) \otimes (a_3 \otimes a_4), 
\  
((a_1 \otimes a_2) \otimes a_3) \otimes a_4,
\  
a_1 \otimes ((I \otimes (a_2 \otimes I)) \otimes (a_3 \otimes a_4)).
\]
The classical definition is `biased' towards arities $2$ and $0$: it gives
them a special status.  Here we eliminate the bias.  

The coherence data for an unbiased monoidal category consists of
isomorphisms such as
\[
((a_1 \otimes a_2 \otimes a_3) \otimes (a_4 \otimes (a_5)) \otimes a_6)
\goiso
(a_1 \otimes (a_2 \otimes a_3 \otimes a_4 \otimes a_5) \otimes a_6).
\]
In fact, all such isomorphisms can be built up from two families of special
cases (the $\gamma$ and $\iota$ of~\ref{defn:lax-mon-cat} below).  We
require, of course, that the coherence isomorphisms satisfy all sensible
axioms.

We will show in the next section that the unbiased and classical
definitions are equivalent, in a strong sense.  The unbiased definition
seems much the more natural (having, for instance, no devious coherence
axioms), and much more useful for the purposes of theory.  When verifying
that a particular example of a category has a monoidal structure, it is
sometimes easier to use one definition, sometimes the other; but the
equivalence result means that we can take our pick.

Here is the definition of unbiased monoidal category.  It is no extra work
to define at the same time `lax monoidal categories', in which the
coherence maps need not be invertible. 

\begin{defn}	\lbl{defn:lax-mon-cat}
A \demph{lax%
\index{monoidal category!lax}
monoidal category} $A$ (or properly, $(A, \otimes,
\gamma, \iota)$) consists of 
\begin{itemize}
\item a category $A$
\item for each $n\in \nat$, a functor $\otimes_n: A^n \go A$,%
\glo{otimesn}
called \demph{$n$-fold%
\index{n-fold@$n$-fold!tensor}
tensor} and written
\[
(a_1, \ldots, a_n) \goesto (a_1 \otimes \cdots \otimes a_n)
\glo{otimesninfix}
\]
\item for each $n, k_1, \ldots, k_n \in \nat$ and double sequence $((a_1^1,
\ldots, a_1^{k_1}), \ldots, (a_n^1, \ldots, a_n^{k_n}))$ of objects of $A$,
a map
\[
\begin{array}{rl}
\gamma_{((a_1^1, \ldots, a_1^{k_1}), \ldots, (a_n^1, \ldots,
a_n^{k_n}))}:&
((a_1^1 \otimes \cdots \otimes a_1^{k_1}) \otimes \cdots \otimes
(a_n^1 \otimes \cdots \otimes a_n^{k_n}))\\
&\go
(a_1^1 \otimes\cdots\otimes a_1^{k_1} \otimes\cdots\otimes a_n^1
\otimes\cdots\otimes a_n^{k_n})
\end{array}
\glo{gammalaxmoncat}
\]
\item for each object $a$ of $A$, a map
\[
\iota_a: a \go (a),
\glo{iotalaxmoncat}
\]
\end{itemize}
with the following properties:
\begin{itemize}
\item
$\gamma_{((a_1^1, \ldots, a_1^{k_1}), \ldots, (a_n^1, \ldots,
a_n^{k_n}))}$ is natural in each of the $a_i^j$'s, and $\iota_a$ is natural
in $a$

\item
associativity: for any $n, m_p, k_p^q \in \nat$ and triple sequence 
$((( a_{p,q,r} )_{r=1}^{k_p^q} )_{q=1}^{m_p} )_{p=1}^{n}$ of objects, the
diagram
\[
\begin{diagram}[tight,width=4em,shortfall=2em,noPS] %,leftflush]
&&
\begin{scriptarray}
\scriptstyle 
  (
      ( 
          ( a_{1, 1, 1} \otimes\cdots\otimes a_{1, 1, k_1^1} )
          \otimes\cdots\otimes
          ( a_{1, m_1, 1} \otimes\cdots\otimes a_{1, m_1, k_1^{m_1}} )
      )
      \otimes\cdots
\\
\scriptstyle
\otimes
      ( 
          ( a_{n, 1, 1} \otimes\cdots\otimes a_{n, 1, k_n^1} )
          \otimes\cdots\otimes
          ( a_{n, m_n, 1} \otimes\cdots\otimes a_{n, m_n, k_n^{m_n}} ) 
      ) 
  )
\end{scriptarray}
& & \\
&\ldTo<{( \gamma_{D_1} \otimes\cdots\otimes \gamma_{D_n} )}
& &\rdTo>{\gamma_{D'}} & \\
\begin{scriptarray}
\scriptstyle
( 
    ( a_{1, 1, 1} \otimes\cdots\otimes a_{1, m_1, k_1^{m_1}} )
    \otimes\cdots
\\
\scriptstyle
\otimes
    ( a_{n, 1, 1} \otimes\cdots\otimes a_{n, m_n, k_n^{m_n}} )
)
\end{scriptarray}
& & & & 
\begin{scriptarray}
\scriptstyle
(
    ( a_{1, 1, 1} \otimes\cdots\otimes a_{1, 1, k_1^1} )
    \otimes\cdots
\\
\scriptstyle
\otimes
    ( a_{n, m_n, 1} \otimes\cdots\otimes a_{n, m_n, k_n^{m_n}} )
) 
\end{scriptarray}
\\
&\rdTo<{\gamma_D} & &\ldTo>{\gamma_{D''}} & \\
& &
\scriptstyle
(
   a_{1, 1, 1} \otimes\cdots\otimes a_{n, m_n, k_n^{m_n}} 
) 
& & \\
\end{diagram}
\]
commutes, where the double sequences $D_p, D, D', D''$ are
\begin{eqnarray*}
D_p	&=& 
\bftuple{\bftuple{a_{p, 1, 1}}{a_{p, 1, k_p^1}}}{\bftuple{a_{p, m_p,
1}}{a_{p, m_p, k_p^{m_p}}}},	\\
D	&=&	
\bftuple{\bftuple{a_{1, 1, 1}}{a_{1, m_1, k_1^{m_1}}}}{\bftuple{a_{n, 1,
1}}{a_{n, m_n, k_n^{m_n}}}},	\\
D'	&=&	
(
\bftuple{( a_{1, 1, 1} \otimes\cdots\otimes a_{1, 1, k_1^1} )}{(
a_{1, m_1, 1} \otimes\cdots\otimes a_{1, m_1, k_1^{m_1}} )}
, \ldots,	\\
&&
\bftuple{( a_{n, 1, 1} \otimes\cdots\otimes a_{n, 1, k_n^1} 
)}{( a_{n, m_n, 1} \otimes\cdots\otimes a_{n, m_n, k_n^{m_n}} )} ),	\\
D''	&=&
\bftuple{\bftuple{a_{1, 1, 1}}{a_{1, 1, k_1^1}}}{\bftuple{a_{n, m_n,
1}}{a_{n, m_n, k_n^{m_n}}}}
\end{eqnarray*}
\item
identity: for any $n\in\nat$ and sequence \bftuple{a_1}{a_n} of objects, the
diagrams
\begin{eqnarray*}
\begin{diagram}[scriptlabels,size=2em]
( a_1 \otimes\cdots\otimes a_n )	
&\rTo^{( \iota_{a_1} \otimes\cdots\otimes \iota_{a_n} )}
&( ( a_1 ) \otimes\cdots\otimes ( a_n ) )	\\
&\rdTo<{1}	
&\dTo>{\gamma_{\bftuple{(a_1)}{(a_n)}}}	\\
&& ( a_1 \otimes\cdots\otimes a_n )	\\
\end{diagram}
\\ 
\begin{diagram}[scriptlabels,size=2em]
( ( a_1 \otimes\cdots\otimes a_n ) )
&\lTo^{\iota_{( a_1 \otimes\cdots\otimes a_n )}}
&( a_1 \otimes\cdots\otimes a_n )	\\
\dTo<{\gamma_{(\bftuple{a_1}{a_n})}}
&\ldTo>{1}	& \\
( a_1 \otimes\cdots\otimes a_n ) & & \\
\end{diagram}
\end{eqnarray*}
commute.
\end{itemize}
An \demph{unbiased%
\index{monoidal category!unbiased}
monoidal category} is a lax monoidal category $(A,
\otimes, \gamma, \iota)$ in which each $\gamma_{((a_1^1, \ldots,
a_1^{k_1}), \ldots, (a_n^1, \ldots, a_n^{k_n}))}$ and each $\iota_a$ is
an isomorphism.  An \demph{unbiased strict monoidal category} is a lax
monoidal category $(A, \otimes, \gamma, \iota)$ in which each of the
$\gamma$'s and $\iota$'s is an identity map.
\end{defn}

\begin{remarks}{rmks:u-moncat}
\item
The bark of the associativity axiom is far worse than its bite.  All it
says is that any two ways of removing brackets are equal: for
instance, that the diagram
\begin{diagram}[width=4em,tight]%[leftflush=2.5em,noPS]
 & &
\begin{array}{l}
  \{[(a \otimes b) \otimes (c \otimes d)]
  \\
  \otimes [(e \otimes f) \otimes (g \otimes h)]\} 
\end{array}
& & \\ 
 &\ldTo<{\{\gamma \otimes \gamma\}} & &\rdTo>{\gamma} & \\ 
\begin{array}{l}
  \{[a \otimes b \otimes c \otimes d] 
  \\
  \otimes [e \otimes f \otimes g \otimes h]\} 
\end{array}
& & & &
\begin{array}{l}
  \{(a \otimes b) \otimes (c \otimes d)
  \\
  \otimes (e \otimes f) \otimes (g \otimes h)\} 
\end{array}
\\ 
 &\rdTo<{\gamma} & &\ldTo>{\gamma} & \\
 & &\{a \otimes b \otimes c \otimes d \otimes e \otimes f \otimes g \otimes
 h \} & & \\ 
\end{diagram}
commutes.  This is exactly the role of the associativity axiom for a monad
such as `free semigroup' on \Set, as observed in Mac Lane \cite[VI.4, after
Proposition 1]{MacCWM}.

\item	\lbl{rmk:axioms-obvious}
The coherence%
\index{coherence!axioms}
axioms for an unbiased monoidal category are `canonical' and
rather obvious, in contrast to those for classical monoidal categories;
they are the same shape as the diagrams expressing the associativity and
unit axioms for a monoid or monad~(\ref{defn:monad}). 

\item	\lbl{rmk:strict-unbiased}
In an unbiased \emph{strict} monoidal category, the coherence axioms
(naturality, associativity and identity) hold automatically.  Clearly,
unbiased strict monoidal categories are in one-to-one correspondence
with ordinary strict monoidal categories.

\end{remarks}

We have given a completely explicit definition of unbiased monoidal
category, but a more abstract version is possible.  First recall that if
\cat{C} is a strict 2-category then there is a notion of a \demph{strict
2-monad}%
\lbl{p:defn-2-monad}\index{two-monad@2-monad} 
$(T, \mu, \eta)$ on \cat{C}, and there are notions of \demph{strict},
\demph{weak} and \demph{lax algebra}%
\index{algebra!two-monad@for 2-monad}
for such a 2-monad.  (See Kelly and
Street~\cite{KSRE2}, for instance; terminology varies between authors.)  In
particular, `free strict monoidal category' is a strict 2-monad $(\blank^*,
\mu, \eta)$ on $\Cat$; the functor $\blank^*$ is given on objects of \Cat\
by
\[
A^* = \coprod_{n\in\nat} A^n.
\]
A (small) unbiased monoidal category is precisely a weak algebra for
this 2-monad.  Explicitly, this says that an unbiased monoidal
category consists of a category $A$ together with a functor $\otimes: A^*
\go A$ and natural isomorphisms
\begin{equation}
\label{eq:coh-transfs}
\begin{diagram}
A^{**}			&\rTo^{\mu_A}	&A^*		\\
\dTo<{\otimes^*}	&\nent \gamma	&\dTo>{\otimes}	\\
A^*			&\rTo_{\otimes}	&A		\\
\end{diagram}
\diagspace
\begin{diagram}[size=1.5em]
A	&		&\rTo^{\eta_A}	&		&A^*		\\
	&\rdTo(4,4)<{1}	&		&\nent \iota	&		\\
	&		&		&		&\dTo>{\otimes}	\\
	& 		&		&		&		\\
	&		&		&		&A		\\
\end{diagram}
\end{equation}
satisfying associativity and identity axioms: the diagrams
\[
\begin{diagram}[size=2em]
\otimes\of\otimes^{*}\of\otimes^{**}	&\rTo^{\gamma*1}&\otimes\of\mu_{A}\of\otimes^{**}
&\rEquals	&\otimes\of\otimes^{*}\of\mu_{A^*}	\\
\dTo<{1*\gamma^*}	&	&	&	&\dTo>{\gamma*1}\\
\otimes\of\otimes^{*}\of\mu^*_{A}	&\rTo_{\gamma *1}
&\otimes\of\mu_{A}\of\mu^*_{A}	&\rEquals
&\otimes\of\mu_{A}\of\mu_{A^*}\\
\end{diagram}
\]
\[
\begin{diagram}[size=2em]
\otimes\of 1_{A}^*		&\rTo^{1*\iota^*}	&
\otimes\of\otimes^{*}\of\eta_{A}^*			\\
				&\rdTo<1		&
\dTo>{\gamma*1}						\\
				&			&
\otimes\of\mu_{A}\of\eta^*_{A}				\\
\end{diagram}
\diagspace
\begin{diagram}
\otimes\of\otimes^{*}\of\eta_{A^*}	&\rEquals	&
\otimes\of\eta_{A}\of\otimes	&\lTo^{\iota*1}	&1_A \of\otimes	\\
\dTo<{\gamma*1}				&		&
				&\ldTo(4,2)>1	&		\\
\otimes\of\mu_{A}\of\eta_{A^*}		&		&
				&		&		\\
\end{diagram}
\]
commute.  This may easily be verified.  Similarly, a lax monoidal category
is precisely a lax algebra for the 2-monad (for `lax' means that the
natural transformations $\gamma$ and $\iota$ are no longer required to be
isomorphisms) and an unbiased strict monoidal category is precisely a
strict algebra (for `strict' means that $\gamma$ and $\iota$ are required
to be identities, so that the diagrams~\bref{eq:coh-transfs} containing
them commute).

A different abstract way of defining unbiased monoidal category will be
explored in~\ref{sec:alg-notions}. 

The next step is to define maps between (lax and) unbiased monoidal
categories.  Again, we could use the language of 2-monads to do this, but
opt instead for an explicit definition.

\begin{defn}	\lbl{defn:u-lax-mon-ftr} 
Let $A$ and $A'$ be lax monoidal categories.  Write $\otimes$ for the
tensor and $\gamma$ and $\iota$ for the coherence maps in both categories.
A \demph{lax monoidal%
\index{monoidal functor!unbiased}
functor} $(P, \pi): A \go A'$ consists of
\begin{itemize}
\item
a functor $P: A \go A'$
\item
for each $n\in\nat$ and sequence $a_1, \ldots, a_n$ of objects of $A$, a
map 
\[
\pi_{a_1, \ldots, a_n}:
(Pa_1 \otimes\cdots\otimes Pa_n)
\go
P(a_1 \otimes\cdots\otimes a_n),
\]
\end{itemize}
such that
\begin{itemize}
\item
$\pi_{a_1, \ldots, a_n}$ is natural in each $a_i$
\item
for each $n, k_i \in \nat$ and double sequence
\bftuple{\bftuple{a_1^1}{a_1^{k_1}}}{\bftuple{a_n^1}{a_n^{k_n}}} of objects
of $A$, the diagram
\[
\begin{diagram}[scriptlabels,width=2em]
\begin{scriptarray}
  \scriptstyle
  (
     ( Pa_1^1 \otimes\cdots\otimes Pa_1^{k_1} )
     \otimes\cdots
  \\
  \scriptstyle
     \otimes 
     ( Pa_n^1 \otimes\cdots\otimes Pa_n^{k_n} )
  )
\end{scriptarray}
&\rTo^{\gamma_{\bftuple{\bftuple{Pa_1^1}{Pa_1^{k_1}}}{\bftuple{Pa_n^1}{Pa_n^{k_n}}}}} 
&
\scriptstyle
( Pa_1^1 \otimes\cdots\otimes Pa_n^{k_n} )	\\
\dTo>{( \pi_{\bftuple{a_1^1}{a_1^{k_1}}} \otimes\cdots\otimes
        \pi_{\bftuple{a_n^1}{a_n^{k_n}}}  )}
& & \\
\begin{scriptarray}
  \scriptstyle
  (
      P ( a_1^1 \otimes\cdots\otimes a_1^{k_1} )
      \otimes\cdots
  \\
  \scriptstyle
  \otimes
      P ( a_n^1 \otimes\cdots\otimes a_n^{k_n} )
  )
\end{scriptarray}
&&\dTo<{\pi_{\bftuple{a_1^1}{a_n^{k_n}}}}	\\
\dTo>{\pi_{\bftuple{( a_1^1 \otimes\cdots\otimes a_1^{k_1} )}{( a_n^1
\otimes\cdots\otimes a_n^{k_n} )}}} 
& & \\
\begin{scriptarray}
  \scriptstyle
  P (
       ( a_1^1 \otimes\cdots\otimes a_1^{k_1} )
       \otimes\cdots
  \\
  \scriptstyle
  \otimes
       ( a_n^1 \otimes\cdots\otimes a_n^{k_n} )
    )
\end{scriptarray}
&\rTo_{P\gamma_{\bftuple{\bftuple{a_1^1}{a_1^{k_1}}}{\bftuple{a_n^1}{a_n^{k_n}}}}}
&
\scriptstyle
P ( a_1^1 \otimes\cdots\otimes a_n^{k_n} ) \\
\end{diagram}
\]
commutes
\item
for each 1-cell $a$, the diagram
\begin{diagram}[size=2em]
Pa	&\rTo^{\iota_{Pa}}	&( Pa )	\\
\dEquals&			&\dTo>{\pi_a}	\\
Pa	&\rTo_{P\iota_a}	&P( a )	\\
\end{diagram}
commutes.
\end{itemize}
A \demph{weak monoidal functor} is a
lax monoidal functor $(P, \pi)$ for which each of the maps
$\pi_{a_1, \ldots, a_n}$ is an isomorphism.  A \demph{strict monoidal
functor} is a lax monoidal functor $(P, \pi)$ for which each of the maps
$\pi_{a_1, \ldots, a_n}$ is an identity map (in which case $P$ preserves
composites and identities strictly).
\end{defn}

We remarked in~\ref{rmks:u-moncat}\bref{rmk:axioms-obvious} that the
coherence%
\index{coherence!axioms}
axioms for an unbiased monoidal category were rather obvious,
having the shape of the axioms for a monoid or monad.  Perhaps the
coherence axioms for an unbiased lax monoidal functor are a little less
obvious; they are, however, the same shape as the axioms for a lax map of
monads ($=$ monad functor, p.~\pageref{p:lax-map-of-monads}), and in any
case seem `canonical'.

Lax monoidal functors can be composed in the evident way, and the weak and
strict versions are closed under this composition.  There are also the
evident identities.  So we obtain $3\times 3 = 9$ possible categories,
choosing one of `strict', `weak' or `lax' for both the objects and the
maps.  In notation explained in a moment, the inclusions of subcategories
are as follows:
\[
\begin{diagram}[height=1.2em]
\fcat{LaxMonCat}_\mr{str}	&\sub	&\fcat{LaxMonCat}_\mr{wk}	&\sub
&\fcat{LaxMonCat}_\mr{lax} \\
\rotsub	&	&\rotsub	&	&\rotsub	\\	
\UMCstr	&\sub	&\UMCwk		&\sub	&\UMClax	\\
\rotsub	&	&\rotsub	&	&\rotsub	\\	
\fcat{StrMonCat}_\mr{str}	&\sub	&\fcat{StrMonCat}_\mr{wk}	&\sub
&\fcat{StrMonCat}_\mr{lax}. \\
\end{diagram}%
\glo{ninemoncat}%
\index{monoidal category!nine categories of}%
\]
For all three categories in the bottom (respectively, middle or top) row,
the objects are small strict (respectively, unbiased or lax) monoidal
categories.  For all three categories in the left-hand (respectively,
middle or right-hand) column, the maps are strict (respectively, weak or
lax) monoidal functors.  It is easy to check that the three categories in
the bottom row are isomorphic to the corresponding three categories in the
classical definition.

Of the nine categories, the three on the bottom-left to top-right diagonal
are the most conceptually natural: a level of strictness has been chosen
and stuck to.  In this chapter our focus is on the middle entry, $\UMCwk$,
where everything is weak.

Note that this $3\times 3$ picture does not appear in the classical,%
\index{monoidal category!unbiased vs. classical@unbiased \vs.\ classical}
`biased', approach to monoidal categories.  There the top row is obscured,
as there is no very satisfactory way to laxify%
\index{monoidal category!lax}
the classical definition of
monoidal category.  Admittedly it is possible to drop the condition that
the associativity maps $(a\otimes b) \otimes c \go a\otimes (b\otimes c)$
and unit maps $a\otimes 1 \go a \og 1\otimes a$ are isomorphisms (as
Borceux%
\index{Borceux, Francis}
does in his \cite{Borx1}, just after Definition 7.7.1), but somehow
this does not seem quite right.

To complete the picture, and to make possible the definition of equivalence
of unbiased monoidal categories, we define transformations.

\begin{defn}
Let $(P, \pi), (Q, \chi): A \go A'$ be lax monoidal functors between lax
monoidal categories.  A \demph{monoidal transformation}%
\index{monoidal transformation!unbiased}
$\sigma: (P, \pi)
\go (Q, \chi)$ is a natural transformation
\[
A \ctwomult{P}{Q}{\sigma} A'
\]
such that for all $a_1, \ldots, a_n \in A$, the diagram
\begin{diagram}[size=2em]
(Pa_1 \otimes \cdots \otimes Pa_n)	&
\rTo^{\pi_{a_1, \ldots, a_n}}		&
P(a_1 \otimes\cdots\otimes a_n)		\\
\dTo<{(\sigma_{a_1} \otimes\cdots\otimes \sigma_{a_n})}	&	&
\dTo>{\sigma_{(a_1 \otimes\cdots\otimes a_n)}}	\\
(Qa_1 \otimes \cdots \otimes Qa_n)	&
\rTo_{\chi_{a_1, \ldots, a_n}}		&
Q(a_1 \otimes\cdots\otimes a_n)		\\
\end{diagram}
commutes.
\end{defn}
(This time, there is only one possible level of strictness.)

Monoidal transformations can be composed in the expected ways, so that the
nine categories above become strict 2-categories.  In particular, $\UMCwk$
is a 2-category, so~(\ref{propn:bicat-eqv-eqv}) there is a notion of
equivalence of unbiased monoidal categories.  Explicitly, $A$ and $A'$ are
\demph{equivalent}%
\index{equivalence!unbiased monoidal categories@of unbiased monoidal categories} 
if there exist weak monoidal functors and invertible
monoidal transformations
\[
\begin{diagram}
A &\pile{\rTo^{(P, \pi)} \\ \lTo_{(Q, \chi)}} &A',
\end{diagram}
\diagspace
A \ctwomult{1}{(Q, \chi) \of (P, \pi)}{\eta} A,
\diagspace
A' \ctwomult{(P, \pi) \of (Q, \chi)}{1}{\epsln} A',
\]
and it makes no difference if we insist that $((P, \pi), (Q, \chi), \eta,
\epsln)$ forms an adjunction in $\UMCwk$.  As we might expect from the case
of classical monoidal categories~(\ref{propn:mon-eqv-eqv}), there is the
following alternative formulation.
\begin{propn}
Let $A$ and $A'$ be unbiased monoidal categories.  Then $A$ and $A'$ are
equivalent if and only if there exists a weak monoidal functor $(P, \pi): A
\go A'$ whose underlying functor $P$ is full, faithful and essentially
surjective on objects.
\end{propn}

\begin{proof} 
\latin{Mutatis mutandis}, this is the same as the proof
of~\ref{propn:mon-eqv-eqv}. 
\done
\end{proof}

We can now state and prove a coherence%
\index{coherence!monoidal categories@for monoidal categories!unbiased} 
theorem for unbiased
monoidal categories.
\begin{thm}	\lbl{thm:eqv-coh-umc}
Every unbiased monoidal category is equivalent to a strict monoidal category.
\end{thm}
\begin{proof}
Let $A$ be an unbiased monoidal category.  We construct an (unbiased)
strict monoidal category $\st(A)$,%
\glo{strictcover}
the \demph{strict%
\index{cover, strict}
cover} of $A$, and a
weak monoidal functor $(P, \pi): \st(A) \go A$ whose underlying functor $P$
is full, faithful and essentially surjective on objects.  By the last
proposition, this is enough.

An object of $\st(A)$ is a finite sequence $(a_1, \ldots, a_n)$ of objects
of $A$ (with $n\in\nat$).  A map $(a_1, \ldots, a_n) \go (b_1, \ldots,
b_m)$ in $\st(A)$ is a map $(a_1 \otimes\cdots\otimes a_n) \go (b_1
\otimes\cdots\otimes b_m)$ in $A$, and composition and identities in
$\st(A)$ are as in $A$.  The tensor in $\st(A)$ is given on objects by
concatenation: 
\[
((a_1^1, \ldots, a_1^{k_1}) \otimes\cdots\otimes (a_n^1, \ldots,
a_n^{k_n}))
=
(a_1^1, \ldots, a_1^{k_1}, \ldots, a_n^1, \ldots, a_n^{k_n}).
\]
To define the tensor of maps in $\st(A)$, take maps 
\begin{eqnarray*}
(a_1^1 \otimes\cdots\otimes a_1^{k_1})	&
\goby{f_1}				&
(b_1^1 \otimes\cdots\otimes b_1^{l_1})	\\
\vdots	&\vdots	&\vdots	\\
(a_n^1 \otimes\cdots\otimes a_n^{k_n})	&
\goby{f_n}				&
(b_n^1 \otimes\cdots\otimes b_n^{l_n})
\end{eqnarray*}
in $A$; then their tensor product in $\st(A)$ is the composite map
\begin{eqnarray*}
(a_1^1 \otimes\cdots\otimes a_n^{k_n})	&
\goby{\gamma^{-1}}			&
((a_1^1 \otimes\cdots\otimes a_1^{k_1}) \otimes\cdots\otimes (a_n^1
\otimes\cdots\otimes a_n^{k_n}))	\\
					&
\goby{(f_1 \otimes\cdots\otimes f_n)}	&
((b_1^1 \otimes\cdots\otimes b_1^{l_1}) \otimes\cdots\otimes (b_n^1
\otimes\cdots\otimes b_n^{l_n}))	\\
					&
\goby{\gamma}				&
(b_1^1 \otimes\cdots\otimes b_n^{l_n})
\end{eqnarray*}
in $A$.  It is absolutely straightforward and not too arduous to check that
$\st(A)$ with this tensor forms a strict monoidal category.  

The functor $P: \st(A) \go A$ is defined by
$
(a_1, \ldots, a_n) \goesto (a_1 \otimes\cdots\otimes a_n)
$
on objects, and `is the identity' on maps---in other words, performs the
identification
\[
\st(A)((a_1, \ldots, a_n), (b_1, \ldots, b_m))
=
A((a_1 \otimes\cdots\otimes a_n), (b_1 \otimes\cdots\otimes b_m)).
\]
For each double sequence 
\[
D = ((a_1^1, \ldots, a_1^{k_1}), \ldots, (a_n^1, \ldots, a_n^{k_n}))
\]
of objects of $A$, the isomorphism
\[
\pi_D: 
(P(a_1^1, \ldots, a_1^{k_1}) \otimes\cdots\otimes 
P(a_n^1, \ldots, a_n^{k_n}))
\goiso
P((a_1^1 \otimes\cdots\otimes a_1^{k_1}) \otimes\cdots\otimes (a_n^1
\otimes\cdots\otimes a_n^{k_n}))
\]
is simply $\gamma_D$.  It is also quick and straightforward to check that
$(P,\pi)$ is a weak monoidal functor.

$P$ is certainly full and faithful.  Moreover, for each $a\in A$ we have
an isomorphism
\[
\iota_a: a \goiso (a) = Pa,
\]
and this proves that $P$ is essentially surjective on objects.
\done
\end{proof}

As discussed on p.~\pageref{p:coherence-discussion}, coherence%
\index{coherence|(}
theorems
take various forms, usually falling under one of two headings: `all
diagrams commute' or `every weak thing is equivalent to a strict thing'.
The one here is of the latter type.  In~\ref{sec:alg-notions} we will prove
a coherence theorem for unbiased monoidal categories that is essentially of
the former type.

The proof above was adapted from Joyal%
\index{Joyal, Andr\'e}
and Street~\cite[p.~29]{JSBTC}.%
\index{Street, Ross!coherence for monoidal categories@on coherence for monoidal categories}
There the $\st$ construction was done for \emph{classical}%
\index{monoidal category!unbiased vs. classical@unbiased \vs.\ classical}
monoidal
categories $A$, for which the situation is totally different.  What happens
is that the $n$-fold tensor product used in the definition of both $\st(A)$
and $P$ must be replaced by some derived, non-canonical, $n$-fold tensor
product, such as
\[
(a_1, \ldots, a_n) 
\goesto 
a_1 \otimes (a_2 \otimes (a_3 \otimes \cdots \otimes (a_{n-1} \otimes a_n)
\cdots ));
\]
and then in order to define both $\pi$ and the tensor product of maps in
$\st(A)$, it is necessary to use coherence isomorphisms such as
\[
\begin{array}{rl}
&
(a_1 \otimes (a_2 \otimes (a_3 \otimes a_4)))
\ \otimes\ 
(a_5 \otimes (a_6 \otimes a_7))
\\
\goiso &
a_1 \otimes (a_2 \otimes (a_3 \otimes (a_4 \otimes (a_5 \otimes 
(a_6 \otimes a_7))))) .
\end{array}
\]
It would be folly to attempt to define these coherence isomorphisms, and
prove that $\st(A)$ and $(P, \pi)$ have the requisite properties, without
the aid of a coherence theorem for monoidal categories.  That is, the work
required to do this would be of about the same volume and kind as the work
involved in the syntactic proof of the `all diagrams commute' coherence
theorem, so one might as well have proved that coherence theorem anyway.
In contrast, the proof that every \emph{unbiased} monoidal category is
equivalent to a strict one is easy, short, and needs no supporting results.

This does not, however, provide a short cut to proving any kind of
coherence result for classical monoidal categories.  In the next section we
will see that unbiased and classical monoidal categories are essentially
the same, and it then follows from Theorem~\ref{thm:eqv-coh-umc} that every
classical monoidal category is equivalent to a strict one.  However, as
with any serious undertaking involving classical monoidal categories, the
proof that they are the same as unbiased ones is close to impossible
without the use of a coherence theorem.%
\index{coherence|)}

\section{Algebraic notions of monoidal category}
\lbl{sec:alg-notions}

We have already seen two notions of monoidal category: classical and
unbiased.  In one of them there was an $n$-fold tensor product just for
$n\in\{0,2\}$, and in the other there was an $n$-fold tensor product for
all $n\in\nat$.  But what happens if we take a notion of monoidal category
in which there is an $n$-fold tensor product for each $n$ lying in some
other subset of $\nat$?  More generally, what if we allow any number of
different $n$-fold tensors (zero, one, or more) for each value of
$n\in\nat$?  

For instance, we might choose to take a notion of monoidal category in
which there are 6 unit objects, a single 3-fold tensor product, 8 11-fold
tensor products, and $\aleph_4$ 38-fold tensor products.  Just as long as
we add in enough coherence isomorphisms to ensure that any two $n$-fold
tensor products built up from the given ones are canonically isomorphic,
this new notion of monoidal category ought to be essentially the same as
the classical notion.  

This turns out to be the case.  Formally, we start with a `signature'%
\index{signature}
$\Sigma \in \Set^\nat$.  (In the example this was given by $\Sigma(0)=6$,
$\Sigma(1)=\Sigma(2)=0$, $\Sigma(3)=1$, and so on.)  From this we define
the category $\Sigma\hyph\MCwk$ of `$\Sigma$-monoidal%
\index{Sigma-monoidal category@$\Sigma$-monoidal category}%
\index{monoidal category!Sigma-@$\Sigma$-}
categories' and weak
monoidal functors between them.  We then show that, up to equivalence,
$\Sigma\hyph\MCwk$ is independent of the choice of $\Sigma$, assuming only
that $\Sigma$ is large enough that we can build at least one $n$-fold
tensor product for each $n\in\nat$.  We also show that the category $\MCwk$
of classical monoidal categories is isomorphic to $\Sigma\hyph\MCwk$ for a
certain value of $\Sigma$, and that the same goes for the unbiased version
$\UMCwk$ (for a different value of $\Sigma$); it follows that the
classical%
\index{monoidal category!unbiased vs. classical@unbiased \vs.\ classical}
and unbiased definitions are equivalent.

In the introduction to this chapter, I argued that comparing definitions of
monoidal category is important for understanding higher-dimensional
category theory.  Here is a further reason why it is important, specific to
this particular, `algebraic', family of definitions of monoidal category.

Consider the definition of monoid.  Usually a monoid is defined as a set
equipped with a binary operation and a nullary operation, satisfying
associativity and unit equations.  `Weakening'%
\index{weakening!presentation-sensitivity of}
or `categorifying'%
\index{categorification!presentation-sensitivity of}
this
definition, we obtain the classical definition of (weak) monoidal category.
However, this process of categorification is dependent on presentation.%
\index{presentation, sensitivity of categorification to}
That is, we could equally well have defined a monoid as a set equipped with
one $n$-ary operation for each $n\in\nat$, satisfying appropriate
equations, and categorifying \emph{this} gives a different notion of
monoidal category---the unbiased one.  So different presentations of
the same 0-dimensional theory (monoids) give, under this process of
categorification, different 1-dimensional theories (of monoidal category).

Thus, our purpose is to show that in this particular situation, the
presentation-sensitivity of categorification disappears when we work up to
equivalence.

More generally, a fully-developed theory of weak $n$-categories might
include a formal process of weakening, which would take as input a theory
of strict structures and give as output a theory of weak structures.  If
the weakening process depended on how the theory of strict structures was
presented, then we would have to ask whether different presentations always
gave equivalent theories of weak structures.  We might hope so; but who
knows?

Our first task is to define $\Sigma$-monoidal categories and maps between
them, for an arbitrary $\Sigma \in \Set^\nat$.  Here it is in outline.  A
$\Sigma$-monoidal category should be a category $A$ equipped with a tensor
product $\otimes_\sigma: A^n \go A$ for each $\sigma\in\Sigma(n)$.  There
should also be a coherence isomorphism between any pair of derived $n$-fold
tensors: for instance, there should be a specified isomorphism
\[
\otimes_{\sigma_1} (\otimes_{\sigma_2} (a_1, a_2, a_3), a_4, 
\otimes_{\sigma_3} (a_5) ) 
\goiso
\otimes_{\sigma_4} (a_1, a_2, \otimes_{\sigma_5} (a_3, a_4), a_5 )
\]
for any
\[
\sigma_1, \sigma_2 \in \Sigma(3), \sigma_3 \in \Sigma(1),
\sigma_4 \in \Sigma(4), \sigma_5 \in \Sigma(2)
\]
and any objects $a_1, \ldots, a_5$ of $A$.  This isomorphism is naturally
depicted as
\begin{equation}	\label{eq:Sigma-mon-cat-iso}
% \drmk{Picture with one tree $\goiso$ another.  Labelled with $\sigma_i$'s
% and $a_i$'s.}
\begin{centredpic}
\begin{picture}(10,6.2)(-0.5,0)
% bottom layer
\put(4.5,0){\line(0,1){1.5}}
% middle layer
\cell{4.5}{1.5}{c}{\vx}
\put(4.5,1.5){\line(-2,1){3}}
\put(4.5,1.5){\line(0,1){1.5}}
\put(4.5,1.5){\line(2,1){3}}
% top layer
\cell{1.5}{3}{c}{\vx}
\put(1.5,3){\line(-1,1){1.5}}
\put(1.5,3){\line(0,1){1.5}}
\put(1.5,3){\line(1,1){1.5}}
\cell{7.5}{3}{c}{\vx}
\put(7.5,3){\line(0,1){1.5}}
% sigma_i labels
\cell{4.8}{1.5}{tl}{\sigma_1}
\cell{1.2}{3}{r}{\sigma_2}
\cell{7.8}{3}{l}{\sigma_3}
% a_i labels
\cell{0}{4.7}{b}{a_1}
\cell{1.5}{4.7}{b}{a_2}
\cell{3}{4.7}{b}{a_3}
\cell{4.5}{3.2}{b}{a_4}
\cell{7.5}{4.7}{b}{a_5}
\end{picture}
\end{centredpic}
\diagspace
\goiso
\diagspace
\begin{centredpic}
\begin{picture}(9.2,6.2)(0,0)
% bottom layer
\put(4.5,0){\line(0,1){1.5}}
% middle layer
\cell{4.5}{1.5}{c}{\vx}
\put(4.5,1.5){\line(-3,1){4.5}}
\put(4.5,1.5){\line(-1,1){1.5}}
\put(4.5,1.5){\line(1,1){1.5}}
\put(4.5,1.5){\line(3,1){4.5}}
% top layer
\cell{6}{3}{c}{\vx}
\put(6,3){\line(-1,1){1.5}}
\put(6,3){\line(1,1){1.5}}
% sigma_i labels
\cell{4.8}{1.5}{tl}{\sigma_4}
\cell{6.3}{3}{l}{\sigma_5}
% a_i labels
\cell{0}{3.2}{b}{a_1}
\cell{3}{3.2}{b}{a_2}
\cell{4.5}{4.7}{b}{a_3}
\cell{7.5}{4.7}{b}{a_4}
\cell{9}{3.2}{b}{a_5}
\end{picture}
\end{centredpic}.
\end{equation}
Recall from~\ref{sec:om-further} that labelled trees arise as operations in
free operads: the two derived tensor products drawn above are elements of
$(F\Sigma)(5)$, where $F\Sigma$ is the free operad on $\Sigma$.  Hence
$F\Sigma$ is the operad of (derived) tensor operations.  To obtain the
coherence isomorphisms, replace the set $(F\Sigma)(n)$ by the indiscrete
category $I((F\Sigma)(n))$ whose objects are the elements of
$(F\Sigma)(n)$; then the picture above shows a typical map in
$I((F\Sigma)(5))$.

Also recall from p.~\pageref{p:defn-V-Operad} that there is a notion of
`$\cat{V}$-operad' for any symmetric monoidal category \cat{V}, and an
accompanying notion of an algebra (in \cat{V}) for any \cat{V}-operad.  The
categories $I((F\Sigma)(n))$ form a \Cat-operad, an algebra for which
consists of a category $A$, a functor $A^n \go A$ for each derived $n$-fold
tensor operation, and a natural isomorphism between any two functors $A^n
\go A$ so arising, all fitting together coherently.
% (The coherence is encoded in the uniqueness of maps in an indiscrete
% category.)  
This is exactly what we want a $\Sigma$-monoidal category to
be.

Precisely, define the functor
\[
\begin{array}{rrcl}
\blank\hyph\MCwk: 	&\Set^\nat	&\go	&\CAT^\op	\\
			&\Sigma		&\goesto&\Sigma\hyph\MCwk
\end{array}
\]%
\glo{blankMonCatwk}%
to be the composite of the functors
\[
\Set^\nat \goby{F} 
\Set\hyph\Operad \goby{I_*}
\Cat\hyph\Operad \goby{\Alg_\mr{wk}}
\CAT^\op,
\]%
\glo{CatOperad}\glo{Algwk}%
where the terms involved are now defined in turn.
\begin{itemize}
\item $\Set\hyph\Operad$ is the category of ordinary (\Set-)operads,
usually just called $\Operad$~(\ref{defn:plain-opd}).
\item $\Cat\hyph\Operad$ is the category of \Cat-operads
(p.~\pageref{p:defn-V-Operad}).
\item $F$ is the free operad functor (p.~\pageref{p:free-mti-ftr}).
\item $I: \Set \go \Cat$ is the functor assigning to each set the
indiscrete%
\index{category!indiscrete}
category on it (p.~\pageref{p:indiscrete}), and since $I$
preserves products, it induces a functor $I_*: \Set\hyph\Operad \go
\Cat\hyph\Operad$.
\item For a \Cat-operad%
\index{Cat-operad@$\Cat$-operad}
$R$, the (large) category $\Alg_\mr{wk}(R)$
consists of $R$-algebras and \demph{weak
maps} between them.  So by
definition, an object of $\Alg_\mr{wk}(R)$ is a category $A$ together with
a sequence $(R(n) \times A^n \goby{\act{n}} A)_{n\in\nat}$ of functors,
compatible with the composition and identities of the operad $R$, and a map
$A \go A'$ in $\Alg_\mr{wk}(R)$ is a functor $P: A \go A'$ together with a
natural isomorphism
\begin{diagram}
R(n) \times A^n		&\rTo^{\act{n}}	&A	\\
\dTo<{1 \times P^n}	&\nent \pi_n	&\dTo>P	\\
R(n) \times A'^n	&\rTo_{\act{n}}	&A'	\\
\end{diagram}
for each $n\in\nat$, satisfying the coherence axioms in
Fig.~\ref{fig:wk-alg-coh}.  

\begin{figure}
\[
\begin{array}{c}
\begin{array}{rl}
&
\begin{diagram}[scriptlabels]
\begin{scriptarrayc}
  \scriptstyle
  R(n) \times \{ R(k_1) \times\cdots\times R(k_n) \} 
  \\
  \scriptstyle
  \times A^{k_1 + \cdots + k_n}	
\end{scriptarrayc}						&
\rTo^{1 \times \act{k_1} \times\cdots\times \act{k_n}}		&
\scriptstyle
R(n) \times A^n							&
\rTo^{\act{n}}							&
\scriptstyle
A								\\
\dTo<{1\times P^{k_1 + \cdots + k_n}}				&
\nent \scriptstyle
1\times \pi_{k_1} \times\cdots\times \pi_{k_n}		&
\dTo~{1 \times P^n}						&
\nent \scriptstyle
\pi_n							&
\dTo>P								\\
\begin{scriptarrayc}
  \scriptstyle
  R(n) \times \{ R(k_1) \times\cdots\times R(k_n) \} 
  \\
  \scriptstyle
  \times A'^{k_1 + \cdots + k_n}
\end{scriptarrayc}						&
\rTo_{1 \times \act{k_1} \times\cdots\times \act{k_n}}		&
\scriptstyle
R(n) \times A'^n							&
\rTo_{\act{n}}							&
\scriptstyle
A'								\\
\end{diagram}
\\
\\
\\
=
&
\begin{diagram}[scriptlabels]
\begin{scriptarrayc}
  \scriptstyle
  R(n) \times \{ R(k_1) \times\cdots\times R(k_n) \} 
  \\
  \scriptstyle
  \times A^{k_1 + \cdots + k_n}
\end{scriptarrayc}						&
\rTo^{\comp \times 1}						&
\scriptstyle
R(k_1 + \cdots + k_n) \times A^{k_1 + \cdots + k_n}			&
\rTo^{\act{k_1 + \cdots + k_n}}					&
\scriptstyle
A								\\
\dTo<{1 \times P^{k_1 + \cdots + k_n}}				&
\neeq								&
\dTo~{1 \times P^{k_1 + \cdots + k_n}}				&
\nent \scriptstyle
\pi_{k_1 + \cdots + k_n}					&
\dTo>P								\\
\begin{scriptarrayc}
  \scriptstyle
  R(n) \times \{ R(k_1) \times\cdots\times R(k_n) \} 
  \\
  \scriptstyle
  \times A'^{k_1 + \cdots + k_n}
\end{scriptarrayc}						&
\rTo_{\comp \times 1}						&
\scriptstyle
R(k_1 + \cdots + k_n) \times A'^{k_1 + \cdots + k_n}		&
\rTo_{\act{k_1 + \cdots + k_n}}					&
\scriptstyle
A'								\\
\end{diagram}
\end{array}
\\
\\
\\
\\
\begin{diagram}[scriptlabels]
\scriptstyle
1 \times A^1		&
\rTo^\diso	&
\scriptstyle
A	\\
\dTo<{1 \times P^1}	&
\neeq		&
\dTo>P	\\
\scriptstyle
1 \times A'^1		&
\rTo_\diso	&
\scriptstyle
A'	\\
\end{diagram}
\diagspace
=
\diagspace
\begin{diagram}[scriptlabels]
\scriptstyle
1 \times A^1		&
\rTo^{\ids \times 1}	&
\scriptstyle
R(1) \times A^1	&
\rTo^{\act{1}}		&
\scriptstyle
A			\\
\dTo<{1 \times P^1}	&
\neeq			&
\dTo~{1 \times P^1}	&
\nent \scriptstyle
\pi_1		&
\dTo>P			\\
\scriptstyle
1 \times A'^1		&
\rTo_{\ids \times 1}	&
\scriptstyle
R(1) \times A'^1	&
\rTo_{\act{1}}		&
\scriptstyle
A'			\\
\end{diagram}
\end{array}
\]
\caption{Coherence axioms for a weak map of $R$-algebras}
\label{fig:wk-alg-coh}
\end{figure}

\item For a map $H: R \go S$ of \Cat-operads, the induced functor $H^*:
\Alg_\mr{wk}(S) \go \Alg_\mr{wk}(R)$ is defined by composition with $H$: if
$A$ is an $S$-algebra then the resulting $R$-algebra has the same
underlying category and $R$-action given by
\[
R(n) \times A^n \goby{H_n \times 1} S(n) \times A^n \goby{\act{n}} A.
\]
\end{itemize}

In the construction above we used \emph{weak} maps between algebras for a
\Cat-operad, but we could just as well have used \demph{lax maps} (by
dropping the insistence that the natural transformations $\pi_n$ are
isomorphisms) or \demph{strict maps} (by insisting that the $\pi_n$'s are
identities).  The resulting functors $\Cat\hyph\Operad \go \CAT^\op$ are,
of course, called $\Alg_\mr{lax}$ and $\Alg_\mr{str}$,%
\glo{Algstr}
and so we have three functors
\[
\blank\hyph\MClax, \,
\blank\hyph\MCwk, \,
\blank\hyph\MCstr:
\Set^\nat \go \CAT^\op.
\glo{blankMonCatlax}
\]

\begin{defn}	\lbl{defn:Sigma-mon-cat}
Let $\Sigma\in\Set^\nat$.  A \demph{$\Sigma$-monoidal category}%
\index{Sigma-monoidal category@$\Sigma$-monoidal category}%
\index{monoidal category!Sigma-@$\Sigma$-}
is an
object of $\Sigma\hyph\MClax$ (or equivalently, of $\Sigma\hyph\MCwk$ or
$\Sigma\hyph\MCstr$).  A \demph{lax} (respectively, \demph{weak} or
\demph{strict}) \demph{monoidal functor}%
\index{monoidal functor!Sigma-@$\Sigma$-}
between $\Sigma$-monoidal
categories is a map in $\Sigma\hyph\MClax$ (respectively,
$\Sigma\hyph\MCwk$ or $\Sigma\hyph\MCstr$).
\end{defn}

Now we can state the results.  First, the notion of $\Sigma$-monoidal
category really does generalize the notions of unbiased and classical
monoidal category, as intended all along:

\begin{thm}%
\index{coherence!monoidal categories@for monoidal categories!unbiased} 
\lbl{thm:diag-coh-umc}
\thmname{Coherence for unbiased monoidal categories and functors}
Writing $1$ for the terminal object of $\Set^\nat$, there are
isomorphisms of categories
\begin{eqnarray*}
\UMClax		&\iso 	&1\hyph\MClax,	\\
\UMCwk		&\iso  	&1\hyph\MCwk , 	\\
\UMCstr		&\iso	&1\hyph\MCstr.   
\end{eqnarray*}
\end{thm}

\begin{proof}  
See Appendix~\ref{sec:app-UMC}.  It takes almost no calculation to see that
there is a canonical functor $1\hyph\MClax \go \UMClax$.  To see that it
is an isomorphism requires calculations using the coherence axioms for an
unbiased monoidal category.  Restricting to the weak and strict cases is
simple.  \done
\end{proof}

\begin{thm}
\lbl{thm:diag-coh-mc}
\thmname{Coherence%
\index{coherence!monoidal categories@for monoidal categories!classical} 
for classical monoidal categories and functors}
Write $\Sigma_\mr{c}$%
\glo{Sigmac}
for the object of $\Set^\nat$ given by
$\Sigma_\mr{c}(n)=1$ for $n\in\{0,2\}$ and $\Sigma_\mr{c}(n)=\emptyset$
otherwise.  Then there are isomorphisms of categories
\begin{eqnarray*}
\MClax		&\iso 	&\Sigma_\mr{c}\hyph\MClax,	\\
\MCwk		&\iso  	&\Sigma_\mr{c}\hyph\MCwk , 	\\
\MCstr		&\iso	&\Sigma_\mr{c}\hyph\MCstr.   
\end{eqnarray*}
\end{thm}

\begin{proof}  
See Appendix~\ref{sec:app-MC}.  The strategy is just the same as
in~\ref{thm:diag-coh-umc}, except that this time the calculations are in
principle much more tricky (because of the irregularity of the data and
axioms for a classical monoidal category) but in practice can be omitted
(by relying on the coherence theorems of others).  \done
\end{proof}

It took a little work to define `$\Sigma$-monoidal category', but I hope it
will be agreed that it is a completely natural definition, free from
\latin{ad hoc} coherence axioms and correctly embodying the idea of a
monoidal category with as many primitive tensor operations as are specified
by $\Sigma$.  So, for instance, the statement that the objects of $\MCwk$
correspond one-to-one with the objects of $\Sigma_\mr{c}\hyph\MCwk$ says
that the coherence data and axioms in the classical definition of monoidal
category are exactly right.  Were there no such isomorphism, it would be
the coherence data or axioms at fault, not the definition of
$\Sigma$-monoidal category.  This is why~\ref{thm:diag-coh-umc}
and~\ref{thm:diag-coh-mc} are called coherence theorems.

The unbiased coherence theorem~(\ref{thm:diag-coh-umc}) also tells us that
unbiased monoidal categories play a universal role: for any $\Sigma \in
\Set^\nat$, the unique map $\Sigma \go 1$ induces a canonical map from
unbiased monoidal categories to $\Sigma$-monoidal categories,
\[
\UMCwk \iso 1\hyph\MCwk \go \Sigma\hyph\MCwk.  
\]
Concretely, if we are given an unbiased monoidal category $A$ then we can
define a $\Sigma$-monoidal category by taking $\otimes_\sigma$ to be the
$n$-fold tensor $\otimes_n$, for each $\sigma\in \Sigma(n)$.

\begin{thm}	\lbl{thm:irrel-sig}%
\index{irrelevance of signature}\index{signature!irrelevance of}
\thmname{Irrelevance of signature for monoidal categories}
For any plausible $\Sigma, \Sigma' \in \Set^\nat$, there are equivalences
of categories
\[
\Sigma\hyph\MClax \eqv \Sigma'\hyph\MClax,
\diagspace
\Sigma\hyph\MCwk \eqv \Sigma'\hyph\MCwk.
\]
\end{thm}
Here $\Sigma \in \Set^\nat$ is called \demph{plausible}%
\index{plausible}
if $\Sigma(0) \neq
\emptyset$ and for some $n\geq 2$, $\Sigma(n) \neq \emptyset$.  This means
that $n$-fold tensor products can be derived for all $n\in\nat$.  In
contrast, if $\Sigma(0) = \emptyset$ then $(F\Sigma)(0) = \emptyset$, so
$(IF\Sigma)(0) = \emptyset$, so in a $\Sigma$-monoidal category all the
tensor operations are of arity $n \geq 1$---there is no unit object.
Dually, if $\Sigma(n) = \emptyset$ for all $n\geq 2$ then there is no
derived binary tensor.  In these cases we would not expect
$\Sigma$-monoidal categories to be much like ordinary monoidal categories.
So plausibility is an obvious minimal requirement.

\begin{proof}
It is enough to prove the result in the case $\Sigma'=1$.  As can be seen
from the explicit description of free operads on
p.~\pageref{p:free-plain-clauses}, plausibility of $\Sigma$ says exactly
that $(F\Sigma)(n)\neq \emptyset$ for each $n\in\nat$, or equivalently that
there exists a map $1\go UF\Sigma$ in $\Set^\nat$, where $U:
\Set\hyph\Operad \go \Set^\nat$ is the forgetful functor.  By adjointness,
this says that there is a map $F1 \go F\Sigma$ of \Set-operads, giving a
map $I_*F1 \go I_*F\Sigma$ of \Cat-operads.  On the other hand, $1$ is the
terminal object of $\Set^\nat$, so we have maps $I_*F1 \oppairu
I_*F\Sigma$.

In brief, the rest of the proof runs as follows.  $\Cat\hyph\Operad$
naturally has the structure of a 2-category, because \Cat\ does; and if two
objects of $\Cat\hyph\Operad$ are equivalent then so are their images under
both $\Alg_\mr{lax}$ and $\Alg_\mr{wk}$.  By the nature of indiscrete
categories, the existence of maps $I_*F1 \oppairu I_*F\Sigma$ implies that
$I_*F1 \eqv I_*F\Sigma$.  The result follows.  However, the 2-categorical
details are rather tiresome to check and the reader may prefer to avoid
them.  The main reason for including them below is that they reveal why the
theorem holds at the lax and weak levels but not at the strict level.

So, for the conscientious, a \demph{transformation}
\[
R \ctwomult{H}{H'}{\alpha} S
\]
of $\Cat$-operads is a sequence
\[
\left(
R(n) \ctwomult{H_n}{H'_n}{\alpha_n} S(n)
\right)_{n\in\nat}
\]
of natural transformations, such that
\[
\begin{array}{rl}
&
\scriptstyle
R(n) \times R(k_1) \times\cdots\times R(k_n)
\cone{\comp}
R(k_1 + \cdots + k_n)
\ctwomult{H_{k_1 + \cdots + k_n}}%
{H'_{k_1 + \cdots + k_n}}%
{}
S(k_1 + \cdots + k_n)				\\
=	&
\scriptstyle
R(n) \times R(k_1) \times\cdots\times R(k_n)
\ctwomult{H_n \times H_{k_1} \times\cdots H_{k_n}}%
{H'_n \times H'_{k_1} \times\cdots H'_{k_n}}%
{}
S(n) \times S(k_1) \times\cdots\times S(k_n)
\cone{\comp}
S(k_1 + \cdots + k_n),
\end{array}
\]
where the unlabelled 2-cells are respectively $\alpha_{k_1 + \cdots + k_n}$
and $\alpha_n \times \alpha_{k_1} \times\cdots\times \alpha_{k_n}$, and
\[
1 \cone{\ids} R(1) \ctwomult{H_1}{H'_1}{\alpha_1} S(1)
\ =\ 
1 \ctwomult{\ids}{\ids}{1} S(1).
\]
With the evident compositions, $\Cat\hyph\Operad$%
\glo{CatOperad2cat}
becomes a strict
2-category.  The statement on indiscrete categories is easily proved; in
fact, if $R$ is any $\Cat$-operad and $Q$ any $\Set$-operad then there is a
unique transformation between any pair of maps $R \parpair{}{} I_*Q$.

The functor $\Alg_\mr{lax}$ becomes a strict map $\Cat\hyph\Operad \go
\CAT^{\mr{co}\,\op}$ of 2-categories, where the codomain is $\CAT$ with
both the 1-cells and the 2-cells reversed~(\ref{sec:bicats}).  To see this,
let $\alpha$ be a transformation of \Cat-operads as above.  We need to
produce a natural transformation
\[
\Alg_\mr{lax}(R) 
\ctwomultcoop{H^*}{H'^*}{}
\Alg_\mr{lax}(S),
\]
which in turn means producing for each $S$-algebra $A$ a lax map $H'^*(A)
\go H^*(A)$ of $R$-algebras.  The composite natural transformation
\[
R(n) \times A^n
\ctwomult{H_n}{H'_n}{\alpha_n} 
S(n) \times A^n
\cone{\mi{act}_n^A}
A
\]
can be re-drawn as a natural transformation
\begin{eqnarray}
\label{diag:induced-lax-map}
\begin{diagram}
R(n) \times A^n 	&\rTo^{\mi{act}_n^{H'^*(A)}}	&A	\\
\dTo<1			&\nent 				&\dTo>1	\\
R(n) \times A^n 	&\rTo_{\mi{act}_n^{H^*(A)}}	&A,	\\
\end{diagram}
\end{eqnarray}
and this gives the desired lax map.  (To make the necessary distinctions,
superscripts have been added to the $\act{n}$'s naming the algebra
concerned.)  We thus obtain a strict map
\[
\Alg_\mr{lax}: \Cat\hyph\Operad \go \CAT^{\mr{co}\,\op}
\]
of 2-categories.

It follows immediately that if $R$ and $S$ are equivalent \Cat-operads then
$\Alg_\mr{lax}(R)$ and $\Alg_\mr{lax}(S)$ are equivalent categories.  To
see that $\Alg_\mr{wk}$ also preserves equivalence of objects, note that if
$\alpha$ is an invertible transformation of \Cat-operads then the natural
transformation~\bref{diag:induced-lax-map} is also invertible, and so
defines a weak map of algebras.  Put another way, $\Alg_\mr{wk}$ is a map
from the 2-category (\Cat-operads $+$ maps $+$ invertible transformations)
into $\CAT^{\mr{co}\,\op}$.

(However,~\bref{diag:induced-lax-map} only induces a \emph{strict} map
$H'^*(A) \go H^*(A)$ if it is the identity transformation, which means that
if $R$ and $S$ are equivalent \Cat-operads then $\Alg_\mr{str}(R)$ and
$\Alg_\mr{str}(S)$ are not necessarily equivalent categories.)  \done
\end{proof}

\begin{cor}%
\index{monoidal category!unbiased vs. classical@unbiased \vs.\ classical}
There are equivalences of categories
\[
\UMClax \eqv \MClax, 
\diagspace
\UMCwk \eqv \MCwk.
\]
\end{cor}

\begin{proof}
\ref{thm:diag-coh-umc} $+$ \ref{thm:diag-coh-mc} $+$ \ref{thm:irrel-sig}.
\done
\end{proof}

We arrived at the equivalence of unbiased and classical monoidal categories
by a roundabout route,
\[
\UMClax \iso 1\hyph\MClax \eqv \Sigma_\mr{c}\hyph\MClax \iso \MClax,
\]
so it may be useful to consider a direct proof.  Given an unbiased monoidal
category $(A, \otimes, \gamma, \iota)$, we can canonically write down a
classical monoidal category: the underlying category is $A$, the tensor is
$\otimes_2$, the unit is $\otimes_0$, and the associativity and unit
coherence isomorphisms are formed from certain components of $\gamma$ and
$\iota$.  The converse process is non-canonical:
\lbl{p:c-to-u}
we have to choose for each $n\in\nat$ an $n$-fold tensor operation built up
from binary tensor operations and the unit object.  Put another way, we
have to choose for each $n\in\nat$ an $n$-leafed classical tree (where
`classical' means that each vertex has either $2$ or $0$ outgoing edges, as
in~\ref{eg:opd-of-cl-trees}); and since $F\Sigma_\mr{c}$ is exactly the
operad of classical trees, this corresponds to the step of the proof
of~\ref{thm:irrel-sig} where we chose a map $1 \go UF\Sigma$ in
$\Set^\nat$.

All of the results above can be repeated with monoidal transformations
brought into the picture.  Then $\Sigma$-monoidal categories form a strict
2-category, and all the isomorphisms and equivalences of categories become
isomorphisms and equivalences of 2-categories.  That we did not \emph{need}
to mention monoidal transformations in order to prove the equivalence of
the various notions of monoidal category says something about the strength
of our equivalence result.  For suppose we start with a $\Sigma$-monoidal
category, derive from it a $\Sigma'$-monoidal category, and derive from
that a second $\Sigma$-monoidal category.  Then, in our construction, the
two $\Sigma$-monoidal categories are not just equivalent in the 2-category
$\Sigma\hyph\MCwk$, but isomorphic.  So Theorem~\ref{thm:irrel-sig} is `one
level%
\lbl{p:one-level-better}
better' than might be expected.

\section{Non-algebraic notions of monoidal category}
\lbl{sec:non-alg-notions}

Coherence%
\index{coherence!axioms}
axioms have got a bad name for themselves: unmemorable,
unenlightening, and unwieldy, they are often regarded as bureaucracy to be
fought through grimly before getting on to the real business.  Some people
would, therefore, like to create a world where there are no coherence
axioms at all---or anyway, as few as possible.

Whatever the merits of this aspiration, it is a fact that it can be
achieved in some measure; that is, there exist approaches to various higher
categorical structures that involve almost no coherence axioms.  In this
chapter we look at two different such approaches for monoidal categories.
The first exploits the relation between monoidal categories and
multicategories.  I will describe it in some detail.  The second is based
on the idea of the nerve of a category, and has its historical roots in the
homotopy-algebraic structures known as $\Gamma$-spaces.  Since it has less
to do with the main themes of this book, I will explain it more sketchily.

\index{multicategory!underlying|(}%
So, we start by looking at monoidal categories \vs.\
multicategories.  One might argue that multicategories are conceptually
more primitive than monoidal categories: that an operation taking several
inputs and producing one output is a more basic idea than a set whose
elements are ordered tuples.  In any case, the first example of a `tensor%
\index{tensor!ring@over ring}
product' that many of us learned was implicitly introduced via
multicategories---the tensor product $V\otimes W$ of two vector spaces
being characterized as the codomain of a `universal bilinear map' out of
$V, W$.  Now, the monoidal category of vector spaces contains no more or
less information than the multicategory of vector spaces; given either one,
the other can be derived in its entirety.  So, with the thoughts in our
head that multicategories are basic and coherence axioms are bad, we might
hope to `define' a monoidal category as a multicategory in which there are
enough universal maps around.  And this is what we do.

Formally, we have a functor $V$%
\glo{Vmonmti}
assigning to each monoidal category its
underlying multicategory; we want to show that $V$ gives an equivalence
between (monoidal categories) and some subcategory \cat{R} of $\Multicat$;
and we want, moreover, to describe \cat{R}.  In the terms of the previous
paragraph, \cat{R} consists of those multicategories containing `enough
universal maps'.

Before starting we have to make precise something that has so far been left
vague. In Example~\ref{eg:multi-mon}, we defined the underlying
multicategory $C$ of a monoidal category $A$ to have the same objects as
$A$ and maps given by
\[
C(a_1, \ldots, a_n ; a) = A(a_1 \otimes\cdots\otimes a_n, a),
\]
but we deferred the question of what exactly the expression $a_1
\otimes\cdots\otimes a_n$ meant.  We answer it now, and so obtain a precise
definition of the functor $V$.

If $A$ is a strict monoidal category then the expression makes perfect
sense.  If $A$ is not strict then it still makes perfect sense as long as
we have chosen to use unbiased rather%
\index{monoidal category!unbiased vs. classical@unbiased \vs.\ classical}
than classical monoidal categories:
take $a_1 \otimes\cdots\otimes a_n$ to mean $\otimes_n(a_1, \ldots, a_n)$.
Bringing into play the coherence maps of $A$, we can also define
composition in $C$, and so obtain the entire multicategory structure of $C$
without trouble; we arrive at a functor
\[
V: \UMCwk \go \Multicat.
\]

What if we insist on starting from a classical monoidal category?  We can
certainly obtain a multicategory by passing first from classical to
unbiased and then applying the functor $V$ just mentioned.  This passage
amounts to a choice of an $n$-leafed classical%
\index{tree!classical}
tree $\tau_n\in\ctr(n)$ for
each $n\in\nat$ (see p.~\pageref{p:c-to-u}), and the resulting
multicategory $C$ has the same objects as $A$ and maps given by
\[
C(a_1, \ldots, a_n; a) = A(\otimes_{\tau_n}(a_1, \ldots, a_n), a)
\]
where $\otimes_{\tau_n}: A^n \go A$ is `tensor according to the shape of
$\tau_n$' (defined formally in~\ref{sec:app-MC}).  So, any such sequence
$\tau_\bullet = (\tau_n)_{n\in\nat}$ induces a functor
\[
V_{\tau_\bullet}:
\MCwk \go \Multicat.
\]
For instance, if we take one of the two most obvious choices of sequence
$\tau_\bullet$ and write $C = V_{\tau_\bullet}(A)$ as usual, then
\[
C(a_1, \ldots, a_n; a) = A(
a_1 \otimes (a_2 \otimes (a_3 \otimes \cdots \otimes (a_{n-1} \otimes a_n)
\cdots )), a ).
\]
But there is also a way of passing from classical monoidal categories to
multicategories without making any arbitrary choices.  Let $A$ be a
classical monoidal category.  For $n\in\nat$ and $\tau, \tau' \in \ctr(n)$,
let 
\[
A^n \ctwomult{\otimes_\tau}{\otimes_{\tau'}}{\delta_{\tau,\tau'}} A
\]
be the canonical isomorphism whose existence is asserted by the coherence
theorem~(\ref{sec:mon-cats}).  Now define a multicategory $C$ by taking an
object to be, as usual, just an object of $A$, and a map $a_1, \ldots, a_n
\go a$ to be a family $(f_\tau)_{\tau\in\ctr(n)}$ in which $f_\tau \in
A(\otimes_\tau (a_1, \ldots, a_n), a)$ for each $\tau\in\ctr(n)$ and
$f_{\tau'} \of \delta_{\tau,\tau'} = f_\tau$ for all $\tau, \tau' \in
\ctr(n)$.  This yields another functor
\[
V': \MCwk \go \Multicat.
\]
However, a family $(f_\tau)_{\tau\in\ctr(n)}$ as above is entirely
determined by any single $f_\tau$, so the multicategory $C = V'(A)$ just
constructed is isomorphic to the multicategory $V_{\tau_\bullet}(A)$
obtained by choosing a particular sequence $\tau_\bullet$ of trees.  So $V'
\iso V_{\tau_\bullet}$ for any $\tau_\bullet$, and henceforth we write $V'$
or $V_{\tau_\bullet}$ as just $V$.

We now have ways of obtaining a multicategory from either an unbiased or a
classical monoidal category, and it is a triviality to check that these are
compatible with the equivalence between unbiased and classical: the diagram
\begin{diagram}[size=2em]
\UMCwk 	&	&\rTo^{\eqv}	&	&\MCwk	\\
	&\rdTo<V&		&\ldTo>V&	\\
	&	&\Multicat	&	&	\\
\end{diagram}
commutes up to canonical isomorphism.  So we know unambiguously
what it means for a multicategory to be `the underlying multicategory of
some monoidal category', and similarly for maps.  

\begin{figure}
\[
% DOMAIN
%
\begin{centredpic}
\begin{picture}(16,12)(0,-6)
% transistors
\cell{10}{0}{l}{\tusualdotty{\exists !\ovln{f}}}
\cell{2}{4}{l}{\tusual{}}
\cell{2}{-4}{l}{\tusual{}}
% labels on left-hand transistors
\cell{3.8}{4}{bl}{\tnelabel{u(c_1^1, \ldots, c_1^{k_1})}}
\cell{3.8}{-4}{tl}{\tselabel{u(c_n^1, \ldots, c_n^{k_n})}}
% short wires
\cell{14}{0}{l}{\toutputrgt{c}}
\cell{2}{4}{r}{\tinputslft{c_1^1}{c_1^{k_1}}}
\cell{2}{-4}{r}{\tinputslft{c_n^1}{c_n^{k_n}}}
% long wires
\qbezier(10,1.5)(8.5,1.5)(8,2.75)
\qbezier(6,4)(7.5,4)(8,2.75)
\cell{8}{3.2}{l}{\otimes(c_1^1, \ldots, c_1^{k_1})}
\qbezier(10,-1.5)(8.5,-1.5)(8,-2.75)
\qbezier(6,-4)(7.5,-4)(8,-2.75)
\cell{8}{-3.2}{l}{\otimes(c_n^1, \ldots, c_n^{k_n})}
% ellipses
\cell{9.2}{0.3}{c}{\vdots}
\cell{3}{0}{c}{\cdot}
\cell{3}{1}{c}{\cdot}
\cell{3}{-1}{c}{\cdot}
\end{picture}
\end{centredpic}
\mbox{\hspace{2em}}
=
\mbox{\hspace{2em}}
%
% CODOMAIN
%
\begin{centredpic}
\begin{picture}(8,12)(0,-6)
% transistor
\put(2,-6){\line(0,1){12}}
\put(6,0){\line(-2,-3){4}}
\put(6,0){\line(-2,3){4}}
\cell{3.5}{0}{c}{f}
% tags
\cell{2}{4}{r}{\tinputslft{c_1^1}{c_1^{k_1}}}
\cell{2}{-4}{r}{\tinputslft{c_n^1}{c_n^{k_n}}}
\cell{1.2}{0.3}{c}{\vdots}
\cell{6}{0}{l}{\toutputrgt{c}}
\end{picture}
\end{centredpic}
\]
% \hand{50}{1}
\caption{Representation of a multicategory}
\label{fig:repn}
\end{figure}
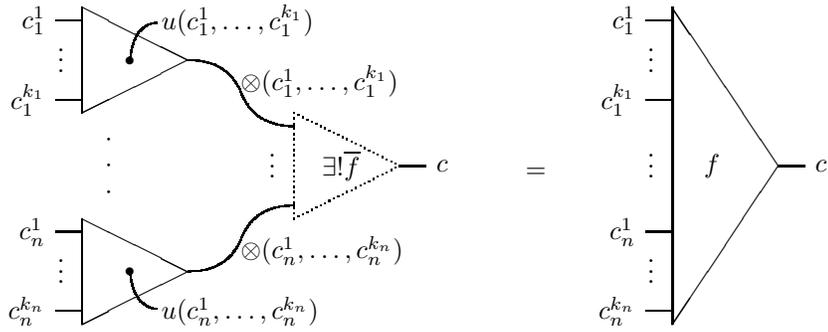

The results for which we hoped, exhibiting monoidal categories as special
multicategories, can be phrased in various different ways.
\begin{defn}	\lbl{defn:repn-multi}
A \demph{representation}%
\index{multicategory!representation of}
of a multicategory $C$ consists of an object
$\otimes(c_1, \ldots, c_n)$ and a map 
\[
u(c_1, \ldots, c_n):
c_1, \ldots, c_n 
\go
% \goby{u(c_1, \ldots, c_n)} 
\otimes(c_1, \ldots, c_n)
\]
for each $n\in\nat$ and $c_1, \ldots, c_n \in C$, with the following
factorization property (Fig.~\ref{fig:repn}): for any objects $c_1^1,
\ldots, c_1^{k_1}, \ldots, c_n^1, \ldots, c_n^{k_n}, c$ and any map $f:
c_1^1, \ldots, c_n^{k_n} \go c$, there is a unique map
\[
\ovln{f}: 
\otimes(c_1^1, \ldots, c_1^{k_1}), \ldots, \otimes(c_n^1, \ldots,
c_n^{k_n})
\go c
\]
such that 
\[
\ovln{f} \of 
(u(c_1^1, \ldots, c_1^{k_1}), \ldots, u(c_n^1, \ldots, c_n^{k_n}))
= f.
\]
A multicategory is \demph{representable}%
\index{multicategory!representable}
if it admits a
representation.
\end{defn}

\begin{defn}	\lbl{defn:pre-univ-simple}
A map $c_1, \ldots, c_n \goby{u} c'$ in a multicategory is
\demph{pre-universal}%
\index{pre-universal!map in multicategory}
if (Fig.~\ref{fig:univ-simple}(a)) for any object
$c$ and map $c_1, \ldots, c_n \goby{f} c$, there is a unique map $c'
\goby{\ovln{f}} c$ such that $\ovln{f} \of u = f$.
\end{defn}
\begin{figure}
\centering
\setlength{\unitlength}{1em}
\begin{picture}(26.5,22)(0,0)
\cell{0}{20}{l}{%
\begin{picture}(14,4)(0,-2)
\cell{2}{0}{r}{\tinputslft{c_1}{c_n}}
\cell{2}{0}{l}{\tusual{u}}
\put(6,0){\line(1,0){2}}
\cell{7}{0.2}{b}{c'}
\cell{8}{0}{l}{\tusualdotty{\exists !\ovln{f}}}
\cell{12}{0}{l}{\toutputrgt{c}}
\end{picture}}
\cell{16}{20}{c}{=}
\cell{18}{20}{l}{%
\begin{picture}(8.5,4)(-0.5,-2)
\cell{2}{0}{r}{\tinputslft{c_1}{c_n}}
\cell{2}{0}{l}{\tusual{f}}
\cell{6}{0}{l}{\toutputrgt{c}}
\end{picture}}
\cell{13.25}{16}{b}{\textrm{(a)}}
\cell{0}{8}{l}{%
\begin{picture}(14,12)(0,-6)
% leftmost tags
\cell{2}{4}{r}{\tinputslft{a_1}{a_p}}
\cell{2}{0}{r}{\tinputslft{c_1}{c_n}}
\cell{2}{-4}{r}{\tinputslft{b_1}{b_q}}
% lefthand transistor
\cell{2}{0}{l}{\tusual{u}}
% joining wires
\qbezier(2,5.5)(3.875,5.5)(4.5,4.75)
\qbezier(7,4)(5.125,4)(4.5,4.75)
\qbezier(2,2.5)(3.875,2.5)(4.5,1.75)
\qbezier(7,1)(5.125,1)(4.5,1.75)
\put(6,0){\line(1,0){2}}
\cell{6.8}{0.1}{b}{c'}
\qbezier(2,-5.5)(3.875,-5.5)(4.5,-4.75)
\qbezier(7,-4)(5.125,-4)(4.5,-4.75)
\qbezier(2,-2.5)(3.875,-2.5)(4.5,-1.75)
\qbezier(7,-1)(5.125,-1)(4.5,-1.75)
% inputs to righthand transistor
\cell{8}{2.5}{r}{\tinputslft{}{}}
\cell{8}{-2.5}{r}{\tinputslft{}{}}
% righthand transistor
\qbezier[45](8,4.5)(8,0)(8,-4.5)
\qbezier[32](8,4.5)(10,2.25)(12,0)
\qbezier[32](8,-4.5)(10,-2.25)(12,0)
\cell{9.7}{0}{c}{\exists !\ovln{f}}
% output tag
\cell{12}{0}{l}{\toutputrgt{c}}
\end{picture}}
\cell{16}{8}{c}{=}
\cell{18}{8}{l}{%
\begin{picture}(8.5,12)(-0.5,-6)
% leftmost tags
\cell{2}{4}{r}{\tinputslft{a_1}{a_p}}
\cell{2}{0}{r}{\tinputslft{c_1}{c_n}}
\cell{2}{-4}{r}{\tinputslft{b_1}{b_q}}
% transistor
\put(2,6){\line(0,-1){12}}
\put(2,6){\line(2,-3){4}}
\put(2,-6){\line(2,3){4}}
\cell{3.7}{0}{c}{f}
% output tag
\cell{6}{0}{l}{\toutputrgt{c}}
\end{picture}}
\cell{13.25}{0}{b}{\textrm{(b)}}
\end{picture}
% \hand{90}{2}
\caption{(a) Pre-universal map, and (b) universal map}
\label{fig:univ-simple}
\end{figure}
\begin{defn}	\lbl{defn:univ-simple}
A map $c_1, \ldots, c_n \goby{u} c'$ in a multicategory is 
\demph{universal}%
\index{universal!map in multicategory}
if (Fig.~\ref{fig:univ-simple}(b)) for any objects
$a_1, \ldots, a_p, b_1, \ldots, b_q, c$ and any map 
\[
a_1, \ldots, a_p, c_1, \ldots, c_n, b_1, \ldots, b_q \goby{f} c,
\]
there is a unique map $a_1, \ldots, a_p, c', b_1, \ldots, b_q
\goby{\ovln{f}} c$ such that $\ovln{f} \of_{p+1} u = f$.
(See~\ref{sec:om-further} for the $\of_{p+1}$ notation.)
\end{defn}

Here is the main result.
\begin{thm}	\lbl{thm:rep-multi}
\begin{enumerate}
\item	\lbl{part:rep-multi-objs}
The following conditions on a multicategory $C$ are equivalent:
\begin{itemize}
\item $C \iso V(A)$ for some monoidal category $A$
\item $C$ is representable
\item every sequence $c_1, \ldots, c_n$ of objects of $C$ is the domain of
some pre-universal map, and the composite of pre-universal maps is
pre-universal 
\item every sequence $c_1, \ldots, c_n$ of objects of $C$ is the domain of
some universal map.
\end{itemize}
Under these equivalent conditions, a map in $C$ is universal if and only if
it is pre-universal.
\item	\lbl{part:rep-multi-maps}
Let $A$ and $A'$ be monoidal categories.  The following conditions on
a map $H: V(A) \go V(A')$ of multicategories are equivalent:
\begin{itemize}
\item $H = V(P, \pi)$ for some weak monoidal functor $(P, \pi): A \go A'$
\item $H$ preserves%
\index{universal!preservation}
universal maps (that is, if $u$ is a universal
map in $V(A)$ then $Hu$ is a universal map in $V(A')$).
\end{itemize}
\end{enumerate}
\end{thm}
The functor $V$ is faithful, and therefore provides an equivalence between
$\UMCwk$ or $\MCwk$ (as you prefer) and the subcategory $\fcat{RepMulti}$
of $\Multicat$ consisting of the representable multicategories and the
universal-preserving maps.

The only subtle point here is that the existence of a pre-universal map for
every given domain is \emph{not} enough to ensure that the multicategory
comes from a monoidal category.  A specific example appears in
Leinster~\cite{FM}, but the point can be explained here in the familiar
context of vector spaces.  Suppose we are aware that for each pair $(X, Y)$
of vector spaces, there is an object $X\otimes Y$ and a bilinear map $X, Y
\goby{u_{X,Y}} X\otimes Y$ with the traditional universal property (which
we are calling pre-universality).  Then it does \emph{not} follow for
purely formal reasons that the tensor product is associative (up to
isomorphism): one has to use some actual properties of vector spaces.
Essentially, one either has to show that trilinear maps of the form
$u_{X\otimes Y, Z} \of (u_{X,Y}, 1_Z)$ and $u_{X, Y\otimes Z} \of (1_X,
u_{Y,Z})$ are pre-universal, or show that maps of the form $u_{X,Y}$ are in
fact universal.

The energetic reader with plenty of time on her hands will have no
difficulty in proving Theorem~\ref{thm:rep-multi}; the main ideas have been
explained and it is just a matter of settling the details.  So in a sense
that is an end to the matter: monoidal categories can be recognized as
multicategories with a certain property, and monoidal functors similarly,
all as hoped for originally.

We can, however, take things further.  With just a little more work than a
direct proof would involve, Theorem~\ref{thm:rep-multi} can be seen as a
special case of a result in the theory of fibrations%
\index{fibration!multicategories@of multicategories}
of multicategories.
This theory is a fairly predictable extension of the theory of fibrations
of ordinary categories, and the result of which~\ref{thm:rep-multi} is a
special case is the multicategorical analogue of a standard result on
categorical fibrations.

The basic theory of fibrations of multicategories is laid out in
Leinster~\cite{FM}, which culminates in the deduction of~\ref{thm:rep-multi}
and some related facts on, for instance, \emph{strict} monoidal categories
as multicategories.  Here is the short story.  For any category $D$, a
fibration (or really, opfibration)
% $C \go D$ 
over $D$ is essentially the same thing as a weak functor $D \go \Cat$.  (We
looked at the case of \emph{discrete} fibrations in~\ref{sec:cats}.)  With
appropriate definitions, a similar statement can be made for
multicategories.  Taking $D$ to be the terminal multicategory $1$, we find
that the unique map $C \go 1$ is a fibration exactly when $C$ is
representable, and that weak functors $1 \go \Cat$ are exactly unbiased
monoidal categories.  (Universal and pre-universal maps in $C$ correspond
to what are usually called cartesian and pre-cartesian maps.)  So a
representable multicategory is essentially the same thing as a monoidal
category.%
\index{multicategory!underlying|)}

\paragraph*{}

We now consider a different non-algebraic notion of monoidal category:
`homotopy monoidal categories'.  The idea can be explained as follows.

Recall that every small category has a nerve, and that this allows
categories to be described as simplicial sets satisfying certain
conditions.  Explicitly, if $n\in\nat$ then let $\upr{n} = \{ 0, 1, \ldots,
n \}$,%
\glo{upr}
and let $\Delta$%
\glo{Delta}\index{simplex category $\Delta$}
be the category whose objects are $\upr{0},
\upr{1}, \upr{2}, \ldots$ and whose morphisms are all order-preserving
functions; so $\Delta$ is equivalent to the category of nonempty finite
totally ordered sets.  A functor $\Delta^\op \go \Set$ is called a
\demph{simplicial%
\index{simplicial set}
set}.  (More generally, a functor $\Delta^\op \go \Eee$
is called a \demph{simplicial%
\index{simplicial object}
object in $\Eee$}.)  Any ordered set $(I,
\leq)$ can be regarded as a category with object-set $I$ and with exactly
one morphism $i \go j$ if $i\leq j$, and none otherwise.  This applies in
particular to the ordered sets \upr{n}, and so we may define the
\demph{nerve}% 
\lbl{p:defn-nerve}\index{nerve!category@of category} 
$NA$ of a small category $A$ as the simplicial set
\[
\begin{array}{rrcl}
NA: 	&\Delta^\op	&\go 		&\Set			\\
	&\upr{n}	&\goesto	&\Cat(\upr{n}, A).
\end{array}
\]
This gives a functor $N: \Cat \go \ftrcat{\Delta^\op}{\Set}$, which turns
out to be full and faithful.  Hence $\Cat$ is equivalent to the full
subcategory of $\ftrcat{\Delta^\op}{\Set}$ whose objects are those
simplicial sets $X$ isomorphic to $NA$ for some small category $A$.  There
are various intrinsic characterizations of such simplicial sets $X$.  We do
not need to think about the general case for now, only the special case of
one-object categories, that is, monoids.

So: let $k, n_1, \ldots, n_k \in \nat$, and for each $j\in \{ 1, \ldots, k
\}$, define a map $\iota_j$ in $\Delta$ by 
\[
\begin{array}{rrcl}
\iota_j: 	&\upr{n_j}	&\go		&\upr{n_1 + \cdots + n_k}\\
		&p		&\goesto	&n_1 + \cdots + n_{j-1} + p.
\end{array}
\]
Given also a simplicial set $X$, let%
\index{Segal, Graeme!map}
\[
\xi_{n_1, \ldots, n_k}: 
X\upr{n_1 + \cdots + n_k} 
\go 
X\upr{n_1} \times\cdots\times X\upr{n_k}
\]
be the map whose $j$th component is $X(\iota_j)$.  Write
\[
\xi^{(k)} = \xi_{1, \ldots, 1} : X\upr{k} \go X\upr{1}^k.
\]
\begin{propn}	\lbl{propn:simp-isos}
The following conditions on a simplicial set $X$ are equivalent:
\begin{enumerate}
\item	\lbl{item:simp-iso-general}
$\xi_{n_1, \ldots, n_k}: 
X\upr{n_1 + \cdots + n_k} 
\go 
X\upr{n_1} \times\cdots\times X\upr{n_k}$
is an isomorphism for all $k, n_1, \ldots, n_k \in \nat$
\item 	\lbl{item:simp-iso-classical}
$\xi_{m,n}: X\upr{m+n} \go X\upr{m} \times X\upr{n}$
is an isomorphism for all $m, n \in \nat$, and the unique map $X\upr{0} \go
1$ is an isomorphism
\item 	\lbl{item:simp-iso-powers}
$\xi^{(k)}: X\upr{k} \go X\upr{1}^k$
is an isomorphism for all $k\in\nat$
\item 	\lbl{item:simp-iso-nerve}
$X \iso NA$ for some monoid $A$.
\end{enumerate}
\end{propn}

\begin{proof} 
Straightforward.  Note that~\bref{item:simp-iso-classical} is
just~\bref{item:simp-iso-general} restricted to $k\in\{0,2\}$, and
similarly that~\bref{item:simp-iso-powers} is
just~\bref{item:simp-iso-general} in the case $n_1 = \cdots = n_k = 1$.
\done
\end{proof}

Monoids are, therefore, the same thing as simplicial sets satisfying any of
the conditions \bref{item:simp-iso-general}--\bref{item:simp-iso-powers}.
The proposition can be generalized: replace $\Set$ by any category $\Eee$
possessing finite products to give a description of monoids in $\Eee$ as
certain simplicial objects in $\Eee$.  In particular, if we take
$\Eee=\Cat$ then we obtain a description of strict monoidal categories as
certain simplicial objects in $\Cat$.  This suggests that a `righteous'
(weak)%
\index{weakening}
notion of monoidal category could be obtained by changing the
isomorphisms in Proposition~\ref{propn:simp-isos} to equivalences.

\begin{propn}	\lbl{propn:simp-eqs}
The following conditions on a functor $X: \Delta^\op \go \Cat$ are
equivalent: 
\begin{enumerate}
\item	\lbl{item:simp-eq-general}
$\xi_{n_1, \ldots, n_k}: 
X\upr{n_1 + \cdots + n_k} 
\go 
X\upr{n_1} \times\cdots\times X\upr{n_k}$
is an equivalence for all $k, n_1, \ldots, n_k \in \nat$
\item 	\lbl{item:simp-eq-classical}
$\xi_{m,n}: X\upr{m+n} \go X\upr{m} \times X\upr{n}$
is an equivalence for all $m, n \in \nat$, and the unique map $X\upr{0} \go
1$ is an equivalence
\item 	\lbl{item:simp-eq-powers}
$\xi^{(k)}: X\upr{k} \go X\upr{1}^k$
is an equivalence for all $k\in\nat$.
\end{enumerate}
\end{propn}

\begin{proof}
As noted above, both the implications
\bref{item:simp-eq-general}$\implies$\bref{item:simp-eq-classical} and 
\bref{item:simp-eq-general}$\implies$\bref{item:simp-eq-powers} are
trivial.  Their converses are straightforward inductions.
\done
\end{proof}

\begin{defn}%
\index{monoidal category!homotopy}
A \demph{homotopy monoidal category} is a functor $X: \Delta^\op \go \Cat$
satisfying the equivalent conditions
\ref{propn:simp-eqs}\bref{item:simp-eq-general}--\bref{item:simp-eq-powers}.
\fcat{HMonCat}%
\glo{HMonCat}
is the category of homotopy monoidal categories and natural
transformations between them.
\end{defn}

A similar definition can be made with \fcat{Top} replacing \Cat\ and
homotopy equivalences replacing categorical equivalences, to give a notion
of `topological monoid%
\index{monoid!topological}
up to homotopy' (Leinster~\cite[\S 4]{UTHM}).  It
can be shown that any loop%
\index{loop space}
space provides an example.  In fact, loop space
theory was where the idea first arose: there, topological monoids up to
homotopy were called `special $\Delta$-spaces'%
\index{Delta-space@$\Delta$-space}%
\index{special}%
\index{simplicial space}
or `special simplicial
spaces' (Segal~\cite{SegCCT},%
\index{Segal, Graeme}
Anderson~\cite{And},
Adams~\cite[p.~63]{Ad}). 

Before I say anything about the comparison with ordinary monoidal
categories, let me explain another route to the notion of homotopy monoidal
category.  

Let \scat{D} be the augmented simplex category~(\ref{eg:str-mon-D}), whose
objects are the (possibly empty) finite totally ordered sets $\lwr{n} =
\{1, \ldots, n\}$.  The fact that $\Delta$ is \scat{D} with the object
\lwr{0} removed is a red herring:%
\index{herring, red}%
\index{simplex category $\Delta$!augmented simplex category@\vs.\ augmented simplex category}% 
\index{augmented simplex category $\scat{D}$!simplex category@\vs.\ simplex category} 
we will not use \emph{this} connection
between $\Delta$ and \scat{D}.  

Now, $\scat{D}$ is the free monoidal category containing a monoid,%
\index{monoid!monoidal category@in monoidal category}
in the
sense that for any monoidal category $(\Eee, \otimes, I)$ (classical, say),
there is an equivalence
\[
\MCwk( (\scat{D},+,0), (\Eee,\otimes,I) )
\eqv
\Mon( \Eee, \otimes, I )
\]
between the category of weak monoidal functors $\scat{D} \go \Eee$ and the
category of monoids in $\Eee$.  For given a weak monoidal functor $\scat{D}
\go \Eee$, the image of any monoid in \scat{D} is a monoid in $\Eee$, and
in particular, the object \lwr{1} of \scat{D} has a unique monoid
structure, giving a monoid in \Eee.  Conversely, given a monoid $A$ in
$\Eee$, there arises a weak monoidal functor $\scat{D} \go \Eee$ sending
\lwr{n} to $A^n$.

Taking $\Eee = \Cat$ describes strict monoidal categories as weak monoidal
functors $(\scat{D},+,0) \go (\Cat,\times,1)$.  Such a weak monoidal
functor is an ordinary functor $W: \scat{D} \go \Cat$ together with
isomorphisms
\begin{equation}	\label{eq:colax-coh-maps}
\omega_{m,n}: W(\lwr{m} + \lwr{n}) \go W\lwr{m} \times W\lwr{n},
\diagspace
\omega_\cdot: W\lwr{0} \go 1
\end{equation}
($m, n \in\nat$) satisfying coherence axioms.  By changing `isomorphisms'
to `equivalences' in the previous sentence, we obtain another notion of
weak monoidal category.  Formally, given monoidal categories $\cat{D}$ and
$\cat{E}$, write $\fcat{MonCat}_\mr{colax}(\cat{D}, \cat{E})$ for the
category of colax monoidal functors $\cat{D} \go \cat{E}$ (as defined
in~\ref{defn:mon-ftr}) and monoidal transformations between them.  The new
`weak monoidal categories' are the objects of the category \fcat{HMonCat'}
defined as follows.
\begin{defn}	\lbl{defn:hmoncat-prime}
\fcat{HMonCat'} is the full subcategory of
\[
\fcat{MonCat}_\mr{colax}((\scat{D},+,0), (\Cat,\times,1))
\]
consisting of the colax monoidal functors $(W, \omega)$ for which each of
the functors~\bref{eq:colax-coh-maps} ($m, n \in \nat$) is an equivalence
of categories.
\end{defn}

You might object that this is useless as a definition of monoidal category,
depending as it does on pre-existing concepts of monoidal category and
colax monoidal functor.  There are several responses.  One is that we could
give an explicit description of what a colax monoidal functor $\scat{D} \go
\Cat$ is, along the lines of the traditional description of cosimplicial
objects by face and degeneracy maps, and this would eliminate the
dependence.  Another is that we will soon show that $\fcat{HMonCat'} \iso
\fcat{HMonCat}$, and $\fcat{HMonCat}$ is defined without mention of
monoidal categories.  A third is that while~\ref{defn:hmoncat-prime} might
not be good as a \emph{definition} of monoidal category, it is a useful
reformulation: for instance, if we change the monoidal category
$(\Cat,\times,1)$ to the monoidal category of chain
complexes and change
categorical equivalence to chain homotopy equivalence, then we obtain a
reasonable notion of homotopy%
\index{homotopy-algebraic structure}%
\index{algebra!differential graded}%
\index{chain complex!homotopy DGA}
differential%
\lbl{p:hty-dgas}
graded algebra.  This exhibits an advantage of the $\scat{D}$ approach over
the $\Delta$ approach: we can use it to discuss homotopy monoids in
monoidal categories where the tensor is not cartesian product.  Much more
on this can be found in my~\cite{UTHM} and~\cite{HAO}.

To show that $\fcat{HMonCat'} \iso \fcat{HMonCat}$, we first establish a
connection between $\Delta$ and $\scat{D}$.

\begin{propn}	\lbl{propn:D-Delta}
Let $\Eee$ be a category with finite products.  Then there is an
isomorphism of categories
\[
\fcat{MonCat}_\mr{colax}((\scat{D},+,0), (\Eee,\times,1))
\iso
\ftrcat{\Delta^\op}{\Eee}.
\]
\end{propn}

\begin{proof}  
This is a special case of a general result on Kleisli%
\index{Kleisli!category}
categories
(Leinster~\cite[3.1.6]{HAO}; beware the different notation).  It can also
be proved directly in the following way.  A functor $W: \scat{D} \go \Eee$
is conventionally depicted as a diagram
\begin{diagram}[width=2em,tight]
\cdot	&
\rTo	&
\cdot	&
\pile{\rTo\\ \lTo\\ \rTo}	&
\cdot	&
\pile{\rTo\\ \lTo\\ \rTo\\ \lTo\\ \rTo}	&
\cdot	&
\cdots
\end{diagram}
of objects and arrows in $\Eee$.  A colax monoidal structure $\omega$ on
$W$ amounts to a pair of maps
\[
W\lwr{m} \og W(\lwr{m} + \lwr{n}) \go W\lwr{n}
\]
for each $m,n\in\nat$, satisfying axioms implying, among other things,
that all of these maps can be built up from the special cases
\[
W\lwr{m} \og W(\lwr{m} + \lwr{1}), 
\diagspace
W(\lwr{1} + \lwr{n}) \go W\lwr{n}.
\]
So a colax monoidal functor $(W, \omega): \scat{D} \go \Eee$ looks like
\begin{diagram}[width=2em,tight]
\cdot	&
\pile{\lGet\\ \rTo\\ \lGet}	&
\cdot	&
\pile{\lGet\\ \rTo\\ \lTo\\ \rTo\\ \lGet}	&
\cdot	&
\pile{\lGet\\ \rTo\\ \lTo\\ \rTo\\ \lTo\\ \rTo\\ \lGet}	&
\cdot	&
\cdots,
\end{diagram}
and this is the conventional picture of a simplicial object in \Eee.
\done
\end{proof}

So if $\Eee$ is a category with ordinary, cartesian, products, simplicial
objects in $\Eee$ are the same as colax monoidal functors from $\scat{D}$
to $\Eee$.  This fails when $\Eee$ is a monoidal category whose tensor
product is not the cartesian product.  It could be argued that in this
situation, it would be better to define a simplicial object in $\Eee$ not
as a functor $\Delta^\op \go \Eee$, but rather as a colax monoidal functor
$\scat{D} \go \Eee$.  For example, it was the colax monoidal version that
made possible the definition of homotopy differential graded algebra
referred to above.

In Proposition~\ref{propn:D-Delta}, if the colax monoidal functor $(W,
\omega)$ corresponds to the simplicial object $X$ then we have
\[
W\lwr{n} = X\upr{n},
\diagspace
\omega_{m,n} = \xi_{m,n}, 
\diagspace
\omega_\cdot = \xi_\cdot
\]
(where $\xi_\cdot: X\upr{0} \go 1$ means $\xi_{n_1, \ldots, n_k}$ in the case
$k=0$).  So by using condition~\bref{item:simp-eq-classical} of
Proposition~\ref{propn:simp-eqs} we obtain:
\begin{cor}
The isomorphism of Proposition~\ref{propn:D-Delta} restricts to an
isomorphism $\fcat{HMonCat'} \iso \fcat{HMonCat}$.
\done
\end{cor}

Our two ways of approaching homotopy monoidal categories, via either nerves
in a finite product category or monoids in a monoidal category, are
therefore equivalent in a strong sense.  The next question is: are they
also equivalent to the standard notion of monoidal category?  I will stop
short of a precise equivalence result, and instead just indicate how to
pass back and forward between homotopy and `ordinary' monoidal categories.

So, let us start from an unbiased%
\index{monoidal category!unbiased}
monoidal category $(A, \otimes, \gamma,
\iota)$ and define from it a homotopy monoidal category $X \in
\fcat{HMonCat}$.  To define $X$ we will need to use unbiased bicategories
(defined in the next section).  For each $n\in\nat$ the ordered set \upr{n}
may be regarded as a category, and so as an unbiased bicategory in which
all 2-cells are identities.  Also, the unbiased monoidal category $A$ may
be regarded as an unbiased bicategory with only one object.  The homotopy
monoidal category $X$ is defined by taking $X\upr{n}$ to be the category
whose objects are all weak functors $\upr{n} \go A$ of unbiased
bicategories and whose maps are transformations of a suitably-chosen kind.
This is just a categorification%
\index{categorification}
of the usual nerve%
\index{nerve!categorified}
construction.

Conversely, let us start with $X \in \fcat{HMonCat}$ and derive from $X$ an
unbiased monoidal category $(A, \otimes, \gamma, \iota)$.  The category $A$
is $X\upr{1}$.  To obtain the rest of the data we first choose for each
$n\in\nat$ a functor $\psi^{(n)}: X\upr{1}^n \go X\upr{n}$ and natural
isomorphisms
\[
\eta^{(n)}: 1 \goiso \xi^{(n)} \of \psi^{(n)},
\diagspace
\epsln^{(n)}: \psi^{(n)} \of \xi^{(n)} \goiso 1
\]
such that $(\psi^{(n)},\xi^{(n)},\eta^{(n)},\epsln^{(n)})$ forms an adjoint
equivalence, which is possible since the functor $\xi^{(n)}$ is an
equivalence~(\ref{propn:eqv-eqv}).  Define $\delta^{(n)}: \upr{1} \go
\upr{n}$ by $\delta^{(n)}(0) = 0$ and $\delta^{(n)}(1) = n$.  Then define
$\otimes_n: A^n \go A$ as the composite
\[
X\upr{1}^n \goby{\psi^{(n)}} X\upr{n} \goby{X\delta^{(n)}} X\upr{1}.
\]
The coherence isomorphisms $\gamma$ and $\iota$ are defined from the
$\eta^{(n)}$'s and $\epsln^{(n)}$'s in a natural way.  The coherence axioms
follow from the fact that we chose \emph{adjoint} equivalences---that is,
they follow from the triangle identities.  So we arrive at an unbiased
monoidal category $(A, \otimes, \gamma, \iota)$.

The processes above determine functors $\UMCwk \oppairu \fcat{HMonCat}$.
It is fairly easy to see that composing one functor with the other does not
yield a functor isomorphic to the identity, either way round.  I believe,
however, that if \fcat{HMonCat} is made into a 2-category in a suitable way
then each composite functor is \emph{equivalent} to the identity.  This
would mean that the notion of homotopy monoidal category is essentially the
same as the other notions of monoidal category that we have discussed.

\paragraph*{}

To summarize the chapter so far: we have formalized the idea of non-strict
monoidal category in various ways, and, with the exception of homotopy
monoidal categories, shown that all the formalizations are equivalent.
Precisely, the following categories are equivalent:%
\index{monoidal category!definitions of}
\begin{quote}
\begin{tabular}{ll}
\MCwk		&(classical monoidal categories)	\\
\UMCwk		&(unbiased monoidal categories)		\\
$\Sigma\hyph\MCwk$&($\Sigma$-monoidal categories, for any plausible
		   $\Sigma$)\\
\fcat{RepMulti}	&(representable multicategories)	
\end{tabular}
\end{quote}
and these equivalences can presumably be extended to equivalences of
2-categories.  The categories
\begin{quote}
\begin{tabular}{ll}
\fcat{HMonCat}	&(homotopy monoidal categories, via simplicial objects)	\\
\fcat{HMonCat'}	&(homotopy monoidal categories, via monoidal functors)	
\end{tabular}
\end{quote}
are isomorphic, and there is reasonable hope that if they are made into
2-categories then they are equivalent to $\MCwk$.  There are still more
notions of monoidal category that we might contemplate---for instance, the
anamonoidal%
\index{anamonoidal category}
categories of Makkai~\cite{MakAAC}---but%
\index{Makkai, Michael}
we leave it at that.

\section{Notions of bicategory}
\lbl{sec:notions-bicat}

Everything that we have done for monoidal categories can also be done for
bicategories.  This is usually at the expense of setting up some slightly
more sophisticated language, which is why things so far have been done for
monoidal categories only.  Here we run through what we have done for
monoidal categories and generalize it to bicategories, noting any wrinkles.

\minihead{Unbiased bicategories}
% \paragraph*{}

\begin{defn}	\lbl{defn:lax-bicat}
A \demph{lax%
\index{bicategory!lax}
bicategory} $B$ (or properly, $(B, \of, \gamma, \iota)$)
consists of
\begin{itemize}
\item a class $B_0$ (often assumed to be a set), whose elements are called
\demph{objects} or \demph{0-cells}%
\index{cell!unbiased bicategory@of unbiased bicategory}
of $B$
\item for each $a,b\in B_0$, a category $B(a,b)$, whose objects are called
\demph{1-cells} and whose morphisms are called \demph{2-cells}
\item for each $n\in \nat$ and $a_0, \ldots, a_n \in B_0$, a functor 
\[
\ofdim{n}: B(a_{n-1}, a_n) \times \cdots \times B(a_0, a_1)
\go
B(a_0, a_n),
\]%
\glo{nfoldcompbicat}%
called \demph{$n$-fold%
\index{n-fold@$n$-fold!composition}
composition} and written
\[
\begin{array}{rcl}
(f_n, \ldots, f_1) 	&\goesto	&(f_n \of\cdots\of f_1)		\\
(\alpha_n, \ldots, \alpha_1)	
			&\goesto
					&(\alpha_n * \cdots * \alpha_1)
\end{array}
\]%
\glo{nfoldcompbicatinfix}\glo{nfoldstar}%
where the $f_i$'s are 1-cells and the $\alpha_i$'s are 2-cells
\item 
a 2-cell
\[
\begin{array}{rl}
\gamma_{((f_n^{k_n}, \ldots, f_n^1), \ldots, (f_1^{k_1}, \ldots,
f_1^1))}: &
((f_n^{k_n} \of \cdots \of f_n^1) \of \cdots \of
(f_1^{k_1} \of \cdots \of f_1^1))\\
&\go 
(f_n^{k_n} \of\cdots\of f_n^1 \of\cdots\of f_1^{k_1}
\of\cdots\of f_1^1)
\end{array}
\]%
\glo{gammabicat}%
for each $n, k_1, \ldots, k_n \in \nat$ and double sequence $((f_n^{k_n},
\ldots, f_n^1), \ldots, (f_1^{k_1}, \ldots, f_1^1))$ of 1-cells for which
the composites above make sense
\item a 2-cell
\[
\iota_f: f \go (f)
\]%
\glo{iotabicat}%
for each 1-cell $f$,
\end{itemize}
satisfying naturality and coherence axioms analogous to those for unbiased
monoidal categories (Definition~\ref{defn:lax-mon-cat}).

A lax bicategory $(B, \of, \gamma, \iota)$ is called an \demph{unbiased%
\index{bicategory!unbiased}
bicategory} (respectively, an \demph{unbiased strict 2-category}) if all
the components of $\gamma$ and $\iota$ are invertible 2-cells
(respectively, identity 2-cells).
\end{defn}

\begin{remarks}{rmks:u-bicat}
\item A lax bicategory with exactly one object is, of course, just a lax
monoidal category, and similarly for the weak and strict versions.
 
\item Unbiased strict 2-categories are in one-to-one correspondence with
ordinary strict 2-categories, easily.
\end{remarks}

As in the case of monoidal categories, there is an abstract version of the
definition of unbiased bicategory phrased in the language of 2-monads.
Previously we took the 2-monad `free strict monoidal category' on $\Cat$;
now we take the 2-monad `free strict 2-category' on $\Cat\hyph\Gph$, the
strict 2-category of $\Cat$-graphs.%
\index{Cat-graph@$\Cat$-graph}
 We met the ordinary category
$\Cat\hyph\Gph$ earlier~(\ref{defn:V-gph}).  The 2-category structure on
$\Cat$ induces a 2-category structure on $\Cat\hyph\Gph$ as follows: given
maps $P, Q: B \go B'$ of $\Cat$-graphs, there are only any 2-cells of the
form
\[
B \ctwomult{P}{Q}{} B'
\]
when $P_0 = Q_0: B_0 \go B'_0$, and in that case a 2-cell $\zeta$ is a
family of natural transformations
\[
\left(P_{a,b} \goby{\zeta_{a,b}} Q_{a,b} \right)_{a,b \in B_0}.
\]
Now, there is a forgetful map $\fcat{Str}\hyph 2\hyph\Cat \go
\Cat\hyph\Gph$ of 2-categories, and this has a left adjoint, so there is an
induced 2-monad $(\blank^*,\mu,\eta)$ on $\Cat\hyph\Gph$.  An unbiased
bicategory is exactly a weak algebra for this 2-monad, and the same applies
in the lax and strict cases.

\begin{defn}
Let $B$ and $B'$ be lax bicategories.  A \demph{lax functor}%
\index{bicategory!unbiased!functor of}%
\index{functor!unbiased bicategories@of unbiased bicategories}
$(P, \pi): B \go B'$ consists of
\begin{itemize}
\item
a function $P_0: B_0 \go B'_0$ (usually just written as $P$)
\item
for each $a,b\in B_0$, a functor $P_{a,b}: B(a,b) \go B'(Pa,Pb)$
\item
for each $n\in\nat$ and composable sequence $f_1, \ldots, f_n$ of 1-cells
of $B$, a 2-cell
\[
\pi_{f_n, \ldots, f_1}:
(Pf_n \of\cdots\of Pf_1)
\go
P(f_n \of\cdots\of f_1),
\]
\end{itemize}
satisfying axioms analogous to those in the definition of lax monoidal
functor~(\ref{defn:u-lax-mon-ftr}).  A \demph{weak functor} (respectively,
\demph{strict functor}) from $B$ to $B'$ is a lax functor $(P, \pi)$ for
which each component of $\pi$ is an invertible 2-cell (respectively,
an identity 2-cell).
\end{defn}

As an example of the benefits of the unbiased approach, consider
`\Hom-functors'.%
\index{Hom-functor}
 For any category $A$ there is a functor
\[
\Hom: A^\op \times A	\go \Set
\]
defined on objects by $(a,b) \goesto A(a,b)$ and on morphisms by
composition---that is, morphisms $a' \goby{f} a$ and $b \goby{g} b'$ in $A$
induce the function
\[
\begin{array}{rrcl}
\Hom(f,g):	&\Hom(a,b)	&\go	&\Hom(a',b'),	\\
		&p		&\goesto&g\of p\of f.
\end{array}
\]
Suppose we want to imitate this construction for bicategories, changing the
category $A$ to a bicategory $B$ and looking for a weak functor $\Hom:
B^\op \times B \go \Cat$.  If we use classical bicategories then we have a
problem: there is no such composite as $g\of p\of f$, and the best we can
do is to choose some substitute such as $(g\of p)\of f$ or $g\of (p\of f)$.
Although we could, say, consistently choose the first option and so arrive
at a weak functor $\Hom$, this is an arbitrary choice.  So there is no
canonical \Hom-functor in the classical world.  In the unbiased world,
however, we can simply take the ternary composite $(g\of p\of f)$, and
everything runs smoothly.

By choosing different strengths of bicategory and of maps between them, we
again obtain $9$ different categories:%
\index{bicategory!nine categories of}
\[
\begin{diagram}[height=1.2em]
\fcat{LaxBicat}_\mr{str}	&\sub	&\fcat{LaxBicat}_\mr{wk}	&\sub
&\fcat{LaxBicat}_\mr{lax} \\
\rotsub	&	&\rotsub	&	&\rotsub	\\	
\UBistr	&\sub	&\UBiwk		&\sub	&\UBilax	\\
\rotsub	&	&\rotsub	&	&\rotsub	\\	
\fcat{Str}\hyph 2\hyph\Cat_\mr{str}	&\sub	&\
\fcat{Str}\hyph 2\hyph\Cat_\mr{wk}	&\sub	&\
\fcat{Str}\hyph 2\hyph\Cat_\mr{lax}.			\\
\end{diagram}
\]%
\glo{ninebicat}%
For instance, $\UBiwk$ is the category of unbiased bicategories and weak
functors.  

Differences from the theory of unbiased monoidal categories emerge when we
try to define transformations between functors between unbiased
bicategories.  This should not come as a surprise given what we already
know in the classical case about transformations and modifications of
bicategories \vs.\ transformations of monoidal
categories~(\ref{eg:mon-cat-bicat-transf}).  More mysteriously, there seems
to be no satisfactorily unbiased way to formulate a definition of
transformation%
\index{transformation!unbiased bicategories@of unbiased bicategories}%
\index{bicategory!unbiased!transformation of}
or modification%
\index{modification!unbiased}%
\index{bicategory!unbiased!modification of}
for unbiased bicategories; it seems that we are forced to grit our teeth
and write down biased-looking definitions.  This done, we obtain a notion
of biequivalence%
\index{biequivalence!unbiased}
of unbiased bicategories.  Just as in the classical
case~(\ref{propn:bieqv-eqv}), biequivalence amounts to the existence of a
weak functor that is essentially surjective on objects and locally an
equivalence.  The $\st$%
\glo{strictcoverbicat}\index{cover, strict}
construction for monoidal
categories~(\ref{thm:eqv-coh-umc}) generalizes without trouble to give

\begin{thm}	%\lbl{thm:eqv-coh-ubicat}%
\index{coherence!bicategories@for bicategories!unbiased}
Every unbiased bicategory is biequivalent to a strict 2-category.  
\done
\end{thm}

\begin{example}
Every topological space $X$ has a fundamental%
\index{fundamental!2-groupoid}
2-groupoid $\Pi_2 X$.%
\glo{fundubicat}
 We
saw how to define $\Pi_2 X$ as a classical bicategory
in~\ref{eg:bicat-Pi}.  Here we consider the unbiased version of $\Pi_2 X$,
in which $n$-fold composition is defined by choosing for each $n\in\nat$ a
reparametrization map $[0,1] \go [0,n]$ (the most obvious choice being
multiplication by $n$).  A 0-cell of the strict%
\index{cover, strict}
cover $\st(\Pi_2 X)$ is
a point of $X$, and a 1-cell is a pair $(n,\gamma)$ where $n\in\nat$ and
$\gamma: [0,n] \go X$ with $\gamma(0) = x$ and $\gamma(n) = y$.  This is
essentially the technique of Moore%
\index{Moore loop}
loops (Adams~\cite[p.~31]{Ad}), used to show that every loop space is
homotopy equivalent to a strict topological monoid.
\end{example}

\minihead{Algebraic notions of bicategory}
% \paragraph*{}

Here we take the very general family of algebraic notions of monoidal
category considered in~\ref{sec:alg-notions} and imitate it for
bicategories.  In other words, we set up a theory of $\Sigma$-bicategories,
where the `signature' $\Sigma$ is a sequence of sets.

Recall that we defined $\Sigma$-monoidal categories in three steps, the
functor $\blank\hyph\MCwk$ being the composite
\[
\Set^\nat \goby{F} 
\Set\hyph\Operad \goby{I_*}
\Cat\hyph\Operad \goby{\Alg_\mr{wk}}
\CAT^\op.
\]
The first two steps create an operad consisting of all the derived tensor
products arising from $\Sigma$ and all the coherence isomorphisms between
them.  The third takes algebras for this operad (or, if you prefer, models
for this theory), which in this case means forming the category of
`monoidal categories' with the kind of products described by the operad.
It is therefore only the third step that we need to change here.

So, given a \Cat-operad $R$, define a category $\fcat{CatAlg}_\mr{wk}(R)$%
\glo{CatAlgwk}
as follows.  An object is a \demph{categorical $R$-algebra},%
\index{categorical algebra for operad}%
\index{Cat-operad@$\Cat$-operad!category over}
 that is, a
\Cat-graph $B$ together with a functor
\[
\act{n}:
R(n) \times B(a_{n-1}, a_n) \times\cdots\times B(a_0, a_1)
\go
B(a_0, a_n)
\]
for each $n\in\nat$ and $a_0, \ldots, a_n \in B_0$, satisfying axioms very
similar to the usual axioms for an algebra for an operad.  (So if $B_0$ has
only one element then $B$ is just an $R$-algebra in the usual sense.)  A
\demph{weak map} $B \go B'$ of categorical $R$-algebras is a map $P: B \go
B'$ of \Cat-graphs together with a natural isomorphism
\begin{diagram}
R(n) \times B(a_{n-1}, a_n) \times\cdots\times B(a_0, a_1)	&
\rTo^{\act{n}}	&B(a_0, a_n)					\\
\dTo<{1 \times P_{a_{n-1}, a_n} \times\cdots\times P_{a_0, a_1}}&
\nent \pi_n	&\dTo>{P_{a_0,a_n}}				\\
R(n) \times B'(a_{n-1}, a_n) \times\cdots\times B'(a_0, a_1)	&
\rTo_{\act{n}}	&B'(Pa_0, Pa_n)					\\
\end{diagram}
for each $n\in\nat$ and $a_0, \ldots, a_n \in B_0$, satisfying axioms like
the ones in the monoidal case (Fig.~\ref{fig:wk-alg-coh},
p.~\pageref{fig:wk-alg-coh}).  This defines a category
$\fcat{CatAlg}_\mr{wk}(R)$.  Defining $\fcat{CatAlg}_\mr{wk}$ on maps in
the only sensible way, we obtain a functor
\[
\fcat{CatAlg}_\mr{wk}: \Cat\hyph\Operad \go \CAT^\op,
\]
and we then define the functor $\blank\hyph\Biwk$%
\glo{blankBi}
as the composite
\[
\Set^\nat \goby{F} 
\Set\hyph\Operad \goby{I_*}
\Cat\hyph\Operad \goby{\fcat{CatAlg}_\mr{wk}}
\CAT^\op.
\]
The lax and strict cases are, of course, done similarly.

\begin{defn}%
\index{Sigma-bicategory@$\Sigma$-bicategory}%
\index{bicategory!Sigma-@$\Sigma$-}
Let $\Sigma\in\Set^\nat$.  A \demph{$\Sigma$-bicategory} is an object of
$\Sigma\hyph\Bilax$ (or equivalently, of $\Sigma\hyph\Biwk$ or
$\Sigma\hyph\Bistr$).  A \demph{lax} (respectively, \demph{weak} or
\demph{strict}) \demph{functor}%
\index{functor!Sigma-bicategories@of $\Sigma$-bicategories}
between $\Sigma$-bicategories is a map in
$\Sigma\hyph\Bilax$ (respectively, $\Sigma\hyph\Biwk$ or
$\Sigma\hyph\Bistr$).
\end{defn}

All of the equivalence results in Section~\ref{sec:alg-notions} go
through.  Hence there are isomorphisms of categories%
\index{coherence!bicategories@for bicategories!unbiased}%
\index{coherence!bicategories@for bicategories!classical}
\[
\fcat{UBicat}_{-} \iso 1\hyph\Bicat_{-},
\diagspace
\Bicat_{-} \iso \Sigma_\mr{c}\hyph\Bicat_{-},
\]
where `$-$' represents any of `$\mr{lax}$', `$\mr{wk}$' or
`$\mr{str}$'.  The proofs are as in the monoidal
case~(\ref{thm:diag-coh-umc},~\ref{thm:diag-coh-mc}) with only cosmetic
changes.  Then, there is the irrelevance%
\index{irrelevance of signature}%
\index{signature!irrelevance of}
of signature theorem, analogous
to~\ref{thm:irrel-sig}:
\[
\Sigma\hyph\Bicat_{-} \eqv \Sigma'\hyph\Bicat_{-}
\]
for all plausible $\Sigma$ and $\Sigma'$, where `$-$' is either
`$\mr{lax}$' or `$\mr{wk}$'.  Again, the proof is essentially unchanged.
As a corollary, unbiased bicategories are the same as classical
bicategories:%
\index{bicategory!unbiased vs. classical@unbiased \vs.\ classical}
\[
\fcat{UBicat}_{-} \eqv \Bicat_{-}
\]
where `$-$' is either `$\mr{lax}$' or `$\mr{wk}$'.  

We saw on p.~\pageref{p:one-level-better} that the equivalence between
unbiased and classical monoidal categories is `one level better' than might
be expected, because it does not refer upwards to transformations.  In the
case of bicategories it is \emph{two} levels better, because modifications
are not needed either.

\minihead{Non-algebraic notions of bicategory}
% \paragraph*{}

I will say much less about these.  

Take representable%
\index{multicategory!representable}
multicategories first.  To formulate a notion of a
bicategory as a multicategory satisfying a representability condition, we
need to use a new kind of multicategory: instead of the arrows looking like
\[
\begin{centredpic}
\begin{picture}(4,8)(-2,0)
\cell{0}{6}{b}{\tinputsvert{a_1}{a_n}}
\cell{0}{6}{t}{\tusualvert{\alpha}}
\cell{0}{2}{t}{\toutputvert{a}}
\end{picture}
\end{centredpic},
\]
they should look like
\[
\setlength{\unitlength}{1mm}
\begin{picture}(36,15)(0,-2)
% Zero-cell marks
\cell{0}{0}{c}{\zmark}
\cell{6}{8}{c}{\zmark}
\cell{36}{0}{c}{\zmark}
% One-cell arrows
\put(0,0){\vector(3,4){6}}
\put(6,8){\vector(3,1){9}}
\put(30,8){\vector(3,-4){6}}
\put(0,0){\vector(1,0){36}}
% Two-cell arrow
\cell{18}{4.5}{c}{\Downarrow}
% \put(18,7){\vector(0,-1){5}}
% Dot-dot-dot
\cell{22}{9.5}{c}{\cdots}
% Labels
\cell{-2.5}{0}{c}{a_0}
\cell{4}{9}{c}{a_1}
\cell{39}{0}{c}{a_n.}
\cell{1}{5}{c}{f_1}
\cell{10}{11.5}{c}{f_2}
\cell{35.5}{5}{c}{f_n}
\cell{18}{-1.5}{c}{g}
\cell{20.5}{4.5}{c}{\alpha}
\end{picture}
\]
We will consider such multicategories later: they are a special kind of
`$\fc$-multicategory' (Example~\ref{eg:fcm-vdisc}), and they also belong to
the world of opetopic structures (Chapter~\ref{ch:opetopic}).  All the
results on monoidal categories as representable multicategories can be
extended unproblematically to bicategories, and the same goes for the
theory of fibrations%
\index{fibration!multicategories@of multicategories}
of multicategories (Leinster~\cite{FM}).

Consider, finally, homotopy monoidal categories.  A homotopy%
\index{bicategory!homotopy}
bicategory can
be defined as a functor $\Delta^\op \go \Cat$%
\index{simplicial object}
satisfying conditions similar
to (but, of course, looser than) those in Proposition~\ref{propn:simp-eqs}.
The new conditions involve pullbacks rather than products, and can be found
by considering nerves%
\index{nerve!category@of category}
of categories in general instead of just nerves of
monoids.  See~\ref{sec:non-alg-defns-n-cat} for further remarks.

\begin{notes}

The notion of unbiased monoidal category has been part of the collective
consciousness for a long while (Kelly~\cite{KelCD},
Hermida~\cite[9.1]{HerRM}).  Around 30 years ago Kelly%
\index{Kelly, Max}
and his
collaborators began investigating 2-monads and 2-dimensional%
\index{algebraic theory!two-dimensional@2-dimensional}
algebraic
theories (see Blackwell, Kelly and Power~\cite{BKP}, for instance), and
they surely knew that unbiased monoidal categories were equivalent to
classical monoidal categories in the way described above, although I have
not been able to find anywhere this is made explicit before my
own~\cite{OHDCT}.  If I am interpreting the (somewhat daunting) literature
correctly, one can also deduce from it some of the more general results
in~\ref{sec:alg-notions} on algebraic notions of monoidal category, but
there are considerable technical issues to be understood before one is in a
position to do so.  Here, in contrast, we have a quick and natural route to
the results.  Essentially, operads take the place of 2-monads as the means
of describing an algebraic theory on $\Cat$.  In order even to \emph{state}
the problem at hand, in any approach, the explicit or implicit use of
operads seems inevitable; an advantage of our approach is that we use
nothing more.

The idea that a monoidal category is a multicategory with enough universal
arrows goes back to Lambek's paper introducing
multicategories~\cite{LamDSCII} and is implicit in the definition of weak
$n$-category proposed by Baez and Dolan~\cite{BDHDA3}, but as far as I know
did not appear in print until Hermida~\cite{HerRM}.%
\index{Hermida, Claudio}

Homotopy monoidal categories were studied in my own~\cite{HAO} paper and
its introductory companion~\cite{UTHM}, in the following wider context:
given an operad $P$ and a monoidal category $\cat{A}$ with a distinguished
class of maps called `equivalences', there is a notion of `homotopy%
\index{homotopy-algebraic structure}%
\index{operad!homotopy-algebra for}
$P$-algebra in $\cat{A}$'.  The case $P = 1$, $\cat{A} = \Cat$ gives
homotopy monoidal categories; the case of homotopy differential graded
algebras was mentioned on p.~\pageref{p:hty-dgas}.  Homotopy monoidal
categories are related to the work of Simpson,%
\index{Simpson, Carlos}
Tamsamani,%
\index{Tamsamani, Zouhair}
To\"en,%
\index{To\"en, Bertrand}
and
Vezzosi%
\index{Vezzosi, Gabriele}
mentioned in the Notes to Chapter~\ref{ch:other-defns}.

\end{notes}

\part{Operads}
\label{part:operads}

\chapter{Generalized Operads and Multicategories: Basics}
\lbl{ch:gom-basics}

\chapterquote{%
Three minutes' thought would suffice to find this out; but thought is
irksome and three minutes is a long time}{%
A.E. Houseman}

\noindent
In a category, an arrow has a single object as its domain and a single
object as its codomain.  In a multicategory, an arrow has a finite sequence
of objects as its domain and a single object as its codomain.  What other
things could we have for the domain of an arrow, keeping a single object as
the codomain?  Could we, for instance, have a tree or a many-dimensional
array of objects?  In different terms (logic or computer science): what can
the input%
\index{input type}
type of an operation be?

In this chapter---the central chapter of the book---we answer these
questions.  We formalize the idea of an input type, and for each input type
we define a corresponding theory of operads and multicategories.  For
instance, the input type might be `finite sequences', and this yields the
theory of ordinary operads and multicategories.  

From now on, operads and multicategories as defined in Chapter~\ref{ch:om}
will be called \demph{plain operads}%
\index{operad!plain}%
\index{plain}
and \demph{plain multicategories}.%
\index{multicategory!plain}
Some mathematicians are used to their operads coming equipped with
symmetric group actions; they should take `plain' as a pun on `planar', to
remind them that these operads do not.

The formal strategy is as follows.  A small category $C$
can be described as consisting of sets and functions%
\index{category!internal}
\[
\begin{slopeydiag}
	&	&C_1	&	&	\\
	&\ldTo<\dom&	&\rdTo>\cod&	\\
C_0	&	&	&	&C_0
\end{slopeydiag}
\diagspace
\begin{array}{c}
C_1 \times_{C_0} C_1 \goby{\comp} C_1,	\\
\\
C_0 \goby{\ids} C_1
\end{array}
\]
satisfying associativity and identity axioms, which can be interpreted as
commutative diagrams in \Set.  Here $C_1 \times_{C_0} C_1$ is a certain
pullback, as explained on p.~\pageref{p:defn-caty-pb}.  Similarly, let $T:
\Set \go \Set$ be the functor sending a set $A$ to the underlying set
$\coprod_{n\in\nat} A^n$ of the free monoid on $A$: then a multicategory
can be described as consisting of sets and functions
\[
\begin{slopeydiag}
	& &C_1 & & \\ &\ldTo<\dom& &\rdTo>\cod& \\ TC_0 & & & &C_0
\end{slopeydiag}
\diagspace
\begin{array}{c}
C_1 \times_{TC_0} TC_1 \goby{\comp} C_1,\\
\\
C_0 \goby{\ids} C_1
\end{array}
\]
satisfying associativity and identity axioms (expressed using the monad
structure on $T$).  Here $C_0$ is the set of all objects, $C_1$ is the set
of all arrows, $\dom$ assigns to an arrow the sequence of objects that is
its domain, and $\cod$ assigns to an arrow the single object that is its
codomain.  The crucial point is that this formalism works for any monad $T$
on any category $\Eee$, as long as $\Eee$ and $T$ satisfy some simple
conditions concerning pullbacks.  This gives a definition of
$T$-multicategory for any such $\Eee$ and $T$; in the terms above, the
pair $(\Eee,T)$ is the `input type'.  So when $T$ is the identity monad on
$\Set$, a $T$-multicategory is an ordinary category, and when $T$ is the
free-monoid monad on $\Set$, a $T$-multicategory is a plain multicategory.

We also define $T$-operads.  As in the plain case, these are simply
 $T$-multicategories $C$ with only one object---or formally, those in which
 $C_0$ is a terminal object of $\Eee$.

There is a canonical notion of an algebra for a $T$-multicategory.  Like
their plain counterparts, generalized operads and multicategories can be
regarded as algebraic theories%
\index{algebraic theory!generalized multicategory as}
(single- and multi-sorted, respectively);
algebras are the accompanying notion of model.  

We start~(\ref{sec:cart-monads}) by describing the simple conditions on
$\Eee$ and $T$ needed to make the definitions work.  Next~(\ref{sec:om})
are the definitions of $T$-multicategory and $T$-operad, and
then~(\ref{sec:algs}) the definition of algebra for a $T$-multicategory (or
$T$-operad).  There are many examples throughout, but some of the most
important ones are done only very briefly; we do them in detail in later
chapters.

\section{Cartesian monads}
\lbl{sec:cart-monads}

In this section we introduce the conditions required of a monad
$(T,\mu,\eta)$ on a category $\Eee$ in order that we may define the notions
of $T$-multicategory and $T$-operad.
\begin{defn}	\lbl{defn:cartesian}
\begin{enumerate}
\item A category $\Eee$ is \demph{cartesian}%
\index{category!cartesian}%
\index{cartesian}
if it has all pullbacks.
\item A functor $\Eee \goby{T} \cat{F}$ is \demph{cartesian}%
\index{functor!cartesian}
if it preserves
pullbacks. 
\item 	\lbl{part:cart-transf}
A natural transformation $\Eee \ctwo{S}{T}{\alpha} \cat{F}$ is
\demph{cartesian}%
\index{transformation!cartesian}
if for each map $A \goby{f} B$ in $\Eee$, the naturality
square
\begin{diagram}[size=2em]
SA		&\rTo^{Sf}	&SB		\\
\dTo<{\alpha_A}	&		&\dTo>{\alpha_B}\\
TA		&\rTo_{Tf}	&TB		\\
\end{diagram}
is a pullback.
\item A monad $(T, T^2\goby{\mu}T, 1\goby{\eta}T)$%
\glo{Tmonad}
on a category $\Eee$ is
\demph{cartesian}%
\index{monad!cartesian}
if the category $\Eee$, the functor $T$, and the natural
transformations $\mu$ and $\eta$ are all cartesian.
\end{enumerate}
\end{defn}

\begin{remarks}{rmks:cart-defns}
\item All of our examples of cartesian categories will have a terminal%
\index{terminal object}
object, hence all finite limits.
\item When the category $\cat{E}$ has a terminal object, a (necessary and)
sufficient condition for the natural transformation $\alpha$
of~\bref{part:cart-transf} to be cartesian is that for each object $A$ of
$\Eee$, the naturality square for the unique map $A \go 1$ is a pullback.  
\item Cartesian categories, cartesian functors and cartesian natural
transformations form a sub-2-category \fcat{CartCat} of \Cat, and a
cartesian monad is exactly a monad in \fcat{CartCat}.  (See
p.~\pageref{p:defn-monad-in-bicaty} for the definition of monad in a
2-category.)
\item As is customary, we often write $T$ to mean the whole monad
$(T,\mu,\eta)$. 
\end{remarks}

The rest of the section is examples.

\begin{example}
The identity monad on a cartesian category is cartesian.
\end{example}

\begin{example}		\lbl{eg:mon-monoids}
Let $\Eee = \Set$ and let $T$ be the free-monoid%
\index{monoid!free}
monad on \Eee.  Certainly
$\Eee$ is cartesian.  An easy calculation shows that the monad $T$ is
cartesian too: Leinster~\cite[1.4(ii)]{GOM}.
\end{example}

\begin{example}		\lbl{eg:comm-not-cart} 
A non-example.  Let $\Eee=\Set$ and let $(T, \mu, \eta)$ be the free
commutative%
\index{monoid!commutative!free}
monoid monad.  This fails to be cartesian on two counts: $\mu$
is not cartesian (for instance, its naturality square at the unique map
$2\go 1$ is not a pullback), and the functor $T$ does not preserve
pullbacks.  Let us show the latter in detail, using an argument of
Weber~\cite[2.7.2]{Web}.%
\index{Weber, Mark}
 First note that for any set $A$,
\[
TA = \coprod_{n\in\nat} A^n / S_n
\]
where $S_n$ is the $n$th symmetric group acting on $A^n$ in the natural
way.  Write $[a_1, \ldots, a_n]$ for the equivalence class of $(a_1,
\ldots, a_n) \in A^n$ under this action, so that $[a_1, \ldots, a_n] =
[b_1, \ldots, b_n]$ if and only if there exists $\sigma \in S_n$ such that
$b_i = a_{\sigma i}$ for all $i$.  Now let $x, x', y, y', z$ be distinct
formal symbols and consider applying $T$ to the pullback square
\[
\begin{diagram}[size=2em]
\{(x,y), (x,y'), (x,y'), (x',y')\}	&\rTo	&\{y,y'\}	\\
\dTo					&	&\dTo		\\
\{x,x'\}				&\rTo	&\{z\}.		\\
\end{diagram}
\]
There are distinct elements
\[
[(x,y), (x',y')],
\diagspace
[(x,y'), (x',y)]
\]
of $T \{(x,y), (x,y'), (x,y'), (x',y')\}$, yet both map to $[x,x']
\in T\{x,x'\}$ and to $[y,y'] \in T\{y,y'\}$.  Hence the image of the
square under $T$ is not a pullback.
\end{example}

\begin{example}		\lbl{eg:mon-CJ}%
\index{algebraic theory|(}%
\index{strongly regular theory|(}
Any algebraic theory gives rise to a monad on $\Set$ (its free algebra
monad), and if the theory is strongly regular in the sense
of~\ref{eg:opd-sr} then the monad is cartesian.  Carboni and Johnstone
proved this first~\cite{CJ}; an alternative proof is in
Appendix~\ref{app:special-cart}.  The result implies that the free monoid
monad~(\ref{eg:mon-monoids}) is cartesian, and, as suspected
in~\ref{eg:opd-sr}, that the theory of commutative
monoids~(\ref{eg:comm-not-cart}) is not strongly regular.  Further examples
appear below.
\end{example}

\begin{example}	\lbl{eg:mon-mon-with-inv}
Let $\Eee=\Set$ and let $T$ be the monad corresponding to the theory of
monoids with involution.  By definition, a \demph{monoid with involution}%
\index{monoid!involution@with involution}
is a monoid equipped with an endomorphism whose composite with itself is
the identity; so this theory is defined by the usual operations and
equations for the theory of monoids together with a unary operation
$\blank^\circ$ satisfying
\[
x^{\circ\circ} = x,
\diagspace
(x\cdot y)^\circ = x^\circ \cdot y^\circ, 
\diagspace
1^\circ = 1.
\]
These equations are strongly regular, so the monad $T$ is cartesian.  Any
abelian group has an underlying monoid with involution, given by $x^\circ =
x^{-1}$.   
\end{example}

\begin{example}		\lbl{eg:mon-mon-with-anti-inv}
More interestingly, let $\Eee=\Set$ and let $T$ be the monad for
\demph{monoids with anti-involution},%
\index{monoid!anti-involution@with anti-involution}
where an anti-involution on a monoid
is a unary operation $\blank^\circ$ satisfying
\[
x^{\circ\circ} = x,
\diagspace
(x\cdot y)^\circ = y^\circ \cdot x^\circ, 
\diagspace
1^\circ = 1.
\]
Any group whatsoever has an underlying monoid with anti-involution, again
given by $x^\circ = x^{-1}$.  Now, this is \emph{not} a strongly regular
presentation of the theory, and in fact there is no strongly regular
presentation, but $T$ \emph{is} a cartesian monad.  So not all cartesian
monads on $\Set$ arise from strongly regular theories.  These assertions
are proved in Example~\ref{eg:fam-rep-not-opdc}; that $T$ is cartesian was
pointed out to me by Peter Johnstone.%
\index{Johnstone, Peter}
\end{example}

\begin{example}		\lbl{eg:mon-exceptions}
Let $\Eee = \Set$ and fix a set $S$.  The endofunctor $S + \dashbk$%
\index{coproduct!monad from}
on
$\Eee$ has a natural monad structure, and the monad is cartesian, corresponding
to the algebraic theory consisting only of one constant for each member of
$S$.  In particular, if $S=1$ then this is the theory of pointed sets.
\end{example}

\begin{example}		\lbl{eg:mon-action}%
\index{monoid!action of}%
\index{action!monoid@of monoid}
Algebraic theories generated by just unary operations correspond to monads
on \Set\ of the form $M\times\dashbk$, where $M$ is a monoid and the unit
and multiplication of the monad are given by those of the monoid.  Since
any equation formed from unary operations and a single variable is strongly
regular, the monad $M\times\dashbk$ on \Set\ is always cartesian.
\end{example}

\begin{example}		\lbl{eg:mon-free-cl-opd}%
\index{operad!free}
Let $\Eee = \Set^\nat$ and let $T$ be the monad `free plain operad' on
$\Eee$.  As we saw in~\ref{sec:om-further}, the functor $T$ forms trees%
\index{tree!vertices labelled@with vertices labelled}
with labelled vertices: for instance, if $A \in \Eee$, $a_1 \in A(3)$,
$a_2 \in A(1)$, and $a_3 \in A(2)$, then
\[
% \drmk{tree corr to tensor product }
% \otimes_{a_1}(\otimes_{a_2}(-), -, \otimes_{a_3}(-,-))
\setlength{\unitlength}{1.5em}
\begin{picture}(6,3)(-1.5,0)
% 0th layer
\put(1.5,0){\line(0,1){1}}
% 1st layer
\cell{1.5}{1}{c}{\vx}
\put(1.5,1){\line(-3,2){1.5}}
\put(1.5,1){\line(0,1){1}}
\put(1.5,1){\line(3,2){1.5}}
\cell{0}{2}{c}{\vx}
\cell{3}{2}{c}{\vx}
% top layer
\put(0,2){\line(0,1){1}}
\put(3,2){\line(-1,1){1}}
\put(3,2){\line(1,1){1}}
% labels
\cell{1.8}{1}{tl}{a_1}
\cell{-0.3}{2}{tr}{a_2}
\cell{3.3}{2}{tl}{a_3}
\end{picture}
\]
is an element of $(TA)(4)$.  That $T$ is cartesian follows from theory we
develop later~(\ref{eg:free-cl-opd-cart}).  
\end{example}

\begin{example}		\lbl{eg:mon-tree}
Consider the finitary algebraic theory on $\Set$ generated by one $n$-ary
operation for each $n\in\nat$ and no equations. This theory is strongly
regular, so the induced monad $(T, \mu, \eta)$ on $\Set$ is cartesian.

Given a set $A$, the set $TA$ can be described inductively by
\begin{itemize}
\item if $a\in A$ then $a\in TA$
\item if $t_1, \ldots, t_n \in TA$ then $(t_1, \ldots, t_n) \in TA$.
\end{itemize}
We have already looked at the case $A=1$, where $TA$ is the set of
unlabelled trees~(\ref{eg:opd-of-trees}).  Similarly, an
element of $TA$ can be drawn as a tree%
\index{tree!leaves labelled@with leaves labelled}
whose leaves are labelled by
elements of $A$: for instance,
\[
((a_1, a_2, ()), a_3, (a_4, a_5)) \in TA
\]
is drawn as
\[
\begin{array}{c}
\setlength{\unitlength}{1.5em}
\begin{picture}(7,4.2)(-0.5,0)
% 0th layer
\put(3,0){\line(0,1){1}}
% 1st layer
\cell{3}{1}{c}{\vx}
\put(3,1){\line(-2,1){2}}
\put(3,1){\line(0,1){1}}
\put(3,1){\line(2,1){2}}
% 2nd layer
\cell{1}{2}{c}{\vx}
\put(1,2){\line(-1,1){1}}
\put(1,2){\line(0,1){1}}
\put(1,2){\line(1,1){1}}
\cell{5}{2}{c}{\vx}
\put(5,2){\line(-1,1){1}}
\put(5,2){\line(1,1){1}}
% invisible top layer
\cell{2}{3}{c}{\vx}
% labels
\cell{0}{3.2}{b}{a_1}
\cell{1}{3.2}{b}{a_2}
\cell{3}{2.2}{b}{a_3}
\cell{4}{3.2}{b}{a_4}
\cell{6}{3.2}{b}{a_5}
\end{picture}
\end{array}.
\]
(Contrast this with the previous example, where the \emph{vertices} were
labelled.)  The unit $A\go TA$ is 
\[
a 
\goesto
\begin{array}{c}
\setlength{\unitlength}{1.5em}
\begin{picture}(1,1.6)(-0.5,0)
\put(0,0){\line(0,1){1}}
\cell{0}{1.2}{b}{a}
\end{picture}
\end{array},
\]
and multiplication $T^2 A \go TA$ takes a tree whose leaves are labelled by
elements of $TA$ (for instance,
\[
\begin{array}{c}
\setlength{\unitlength}{1.5em}
\begin{picture}(4,3.7)(-0.5,0)
% 0th layer
\put(2,0){\line(0,1){1}}
% 1st layer
\cell{2}{1}{c}{\vx}
\put(2,1){\line(-1,1){1}}
\put(2,1){\line(1,1){1}}
% 2nd layer
\cell{1}{2}{c}{\vx}
\put(1,2){\line(-1,1){1}}
\put(1,2){\line(1,1){1}}
% invisible 3rd layer
\cell{2}{3}{c}{\vx}
% labels
\cell{0}{3.2}{b}{t_1}
\cell{3}{2.2}{b}{t_2}
\end{picture}
\end{array},
\]
where
\[
t_1=
\begin{array}{c}
\setlength{\unitlength}{1.5em}
\begin{picture}(3,3.7)(-0.5,0)
% 0th layer
\put(1,0){\line(0,1){1}}
% 1st layer
\cell{1}{1}{c}{\vx}
\put(1,1){\line(-1,1){1}}
\put(1,1){\line(1,1){1}}
% 2nd layer
\cell{2}{2}{c}{\vx}
\put(2,2){\line(0,1){1}}
% labels
\cell{0}{2.2}{b}{a_1}
\cell{2}{3.2}{b}{a_2}
\end{picture}
\end{array},
\diagspace
t_2=
\begin{array}{c}
\setlength{\unitlength}{1.5em}
\begin{picture}(3,3.7)(-0.5,0)
% 0th layer
\put(1,0){\line(0,1){1}}
% 1st layer
\cell{1}{1}{c}{\vx}
\put(1,1){\line(0,1){1}}
% 2nd layer
\cell{1}{2}{c}{\vx}
\put(1,2){\line(-1,1){1}}
\put(1,2){\line(0,1){1}}
\put(1,2){\line(1,1){1}}
% 3rd layer (invisible)
\cell{1}{3}{c}{\vx}
% labels
\cell{0}{3.2}{b}{a_3}
\cell{2}{3.2}{b}{a_4}
\end{picture}
\end{array}
)
\]
and expands the labels to produce a tree whose leaves are labelled by
elements of $A$ (here,
\[
\begin{array}{c}
\setlength{\unitlength}{1.5em}
\begin{picture}(7,5.7)(-0.5,0)
% 0th layer
\put(3.5,0){\line(0,1){1}}
% 1st layer
\cell{3.5}{1}{c}{\vx}
\put(3.5,1){\line(-3,2){1.5}}
\put(3.5,1){\line(3,2){1.5}}
% 2nd layer
\cell{2}{2}{c}{\vx}
\put(2,2){\line(-1,1){1}}
\put(2,2){\line(1,1){1}}
\cell{5}{2}{c}{\vx}
\put(5,2){\line(0,1){1}}
% 3rd layer
\cell{1}{3}{c}{\vx}
\put(1,3){\line(-1,1){1}}
\put(1,3){\line(1,1){1}}
\cell{3}{3}{c}{\vx}
\cell{5}{3}{c}{\vx}
\put(5,3){\line(-1,1){1}}
\put(5,3){\line(0,1){1}}
\put(5,3){\line(1,1){1}}
% 4th layer
\cell{2}{4}{c}{\vx}
\put(2,4){\line(0,1){1}}
\cell{5}{4}{c}{\vx}
% labels
\cell{0}{4.2}{b}{a_1}
\cell{2}{5.2}{b}{a_2}
\cell{4}{4.2}{b}{a_3}
\cell{6}{4.2}{b}{a_4}
\end{picture}
\end{array}
).
\]
\end{example}

\begin{example}%
\index{algebraic theory!free}
Similar statements can be made for any free theory.  That is, given any
sequence $(\Sigma(n))_{n\in \nat}$ of sets, the finitary algebraic theory
on \Set\ generated by one $n$-ary operation for each $\sigma \in
\Sigma(n)$, and no equations, is strongly regular, so the induced monad
$(T,\mu,\eta)$ on $\Set$ is cartesian.  (The previous example was the case
$\Sigma(n)=1$ for all $n$.)  Given a set $A$, an element of $TA$ can be
drawn as a tree%
\index{tree!leaves and vertices labelled@with leaves and vertices labelled}
in which the leaves are labelled by elements of $A$ and the
vertices with $n$ branches coming up out of them are labelled by elements
of $\Sigma(n)$: in other words, as a diagram like the domain or codomain
of~\bref{eq:Sigma-mon-cat-iso} (p.~\pageref{eq:Sigma-mon-cat-iso}).%
\index{algebraic theory|)}%
\index{strongly regular theory|)}
\end{example}

\begin{example}		\lbl{eg:mon-free-topo-monoid}%
\index{monoid!topological!free}
The monad `free topological monoid' on $\Top$, whose functor part sends a
space $A$ to the disjoint union of cartesian powers
$\coprod_{n\in\nat}A^n$, is cartesian (by direct calculation).
\end{example}

\begin{example}		\lbl{eg:mon-free-str-mon-cat}%
\index{monoidal category!strict!free}
Similarly, the monad `free strict monoidal category' ($=$ free
monoid) on $\Cat$ is cartesian.  In fact, the free monoid monad on a
category $\Eee$ is always cartesian provided that $\Eee$ satisfies the
(co)limit conditions necessary to ensure that the usual free monoid
construction $A \goesto \coprod_{n\in\nat}A^n$ works in $\Eee$.
\end{example}

\begin{example}		\lbl{eg:mnd-sym-mon}%
\index{monoidal category!symmetric!free strict}
The monad `free symmetric strict monoidal category' on \Cat\ is also
cartesian, as a (lengthy) direct calculation shows.  If we had taken strict
symmetric monoidal categories instead---in other words, insisted that the
symmetries were strict---then we would just be looking at commutative
monoids in \Cat, and the monad would fail to be cartesian by the argument
of Example~\ref{eg:comm-not-cart}.  So weakening%
\index{weakening}%
\index{categorification}
the equality $x \otimes y
= y\otimes x$ to an isomorphism makes a bad monad good.
\end{example}

\begin{example}		\lbl{eg:fc-mnd}
Let $\scat{H}$ be the category $(0 \parpair{\sigma}{\tau} 1)$, so that
\ftrcat{\scat{H}^\op}{\Set} is the category of directed graphs.%
\index{graph!directed}
 The
forgetful functor $\Cat \go \ftrcat{\scat{H}^\op}{\Set}$ has a left adjoint
and therefore induces a monad $\fc$%
\glo{fc}
(`free category')%
\index{category!free (fc)@free ($\fc$)}
on
$\ftrcat{\scat{H}^\op}{\Set}$.  It follows from later
theory~(\ref{eg:fc-cart}) that $\fc$ is cartesian.
\end{example}

\begin{example}		\lbl{eg:glob-mnd}
In Part~\ref{part:n-categories} we consider the free strict
$\omega$-category%
\index{omega-category@$\omega$-category!strict!free}
monad on the category of globular%
\index{globular set}
sets, and, for
$n\in\nat$, the free strict $n$-category%
\index{n-category@$n$-category!strict!free}
monad on the category of
$n$-globular sets.  All of these monads are cartesian.  The previous
example is the case $n=1$.
\end{example}

\begin{example}		\lbl{eg:cubical-mnd}
In~\ref{sec:cl-strict} we looked at strict double%
\index{double category!strict!free}
categories and, more
generally, at strict $n$-tuple categories.%
\index{n-tuple category@$n$-tuple category!strict!free}
 In particular, we saw that a
strict $n$-tuple category could be described as a functor $(\scat{H}^n)^\op
\go \Set$ with extra structure, where $\scat{H}$ is as in~\ref{eg:fc-mnd}.
The forgetful functor from the category of strict $n$-tuple categories and
strict maps between them to the functor category
$\ftrcat{(\scat{H}^n)^\op}{\Set}$ has a left adjoint, the adjunction is
monadic, and the induced monad on $\ftrcat{(\scat{H}^n)^\op}{\Set}$ is
cartesian.  This can be shown by a similar method to that used for
strict $n$-categories in Appendix~\ref{app:free-strict}.
\end{example}

Alert readers may have noticed that nearly every one of the above examples 
of a cartesian monad on $\Set$ is, in fact, the free-algebra monad for a
certain plain operad.  This is no coincidence, as we discover
in~\ref{sec:alt-app}.

\section{Operads and multicategories}
\lbl{sec:om}

We now define `$T$-multicategory' and `$T$-operad', for any cartesian monad
$T$ on a cartesian category $\Eee$.  That is, for each such $\Eee$ and $T$
we define a category $T\hyph\Multicat$ of $T$-multicategories and a full
subcategory $T\hyph\Operad$ consisting of the $T$-operads.  In the case
$\Eee=\Set$ and $T=\textrm{(free monoid)}$ these are, respectively,
$\Multicat$ and $\Operad$; in the case $\Eee=\Set$ and $T=\id$ they are
$\Cat$ and $\fcat{Monoid}$.%
\glo{Monoid}

The strategy for making these definitions is as described in the
introduction to this chapter, dressed up a little: instead of handling the
data and axioms for a $T$-multicategory directly, we introduce a bicategory
$\Sp{\Eee}{T}$%
\glo{Sp}
and define a $T$-multicategory as a monad in $\Sp{\Eee}{T}$.
This amounts to the same thing, as we shall see.  (All we are doing is
generalizing the description of a small category as a monad in the
bicategory%
\index{bicategory!spans@of spans}%
\index{span}
of spans: B\'enabou~\cite[5.4.3]{Ben}.)%
\index{Benabou, Jean@B\'enabou, Jean}

\begin{defn}	\lbl{defn:T-spans} 
For any cartesian monad $(T,\mu,\eta)$ on a cartesian category \Eee, the
bicategory $\Sp{\Eee}{T}$ is defined as follows:
\begin{description}
\item[0-cells] are objects $E$ of \Eee
\item[1-cells $E \go E'$] are diagrams
\begin{slopeydiag}
	&	&M	&	&	\\
	&\ldTo<d&	&\rdTo>c&	\\
TE	&	&	&	&E'	\\
\end{slopeydiag}
in \Eee
\item[2-cells $(M,d,c) \go (N,q,p)$] are maps $M \go N$ in $\Eee$ such that
\begin{slopeydiag}
	&	&M	&	&	\\
	&\ldTo<d&	&\rdTo>c&	\\
TE	&	&\dTo	&	&E'	\\
	&\luTo<q&	&\ruTo>p&	\\
	&	&N	&	&	\\
\end{slopeydiag}
commutes
\item[1-cell composition:] the composite of $1$-cells
\[
\begin{slopeydiag}
	&	&M	&	&	\\
	&\ldTo<{d}&	&\rdTo>{c}&	\\
TE	&	&	&	&E'	\\
\end{slopeydiag}
\diagspace
\begin{slopeydiag}
	&	&M'	&	&	\\
	&\ldTo<{d'}&	&\rdTo>{c'}&	\\
TE'	&	&	&	&E''	\\
\end{slopeydiag}
\]
is given by composing along the upper slopes of the diagram
\[
\begin{slopeydiag}
   &       &   &       &   &       &M'\of M\Spbk& &   &       &   \\
   &       &   &       &   &\ldTo  &      &\rdTo  &   &       &   \\
   &       &   &       &TM &       &      &       &M' &       &   \\
   &       &   &\ldTo<{Td}&&\rdTo>{Tc}&   &\ldTo<{d'}&&\rdTo>{c'}&\\
   &       &T^2 E&     &   &       &TE'   &       &   &       &E''\\
   &\ldTo<{\mu_E}&&    &   &       &      &       &   &       &   \\
TE &       &   &       &   &       &      &       &   &       &   \\
\end{slopeydiag}
\]
in $\Eee$, where the right-angle mark in the top square indicates that the
square is a pullback and we assume from now on that a particular
choice%
\index{pullback!choice of}
of pullbacks in $\Eee$ has been made

\item[1-cell identities:] the identity on $E$ is
\begin{slopeydiag}
	&	&E	&	&	\\
	&\ldTo<{\eta_E}&&\rdTo>{1}&	\\
TE	&	&	&	&E	\\
\end{slopeydiag}
\item[2-cell compositions and identities] are defined in the evident way 
\item[coherence 2-cells:] the associativity and unit 2-cells are defined
using the universal property of pullback.
\end{description}
\end{defn}

Since the choice of pullbacks in $\Eee$ was arbitrary, it is inevitable
that composition of $1$-cells in $\Sp{\Eee}{T}$ does not obey strict
associativity or unit laws.  That it obeys them up to isomorphism is a
consequence of $T$ being cartesian.  Changing the choice%
\index{pullback!choice of}
of pullbacks in
$\Eee$ only changes the bicategory $\Sp{\Eee}{T}$ up to isomorphism (in the
category of bicategories and weak functors): see
p.~\pageref{p:change-of-shape-rmks}.

Here is the most important definition in this book.  It is due to
Burroni~\cite{Bur},%
\index{Burroni, Albert}
and in the form presented here uses the notion of monad
in a bicategory (p.~\pageref{p:defn-monad-in-bicaty}).
\begin{defn}		\lbl{defn:T-multicat}
Let $T$ be a cartesian monad on a cartesian category $\Eee$.  A
\demph{$T$-multicategory}%
\index{generalized multicategory}
is a monad in the bicategory $\Sp{\Eee}{T}$.
\end{defn}

A $T$-multicategory $C$ therefore consists of a diagram
\[
\begin{slopeydiag}
	&		&C_1		&		&	\\
	&\ldTo<\dom	&		&\rdTo>\cod	&	\\
TC_0	&		&		&		&C_0	\\
\end{slopeydiag}
\glo{C0genmti}\glo{C1genmti}\glo{domgenmti}\glo{codgenmti}
\]
in $\Eee$ together with maps 
\[
C_{1}\of C_{1} = C_1 \times_{TC_0} TC_1 \goby{\comp} C_{1},
\diagspace
C_{0} \goby{\ids} C_{1}
\]%
\glo{compgenmti}\glo{idsgenmti}%
satisfying associativity and identity axioms---exactly as promised at the
start of the chapter.

\begin{defn}		\lbl{defn:T-operad}
Let $T$ be a cartesian monad on a cartesian category $\Eee$.  A
\demph{$T$-operad}%
\index{generalized operad}
is a $T$-multicategory $C$ such that $C_0$ is a terminal
object of $\Eee$.
\end{defn}
Just as a $T$-multicategory is a generalized category, a $T$-operad is a
generalized monoid.  Explicitly, a monoid in a category $\Eee$ with finite
limits consists of an object $M$ of $\Eee$
together with maps 
\[
M \times M \goby{\mult} M,
\diagspace
1 \goby{\unit} M
\]
in $\Eee$, satisfying associativity and identity axioms.  Take a cartesian
monad $T$ on $\Eee$: then a $T$-operad consists of an object of
$\Eee$ over $T1$, say $P \goby{d} T1$, together with maps
\[
P \times_{T1} TP \goby{\comp} P,
\diagspace
1 \goby{\ids} P
\]
over $T1$ in $\Eee$, again satisfying associativity and identity axioms.

A map of $T$-multicategories is a map of the underlying `graphs' 
preserving composition and identities, as follows.
\begin{defn}	\lbl{defn:T-graph}
Let $T$ be a cartesian monad on a cartesian category $\Eee$. A
\demph{$T$-graph}%
\index{generalized graph ($T$-graph)}%
\index{T-graph@$T$-graph}\index{graph!generalized ($T$-)}
is a diagram
\[
\begin{slopeydiag}
	&		&C_1		&		&	\\
	&\ldTo<\dom	&		&\rdTo>\cod	&	\\
TC_0	&		&		&		&C_0	\\
\end{slopeydiag}
\]
in $\Eee$ (that is, an endomorphism 1-cell in $\Sp{\Eee}{T}$).  A \demph{map
$C \goby{f} C'$ of $T$-graphs} is a pair $(C_0 \goby{f_0} C'_0,\ 
C_1 \goby{f_1} C'_1)$ of maps in $\Eee$ such that
\[
\begin{slopeydiag}
	&	&C_1	&	&	\\
	&\ldTo	&	&\rdTo	&	\\
TC_{0}	&	&\dTo>{f_1}&	&C_0	\\
\dTo<{Tf_{0}}&	&C'_1&	&\dTo>{f_0}\\
	&\ldTo	&	&\rdTo	&	\\
TC'_0&	&	&	&C'_0\\
\end{slopeydiag}
\]
commutes.  The category of $T$-graphs is written $T\hyph\Graph$.%
\glo{TGraph}
\end{defn}

This definition uses two different notions of a map between objects of
\Eee: on the one hand, genuine maps in \Eee, and on the other, spans ($=$
1-cells of $\Sp{\Eee}{T}$).  In Chapter~\ref{ch:fcm} we will integrate the
objects of $\Eee$ and these two different kinds of map into a single
structure, an `\fc-multicategory'.

\begin{defn}	\lbl{defn:multifunctor}%
\index{generalized multicategory!map of}
A \demph{map $C \goby{f} C'$ of $T$-multicategories} is a map $f$ of
their underlying graphs such that the diagrams 
\[
\begin{diagram}[size=2em]
C_0		&\rTo^{\ids}		&C_1		\\	
\dTo<{f_0}	&			&\dTo>{f_1}	\\
C'_0	&\rTo_{\ids}		&C'_1	\\
\end{diagram}
\diagspace
\begin{diagram}[size=2em]
C_1\of C_1	&\rTo^{\comp}		&C_1		\\
\dTo<{f_1 *f_1}	&			&\dTo>{f_1}	\\
C'_1\of C'_1	&\rTo_{\comp}		&C'_1		\\
\end{diagram}
\]
commute, where $f_1 * f_1$ is the evident map induced by two copies of $C_1
\goby{f_1} C'_1$.  The category of $T$-multicategories and maps between
them is written $T\hyph\Multicat$.%
\glo{TMulticat}
 The full subcategory consisting of
$T$-operads is written $T\hyph\Operad$.%
\glo{TOperad}
\end{defn}

So for any cartesian monad $T$, we have categories and functors
\[
T\hyph\Operad \rIncl 
T\hyph\Multicat \goby{\textrm{\footnotesize forgetful}}
T\hyph\Graph.
\]
When the extra clarity is needed, we will refer to $T$-multicategories as
\demph{\Cartpr-multicategories}%
\index{ET-multicategory etc.@`$\Cartpr$-multicategory' etc.}
and to the category they form as
$\Cartpr\hyph\Multicat$; similarly for operads and graphs. 

We now look at some examples of generalized multicategories.  In many of
the most interesting ones the input of each operation/arrow forms quite a
complicated shape, such as a diagram of pasted-together higher-dimensional
cells.  These examples are only described briefly here, with proper
discussions postponed to later chapters.

We start with the two motivating cases.

\begin{example}
Let $T$ be the identity monad on $\Eee=\Set$.  Then $\Sp{\Eee}{T}$ is what
is usually called the `bicategory%
\index{bicategory!spans@of spans}%
\index{span}
of spans' (B\'enabou~\cite[2.6]{Ben}), and a
monad in $\Sp{\Eee}{T}$ is just a small category.  So
\[
(\Set,\id)\hyph\Multicat 
\eqv
\Cat,
\diagspace
(\Set,\id)\hyph\Operad
\eqv
\fcat{Monoid}.
\]
More generally, if $\Eee$ is any cartesian category then
$(\Eee, \id)$-multicategories are categories%
\index{category!internal}
in $\Eee$ and
$(\Eee, \id)$-operads are monoids%
\index{monoid!monoidal category@in monoidal category}
in $\Eee$.
\end{example}

\begin{example}	\lbl{eg:cl-oms-are-gen}%
\index{monoid!free}
Let $T$ be the free monoid monad on the category $\Eee$ of sets.  A
$T$-graph
\[
TC_0 = \coprod_{n\in\nat}C_0^n \ogby{\dom} C_1 \goby{\cod} C_0
\]
amounts to a set $C_0$ `of objects' together with a set $C(a_1, \ldots,
a_n; a)$ for each $n\geq 0$ and $a_1, \ldots, a_n, a \in C_0$.  
The composite 1-cell
\[
\begin{slopeydiag}
           &         &C_1 \of C_1      &          &           \\
           &\ldTo<{\dom'}&             &\rdTo>{\cod'}&        \\
TC_0       &         &                 &          &C_0        \\
\end{slopeydiag}
\]
in $\Sp{\Eee}{T}$ is as follows: the set $C_1 \of C_1$ at the apex is 
\[
C_1 \times_{TC_0} TC_1	
=	
\coprod_{%
\begin{scriptarray}
\scriptstyle
n, k_1, \ldots, k_n \in\nat,	\\
\scriptstyle
a_i^j, a_i, a \in C_0
\end{scriptarray}}	
\begin{array}[t]{l}
C(a_1, \ldots, a_n; a) \times
C(a_1^1, \ldots, a_1^{k_1}; a_1) \times\cdots
\\
\times C(a_n^1, \ldots, a_n^{k_n}; a_n),
\end{array}
\]
an element of which looks like the left-hand side of
Fig.~\ref{fig:multi-comp} (p.~\pageref{fig:multi-comp}); the function
$\dom'$ sends this element to 
\[
(a_1^1, \ldots, a_1^{k_1}, \ldots, a_n^1, \ldots, a_n^{k_n}) \in TC_0,
\]
and $\cod'$ sends it to $a$.  The identity 1-cell on $C_0$ in
$\Sp{\Eee}{T}$ is 
\[
\begin{slopeydiag}
           &         &C_0           &          &           \\
           &\ldTo<{\eta_{C_0}}&     &\rdTo>{1} &           \\
TC_0       &         &              &          &C_0,       \\
\end{slopeydiag}
\]
and $\eta_{C_0}$ sends $a\in C_0$ to $(a) \in TC_0$.  A
$T$-multicategory structure on a $T$-graph $C$ therefore consists of a
function $\comp$ as in Fig.~\ref{fig:multi-comp} and a function $\ids$
assigning to each object $a\in C_0$ an `identity' element $1_a \in C(a;a)$,
obeying associativity and identity laws.  So a $T$-multicategory is just a
plain%
\index{multicategory!generalized multicategory@as generalized multicategory}
multicategory; indeed, there are equivalences of categories
\[
T\hyph\Multicat \eqv \Multicat, 
\diagspace
T\hyph\Operad \eqv \Operad.
\index{operad!generalized operad@as generalized operad}
\]
\end{example}

\begin{example}	\lbl{eg:sym-ops}
Suppose we try to realize symmetric operads~(\ref{defn:sym-mti})%
\index{operad!symmetric vs. generalized@symmetric \vs.\ generalized}%
\index{multicategory!symmetric vs. generalized@symmetric \vs.\ generalized!multicategory for multicategories}
as
$T$-operads for some $T$.  A first attempt might be to take the free
commutative monoid monad $T$ on $\Set$. But this is both misguided and
doomed to failure: misguided because if $P$ is a symmetric operad then the
maps $\dashbk\cdot\sigma: P(n) \go P(n)$ coming from permutations
$\sigma\in S_n$ are only isomorphisms, not%
\index{weakening}
identities; and doomed to
failure because $T$ is not cartesian~(\ref{eg:comm-not-cart}), which
prevents us from making a definition of $T$-operad---in particular, from
expressing associativity and identity laws.  A better idea is to take the
free symmetric%
\index{monoidal category!symmetric!free strict}
strict monoidal category monad on $\Cat$, thus replacing
identities by isomorphisms: see~\ref{eg:mti-sym} below.
\end{example}

\begin{example}	\lbl{eg:mti-inv}
Let $\Eee=\Set$ and let $T$ be the free monoid-with-involution%
\index{monoid!involution@with involution}
monad on
$\Eee$, as in~\ref{eg:mon-mon-with-inv}.  A $T$-multicategory looks
like a plain multicategory except that the arrows are of the form
\[
a_1^{\sigma_1}, \ldots, a_n^{\sigma_n} \go a
\]
where the $a_i$'s and $a$ are objects and $\sigma_i \in \{-1,+1\}$.  (It is
sometimes convenient to write $x^{-1}$ instead of $x^\circ$ and $x^{+1}$
instead of $x$.)  For instance:
\begin{enumerate}
\item \lbl{eg:mti-inv-Cat}
There is a large $T$-multicategory $\Cat$%
\index{category!generalized multicategory of}
whose objects are all
small categories and in which a map
\[
A_1^{\sigma_1}, \ldots, A_n^{\sigma_n} \go A
\]
is a functor 
\[
A_1^{\sigma_1} \times\cdots\times A_n^{\sigma_n} \go A
\]
where 
\[
A_i^{\sigma_i} = 
\left\{
\begin{array}{ll}
A_i^\op	&\textrm{if } \sigma_i = -1	\\
A_i	&\textrm{if } \sigma_i = +1.	\\
\end{array}
\right.
\]
Composition in $\Cat$ is usual composition of functors, taking opposites
where necessary.  That $\Cat$ does form a $T$-multicategory (and not just a
plain multicategory) is a statement about the behaviour of contravariance
with respect to products and functors.  A `calculus of substitution' of
this kind was envisaged by Kelly%
\index{Kelly, Max}
in the introduction to his~\cite{KelMVFI}
paper. 
\item 	\lbl{eg:mti-inv-gen}
More generally, suppose that $(\cat{A},\otimes,I)$ is a monoidal
category and $\blank^\circ: \cat{A} \go \cat{A}$ a functor for which there
are coherent natural isomorphisms
\[
(A^\circ)^\circ \iso A,
\diagspace
(A \otimes B)^\circ \iso A^\circ \otimes B^\circ,
\diagspace
I^\circ \iso I
\]
($A, B \in \cat{A}$): then we obtain a $T$-multicategory in the same way as
we did for $\Cat$ above.  (Beware that although a typical duality%
\index{duality!operator}
operator
$\blank^\circ$, such as that for duals of finite-dimensional vector spaces,
does satisfy the three displayed isomorphisms, it is a \emph{contravariant}
functor on $\cat{A}$, so does not give a $T$-multicategory.)
\end{enumerate}
\end{example}

\begin{example}	\lbl{eg:mti-anti-inv}%
\index{monoid!anti-involution@with anti-involution}
Now consider the monad $T$ for monoids with \emph{anti}-involution, as
in~\ref{eg:mon-mon-with-anti-inv}.  $T$-multicategories are the
same as in the previous example except that substitution reverses order.
For instance:
\begin{enumerate}
\item Example~\ref{eg:mti-inv}\bref{eg:mti-inv-Cat} is also an example of a
$T$-multicategory for the present $T$, since $A \times B \iso B \times A$
for categories $A$ and $B$.
\item Example~\ref{eg:mti-inv}\bref{eg:mti-inv-gen} can also be repeated,
except that now we require $\blank^\circ$ to reverse, rather than preserve,
the order of the tensored factors:
\[
(A \otimes B)^\circ \iso (B\otimes A)^\circ
\]
naturally in $A, B \in \cat{A}$.  
\item Loop%
\index{loop space!generalized multicategory from}
spaces give an example.  Fix a space $X$ with a basepoint $x$.
In~\ref{eg:mon-cat-loops} we met the monoidal category $\cat{A}$ whose
objects are loops based at $x$ and whose tensor $\otimes$ is concatenation
of loops.  There is a functor $\blank^\circ: \cat{A} \go \cat{A}$ sending a
loop to the same loop run backwards, and this satisfies
\[
(\gamma^\circ)^\circ = \gamma,
\diagspace
(\gamma \otimes \delta)^\circ = \delta^\circ \otimes \gamma^\circ,
\diagspace
\mr{const}_x^\circ = \mr{const}_x
\]
for all loops $\gamma$ and $\delta$.  So there is a resulting
$T$-multicategory whose objects are loops and whose maps encode all the
information about concatenation of loops, homotopy classes of homotopies
between loops, and reversal%
\index{invertibility}
of loops.
\end{enumerate}
\end{example}

\begin{example}	\lbl{eg:mti-exceptions}
Let $\Eee = \Set$ and let $T$ be the monad $1 + \dashbk$%
\index{coproduct!monad from}%
\index{pointed set}
of~\ref{eg:mon-exceptions}.  A $T$-graph is a diagram $1+C_0 \ogby{\dom}
C_1 \goby{\cod} C_0$ of sets and functions.  If we regard $1 + C_0$ as a
subset of the free monoid $\coprod_{n\in\nat} C_0^n$ on $C_0$ and recall
Example~\ref{eg:cl-oms-are-gen} then it is clear that a $T$-multicategory
is exactly a plain multicategory in which all arrows are either unary or
nullary.  So it is natural to draw an arrow $\theta$ of $C$ as either
\[
\begin{centredpic}
\begin{picture}(8,4)(0,-2)
\cell{2}{0}{r}{\tinputlft{a}}
\cell{2}{0}{l}{\tusual{\theta}}
\cell{6}{0}{l}{\toutputrgt{b}}
\end{picture}
\end{centredpic}
\diagspace
\textrm{or}
\diagspace
\begin{centredpic}
\begin{picture}(6,4)(0,-2)
\cell{0}{0}{l}{\tusual{\theta}}
\cell{4}{0}{l}{\toutputrgt{b}}
\end{picture}
\end{centredpic}
\]
where in the first case $\dom(\theta) = a \in C_0$, in the second
$\dom(\theta)$ is the unique element of $1$, and in both $\cod(\theta) =
b$.  The unary arrows form a category $D$, and the nullary arrows define a
functor $Y: D \go \Set$ in which $Y(b)$ is the set of nullary arrows with
codomain $b$.  So a $T$-multicategory is the same thing as a small category
$D$ together with a functor $Y: D \go \Set$, and in particular, a
$T$-operad is a monoid acting%
\index{action!monoid@of monoid}%
\index{monoid!action of}
on a set.  Similarly, if $S$ is any set then
a $(\Set, S + \dashbk)$-multicategory is a small category $D$ together with
an $S$-indexed family of functors $(D \goby{Y_s} \Set)_{s \in S}$.

Another way of putting this is that, when $T = 1 + \dashbk$, a
$T$-multicategory is a discrete opfibration%
\index{fibration!discrete opfibration}
(p.~\pageref{p:defn-cl-d-opfib}).  In fact, $T\hyph\Multicat$ is equivalent
to the category whose objects are discrete opfibrations between small
categories and whose morphisms are commutative squares.
\end{example}

\begin{example}	\lbl{eg:semigp-mti}
As a kind of dual to the last example, let $\Eee=\Set$ and let $T$ be the
free semigroup%
\index{semigroup!free}
monad, $TA = \coprod_{n\geq 1} A^n$.  Then a
$T$-multicategory is exactly a plain multicategory with no nullary%
\index{nullary!arrow}
arrows.  In particular, a $T$-operad $P$ is a family $(P(n))_{n\geq 1}$ of
sets, indexed over positive numbers, equipped with composition and
identities of the usual kind.  Some authors prefer to exclude the
possibility of nullary operations: see the Notes to Chapter~\ref{ch:om}.
\end{example}

\begin{example}	\lbl{eg:M-times-mti}
Fix a monoid%
\index{monoid!action of}%
\index{action!monoid@of monoid}
$M$ and let $T$ be the monad $M\times\dashbk$ on $\Set$, as
in~\ref{eg:mon-action}.  Then a $T$-graph is a diagram $M\times C_0 \og C_1
\go C_0$, and by projecting onto the two factors of $M \times C_0$ we find
that $T\hyph\Multicat$ is isomorphic to the category $\Cat/M$ of categories
over $M$.  (Here we regard a monoid as a one-object category:
p.~\pageref{p:degen-cat-monoid}.)  In particular, $T\hyph\Operad \iso
\fcat{Monoid}/M$.

For a specific example, let $M$ be the (large) monoid of all cardinals%
\index{cardinals}
under multiplication.  Let $C$ be the (large) category of fields%
\index{fields, category of}
and
homomorphisms between them (which are, of course, all injective).  By
taking degrees%
\index{degree of extension}
of extensions we obtain a functor $\pi: C \go M$, making
fields into an $(M\times\dashbk)$-multicategory.
\end{example}

\begin{example}	\lbl{eg:mti-free-cl-opd}%
\index{operad!free}
Let $T$ be the free plain operad monad on $\Eee = \Set^\nat$, as
in~\ref{eg:mon-free-cl-opd}.  In a $T$-multicategory the objects form a
graded set $(C_0(n))_{n\in\nat}$ and the arrows look like%
\index{tree!vertices labelled@with vertices labelled}
\[
\setlength{\unitlength}{1.5em}
\begin{picture}(14.5,3)(0,-1.5)
\cell{0}{0}{l}{%
\begin{picture}(6,3)(-1.5,0)
% 0th layer
\put(1.5,0){\line(0,1){1}}
% 1st layer
\cell{1.5}{1}{c}{\vx}
\put(1.5,1){\line(-3,2){1.5}}
\put(1.5,1){\line(0,1){1}}
\put(1.5,1){\line(3,2){1.5}}
\cell{0}{2}{c}{\vx}
\cell{3}{2}{c}{\vx}
% 2nd layer
\put(0,2){\line(0,1){1}}
\put(3,2){\line(-1,1){1}}
\put(3,2){\line(1,1){1}}
% labels
\cell{1.8}{1}{tl}{a_1}
\cell{-0.3}{2}{tr}{a_2}
\cell{3.3}{2}{tl}{a_3}
\end{picture}}
\put(6.5,0){\vector(1,0){3}}
\cell{8}{0.2}{b}{\theta}
\cell{10.5}{0}{l}{%
\begin{picture}(4,2)(-2,0)
% lower layer
\put(0,0){\line(0,1){1}}
% upper layer
\cell{0}{1}{c}{\vx}
\put(0,1){\line(-2,1){2}}
\put(0,1){\line(-2,3){0.66666667}}
\put(0,1){\line(2,3){0.66666667}}
\put(0,1){\line(2,1){2}}
% label
\cell{0.3}{1}{tl}{a}
\end{picture}}
\end{picture}
\]
($a_1 \in C_0(3)$, $a_2 \in C_0(1)$, $a_3 \in C_0(2)$, $a \in C_0(4)$),
where the tree in the codomain is always the corolla with the same number
of leaves as the tree in the domain.  A typical example of composition is
that arrows
\[
\setlength{\unitlength}{1.5em}
\begin{picture}(17,8)(0,-4.5)
% 
% 1ST INPUT TREE
% 
\cell{0}{-1.5}{tl}{%
\begin{picture}(3,3)(0,0)
% 0th layer
\put(1,0){\line(0,1){1}}
% 1st layer
\cell{1}{1}{c}{\vx}
\put(1,1){\line(-1,1){1}}
\put(1,1){\line(1,1){1}}
% 2nd layer
\cell{2}{2}{c}{\vx}
\put(2,2){\line(-1,1){1}}
\put(2,2){\line(1,1){1}}
% labels
\cell{1.2}{1.2}{tl}{a_1^1}
\cell{2.2}{2.2}{tl}{a_1^2}
\end{picture}}
%
% 2ND INPUT TREE
% 
\cell{0}{0.5}{bl}{%
\begin{picture}(3,3)(-1,0)
% 0th layer
\put(1,0){\line(0,1){1}}
% 1st layer
\cell{1}{1}{c}{\vx}
\put(1,1){\line(-1,1){1}}
\put(1,1){\line(0,1){1}}
\put(1,1){\line(1,1){1}}
% 2nd layer
\cell{0}{2}{c}{\vx}
\put(0,2){\line(0,1){1}}
\cell{1}{2}{c}{\vx}
% labels
\cell{1.2}{1.2}{tl}{a_2^1}
\cell{-0.2}{2}{r}{a_2^2}
\cell{1}{2.2}{b}{a_2^3}
\end{picture}}
% 
% INTERMEDIATE TREE
% 
\cell{6}{0}{l}{%
\begin{picture}(4,3)(0,0)
% 0th layer
\put(1.5,0){\line(0,1){1}}
% 1st layer
\cell{1.5}{1}{c}{\vx}
\put(1.5,1){\line(-3,2){1.5}}
\put(1.5,1){\line(0,1){1}}
\put(1.5,1){\line(3,2){1.5}}
\cell{3}{2}{c}{\vx}
% 2nd layer
\put(3,2){\line(-1,1){1}}
\put(3,2){\line(1,1){1}}
% labels
\cell{1.8}{1}{tl}{a_1}
\cell{3.3}{2}{tl}{a_2}
\end{picture}}
% 
% OUTPUT TREE
% 
\cell{13}{0}{l}{%
\begin{picture}(4,2)(-2,0)
% lower layer
\put(0,0){\line(0,1){1}}
% upper layer
\cell{0}{1}{c}{\vx}
\put(0,1){\line(-2,1){2}}
\put(0,1){\line(-2,3){0.66666667}}
\put(0,1){\line(2,3){0.66666667}}
\put(0,1){\line(2,1){2}}
% label
\cell{0.3}{1}{tl}{a}
\end{picture}}
% 
% ARROWS
%
% theta_1
\qbezier(3.5,-3)(5.4,-3)(7.2,-0.8)
\put(7.2,-0.8){\vector(3,4){0}}
\cell{5.4}{-3}{c}{\theta_1}
% theta_2
\qbezier(3.5,2)(5.5,2.5)(8.6,0.7)
\put(8.6,0.7){\vector(3,-2){0}}
\cell{5.5}{2.5}{c}{\theta_2}
% theta
\put(10.5,0){\vector(1,0){3.5}}
\cell{12}{0.2}{b}{\theta}
\end{picture}
% \hand{60}{12a}
\]
compose to give a single arrow
\[
\begin{array}{c}
\setlength{\unitlength}{1.5em}
\begin{picture}(14,5)(0,-2.5)
% 
% DOMAIN
%
\cell{0}{0}{l}{%
\begin{picture}(4,5)(0,0)
% 0th layer
\put(1,0){\line(0,1){1}}
% 1st layer
\cell{1}{1}{c}{\vx}
\put(1,1){\line(-1,1){1}}
\put(1,1){\line(1,1){1}}
% 2nd layer
\cell{2}{2}{c}{\vx}
\put(2,2){\line(-1,1){1}}
\put(2,2){\line(1,1){1}}
% 3rd layer
\cell{3}{3}{c}{\vx}
\put(3,3){\line(-1,1){1}}
\put(3,3){\line(0,1){1}}
\put(3,3){\line(1,1){1}}
% 4th layer
\cell{2}{4}{c}{\vx}
\put(2,4){\line(0,1){1}}
\cell{3}{4}{c}{\vx}
% labels
\cell{1.2}{1.2}{tl}{a_1^1}
\cell{2.2}{2.2}{tl}{a_1^2}
\cell{3.2}{3.2}{tl}{a_2^1}
\cell{1.8}{4}{r}{a_2^2}
\cell{3}{4.2}{b}{a_2^3}
\end{picture}}
\put(5,0){\vector(1,0){4}}
\cell{7}{0.2}{b}{\theta \of (\theta_1, \theta_2)}
%
% CODOMAIN
%
\cell{10}{0}{l}{%
\begin{picture}(4,2)(-2,0)
% lower layer
\put(0,0){\line(0,1){1}}
% upper layer
\cell{0}{1}{c}{\vx}
\put(0,1){\line(-2,1){2}}
\put(0,1){\line(-2,3){0.66666667}}
\put(0,1){\line(2,3){0.66666667}}
\put(0,1){\line(2,1){2}}
% label
\cell{0.3}{1}{tl}{a}
\end{picture}}
\end{picture}
\end{array}.
% \hand{35}{12b}.
\]
In the case of $T$-operads the labels $a_i$ vanish, so a $T$-operad
consists of a family of sets $(P(\tau))_{\mathrm{trees\ }\tau}$ equipped
with composition and identities.  

We have seen that when $T$ is the identity monad on \Set, a $T$-operad is
exactly a monoid.  We have seen that when $T$ is the free monoid monad, a
$T$-operad is exactly a plain operad.  We have seen that when $T$ is the
free plain operad monad, a $T$-operad is as just described.  This process
can be iterated indefinitely, producing the shapes called `opetopes';%
\index{opetope}
that
is the subject of Chapter~\ref{ch:opetopic}.  The $T$-operads of the
present example are described in~\ref{sec:opetopes} under the name of
`$T_2$-operads'.
\end{example}

\begin{example}		\lbl{eg:tree-mti}
A different example involving trees%
\index{tree!leaves labelled@with leaves labelled}
takes $\Eee$ to be $\Set$ and $T$ to be
the free algebraic theory%
\index{algebraic theory!free}
of~\ref{eg:mon-tree}.  In this context labels
appear on leaves rather than vertices, and trees are amalgamated by grafting
leaves to roots rather than by substituting trees into vertices.  A
$T$-multicategory consists of a set $C_0$ of objects and hom-sets like
\[
C
\left(
\begin{array}{c}
\setlength{\unitlength}{1.5em}
\begin{picture}(4.5,4.6)(-0.25,-0.6)
% 0th layer
\put(2.5,0){\line(0,1){1}}
% 1st layer
\cell{2.5}{1}{c}{\vx}
\put(2.5,1){\line(-3,2){1.5}}
\put(2.5,1){\line(3,2){1.5}}
% 2nd layer
\cell{1}{2}{c}{\vx}
\put(1,2){\line(-1,1){1}}
\put(1,2){\line(1,1){1}}
% 3rd layer
\cell{2}{3}{c}{\vx}
% labels
\cell{0}{3.2}{b}{a_1}
\cell{4}{2.2}{b}{a_2}
\cell{2.5}{-0.2}{t}{a}
\end{picture}
\end{array}
\right)
\]
($a_{1}, a_{2}, a \in C_0$), together with an identity 
\[
1_a 
\in 
C
\left(
\begin{array}{c}
\setlength{\unitlength}{1.5em}
\begin{picture}(1,2)(-0.5,-0.5)
\put(0,0){\line(0,1){1}}
\cell{0}{-0.2}{t}{a}
\cell{0}{1.2}{b}{a}
\end{picture}
\end{array}
\right)
% $\utree$ labelled $a$ at top and bottom
\]
for each $a\in C_0$ and composition functions like
\begin{eqnarray*}
C
\left( 
\begin{array}{c}
\setlength{\unitlength}{1.5em}
\begin{picture}(4.5,4.6)(-0.25,-0.6)
% 0th layer
\put(2.5,0){\line(0,1){1}}
% 1st layer
\cell{2.5}{1}{c}{\vx}
\put(2.5,1){\line(-3,2){1.5}}
\put(2.5,1){\line(3,2){1.5}}
% 2nd layer
\cell{1}{2}{c}{\vx}
\put(1,2){\line(-1,1){1}}
\put(1,2){\line(1,1){1}}
% 3rd layer
\cell{2}{3}{c}{\vx}
% labels
\cell{0}{3.2}{b}{a_1}
\cell{4}{2.2}{b}{a_2}
\cell{2.5}{-0.2}{t}{a}
\end{picture}
\end{array}
\right)
\times
% 
% \left\{	
C
\left( 
\begin{array}{c}
\setlength{\unitlength}{1.5em}
\begin{picture}(2.5,4.8)(-0.25,-0.7)
% 0th layer
\put(1,0){\line(0,1){1}}
% 1st layer
\cell{1}{1}{c}{\vx}
\put(1,1){\line(-1,1){1}}
\put(1,1){\line(1,1){1}}
% 2nd layer
\cell{2}{2}{c}{\vx}
\put(2,2){\line(0,1){1}}
% labels
\cell{0}{2.2}{b}{a_1^1}
\cell{2}{3.2}{b}{a_1^2}
\cell{1}{-0.2}{t}{a_1}
\end{picture}
\end{array}
\right)
\times
C
\left( 
\begin{array}{c}
\setlength{\unitlength}{1.5em}
\begin{picture}(2.5,4.8)(-0.25,-0.7)
% 0th layer
\put(1,0){\line(0,1){1}}
% 1st layer
\cell{1}{1}{c}{\vx}
\put(1,1){\line(0,1){1}}
% 2nd layer
\cell{1}{2}{c}{\vx}
\put(1,2){\line(-1,1){1}}
\put(1,2){\line(0,1){1}}
\put(1,2){\line(1,1){1}}
% 3rd layer (invisible)
\cell{1}{3}{c}{\vx}
% labels
\cell{0}{3.2}{b}{a_2^1}
\cell{2}{3.2}{b}{a_2^2}
\cell{1}{-0.2}{t}{a_2}
\end{picture}
\end{array}
\right)
% \right\}					
\\
\go
C
\left( 
\begin{array}{c}
\setlength{\unitlength}{1.5em}
\begin{picture}(6.5,6.7)(-0.25,-0.6)
% 0th layer
\put(3.5,0){\line(0,1){1}}
% 1st layer
\cell{3.5}{1}{c}{\vx}
\put(3.5,1){\line(-3,2){1.5}}
\put(3.5,1){\line(3,2){1.5}}
% 2nd layer
\cell{2}{2}{c}{\vx}
\put(2,2){\line(-1,1){1}}
\put(2,2){\line(1,1){1}}
\cell{5}{2}{c}{\vx}
\put(5,2){\line(0,1){1}}
% 3rd layer
\cell{1}{3}{c}{\vx}
\put(1,3){\line(-1,1){1}}
\put(1,3){\line(1,1){1}}
\cell{3}{3}{c}{\vx}
\cell{5}{3}{c}{\vx}
\put(5,3){\line(-1,1){1}}
\put(5,3){\line(0,1){1}}
\put(5,3){\line(1,1){1}}
% 4th layer
\cell{2}{4}{c}{\vx}
\put(2,4){\line(0,1){1}}
\cell{5}{4}{c}{\vx}
% labels
\cell{0}{4.2}{b}{a_1^1}
\cell{2}{5.2}{b}{a_1^2}
\cell{4}{4.2}{b}{a_2^1}
\cell{6}{4.2}{b}{a_2^2}
\cell{3.5}{-0.2}{t}{a}
\end{picture}
\end{array}
\right)
\end{eqnarray*}
($a_i^j, a_i, a \in C_0$).  Put another way, a $T$-multicategory $C$ is
a plain multicategory in which the hom-sets are graded by trees: to each $a_1,
\ldots, a_n, a \in C_0$ and $\tau \in \tr(n)$ there is associated a set
$C_\tau(a_1, \ldots, a_n; a)$, composition consists of functions 
\begin{eqnarray*}
C_\tau(a_1, \ldots, a_n; a)
\times
C_{\tau_1}(a_1^1, \ldots, a_1^{k_1}; a_1)
\times\cdots\times
C_{\tau_n}(a_n^1, \ldots, a_n^{k_n}; a_n)\\
\go
C_{\tau\sof(\tau_1, \ldots, \tau_n)}(a_1^1, \ldots, a_n^{k_n}; a),
\end{eqnarray*}
and the identity on $a\in C_0$ is an element of $C_\utree(a;a)$.  In
particular, a $T$-operad is a family $(P(\tau))_{\mathrm{trees\ }\tau}$ of
sets together with compositions
\[
P(\tau) \times P(\tau_1) \times\cdots\times P(\tau_n) 
\go
P(\tau\sof(\tau_1, \ldots, \tau_n))
\]
and an identity element of $P(\utree)$.  
$T$-multicategories are a simplified version of the `relaxed%
\index{multicategory!relaxed}%
\index{relaxed multicategory}
multicategories' mentioned in~\ref{eg:relaxed}.
\end{example}

\begin{example}	\lbl{eg:mti-Top}
Recall that given any symmetric monoidal category $\cat{V}$, there is a
notion of `plain multicategory enriched%
\index{enrichment!plain multicategory@of plain multicategory!symmetric monoidal category@in symmetric monoidal category}
in $\cat{V}$'
(p.~\pageref{p:sym-enr-mti}), the one-object version of which is `plain
operad in $\cat{V}$' (p.~\pageref{p:defn-V-Operad}).  If $\Eee=\Top$ and
$T$ is the free topological%
\index{monoid!topological!free}
monoid monad~(\ref{eg:mon-free-topo-monoid})
then a $T$-operad is precisely an operad in $\Top$.%
\index{operad!topological}
 However,
$T$-multicategories are \emph{not} the same thing as multicategories
enriched in $\Top$, as in a $T$-multicategory $C$ there is a topology on
the set of objects.  A multicategory enriched in $\Top$ is a
$T$-multicategory in which the set of objects has the discrete topology.
This difference should not found be surprising or disappointing: it
exhibits the tension between internal%
\index{internal!enriched@\vs.\ enriched}%
\index{enrichment!internal@\vs.\ internal}
and enriched category theory,
previously discussed in~\ref{sec:cl-enr}.
\end{example}

\begin{example}	\lbl{eg:mti-Cat}
Similarly, if $\Eee=\fcat{Cat}$ and $T$ is the free strict monoidal
category monad~(\ref{eg:mon-free-str-mon-cat}) then a $T$-operad is exactly
a $\Cat$-operad.%
\index{Cat-operad@$\Cat$-operad!generalized operad@as generalized operad}
 We saw some examples of these in~\ref{sec:alg-notions}:
the operads `$F\Sigma$' determining different theories of monoidal
category.  Another example is the $\Cat$-operad $(\TR(n))_{n\in\nat}$,
where $\TR(n)$ is the category of $n$-leafed trees and maps between them,
defined and discussed in~\ref{sec:trees}.
\end{example}

\begin{example}		\lbl{eg:mti-sym}%
\index{multicategory!symmetric vs. generalized@symmetric \vs.\ generalized}
Let $T$ be the free symmetric%
\index{monoidal category!symmetric!free strict}
strict monoidal category monad on $\Cat$, as
in~\ref{eg:mnd-sym-mon}.  Any symmetric multicategory $A$ gives rise to a
$T$-multicategory $C$ as follows.  The category $C_0$ is discrete, with the
same objects as $A$.  The objects of the category $C_1$ are the arrows of
$A$; the arrows of $C_1$ are of the form
\[
(a_{\sigma 1}, \ldots, a_{\sigma n} \goby{\theta\cdot\sigma} a)
\diagspace
\goby{\sigma}
\diagspace
(a_1, \ldots, a_n \goby{\theta} a)
\]
---in other words, an arrow $\phi \go \theta$ is a permutation $\sigma$
such that $\phi = \theta\cdot\sigma$.  The category $TC_0$ has objects all
finite sequences $(a_1, \ldots, a_n)$ of objects of $A$, and arrows of the
form
\[
(a_{\sigma 1}, \ldots, a_{\sigma n})
\diagspace
\goby{\sigma}
\diagspace
(a_1, \ldots, a_n).
\]
The rest of the structure of $C$ is obvious. 

This defines a full and faithful functor
\[
\fcat{SymMulticat} \rIncl T\hyph\Multicat,
\]
so symmetric multicategories could equivalently be defined as
$T$-multicategories with certain properties.  In particular, symmetric
operads are special $T$-operads.  See also the comments on
p.~\pageref{p:enhanced} on `enhanced symmetric multicategories', where the
category $C_0$ is not required to be discrete.
\end{example}

\begin{example}	\lbl{eg:mti-fc}
The free category monad $\fc$ on the category of directed
graphs~(\ref{eg:fc-mnd}) gives rise to a notion of $\fc$-multicategory.%
\index{fc-multicategory@$\fc$-multicategory}
This is the subject of Chapter~\ref{ch:fcm}. 
\end{example}

\begin{example}	\lbl{eg:glob-ops}
Let $T$ be the free strict $\omega$-category monad on the category $\Eee$
of globular sets~(\ref{eg:glob-mnd}).  Then $T$-operads are exactly
globular%
\index{globular operad}
operads, which we
study in depth in Part~\ref{part:n-categories} and which can be used to
specify different theories of $\omega$-category.
In~\ref{sec:non-alg-defns-n-cat} we consider briefly the more general
$T$-multicategories.
\end{example}

\begin{example}%
\index{n-tuple category@$n$-tuple category!strict!free}
The cubical analogue of the previous example takes $T$ to be the free
strict $n$-tuple category monad of~\ref{eg:cubical-mnd}.  A weak%
\index{n-tuple category@$n$-tuple category!weak}
$n$-tuple
category might be defined as an algebra for a certain $T$-operad.  We look
at weak double categories in~\ref{sec:wk-dbl}, but otherwise do not pursue
the cubical case.
\end{example}

\begin{example}	\lbl{eg:multi-alg}
Let $T$ be a cartesian monad on a cartesian category $\Eee$, let
$X\in\Eee$, and let $h: TX \go X$ be a map in $\Eee$.  The $T$-graph
$(\graph{TX}{X}{1}{h})$ can be given the structure of a $T$-multicategory
in at most one way, and such a structure exists if and only if $TX \goby{h}
X$ is an algebra for the monad $T$.  (If it does exist then $\comp=\mu_X$
and $\ids=\eta_X$.)  This defines a full and faithful functor
\[
\blank^{+}: \Eee^T \go T\hyph\Multicat
\]%
\glo{plusalgtomulti}%
\index{plus construction $\blank^+$}%
\index{generalized multicategory!algebra@from algebra}%
turning algebras into multicategories.  In~\ref{sec:non-alg-defns-n-cat} we
will use the idea that a $T$-multicategory is a generalized $T$-algebra to
formulate a notion of weak $n$-category.
\end{example}

In the situation of ordinary categories rather than generalized
multicategories, not only is there the concept of a functor between
categories, but also there are the concepts of a module between categories
(p.~\pageref{p:defn-cat-module}) and a natural transformation between
functors.  The same goes for plain multicategories, as we saw
in~\ref{sec:om-further}.  In fact, both of these concepts make sense for
$T$-multicategories in general.  We meet the definitions and see the
connection between them in~\ref{sec:mmm}.

$\Eee$ and $T$ have so far been regarded as fixed.  But we would
expect some kind of functoriality:%
\lbl{p:change-of-shape-rmks}
if $T'$ is another cartesian monad on another category $\Eee'$ then a `map'
from $(\Eee,T)$ to $(\Eee',T')$ should induce a functor
\[
\Cartpr\hyph\Multicat \go
(\Eee', T')\hyph\Multicat. 
\]
We prove this in~\ref{sec:change}.  That the category
$\Cartpr\hyph\Multicat$ is independent (up to isomorphism) of the choice%
\index{pullback!choice of}
of
pullbacks in $\Eee$ follows from this functoriality by considering the
identity map from ($\Eee$ with one choice of pullbacks) to ($\Eee$ with
another choice of pullbacks).

\section{Algebras}
\lbl{sec:algs}

Theories have models, groups have representations, categories have
set-valued functors, and multicategories have algebras.  Here we meet
algebras for generalized (operads and) multicategories.  There are several
ways of framing the definition.  I have chosen the one that seems most
useful in practice; two alternatives are discussed in~\ref{sec:alg-fibs}
and~\ref{sec:endos}.

Let us begin by considering algebras for a plain multicategory%
\index{multicategory!algebra for}
$C$.  These
are maps from $C$ into the multicategory of sets, but this is not much use
for generalization as there is not necessarily a sensible $T$-multicategory
of sets for arbitrary cartesian $T$.  However, as we saw in
Chapter~\ref{ch:om}, algebras for plain multicategories can be described
without explicit reference to the multicategory of sets.  In the special
case of plain operads $P$, an algebra is a set $X$ together with a map
\[
\coprod_{n\in\nat} P(n) \times X^n \go X
\]
(usually written $(\theta, x_1, \ldots, x_n) \goesto \ovln{\theta}(x_1,
\ldots, x_n)$), satisfying axioms expressing compatibility with composition
and identities in $P$.  Let $T_P$ be the endofunctor $X
\goesto \coprod_{n\in\nat} P(n) \times X^n$ on \Set.  Then the composition
and identity of the operad $P$ induce a monad structure on $T_P$:
the multiplication has components
\begin{eqnarray*}
T_P^2 X		&
\iso	&
\coprod_{n, k_1, \ldots, k_n \in \nat} P(n) \times P(k_1)
\times\cdots\times P(k_n) \times X^{k_1 + \cdots + k_n}	\\
	&
\goby{\coprod \comp \times 1}	&
\coprod_{m\in\nat} P(m) \times X^m
=
T_P X
\end{eqnarray*}
and the unit has components
\[
X \iso
1 \times X^1
\goby{\ids \times 1}
P(1) \times X^1
\rIncl
\coprod_{m\in\nat} P(m) \times X^m
= 
T_P X.
\]
An algebra for the operad $P$ is precisely an algebra for the monad $T_P$.
More generally, an algebra for a plain multicategory $C$ is a family
$(X(a))_{a\in C_0}$ of sets together with a map
\[
\coprod_{n\in\nat, a_1, \ldots, a_n \in C_0} 
C(a_1, \ldots, a_n; a) \times X(a_1) \times \cdots \times X(a_n)
\go 
X(a)
\]
for each $a\in C_0$, satisfying axioms, and again the endofunctor
\[
\label{p:cl-endoftr}
T_C: 
(X(a))_{a\in C_0}
\goesto
\left(
\coprod_{a_1, \ldots, a_n \in C_0} 
C(a_1, \ldots, a_n; a) \times X(a_1) \times \cdots \times X(a_n)
\right)_{a\in C_0}
\]
on $\Set^{C_0}$ naturally has the structure of a monad, an algebra for
which is exactly a $C$-algebra.  Since there is an equivalence of
categories $\Set^{C_0} \eqv \Set/C_0$ (p.~\pageref{eq:Set-slice-power}), we
have succeeded in expressing the definition of an algebra for a plain
operad or multicategory in a completely internal way---that is, completely
in terms of the objects and arrows of the category $\Set$ and the
free-monoid monad on $\Set$.  The strategy for arbitrary cartesian $\Eee$
and $T$ is now clear: an algebra for a $T$-multicategory $C$ should be
defined as an algebra for a certain monad $T_C$ on $\Eee/C_0$.

\begin{figure}
\hspace{10em}\epsfig{file=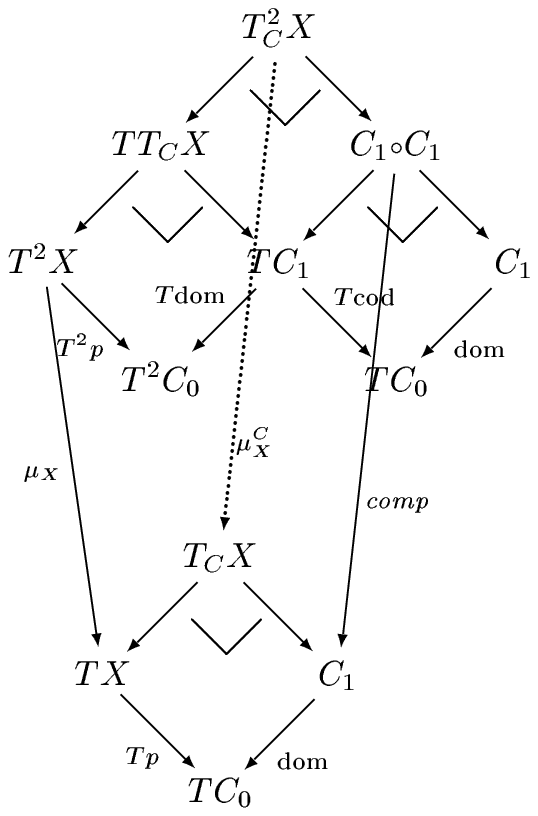}
\caption{Definition of $\mu^C_X$}
\label{fig:mu-diag}
\end{figure}

So, let $(T,\mu,\eta)$ be a cartesian monad on a cartesian category $\Eee$,
and let $C$ be a $T$-multicategory.  If $X = (X \goby{p} C_0)$ is an object
of $\Eee/C_0$ then let $T_C X = (T_C X \goby{p'} C_0)$ be the boxed
composite in the diagram
\[
\setlength{\unitlength}{1em}
\begin{picture}(13,9.5)
\cell{0}{0}{bl}{%
\begin{slopeydiag}
	&	&T_C X\Spbk&	&	&	&	\\
	&\ldTo	&	&\rdTo	&	&	&	\\
TX	&	&	&	&C_1	&	&	\\
	&\rdTo<{Tp}&	&\ldTo>\dom&	&\rdTo>\cod&	\\
	&	&TC_0	&	&	&	&C_0.	\\
\end{slopeydiag}}
\put(10.9,-1){\line(1,1){2}}
\put(2.4,7.5){\line(1,1){2}}
\put(10.9,-1){\line(-1,1){8.5}}
\put(12.9,1){\line(-1,1){8.5}}
\end{picture}
\]
(So the object $T_C X$ of $\Eee$ is defined as a pullback.)  This gives a
functor $T_C: \Eee/C_0 \go \Eee/C_0$.  The unit $\eta^C_X: X \go T_C X$
is the unique map making
\[
\begin{slopeydiag}
	&			&X		&	&	\\
	&\ldTo(2,5)<{\eta_X}	&		&\rdTo>{p}&	\\
	&			&\dGet>{\eta^C_X}&	&C_0	\\
	&			&T_C X\Spbk	&	&\dTo>{\ids}\\
	&\ldTo			&		&\rdTo	&	\\
TX	&			&		&	&C_1	\\
	&\rdTo<{Tp}		&		&\ldTo>{\dom}&	\\
	&			&TC_0		&	&	\\
\end{slopeydiag}
\]
commute.  Similarly, the multiplication $\mu^C_X: T_C^2 X \go T_C X$ is
the unique map making the diagram in Fig.~\ref{fig:mu-diag}
commute.  (That we do have pullbacks as in the top half of the diagram is
an easy consequence of the definition of $T_C$, using the Pasting
Lemma,~\ref{lemma:pasting}.)  This defines a monad $(T_C, \mu^C,
\eta^C)$%
\glo{inducedmonad}
on the category $\Eee/C_0$. 

\begin{defn}%
\index{algebra!generalized multicategory@for generalized multicategory}%
\index{generalized multicategory!algebra for}
Let $T$ be a cartesian monad on a cartesian category $\Eee$, and let $C$ be
a $T$-multicategory.  Then the category $\Alg(C)$%
\glo{Alggenmti}
of \demph{algebras for
$C$} is the category $(\Eee/C_0)^{T_C}$ of algebras for the monad $T_C$ on
$\Eee/C_0$.
\end{defn}

For a more abstract derivation of the induced monad $T_C$, note that there
is a weak functor $\Sp{\Eee}{T} \go \CAT$ sending a 0-cell $E$ to the
category $\Eee/E$%
\index{category!slice}
and defined on 1- and 2-cells by pullback.  Under this
weak functor, any monad in \Sp{\Eee}{T} gives rise to a monad in $\CAT$;
thus, a $T$-multicategory $C$ gives rise to a monad $(T_C, \mu^C, \eta^C)$
on $\Eee/C_0$.

\begin{example}	\lbl{eg:alg-id}
Let $T$ be the identity monad on $\Eee=\Set$ and let $C$ be a
$T$-multicategory, that is, a small category.  Given a set $X\goby{p}C_0$
over $C_0$, write $X(a)$ for the fibre $p^{-1}\{a\}$ over $a\in C_0$.  Then
\begin{eqnarray*}
(T_C X)(a) 	&= 	&\coprod_{b\in C_0} C(b,a) \times X(b), 	\\
(T_C^2 X)(a) 	&= 	&\coprod_{c,b\in C_0} C(c,b) \times C(b,a) \times 
			 X(c).
\end{eqnarray*}
The multiplication map $\mu^C_X: (T_C^2 X)(a) \go (T_C X)(a)$ is given by
composition in $C$, and similarly $\eta^C_X: X(a) \go (T_C X)(a)$ is given
by the identity $1_a$ in $C$.  A $C$-algebra is, therefore, a family
$(X(a))_{a\in C_0}$ of sets together with a function
\begin{equation}	\label{eq:cat-action}
\coprod_{b\in C_0} C(b,a) \times X(b) 
\go
X(a)
\end{equation}
for each $a\in C_0$, compatible with composition and identities in $C$.  So
$\Alg(C) \eqv \ftrcat{C}{\Set}$.

Another way to put this is that for any category $C$, the forgetful functor
$\ftrcat{C}{\Set} \go \ftrcat{C_0}{\Set}$ is monadic and the induced monad
on $\ftrcat{C_0}{\Set} \eqv \Set/C_0$ is $T_C$.  When $C$ is a one-object
category, regarded as a monoid $M$, the resulting monad $T_C$ on $\Set$ is
$M\times\dashbk$, and the category of algebras is the category of left
$M$-sets.%
\index{monoid!action of}%
\index{action!monoid@of monoid}

More generally, let $T$ be the identity monad on any cartesian category
$\Eee$, so that an \Cartpr-multicategory $C$ is a category in $\Eee$: then
$\Alg(C)$ is the category of `left $C$-objects' or `diagrams%
\index{diagram on internal category}
on $C$' as
defined in, for instance, Mac Lane and Moerdijk~\cite[V.7]{MM}.
\end{example}

\begin{example}	\lbl{eg:alg-cl-is-gen} 
Let $T$ be the free monoid monad on the category $\Eee = \Set$ and let $C$
be a $T$-multicategory, that is, a plain multicategory.%
\index{multicategory!algebra for}
 A $C$-algebra 
consists of a family $(X(a))_{a\in C_0}$ of sets together with a function
$(T_C X)(a) \go X(a)$ for each $a \in C_0$, satisfying compatibility
axioms.  By definition of $T_C X$,
\[
(T_C X)(a) =
\coprod_{a_1, \ldots, a_n \in C_0} 
C(a_1, \ldots, a_n; a) \times X(a_1) \times \cdots \times X(a_n),
\]
exactly as on p.~\pageref{p:cl-endoftr}, and we find that the category of
algebras for the $T$-multicategory $C$ is indeed equivalent to the category
of algebras for the plain multicategory $C$.
\end{example}

\begin{example}%
\index{strongly regular theory}
Any strongly regular algebraic theory~(\ref{eg:mon-CJ}) can be described by
a plain operad.  That is, if a monad on $\Set$ arises from a strongly
regular theory then it is isomorphic to the monad $T_P$ arising from some
plain operad $P$; in fact, the converse holds too.  The proofs are
in~\ref{sec:opds-alg-thys}.
\end{example}

\begin{example}	\lbl{eg:alg-exceptions}
When $\Eee=\Set$ and $T= 1+\dashbk$,%
\index{coproduct!monad from}%
\index{pointed set}
as in~\ref{eg:mti-exceptions}, a
$T$-multicategory is an ordinary category $D$ together with a functor
$D\goby{Y}\Set$.  A $(D,Y)$-algebra consists of a functor $D \goby{X} \Set$
together with a natural transformation
\[
D \ctwomult{Y}{X}{} \Set.
\] 
\end{example}

\begin{example}
Let $T$ be the free semigroup%
\index{semigroup!free}
monad on $\Set$, so that a $T$-multicategory
is a plain multicategory with no nullary%
\index{nullary!arrow}
arrows~(\ref{eg:semigp-mti}).
Then an algebra for a $T$-multicategory $C$ is exactly an algebra for the
underlying plain multicategory.  In fact, plain multicategories can be
identified with pairs $(C, X)$ where $C$ is a $T$-multicategory and $X$ is
a $C$-algebra, by making elements of $X(a)$ correspond to nullary arrows
into $a$.
\end{example}

\begin{example}	\lbl{eg:M-times-alg}
Let $M$ be a monoid%
\index{monoid!action of}%
\index{action!monoid@of monoid}
and let $\Cartpr = (\Set, M\times\dashbk)$, so that
a $T$-mul\-ti\-cat\-e\-gory is a category $C$ together with a functor $C
\goby{\phi} M$~(\ref{eg:M-times-mti}).  Then the category of algebras for
$(C, \phi)$ is simply \ftrcat{C}{\Set}, regardless of what $\phi$ is.
This can be seen by working out $T_C$ explicitly, but we will be able to
understand the situation better after we have seen an alternative way of
defining generalized multicategories and their algebras~(\ref{sec:alt-app},
especially the remarks after Corollary~\ref{cor:TC-mti}.)
\end{example}

\begin{example}		\lbl{eg:alg-terminal}
Let $T$ be any cartesian monad on any category $\Eee$ with finite limits.
Then $(T1 \ogby{1} T1 \goby{!} 1)$ is the terminal $T$-graph, and carries a
unique multicategory structure, so is also the terminal%
\index{generalized multicategory!terminal}
$T$-multicategory.
The induced monad on $\Eee/1 \iso \Eee$ is, inevitably, just $(T, \mu,
\eta)$, and so an algebra for the terminal $T$-multicategory is just an
algebra for $T$.
\end{example}

\begin{example}	\lbl{eg:alg-free-cl-opd}%
\index{operad!free}
Let $T$ be the free plain operad monad on the category $\Eee=\Set^\nat$ of
sequences of sets (\ref{eg:mon-free-cl-opd},~\ref{eg:mti-free-cl-opd}).
Then by the previous example, an algebra for the terminal $T$-multicategory
is precisely a plain operad.  Compare Example~\ref{eg:sym-multi-for-opds},
where we defined a \emph{symmetric}%
\index{multicategory!symmetric vs. generalized@symmetric \vs.\ generalized}%
\index{generalized multicategory!operads@for operads}
multicategory whose algebras were plain operads.  This is a minor theme of
this book: objects equipped with symmetries are replaced by objects
equipped with a more refined geometrical structure.
\end{example}

\begin{example}%
\index{algebraic theory!free}
Let $T$ be the free algebraic theory on one operation of each
arity~(\ref{eg:tree-mti}).  Let $P$ be a $T$-operad.  Then a $P$-algebra
structure on a set $X$ consists of a function $X^n \go X$ for each
$n$-leafed tree%
\index{tree!leaves labelled@with leaves labelled}
$\tau$ and each element of $P(\tau)$, satisfying axioms
expressing compatibility with the composition and identity of $P$.
\end{example}

\begin{example}
If $T$ is the free monoid monad on $\Top$%
\index{monoid!topological!free}
or $\Cat$,%
\index{monoidal category!strict!free}
as in~\ref{eg:mti-Top}
or~\ref{eg:mti-Cat}, then an algebra for a $T$-operad $P$ is an algebra in
the usual sense, that is, a space or category $X$ with continuous or
functorial actions%
\index{action!topological monoid@of topological monoid}%
\index{action!monoidal category@of monoidal category}
by the $P(n)$'s.
\end{example}

\begin{example}	\lbl{eg:alg-fc-path}
Historically, one of the most important plain operads has been the little%
\index{operad!little intervals}
intervals operad, that is, the little 1-disks operad
of~\ref{eg:opd-little-disks}, whose algebras are roughly speaking the same
thing as loop spaces.  A plain operad is an operad for the free monoid
monad, and a monoid is a category in which all arrows begin and end at the
same point.  If we are interested in paths%
\index{path!loop@\vs.\ loop}
rather than loops on a basepoint
then it makes sense to replace monoids by arbitrary categories.  Indeed, we
show in~\ref{eg:fc-Trimble} that if $\fc$ is the free category monad on the
category $\Eee$ of directed graphs then there is a certain \fc-operad%
\index{fc-operad@$\fc$-operad}
$P$
such that the paths in any fixed space naturally form a $P$-algebra.  This
solves the language problem posed in~\ref{eg:opd-Trimble}.
\end{example}

\begin{example}
Let $T$ be the free strict $\omega$-category monad on the category $\Eee$
of globular sets, so that a $T$-operad is a `globular%
\index{globular operad}
operad'~(\ref{eg:glob-ops}).  In Chapter~\ref{ch:a-defn} we construct a
certain operad $L$, the initial `operad-with-contraction', and define a
weak $\omega$-category to be an $L$-algebra.  In~\ref{sec:alg-defns-n-cat}
we consider some other possible definitions%
\index{omega-category@$\omega$-category!definitions of}%
\index{n-category@$n$-category!definitions of}
of weak $\omega$-category, some
of which are also of the form `a weak $\omega$-category is a $P$-algebra'
for different choices of globular operad $P$.
\end{example}

\begin{example}		\lbl{eg:alg-to-multi}%
\index{monad!sliced by algebra}%
\index{slice!monad by algebra@of monad by algebra}
Let $T$ be any monad on any category $\Eee$ and let $h= (TX \goby{h} X)$ be
any $T$-algebra.  Then there is a monad $T/h$%
\glo{sliceofmonadbyalg}
on $\Eee/X$ whose functor
part acts on objects by
\[
\left(
\vslob{Y}{p}{X}
\right)
\diagspace 
\goesto 
\diagspace
\left(
\begin{diagram}[height=1.5em]
TY	\\ \dTo>{Tp}	\\ TX	\\ \dTo>{h}	\\ X
\end{diagram}
\right).
\]
Algebras for $T/h$ are just $T$-algebras over $h$; precisely, $\Eee^{T/h}
\iso \Eee^T/h$.

Now recall from~\ref{eg:multi-alg} that when $T$ is a cartesian monad on a
cartesian category $\Eee$, any $T$-algebra $(TX \goby{h} X)$ defines a
$T$-multicategory
\[
h^+ = (TX \ogby{1} TX \goby{h} X).
\]%
\index{plus construction $\blank^+$}%
This induces a monad $T_{h^+}$ on $\Eee/X$.  So starting with $\Eee$, $T$,
and $(TX \goby{h} X)$, we obtain the two monads $T/h$ and $T_{h^+}$ on
$\Eee/X$; they are, inevitably, isomorphic.  So $\Alg(h^+) \iso \Eee^T /
h$.  Example~\ref{eg:alg-terminal} is the special case where $h$ is the
terminal algebra.
\end{example}

\begin{example}
For any $T$-multicategory $C$, the object $(C_1 \goby{\cod} C_0)$ of
$\Eee/C_0$ naturally has the structure of a $T$-algebra.  When
$T$ is the identity monad on $\Set$, so that $C$ is a small category, this
algebra is the functor
\[
\begin{array}{rcl}
C	&\go		&\Set,	\\
a	&\goesto	&\coprod_{a'\in C_0} C(a',a)	
\end{array}
\]
sometimes called the Cayley%
\index{Cayley representation}%
\index{representation theorem}
representation of $C$.  
\end{example}

In this section we have seen how to associate to each $T$-multicategory $C$
a category $\Alg(C)$.  We would expect some kind of functoriality in $C$.
When $T$ is the identity monad on $\Set$, a functor $C \go C'$ between
($T$-multi)categories induces a functor in the opposite direction,
\[
\Alg(C') \eqv \ftrcat{C'}{\Set} \go \ftrcat{C}{\Set} \eqv \Alg(C).
\]
The same holds when $T$ is the free monoid monad on $\Set$, viewing
$C$-algebras as multicategory maps $C\go \Set$ as in the original
definition,~\ref{defn:alg-multi}.

In fact, the construction works for an arbitrary cartesian monad $T$: any
map $f: C \go C'$ of $T$-multicategories induces a functor $\Alg(C') \go
\Alg(C)$.%
\lbl{p:Alg-functorial}
First, we have the functor $f_0^*: \Eee/C'_0 \go \Eee/C_0$ defined by
pullback along $f_0: C_0 \go C'_0$.  Then there is a naturally-arising
natural transformation
\[
\begin{diagram}
\Eee/C'_0		&\rTo^{T_{C'}}	&\Eee/C'_0	\\
\dTo<{f_0^*}		&\nent \phi		&\dTo>{f_0^*}		\\
\Eee/C_0		&\rTo_{T_C}		&\Eee/C_0,		\\
\end{diagram}
\]
which the reader will easily be able to determine.  This $\phi$ is
compatible with the monad structures on $T_{C'}$ and $T_C$: in the
terminology defined in~\ref{sec:more-monads}, $(f_0^*, \phi)$ is a `lax map
of monads' from $T_{C'}$ to $T_C$.  It follows that there is an induced
functor%
\index{algebra!generalized multicategory@for generalized multicategory!change of shape (induced functor)}
on the categories of algebras for these monads, that is, from
$\Alg(C')$ to $\Alg(C)$.

This construction defines a map
\[
\Alg: (T\hyph\Multicat)^\op \go \CAT.
\]
Since the induced functors are defined by pullback, it is inevitable that
this map does not preserve composites and identities strictly, but only up
to coherent isomorphism.  Precisely, it is a weak functor from the
2-category $(T\hyph\Multicat)^\op$ whose only 2-cells are identities to the
2-category $\CAT$.  If we also bring into play transformations between
$T$-multicategories~(\ref{sec:mmm}), then $T\hyph\Multicat$ becomes a
strict 2-category and $\Alg$ a weak functor between 2-categories.

\begin{notes}

$T$-multicategories were introduced by Burroni%
\index{Burroni, Albert}
in his~\cite{Bur} paper,
where they went by the name of $T$-categories.  He showed how to define
them for \emph{any} monad $T$, though he concentrated on cartesian $T$ (and
it is not clear that $T$-categories are useful outside this case).  The
basic idea has been independently rediscovered on at least two occasions:
by Hermida~\cite{HerRM}%
\index{Hermida, Claudio}
and by Leinster~\cite{GOM}.  The notion of an
algebra for a $T$-multicategory seems not to have appeared before the
latter paper.

The shape of this chapter is typical of much of this text: while the
formalism is quite simple (in this case, the definition of
$T$-multicategory and of algebra), it can take a long time to see what it
means concretely in particular cases of interest.  Indeed, in this chapter
we have restricted ourselves to the simpler instances of $T$, leaving some
of the more advanced examples to chapters where they can be explored at
greater leisure (Ch.~\ref{ch:fcm}, \ref{ch:opetopic}, \ref{ch:globular}).

Hermida called $\Sp{\Eee}{T}$ the `Kleisli%
\index{Kleisli!bicategory of spans}
bicategory of spans'
in~\cite{HerRM}; the formal similarity between the definition of
$\Sp{\Eee}{T}$ and the usual construction of a Kleisli category is evident.

Dmitry Roytenberg%
\index{Roytenberg, Dmitry}
suggested to me that something like
Example~\ref{eg:alg-fc-path} ought to exist.

\end{notes}

\chapter{Example: $\mathbf{fc}$-Multicategories}
\lbl{ch:fcm}

\chapterquote{%
A lot of people are afraid of heights.  Not me.  I'm afraid of widths.}{%
Steven Wright}

\noindent
The generalized multicategories that we are interested in typically have some
geometry to them.  They are often `higher-dimensional' in some sense.  In
this chapter we study a 2-dimensional example, \fc-multicategories, which
are $T$-multicategories when $T$ is the free category monad $\fc$ on the
category of directed graphs.

This case is interesting for a variety of reasons.  First,
\fc-multicategories turn out to encompass a wide range of familiar
2-dimensional structures, including bicategories, double categories,
monoidal categories and plain multicategories.  Second, there are two
well-known ideas for which \fc-multicategories provide a cleaner and more
general context than is traditional: the `bimodules construction' (usually
done on bicategories) and the enrichment of categories (usually done in
monoidal categories).  Third, these 2-dimensional structures are the second
rung on an infinite ladder of higher-dimensional structures (the first rung
being ordinary categories), and give us clues about the behaviour of the
more difficult, less easily visualized higher rungs.

We start~(\ref{sec:fcm}) by unwinding the definition of \fc-multicategory
to give a completely elementary description.  As mentioned above, various
familiar structures arise as special kinds of \fc-multicategory; we show
how this happens.  

A nearly-familiar structure that arises as a special kind of
\fc-multicategory is the `weak double category', in which horizontal
composition only obeys associativity and unit laws up to coherent
isomorphism.  We define these in~\ref{sec:wk-dbl} and give examples, of
which there are many natural and pleasing ones.  

Section~\ref{sec:mmm} is on the `bimodules%
\index{bimodules construction}
construction' or, as we prefer
to call it, the `monads%
\index{monads construction}
construction'.  This takes an \fc-multicategory $C$
as input and produces a new \fc-multicategory $\Mon(C)$ as output.  For
example, if $C$ is the \fc-multicategory of sets (and functions, and spans)
then $\Mon(C)$ is the \fc-multicategory of categories (and functors, and
modules); if $C$ is abelian groups (and homomorphisms) then $\Mon(C)$ is
rings (and homomorphisms, and modules).  We wait until~\ref{sec:enr-mtis}
to see what a `category enriched in an \fc-multicategory' is.

\section{$\mathbf{fc}$-multicategories}
\lbl{sec:fcm}

Let $\scat{H}$ be the category $(0 \parpair{\sigma}{\tau} 1)$ and $\Eee
= \ftrcat{\scat{H}^\op}{\Set}$.  Then $\Eee$ is the category of directed
graphs, and there is a forgetful functor $U: \Cat \go \Eee$.  This has a
left adjoint $F$: if $E \in \Eee$ then objects of $FE$ are vertices of $E$,
arrows in $FE$ are strings 
\[
x_0 \goby{p_1} x_1 \goby{p_2} \ \cdots\ \goby{p_n} x_n
\]
of edges in $E$ (with $n\geq 0$), and composition is concatenation.  The
adjunction induces a monad $T=U\of F$ on $\Eee$, which can be shown to be
cartesian either by direct calculation or by applying some general
theory~(\ref{eg:fc-cart}).  We write $T = \fc$,%
\glo{fcprecise}%
\index{category!free (fc)@free ($\fc$)}
for `free category', and so
we have the notion of an \fc-multicategory.

What is an \fc-multicategory,%
\index{fc-multicategory@$\fc$-multicategory}
explicitly?  An \fc-graph%
\index{fc-graph@$\fc$-graph}
$C$ is a diagram
\[
\begin{diagram}[height=2em,noPS]
	&	&C_1=(C_{11}&\pile{\rTo\\ \rTo}	&C_{10})&	&	\\
	&	&\ldTo	&			&\rdTo	&	&	\\
\fc(C_0)=(C'_{01}&\pile{\rTo\\ \rTo}&C_{00})&		&
C_0=(C_{01} &\pile{\rTo\\ \rTo}&C_{00})\\
\end{diagram}
\]
where $C_1$ and $C_0$ are directed graphs and the diagonal arrows are maps
of graphs, the $C_{ij}$'s are sets and the horizontal arrows are functions,
and $C'_{01}$ is the set of strings of edges in $C_0$.  Think of elements
of $C_{00}$ as \demph{objects} or \demph{0-cells},%
\index{cell!fc-multicategory@of $\fc$-multicategory}
elements of $C_{01}$ as
\demph{horizontal 1-cells}, elements of $C_{10}$ as \demph{vertical
1-cells}, and elements of $C_{11}$ as \demph{2-cells}, as in the picture
\begin{equation}	\label{diag:fcm-two-cell}
\begin{fcdiagram}
a_0	&\rTo^{m_1}	&a_1	&\rTo^{m_2}	&\ 	&\cdots	
&\ 	&\rTo^{m_n}	&a_n	\\
\dTo<{f}&		&	&		&\Downarrow\tcs{\theta}&
&	&		&\dTo>{f'}\\
a	&		&	&		&\rTo_{m}	&	
&	&		&a'	\\
\end{fcdiagram}
\end{equation}
($n\geq 0$, $a_i, a, a' \in C_{00}$, $m_i, m \in C_{01}$, $f, f' \in C_{10}$,
$\theta\in C_{11}$).  An \fc-multicategory
structure on the \fc-graph $C$
amounts to composition and identities of two types.  First, the directed
graph 
\[
\begin{slopeydiag}
	&	&C_{10}	&	&	\\
	&\ldTo	&	&\rdTo	&	\\
C_{00}	&	&	&	&C_{00}	\\
\end{slopeydiag}
\]
has the structure of a category; in other words, vertical 1-cells can be
composed and there is an identity vertical 1-cell on each 0-cell.  Second,
there is a composition function for 2-cells,
\begin{equation}	\label{eq:pasted-two-cells}
\begin{array}{l}
\begin{diagram}[width=.5em,height=2em]
\blob&\rTo^{m_1^1}&\cdots&\rTo^{m_1^{k_1}}&
\blob&\rTo^{m_2^1}&\cdots&\rTo^{m_2^{k_2}}&\blob&
\ &\cdots&\ &
\blob&\rTo^{m_n^1}&\cdots&\rTo^{m_n^{k_n}}&\blob\\
\dTo<{f_0}&&\Downarrow\tcs{\theta_1}&&
\dTo<{f_1}&&\Downarrow\tcs{\theta_2}&&\dTo&
\ &\cdots&\ &
\dTo&&\Downarrow\tcs{\theta_n}&&\dTo>{f_n}\\
\blob&&\rTo_{m_1}&&
\blob&&\rTo_{m_2}&&\blob&
\ &\cdots&\ &
\blob&&\rTo_{m_n}&&\blob\\
\dTo<{f}&&&&&&&\Downarrow\tcs{\theta} &&&&&&&&&\dTo>{f'}\\
\blob&&&&&&&&\rTo_{m}&&&&&&&&\blob\\
\end{diagram}
\\
\\
\goesto\diagspace
\begin{diagram}[width=.5em,height=2em]
\blob&\rTo^{m_1^1}&\ &&
&&&&\cdots&
&&&
&&\ &\rTo^{m_n^{k_n}}&\blob\\
\dTo<{f\of f_0}&&&&&&&&
\Downarrow\tcs{\theta\of\tuple{\theta_1}{\theta_2}{\theta_n}}
&&&&&&&&\dTo>{f'\of f_n}\\
\blob&&&&&&&&\rTo_{m}&&&&&&&&\blob\\
\end{diagram}
\end{array}
\end{equation}
($n\geq 0, k_i\geq 0$, with \blob's representing objects), and a function
assigning an identity 2-cell to each horizontal 1-cell,
\[
\begin{fcdiagram}
a&\rTo^{m}&a'\\
\end{fcdiagram}
\diagspace\goesto\diagspace
\begin{fcdiagram}
a&\rTo^{m}&a'\\
\dTo<{1_a}&\Downarrow \tcs{1_m}&\dTo>{1_{a'}}\\
a&\rTo_{m}&a'.\\
\end{fcdiagram}
\]
The composition and identities obey associativity and identity laws, which
ensure that any diagram of pasted-together 2-cells with a rectangular
boundary has a well-defined composite.

The pictures in the nullary%
\index{nullary!composite}%
\index{fc-multicategory@$\fc$-multicategory!nullary composite in}
case are worth a short comment.
When $n=0$, the 2-cell of diagram~\bref{diag:fcm-two-cell} is drawn as
\[
\begin{fcdiagram}
a_0		&\rEquals		&a_0		\\
\dTo<{f}	&\Downarrow\tcs{\theta}	&\dTo>{f'}	\\
a		&\rTo_{m}		&a',		\\
\end{fcdiagram}
\]
and the diagram of pasted-together 2-cells in the domain
of~\bref{eq:pasted-two-cells} is drawn as
\[
\begin{fcdiagram}
b_0		&\rEquals		&b_0		\\
\dTo<{f_0}	&=			&\dTo>{f_0}	\\
a_0		&\rEquals		&a_0		\\
\dTo<{f}	&\Downarrow\tcs{\theta}	&\dTo>{f'}	\\
a		&\rTo_{m}		&a'.		\\
\end{fcdiagram}
\]
The composite%
\lbl{p:null-notation}
of this last diagram will be written as $\theta\of f_0$.

So in completely elementary terms, an \fc-multicategory consists of
\begin{itemize}
\item a set of objects
\item for each pair $(a,a')$ of objects, a set of vertical 1-cells 
$
\begin{fcdiagram}
a	\\
\dTo	\\
a'	\\
\end{fcdiagram}
$
\item for each pair $(a,a')$ of objects, a set of horizontal 1-cells
$a \go a'$
\item for each $a_0, \ldots, a_n, a, a', m_1, \ldots, m_n, m, f, f'$ as
in~\bref{diag:fcm-two-cell}, a set of 2-cells $\theta$
\item composition and identity functions for vertical 1-cells, as described
above
\item composition and identity functions for 2-cells, as described above,
\end{itemize}
satisfying associativity and identity axioms.  Having given this elementary
description I will feel free to refer to large \fc-multicategories, in
which the collection of cells is a proper class.

\begin{example}	\lbl{eg:fcm-Ring}
There is an \fc-multicategory \fcat{Ring}%
\glo{fcRing}%
\index{fc-multicategory@$\fc$-multicategory!rings@of rings}
in which
\begin{itemize}
\item 0-cells are (not necessarily commutative) rings
\item vertical 1-cells are ring homomorphisms
\item a horizontal 1-cell $A \go A'$ is an
  $(A',A)$-bimodule~(\ref{eg:bicat-Ring})%
\index{module!bimodule over rings}

\item a 2-cell  
\[
\begin{fcdiagram}
A_0	&\rTo^{M_1}	&A_1	&\rTo^{M_2}	&\ 	&\cdots	
&\ 	&\rTo^{M_n}	&A_n	\\
\dTo<{f}&		&	&		&\Downarrow\tcs{\theta}&
&	&		&\dTo>{f'}\\
A	&		&	&		&\rTo_{M}	&	
&	&		&A'	\\
\end{fcdiagram}
\]
is an abelian group homomorphism
\[
\theta: 
M_n \otimes_{A_{n-1}} M_{n-1} \otimes_{A_{n-2}} \cdots
\otimes_{A_1} M_1
\go 
M
\]
satisfying
\[
\theta(a_n\cdot m_n \otimes m_{n-1} \otimes\cdots\otimes
m_2 \otimes m_1\cdot a_1)
=
f'(a_n) \cdot 
\theta(m_n \otimes\cdots\otimes m_1)
\cdot f(a_1)
\]
(that is, $\theta$ is a homomorphism of $(A_n, A_0)$-bimodules if $M$ is
given an $(A_n, A_0)$-bimodule structure via $f$ and $f'$)
\item composition and identities are defined in the evident way.
\end{itemize}
Thus, rings, homomorphisms of rings, modules over rings, homomorphisms of
modules, and tensor products of modules are integrated into a single
structure.  The category formed by the objects and vertical 1-cells is the
ordinary category of rings and homomorphisms.  When the distinction needs
making, I will write $\fcat{Ring}_1$ for this `1-dimensional' category and
$\fcat{Ring}_2$%
\glo{fcRing2}
for the `2-dimensional' \fc-multicategory.  Similar
notation is extended to similar examples.
\end{example}

Many familiar 2-dimensional%
\index{algebraic theory!two-dimensional@2-dimensional}
structures are degenerate \fc-multicategories.
The following examples demonstrate this; Fig.~\ref{fig:degens} is a
summary.
\begin{figure}\protect\small
\begin{tabular}{llll}	
			&\emph{No horizontal}	&\emph{Weak horizontal}	&
\emph{Strict horizontal}\\
			&\emph{composition}	&\emph{composition}	&
\emph{composition}	\\
\\
\emph{No degeneracy}	&\fc-multicategory	&Weak double		&
Strict double		\\
			&			&category		&
category		\\
\\
\emph{All vertical 1-cells}&Vertically discrete	&Bicategory		&
Strict			\\
\emph{are identities}	&\fc-multicategory	&			&
2-category		\\
\\
\emph{Only one object and}&Plain		&Monoidal		&
Strict monoidal \\
\emph{one vertical 1-cell}&multicategory	&category		&
category
\end{tabular}%
\index{fc-multicategory@$\fc$-multicategory!degenerate}
\caption{Some of the possible degeneracies of an \fc-multicategory. The
left-hand column refers to degeneracies in the category formed by the
objects and vertical 1-cells. The top row says whether the
\fc-multicategory structure arises from a composition rule for horizontal
1-cells. See Examples~\ref{eg:fcm-str-dbl}--\ref{eg:fcm-mon-cat}}
\label{fig:degens}
\end{figure}

\begin{example}	\lbl{eg:fcm-str-dbl}%
\index{double category!strict!fc-multicategory@as $\fc$-multicategory}%
\index{fc-multicategory@$\fc$-multicategory!strict double category@from strict double category}
Any strict double category gives rise to an \fc-multicategory, in which a
2-cell as in~\bref{diag:fcm-two-cell} is a 2-cell
\begin{diagram}[height=2em]
a_0	&\rTo^{m_n \of \cdots \of m_1}	&a_n\\
\dTo<{f}&\Downarrow\tcs{\theta}		&\dTo>{f'}\\
a	&\rTo_{m}			&a'\\
\end{diagram}
in the double category.
\end{example}

\begin{example}	\lbl{eg:fcm-wk-dbl}
The last example works just as well if we start with a double category in
which the composition of horizontal 1-cells only obeys weak laws.  We call
these `weak%
\index{fc-multicategory@$\fc$-multicategory!weak double category@from weak double category}%
\index{double category!weak!fc-multicategory@as $\fc$-multicategory}
double categories' and discuss them in detail in the next
section; a typical example is \fcat{Ring}~(\ref{eg:fcm-Ring}).  A similar
example, and really the archetypal weak double category, is $\Set_2$,%
\glo{fcSet2}%
\index{fc-multicategory@$\fc$-multicategory!sets@of sets}
defined as follows:
\begin{itemize}
\item objects are sets
\item vertical 1-cells are functions (and the category formed by objects
and vertical 1-cells is the ordinary category $\Set_1$ of sets and
functions; see~\ref{eg:fcm-Ring} for the notation)
\item horizontal 1-cells are spans:%
\index{span}
that is, a horizontal 1-cell $A \go A'$
is a diagram
\[
\begin{slopeydiag}
	&	&M	&	&	\\
	&\ldTo	&	&\rdTo	&	\\
A	&	&	&	&A'	\\
\end{slopeydiag}
\]
of sets and functions
\item a 2-cell inside
\[
\begin{slopeydiag}
	&	&M	&	&	\\
	&\ldTo	&	&\rdTo	&	\\
A	&	&	&	&A'	\\
\dTo<f	&	&N	&	&\dTo>{f'}\\
	&\ldTo	&	&\rdTo	&	\\
B	&	&	&	&B'	\\
\end{slopeydiag}
\]
is a function $\theta: M \go N$ making the diagram commute
\item horizontal composition is by pullback.
\end{itemize}
That every weak double category has an underlying \fc-multicategory is
exactly analogous to every weak monoidal category having an underlying
plain multicategory.
\end{example}

\begin{example}	\lbl{eg:fcm-vdisc}
Consider $\fc$-multicategories $C$ in which all vertical 1-cells are
identities.  This means that the category formed by the objects and
vertical 1-cells is discrete, so we call $C$ \demph{vertically discrete}.%
\index{fc-multicategory@$\fc$-multicategory!vertically discrete}
A vertically discrete \fc-multicategory consists of some objects $a, a',
\ldots$, some 1-cells $m, m',\ldots$, and some 2-cells looking like
\begin{equation}	\label{diag:fc-ope-two-cell}
\begin{diagram}[size=1.5em,tight,scriptlabels]
		&	&	&a_2	&\ldots	&	&	&	\\
		&a_1	&\ruEdge(2,1)^{m_2}&&	&	&a_{n-1}&	\\
\ruEdge(1,2)<{m_1}&	&	&	&\Downarrow \tcs{\theta}&&
&\rdEdge(1,2)>{m_{n}}							\\
a_0		&	&	&\rEdge_{m}&	&	&	&a_n,	\\
\end{diagram}
\end{equation}
together with a composition function
\[
\begin{array}{c}
\setlength{\unitlength}{1mm}
\begin{picture}(49,25)(0,-2)
\cell{0}{0}{bl}{\epsfig{file=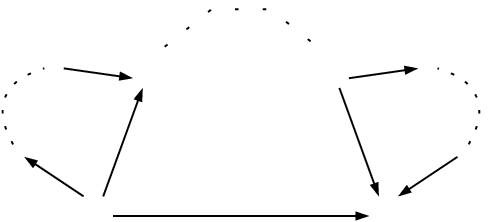}}
% 0-cell labels
\cell{10}{1.5}{t}{a_0}
\cell{15.5}{15.5}{t}{a_1}
\cell{35}{15.5}{t}{a_{n-1}}
\cell{40}{1.5}{t}{a_n}
% 1-cell labels
\cell{25}{0}{t}{\scriptstyle m}
\cell{15}{8}{c}{\scriptstyle m_1}
\cell{34.5}{8}{c}{\scriptstyle m_n}
\cell{5}{2}{c}{\scriptstyle m_1^1}
\cell{10}{17.5}{c}{\scriptstyle m_1^{k_1}}
\cell{39}{17.5}{c}{\scriptstyle m_n^1}
\cell{46.5}{4}{c}{\scriptstyle m_n^{k_n}}
% 2-cell labels
\cell{25}{7}{c}{\Downarrow\theta}
\da{8}{8}{68}
\cell{9}{11}{c}{\theta_1}
\da{42}{8}{-68}
\cell{41}{11}{c}{\theta_n}
\end{picture}
\end{array}
\ \goesto \ 
\begin{array}{c}
\setlength{\unitlength}{1mm}
\begin{picture}(49,25)(0,-2)
\cell{0}{0}{bl}{\epsfig{file=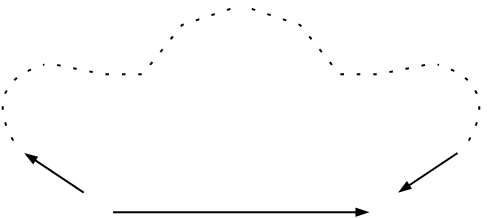}}
% 0-cell labels
\cell{10}{1.5}{t}{a_0}
\cell{40}{1.5}{t}{a_n}
% 1-cell labels
\cell{25}{0}{t}{\scriptstyle m}
\cell{5}{2}{c}{\scriptstyle m_1^1}
\cell{46.5}{4}{c}{\scriptstyle m_n^{k_n}}
% 2-cell labels
\cell{25}{7}{c}{\Downarrow \theta \of (\theta_1, \ldots, \theta_n)}
\end{picture}
\end{array}
\]
and an identity function
\[
\begin{array}{c}
\setlength{\unitlength}{1mm}
\begin{picture}(33,5)
\cell{3}{2}{l}{\epsfig{file=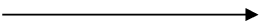}}
\cell{2.5}{1.3}{br}{a}
\cell{30}{1.3}{bl}{a'}
\cell{16}{3}{b}{\scriptstyle m}
\end{picture}
\end{array}
\diagspace\goesto\diagspace
\begin{array}{c}
\setlength{\unitlength}{1mm}
\begin{picture}(33,13)(0,1)
\cell{3}{3}{bl}{\epsfig{file=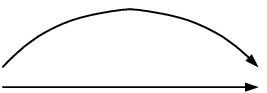}}
\cell{2.5}{3.5}{br}{a}
\cell{30}{3.5}{bl}{a'}
\cell{16}{2.5}{t}{\scriptstyle m}
\cell{16}{12.5}{b}{\scriptstyle m}
\cell{16}{7}{c}{\Downarrow 1_m}
\end{picture}
\end{array}
\]
obeying the inevitable associativity and identity laws.  This leads us into
the world of opetopes%
\index{opetope}
(Chapter~\ref{ch:opetopic}).
\end{example}

\begin{example}	\lbl{eg:fcm-bicat}
A bicategory%
\index{bicategory!fc-multicategory@as $\fc$-multicategory}%
\index{fc-multicategory@$\fc$-multicategory!bicategory@from bicategory}
is a weak double category in which the only vertical 1-cells
are identities, and so gives rise to an \fc-multicategory.  Explicitly, if
$B$ is a bicategory then there is a vertically discrete $\fc$-multicategory
whose objects are those of $B$, whose horizontal 1-cells are the 1-cells of
$B$, and whose 2-cells~\bref{diag:fc-ope-two-cell} are 2-cells
\[
a_0 
\ctwomult{\scriptstyle (m_n \sof\cdots\sof m_1)}%
{\scriptstyle m}{\tcs{\theta}} 
a_n
\]
in $B$.
\end{example}

\begin{example}	\lbl{eg:fcm-cl-mti}%
\index{multicategory!fc-multicategory@as $\fc$-multicategory}%
\index{fc-multicategory@$\fc$-multicategory!plain multicategory@from plain multicategory}
Plain multicategories are the same as \fc-multicategories with only
one object and one vertical 1-cell.  If the plain multicategory is
called $M$ then we call the corresponding \fc-multicategory $\Sigma M$,%
\glo{Sigmamtifc}
the \demph{suspension}%
\index{suspension!plain multicategory@of plain multicategory}%
\index{multicategory!suspension of}
of $M$.  Horizontal 1-cells in $\Sigma M$ are objects
of $M$, and 2-cells
\begin{equation}	\label{diag:fcm-cl-mti-cell}
\begin{diagram}[size=1.5em,tight,abut]
		&	&	&\	&\ldots	&	&	&	\\
		&\bullet&\ruTo(2,1)^{m_2}&&	&	&\	&	\\
\ruTo(1,2)<{m_1}&	&	&	&\Downarrow \tcs{\theta}&&	
&\rdTo(1,2)>{m_n}							\\
\bullet		&	&	&\rTo_{m}&	&	&	&\bullet\\
\end{diagram}
\end{equation}
in $\Sigma M$ are maps
\[
m_1, \ldots, m_n \goby{\theta} m
\]
in $M$.
\end{example}

\begin{example}	\lbl{eg:fcm-cl-opd}
In particular, plain operads are \fc-multicategories in which there is only
one object, one vertical 1-cell and one horizontal 1-cell.  An \fc-operad%
\index{fc-operad@$\fc$-operad}
is an \fc-multicategory in which there is only one object and one
horizontal 1-cell, so a plain operad is a special kind of \fc-operad.

If we are going to take the suspension%
\index{suspension!plain operad@of plain operad}%
\index{operad!suspension of}
idea seriously then we should write
\[
\Sigma: \Operad \rIncl \Multicat
\]
for the inclusion of operads as one-object multicategories.  We also
have~(\ref{eg:fcm-cl-mti}) the suspension map
\[
\Sigma: \Multicat \rIncl \fc\hyph\Multicat,
\]
so the inclusion of operads into \fc-multicategories should be written
$\Sigma^2$: double suspension.  Of course, we usually leave the first
inclusion nameless.
\end{example}

\begin{example}	\lbl{eg:fcm-mon-cat}
As a special case of both~\ref{eg:fcm-bicat} and~\ref{eg:fcm-cl-mti}, any
monoidal%
\index{monoidal category!fc-multicategory@as $\fc$-multicategory}%
\index{fc-multicategory@$\fc$-multicategory!monoidal category@from monoidal category}
category $M$ gives rise to an \fc-multicategory $\Sigma M$ in
which there is one object and one vertical 1-cell.  Horizontal 1-cells are
objects of $M$, and 2-cells~\bref{diag:fcm-cl-mti-cell} are maps 
\[
\theta: (m_n \otimes\cdots\otimes m_1) \go m
\]
in $M$.  
\end{example}

\begin{example}	\lbl{eg:fcm-matrices}
Here is a family of \fc-multicategories that are not usually degenerate%
\index{fc-multicategory@$\fc$-multicategory!non-degenerate}
in
any of the ways listed above.  Let $V$ be a plain multicategory and define
the \fc-multicategory $\Set[V]$ as follows:
\begin{itemize}
\item objects are sets
\item vertical 1-cells are functions
\item a horizontal 1-cell $A \go A'$ is a family $(m_{a,a'})_{a\in A, a'\in
A'}$ of objects of $V$
\item a 2-cell
\[
\begin{fcdiagram}
A_0	&\rTo^{m^1}	&A_1	&\rTo^{m^2}	&\ 	&\cdots	
&\ 	&\rTo^{m^n}	&A_n	\\
\dTo<{f}&		&	&		&\Downarrow\tcs{\theta}&
&	&		&\dTo>{f'}\\
A	&		&	&		&\rTo_{m}	&	
&	&		&A'	\\
\end{fcdiagram}
\]
is a family
\[
\left(
m^1_{a_0, a_1}, \ldots, m^n_{a_{n-1}, a_n}
\goby{\theta_{a_0, \ldots, a_n}}
m_{f(a_0), f'(a_n)}
\right)_{a_0 \in A_0, \ldots, a_n \in A_n}
\]
of maps in $V$
\item composition and identities are obvious. 
\end{itemize}
If $V=\Set$ then $\Set[V]$ is the \fc-multicategory $\Set_2$
of~\ref{eg:fcm-wk-dbl}.  But $\Set[V]$ is not (the underlying
\fc-multicategory of) a weak double category unless $V$ happens to be (the
underlying plain multicategory of) a monoidal category, since $\Set[V]$
does not usually have `enough universal 2-cells'.  Compare the results on
representable%
\index{fc-multicategory@$\fc$-multicategory!representable}
multicategories in~\ref{sec:non-alg-notions}.
\end{example}

We will not be much concerned with algebras%
\index{algebra!for fc-multicategory@for $\fc$-multicategory}%
\index{fc-multicategory@$\fc$-multicategory!algebra for}
for \fc-multicategories, but
let us look at them briefly.  Given an object $E$ of $\Eee$, that is, a
directed graph, an object $(X\go E)$ of $\Eee/E$ amounts to a set $X(a)$
for each vertex $a$ of $E$ and a span 
\[
\begin{slopeydiag}
	&	&X(m)	&	&	\\
	&\ldTo	&	&\rdTo	&	\\
X(a)	&	&	&	&X(b)	\\
\end{slopeydiag}
\]
for each edge $a \goby{m} b$ of $E$.  An algebra for an \fc-multicategory
$C$ is an object $(X \go C_0)$ of $\Eee/C_0$ together with an action of $C$
on $X$, and so a $C$-algebra is just an \fc-multicategory map from $C$ into
the (large) \fc-multicategory $\Set_2$ defined in~\ref{eg:fcm-wk-dbl}.

There are some interesting examples of algebras when $C$ is the
\fc-multicategory coming from a plain operad $P$: in the notation
of~\ref{eg:fcm-cl-opd}, $C = \Sigma P$.  A $(\Sigma P)$-algebra is called a
\demph{categorical $P$-algebra}%
\index{categorical algebra for operad}%
or a \demph{$P$-category},%
\lbl{p:P-category}%
\index{operad!category over}
and consists of a directed graph $X$ together with a function
\begin{equation}	\label{eq:P-cat-action}
P(n) \times X(x_{n-1}, x_n) \times \cdots \times X(x_0, x_1)
\go
X(x_0, x_n)
\end{equation}
for each $n\in\nat$ and sequence $x_0, \ldots, x_n$ of vertices of $X$,
satisfying the evident axioms.  (Here $X(x,x')$ denotes the set of edges
from $x$ to $x'$.)  If $P=1$ then a $P$-category is exactly a category;
there are some less trivial examples too.

\begin{example}	\lbl{eg:fc-Trimble}%
\index{operad!path reparametrizations@of path reparametrizations}
Let $E$ be the operad in which 
\[
E(n) =
\{
\textrm{endpoint-preserving continuous maps }
[0,1] \go [0,n]
\},
\]
as in~\ref{eg:opd-Trimble}, so that any loop space is naturally an
$E$-algebra.  Then any path-space is naturally an $E$-category.  In other
words, take a topological space $Y$ and let $X$ be the graph whose vertices
are points of $Y$ and whose edges are continuous maps from $[0,1]$ into
$Y$: then $X$ is naturally a $E$-category.  
\end{example}

\begin{example}	\lbl{eg:fc-A-infty}%
\index{A-@$A_\infty$-!algebra}%
\index{A-@$A_\infty$-!category}
Similarly, taking $A_\infty$ to be the operad of chain complexes whose
algebras are $A_\infty$-algebras (p.~\pageref{p:A-infty-alg}) gives us the
standard notion of an \demph{$A_\infty$-category}.  For this to make sense,
we must work in a world where everything is enriched in chain complexes: so
$X$ consists of a set $X_0$ of objects (or vertices) together with a chain
complex $X(x,x')$ for each pair $(x,x')$ of objects, and
in~\bref{eq:P-cat-action} the $\times$'s become $\otimes$'s.  We do not
meet enriched generalized multicategories properly
until~\ref{sec:enr-mtis}, but it is clear how things should work in this
particular situation.
\end{example}

\section{Weak double categories}
\lbl{sec:wk-dbl}

Generalizing both bicategories and strict double categories are `weak
double categories', introduced informally in~\ref{sec:fcm}.  Recall that
these are only weak in the horizontal direction: vertical composition still
obeys strict laws.  We do not consider double categories weak in both
directions.

The definition of weak double category is a cross between that of strict
double category and that of unbiased bicategory.
\begin{defn}%
\index{double category!weak}
A \demph{weak double category} $D$ consists of some data subject
to some axioms.  The data is:
\begin{itemize}
\item A diagram
\[
\begin{slopeydiag}
	&	&D_1	&	&	\\
	&\ldTo<\dom&	&\rdTo>\cod&	\\
D_0	&	&	&	&D_0	\\
\end{slopeydiag}
\]
of categories and functors.  The objects of $D_0$ are called the
\demph{0-cells}%
\index{cell!weak double category@of weak double category}
or \demph{objects} of $D$, the maps in $D_0$ are the
\demph{vertical 1-cells} of $D$, the objects of $D_1$ are the
\demph{horizontal 1-cells} of $D$, and the maps in $D_1$ are the
\demph{2-cells} of $D$, as in the picture
\begin{equation}	\label{diag:wk-dbl-two-cell}
\begin{fcdiagram}
a	&\rTo^m			&a'	\\
\dTo<f	&\Downarrow\tcs{\theta}	&\dTo>{f'}\\
b	&\rTo_p			&b'	\\
\end{fcdiagram}
\end{equation}
where $a \goby{f} b$, $a' \goby{f'} b'$ are maps in $D_0$ and $m
\goby{\theta} p$ is a map in $D_1$, with $\dom(m) = a$, $\dom(\theta)=f$,
and so on. 
\item For each $n\geq 0$, a functor $\comp^{(n)}: D_1^{(n)} \go D_1$ such
that the diagram
\[
\begin{slopeydiag}
	&	&D_1^{(n)}		&	&	\\
	&\ldTo	&			&\rdTo	&	\\
D_0	&	&\dTo~{\comp^{(n)}}	&	&D_0	\\
	&\luTo<\dom&			&\ruTo>\cod&	\\
	&	&D_1			&	&	\\
\end{slopeydiag}
\]
commutes, where $D_1^{(n)}$ is the limit of the diagram
\[
\begin{slopeydiag}
	&	&D_1	&	&	&	&D_1	&	&
	&	&	&	&D_1	&	&	\\
	&\ldTo<\dom&	&\rdTo>\cod&	&\ldTo<\dom&	&\rdTo>\cod&
	&	&	&\ldTo<\dom&	&\rdTo>\cod&	\\
D_0	&	&	&	&D_0	&	&	&	&
\ 	&\cdots	&\ 	&	&	&	&D_0	\\
\end{slopeydiag}
\]
containing $n$ copies of $D_1$.  The functor $\comp^{(n)}$ is called
\demph{$n$-fold horizontal composition}%
\index{n-fold@$n$-fold!horizontal composition}
and written
\[
\begin{array}{rl}
&
\begin{fcdiagram}
a_0		&\rTo^{m_1}		&a_1		&\rTo^{m_2}	&
\cdots 	&\rTo^{m_n}		&a_n	\\
\dTo<{f_0}	&\Downarrow\tcs{\theta_1}&\dTo<{f_1}	&\Downarrow\tcs{\theta_2}&
	&\Downarrow\tcs{\theta_n}	&\dTo>{f_n}\\
b_0		&\rTo_{p_1}		&b_1		&\rTo_{p_2}	&
\cdots 	&\rTo_{p_n}		&b_n	\\
\end{fcdiagram}
\\
\\
\goesto	&
\begin{fcdiagram}
a_0 		&\rTo^{(m_n \of\cdots\of m_1)} 	&a_n	\\
\dTo<{f_0}	&\Downarrow\tcs{(\theta_n * \cdots * \theta_1)}&\dTo>{f_n}\\
b_0 		&\rTo_{(p_n \of\cdots\of p_1)} 	&b_n.	\\
\end{fcdiagram}
\end{array}
\]
\item For each double sequence 
\[
\mathbf{m} = 
((m_1^1, \ldots, m_1^{k_1}), \ldots, (m_n^1, \ldots, m_n^{k_n}))
\]
of horizontal 1-cells such that the composites below make sense, an
invertible 2-cell 
\[
\begin{fcdiagram}
\blob	&\rTo^{((m_n^{k_n} \of\cdots\of m_n^1) \of\cdots\of (m_1^{k_1}
\of\cdots\of m_1^1))}	&\blob	\\
\dTo<1	&\Downarrow \tcs{\gamma_{\mathbf{m}}}	&\dTo>1	\\
\blob 	&\rTo_{(m_n^{k_n} \of\cdots\of m_1^1)}	&\blob	\\
\end{fcdiagram}
\]
(where `invertible' refers to vertical composition).
\item For each horizontal 1-cell $a \goby{m} a'$, an invertible 2-cell
\[
\begin{fcdiagram}
a	&\rTo^m			&a'	\\
\dTo<1	&\Downarrow \tcs{\iota_m}	&\dTo>1	\\
a	&\rTo_{(m)}		&a'.	\\	
\end{fcdiagram}
\]
\end{itemize}
The axioms are:
\begin{itemize}
\item $\gamma_{\mathbf{m}}$ is natural in each of the $m_i^j$'s, and
$\iota_m$ is natural in $m$.  In the case of $\iota$, this means that for
each 2-cell $\theta$ as in~\bref{diag:wk-dbl-two-cell} we have
\[
\begin{fcdiagram}
a	&\rTo^m			&a'	\\
\dTo<f	&\Downarrow\tcs{\theta}	&\dTo>{f'}\\
b	&\rTo^p			&b'	\\
\dTo<1	&\Downarrow\tcs{\iota_p}&\dTo>1	\\
b	&\rTo_{(p)}		&b'	\\
\end{fcdiagram}
\diagspace
=
\diagspace
\begin{fcdiagram}
a	&\rTo^m			&a'	\\
\dTo<1	&\Downarrow\tcs{\iota_m}	&\dTo>1\\
a	&\rTo^{(m)}		&a'	\\
\dTo<f	&\Downarrow\tcs{(\theta)}	&\dTo>{f'}\\
b	&\rTo_{(p)}		&b',	\\
\end{fcdiagram}
\]
and similarly for $\gamma$.  
\item $\gamma$ and $\iota$ satisfy associativity and identity coherence
axioms analogous to those in the definition~(\ref{defn:lax-mon-cat}) of
lax monoidal category.  
\end{itemize}
\end{defn}

It is clear how to define lax maps between weak double categories, again
working by analogy with unbiased monoidal categories or unbiased
bicategories, and there is a forgetful functor
\[
(\textrm{weak double categories and lax maps})
\go
\fc\hyph\Multicat.
\]

\begin{example}	\lbl{eg:wk-dbl-degen}%
\index{double category!weak!degenerate}
There are several degenerate cases.  A weak double category in which the
coherence cells $\gamma$ and $\iota$ are all identities is exactly a strict
double category.  A weak double category whose only vertical 1-cells are
identities is exactly an (unbiased) bicategory.  A weak double category
whose only horizontal 1-cells are identities is exactly a strict
2-category. 
\end{example}

\begin{example}
The \fc-multicategory $\fcat{Ring}_2$%
\index{double category!weak!rings@of rings}
of~\ref{eg:fcm-Ring} is a weak double
category, horizontal composition being tensor of modules.
\end{example}

\begin{example}	\lbl{eg:wk-dbl-T-span}
The \fc-multicategory $\Set_2$%
\index{double category!weak!sets@of sets}
of~\ref{eg:fcm-wk-dbl} is also a weak double
category.  It is formed from sets, functions, and spans,%
\index{span}
which is
reminiscent of the definition of $T$-multicategories and maps between them
in~\ref{sec:om}.  Indeed, let $T$ be a cartesian monad on a cartesian
category $\Eee$, and define a weak double category as follows:
\begin{itemize}
\item objects are objects of $\Eee$
\item vertical 1-cells are maps in $\Eee$
\item horizontal 1-cells $E \go E'$ are diagrams
\[
\begin{slopeydiag}
	&	&M	&	&	\\
	&\ldTo<d&	&\rdTo>c&	\\
TE	&	&	&	&E'	\\
\end{slopeydiag}
\]
in $\Eee$
\item 2-cells $\theta$ are commutative diagrams
\[
\begin{slopeydiag}
	&	&L	&	&	\\
	&\ldTo	&	&\rdTo	&	\\
TD	&	&\dTo>\theta&	&D'	\\
\dTo<{Tf}&	&M	&	&\dTo>{f'}\\
	&\ldTo	&	&\rdTo	&	\\
TE	&	&	&	&E'	\\
\end{slopeydiag}
\]
\item vertical composition is as in $\Eee$
\item horizontal composition is by pullback: the composite of horizontal
1-cells
\[
\begin{slopeydiag}
	&	&M_1	&	&	\\
	&\ldTo	&	&\rdTo	&	\\
TE_0	&	&	&	&E_1,	\\
\end{slopeydiag}
\diagspace
\begin{slopeydiag}
	&	&M_2	&	&	\\
	&\ldTo	&	&\rdTo	&	\\
TE_1	&	&	&	&E_2,	\\
\end{slopeydiag}
\diagspace \ldots, \diagspace 
\begin{slopeydiag}
	&	&M_n	&	&	\\
	&\ldTo	&	&\rdTo	&	\\
TE_{n-1}&	&	&	&E_n	\\
\end{slopeydiag}
\]
is the limit of the diagram
\[
\begin{slopeydiag}
&&	&	&T^{n-1}M_1&	&	&	&T^{n-2}M_2&	&
	&	&	&	&M_n	&	&	\\
&&	&\ldTo	&	&\rdTo	&	&\ldTo	&	&\rdTo	&
	&	&	&\ldTo	&	&\rdTo	&	\\
&&T^n E_0&	&	&	&T^{n-1} E_1&	&	&	&
\ 	&\cdots	&\ 	&	&	&	&E_n	\\
&\ldTo<{\mu^{(n)}_{E_0}}&&&&	&	&	&	&	&
	&	&	&	&	&	&	\\
TE_0	&	&	&	&	&	&	&	&
	&	&	&	&	&	&	\\
\end{slopeydiag}
\]
where $\mu^{(n)}$ is the $n$-fold multiplication of the monad $T$; 
horizontal composition of 2-cells works similarly.
\end{itemize}
Any weak double category%
\index{double category!weak!bicategory@\vs.\ bicategory}%
\index{bicategory!weak double category@\vs.\ weak double category}
yields a bicategory by discarding all the vertical
1-cells except for the identities, and applying this to the weak double
category just defined yields the bicategory $\Sp{\Eee}{T}$
(Definition~\ref{defn:T-spans}).  It is therefore reasonable to call our
weak double category $\Sp{\Eee}{T}$%
\glo{Spwkdbl}
too.  In particular, $\Sp{\Set}{\id} =
\Set_2$.
\end{example}

\begin{example}	\lbl{eg:wk-dbl-Cat}
Categories, functors and modules%
\index{module!categories@over categories}
are integrated in the following weak
double category $\Cat_2$:%
\glo{wkdblCat2}%
\index{double category!weak!categories@of categories}%
\index{category!weak double category of}
\begin{itemize}
\item objects are small categories
\item vertical 1-cells are functors
\item horizontal 1-cells are modules~(\ref{sec:om-further}): that is, a
horizontal 1-cell $A \go A'$ is a functor $A^\op \times A' \go \Set$
\item a 2-cell
\[
\begin{fcdiagram}
A	&\rTo^M		&A'	\\
\dTo<F	&\Downarrow	&\dTo>{F'}\\
B	&\rTo_P		&B'	\\
\end{fcdiagram}
\]
is a natural transformation
\[
\begin{diagram}[height=1.2em]
A^\op \times A'		&		&	&	\\
			&\rdTo(3,2)>M	&	&	\\
\dTo<{F\times F'}	&\sent		&	&\Set	\\
			&		&\ruTo(3,2)>P	&	\\
B^\op \times B'		&		&	&	\\
\end{diagram}
\]
\item vertical composition is ordinary composition of functors and
natural transformations
\item horizontal composition is tensor of modules: the composite $(M_n
\otimes \cdots \otimes M_1)$ of
\[
\begin{fcdiagram}
A_0	&\rTo^{M_1}	&A_1		&\rTo^{M_2}	&
\ 	&\cdots	&\ 	&\rTo^{M_n} 	&A_n,	\\
\end{fcdiagram}
\]
is given by the coend formula
\[
(M_n \otimes\cdots\otimes M_1)(a_0, a_n)
=
\int^{a_1, \ldots, a_{n-1}}
M_n(a_{n-1}, a_n) \times\cdots\times M_1(a_0, a_1)
\glo{coend}
\]
($a_0 \in A_0, a_n \in A_n$), or when $n=0$ by $I_A(a,a') = A(a,a')$;%
\glo{idcatmodule}
similarly 2-cells.  (A coend%
\index{coend}
is a colimit of sorts: Mac
Lane~\cite[Ch.~IX]{MacCWM}.)
\end{itemize}
The weak double category $\Cat_2$ incorporates not only categories,
functors and modules, but also natural transformations: for by a Yoneda%
\index{Yoneda!Lemma}
argument, a 2-cell
\[
\begin{fcdiagram}
A	&\rTo^{I_A}	&A	\\
\dTo<F	&\Downarrow	&\dTo>{F'}\\
B	&\rTo_{I_B}	&B	\\
\end{fcdiagram}
\]
amounts to a natural transformation $F \go F'$.
\end{example}

In the next section we will see that $T$-multicategories and maps,
transformations and modules between them form an \fc-multicategory, but not
in general a weak double category.

\begin{example}	\lbl{eg:wk-dbl-Monoid}%
\index{double category!weak!monoids@of monoids}
There is a weak double category $\fcat{Monoid}_2$ made up of monoids,
homomorphisms of monoids, bimodules over monoids, and maps between them.
This is just like $\fcat{Ring}_2$ but with all the additive structure
removed.  Alternatively, it is the substructure of $\Cat_2$ whose 0-cells
are just the one-object categories, and with all 1- and 2-cells between
them.
\end{example}

\begin{example}	\lbl{eg:wk-dbl-cospan}
We have considered using spans as horizontal 1-cells; we can also use
cospans.%
\index{cospan}
 Thus, if $\Eee$ is any category in which all pushouts exist then
there is a weak double category whose objects and vertical 1-cells form the
category $\Eee$, whose horizontal 1-cells $A \go A'$ are diagrams
\[
\begin{slopeydiag}
	&	&M	&	&	\\
	&\ruTo	&	&\luTo	&	\\
A	&	&	&	&A'	\\
\end{slopeydiag}
\]
in $\Eee$, and with the rest of the structure defined in the obvious way.  

For instance, let $\Eee$ be the category of topological%
\index{double category!weak!topological spaces@of topological spaces}
spaces and
embeddings: then $M$ is a space containing copies of $A$ and $A'$ as
subspaces, and horizontal composition in the weak double category is gluing
along subspaces. 
\end{example}

\begin{example}%
\index{double category!weak!manifolds@of manifolds}%
\index{manifold!weak double category of}
Adapting the previous example slightly, we obtain a weak double category
$n\hyph\fcat{Mfd}_2$ of $n$-manifolds and cobordisms,%
\index{cobordism}
for each $n\geq 0$.
An object is an $n$-dimensional manifold (topological, say, and without
boundary).  A vertical 1-cell is a continuous map.  A horizontal 1-cell $A
\go A'$ is an $(n+1)$-manifold $M$ together with a homeomorphism $h$
between the boundary $\bdry M$ of $M$ and the disjoint union $A \disjt A'$.
A 2-cell
\[
\begin{fcdiagram}
A	&\rTo^{(M,h)}	&A'		\\
\dTo<f	&\Downarrow	&\dTo>{f'}	\\
B	&\rTo_{(P,k)}	&B'		\\
\end{fcdiagram}
\]
is a continuous map $\theta: M \go P$ making the diagram
\[
\begin{diagram}[size=2em]
\bdry M		&\rTo_{\diso}^h	&A \disjt A'		\\
\dTo<\theta	&		&\dTo>{f\disjt f'}	\\
\bdry P		&\rTo^{\diso}_k	&B \disjt B'		\\
\end{diagram}
\]
commute.  Vertical composition is composition of functions; horizontal
composition is by gluing.  

Fig.~\ref{fig:cobordisms} 
\begin{figure}
% \[
\setlength{\unitlength}{1mm}
% \begin{array}{c}
\centering
\begin{picture}(102,46)(-1,0)
% cobordisms
\cell{-1}{0}{bl}{\epsfig{file=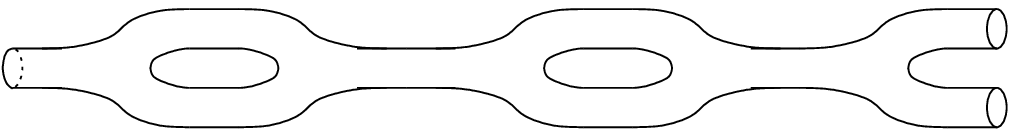}}
\cell{-1}{26}{bl}{\epsfig{file=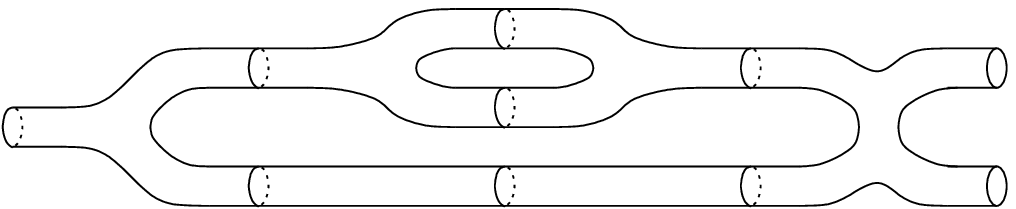}}
% arrows
\put(0,30){\vector(0,-1){20}}
\cell{2}{19}{l}{\cong}
\put(100,24){\vector(0,-1){10}}
\cell{98}{19}{r}{\cong}
\cell{49}{19}{l}{\Downarrow\ \cong}
\end{picture}
% \end{array}
% \]
% \hand{55}{15}
\caption{2-cell in the \fc-multicategory $1\hyph\fcat{Mfd}_2$, with 4
horizontal 1-cells along the top row}
\label{fig:cobordisms}
\end{figure}
shows a 2-cell in the underlying \fc-multicategory of $1\hyph\fcat{Mfd}_2$,
in a case where all the continuous maps involved are homeomorphisms.

Similar constructions can be made for other types of manifold: oriented,
smooth, holomorphic, \ldots.  There is a slight problem with identities for
horizontal composition, as the horizontal identity on $A$ `ought' to be
just $A$ itself (compare~\ref{eg:wk-dbl-cospan}) but this is not usually
counted as an $(n+1)$-manifold with boundary.  Whether or not we can fix
this, there is certainly an underlying $\fc$-multicategory: the problem is
purely one of representability.
\end{example}

\section{Monads, monoids and modules}
\lbl{sec:mmm}

Modules have loomed large in our examples of \fc-multicategories, appearing
as horizontal 1-cells.  Monoids and their cousins (such as rings and
categories) have also been prominent, appearing as 0-cells.  Here we show
how any \fc-multicategory $C$ gives rise to a new \fc-multicategory
$\Mon(C)$, whose 0-cells are monoids/monads in $C$ and whose horizontal
1-cells are modules between them.  

As I shall explain, this construction has traditionally been carried out in
a narrower context than \fc-multicategories, which has meant working under
certain technical restrictions.  If we expand to the wider context of
\fc-multicategories then the technicalities vanish.  

For an example of the traditional construction, start with the monoidal
category $(\Ab, \otimes, \integers)$ of abelian groups.  A monoid therein
is just a ring.  Given two monoids $A$, $B$ in a monoidal category, a
\demph{$(B,A)$-module}%
\index{module!monoids@over monoids}
is an object $M$ equipped with compatible left and
right actions
\[
B \otimes M \goby{\act{\mi{L}}} M,
\diagspace
M \otimes A \goby{\act{\mi{R}}} M,
\]
and there is an obvious notion of map between $(B,A)$-modules.  In this
particular case these are (bi)modules and their homomorphisms, in the usual
sense.  So we have almost arrived at the bicategory of~\ref{eg:bicat-Ring}:
0-cells are rings, 1-cells are modules, and 2-cells are maps of modules.
The only missing ingredient is the tensor product of modules.  If $M$ is a
$(B,A)$-module and $N$ a $(C,B)$-module then $N \otimes_B M$ is a quotient
of the abelian group tensor $N\otimes M$; categorically, it is the
(reflexive) coequalizer%
\index{coequalizer}
\[
\begin{diagram}[scriptlabels]
N \otimes B \otimes M	&
\pile{	\rTo^{\act{\mi{R}} \otimes 1_M}\\
	\rTo_{1_N \otimes \act{\mi{L}}} }	&
N\otimes M	&
\rTo	&
N \otimes_B M.	\\
\end{diagram}
\]
The abelian group $N\otimes M$ acquires a $(C,A)$-module structure just as
long as the endofunctors $C\otimes \dashbk$ and $\dashbk\otimes A$ of $\Ab$
preserve (reflexive) coequalizers, which they do.  The rest of the
bicategory structure comes easily.  So from the monoidal category $(\Ab,
\otimes, \integers)$ we have derived the bicategory of rings, modules and
maps of modules.  

In the same way, any monoidal category $(\cat{A}, \otimes, I)$ gives rise
to a bicategory $\Mon(\cat{A})$ of monoids and modules in $\cat{A}$, as
long as $\cat{A}$ has reflexive coequalizers and these are preserved by the
functors $A \otimes \dashbk$ and $\dashbk \otimes A$ for all $A \in
\cat{A}$.  Generalizing, any bicategory $\cat{B}$ gives rise to a
bicategory $\Mon(\cat{B})$ of monads and modules in $\cat{B}$, as long as
$\cat{B}$ locally has reflexive coequalizers and these are preserved by the
functors $f\of\dashbk$ and $\dashbk\of f$ for all 1-cells $f$.

There are two unsatisfactory aspects to this construction.  One is the
necessity of the coequalizer conditions.  The other is that homomorphisms
of monoids are conspicuous by their absence: for instance, ring
homomorphisms are not an explicit part of $\Mon(\Ab)$ (although they can be
recovered as those 1-cells of $\Mon(\Ab)$ that have a right adjoint).  The
following definition for \fc-multicategories solves both problems.

\begin{defn}
Let $C$ be an \fc-multicategory. The \fc-multicategory $\Mon(C)$%
\glo{fcmtiMon}%
\index{monads construction}
is defined as
follows.
\begin{itemize}
\item
A 0-cell of $\Mon(C)$ is an \fc-multicategory map $1\go C$.  That is, it is
a 0-cell $a$ of $C$ together with a horizontal 1-cell $a\goby{t}a$ and
2-cells
\[
\begin{fcdiagram}
a	&\rTo^{t}	&a			&\rTo^{t}	&a	\\
\dTo<{1}&		&\Downarrow\tcs{\mu}	&		&\dTo>{1}\\
a	&		&\rTo_{t}		&		&a	\\
\end{fcdiagram}
\diagspace
\begin{fcdiagram}
a	&\rEquals		&a	\\
\dTo<{1}&\Downarrow\tcs{\eta}	&\dTo>{1}\\
a	&\rTo_{t}		&a	\\
\end{fcdiagram}
\]
satisfying the usual axioms for a monad, $\mu\of\pr{\mu}{1_t} =
\mu\of\pr{1_t}{\mu}$ and $\mu\of\pr{\eta}{1_t} = 1_t = \mu\of\pr{1_t}{\eta}$.

\item
A vertical 1-cell
\[
\begin{fcdiagram}
(a,t,\mu,\eta)	\\
\dTo		\\
(\hat{a},\hat{t},\hat{\mu},\hat{\eta})\\
\end{fcdiagram}
\]
in $\Mon(C)$ is a vertical 1-cell
\[
\begin{fcdiagram}
a	\\
\dTo>f	\\
\hat{a}	\\
\end{fcdiagram}
\]
in $C$ together with a 2-cell
\[
\begin{fcdiagram}
a		&\rTo^{t}		&a		\\
\dTo<{f}	&\Downarrow\tcs{\omega}	&\dTo>{f}	\\
\hat{a}		&\rTo_{\hat{t}}		&\hat{a}	\\
\end{fcdiagram}
\]
such that $\omega\of\mu = \hat{\mu}\of\pr{\omega}{\omega}$ and
$\omega\of\eta = \hat{\eta}\of f$.  (The notation $\hat{\eta} \of f$ is
explained on p.~\pageref{p:null-notation}.)

\item
A horizontal 1-cell $(a,t,\mu,\eta) \rTo (a',t',\mu',\eta')$ is a
horizontal 1-cell $a\goby{m}a'$ in $C$ together with 2-cells
\[
\begin{fcdiagram}
a	&\rTo^{t}	&a			&\rTo^{m}	&a'	\\
\dTo<{1}&		&\Downarrow\tcs{\chi}	&		&\dTo>{1}\\
a	&		&\rTo_{m}		&		&a'	\\
\end{fcdiagram}
\diagspace
\begin{fcdiagram}
a	&\rTo^{m}	&a'			&\rTo^{t'}	&a'	\\
\dTo<{1}&		&\Downarrow\tcs{\chi'}	&		&\dTo>{1}\\
a	&		&\rTo_{m}		&		&a'	\\
\end{fcdiagram}
\]
satisfying the usual module axioms, $\chi\of\pr{\mu}{1_m} =
\chi\of\pr{1_t}{\chi}$, $\chi\of\pr{\eta}{1_m}=1_m$, and dually for
$\chi'$, and the `commuting actions' axiom, $\chi'\of\pr{\chi}{1_{t'}} =
\chi\of\pr{1_t}{\chi'}$.

\item
A 2-cell
\[
\begin{fcdiagram}
t_0	&\rTo^{m_1}	&t_1	&\rTo^{m_2}	&\ 	&\cdots	
&\ 	&\rTo^{m_{n}}	&t_n	\\
\dTo<{f}&		&	&		&\Downarrow	&	
&	&		&\dTo>{f'}\\
t	&		&	&		&\rTo_{m}	&	
&	&		&t'	\\
\end{fcdiagram}
\]
in $\Mon(C)$, where $t$ stands for $(a,t,\mu,\eta)$, $m$ for
$(m, \chi, \chi')$, $f$ for \pr{f}{\omega}, and so on, is a 2-cell
\[
\begin{fcdiagram}
a_0	&\rTo^{m_1}	&a_1	&\rTo^{m_2}	&\ 	&\cdots	
&\ 	&\rTo^{m_{n}}	&a_n	\\
\dTo<{f}&		&	&		&\Downarrow\tcs{\theta}&	
&	&		&\dTo>{f'}\\
a	&		&	&		&\rTo_{m}	&	
&	&		&a'	\\
\end{fcdiagram}
\]
in $C$ satisfying the `external equivariance' axioms
\begin{eqnarray*}
\theta\of(\chi_1,\range{1_{m_2}}{1_{m_n}}) 		&=&
\chi\of\pr{\omega}{\theta}				\\
\theta\of(\range{1_{m_1}}{1_{m_{n-1}}},\chi'_n)	&=&
\chi'\of\pr{\theta}{\omega'}
\end{eqnarray*}
and the `internal equivariance' axioms
\begin{eqnarray*}
\lefteqn{\theta\of(\range{1_{m_1}}{1_{m_{i-2}}}, \chi'_{i-1}, 1_{m_{i}},
\range{1_{m_{i+1}}}{1_{m_n}})		=}				\\
&&\theta\of(\range{1_{m_1}}{1_{m_{i-2}}}, 1_{m_{i-1}}, \chi_i,
\range{1_{m_{i+1}}}{1_{m_n}}) 
\end{eqnarray*}
for $2\leq i\leq n$.

\item
Composition and identities for both 2-cells and vertical 1-cells in
$\Mon(C)$ are just composition and identities in $C$.
\end{itemize}
\end{defn}
With the obvious definition on maps of \fc-multicategories, this gives a
functor 
\[
\Mon: \fc\hyph\Multicat \go \fc\hyph\Multicat.
\]

\begin{example}%
\index{double category!weak!bicategory@\vs.\ bicategory}%
\index{bicategory!weak double category@\vs.\ weak double category}
Our new $\Mon$ generalizes the traditional $\Mon$ in the following sense.
Let $B$ be a bicategory satisfying the conditions on local reflexive
coequalizers mentioned above, so that it is possible to construct the
bicategory $\Mon(B)$ in the traditional way.  Let $C$ be the
\fc-multicategory corresponding to $B$ (Example~\ref{eg:fcm-bicat}), with
only trivial vertical 1-cells.  Then a 0-cell of $\Mon(C)$ is a monad in
$B$, a horizontal 1-cell $t \go t'$ is a $(t',t)$-module in $B$, and a
2-cell of the form 
\[
\begin{fcdiagram}
t_0	&\rTo^{m_1}	&t_1	&\rTo^{m_2}	&\ 	&\cdots	
&\ 	&\rTo^{m_{n}}	&t_n	\\
\dTo<1	&		&	&		&\Downarrow	&	
&	&		&\dTo>1	\\
t_0	&		&	&		&\rTo_m		&	
&	&		&t_n	\\
\end{fcdiagram}
\]
is a map
\[
m_n \otimes_{t_{n-1}} \cdots 
\otimes_{t_2} m_2 \otimes_{t_1} m_1
\go
m
\]
of $(t_n, t_0)$-modules---that is, a 2-cell in $\Mon(B)$.  So if we discard
the non-identity vertical 1-cells of $\Mon(C)$ then we obtain the
\fc-multicategory corresponding to $\Mon(B)$.
\end{example}

\begin{example}%
\index{fc-multicategory@$\fc$-multicategory!rings@of rings}
Take the monoidal category $(\Ab,\otimes,\integers)$ and the corresponding
\fc-multicategory $\Sigma\Ab$ (Example~\ref{eg:fcm-mon-cat}).  Then
$\Mon(\Sigma\Ab)$ is $\fcat{Ring}_2$, the \fc-multicategory of
Example~\ref{eg:fcm-Ring} in which 0-cells are rings, vertical 1-cells are
homomorphisms of rings, horizontal 1-cells are modules%
\index{module!bimodule over rings}
between rings, and
2-cells are module maps of a suitable kind.  
\end{example}

\begin{example}%
\index{fc-multicategory@$\fc$-multicategory!monoids@of monoids}
Similarly, if we take the monoidal category $(\Set,\times,1)$ then
$\Mon(\Sigma\Set)$ is (the underlying \fc-multicategory of) the weak double
category $\fcat{Monoid}_2$ (Example~\ref{eg:wk-dbl-Monoid}).
\end{example}

\begin{example}
Take the \fc-multicategory $\Set_2$ of sets, functions and
spans~(\ref{eg:fcm-wk-dbl}).  Then a 0-cell of $\Mon(\Set_2)$ is a monad in
$\Set_2$, that is, a small category.%
\index{fc-multicategory@$\fc$-multicategory!categories@of categories}%
\index{double category!weak!categories@of categories}%
\index{category!weak double category of}
 A vertical 1-cell is a functor; a
horizontal 1-cell is a module.  A 2-cell
\[
\begin{fcdiagram}
A_0	&\rTo^{M_1}	&A_1	&\rTo^{M_2}	&\ 	&\cdots	
&\ 	&\rTo^{M_n}	&A_n	\\
\dTo<F	&		&	&		&\Downarrow\tcs{\theta}&
&	&		&\dTo>{F'}\\
A	&		&	&		&\rTo_{M}	&	
&	&		&A'	\\
\end{fcdiagram}
\]
consists of a function
\[
\theta_{a_0, \ldots, a_n}: 
M_n(a_{n-1}, a_n) \times\cdots\times M_1(a_0, a_1) 
\go
M(Fa_0, F'a_n)
\]
for each $a_0 \in A_0, \ldots, a_n \in A_n$, such that this family is
natural in the $a_i$'s.  So $\Mon(\Set_2)$ is $\Cat_2$, the weak double
category of~\ref{eg:wk-dbl-Cat}.
\end{example}

This last example leads us to definitions of module and transformation for
generalized multicategories.  Let $T$ be a cartesian monad on a cartesian
category $\Eee$.  Recall the \fc-multicategory $\Sp{\Eee}{T}$
of~\ref{eg:wk-dbl-T-span}, and consider the \fc-multicategory
$\Mon(\Sp{\Eee}{T})$.  A 0-cell of $\Mon(\Sp{\Eee}{T})$ is a monad in
$\Sp{\Eee}{T}$, that is, a $T$-multicategory, so we write $T\hyph\Multicat$
for $\Mon(\Sp{\Eee}{T})$.  A vertical 1-cell is a map of
$T$-multicategories.  A horizontal 1-cell $A \go B$ will be called a
\demph{$(B,A)$-module}.%
\index{module!generalized multicategories@over generalized multicategories}%
\index{generalized multicategory!module over}
 To see what modules are explicitly, note that the
underlying graph
\[
\begin{slopeydiag}
	&	&B_1	&	&	\\
	&\ldTo<\dom&	&\rdTo>\cod&	\\
TB_0	&	&	&	&B_0	\\
\end{slopeydiag}
\]
of $B$ is a horizontal 1-cell in $\Sp{\Eee}{T}$, composition with which
defines an endofunctor $B\of\dashbk$ of the category $\Sp{\Eee}{T} (A_0,
B_0)$ of diagrams of the form  
\begin{equation}	\label{diag:module-span}
\begin{slopeydiag}
	&	&M	&	&	\\
	&\ldTo	&	&\rdTo	&	\\
TA_0	&	&	&	&B_0.	\\
\end{slopeydiag}
\end{equation}
Moreover, the multicategory structure of $B$ gives $B\of\dashbk$ the
structure of a monad on $\Sp{\Eee}{T} (A_0, B_0)$.  Dually, the
multicategory $A$ also induces a monad $\dashbk\of A$ on $\Sp{\Eee}{T}
(A_0, B_0)$, and there is a canonical isomorphism $(B\of\dashbk)\of A \iso
B\of(\dashbk\of A)$ compatible with the monad structures (formally, a
distributive law:~\ref{defn:distrib-law}).  So a $(B,A)$-module is a
diagram~\bref{diag:module-span} in $\Eee$ together with maps $\act{\mi{L}}:
B \of M \go M$ and $\act{\mi{R}}: M \of A \go M$ such that
$(M,\act{\mi{L}})$ and $(M,\act{\mi{R}})$ are algebras for the monads
$B\of\dashbk$ and $\dashbk\of A$ respectively, and the following diagram
commutes:
\[
\begin{diagram}[size=2em,scriptlabels]
	&(B\of M)\of A 	&\iso 	&B\of (M\of A)	&	\\
\ldTo(1,2)<{\act{\mi{L}}*1_A}&&	&		&
\rdTo(1,2)>{1_B * \act{\mi{R}}}\\
M\of A	&		&	&		&B\of M	\\
	&\rdTo<{\act{\mi{R}}}&	&\ldTo>{\act{\mi{L}}}&	\\
	&		&M.	&		&	\\
\end{diagram}
\]

This definition generalizes the definition for plain
multicategories~(\ref{sec:om-further}).  Just as there, it also makes sense
to define left and right modules for a $T$-multicategory: do not take a
whole span $M$, but only the right-hand or left-hand half.  A \demph{left
$B$-module} is nothing other than a $B$-algebra.  A \demph{right $A$-module}
(or \demph{$A$-coalgebra})%
\index{coalgebra}
is an object $M$ of $\Eee/TA_0$ (thought of as a
span
\[
\begin{slopeydiag}
	&	&M	&	&	\\
	&\ldTo	&	&\rdGet	&	\\
TA_0	&	&	&	&\cdot	\\
\end{slopeydiag}
\]
with a phantom right leg) together with a map $h: M\of A \go M$ making $M$
into an algebra for the monad $\dashbk\of A$ on $\Eee/TA_0$.  

In general, $T\hyph\Multicat$ is not a weak double category: there is no
tensor%
\index{tensor!generalized multicategory@over generalized multicategory}
product of modules.  For to tensor modules we would need $\Eee$ to
possess reflexive coequalizers and $T$ to preserve them, and although the
requirement on $\Eee$ is satisfied in almost all of the examples that we
consider, the requirement on $T$ is not---for instance, when $T=\fc$.

There is, however, an `identity' horizontal 1-cell on each object.  Given a
$T$-multicategory $A$, this is the $(A,A)$-module $I_A$ whose underlying
span is
\[
\begin{slopeydiag}
	&	&A_1	&	&	\\
	&\ldTo<\dom&	&\rdTo>\cod&	\\
TA_0	&	&	&	&A_0	\\
\end{slopeydiag}
\]
and whose actions are both composition in $A$.  Recalling the remarks on
natural transformations in Example~\ref{eg:wk-dbl-Cat}, we make the
following definition.  Let $A \parpair{f}{f'} B$ be maps of
$T$-multicategories.  A \demph{transformation}%
\index{transformation!generalized multicategories@of generalized multicategories}%
\index{generalized multicategory!transformation of}
\begin{equation}	\label{diag:T-transf}
A 
\ctwomult{\scriptstyle f}{\scriptstyle f'}{\scriptstyle \alpha} 
B
\end{equation}
of $T$-multicategories is a 2-cell 
\[
\begin{fcdiagram}
A	&\rTo^{I_A}		&A	\\
\dTo<f	&\Downarrow\tcs{\alpha}	&\dTo>{f'}\\
B	&\rTo_{I_B}		&B
\end{fcdiagram}
\]
in $T\hyph\Multicat$.  Explicitly, a transformation $f \go f'$ is a map
$\alpha: A_1 \go B_1$ such that the diagrams 
\begin{equation}	\label{eqs:T-transf}
\begin{slopeydiag}
	&		&A_1	&		&	\\
	&\ldTo^\dom	&	&\rdTo^\cod	&	\\
TA_0	&		&\dTo>\alpha&		&A_0	\\
\dTo<{Tf_0}&		&B_1	&		&\dTo>{f'_0}\\
	&\ldTo^\dom	&	&\rdTo^\cod	&	\\
TB_0	&		&	&		&B_0	\\
\end{slopeydiag}
\diagspace
\diagspace
%
% \begin{slopeydiag}
\begin{diagram}[size=2em,scriptlabels]
	&		&A_1\of A_1	&		&	\\
	&\ldTo<{\alpha*f_1}&\dTo~\comp	&\rdTo>{f'_1*\alpha}&	\\
B_1\of B_1&		&A_1		&		&B_1\of B_1\\
	&\rdTo<\comp	&\dTo~\alpha	&\ldTo>\comp	&	\\
	&		&B_1		&		&	\\
\end{diagram}
% \end{slopeydiag}
\end{equation}
commute, where the map $\alpha * f_1: A_1 \of A_1 \go B_1\of B_1$ is
induced by the maps $\alpha: A_1 \go B_1$ and $f_1: A_1 \go B_1$,
% via the definitions of $A_1 \of A_1$ and $B_1 \of B_1$ as pullbacks, 
and similarly $f'_1 * \alpha$.

This may come as a surprise: if $A$ and $B$ are ordinary categories and
$\alpha$ a natural transformation as in~\bref{diag:T-transf} then we would
usually think of $\alpha$ as a map $A_0 \go B_1$, not $A_1 \go B_1$.
Observe, however, that $\alpha$ assigns to each map $a \goby{\theta} a'$ in
$A$ a map $fa \goby{\alpha_\theta} f'a'$ in $B$, the diagonal of the
naturality square for $\alpha$ at $\theta$.  Since $\alpha_a$ can be
recovered as $\alpha_{1_a}$ for each $a\in A$, it must be possible to
define a natural transformation~\bref{diag:T-transf} as a map $A_1 \go B_1$
satisfying axioms.  The same goes for transformations of plain
multicategories, where now $\alpha_\theta$ is either side of the equation
in Fig.~\ref{fig:cl-transf-axiom} (p.~\pageref{fig:cl-transf-axiom}).  For
generalized multicategories we also have the choice between viewing
transformations as maps $A_0 \go B_1$ or as maps $A_1 \go B_1$; the axioms
for the $A_0 \go B_1$ version can easily be worked out (or looked up:
Leinster \cite[1.1.1]{GECM}).  

For general cartesian $T$, suppose that we discard from the
\fc-multicategory $T\hyph\Multicat$ all of the horizontal 1-cells except
for the `identities' (those of the form $I_A$).  Then we obtain a weak
double category whose only horizontal 1-cells are identities---in other
words~(\ref{eg:wk-dbl-degen}), a strict 2-category.  Concretely:
transformations can be composed vertically and horizontally, and this makes
$T\hyph\Multicat$%
\glo{TMulticat2cat}
into a strict 2-category.

\begin{notes}

Many \fc-multicategories are weak double categories.  Many
\fc-multicategories also have the property that every vertical 1-cell
$f: a \go b$ gives rise canonically to cells
\[
\begin{fcdiagram}
a	&\rEquals	&a	\\
\dTo<1	&\Downarrow	&\dTo>f	\\
a	&\rTo		&b\makebox[0em][l]{,}	\\	
\end{fcdiagram}
\diagspace
\begin{fcdiagram}
a	&\rEquals	&a	\\
\dTo<f	&\Downarrow	&\dTo>1	\\
b	&\rTo		&a
\end{fcdiagram}
\]
(as is familiar in $\fcat{Ring}_2$, for instance).  These two properties
combine to say that any diagram
\[
\begin{fcdiagram}
a_0	&\rTo^{m_1}	&a_1	&\rTo^{m_2}	&\ 	&\cdots	
&\ 	&\rTo^{m_n}	&a_n	\\
\dTo<{f}&		&	&		&	&
&	&		&\dTo>{f'}\\
a	&		&	&		&	&	
&	&		&a'	\\
\end{fcdiagram}
\]
has a 2-cell filling it in (canonically, up to isomorphism): the weak
double category structure provides a fill-in when $f$ and $f'$ are both
identities, the vertical-to-horizontal property provides a fill-in when
$n=0$ and one of $f$ or $f'$ is an identity, and the general case can be
built up from these special cases.  The exact situation is still unclear,
but there are obvious similarities between this, the representations of
plain multicategories in~\ref{sec:non-alg-notions}, and the universal
opetopic cells in~\ref{sec:ope-n}.

The $\fc$-multicategory $\Set[V]$ of~\ref{eg:fcm-matrices} appears in a
different guise in Carboni,%
\index{Carboni, Aurelio}
Kasangian%
\index{Kasangian, Stefano}
and Walters~\cite{CKW},%
\index{Walters, Robert}
under the
name of `matrices'.%
\index{matrix}
 Weak double categories have also been studied by
Grandis%
\index{Grandis, Marco}
and Par\'e~\cite{GPLDC, GPADC}.%
\index{Pare, Robert@Par\'e, Robert}

\end{notes}

\chapter{Generalized Operads and Multicategories: Further Theory}
\lbl{ch:gom-further}

\chapterquote{%
The last paragraph plays on the postmodern fondness for
`multidimensionality' and `nonlinearity' by inventing a nonexistent field:
`multidimensional (nonlinear) logic'}{%
Sokal and Bricmont~\cite{SB}}

\noindent
This chapter is an assortment of topics in the theory of
$T$-multicategories.  Some are included because they answer natural
questions, some because they connect to established concepts for classical
operads, and some because we will need them later.  The reader who wants to
get on to geometrically interesting structures should skip this chapter and
come back later if necessary; there are no pictures here.

In~\ref{sec:more-monads} we `recall' some categorical language: maps
between monads, mates under adjunctions, and distributive laws.  This
language has nothing intrinsically to do with generalized multicategories,
but it will be efficient to use it in some parts of this chapter
(\ref{sec:alt-app},~\ref{sec:change}) and later chapters.  

The first three proper sections each recast one of the principal
definitions.  In~\ref{sec:alt-app} we find an alternative definition of
generalized multicategory, which amounts to characterizing a
$T$-multicategory $C$ by its free-algebra monad $T_C$ plus one small extra
piece of data.  Both~\ref{sec:alg-fibs} and~\ref{sec:endos} are
alternative ways of defining an algebra for a $T$-multicategory.  The
former generalizes the categorical fact that $\Set$-valued functors can be
described as discrete fibrations, and the latter the fact that a classical
operad-algebra is a map into an endomorphism operad (often taken as the
\emph{definition} of algebra by `working mathematicians' using operads).

The next two sections are also generalizations.  The free $T$-multicategory
on a $T$-graph is discussed in~\ref{sec:free-mti}.  Abstract as this may
seem, it is crucial to the way in which geometry arises spontaneously from
category theory: witness the planar trees of~\ref{sec:om-further} and the
opetopes of Chapter~\ref{ch:opetopic}.  Then we make a definition to
complete the phrase: `plain multicategories are to monoidal categories as
$T$-multicategories are to \ldots'; we call these things `$T$-structured
categories'~(\ref{sec:struc}).

So far all is generalization, but the final two sections are genuinely new.
A choice of cartesian monad $T$ specifies the type of input shape that the
operations in a $T$-multicategory will have, so changing $T$ amounts to a
change of shape.  In~\ref{sec:change} we show how to translate between
different shapes, in other words, how a relation between two monads $T$ and
$T'$ induces a relation between the classes of $T$- and
$T'$-multicategories.  Finally, in~\ref{sec:enr-mtis} we take a short look
at enrichment of generalized multicategories.  As mentioned in the
Introduction, there is much more to this than one might guess, and in
particular there is more than we have room for; this is just a taste.

\section{More on monads}
\lbl{sec:more-monads}%
\index{monad|(}

To compare approaches to higher categorical structures using different
shapes or of different dimensions we will need a notion of map $(\Eee, T)
\go (\Eee', T')$ between monads.  Actually, there are various such notions:
lax, colax, weak and strict.  Such comparisons lead to functors between
categories of structures, and to discuss adjunctions between such functors
it will be convenient to use the language of mates (an Australian
creation, of course).  It will also make certain later proofs
(\ref{thm:wk-2},~\ref{propn:free-enr}) easier if we know a little about
distributive laws, which are recipes for gluing together two monads $T$,
$T'$ on the same category to give a monad structure on the composite
functor $T'\of T$.

Nothing here is new.  I learned this material from Street~\cite{StrFTM}%
\index{Street, Ross}
and
Kelly%
\index{Kelly, Max}
and Street~\cite{KSRE2},%
\index{Street, Ross}
although I have changed some terminology.
Distributive laws were introduced by Beck~\cite{Beck}.%
\index{Beck, Jon}
 Where Street
discusses monads in an arbitrary 2-category $\cat{V}$, we stick to the case
$\cat{V} = \CAT$, since that is all we need.

Let $T = (T, \mu, \eta)$ be a monad on a category $\Eee$ and $T' = (T',
\mu', \eta')$ a monad on a category $\Eee'$.  A \demph{lax map of monads}%
\index{monad!map of|(}%
\lbl{p:lax-map-of-monads}
$(\Eee, T) \go (\Eee', T')$ is a functor $Q: \Eee \go \Eee'$ together with
a natural transformation
\[
\begin{diagram}
\Eee		&\rTo^T		&\Eee		\\
\dTo<Q		&\nent\tcs{\psi}&\dTo>Q		\\
\Eee'		&\rTo_{T'}	&\Eee'		\\
\end{diagram}
\]
making the diagrams 
\[
\begin{diagram}[height=2em]
T'^2 Q		&\rTo^{T'\psi}	&T'QT		&\rTo^{\psi T}	&QT^2	\\
\dTo<{\mu' Q}	&		&		&		&\dTo>{Q\mu}\\
T'Q		&		&\rTo_{\psi}	&		&QT	\\
\end{diagram}
\diagspace
\begin{diagram}[height=2em]
Q		&\rEquals	&Q		\\
\dTo<{\eta' Q}	&		&\dTo>{Q \eta}	\\
T'Q		&\rTo_{\psi}	&QT		\\
\end{diagram}
\]
commute.  If $(\Eee, T) \goby{(\twid{Q}, \twid{\psi})} (\Eee', T')$ is
another lax map of monads then a \demph{transformation}%
\index{monad!transformation of}%
\index{transformation!monads@of monads}
$(Q, \psi) \go
(\twid{Q}, \twid{\psi})$ is a natural transformation $Q \goby{\alpha}
\twid{Q}$ such that 
\[
\begin{diagram}[size=2em]
T'Q		&\rTo^\psi		&QT		\\
\dTo<{T'\alpha}	&			&\dTo>{\alpha T}\\
T'\twid{Q}	&\rTo_{\twid{\psi}}	&\twid{Q}T	\\
\end{diagram}
\]
commutes.  There is a strict 2-category%
\index{monad!two-category of@2-category of}
$\fcat{Mnd}_\mr{lax}$%
\glo{Mndlax}
whose 0-cells
are pairs \pr{\Eee}{T}, whose 1-cells are lax maps of monad, and whose
2-cells are transformations.

Dually, if $T$ and $T'$ are monads on categories $\Eee$ and $\Eee'$
respectively then a \demph{colax map of monads} $(\Eee, T) \go (\Eee', T')$
consists of a functor $P: \Eee \go \Eee'$ together with a natural
transformation 
\[
\begin{diagram}
\Eee		&\rTo^T		&\Eee		\\
\dTo<P		&\swnt\tcs{\phi}&\dTo>P		\\
\Eee'		&\rTo_{T'}	&\Eee'		\\
\end{diagram}
\]
satisfying axioms dual to those for lax maps.  With the accompanying notion
of \demph{transformation} between colax maps of monads, we obtain another
strict 2-category $\fcat{Mnd}_\mr{colax}$.

A \demph{weak map of monads} is a lax map $(Q, \psi)$ of monads in which
$\psi$ is an isomorphism (or equivalently, a colax map $(P, \phi)$ in which
$\phi$ is an isomorphism), and a \demph{strict map of monads}%
\index{monad!map of|)}
is a lax map
$(Q, \psi)$ in which $\psi$ is the identity (and so $T'Q = QT$).

Often $\Eee=\Eee'$ and the functor $\Eee \go \Eee'$ is the identity.  A
natural transformation $\psi: T' \go T$ \demph{commutes with the monad
structures}%
\index{commutes with monad structures}
if $(\id, \psi)$ is a lax map of monads $(\Eee,T) \go
(\Eee,T')$ or, equivalently, a colax map $(\Eee,T') \go (\Eee,T)$.

A crucial property of lax maps of monads is that they induce maps between
categories of algebras:%
\index{algebra!monad@for monad!induced functor}
$(Q, \psi): (\Eee, T) \go (\Eee', T')$ induces the
functor
\[
\begin{array}{rrcl}
Q_* = (Q, \psi)_*: 	&\Eee^T 	&\go 	&\Eee'^{T'},	\\%
\glo{mndmapindftr}
&&\\
			&\bktdvslob{TX}{h}{X}
			&\goesto&
\left(
\begin{diagram}[height=1.5em,scriptlabels]
T'QX		\\
\dTo>{\psi_X}	\\
QTX		\\
\dTo>{Qh}	\\
QX		\\
\end{diagram}
\right).
\end{array}
\]
(Dually, but less usefully for us, a colax map of monads $(P, \phi): (\Eee,
T) \go (\Eee', T')$ induces a functor $\Eee_T \go \Eee'_{T'}$ between
Kleisli categories.)  In fact we have:
\begin{lemma}	\lbl{lemma:lax-map-mnds-is-ftr}
Let $T = (T, \mu, \eta)$ and $T' = (T',\mu',\eta')$ be monads on categories
$\Eee$ and $\Eee'$, respectively.  Then there is a one-to-one
correspondence between lax maps of monads $(\Eee,T) \go (\Eee',T')$ and
pairs $(Q,R)$ of functors such that the square 
\begin{equation}	\label{diag:lax-map-alt}
\begin{diagram}[size=2em]
\Eee^T	&\rTo^R	&\Eee'^{T'}	\\
\dTo<{\mr{forgetful}}	&	&
\dTo>{\mr{forgetful}}	\\
\Eee	&\rTo_Q	&\Eee'		\\
\end{diagram}
\end{equation}
commutes, with a lax map $(Q,\psi)$ corresponding to the pair $(Q, Q_*)$.
\end{lemma}
\begin{proof}
Given $(Q, \psi)$, the square~\bref{diag:lax-map-alt} with $R=Q_*$ plainly
commutes.  Conversely, take a pair $(Q, R)$ such
that~\bref{diag:lax-map-alt} commutes.  For each $X \in \Eee$, write $(T'QTX
\goby{\chi_X} QTX)$ for the image under $R$ of the free algebra $(T^2 X
\goby{\mu_X} TX)$, then put
\[
\psi_X = (T'QX \goby{T'Q\eta_X} T'QTX \goby{\chi_X} QTX).
\]
This defines a lax map of monads $(Q,\psi)$.  It is easily checked that the
two processes described are mutually inverse.
\done 
\end{proof}

Lax and colax maps can be related as `mates'.  Suppose we have an
adjunction
\[
\begin{diagram}[height=2em]
\cat{D}		\\
\uTo<P \ladj \dTo>Q	\\
\cat{D}'	\\
\end{diagram}
\]
and functors $\cat{D} \goby{T} \cat{D}$, $\cat{D}' \goby{T'} \cat{D}'$.
Then there is a one-to-one correspondence between natural transformations
$\phi$ and natural transformations $\psi$ with domains and codomains as
shown:
\[
\begin{diagram}
\cat{D} 	&\rTo^T		&\cat{D}	\\
\uTo<P		&\nwnt \tcs{\phi}&\uTo>P		\\
\cat{D}'	&\rTo_{T'}	&\cat{D}'	\\
\end{diagram}
\diagspace
\begin{diagram}
\cat{D} 	&\rTo^T			&\cat{D}	\\
\dTo<Q		&\nent \tcs{\psi}	&\dTo>Q		\\
\cat{D}'	&\rTo_{T'}		&\cat{D}'.	\\
\end{diagram}
\]
This is given by 
\begin{eqnarray*}
\psi	&=	&
\left(
T'Q \goby{\gamma T'Q} QPT'Q \goby{Q\phi Q} QTPQ \goby{QT\delta} QT
\right),	\\
\phi	&=	&
\left(
PT' \goby{PT'\gamma} PT'QP \goby{P\psi P} PQTP \goby{\delta TP} TP
\right)	
\end{eqnarray*}
where $\gamma$ and $\delta$ are the unit and counit of the adjunction.  We
call $\psi$ the \demph{mate}%
\index{mate}
of $\phi$ and write $\psi = \ovln{\phi}$;%
\glo{mate}
dually, we call $\phi$ the \demph{mate} of $\psi$ and write $\phi =
\ovln{\psi}$.  The world of mates is strictly monogamous: everybody has
exactly one mate, and your mate's mate is you ($\ovln{\ovln{\phi}} = \phi$,
$\ovln{\ovln{\psi}} = \psi$).  All imaginable statements about mates are
true.  In particular, if $T$ and $T'$ have the structure of monads then
$(P, \phi)$ is a colax%
\lbl{p:colax-lax-mate}
map of monads if and only if $(Q, \ovln{\phi})$ is a lax map of monads.

We finish by showing how to glue monads together.  Given two monads $(S,
\mu, \eta)$ and $(S', \mu', \eta')$ on the same category $\cat{C}$, how can
we give the composite functor $S' \of S$ the structure $(\widehat{\mu},
\widehat{\eta})$ of a monad on $\cat{C}$?  The unit is easy---
\begin{equation}	\label{eq:distrib-unit}
\widehat{\eta}
=
\left(
1
\goby{\eta' * \eta}
S' \of S
\right)
\end{equation}
---but for the multiplication we need some extra data---
\begin{equation}	\label{eq:distrib-mult}
\widehat{\mu}
=
\left(
S' \of S \of S' \of S
\goby{S' \sof ? \sof S}
S' \of S' \of S \of S
\goby{\mu' * \mu}
S' \of S
\right)
\end{equation}
---and that is provided by a distributive law.  
\begin{defn}	\lbl{defn:distrib-law}
Let $S$ and $S'$ be monads on the same category.  A \demph{distributive
law}%
\index{distributive law}
$\lambda: S \of S' \go S' \of S$ is a natural transformation such that
$(S', \lambda)$ is a lax map of monads $S \go S$ and $(S, \lambda)$ is a
colax map of monads $S' \go S'$.  
\end{defn}

\begin{lemma}	\lbl{lemma:distrib-gives-monad}
Let $\lambda: S \of S' \go S' \of S$ be a distributive law between monads
$(S, \mu, \eta)$ and $(S', \mu', \eta')$ on a category $\cat{C}$.  Then the
formulas~\bref{eq:distrib-unit} and~\bref{eq:distrib-mult} (with $? =
\lambda$) define a monad structure $(\widehat{\mu}, \widehat{\eta})$ on the
functor $S' \of S$.  If the monads $S$ and $S'$ and the transformation
$\lambda$ are all cartesian then so too is the monad $S' \of S$.  \done
\end{lemma}

The distributive law $\lambda$ determines a lax map of monads $(S',
\lambda): S \go S$, hence a functor $\twid{S}: \cat{C}^S \go \cat{C}^S$.
More incisively, we have the following.
\begin{lemma}	\lbl{lemma:distrib-corr}
Let $S$ and $S'$ be monads on a category $\cat{C}$.  Then there is a
one-to-one correspondence between distributive laws $\lambda: S \of S' \go
S' \of S$ and monads $\twid{S}$ on $\cat{C}^S$ such that the forgetful
functor $U: \cat{C}^S \go \cat{C}$ is a strict map of monads $\twid{S} \go
S'$:
\[
\begin{diagram}[size=2em]
\cat{C}^S	&\rTo^{\twid{S}}	&\cat{C}^S	\\
\dTo<U		&			&\dTo>U		\\
\cat{C}		&\rTo_{S'}		&\cat{C}.	\\
\end{diagram}
\]
If the monad $S'$ is cartesian then so too is the monad $\twid{S}$.
\end{lemma}
\begin{proof}
Straightforward, using Lemma~\ref{lemma:lax-map-mnds-is-ftr}.  
\done
\end{proof}
A distributive law $S\of S' \go S'\of S$ therefore gives two new categories
of algebras, $\cat{C}^{S' \of S}$ and $(\cat{C}^S)^{\twid{S}}$; but they
are isomorphic by general principles of coherence, or more rigorously by
\begin{lemma}	\lbl{lemma:distrib-iso-algs}
Let $S$ and $S'$ be monads on a category $\cat{C}$, let $\lambda: S \of S'
\go S' \of S$ be a distributive law, and let $\twid{S}$ be the
corresponding monad on $\cat{C}^S$.  Then there is a canonical natural
transformation
\[
\begin{diagram}
\cat{C}^S	&\rTo^{\twid{S}}	&\cat{C}^S	\\
\dTo<U		&\nent\tcs{\psi}	&\dTo>U		\\
\cat{C}		&\rTo_{S'\of S}		&\cat{C}.	\\
\end{diagram}
\]
This makes $(U, \psi)$ into a lax map of monads $\twid{S} \go S' \of S$,
and the induced functor 
$
(U, \psi)_*: (\cat{C}^S)^{\twid{S}} \go \cat{C}^{S'\of S}
$
is an isomorphism of categories.  
\end{lemma}

\begin{proof}
The transformation $\psi$ is $S'U\epsln$, where $\epsln$ is the counit of
the free-forgetful adjunction $F\ladj U$ for $S$-algebras.  Algebras for
both $\twid{S}$ and $S' \of S$ can be described as triples $(X, h, h')$
where $X \in \cat{C}$, $h$ and $h'$ are respectively $S$-algebra and
$S'$-algebra structures on $X$, and the following diagram commutes:
\[
\begin{diagram}[size=1.5em]
		&		&SS'X		&\rTo^{Sh'}	&SX	\\
		&\ldTo<{\lambda_X}&		&		&	\\
S'SX		&		&		&		&\dTo>h	\\
\dTo<{S'h}	&		&		&		&	\\
S'X		&		&\rTo_{h'}	&		&X.	\\
\end{diagram}
\]
The details of the proof are, again, straightforward.
\done
\end{proof}%
\index{monad|)}

\section{Multicategories via monads}
\lbl{sec:alt-app}%
\index{generalized multicategory!equivalent definitions of|(}

Operads are meant to be regarded as algebraic%
\index{algebraic theory}
theories of a special kind.
Monads are meant to be regarded as algebraic theories of a general kind.
It is therefore natural to ask whether operads can be re-defined as monads
satisfying certain conditions.  

We show here that the answer is nearly `yes': for any cartesian monad $T$
on a cartesian category $\Eee$, a $T$-operad is the same thing as a
cartesian monad $S$ on $\Eee$ together with a cartesian natural
transformation $S \go T$ commuting with the monad structures.  (The
transformation really is necessary, by the results of
Appendix~\ref{app:special-cart}.)  An algebra for a $T$-operad is then just
an algebra for the corresponding monad $S$.  More generally, a version
holds for $T$-multicategories, to be regarded as \emph{many-sorted}
algebraic theories.

We will need to know a little about slice%
\index{slice!category}
categories.  If $E$ is an object
of a category $\Eee$ then the forgetful functor $U_E: \Eee/E \go \Eee$%
\glo{fgtslice}
creates pullbacks in the following strict sense: if
\[
\begin{diagram}[size=2em]
P	&\rTo	&X	\\
\dTo	&	&\dTo	\\
Y	&\rTo	&Z	\\
\end{diagram}
\]
is a pullback square in $\Eee$ and $(Z \go E)$ is a map in $\Eee$ then the
evident square
\[
\begin{diagram}
\bktdvslob{P}{}{E}	&\rTo	&\bktdvslob{X}{}{E}	\\
\dTo			&	&\dTo			\\
\bktdvslob{Y}{}{E}	&\rTo	&\bktdvslob{Z}{}{E}	\\
\end{diagram}
\]
in $\Eee/E$ is also a pullback.  Hence $U_E$ reflects pullbacks, and if
$\Eee$ is cartesian then the category $\Eee/E$ and the functor $U_E$ are
also cartesian.

\begin{propn}	\lbl{propn:ind-monad-cart}
Let $T$ be a cartesian monad on a cartesian category $\Eee$, and let $C$ be
a $T$-multicategory.  Then the induced monad $T_C$ on $\Eee/C_0$ is
cartesian, and there is a cartesian natural transformation
\[
\begin{diagram}
\Eee/C_0	&\rTo^{T_C}	&\Eee/C_0	\\
\dTo<{U_{C_0}}	&\swnt \tcs{\pi^C}	&\dTo>{U_{C_0}}	\\
\Eee		&\rTo_T		&\Eee		\\	
\end{diagram}
\]
such that $(U_{C_0}, \pi^C)$ is a colax map of monads $(\Eee/C_0, T_C) \go
(\Eee, T)$.  
\end{propn}
\begin{proof}
For $X = (X \goby{p} C_0) \in \Eee/C_0$, let $\pi^C_X$ be the map in the
diagram
\[
\begin{slopeydiag}
	&	&T_C X\Spbk&	&	&	&	\\
	&\ldTo<{\pi^C_X}&	&\rdTo	&	&	&	\\
TX	&	&	&	&C_1	&	&	\\
	&\rdTo<{Tp}&	&\ldTo>\dom&	&\rdTo>\cod&	\\
	&	&TC_0	&	&	&	&C_0	\\
\end{slopeydiag}
\]
defining $T_C X$.  It is easy to check that $\pi^C$ is natural and that
$(U_{C_0}, \pi^C)$ forms a colax map of monads.  Using the Pasting
Lemma~(\ref{lemma:pasting}) it is also easy to check that $\pi^C$ is a
cartesian natural transformation, and from this that $T_C$ is a cartesian
monad.  \done
\end{proof}

The proposition says that any $(\Eee,T)$-multicategory $C$ gives rise to a
triple $(E, S, \pi) = (C_0, T_C, \pi^C)$ where
\begin{enumerate}
\item 	\lbl{item:E}
$E$ is an object of $\Eee$
\item 
$S$ is a cartesian monad on $\Eee/E$
\item 	\lbl{item:pi}
$\pi$ is a cartesian natural transformation 
\[
\begin{diagram}
\Eee/E		&\rTo^{S}	&\Eee/E		\\
\dTo<{U_E}	&\swnt \tcs{\pi}&\dTo>{U_E}	\\
\Eee		&\rTo_T		&\Eee		\\
\end{diagram}
\]
such that $(U_E, \pi)$ is a colax map of monads $S \go T$.
\end{enumerate}
It turns out that this captures exactly what a $T$-multicategory is: every
triple $(E, S, \pi)$ satisfying these three conditions arises from a
$T$-multicategory, and the whole multicategory structure of $C$ can be
recovered from the associated triple $(C_0, T_C, \pi^C)$.  We will prove
this in a moment.  In the case of $T$-operads, assuming that $\Eee$ has a
terminal object $1$, this says that a $T$-operad is a pair $(S, \pi)$ where
$S$ is a cartesian monad on $\Eee$ and $\pi: S \go T$ is a cartesian
natural transformation commuting with the monad structures---just as
promised in the introduction to this section.

In particular, the monad $T_P$ on $\Set$ induced by a plain operad $P$ is
always cartesian, which gives a large class of cartesian monads on $\Set$.
(Not all cartesian monads on $\Set$ are of this type, but all strongly
regular theories are: see Appendix~\ref{app:special-cart}.)  The natural
transformation $\pi^P$ induces the obvious functor from $\fcat{Monoid}$ to
$\Alg(P)$.

We now give the alternative definition of $T$-multicategory.  Notation: if
$E \goby{e} \twid{E}$ is a map in a category $\Eee$ then $e_!$%
\glo{compbang}
is the
functor $\Eee/E \go \Eee/\twid{E}$ given by composition with $e$.

\begin{defn}	\lbl{defn:T-mti-colax}
Let $T$ be a cartesian monad on a cartesian category $\Eee$.  Define a
category $T\hyph\Multicat'$ as follows:
\begin{description}
\item[objects] are triples $(E, S, \pi)$ as
in~\bref{item:E}--\bref{item:pi} above
\item[maps] $(E, S, \pi) \go (\twid{E}, \twid{S}, \twid{\pi})$ are pairs
$(e, \phi)$ where $e: E \go \twid{E}$ is a map in $\Eee$ and 
\[
\begin{diagram}
\Eee/E		&\rTo^S			&\Eee/E		\\
\dTo<{e_!}	&\swnt \tcs{\phi}	&\dTo>{e_!}	\\
\Eee/\twid{E}	&\rTo_{\twid{S}}	&\Eee/\twid{E}	\\
\end{diagram}
\]
is a natural transformation such that $(e_!, \phi)$ is a colax map of
monads $S \go \twid{S}$ and 
\begin{equation}	\label{eq:colax-compat}
\begin{diagram}
\Eee/E			&\rTo^S			&\Eee/E		\\
\dTo<{e_!}		&\swnt \tcs{\phi}	&\dTo>{e_!}	\\
\Eee/\twid{E}		&\rTo~{\twid{S}}	&\Eee/\twid{E}	\\
\dTo<{U_{\twid{E}}}	&\swnt \tcs{\twid{\pi}}	&\dTo>{U_{\twid{E}}}\\
\Eee			&\rTo_T			&\Eee		\\
\end{diagram}
\diagspace
=
\diagspace
\begin{diagram}
\Eee/E			&\rTo^S			&\Eee/E		\\
\dTo<{U_E}		&\swnt \tcs{\pi}	&\dTo>{U_E}	\\
\Eee			&\rTo_T			&\Eee.		\\
\end{diagram}
\end{equation}
\end{description}
\end{defn}
The natural transformation $\phi$ in the definition of map is automatically
cartesian: this follows from equation~\bref{eq:colax-compat}, the fact that
$\pi$ and $\twid{\pi}$ are cartesian, the fact that $U_{\twid{E}}$ reflects
pullbacks, and the Pasting Lemma~(\ref{lemma:pasting}).  The functor $e_!$
is also cartesian.

\begin{propn}	\lbl{propn:T-mti-colax}
For any cartesian monad $T$ on a cartesian category $\Eee$, there is an
equivalence of categories
\[
T\hyph\Multicat \eqv T\hyph\Multicat'.
\]
\end{propn}

\begin{proof}
A $T$-multicategory $C$ gives rise to an object $(C_0, T_C, \pi^C)$ of
$T\hyph\Multicat'$, as in Proposition~\ref{propn:ind-monad-cart}.
Conversely, take an object $(E, S, \pi)$ of $T\hyph\Multicat'$, and define
\begin{eqnarray*}
C_0			&=	&E,					\\
(C_1 \goby{\cod} C_0)	&=	&S(C_0 \goby{1} C_0),			\\
\dom			&= 	&\pi_{1_{C_0}}: C_1 \go TC_0.
% \dom			&= 	&\pi_{(C_0 \goby{1} C_0)}: C_1 \go TC_0.
\end{eqnarray*}
This specifies a $T$-graph $C$, and there is a $T$-multicategory structure
on $C$ given by
\[
\comp = \mu_{1_{C_0}},
\diagspace
\ids = \eta_{1_{C_0}}.
\]
(It takes a little work to see that these make sense.)  The associativity
and identity axioms for the multicategory follow from the coherence axioms
for the colax map of monads $(U_E, \pi)$.  It is straightforward to check
that this extends to an equivalence of categories.
\done
\end{proof}

\begin{cor}	\lbl{cor:T-opd-colax}
Let $T$ be a cartesian monad on a category $\Eee$ with finite
limits.  Then $T\hyph\Operad$ is equivalent to the category
$T\hyph\Operad'$ in which
\begin{itemize}
\item an object is a cartesian monad $S$ on $\Eee$ together with a
cartesian natural transformation $S \go T$ commuting with the monad
structures
\item a map $(S, \pi) \go (\twid{S}, \twid{\pi})$ is a cartesian natural
transformation $\phi: S \go \twid{S}$ commuting with the monad structures
and satisfying $\twid{\pi} \of \phi = \pi$.
\end{itemize}
\end{cor}
\begin{proof}
Restrict the proof of~\ref{propn:T-mti-colax} to the case $C_0=E=1$.
\done
\end{proof}

Since the monad $T_C$ arising from a $T$-multicategory $C$ is cartesian, it
makes sense to ask what $T_C$-multicategories are.  The answer is simple:
\begin{cor}	\lbl{cor:TC-mti}
Let $T$ be a cartesian monad on a cartesian category $\Eee$, and let $C$ be
a $T$-multicategory.  Then there is an equivalence of categories
\[
T_C\hyph\Multicat \eqv T\hyph\Multicat/C,
\]
and if a $T_C$-multicategory $D$ corresponds to a $T$-multicategory
$\ovln{D}$ over $C$ then $\Alg(D) \iso \Alg(\ovln{D})$.
\end{cor}

\begin{proof}
For the first part, it is enough to show that 
\[
T_C\hyph\Multicat' 
\eqv
T\hyph\Multicat'/(C_0, T_C, \pi^C).
\]
An object of the right-hand side is an object $(E, S, \pi)$ of
$T\hyph\Multicat'$ together with a map $E \goby{e} C_0$ and a cartesian
natural transformation $\phi$ such that $(e_!, \phi)$ is a colax map of
monads and $\pi$ is the pasting of $\phi$ and $\pi^C$.  In other words, it
is just an object $(E \goby{e} C_0, S, \phi)$ of $T_C\hyph\Multicat'$.  The
first part follows; and for the second part, $\Alg(D) \iso (\Eee/E)^S \iso
\Alg(\ovln{D})$. \done
\end{proof}

This explains many of the examples in~\ref{sec:om} and~\ref{sec:algs}
(generalized multicategories and their algebras).  Take, for
instance,~\ref{eg:M-times-mti} and~\ref{eg:M-times-alg}, where we fixed a
monoid $M$ and considered $(M\times\dashbk)$-multicategories and their
algebras.  Write $T$ for the free monoid monad on $\Set$, and let $P$ be
the $T$-operad ($=$ plain operad) with $P(1)=M$ and $P(n)=\emptyset$
for $n\neq 1$.  Then $T_P = (M\times\dashbk)$, so by
Corollary~\ref{cor:TC-mti} an $(M\times\dashbk)$-multicategory is a
plain multicategory over $P$.  Evidently a plain multicategory over
$P$ can only have unary arrows, so in fact we have
\[
(M\times\dashbk)\hyph\Multicat \eqv \Cat/M.  
\]
Moreover, the second part of~\ref{cor:TC-mti} tells us that if an
$(M\times\dashbk)$-multicategory corresponds to an object $(C \goby{\phi}
M)$ of $\Cat/M$ then its category of algebras is just $\Alg(C) \eqv
\ftrcat{C}{\Set}$, as claimed in~\ref{eg:M-times-alg}.

This reformulation of $T$-multicategories in terms of monads and colax maps
between them has a dual, using \emph{lax} maps.  As we saw
in~\ref{sec:more-monads}, colax and lax maps can be related using mates.
The adjunctions involved are
\begin{equation}	\label{eq:two-adjns}
\begin{diagram}[height=2em]
\Eee		\\
\uTo<{U_E} \ladj \dTo>{(\dashbk\times E)}	\\
\Eee/E,		\\
\end{diagram}
\diagspace \diagspace
\begin{diagram}[height=2em]
\Eee/\twid{E}		\\
\uTo<{e_!} \ladj \dTo>{e^*}	\\
\Eee/E.			\\
\end{diagram}%
\index{slice!category}
\end{equation}
In the first adjunction $E$ is an object of a category $\Eee$ with finite
limits, and the right adjoint to the forgetful functor $U_E$ sends $X \in
\Eee$ to $(X\times E \goby{\mr{pr}_2} E) \in \Eee/E$.  In the second $e: E
\go \twid{E}$ is a map in $\Eee$, and the right adjoint to $e_!$ is the
functor $e^*$ defined by pullback along $e$.  Actually, the first
is just the second in the case $\twid{E} = 1$.

\begin{defn}
Let $T$ be a cartesian monad on a category $\Eee$ with finite limits.  The
category $T\hyph\Multicat''$ is defined as follows:
\begin{description}
\item[objects] are triples $(E, S, \rho)$ where $E\in \Eee$, $S$ is a
cartesian monad on $\Eee/E$, and $\rho$ is a cartesian natural
transformation such that $(\dashbk\times E, \rho)$ is a lax map of monads
$T \go S$
\item[maps] $(E, S, \rho) \go (\twid{E}, \twid{S}, \twid{\rho})$ are pairs
$(e, \psi)$ where $e: E \go \twid{E}$ in $\Eee$ and $\psi$ is a cartesian
natural transformation such that $(e^*, \psi)$ is a lax map of monads
$\twid{S} \go S$ and an equation dual to~\bref{eq:colax-compat} in
Definition~\ref{defn:T-mti-colax} holds.
\end{description}
\end{defn}

\begin{propn}	%\lbl{propn:T-mti-lax}
For any cartesian monad $T$ on a category $\Eee$ with finite limits, there
is an isomorphism of categories
\[
T\hyph\Multicat'' \iso T\hyph\Multicat'.
\]
\end{propn}

\begin{proof}
Just take mates%
\index{mate}
throughout.  All the categories, functors and natural
transformations involved in the adjunctions~\bref{eq:two-adjns} are
cartesian, so under these adjunctions, the mate of a cartesian natural
transformation is also cartesian.  \done
\end{proof}

It follows that $T\hyph\Multicat'' \eqv T\hyph\Multicat$.  So given any
$T$-multicategory $C$, there is a corresponding lax map of monads
$(\dashbk\times C_0, \rho^C): T \go T_C$, and this induces a functor
$\Eee^T \go \Alg(C)$.  For instance, any monoid $M$ yields an algebra
for any plain multicategory $C$; concretely, this algebra $X$ is given
by $X(a)=M$ for all objects $a$ of $C$ and by $X(\theta) = $ ($n$-fold
multiplication) for all $n$-ary maps $\theta$ in $C$.  Similarly, any map
$C \go \twid{C}$ of $T$-multicategories corresponds to a lax map of monads
$T_{\twid{C}} \go T_C$, and this induces the functor $\Alg(\twid{C}) \go
\Alg(C)$ that we constructed directly at the end of
Chapter~\ref{ch:gom-basics}.%
\index{generalized multicategory!equivalent definitions of|)}

\section{Algebras via fibrations}
\lbl{sec:alg-fibs}%
\index{generalized multicategory!algebra for|(}%
\index{algebra!generalized multicategory@for generalized multicategory|(}

For a small category $C$, the functor category \ftrcat{C}{\Set} is
equivalent to the category of discrete opfibrations over $C$, as we saw
in~\ref{sec:cats}.  Here we extend this result from categories to
generalized multicategories.

By definition (p.~\pageref{p:defn-cl-d-opfib}), a functor $g: D \go C$
between ordinary categories is a discrete opfibration if for each object
$b$ of $D$ and arrow $g(b) \goby{\theta} a$ in $C$, there is a unique arrow
$b \goby{\chi} b'$ in $D$ such that $g(\chi) = \theta$.  Another way of
saying this is that in the diagram
\[
\begin{slopeydiag}
	&		&D_1		&		&	\\
	&\ldTo<{\dom}	&		&\rdTo>{\cod}	&	\\
D_0	&		&\dTo>{g_1}	&		&D_0	\\
\dTo<{g_0}&		&C_1		&		&\dTo>{g_0}\\
	&\ldTo<{\dom}	&		&\rdTo>{\cod}	&	\\
C_0	&		&		&		&C_0	\\
\end{slopeydiag}
\]
depicting $g$, the left-hand `square' is a pullback. 

Generalizing to any cartesian monad $T$ on any cartesian category $\Eee$,
let us call a map $D\goby{g}C$ of $T$-multicategories a \demph{discrete
opfibration}%
\index{fibration!discrete opfibration}
if the square
\begin{diagram}[size=2em,scriptlabels]
TD_{0}		&\lTo^{\dom}	&D_1		\\
\dTo<{Tg_0}	&		&\dTo>{g_1}	\\
TC_{0}		&\lTo^{\dom}	&C_1		\\
\end{diagram}
is a pullback.  We obtain, for any $T$-multicategory $C$, the category
$\fcat{DOpfib}(C)$%
\glo{DOpfibgen}
of discrete opfibrations over $C$: an object is a
discrete opfibration with codomain $C$, and a map from $(D \goby{g} C)$ to
$(D' \goby{g'} C)$ is a map $D \goby{f} D'$ of $T$-multicategories such
that $g' \of f = g$.  Such an $f$ is automatically a discrete opfibration too,
by the Pasting Lemma~(\ref{lemma:pasting}).

\begin{thm}	\lbl{thm:alt-alg}
Let $T$ be a cartesian monad on a cartesian category $\Eee$, and let $C$ be
a $T$-multicategory. Then there is an equivalence of categories
\[
\fcat{DOpfib}(C) \eqv \Alg(C).
\]
\end{thm}

\begin{proof}
A $C$-algebra is an algebra for the monad $T_C$ on $\Eee/C_0$, which sends
an object $X = (X\goby{p} C_0)$ of $\Eee/C_0$ to the boxed composite in the
diagram 
\[
\setlength{\unitlength}{1em}
\begin{picture}(13,9.5)
\cell{0}{0}{bl}{%
\begin{slopeydiag}
	&	&T_C X\Spbk&	&	&	&	\\
	&\ldTo<{\pi_X}&	&\rdTo>{\nu_X}&	&	&	\\
TX	&	&	&	&C_1	&	&	\\
	&\rdTo<{Tp}&	&\ldTo>\dom&	&\rdTo>\cod&	\\
	&	&TC_0	&	&	&	&C_0,	\\
\end{slopeydiag}}
\put(10.9,-1){\line(1,1){2}}
\put(2.4,7.5){\line(1,1){2}}
\put(10.9,-1){\line(-1,1){8.5}}
\put(12.9,1){\line(-1,1){8.5}}
\end{picture}
\]
and therefore consists of an object $(X\goby{p} C_0)$ of $\Eee/C_0$
together with a map $h: T_C X \go X$ over $C_0$, satisfying axioms.

So, given a $C$-algebra $(X \goby{p} C_0, h)$ we obtain a commutative
diagram
\[
\begin{slopeydiag}
	&		&T_C X		&		&	\\
	&\ldTo<{\pi_X}	&		&\rdTo>{h}	&	\\
TX	&		&\dTo>{\nu_X}	&		&X	\\
\dTo<{Tp}&		&C_1		&		&\dTo>{p}\\
	&\ldTo<{\dom}	&		&\rdTo>{\cod}	&	\\
TC_0	&		&		&		&C_0,	\\
\end{slopeydiag}
\]
the left-hand half of which is a pullback square.  The top part of the
diagram defines a $T$-graph $D$, and there is a map $g: D \go C$ defined by
$g_0=p$ and $g_1=\nu_X$.  With some calculation we see that $D$ is
naturally a $T$-multicategory and $g$ a map of $T$-multicategories.
So we have constructed from the $C$-algebra $X$ a discrete
opfibration over $C$.

This defines a functor from $\Alg(C)$ to $\fcat{DOpfib}(C)$, which is
easily checked to be full, faithful and essentially surjective on objects.
\done
\end{proof}

Let us look more closely at the $T$-multicategory $D$ corresponding to a
$C$-algebra $h=(X \goby{p} C_0, h)$.  Generalizing the terminology for
ordinary categories (p.~\pageref{p:defn-caty-elts}), we call $D$ the
\demph{multicategory of elements}%
\index{multicategory!elements@of elements}%
\index{generalized multicategory!elements@of elements}%
\index{slice!generalized multicategory by algebra@of generalized multicategory by algebra}
of $h$ and write $D = C/h$.%
\glo{genmtielts}
\begin{propn}	\lbl{propn:slice-multicat-fib}
Let $T$ be a cartesian monad on a cartesian category $\Eee$, let $C$ be a
$T$-multicategory, and let $h$ be a $C$-algebra.  Then there is an
isomorphism of categories
\[
\fcat{DOpfib}(C/h) \iso \fcat{DOpfib}(C)/h.
\]
\end{propn}
\begin{proof}
Follows from the definition of $\fcat{DOpfib}$, using the observation above
that maps in $\fcat{DOpfib}(C)$ are automatically discrete opfibrations.
\done
\end{proof}
Hence $\Alg(C/h) \eqv \Alg(C)/h$, generalizing
Proposition~\ref{propn:pshf-slice}.  In fact, this equivalence is an
isomorphism.  To see this, recall from~\ref{eg:alg-to-multi} the process of
slicing a monad by an algebra: for any monad $S$ on a category \cat{F} and
any $S$-algebra $k$, there is a monad $S/k$ on \cat{D} with the property
that $\cat{F}^{S/k} \iso \cat{F}^S/k$.
\begin{propn}	\lbl{propn:slice-multicat}
Let $T$ be a cartesian monad on a cartesian category $\Eee$, let $C$ be a
$T$-multicategory, and let $h$ be a $C$-algebra.  Then there is an
isomorphism of monads $T_{C/h} \iso T_C /h$ and an isomorphism of
categories $\Alg(C/h) \iso \Alg(C)/h$. 
\end{propn}
\begin{proof}
The first assertion is easily verified, and the second follows immediately.
\done
\end{proof}

As an example, let $C$ be the terminal%
\index{generalized multicategory!terminal}
$T$-multicategory $1$.  We have $T_1
\iso T$ and so $\Alg(1) \iso \Alg(T)$~(\ref{eg:alg-terminal}).  Given a
$T$-algebra $h = (TX \goby{h} X)$, we therefore obtain a $T$-multicategory
$1/h$.  Plausibly enough, this is the same as the $T$-multicategory $h^+$
of~\ref{eg:multi-alg}, with graph
\[
TX \ogby{1} TX \goby{h} X.
\]
So by the results above, $T_{h^+} = T_{1/h} \iso T/h$ and $\Alg(h^+) =
\Alg(1/h) \iso \Alg(T)/h$; compare~\ref{eg:alg-to-multi}.

We could also define \demph{opalgebras}%
\lbl{p:opalgebras}%
\index{opalgebra}
for a $T$-multicategory $C$ as discrete%
\index{fibration!discrete}
 brations over $C$: that is, as
maps $D \goby{g} C$ of $T$-multicategories such that the right-hand square
\[
\begin{diagram}[size=2em,scriptlabels]
D_1 		&\rTo^\cod	&D_0		\\
\dTo<{g_1}	&		&\dTo>{g_0}	\\
C_1		&\rTo^\cod	&C_0		\\
\end{diagram}
\]
of the diagram depicting $g$ is a pullback.  In the case of ordinary
categories $C$, an opalgebra for $C$ is a functor $C^\op \go \Set$.  In the
case of plain multicategories $C$, an opalgebra is a family $(X(a))_{a \in
C_0}$ of sets together with a function
\[
X(a) \go X(a_1) \times\cdots\times X(a_n)
\]
for each map $a_1, \ldots, a_n \go a$ in $C$, satisfying the obvious
axioms.  Note that these are different from the `coalgebras' or `right
modules' mentioned in~\ref{sec:om-further} and~\ref{sec:mmm}; we do not
discuss them any further.

\section{Algebras via endomorphisms}
\lbl{sec:endos}

An action of a monoid on a set is a homomorphism from the monoid to the
monoid of endomorphisms of the set.  A representation of a Lie algebra is a
homomorphism from it into the Lie algebra of endomorphisms of some vector
space.  An algebra for a plain operad is often defined as a map from it
into the operad of endomorphisms of some set~(\ref{eg:opd-End}).  Here we
show that algebras for generalized multicategories can be described in the
same way, assuming some mild properties of the base category $\Eee$.

First recall from~\ref{eg:alg-multi-End} what happens for plain
multicategories: given any family $(X(a))_{a\in E}$ of sets, there is an
associated plain multicategory $\END(X)$ with object-set $E$ and with
\begin{equation}	\label{eq:cl-endo}
(\END(X))(a_1, \ldots, a_n; a)
=
\Set (X(a_1) \times\cdots\times X(a_n), X(a) ),
\end{equation}
and if $C$ is a plain multicategory with object-set $C_0$ then a
$C$-algebra amounts to a family $(X(a))_{a\in C_0}$ of sets together with a
map $C \go \END(X)$ of multicategories leaving the objects unchanged.

To extend this to generalized multicategories we need to rephrase the
definition of $\END(X)$.  Let $T$ be the free monoid monad on $\Set$.
Recall that given a set $E$, a family $(X(a))_{a\in E}$ amounts to an
object $X \goby{p} E$ of $\Set/E$.  Then note that $X(a_1)
\times\cdots\times X(a_n)$ is the fibre over $(a_1, \ldots, a_n)$ in the
map $TX \goby{Tp} TE$, or equivalently that it is the fibre over $((a_1,
\ldots, a_n), a)$ in the map $TX \times E \goby{Tp \times 1} TE \times E$.
On the other hand, $X(a)$ is the fibre over $((a_1, \ldots, a_n), a)$ in
the map $TE \times X \goby{1 \times p} TE \times E$.  So if we define
objects
\begin{equation}	\label{eq:endo-graphs}
G_1(X) = \bktdvslob{TX \times E}{Tp \times 1}{TE \times E},
\diagspace
G_2(X) = \bktdvslob{TE \times X}{1 \times p}{TE \times E}
\end{equation}
of the category $\Set/(TE \times E)$, then~\bref{eq:cl-endo} says that the
underlying $T$-graph of the multicategory $\END(X)$ is the exponential
$G_2(X)^{G_1(X)}$.

It is now clear what the definition of endomorphism%
\index{endomorphism!generalized multicategory}%
\index{generalized multicategory!endomorphism}
multicategory in the
general case must be.  Let $T$ be a cartesian monad on a cartesian category
$\Eee$.  Assume further that $\Eee$ is \demph{locally%
\index{locally cartesian closed}
cartesian closed}:
for each object $D$ of $\Eee$, the slice category $\Eee/D$ is cartesian
closed (has exponentials).  This is true when $\Eee$ is a presheaf
category, as in the majority of our examples.  (For the definition of
cartesian closed, see, for instance, Mac Lane~\cite[IV.6]{MacCWM}.  For a
more full account, including a proof that presheaf categories are cartesian
closed, see Mac Lane and Moerdijk~\cite[I.6]{MM}; our
Proposition~\ref{propn:pshf-slice} then implies that each slice is
cartesian closed.)  Given $E\in\Eee$, define functors $G_1, G_2: \Eee/E \go
\Eee/(TE \times E)$ by the formulas of~\bref{eq:endo-graphs} above.  

The short story is that for any $X \in \Eee/E$, there is a natural
$T$-multicategory structure on the $T$-graph $\END(X) = G_2(X)^{G_1(X)}$,%
\glo{Endgenmti}
and that if $C$ is any $T$-multicategory then a $C$-algebra amounts to an
object $X$ of $\Eee/C_0$ together with a map $C \go
\END(X)$ of $T$-multicategories fixing the objects.

Here is the long story.  Given $T$ and $\Eee$ as above and $E
\in \Eee$, define a functor
\[
\begin{array}{rrcl}%
\glo{Homgenmti}
\HOM:	&(\Eee/E)^\op \times \Eee/E	&\go	&\Eee/(TE \times E)	\\
	&(X, Y)				&\goesto&G_2(Y)^{G_1(X)}.	
\end{array}
\]
% where the image is an exponential in the category $\Eee/(TE \times E)$.
Consider also the functor
\[
\Eee/(TE \times E) \times \Eee/E
\go
\Eee/E
\]
sending a pair
\[
C = 
\left(
\begin{diagram}[size=1.7em,noPS]
	&	&C_1	&	&	\\
	&\ldTo<\dom&	&\rdTo>\cod&	\\
TE	&	&	&	&E	\\
\end{diagram}
\right),
\diagspace
X =
\bktdvslob{X}{p}{E}
\]
to the boxed diagonal in the pullback diagram
\[
\setlength{\unitlength}{1em}
\begin{picture}(13,9.5)
\cell{0}{0}{bl}{%
\begin{slopeydiag}
	&	&C\of X \Spbk &	&	&	&	\\
	&\ldTo<{\pi_X}&	&\rdTo>{\nu_X}&	&	&	\\
TX	&	&	&	&C_1	&	&	\\
	&\rdTo<{Tp}&	&\ldTo>\dom&	&\rdTo>\cod&	\\
	&	&TE	&	&	&	&E.	\\
\end{slopeydiag}}
\put(10.9,-1){\line(1,1){2}}
\put(2.4,7.5){\line(1,1){2}}
\put(10.9,-1){\line(-1,1){8.5}}
\put(12.9,1){\line(-1,1){8.5}}
\end{picture}
\]
This functor is written more shortly as $(C, X) \goesto C\of X$; of course,
when $C$ has the structure of a $T$-multicategory, we usually write $C\of
\dashbk$ as $T_C$.

\begin{propn}	\lbl{propn:endo-hom-adjn}
Let $T$ be a cartesian monad on a cartesian, locally cartesian closed
category $\Eee$, and let $E$ be an object of $\Eee$.  Then there is an
isomorphism
\begin{equation}	\label{eq:Hom-adjn}
\frac{\Eee}{TE \times E} (C, \HOM(X,Y))
\ \iso\ 
\frac{\Eee}{E} (C\of X, Y)
% (\Eee/(TE \times E)) (C, \HOM(X,Y))
% \iso
% (\Eee/E) (C\of X, Y)
\end{equation}
natural in $C \in \Eee/(TE \times E)$ and $X, Y \in \Eee/E$.
\end{propn}

\begin{proof}
Write $X = (X \goby{p} E)$ and $Y = (Y \goby{q} E)$.  Product in $\Eee/(TE
\times E)$ is pullback over $TE \times E$ in $\Eee$, so the left-hand side
of~\bref{eq:Hom-adjn} is naturally isomorphic to
\[
\frac{\Eee}{TE \times E} (C \times_{TE \times E} G_1(X), G_2(Y)),
\]
where $C \times_{TE \times E} G_1(X) \go TE \times E$ is the diagonal of
the pullback square
\[
\begin{diagram}[size=2em]
\SEpbk C \times_{TE \times E} G_1(X) 	&\rTo	&C_1			\\
\dTo				&		&\dTo>{(\dom,\cod)}	\\
TX \times E			&\rTo_{Tp\times 1}	&TE \times E.	\\
\end{diagram}
\]
But by an easy calculation, we also have a pullback square
\[
\begin{diagram}[size=2em]
\SEpbk C\of X		&\rTo^{\nu_X}		&C_1			\\
\dTo<{(\pi_X, \cod\of \nu_X)}&			&\dTo>{(\dom,\cod)}	\\
TX \times E		&\rTo_{Tp\times 1}	&TE \times E,		\\
\end{diagram}
\]
so in fact an element of the left-hand side of~\bref{eq:Hom-adjn} is a map
$C \of X \go G_2(Y)$ in $\Eee/(TE \times E)$.  This is a map $C \of X \go
TE \times Y$ in $\Eee$ such that
\[
\begin{slopeydiag}
C\of X	&	&\rTo	&	&TE \times Y	\\
	&\rdTo<{(\dom\sof\nu_X, \cod\sof\nu_X)} 
		&	&\ldTo>{1\times q}&	\\
	&	&TE \times E&	&		\\
\end{slopeydiag}
\]
commutes, and this in turn is a map $C\of X \go Y$ in $\Eee/E$.
\done
\end{proof}

Next observe that $\Eee/(TE \times E)$%
\lbl{p:slice-monoidal}
is naturally a monoidal category: it is the full sub-bicategory of
$\Sp{\Eee}{T}$ whose only object is $E$.  Tensor product of objects
of $\Eee/(TE \times E)$ is composition $\of$ of 1-cells in $\Sp{\Eee}{T}$,
and a monoid in $\Eee/(TE \times E)$ is a $T$-multicategory $C$ with
$C_0=E$.  The functor
\[
\begin{array}{rcl}
\Eee/(TE \times E) \,\times\, \Eee/E 	&\go		&\Eee/E	\\
(C,X)					&\goesto	&C \of X\\
\end{array}
\]
then becomes an action of the monoidal category $\Eee/(TE \times E)$ on the
category $\Eee/E$, in the sense of~\ref{eg:mon-cat-action}: there are
coherent natural isomorphisms
\[
D \of (C \of X) \goiso (D\of C)\of X,
\diagspace
X \goiso 1_E \of X
\]
for $C, D \in \Eee/(TE \times E)$, $X \in \Eee/E$. 

\begin{propn}
Let $T$, $\Eee$ and $E$ be as in Proposition~\ref{propn:endo-hom-adjn}.
For each $X \in \Eee/E$, the $T$-graph $\END(X) = \HOM(X,X)$ naturally has
the structure of a $T$-multicategory.
\end{propn}

\begin{proof}
We have to define a composition map $\END(X) \of \END(X) \go \END(X)$.
First let $\mr{ev}_X: \END(X) \of X \go X$ be the map corresponding under
Proposition~\ref{propn:endo-hom-adjn} to the identity $\END(X) \go
\HOM(X,X)$.  Then define composition to be the map corresponding
under~\ref{propn:endo-hom-adjn} to the composite
\[
\END(X) \of \END(X) \of X
\goby{1 * \mr{ev}_X}
\END(X) \of X
\goby{\mr{ev}_X}
X.
\]
The definition of identities is similar but easier.  The associativity and
identity axioms follow from the axioms for an action of a monoidal
category~(\ref{eg:mon-cat-action}).  
\done
\end{proof}

We can now express the alternative definition of algebra.  Given $\Eee$ and
$T$ as above and a $T$-multicategory $C$, let $\Alg'(C)$ be the category in
which
\begin{description}
\item[objects] are pairs $(X, h)$ where $X \in \Eee/C_0$ and $h:
C \go \END(X)$ is a homomorphism of monoids in $\Eee/(TC_0 \times C_0)$
\item[maps] $(X, h) \go (Y, k)$ are maps $f: X \go Y$ in
$\Eee/C_0$ such that 
\[
\begin{diagram}[size=2em]
C		&\rTo^{h}		&\HOM(X,X)	\\
\dTo<{k}	&			&\dTo>{\HOM(1,f)}	\\
\HOM(Y,Y)       &\rTo_{\HOM(f,1)}       &\HOM(X,Y)		\\
\end{diagram}
\]
commutes.
\end{description}
(A homomorphism between monoids in $\Eee/(TC_0 \times C_0)$ is just a map
$f$ between the corresponding multicategories such that $f_0: C_0 \go C_0$
is the identity.)

% This definition is indeed the same as the usual definition of algebra:
%
\begin{thm}
Let $T$ be a cartesian monad on a cartesian, locally cartesian closed
category $\Eee$.  Let $C$ be a $T$-multicategory.  Then there is an
isomorphism of categories $\Alg'(C) \iso \Alg(C)$.
\end{thm}

\begin{proof}
Let $X \in \Eee/C_0$.  Proposition~\ref{propn:endo-hom-adjn} in the case
$Y=X$ gives a bijection between $T$-graph maps $h: C \go \END(X)$ and maps
$\ovln{h}: T_C(X) = C\of X \go X$ in $\Eee/C_0$.  Under this
correspondence, $h$ is a homomorphism of monoids if and only if $\ovln{h}$
is an algebra structure on $X$.  So we have a bijection between the objects
of $\Alg'(C)$ and those of $\Alg(C)$.  The remaining checks are
straightforward.  
\done
\end{proof}%
\index{generalized multicategory!algebra for|)}%
\index{algebra!generalized multicategory@for generalized multicategory|)}

\section{Free multicategories}
\lbl{sec:free-mti}%
\index{generalized multicategory!free|(}

Any directed graph freely generates a category: objects are vertices and
maps are chains of edges.  More generally, any `graph' in which each `edge'
has a finite sequence of inputs and a single output freely generates a
plain multicategory, as explained in~\ref{sec:om-further}: objects are
vertices and maps are trees of edges.  In this short section we extend this
to generalized multicategories.

The construction for plain multicategories involves an infinite recursive
process, so we cannot hope to generalize to arbitrary cartesian $\Eee$ and
$T$---after all, the category $\Eee$ being cartesian only means that it
admits certain finite limits.  There is, however, a class of cartesian
categories $\Eee$ and a class of cartesian monads $T$ for which the free
$\Cartpr$-multicategory construction is possible, the so-called
\demph{suitable}%
\index{category!suitable}%
\index{suitable}%
\index{monad!suitable}
categories and monads.  The definition of suitability is
quite complicated, but fortunately can be treated as a black box: all the
properties of suitable monads that we need are stated in this section, and
the details are confined to Appendix~\ref{app:free-mti}.

First, we have a good stock of suitable categories and monads:
\begin{thm}	\lbl{thm:free-gen}
Any presheaf category is suitable.  Any finitary cartesian monad on a
cartesian category is suitable.
\end{thm}
A functor is said to be \demph{finitary}%
\index{finitary}%
\index{functor!finitary}
if it preserves filtered colimits%
\index{colimit!filtered}
(themselves defined in Mac Lane~\cite[IX.1]{MacCWM}); a monad $(T, \mu,
\eta)$ is said to be \demph{finitary}%
\index{monad!finitary}
if the functor $T$ is finitary.
In almost all of the examples in this book, $T$ is a finitary monad on a
presheaf category.

Second, suitability is a sufficient condition for the existence of free
multicategories:
\begin{thm}	\lbl{thm:free-main}
Let $T$ be a suitable monad on a suitable category $\Eee$.  Then the
forgetful functor
\[
\Cartpr\hyph\Multicat \go \Eee^+ = \Cartpr\hyph\Graph%
\glo{plusofcaty}
\]
has a left adjoint, the adjunction is monadic, and if $T^+$%
\glo{plusofmonad}
is the induced
monad on $\Eee^+$ then both $T^+$ and $\Eee^+$ are suitable.
\end{thm}

\begin{example}	\lbl{eg:fc-cart}
The category of sets and the identity monad are
suitable~(\ref{thm:free-gen}).  In this case Theorem~\ref{thm:free-main}
tells us that there is a free category monad $\fc$ on the category of
directed graphs, and that it is suitable.  In particular it is cartesian, so
it makes sense to talk about $\fc$-multicategories,%
\index{fc-multicategory@$\fc$-multicategory}%
\index{category!free (fc)@free ($\fc$)}
as we did in
Chapter~\ref{ch:fcm}.
\end{example}

Taking the free category on a directed graph leaves the set of objects
(vertices) unchanged, and the corresponding fact for generalized
multicategories is expressed in a variant of the theorem.  Notation: if $E$
is an object of $\Eee$ then $\Cartpr\hyph\Multicat_E$%
\glo{Multicatfixedobjs}
is the subcategory of
$\Cartpr\hyph\Multicat$ whose objects $C$ satisfy $C_0=E$ and whose
morphisms $f$ satisfy $f_0=1_E$.  Observe that $\Eee/(TE \times E)$ is the
category of $T$-graphs with fixed object-of-objects $E$.
\begin{thm}	\lbl{thm:free-fixed}
Let $T$ be a suitable monad on a suitable category $\Eee$, and let $E \in
\Eee$.  Then the forgetful functor
\[
\Cartpr\hyph\Multicat_E \go \Eee^+_E = \Eee/(TE \times E)%
\glo{pluscatyfo}
\]
has a left adjoint, the adjunction is monadic, and if $T^+_E$%
\glo{plusmonadfo}
is the
induced monad on $\Eee^+_E$ then both $T^+_E$ and $\Eee^+_E$ are suitable.
\end{thm}

\begin{example}	\lbl{eg:free-cl-opd-cart}
The free monoid monad on the category of sets is suitable,
by~\ref{thm:free-gen}, hence the free%
\index{operad!free}
plain operad monad on $\Set/\nat$ is
also suitable, by~\ref{thm:free-fixed}.  In particular it is cartesian, as
claimed in~\ref{eg:mon-free-cl-opd}.  
\end{example}

In both examples the theorems were used to establish that $T^+$ and
$\Eee^+$, or $T^+_E$ and $\Eee^+_E$, were cartesian (rather than suitable).
We will use the full iterative strength in~\ref{sec:opetopes} to construct
the `opetopes'.

For technical purposes later on, we will need a refined version of these
results.  Wide pullbacks are defined on p.~\pageref{p:defn-wide-pb}.
\begin{propn}	\lbl{propn:free-refined}
If $\Eee$ is a presheaf category and the functor $T$ preserves wide
pullbacks then the same is true of $\Eee^+$ and $T^+$ in
Theorem~\ref{thm:free-main}, and of $\Eee^+_E$ and $T^+_E$ in
Theorem~\ref{thm:free-fixed}.  Moreover, if $T$ is finitary then so are
$T^+$ and $T^+_E$.
\end{propn}%
\index{generalized multicategory!free|)}

\section{Structured categories}
\lbl{sec:struc}%
\index{monoidal category!multicategory@\vs.\ multicategory|(}%
\index{structured category!generalized multicategory@\vs.\ generalized multicategory|(}%
\index{generalized multicategory!structured category@\vs.\ structured category|(}

Any monoidal category has an underlying plain multicategory.  Here we
meet `$T$-structured categories', for any cartesian monad $T$, which
bear the same relation to $T$-multicategories as strict monoidal categories
do to plain multicategories.  At the end we briefly consider the non-strict
case.

A strict monoidal category is a monoid in $\Cat$, or, equivalently, a
category in $\fcat{Monoid}$.  This makes sense because the category
$\fcat{Monoid}$ is cartesian.  More generally, if $T$ is a cartesian monad
on a cartesian category $\Eee$ then the category $\Eee^T$ of algebras is
also cartesian (since the forgetful functor $\Eee^T \go \Eee$ creates
limits), so the following definition makes sense:
\begin{defn}
Let $T$ be a cartesian monad on a cartesian category $\Eee$.  Then a
\demph{$T$-structured category}%
\index{structured category}
is a category in $\Eee^T$, and we write
$T\hyph\Struc$%
\glo{Struc}
or $\Cartpr\hyph\Struc$ for the category $\Cat(\Eee^T)$ of
$T$-structured categories.
\end{defn}
A $T$-structured category is, incidentally, a generalized multicategory:
\[
T\hyph\Struc \iso (\Eee^T, \id)\hyph\Multicat
\]
where $\id$ is the identity monad.

\begin{example}
For any cartesian category $\Eee$ we have 
\[
(\Eee,\id)\hyph\Struc \iso (\Eee,\id)\hyph\Multicat \iso \Cat(\Eee).
\]
\end{example}

\begin{example}
If $T$ is the free monoid monad on the category $\Eee$ of sets then
$T\hyph\Struc$ is the category $\fcat{StrMonCat}_\mr{str}$ of strict%
\index{monoidal category!strict}
monoidal categories and strict monoidal functors.
\end{example}

For an alternative definition, lift $T$ to a monad $\Cat(T)$ on
$\Cat(\Eee)$, then define a $T$-structured category as an algebra for
$\Cat(T)$.  This is equivalent; more precisely, there is an isomorphism of
categories
\[
\Cat(\Eee^T) \iso \Cat(\Eee)^{\Cat(T)}.
\]
In the plain case $\Cat(T)$ is the free strict monoidal category monad on
$\Cat(\Set) = \Cat$, so algebras for $\Cat(T)$ are certainly the same as
$T$-structured categories.

\begin{example}	\lbl{eg:struc-sr}
Let $T$ be the monad on $\Set$ corresponding to a strongly regular
algebraic theory, as in~\ref{eg:mon-CJ}.  It makes sense to take models of
such a theory in any category possessing finite products.  A $T$-structured
category is an algebra for $\Cat(T)$, which is merely a model of the theory
in $\Cat$.
\end{example}

\begin{example}	\lbl{eg:struc-pointed}
A specific instance is the monad $T = (1 + \dashbk)$ on $\Set$
corresponding to the theory of pointed sets~(\ref{eg:mon-exceptions}).
Then a $T$-structured category is a category $A$ together with a functor
from the terminal category $1$ into $A$; in other words, it is a category
$A$ with a distinguished object.
\end{example}

\begin{example}	\lbl{eg:struc-fin-lims}
The principle described in~\ref{eg:struc-sr} for finite product theories
holds equally for finite limit theories.  For instance, if $T$ is the free
plain operad monad on $\Eee = \Set^\nat$ (as in~\ref{eg:mon-free-cl-opd})
then a $T$-structured category is an operad in $\Cat$, that is, a
$\Cat$-operad (p.~\pageref{p:defn-V-Operad}).  So we now have three
descriptions of $\Cat$-operads:%
\index{Cat-operad@$\Cat$-operad!three definitions of}
as operads in $\Cat$, as $S$-operads where
$S$ is the free strict monoidal category monad on
$\Cat$~(\ref{eg:mti-Cat}), and as $T$-structured categories.
\end{example}

\begin{example}
When $T$ is the free strict $\omega$-category%
\index{omega-category@$\omega$-category!strict!free}
monad on the category of
globular sets~(\ref{eg:glob-mnd}), a $T$-structured category is what has
been called a `strict monoidal%
\index{monoidal globular category}
globular category' (Street~\cite[\S
1]{StrRMB} or Batanin~\cite[\S 2]{BatMGC}).%
\index{Batanin, Michael}
\end{example}

\begin{example}
Take the free category%
\index{category!free (fc)@free ($\fc$)}
monad $\fc$ on the category $\Eee$ of directed
graphs (Chapter~\ref{ch:fcm}).  Then an $\fc$-structured category is a
category in $\Eee^\fc \iso \Cat$, that is, a strict double%
\index{double category!strict}
category.
\end{example}

Every strict monoidal category $A$ has an underlying plain multicategory
$UA$, and every plain multicategory $C$ generates a free strict monoidal
category $FC$, giving an adjunction $F\ladj U$~(\ref{sec:om-further}).  The
same applies for generalized multicategories.  Any $T$-structured category
$A$ has an underlying $T$-multicategory $UA$, whose graph is given by
composing along the upper slopes of the pullback diagram
\[
\begin{slopeydiag}
	&		&(UA)_1\Spbk&	&	&	&	\\
	&\ldTo		&	&\rdTo	&	&	&	\\
TA_0	&		&	&	&A_1	&	&	\\
	&\rdTo<{h_0}	&	&\ldTo<\dom&	&\rdTo>\cod&	\\
	&		&A_0	&	&	&	&A_0	\\
\end{slopeydiag}
\]
in which $h_0: TA_0 \go A_0$ is the $T$-algebra structure on $A_0$.
Conversely, the free $T$-structured category $FC$ on a $T$-multicategory
$C$ has graph
\[
\begin{slopeydiag}
	&	&	&	&TC_1	&	&	&	&	\\
	&	&	&\ldTo<{T\dom}&	&\rdTo(4,4)>{T\cod}&&	&	\\
	&	&T^2 C_0&	&	&	&	&	&	\\
	&\ldTo<{\mu_{C_0}}&&	&	&	&	&	&	\\
TC_0	&	&	&	&	&	&	&	&TC_0	\\
\end{slopeydiag}
\]
and the $T$-algebra structure on $(FC)_i = TC_i$ is $\mu_{C_i}$ ($i=0,1$).
We then have the desired adjunction
\begin{equation}	\label{eq:struc-mti-adjn}
\begin{diagram}[height=2em]
T\hyph\Struc		\\
\uTo<F \ladj \dTo>U	\\
T\hyph\Multicat.	\\
\end{diagram}
\end{equation}

\begin{example}	\lbl{eg:free-struc-D}
Consider once more the free monoid monad $T$ on $\Eee=\Set$.  Take the
terminal plain multicategory $1$, which has graph
\[
\nat \ogby{1} \nat \goby{!} 1.
\]
Then $F1$ is a strict monoidal category with graph
\[
\nat \ogby{+} T\nat \goby{T!} \nat.
\]
The objects of $F1$ are the natural numbers and a map $m \go n$ in $F1$ is
a sequence \bftuple{m_1}{m_n} of natural numbers such that $m_1 +\cdots+
m_n = m$.  So $F1$ is the strict monoidal category $\scat{D}$%
\index{augmented simplex category $\scat{D}$}
of (possibly
empty) finite totally ordered sets, with addition as tensor and $0$ as
unit.  This is also suggested by diagram~\bref{diag:arrows-in-mon-cat}
(p.~\pageref{diag:arrows-in-mon-cat}).
\end{example}

\begin{example}	\lbl{eg:struc-disc}
Any object $K$ of a cartesian category $\cat{K}$ generates a category $DK$%
\glo{discintcat}
in $\cat{K}$, the \demph{discrete%
\index{category!discrete}
category} on $K$, uniquely determined by
its underlying graph
\[
K \ogby{1} K \goby{1} K.
\]
In particular, if $h = (TX \goby{h} X)$ is an algebra for some cartesian
monad $T$ then $Dh$ is a $T$-structured category and $UDh$ is a
$T$-multicategory with graph
\[
TX \ogby{1} TX \goby{h} X.
\]
So $UDh$ is the $T$-multicategory $h^+$%
\index{plus construction $\blank^+$}
discussed in
Example~\ref{eg:multi-alg}, and the triangle of functors
\[
\begin{diagram}[height=2em]
\Eee^T	&\rTo^{D}	&T\hyph\Struc	\\
	&\rdTo<{\blank^+}	&\dTo>{U}	\\
	&			&T\hyph\Multicat\\
\end{diagram}
\]
commutes up to natural isomorphism.  
\end{example}

I have tried to emphasize that multicategories are truly different from
monoidal categories (as well as providing a more natural language in many
situations).  This is borne out by the fact that the functor $U:
\fcat{StrMonCat}_\mr{str} \go \Multicat$ is far from an equivalence: it
is faithful, but neither full~(\ref{eg:map-mti-mon}) nor
essentially surjective on objects (p.~\pageref{p:not-ESO}).  

There is, however, a `representation%
\index{representation theorem}
theorem' saying that every
$T$-multicategory embeds fully in some $T$-structured category.  In
particular, every plain multicategory is a full sub-multicategory of the
underlying multicategory of some strict monoidal
category~(\ref{eg:multi-some-of-mon}).  I do not know of any use for the
theorem; it does not reduce the study of $T$-multicategories to the study
of $T$-structured categories any more than the Cayley%
\index{Cayley representation}
Representation
Theorem reduces the study of finite groups to the study of the symmetric
groups.  I therefore leave the proof as an exercise.  The precise statement
is as follows.
\begin{defn}
Let $T$ be a cartesian monad on a cartesian category $\Eee$.  A map $f: C
\go C'$ of $T$-multicategories is \demph{full%
\index{full and faithful}
and faithful} if the
square
\[
\begin{diagram}[size=2em]
C_1		&\rTo^{(\dom,\cod)}	&TC_0 \times C_0	\\
\dTo<{f_1}	&			&\dTo>{Tf_0 \times f_0}	\\
C'_1		&\rTo_{(\dom,\cod)}	&TC'_0 \times C'_0\\
\end{diagram}
\]
is a pullback.
\end{defn}
\begin{propn}	\lbl{propn:Cayley-multi}
Let $T$ be a cartesian monad on a cartesian category $\Eee$, and consider
the adjunction of~\bref{eq:struc-mti-adjn}
(p.~\pageref{eq:struc-mti-adjn}).  For each $T$-multicategory $C$, the unit
map $C \go UFC$ is full and faithful.  
\done
\end{propn}

We finish with two miscellaneous thoughts.

First, we have been considering generalizations of plain multicategories,
which are structures whose operations are `many in, one out'.  But
in~\ref{sec:om-further} we also considered PROs
(and their symmetric
cousins, PROPs), whose operations are `many%
\index{many in, many out}
in, many out'.  A PRO consists
of a set $S$ (the objects, or `colours', often taken to have only one
element) and a strict monoidal category whose underlying monoid of objects
is the free monoid on $S$.  The generalization to arbitrary cartesian
monads $T$ on cartesian categories $\Eee$ is clear: a \demph{$T$-PRO}%
\index{PRO!generalized}
should be defined as a pair $(S,A)$ where $S\in\Eee$ and $A$ is a
$T$-structured category whose underlying $T$-algebra of objects is the free
$T$-algebra on $S$.  We will not do anything with this definition, but
see~\ref{sec:many} for further discussion of `many in, many out'.

Second, if strict monoidal categories generalize to $T$-structured
categories, what do weak monoidal categories generalize to?  One answer
comes from realizing that the category $\Cat(\Eee)$ has the structure of a
strict 2-category and the monad $\Cat(T)$ the structure of a strict
2-monad.  We can then define a \demph{weak $T$-structured%
\index{structured category!weak}
category} to be a
weak algebra for this 2-monad, and indeed do the same in the lax case.  In
particular, if $T$ is the free plain operad monad on
$\Set^\nat$~(\ref{eg:struc-fin-lims}) then a weak $T$-structured category is
like a $\Cat$-operad,%
\index{Cat-operad@$\Cat$-operad!weak}
but with the operadic composition only obeying
associativity and unit laws up to coherent isomorphism.  For example, let
$P(n)$ be the category of Riemann%
\index{Riemann surface}%
\index{manifold!operad of}
surfaces whose boundaries are identified
with the disjoint union of $(n+1)$ copies of $S^1$, and define composition
by gluing: then $P$ forms a weak $T$-structured category.
Compare~\ref{eg:opd-Riemann} where, not having available the refined
language of generalized multicategories, we had to quotient out by
isomorphism and so lost (for instance) any information about automorphisms
of the objects.  We do not, however, pursue weak structured categories any
further in this book.%
\index{monoidal category!multicategory@\vs.\ multicategory|)}
\index{structured category!generalized multicategory@\vs.\ generalized multicategory|)}%
\index{generalized multicategory!structured category@\vs.\ structured category|)}

\section{Change of shape}
\lbl{sec:change}%
\index{change of shape|(}%
\index{generalized multicategory!change of shape|(}

To do higher-dimensional category theory we are going to want to move%
\index{n-category@$n$-category!definitions of!comparison}
between $n$-categories and $(n+1)$-categories and $\omega$-categories,
between globular and cubical and simplicial structures, and so on.  In
Chapter~\ref{ch:a-defn} we will see that a weak $n$-category can be defined
as an algebra for a certain $\gm{n}$-operad, where $\gm{n}$ is the free
strict $n$-category monad on the category of $n$-globular sets.  So if we
want to be able to relate $n$-categories to $(n+1)$-categories, for
instance, then we will need some way of relating $\gm{n}$-operads to
$\gm{n+1}$-operads and some way of relating their algebras.  In this
section we set up the supporting theory: in other words, we show what
happens to $T$-multicategories and their algebras as the monad $T$
varies.

Formally, we expect the assignment $\Cartpr \goesto \Cartpr\hyph\Multicat$
to be functorial in some way.  We saw in~\ref{sec:more-monads} that there
are notions of lax, colax, weak, and strict maps of monads, some of which
are special cases of others, and we will soon see that a map $\Cartpr \go 
(\Eee', T')$ of any one of these types induces a functor
\[
\Cartpr\hyph\Multicat \go (\Eee',T')\hyph\Multicat.
\]

First we need some terminology.  A lax map of monads $(Q, \psi): \Cartpr
\go (\Eee',T')$ is \demph{cartesian}%
\index{monad!map of!cartesian}%
\index{cartesian}
if the functor $Q$ is cartesian (but
note that the natural transformation $\psi$ need not be cartesian).
Cartesian monads, cartesian lax maps of monads, and transformations form a
sub-2-category $\fcat{CartMnd}_\mr{lax}$%
\glo{CartMndlax}
of $\fcat{Mnd}_\mr{lax}$.  Dually, a colax map of monads $(P, \phi)$ is
\demph{cartesian} if the functor $P$ \emph{and} the natural transformation
$\phi$ are cartesian.  Cartesian monads, cartesian colax maps of monads,
and transformations form a sub-2-category $\fcat{CartMnd}_\mr{colax}$ of
$\fcat{Mnd}_\mr{colax}$.

These definitions appear haphazard, with natural transformations required
to be cartesian, or not, at random.  I can justify this only pragmatically:
they are the conditions required to make the following constructions work.

The main constructions are as follows.  Let $(Q, \psi): \Cartpr \go
(\Eee',T')$ be a cartesian lax map of cartesian monads.  Then there is an
induced functor
\[
Q_* = (Q, \psi)_*: 
\Cartpr\hyph\Multicat \go (\Eee',T')\hyph\Multicat%
\glo{laxmtiindftr}
\]
sending an \Cartpr-multicategory $C$ to the $(\Eee',T')$-multicategory
$Q_* C$ whose underlying graph is given by composing along the upper
slopes of the pullback diagram
\[
\begin{diagram}[width=1.7em,height=1.7em,scriptlabels,noPS]
	&	&(Q_* C)_1\Spbk&	&	&	&	\\
	&\ldTo	&		&\rdTo	&	&	&	\\
T'Q C_0 &	&		&	&QC_1	&	&	\\
	&\rdTo<{\psi_{C_0}}&	&\ldTo<{Q\dom}&	&\rdTo>{Q\cod}&	\\
	&	&QT C_0		&	&	&	&QC_0=(Q_* C)_0\\
\end{diagram}
\]
and whose composition and identities are defined in an evident way.
Dually, let $(P, \phi): \Cartpr \go (\Eee',T')$ be a cartesian colax map of
cartesian monads.  Then there is an induced functor
\[
P_* = (P, \phi)_*: 
\Cartpr\hyph\Multicat \go (\Eee',T')\hyph\Multicat%
\glo{colaxmtiindftr}
\]
sending an \Cartpr-multicategory $C$ to the $(\Eee',T')$-multicategory
$P_* C$ with underlying graph
\[
\begin{diagram}[width=1.7em,height=1.7em,scriptlabels,noPS]
	&	&	&	&PC_1 = (P_* C)_1&&&	&	\\
	&	&	&\ldTo<{P\dom}&	&\rdTo(4,4)>{P\cod}&&	&	\\
	&	&PTC_0	&	&	&	&	&	&	\\
	&\ldTo<{\phi_{C_0}}&&	&	&	&	&	&	\\
T'PC_0	&	&	&	&	&	&	&	&
PC_0=(P_* C)_0.	\\
\end{diagram}
\]
Note that $(Q_* C)_0 = Q(C_0)$ and $(P_* C)_0 = P(C_0)$, so if $Q$ or $P$
preserves all finite limits then $Q_*$ or $P_*$ restricts to a functor
\label{p:change-operads}
\[
\Cartpr\hyph\Operad \go (\Eee',T')\hyph\Operad
\]
between categories of operads.

After filling in all the details we obtain two maps of strict 2-categories,
\[
\fcat{CartMnd_\mr{lax}} \go \fcat{CAT}, 
\diagspace
\fcat{CartMnd_\mr{colax}} \go \fcat{CAT},
\]
both defined on objects by $\Cartpr \goesto \Cartpr\hyph\Multicat$.  The
first is defined using pullbacks, so is only a weak functor; the second is
strict.  These dual functors agree where they intersect: the 2-categories
$\fcat{CartMnd}_\mr{lax}$ and $\fcat{CartMnd}_\mr{colax}$ `intersect' in
the 2-category $\fcat{CartMnd}_\mr{wk}$%
\glo{CartMndwk}
of cartesian monads, cartesian weak
maps of monads, and transformations, and the square
\[
\begin{diagram}[height=2em,width=5em]
\fcat{CartMnd}_\mr{wk}	&\rIncl		&\fcat{CartMnd}_\mr{colax}	\\
\dIncl			&		&\dTo				\\
\fcat{CartMnd}_\mr{lax}	&\rTo		&\CAT				\\
\end{diagram}
\]
commutes up to natural isomorphism.  (A \demph{cartesian%
\index{monad!map of!cartesian}
weak map} of
monads is a cartesian lax map $(Q, \psi)$ in which $\psi$ is a natural
isomorphism; then $\psi$ is automatically cartesian.)

\begin{example}	\lbl{eg:ind-ftr-int-cats}
Let $Q: \cat{D} \go \cat{D'}$ be a cartesian functor between cartesian
categories.  Then $(Q, \id): (\cat{D}, \id) \go (\cat{D'}, \id)$ is a
strict map of monads and so induces (unambiguously) a functor between the
categories of multicategories, 
\[
\left( 
(\cat{D},\id)\hyph\Multicat \goby{Q_*} (\cat{D'},\id)\hyph\Multicat
\right)
= 
\left( 
\Cat(\cat{D}) \goby{Q_*} \Cat(\cat{D'})
\right).
\]
This is the usual induced functor between categories of internal%
\index{category!internal}
categories.
\end{example}

There is another way in which the two processes are compatible, involving
adjunctions.  Suppose we have a diagram
\begin{equation}	\label{diag:adjt-change}
\begin{diagram}[height=2em]
\Cartpr			\\
\uTo<{\pr{P}{\phi}}	
\ladj			
\dTo>{\pr{Q}{\psi}}	\\
(\Eee',T')		\\
\end{diagram}
\end{equation}
in which $(P, \phi)$ is a cartesian colax map of cartesian monads, $(Q,
\psi)$ is a cartesian lax map, and there is an adjunction%
\index{monad!map of!adjunction between}
of functors
$P\ladj Q$ under which $\phi$ and $\psi$ are mates~(\ref{sec:more-monads}).
Then, as may be checked, there arises an adjunction between the functors
$P_*$ and $Q_*$ constructed above:
\begin{equation}   \label{diag:mti-adjn}
\begin{diagram}[height=2em]
\Cartpr\hyph\Multicat	\\
\uTo<{P_*}		
\ladj			
\dTo>{Q_*}		\\
(\Eee',T')\hyph\Multicat.\\
\end{diagram}
\end{equation}

\begin{example}	\lbl{eg:change-struc-adjn}
Let $T$ be a cartesian monad on a cartesian category $\Eee$.  Then there is
a diagram
\[
\begin{diagram}[height=2em]
(\Eee^T, \id)		\\
\uTo<{(F, \nu)}		
\ladj			
\dTo>{(U, \epsln)}	\\
(\Eee,T)		\\
\end{diagram}
\]
of the form~\bref{diag:adjt-change}, in which $F$ and $U$ are the free and
forgetful functors and $\nu$ and $\epsln$ are certain canonical natural
transformations.  This gives rise by the process just described to an
adjunction
\[
\begin{diagram}[height=2em]
(\Eee,T)\hyph\Struc	\\
\uTo<{F_*}		
\ladj			
\dTo>{U_*}		\\
(\Eee,T)\hyph\Multicat,	\\
\end{diagram}%
\index{structured category!generalized multicategory@\vs.\ generalized multicategory}%
\index{generalized multicategory!structured category@\vs.\ structured category}
\]
none other than the adjunction that was the subject of the previous
section.
\end{example}

\begin{example}   \lbl{eg:cart-adjts}
Let $(P, \phi): (\Eee, T) \go (\Eee', T')$ be a cartesian colax map of
cartesian monads, and suppose that the functor $P$ has a right adjoint $Q$.
Then $P_*$ has a right adjoint too: for taking the mate $\psi =
\ovln{\phi}$ gives $Q$ the structure of a cartesian lax map of monads,
% (p.~\pageref{p:colax-lax-mate}), 
leading to an adjunction $P_* \ladj Q_*$.
\end{example}

We have considered change of shape for $T$-multicategories; let us now do
the same for $T$-algebras%
\index{algebra!monad@for monad!induced functor}
and $T$-structured categories, concentrating on
lax rather than colax maps.  If $(Q, \psi): (\Eee, T) \go (\Eee', T')$ is a
cartesian lax map of cartesian monads then the induced functor $\Eee^T \go
\Eee'^{T'}$ is also cartesian, so by Example~\ref{eg:ind-ftr-int-cats}
induces in turn a functor $\Cat(\Eee^T) \go \Cat(\Eee'^{T'})$---in other
words, induces a functor on structured%
\index{structured category!change of shape}
categories,
\[
Q_*: \Cartpr\hyph\Struc \go (\Eee',T')\hyph\Struc.
\]
The change-of-shape processes for algebras, structured categories and
multicategories are compatible: for any cartesian lax map $(Q, \psi):
\Cartpr \go (\Eee',T')$ of cartesian monads, the diagram
\[
\begin{diagram}[height=2em]
\Eee^T			&\rTo^{D}	&
(\Eee,T)\hyph\Struc	&\rTo^{U_*}		&
(\Eee,T)\hyph\Multicat	\\
\dTo<{Q_*}		&			&
\dTo>{Q_*}		&			&
\dTo>{Q_*}		\\
\Eee'^{T'}		&\rTo_{D}	&
(\Eee',T')\hyph\Struc	&\rTo_{U_*}		&
(\Eee',T')\hyph\Multicat\\
\end{diagram}
\]
commutes up to natural isomorphism.  Here $D$ is the discrete category
functor defined in~\ref{eg:struc-disc}, and the commutativity of the
left-hand square can be calculated directly.  The $U_*$'s are the forgetful
functors, which, as we saw in~\ref{eg:change-struc-adjn}, are induced by
lax maps $(U,\epsln)$.  To see that the right-hand square commutes, it is
enough to see that the square of lax maps
\[
\begin{diagram}[height=2em]
(\Eee^T, \id)		&\rTo^{(U,\epsln)}	&(\Eee,T)	\\
\dTo<{(Q_*,\id)}	&			&\dTo>{(Q,\psi)}\\
(\Eee'^{T'}, \id)	&\rTo_{(U,\epsln)}	&(\Eee',T')	\\
\end{diagram}
\]
commutes, and this is straightforward.

Finally, we answer the question posed in the introduction to the section:
how do algebras%
\index{algebra!generalized multicategory@for generalized multicategory!change of shape (induced functor)}
for multicategories behave under change of shape?  An
algebra for an \Cartpr-multicategory $C$ is an object $X$ over $C_0$ acted
on by $C$, so if $(Q, \psi): \Cartpr \go (\Eee',T')$ is a cartesian lax or
colax map then we might hope that $QX$, an object over $QC_0$, would be
acted on by $Q_* C$.  In other words, we might hope for a functor from
$\Alg(C)$ to $\Alg(Q_* C)$.  Such a functor does indeed exist, in both the
lax and colax cases.  

First take a cartesian lax map of cartesian monads, $(Q, \psi): \Cartpr \go
(\Eee',T')$.  For any $T$-multicategory $C$ there is an induced lax map of
monads
\[
(\Eee/C_0, T_C) \go (\Eee'/QC_0, T'_{Q_* C}), 
\]
comprised of the functor $\Eee/C_0 \go \Eee'/QC_0$ induced by $Q$ and a
natural transformation $\psi$ that may easily be determined.  This in turn
induces a functor on categories of algebras:
\[
\begin{array}{rcl}
\Alg(C)		&\go	&\Alg(Q_* C),	\\
&&\\
\left(
\begin{diagram}[height=1.5em,scriptlabels]
T_C X	\\
\dTo>h	\\
X	\\
\end{diagram}
\right)						&
\goesto						&
\left(
\begin{diagram}[height=1.5em,scriptlabels]
T'_{Q_* C} Q X	\\
\dTo>{\psi^C_X}		\\
QT_C X			\\
\dTo>{Qh}		\\
QX			\\
\end{diagram}
\right)
.
\end{array}
\]

Now take a cartesian colax map of cartesian monads, $(P, \phi): \Cartpr \go
(\Eee',T')$.  For any $T$-multicategory $C$ there is an induced \emph{weak}
map of monads 
\[
(\Eee/C_0, T_C) \go (\Eee'/PC_0, T'_{P_* C}), 
\]
which amounts to saying that if $X = (X \goby{p} C_0) \in \Eee/C_0$ then
$T'_{P_* C} P X \iso P T_C X$ canonically; and indeed, we have a diagram
\[
\begin{diagram}[width=1.7em,height=1.7em,scriptlabels,noPS]
	&	&	&	&\Spbk PT_C X&	&	&	&	\\
	&	&	&\ldTo	&	&\rdTo	&	&	&	\\
	&	&\Spbk PTX&	&	&	&PC_1	&	&	\\
	&\ldTo<{\phi_X}&&\rdTo>{PTp}&	&\ldTo>{P\dom}&	&\rdTo>{P\cod}&	\\
T'PX	&	&	&	&PTC_0	&	&	&	&PC_0	\\
	&\rdTo<{T'Pp}&	&\ldTo>{\phi_{C_0}}&&	&	&	&	\\
	&	&T'PC_0,&	&	&	&	&	&	\\
\end{diagram}
\]
giving the isomorphism required.  

\begin{example}	\lbl{eg:change-mon-multi}
In~\ref{sec:om-further} we considered the adjunction between plain
multicategories%
\index{monoidal category!multicategory@\vs.\ multicategory}
and strict monoidal categories, and
in~\ref{eg:change-struc-adjn} we saw that it is induced by certain maps of
monads: 
\[
\begin{diagram}[height=2em]
(\fcat{Monoid}, \id)			\\
\uTo<{(F,\nu)} \ladj \dTo>{(U,\epsln)}	\\
(\Set, \textrm{free monoid})		\\
\end{diagram}
\diagspace
\goesto
\diagspace
\begin{diagram}[height=2em]
\fcat{StrMonCat}_\mr{str}	\\
\uTo<{F_*} \ladj \dTo>{U_*}	\\
\Multicat.			\\
\end{diagram}
\]
The lax map $(U, \epsln)$ induces a functor $\Alg(A) \go \Alg(U_* A)$ for
each strict monoidal category $A$.  When $A$ is regarded as a
$(\fcat{Monoid}, \id)$-multicategory, an $A$-algebra is a lax monoidal
functor $A \go \Set$.  This is the same thing as an algebra for the
underlying multicategory $U_* A$~(\ref{eg:alg-mon}), and the induced
functor is in fact an isomorphism.

Conversely, the colax map $(F, \nu)$ induces a functor $\Alg(C) \go
\Alg(F_* C)$ for each multicategory $C$, whose explicit form is left as an
exercise.
\end{example}

It is no coincidence that the functors induced by $U$ in this example are
isomorphisms.  Roughly speaking, this is because $U$ is monadic:
\begin{propn}	\lbl{propn:shape-distrib}
Let $\cat{C}$ be a cartesian category, let $S$ and $S'$ be cartesian monads
on $\cat{C}$, and let $\lambda: S\of S' \go S'\of S$ be a cartesian
distributive law.  Write
\[
\begin{diagram}
\cat{C}^S	&\rTo^{\twid{S}}	&\cat{C}^S	\\
\dTo<U		&\nent\tcs{\psi}		&\dTo>U		\\
\cat{C}		&\rTo_{S'\of S}		&\cat{C}	\\
\end{diagram}
\]
for the induced lax map of monads, as in~\ref{lemma:distrib-iso-algs}.
Then for any $\twid{S}$-multicategory $C$, there is an isomorphism of
categories $\Alg(C) \iso \Alg((U, \psi)_* C)$.
\end{propn}
The monads $\twid{S}$ and $S'\of S$ are cartesian
(\ref{lemma:distrib-gives-monad},~\ref{lemma:distrib-corr}), so it does
make sense to talk about multicategories for them.  
\begin{proof}
We just prove this in the case where $C$ is an $\twid{S}$-operad $O$, since
that is all we will need later and the proof is a little easier.

As well as the lax map shown, we have a strict map of monads
\[
\begin{diagram}
\cat{C}^S	&\rTo^{\twid{S}}	&\cat{C}^S	\\
\dTo<U		&\neeq			&\dTo>U		\\
\cat{C}		&\rTo_{S'}		&\cat{C}.	\\
\end{diagram}
\]
So we have an $S'$-operad $(U, \id)_* O$, hence a monad $S'_{(U, \id)_*
O}$ on $\cat{C}$, of which $\twid{S}_O$ is a lift to $\cat{C}^S$.  By
Lemma~\ref{lemma:distrib-corr}, there is a corresponding distributive law
\[
S \of S'_{(U, \id)_* O} 
\go 
S'_{(U, \id)_* O} \of S. 
\]
This gives $S'_{(U, \id)_* O} \of S$ the structure of a monad on $\cat{C}$;
but it can be checked that this monad is exactly $(S' \of S)_{(U, \psi)_*
O}$, so by Lemma~\ref{lemma:distrib-iso-algs},
\[
\Alg(O)
=
(\cat{C}^S)^{\twid{S}_O}
\iso
\cat{C}^{S'_{(U, \id)_* O} \of S}
\iso
\cat{C}^{(S' \of S)_{(U, \psi)_* O}}
=
\Alg((U, \psi)_* O),
\]
as required.
\done
\end{proof}
The forgetful functor in Example~\ref{eg:change-mon-multi} is the case
$\cat{C} = \Set$, $S = (\textrm{free monoid})$, $S' = \id$.  We will use
the Proposition when comparing definitions of weak 2-category
in~\ref{sec:wk-2}.%
\index{change of shape|)}%
\index{generalized multicategory!change of shape|)}

\section{Enrichment}	\lbl{sec:enr-mtis}

We finish with a brief look at a topic too large to fit in this book.  Its
slogan is `what can we enrich in?'  

Take, for example, $\Ab$-categories:%
\index{Ab-category@$\Ab$-category}
categories enriched in (or `over')
abelian groups.  The simplest definition of an $\Ab$-category is as a class
$C_0$ of objects together with an abelian group $C(a,b)$ for each $a, b \in
C_0$, a bilinear composition function
\[
C(a,b), C(b,c) \go C(a,c)
\]
for each $a,b,c \in C_0$, and an identity $1_a \in C(a,a)$ for each $a \in
C_0$, satisfying associativity and identity axioms.  You \emph{could}
express composition as a linear map out of a tensor product, but this would
be an irrelevant%
\index{monoidal category!multicategory@\vs.\ multicategory}
elaboration; put another way, the first definition of
$\Ab$-category can be understood by someone who knows what a multilinear
map is but has not yet learned about tensor products.

More generally, if $V$ is a plain multicategory then there is an evident
definition of $V$-enriched%
\index{enrichment!category@of category!plain multicategory@in plain multicategory}
category.  Classically one enriches categories
in monoidal categories (as in~\ref{sec:cl-enr}), but this is an unnaturally
narrow setting; in the terminology of~\ref{sec:non-alg-notions},
representability%
\index{multicategory!representable}
of the multicategory is an irrelevance.

More generally still, suppose we have some type of categorical
structure---`widgets',
say.  Then the question is: for what types of
structure $V$ can we make a definition of `$V$-enriched%
\index{enrichment}
widget'?  In the
previous paragraphs widgets were categories, and we saw that $V$ could be a
plain multicategory.  (In fact, that is not all $V$ can be, as we will soon
discover.)  In my~\cite{GECM} paper, the question is answered when widgets
are $T$-multicategories for almost any cartesian monad $T$.  The general
definition of enriched $T$-multicategory is short, simple, and given below,
but unwinding its implications takes more space than we have; hence the
following sketch.

Let $T$ be a monad on a category $\Eee$, and suppose that both $T$ and
$\Eee$ are suitable~(\ref{sec:free-mti}), so that there is a cartesian free
$T$-multicategory monad $T^+$ on the category $\Eee^+$ of $T$-graphs.  For
any object $C_0$ of $\Eee$, there is a unique $T$-multicategory structure
on the $T$-graph
\[
\begin{slopeydiag}
	&	&TC_0 \times C_0	&		&	\\
	&\ldTo<{\mr{pr}_1}&		&\rdTo>{\mr{pr}_2}&	\\
TC_0	&	&			&		&C_0,	\\
\end{slopeydiag}
\]
and we write this $T$-multicategory as $IC_0$,%
\glo{indiscgenmti}
the \demph{indiscrete}%
\index{generalized multicategory!indiscrete}
$T$-multicategory on $C_0$.  Then $IC_0$ is a $T^+$-algebra, so
by~\ref{eg:multi-alg} gives rise to a $T^+$-multicategory $(IC_0)^+$ whose
domain map is the identity.
\begin{defn}%
\index{enrichment!generalized multicategory@of generalized multicategory}%
\index{generalized multicategory!enriched}
Let $T$ be a suitable monad on a suitable category $\Eee$ and let $V$ be a
$T^+$-multicategory.  A \demph{$V$-enriched $T$-multicategory} is 
an object $C_0$ of $\Eee$ together with a map $(IC_0)^+ \go V$ of
$T^+$-multicategories. 
\end{defn}

\begin{example}
The most basic example is when $T$ is the identity monad on the category
$\Eee$ of sets.  Then $\Eee^+$ is the category of directed graphs, $T^+$ is
the free category monad $\fc$, and $T^+$-multicategories are the
$\fc$-multicategories of Chapter~\ref{ch:fcm}.  We therefore have a notion
of `$V$-enriched%
\index{enrichment!category@of category!fc-multicategory@in $\fc$-multicategory}
category' for any $\fc$-multicategory%
\index{fc-multicategory@$\fc$-multicategory!enrichment in}
$V$.  In the special
case that $V$ is a monoidal category~(\ref{eg:fcm-mon-cat}), we recover the
standard definition of enrichment.  More generally, if $V$ is a
bicategory~(\ref{eg:fcm-bicat}) then we recover the less well-known
definition of category enriched%
\index{enrichment!category@of category!bicategory@in bicategory}%
\index{bicategory!enrichment in}
in a bicategory (Walters~\cite{Wal}),
and if $V$ is a plain multicategory~(\ref{eg:fcm-cl-mti}) then we recover
the definition of category enriched in a plain multicategory.  For more on
enriched categories in this broad sense, see my~\cite{GEC}.
\end{example}

\begin{example}
The next most basic example is when $T$ is the free monoid monad on the
category $\Eee$ of sets.  Theorem~\ref{thm:sm-opetopic} tells us that any
symmetric monoidal category gives rise canonically to what is there called
a $T_2$-multicategory.  By definition, $T_2$ is the free $T$-operad monad
and $T^+$ the free $T$-multicategory monad, so a $T_2$-multicategory is a
special kind of $T^+$-multicategory.  Hence any symmetric monoidal category
gives rise canonically to a $T^+$-multicategory, giving us a definition of
plain multicategory%
\index{enrichment!plain multicategory@of plain multicategory!symmetric monoidal category@in symmetric monoidal category}
enriched in a symmetric monoidal category, and in
particular, of plain operad%
\index{operad!symmetric monoidal category@in symmetric monoidal category}
in a symmetric monoidal category.  These are
exactly the usual definitions (pp.~\pageref{p:sym-enr-mti},
\pageref{p:defn-V-Operad}).
\end{example}

\begin{example}	\lbl{eg:relaxed}
Borcherds~\cite{Borch}%
\index{Borcherds, Richard}
introduced certain structures called `relaxed%
\index{relaxed multicategory}%
\index{multicategory!relaxed}
multilinear categories' in his definition of vertex%
\index{vertex!algebra}
algebras over a vertex
group, and Soibelman~\cite{SoiMTC, SoiMBC}%
\index{Soibelman, Yan}
defined the same structures
independently in his work on quantum affine algebras.  As explained by
Borcherds, they can be regarded as categorical structures in which the maps
have singularities%
\index{singularity}
whose severity is measured by trees.  They also arise
completely naturally in the theory of enrichment: if $T$ is the free monoid
monad on the category $\Eee$ of sets, as in the previous example, then
there is a certain canonical $T^+$-multicategory $V$ such that
$V$-enriched%
\index{enrichment!plain multicategory@of plain multicategory!T-Zmulticategory@in $T^+$-multicategory}
$T$-multicategories are precisely relaxed multilinear categories.  See
Leinster~\cite[Ch.~4]{GECM} for details.
\end{example}

\begin{example}
In the next chapter we introduce the sequence $(T_n)_{n\in\nat}$ of
`opetopic'%
\index{opetopic!monad}
monads.  By definition, $T_n$ is the free $T_{n-1}$-operad
monad, and this means that there is a notion of $T_{n-1}$-multicategory
enriched in a $T_n$-multicategory.  More vaguely, a $T_n$-multicategory is
naturally regarded as an $(n+1)$-dimensional structure, so $n$-dimensional
structures can be enriched in $(n+1)$-dimensional structures.
\end{example}

\begin{notes}

Most parts of this chapter have appeared before \cite[\S 4]{GOM},
\cite[Ch.~3]{OHDCT}.  The thought that an operad is a cartesian monad
equipped with a cartesian natural transformation down to the free monoid
monad (\ref{sec:alt-app}) is closely related to Kelly's%
\index{Kelly, Max}
idea of a
`club'~\cite{KelCD, KelCDC}.%
\index{club}
 See Snydal~\cite{SnyEBG, SnyRMC}%
\index{Snydal, Craig}
for more on
the relaxed multicategory definition of vertex algebra.

\end{notes}

\chapter{Opetopes}
\lbl{ch:opetopic}

\chapterquote{%
John Dee [\ldots] summoned angels of dubious celestial provenance by
invoking names like Zizop, Zchis, Esiasch, Od and Iaod}{%
Eco~\cite{Eco}}

\noindent
Operads lead inescapably into geometry.  In this chapter we see that as
soon as the notion of generalized operad is formulated, the notion of
opetope is unavoidable.  Opetopes are something like simplices: they are a
completely canonical family of polytopes, as pervasive in
higher-dimensional algebra as simplices are in geometry.

In~\ref{sec:opetopes} opetopes are defined and their geometric
representation explained.  Intertwined with the definition of opetope is
the definition of a certain sequence $(T_n)_{n\in\nat}$ of cartesian
monads; we also look at $T_n$-multicategories and their relation to
symmetric multicategories.  Section~\ref{sec:pds} is formally about
$T_n$-structured categories (in the sense of~\ref{sec:struc}); translated
into geometry, this means diagrams of opetopes pasted together.

As we shall see, there is a category of $n$-dimensional pasting diagrams
for each natural number $n$.  When $n=1$ this is the category $\scat{D}$%
\index{augmented simplex category $\scat{D}$}
of
finite totally ordered sets (the augmented simplex category); when $n=2$ it
is a category of trees, as found in parts of quantum algebra.
In~\ref{sec:trees} we analyse the category of trees in some detail.  Some
of this is analogous to some of the standard analysis of $\scat{D}$: see
for instance Mac Lane~\cite[VII.5]{MacCWM}, where the monics, epics,
primitive face and degeneracy maps, and standard factorization properties
are all worked out, and the universal role of $\scat{D}$ (as the free
monoidal category containing a monoid) is established.

Opetopes were invented by Baez%
\index{Baez, John}
and Dolan~\cite{BDHDA3}%
\index{Dolan, James}
so that they could
frame a definition of weak $n$-category.  The strategy is simple: opetopes
together with face maps form a category, a presheaf on that category is
called an opetopic set, and a weak $n$-category is an opetopic set with
certain properties.  This is like both the definition of Kan complex (a
simplicial set with horn-filling properties) and the definition of
representable multicategory (a multicategory with universality
properties,~\ref{defn:repn-multi}).  Opetopic sets are discussed
in~\ref{sec:ope-sets}, and opetopic definitions of weak $n$-category
in~\ref{sec:ope-n}.

We finish~(\ref{sec:many}) with a short section on the `many%
\index{many in, many out}
in, many out'
approach to higher-dimensional category theory.  Opetopes arise from the
concept of operations with many inputs and a single output.  We could start
instead with the concept of operations having both many inputs and many
outputs, and try to construct shapes analogous to opetopes.  But this turns
out to be impossible, as we see.  (Formally, 3-computads do not form a
presheaf category.)  This does not mean that it is hopeless to try to
develop a many-in, many-out framework for higher categorical structures,
but it does mean that such a framework would be qualitatively different
from the simplicial, cubical, globular, and opetopic frameworks, in each of
which there is a genuine category of shapes.

\section{Opetopes}
\lbl{sec:opetopes}

The following table shows some of the types of generalized operad that we
have met.  As usual, $T$ is a cartesian monad on a cartesian category
$\Eee$.
\begin{center}
\begin{tabular}{lll}
$\Eee$		&$T$			&$T$-operads			\\
		&			&				\\
$\Set$		&identity		&monoids			\\
$\Set$		&free monoid		&plain operads		\\
$\Set^\nat \eqv \Set/\nat$	&
		 free plain operad	&
			(see~\ref{eg:mti-free-cl-opd})
\end{tabular}
\end{center}
Each row is generated automatically from the last by taking $T$ to be
`free%
\index{generalized operad!free}
operad of the type in the last row'.  Technicalities aside, it is clear
that this table can be continued indefinitely.

This gives an absolutely fundamental infinite sequence of categories
$\Eee$, monads $T$, and types of generalized operad.  The bulk of this
chapter consists of working out what they look like.  It is plain
from the first few rows that there will be some geometrical content---for
instance, in the third row, the operations in a $T$-operad are indexed by
trees (Example~\ref{eg:mti-free-cl-opd}).

The table can be compared with that of Ginzburg%
\index{Ginzburg, Victor}
and
Kapranov~\cite[p.~204]{GiKa}.%
\index{Kapranov, Mikhail}
 They include columns marked `Geometry'
(whose entries are `vector bundles', `manifolds', and `?(moduli spaces)'),
`Linear Physics', and `Nonlinear Physics'.  Their table, like the one
above, has only three rows.  Here we show how to continue forever---in our
columns, at least.

To do this formally we need to recall the results on free operads
in~\ref{sec:free-mti}.  There we met the concept of `suitable' categories
and monads and saw that they provided a good context for free operads.
Specifically, let $T$ be a suitable monad on a suitable category $\Eee$
with a terminal object.  Then Theorem~\ref{thm:free-fixed} tells us that
the forgetful functor
\[
\Cartpr\hyph\Operad \go \Eee/T1
\]
sending a $T$-operad to its underlying $T$-graph has a left adjoint, that
the induced monad (`free $T$-operad') on $\Eee/T1$ is suitable, and that
the category $\Eee/T1$ is suitable.  Trivially, $\Eee/T1$ has a terminal
object.  Also, Theorem~\ref{thm:free-gen} tells us that the category of
sets and the identity monad on it are both suitable.  This makes possible:
\begin{defn}	\lbl{defn:opetopic-monads}
For each $n\in\nat$, the suitable category $\Eee_n$%
\glo{nthopecat}
and the suitable
monad $T_n$%
\glo{nthopemonad}%
\index{opetopic!monad}
on $\Eee_n$ are defined inductively by
\begin{itemize}
\item $\Eee_0 = \Set$ and $T_0 = \id$
\item $\Eee_{n+1} = \Eee_n / T_n 1$ and
$T_{n+1} = (\textrm{free } T_n \textrm{-operad}).$ 
\end{itemize}
\end{defn}

In fact, there is for each $n\in\nat$ a set $O_n$%
\glo{nopes}
such that $\Eee_n$ is
canonically isomorphic to $\Set/O_n$.  (Recall also from
p.~\pageref{eq:Set-slice-power} that $\Set/O_n \eqv \Set^{O_n}$.)  First,
$O_0 = 1$.  Now suppose, inductively, that $n\geq 0$ and $\Eee_n \iso
\Set/O_n$.  The terminal object of $\Set/O_n$ is $(O_n \goby{1} O_n)$, so
if we put
\[
\bktdvslob{O_{n+1}}{t}{O_n} = T_n \bktdvslob{O_n}{1}{O_n}
\]
then we have
\[
\Eee_{n+1} = 
\frac{\Eee_n}{T_n 1} \iso
\frac{\Set/O_n}{O_{n+1} \goby{t} O_n} \iso
\frac{\Set}{O_{n+1}}.%
\glo{opetgt}
\]
This gives an infinite sequence of sets and functions
\[
\cdots 
\ 
\goby{t} 
O_{n+1} 
\goby{t} 
O_n 
\goby{t} 
\ 
\cdots
\ 
\goby{t} 
O_1 
\goby{t} 
O_0.
\]
The first few steps of the iteration are:
\begin{center}
\begin{tabular}{lllll}
$n$	&$O_n$		&$\Eee_n$	&$T_n$		&$T_n$-operads	\\
	&		&		&		&		\\
0	&$1$		&$\Set$		&identity	&monoids	\\
1	&$1$		&$\Set$		&free monoid	&
					plain operads		\\
2	&$\nat$&
			$\Set/\nat$	&free plain operad&
				(see~\ref{eg:mti-free-cl-opd} and below)\\
3	&$\{\textrm{trees}\}$&
		$\Set/\{\textrm{trees}\}$
\end{tabular}
\end{center}

We call $O_n$ the set of \demph{$n$-dimensional opetopes},%
\index{opetope}
or
\demph{$n$-opetopes}.  `Opetope' is pronounced%
\index{opetope!pronunciation}
with three syllables, as in
`OPEration polyTOPE'; opetopes encode fundamental
information about operations and are naturally represented as polytopes of
the corresponding dimension.  Let us see how this works in low dimensions.

The unique 0-opetope is drawn as a point:
\[
O_0 = \{ \gzero{} \}.
\]
The monad $T_0$ on $\Set/O_0 \iso \Set$ is the identity.  A
$T_0$-multicategory is an ordinary category, with objects $a$ depicted as
labelled points (0-opetopes) and maps $\theta$ as labelled arrows, as
usual:
\[
\gzero{a},
\diagspace
\gfst{a} \topebase{\theta} \glst{b}.
\]
In particular, a $T_0$-operad is a monoid, the elements $\theta$ of which
are drawn as
\[
\gfst{} \topebase{\theta} \glst{}.
\]
The underlying graph structure of a $T_0$-operad is formally a set over
$T_0 1 = O_1$, that is, a family of sets indexed by the elements of $O_1$.
Of course, $O_1$ has only one element and the `underlying graph structure'
of a monoid is just a set.  But since we want to view the elements of a
$T_0$-operad as labels on arrows, we choose to draw the unique element of
$O_1$ as an arrow:
\[
O_1 = \{ \gfst{}\topebase{}\glst{} \}.
\]
So the unique 1-opetope is drawn as a 1-dimensional polytope.

Next we have the monad $T_1 = (\textrm{free monoid})$ on the category
$\Set/O_1 \iso \Set$.  The free monoid $T_1 E$ on a set $E$ has elements
of the form
\[
\gfst{}\topebase{a_1}\gblw{}\topebase{a_2}
\diagspace \cdots \diagspace 
\topebase{a_n}\glst{}
\]
($n\geq 0, a_i \in E$).  A $T_1$-multicategory is a plain multicategory,
and with the diagrams we are using it is natural to draw an object $a$ as
labelling an edge,
\[
\gfst{}\topebase{a}\glst{},
\]
and a map $\theta: a_1, \ldots, a_n \go a$ as labelling a 2-dimensional
region,
\[
\topeq{a_1}{a_2}{a_n}{a}{\Downarrow \theta}.
\]
This is an alternative to our usual picture,
\[
\begin{array}{c}
\setlength{\unitlength}{1em}
\begin{picture}(4,7.5)(-2,0.3)
\cell{0}{6}{b}{\tinputssvert{a_1}{a_2}{a_n}}
\cell{0}{6}{t}{\tusualvert{\theta}}
\cell{0}{2}{t}{\toutputvert{a}}
\end{picture}
\end{array}.
\]
Sometimes we reduce clutter by omitting arrows:
\[
\topeqn{a_1}{a_2}{a_n}{a}{\Downarrow\theta}
\diagspace
\textrm{or}
\diagspace
\topeqn{a_1}{a_2}{a_n}{a}{\theta}.
\]
In this paradigm, composition is drawn as
\[
\begin{array}{c}
\setlength{\unitlength}{1mm}
\begin{picture}(44,22)(-2.5,-2)
\cell{0}{0}{bl}{\epsfig{file=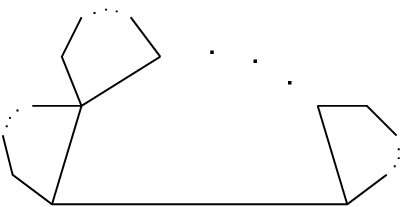}}
% edge labels
\cell{20}{-0.5}{t}{\scriptstyle a}
\cell{9}{5}{c}{\scriptstyle a_1}
\cell{14}{12}{c}{\scriptstyle a_2}
\cell{32}{5}{c}{\scriptstyle a_n}
\cell{2}{-1}{b}{\scriptstyle a_1^1}
\cell{0}{5}{r}{\scriptstyle a_1^2}
\cell{5}{9}{b}{\scriptstyle a_1^{k_1}}
\cell{8}{14}{c}{\scriptstyle a_2^1}
\cell{6}{18.5}{c}{\scriptstyle a_2^2}
\cell{17.5}{18}{c}{\scriptstyle a_2^{k_2}}
\cell{35}{10.5}{b}{\scriptstyle a_n^1}
\cell{39}{8}{bl}{\scriptstyle a_n^2}
\cell{37.5}{3}{tl}{\scriptstyle a_n^{k_n}}
% region labels
\da{20}{5}{0}
\cell{22}{5}{c}{\scriptstyle \theta}
\da{4}{5}{70}
\cell{5}{7}{c}{\scriptstyle \theta_1}
\da{11}{15}{30}
\cell{13}{16}{c}{\scriptstyle \theta_2}
\da{36}{7}{-70}
\cell{37}{5}{c}{\scriptstyle \theta_n}
\end{picture}
\end{array}
\diagspace
\goesto
\diagspace
\begin{array}{c}
\setlength{\unitlength}{1mm}
\begin{picture}(44,22)(-2.5,-2)
\cell{0}{0}{bl}{\epsfig{file=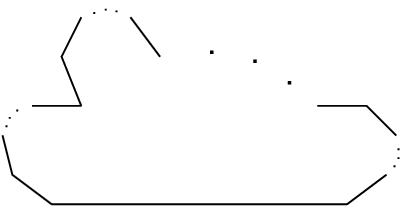}}
% edge labels
\cell{20}{-0.5}{t}{\scriptstyle a}
\cell{2}{-1}{b}{\scriptstyle a_1^1}
\cell{0}{5}{r}{\scriptstyle a_1^2}
\cell{5}{9}{b}{\scriptstyle a_1^{k_1}}
\cell{8}{14}{c}{\scriptstyle a_2^1}
\cell{6}{18.5}{c}{\scriptstyle a_2^2}
\cell{17.5}{18}{c}{\scriptstyle a_2^{k_2}}
\cell{35}{10.5}{b}{\scriptstyle a_n^1}
\cell{39}{8}{bl}{\scriptstyle a_n^2}
\cell{37.5}{3}{tl}{\scriptstyle a_n^{k_n}}
% region labels
\da{15}{5}{0}
\cell{17}{5}{l}{\scriptstyle \theta\sof (\theta_1, \ldots, \theta_n)}
\end{picture}
\end{array}
% \hand{35}{16}
\]
and the function assigning identities as
\[
\topebasen{a}
\diagspace
\goesto
\diagspace
\topean{a}{a}{\Downarrow 1_a}.
\]
A $T_1$-\emph{operad} (plain operad) looks just the same except that the edges
are no longer labelled: so an $n$-ary operation $\theta \in P(n)$ of an
operad $P$ is drawn as
\[
\topeqn{}{}{}{}{\theta}
\]
with $n$ `input' edges along the top and one `output' edge along the
bottom.  The underlying $T_1$-graph of a $T_1$-operad is a set over $T_1 1
= O_2 \iso \nat$, so we draw 2-opetopes as 2-dimensional polytopes:
\[
O_2 = 
\left\{
\topez{}{\Downarrow},
\ 
\topea{}{}{\Downarrow},
\ 
\topeb{}{}{}{\Downarrow},
\ 
\topec{}{}{}{}{\Downarrow},
\ 
\toped{}{}{}{}{}{\Downarrow},
\ 
\ldots
\ \ 
\right\}.
\]

Next we have the monad $T_2 = (\textrm{free plain operad})$ on the category
\[
\Set/O_2 \iso \Set/\nat \eqv \Set^\nat.
\]
An object $E$ of $\Set/O_2$ is a family $(E(\omega))_{\omega\in O_2}$ of
sets, with the elements of $E(\omega)$ regarded as potential labels to be
stuck on the 2-opetope $\omega$.  The free operad $T_2 E$ on $E$ is formed
by pasting together labelled 2-opetopes, with output edges joined to input
edges.  For instance, if $a_1 \in E(3)$, $a_2 \in E(1)$ and $a_3 \in E(2)$
then
\begin{equation}	\label{diag:labelled-two-pd}
\begin{array}{c}
\setlength{\unitlength}{1mm}
\begin{picture}(38,13)
\cell{0}{0}{bl}{\epsfig{file=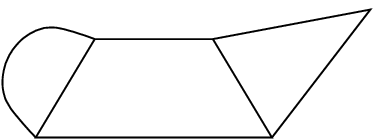}}
\da{15}{5}{0}
\cell{17}{5}{l}{a_1}
\da{4}{6}{60}
\cell{5}{8}{c}{a_2}
\da{27}{8}{-60}
\cell{28.5}{6}{c}{a_3}
\end{picture}
\end{array}
% \hand{25}{17}
\end{equation}
is an element of $(T_2 E)(4)$, and if $b_1 \in E(5)$, $b_2 \in E(0)$, $b_3,
b_4 \in E(3)$ and $b_5 \in E(1)$ then
\begin{equation}	\label{diag:bigger-two-pd}
\begin{array}{c}
\setlength{\unitlength}{1mm}
\begin{picture}(30,22)
\cell{0}{0}{bl}{\epsfig{file=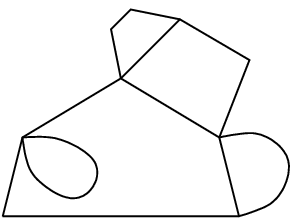}}
\da{15}{5}{0}
\cell{16.5}{5}{l}{b_1}
\da{6.5}{3.7}{60}
\cell{7.5}{5.7}{c}{b_2}
\da{19}{14}{-35}
\cell{20}{14}{tl}{b_3}
\da{13.5}{18}{50}
\cell{15}{19.5}{c}{b_4}
\da{25}{6}{-75}
\cell{26}{3.5}{c}{b_5}
\end{picture}
\end{array}
% \hand{40}{18}
\end{equation}
is an element of $(T_2 E)(8)$.  Previously we represented operations in a
free operad as trees%
\index{tree!vertices labelled@with vertices labelled}
with labelled vertices, so instead
of~\bref{diag:labelled-two-pd} we drew
\[
\setlength{\unitlength}{1em}
\begin{picture}(5,3)(-2.5,0)
% bottom layer
\put(0,0){\line(0,1){1}}
% middle layer
\cell{0}{1}{c}{\vx}
\cell{0.2}{1}{tl}{a_1}
\put(0,1){\line(-3,2){1.5}}
\put(0,1){\line(0,1){1}}
\put(0,1){\line(3,2){1.5}}
% top layer
\cell{-1.5}{2}{c}{\vx}
\cell{-1.7}{2}{tr}{a_2}
\put(-1.5,2){\line(0,1){1}}
\cell{1.5}{2}{c}{\vx}
\cell{1.7}{2}{tl}{a_3}
\put(1.5,2){\line(-1,1){1}}
\put(1.5,2){\line(1,1){1}}
\end{picture}
% \drmk{tree corr to tensor product }
% \otimes_{a_1}(\otimes_{a_2}(-), -, \otimes_{a_3}(-,-))
% \dr{35}{24}
\]
(\ref{sec:om-further} and~\ref{eg:mon-free-cl-opd}).  These pictures are
dual to one another, as Fig.~\ref{fig:two-pd-tree-dual}
\begin{figure}
\centering
\setlength{\unitlength}{1mm}
\begin{picture}(116,20)
\cell{0}{0}{bl}{\epsfig{file=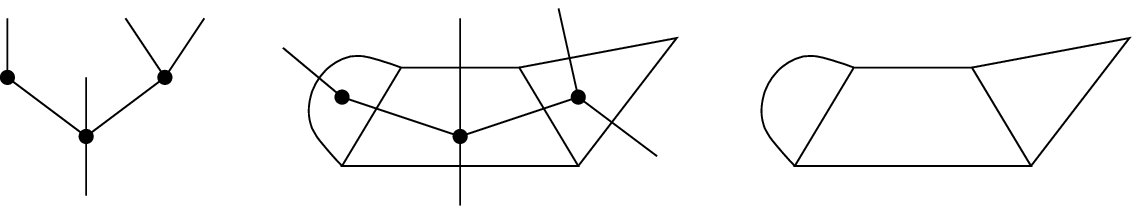}}
\da{93}{9}{0}
\da{81}{11}{60}
\da{105}{11}{-60}
\end{picture}
% \hand{35}{19}
\caption{How a tree corresponds to a diagram of pasted-together 2-opetopes}
\label{fig:two-pd-tree-dual} 
\end{figure}
demonstrates; we return to this correspondence in~\ref{sec:trees}.

We described $T_2$-multicategories in Example~\ref{eg:mti-free-cl-opd} in
terms of trees; we now describe them in terms of opetopes.  The objects of
a $T_2$-multicategory $C$ form a family $(C_0(\omega))_{\omega\in O_2}$ of
sets; put another way, they form a graded set $(C_0(n))_{n\in\nat}$ with,
for instance, $a \in C_0(3)$ drawn as
\[
\topecn{}{}{}{}{\Downarrow a}.
\]
Arrows look like
\begin{equation}	\label{diag:opetopic-three-arrow}
\begin{array}{c}
\setlength{\unitlength}{1mm}
\begin{picture}(38,13)
\cell{0}{0}{bl}{\epsfig{file=simple2pd.eps}}
\da{15}{5}{0}
\cell{17}{5}{l}{a_1}
\da{4}{6}{60}
\cell{5}{8}{c}{a_2}
\da{27}{8}{-60}
\cell{28.5}{6}{c}{a_3}
\end{picture}
\end{array}
\ 
\begin{array}{c}
\theta\\
\epsfig{file=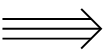}
\end{array}
\ 
\begin{array}{c}
\setlength{\unitlength}{1mm}
\begin{picture}(38,13)
\cell{0}{0}{bl}{\epsfig{file=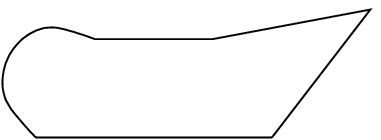}}
\da{15}{5}{0}
\cell{17}{5}{l}{a}
\end{picture}
\end{array}
\end{equation}
($a_1 \in C_0(3)$, $a_2 \in C_0(1)$, $a_3 \in C_0(2)$, $a \in C_0(4)$).
The 2-opetope in the codomain always has the same number of input edges as
the diagram in the domain (4, in this case); here it is drawn irregularly
to make the equality self-evident.  Visualize the whole
of~\bref{diag:opetopic-three-arrow} as a 3-dimensional polytope with one
flat bottom face labelled by $a$ and three curved top faces labelled by
$a_1$, $a_2$ and $a_3$ respectively, and with a label $\theta$ in the
middle.  On a sheet of paper we must settle for 2-dimensional
representations such as the one above.  Composition takes a diagram of
arrows such as
\begin{equation}	\label{diag:2-substn}
\begin{array}{c}
\setlength{\unitlength}{1mm}
\begin{picture}(106,41)
\cell{0}{0}{bl}{\epsfig{file=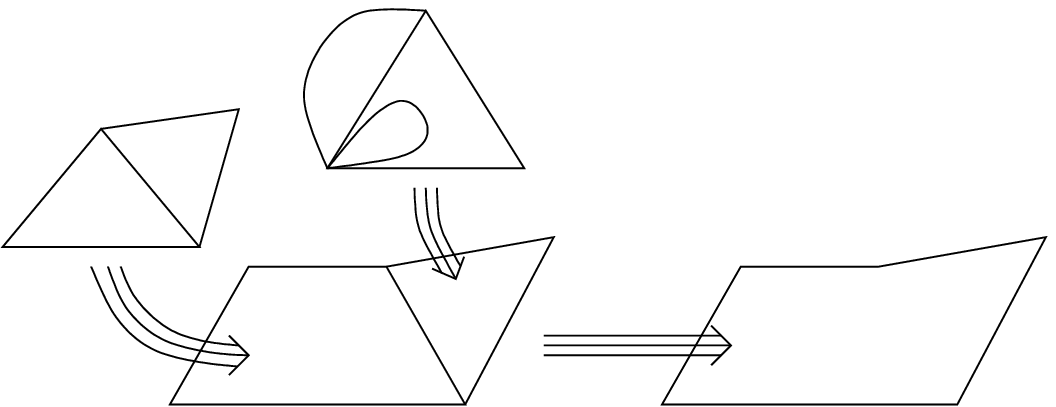}}
% region labels
\cell{32}{6}{c}{a_1}
\cell{47}{10}{c}{a_2}
\cell{10}{20}{c}{a_1^1}
\cell{18}{24}{c}{a_1^2}
\cell{46}{28}{c}{a_2^1}
\cell{40}{28}{c}{a_2^2}
\cell{36}{34}{c}{a_2^3}
\cell{85}{6}{c}{a}
% arrow labels
\cell{18}{10.5}{c}{\theta_1}
\cell{47.5}{19}{c}{\theta_2}
\cell{64}{9.5}{c}{\theta}
\end{picture}
\end{array}
% \hand{50}{20a}
\end{equation}
and produces a single arrow
\begin{equation}	\label{diag:2-substd}
\begin{array}{c}
\setlength{\unitlength}{1mm}
\begin{picture}(39,22)
\cell{0}{0}{bl}{\epsfig{file=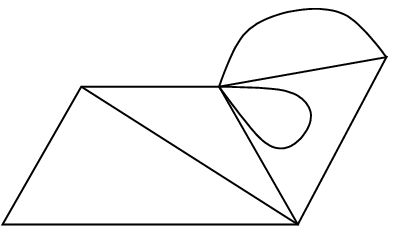}}
\cell{11}{5}{c}{a_1^1}
\cell{21}{10}{c}{a_1^2}
\cell{30}{6}{c}{a_2^1}
\cell{28}{11}{c}{a_2^2}
\cell{30}{19}{c}{a_2^3}
\end{picture}
\end{array}
\ 
\begin{array}{c}
\theta \of (\theta_1, \theta_2)\\
\epsfig{file=threearrow.eps}
\end{array}
\ 
\begin{array}{c}
\setlength{\unitlength}{1mm}
\begin{picture}(39,22)
\cell{0}{0}{bl}{\epsfig{file=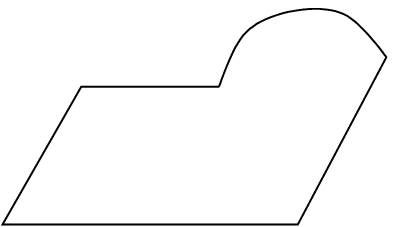}}
\cell{18}{7}{c}{a}
\end{picture}
\end{array}
.
% \hand{25}{20b}.
\end{equation}
(This is exactly the example of~\ref{eg:mti-free-cl-opd}.)  The identity on
an object $a \in C_0(n)$ looks like
\[
\topeqn{}{}{}{}{a} 
\ \ 
\begin{array}{c}
1_a\\
\epsfig{file=threearrow.eps}
\end{array}
\ \ 
\topeqn{}{}{}{}{a}\ , 
\]
where there are $n$ input edges in both the domain and the codomain.  

In particular, a $T_2$-operad consists of a collection of operations such
as
\[
\begin{array}{c}
\setlength{\unitlength}{1mm}
\begin{picture}(38,13)
\cell{0}{0}{bl}{\epsfig{file=simple2pd.eps}}
\da{15}{5}{0}
\da{4}{6}{60}
\da{27}{8}{-60}
\end{picture}
\end{array}
\ 
\begin{array}{c}
\theta\\
\epsfig{file=threearrow.eps}
\end{array}
\ 
\begin{array}{c}
\setlength{\unitlength}{1mm}
\begin{picture}(38,13)
\cell{0}{0}{bl}{\epsfig{file=simple2pdoutline.eps}}
\da{15}{5}{0}
\end{picture}
\end{array}
\]
together with composition and identities as above.  So the elements of $O_3
= T_2(1)$---the 3-opetopes---are thought of as 3-dimensional polytopes, for
instance 
\begin{equation}	\label{diag:three-opetope}
\begin{array}{c}
\setlength{\unitlength}{1mm}
\begin{picture}(38,13)
\cell{0}{0}{bl}{\epsfig{file=simple2pd.eps}}
\da{15}{5}{0}
\da{4}{6}{60}
\da{27}{8}{-60}
\end{picture}
\end{array}
\ 
\begin{array}{c}
\epsfig{file=threearrow.eps}
\end{array}
\ 
\begin{array}{c}
\setlength{\unitlength}{1mm}
\begin{picture}(38,13)
\cell{0}{0}{bl}{\epsfig{file=simple2pdoutline.eps}}
\da{15}{5}{0}
\end{picture}
\end{array}.
\end{equation}
The function $t: O_3 \go O_2$ is `target'; it sends the 3-opetope above to
the 2-opetope with 4 input edges.

This gives a systematic way of portraying opetopes and
$T_n$-multicategories for arbitrary $n$.  In~\ref{sec:ope-sets} we will see
that each $n$-opetope does indeed give rise to an $n$-dimensional
topological space.

Some crude examples of $T_n$-multicategories are provided by symmetric
structures.  First, any commutative%
\index{monoid!commutative!operad from}
monoid $(A,+,0)$ gives rise to a
$T_n$-operad for every $n\in\nat$: an operation of shape $\omega\in O_n$ is
just an element of $A$ (regardless of what $\omega$ is), composition is
$+$, and the identity is $0$.  We have already seen this in the case $n=1$
(plain operads,~\ref{eg:opd-from-comm-mon}).  Formally, let
$\fcat{CommMon}$%
\glo{CommMon}
be the category of commutative monoids and $\Delta: \Set
\go \Set/O_n$%
\glo{Deltadiag}
the functor sending a set $A$ to $(A \times O_n
\goby{\mr{pr}_2} O_n)$; then we have
\begin{thm}	\lbl{thm:cm-ftr-opetopic}
For each $n\in\nat$ there is a canonical functor
\[
\fcat{CommMon} \go T_n\hyph\Operad
\]
making the diagram
\[
\begin{diagram}
\fcat{CommMon}	&\rTo		&T_n\hyph\Operad	\\
\dTo		&		&\dTo			\\
\Set		&\rTo_\Delta	&\Set/O_n		\\
\end{diagram}
\]
commute, where the vertical arrows are the forgetful functors. 
\end{thm}
\begin{proof}
This follows from $T_n$ being a finitary familially representable monad on
a slice of $\Set$: see Example~\ref{eg:cm-ftr-opetopic}.
\done
\end{proof}

\index{multicategory!symmetric vs. generalized@symmetric \vs.\ generalized|(}
More interestingly, any symmetric multicategory $A$ gives rise to a
$T_n$-multicategory $C$ for every $n\in\nat$, canonically up to isomorphism.
The objects of $C$ of shape $\omega\in O_n$ are simply the objects of $A$
(regardless of $\omega$).  The arrows of $C$ are arrows of $A$: for
example, if $n=2$ and $\widehat{a}, \twid{a}, \ldots$ are objects of $A$
then an arrow
\begin{equation}	\label{eq:Tn-sym-arrow}
\begin{array}{c}
\setlength{\unitlength}{1mm}
\begin{picture}(42,19)
\cell{0}{0}{bl}{\epsfig{file=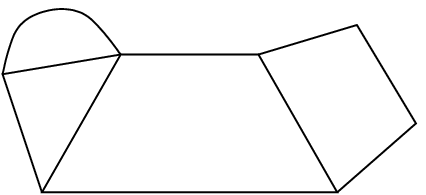}}
\cell{19}{7}{c}{\ovln{a}}
\cell{6}{9}{c}{a'}
\cell{6}{15.5}{c}{\twid{a}}
\cell{34}{9}{c}{\widehat{a}}
\end{picture}
\end{array}
\ 
\begin{array}{c}
\epsfig{file=threearrow.eps}
\end{array}
\ 
\begin{array}{c}
\setlength{\unitlength}{1mm}
\begin{picture}(42,19)
\cell{0}{0}{bl}{\epsfig{file=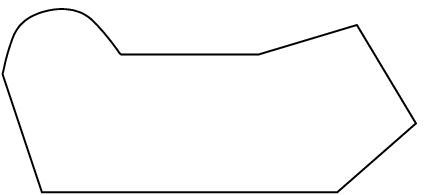}}
\cell{20}{7}{c}{a}
\end{picture}
\end{array}
% \hand{25}{29}
\end{equation}
in $C$ might be defined to be an arrow
\[
a', \widehat{a}, \ovln{a}, \twid{a} \go a
\]%
\lbl{p:sym-ordering}%
in $A$---or indeed, the same but with the four domain objects ordered%
\index{order!non-canonical}%
\index{pasting diagram!opetopic!ordering of}
differently.  There is no canonical ordering of the $2$-opetopes making up
a given $2$-pasting diagram; more precisely, there is no method of ordering
that is stable under substitution of the kind shown in
diagrams~\bref{diag:2-substn} and~\bref{diag:2-substd}.  This is why we
needed to start with a \emph{symmetric} multicategory: composition in $C$
cannot be defined without permuting the lists of objects.  It is also why
$C$ is canonical only up to isomorphism.  The precise statement of the
construction is:
\begin{thm}	\lbl{thm:sm-opetopic}
For each $n\in\nat$ there is a functor 
\[
\fcat{SymMulticat} \go T_n\hyph\Multicat,
\]
canonical up to isomorphism, making the diagram
\[
\begin{diagram}
\fcat{SymMulticat}	&\rTo		&T_n\hyph\Multicat	\\
\dTo<{\blank_0}		&		&\dTo>{\blank_0}	\\
\Set			&\rTo_\Delta	&\Set/O_n		\\
\end{diagram}
\]
commute, where both vertical arrows are the functors assigning to a
multicategory its object of objects.
\end{thm}
\begin{proof}
Again this follows from $T_n$ being a finitary familially representable
monad on a slice of $\Set$: Example~\ref{eg:sm-opetopic}.
\done
\end{proof}

Symmetric%
\lbl{p:sym-discussion}
structures have a ghostly presence throughout this book, hovering just
beyond our world of cartesian monads and generalized multicategories.  Much
of the time we are concerned with labelled cell diagrams such as the domain
of~\bref{eq:Tn-sym-arrow}, and exercise full sensitivity to their geometric
configuration.  By passing to a symmetric multicategory we destroy all the
geometry.

In this sense, symmetric multicategories are the ultimate in crudeness.
There are various situations in mathematics where symmetric multicategories
(or symmetric monoidal categories) are customarily used, but generalized
multicategories provide a more sensitive and more general approach.
Operads in a symmetric monoidal category and categories enriched in a
(symmetric or not) monoidal category are two examples; the more thoughtful
approaches replace symmetric monoidal categories by $T_2$-multicategories
and $\fc$-multicategories, respectively~(\ref{sec:enr-mtis}).  There are
also entire approaches to higher-dimensional category theory based on
symmetric structures, as discussed in~\ref{sec:ope-n}; crude does not mean
ineffective.%
\index{multicategory!symmetric vs. generalized@symmetric \vs.\ generalized|)}

\section{Categories of pasting diagrams}
\lbl{sec:pds}

The principal thing that you can do in a higher-dimensional category is to
take a diagram of cells and form its composite.  Not just any old diagram
will do: it must, for instance, be connected and have the cells oriented
compatibly.  For example, the acceptable diagrams of 1-cells are those of
the form
\begin{equation}	\label{diag:ope-one-pd}
\gfstsu\gonesu\gzersu\gonesu 
\diagspace \cdots \diagspace 
\gonesu\glstsu
\end{equation}
where the number of arrows is a non-negative integer.  Let us call a
composable diagram of pasted-together $n$-cells an `$n$-pasting%
\index{pasting diagram}
diagram'.
For $n\geq 2$, the class of $n$-pasting diagrams depends on the shape of
cells that you have chosen to use in your theory of higher-dimensional
categories: globular, cubical, simplicial, opetopic, \ldots.

Pasting diagrams play an important role in any theory of higher categories.
What distinguishes the opetopic theory is that pasting diagrams are the
same thing as cell shapes of one dimension higher.  For instance, the
1-pasting diagram~\bref{diag:ope-one-pd} with $k$ arrows corresponds to the
2-opetope
\[
\topeq{}{}{}{}{\Downarrow}
\]
with $k$ arrows along the top; we draw them differently, but there is a
natural identification.  So we may conveniently \emph{define} an
\demph{(opetopic) $n$-pasting diagram}%
\index{pasting diagram!opetopic}
to be an $(n+1)$-opetope, for any
$n\in\nat$.  In this section we show how, for each $n$, the $n$-pasting
diagrams form a category $\Pd{n}$---and actually, rather more than just a
category.  

The formal method is as follows.  So far we have looked at $T_n$-operads
and $T_n$-multicategories; now we look at $T_n$-structured%
\index{structured category!opetopic monad@for opetopic monad}
categories.
$\PD{n}$ is defined as the free $T_n$-structured category on the terminal
$T_n$-multicategory, and $\Pd{n}$ as the underlying category of $\PD{n}$.
So, we begin by recalling what $T$-structured categories are in general and
what they look like when $T=T_n$; then we look at free structured
categories; finally, we arrive at the category of $n$-pasting diagrams.
The next section,~\ref{sec:trees}, is a detailed examination of the case
$n=2$: it turns out that $\Pd{2}$ is a category of trees.

Recall from~\ref{sec:struc} that if $T$ is a cartesian monad on a cartesian
category $\Eee$ then a $T$-structured category is an internal category in
$\Eee^T$.  Alternatively, note that $T$ lifts naturally to a monad
$\Cat(T)$ on $\Cat(\Eee)$, and a $T$-structured category is then a
$\Cat(T)$-algebra.

A $T_0$-structured category is a category in $\Set^{T_0} \iso \Set$, that
is, a category.  

A $T_1$-category is a category in $\Set^{T_1} \iso \fcat{Monoid}$, that is,
a strict monoidal category.  (Alternatively, $\Cat(T_1)$ is the free strict
monoidal category monad on $\Cat(\Set) \iso \Cat$, and a $T_1$-structured
category is a $\Cat(T_1)$-algebra.)

A $T_2$-structured category is a category in $(\Set/\nat)^{T_2} \iso
\Operad$.  Alternatively, it is an operad in \Cat, or `\Cat-operad',%
\index{Cat-operad@$\Cat$-operad}
as we
saw in~\ref{eg:struc-fin-lims}.  Diagrammatically, a $T_2$-structured
category $A$ consists of
\begin{itemize}
\label{p:T2-diagrammatic}
\item a set $A_0(k)$ for each $k\in\nat$, with $a \in A_0(k)$ drawn as a
label on the $k$th 2-opetope,
\[
\topeqn{}{}{}{}{a}
\]
\item a set $A(a,b)$ for each $k\in\nat$ and $a,b \in A_0(k)$, with $\theta
\in A(a,b)$ drawn as
\[
\topeqn{}{}{}{}{a} 
\diagspace \goby{\theta} \diagspace
\topeqn{}{}{}{}{b}
\]
\item a function defining composition or `gluing' of objects,
\[
\begin{array}{c}
\setlength{\unitlength}{1mm}
\begin{picture}(41,20)
\cell{0}{0}{bl}{\epsfig{file=opeT1compdom.eps}}
\cell{20}{6}{c}{a}
\cell{4}{6}{c}{a_1}
\cell{11}{15}{c}{a_2}
\cell{37}{6}{c}{a_n}
\end{picture}
\end{array}
\diagspace
\goesto
\diagspace
\begin{array}{c}
\setlength{\unitlength}{1mm}
\begin{picture}(41,20)
\cell{0}{0}{bl}{\epsfig{file=opeT1compcod.eps}}
\cell{20}{6}{c}{a \of (a_1, \ldots, a_n)}
\end{picture}
\end{array}
% \hand{20}{23},
\]
and an identity or `unit' object,
\[
\topean{}{}{1}
\]
\item a function defining composition of arrows,
\[
(a \goby{\theta} b \goby{\phi} c) 
\diagspace\goesto\diagspace
(a \goby{\phi\theta} c),
\]
and an identity arrow $(a \goby{1_a} a)$ on each object $a$
\item a function defining gluing of arrows:
\[
a \goby{\theta} b, \ 
a_1 \goby{\theta_1} b_1, \ \ldots,\  a_k \goby{\theta_k} b_k
\]
give rise to 
\[
a \of (a_1, \ldots, a_k) 
\goby{\theta * (\theta_1, \ldots, \theta_k)}
b \of (b_1, \ldots, b_k),
\]
\end{itemize}
all satisfying the usual kinds of axioms.  The sets $A_0$ and $A_1 =
\coprod_{a,b\in C_0} A(a,b)$ both have the structure of plain operads:
thus, $A$ can be viewed as a category in $\Operad$.  Regrouping the data,
there is for each $n$ a category $A(n)$ whose object-set is $A_0(n)$: thus,
$A$ can also be viewed as an operad in $\Cat$.

An analogous diagrammatic description applies to $T_n$-structured
categories for any $n\in\nat$.

Next recall from~\ref{sec:struc} that there is an adjunction between
$T$-structured categories and $T$-multicategories, which for $T=T_n$ will
be denoted
\[
\begin{diagram}[height=2em]
T_n\hyph\Struc			\\
\uTo<{F_n} \ladj \dTo>{U_n}	\\
T_n\hyph\Multicat.		\\
\end{diagram}
\]
Also recall from~\ref{sec:opetopes} that a $T_n$-multicategory consists
of a set of objects labelling $n$-opetopes, a set of maps whose domains are
labelled $n$-pasting diagrams and whose codomains are labelled single
$n$-opetopes, and functions defining composition and identities.  Now, let
us see what this adjunction looks like.

The functor $U_n$ `forgets how to tensor but remembers multilinear maps'.
When $n=0$ it is the identity, when $n=1$ it sends a strict monoidal
category to its underlying plain multicategory, and when $n=2$ it sends
a \Cat-operad $A$ to the $T_2$-multicategory whose objects are the same as
those of $A$ and whose maps
\[
\begin{array}{c}
\setlength{\unitlength}{1mm}
\begin{picture}(29,21)
\cell{0}{0}{bl}{\epsfig{file=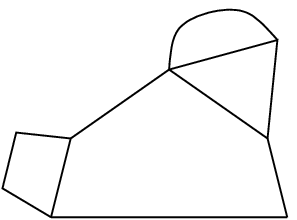}}
\cell{17}{5}{c}{a_1}
\cell{4}{5}{c}{a_2}
\cell{24}{14}{c}{a_3}
\cell{22}{19}{c}{a_4}
\end{picture}
\end{array}
\diagspace
\go
\diagspace
\begin{array}{c}
\setlength{\unitlength}{1mm}
\begin{picture}(29,21)
\cell{0}{0}{bl}{\epsfig{file=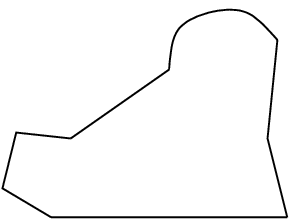}}
\cell{16}{6}{c}{a}
\end{picture}
\end{array}
% \hand{25}{24a}
\]
(for instance) are maps
\[
\begin{array}{c}
\setlength{\unitlength}{1mm}
\begin{picture}(29,21)
\cell{0}{0}{bl}{\epsfig{file=struccod.eps}}
\cell{14}{6}{c}{\scriptstyle a_1 \of (a_2, 1, a_3\of (a_4, 1))}
\end{picture}
\end{array}
\diagspace
\go
\diagspace
\begin{array}{c}
\setlength{\unitlength}{1mm}
\begin{picture}(29,21)
\cell{0}{0}{bl}{\epsfig{file=struccod.eps}}
\cell{16}{6}{c}{a}
\end{picture}
\end{array}
% \hand{25}{24b}
\]
in $A$.

The free functor $F_n$ is formal pasting: if $C$ is a $T_n$-multicategory
then the objects (respectively, arrows) of $F_n C$ are the formal pastings
of objects (respectively, arrows) of $C$.  Trivially, $F_0$ is the
identity.  We described $F_1$ on p.~\pageref{diag:arrows-in-mon-cat}, in a
different diagrammatic style.  

\begin{defn}
Let $n\geq 0$.  The \demph{structured%
\index{structured category!pasting diagrams@of pasting diagrams}%
\index{pasting diagram!opetopic!structured category of}
category of $n$-pasting diagrams},
$\PD{n}$,%
\glo{PDn}
is defined by $\PD{n} = F_n 1$.  In other words, $\PD{n}$ is the
free $T_n$-structured category on the terminal $T_n$-multicategory.
\end{defn}
So $\PD{n}$ is an internal category in $(\Set/O_n)^{T_n}$; its underlying
graph is
\begin{equation}	\label{diag:PD-n-graph}
\begin{slopeydiag}
	&	&T_n^2 1	&	&	\\
	&\ldTo<{\mu_1}&		&\rdTo>{T_n !}&	\\
T_n 1	&	&		&	&T_n 1	\\
\end{slopeydiag}
\end{equation}
where the $T_n$-algebra structures on $T_n 1$ and $T_n^2 1$ are both
components of the multiplication $\mu$ of the monad $T_n$.  

A $T_0$-structured category is just a category, and $\PD{0}$ is the
terminal category (whose object is viewed as the unique 0-pasting diagram
$\gzero{}$).  

A $T_1$-structured category is a strict monoidal category, and we have
already seen in Example~\ref{eg:free-struc-D} that $\PD{1}$ is $\scat{D}$,%
\index{augmented simplex category $\scat{D}$}
the strict monoidal category of (possibly empty) finite totally ordered
sets $\lwr{n} = \{1, \ldots, n\}$; addition is tensor and $\lwr{0}$ is the
unit.  The diagram above is in this case
\[
\begin{slopeydiag}
	&	&\nat^*		&	&	\\
	&\ldTo<+&		&\rdTo	&	\\
\nat	&	&		&	&\nat,	\\
\end{slopeydiag}
\]
where $\nat^*$ is the set of finite sequences of natural numbers and the
right-hand map sends $(m_1, \ldots, m_n)$ to $n$.  

The $\Cat$-operad $\PD{2}$ is described in detail in the next section.

We have been looking at $\PD{n}$, the $T_n$-structured category of
$n$-pasting diagrams, but sometimes it is useful to forget the more
sophisticated structure and pass to the mere category of $n$-pasting
diagrams.  Formally, we have cartesian forgetful functors
\[
(\Set/O_n)^{T_n} \go \Set/O_n \go \Set
\]
and these induce a forgetful functor
\[
T_n \hyph \Struc = \Cat((\Set/O_n)^{T_n})
\go
\Cat(\Set) = \Cat,
\]
making possible the following definition.
\begin{defn}
Let $n\geq 0$.  The \demph{category of $n$-pasting%
\index{pasting diagram!opetopic!category of}
diagrams}, $\Pd{n}$,%
\glo{Pdn}
is
the image of $\PD{n}$ under the forgetful functor $T_n\hyph\Struc \go
\Cat$.  
\end{defn}
For example, $\Pd{0}$ is the terminal category, $\Pd{1}$ is the category
$\scat{D}$ of finite totally ordered sets, and $\Pd{2}$ can be regarded as
a category of trees (see below).  The objects of $\Pd{n}$ really
are the $n$-pasting diagrams: for by~\bref{diag:PD-n-graph}, the object-set of
$\Pd{n}$ is the underlying set of $T_n 1 \in \Set/O_n$, which is the set
$O_{n+1}$ of $(n+1)$-opetopes or $n$-pasting diagrams.

In~\ref{sec:ope-sets} we will consider the category of all opetopes;%
\index{opetope!category of}
beware
that this is quite different from the categories $\Pd{n}$ of $n$-pasting
diagrams.

It is instructive to contemplate the situation for arbitrary cartesian
$\Eee$ and $T$.  We have functors
\[
\begin{diagram}[height=2em,width=4em]
T\hyph\Struc = \Cat(\Eee^T)	&
\rTo^V				&
\Cat(\Eee)			\\
\uTo<F \ladj \dTo>U		&	&	\\
T\hyph\Multicat,			&	&	\\
\end{diagram}
\]
where $V$ is forgetful, and so we have an internal category $VF1$ in
$\Eee$.  The object-of-objects of $VF1$ is $T1$, which may be thought of as
the object of $T$-pasting diagrams; hence $VF1$ may be thought of as the
(internal) category of $T$-pasting%
\index{pasting diagram}
diagrams.  In some situations (such as
when $\Eee = \Eee_n \iso \Set/O_n$) there is an obvious cartesian
`forgetful' functor $\Eee \go \Set$, and then there is an induced functor
$\Cat(\Eee) \go \Cat$, giving a genuine category of $T$-pasting diagrams.
For instance, if $K$ is a set and $T$ is the monad $K + (\dashbk)$ on $\Eee
= \Set$ then this is the category $\scat{P}_K$ defined on
p.~\pageref{p:defn-wide-pb-shape}.

\section{A category of trees}
\lbl{sec:trees}

We have seen that for each $n\geq 0$ there is a category $\Pd{n}$ whose
objects are $n$-pasting diagrams.  We have also seen that 2-pasting
diagrams correspond naturally to trees (Fig.~\ref{fig:two-pd-tree-dual}).
Hence $\Tr = \Pd{2}$%
\glo{Tr}%
\index{tree!category of}
is a category whose objects are trees.  Here we
describe it in some detail.  The definition of a map between trees is
perfectly natural but takes some getting used to; we approach it slowly.

First recall what trees themselves are.  By
definition~(\ref{eg:opd-of-trees}), $\tr$ is the free plain operad on the
terminal object of $\Set^\nat$, and an $n$-leafed tree is an element of
$\tr(n)$.  As we saw, the sets $\tr(n)$ also admit the following recursive
description:
\begin{itemize}
\item $\utree\in\tr(1)$
\item if $n, k_1, \ldots, k_n \in \nat$ and $\tau_1 \in \tr(k_1), \ldots,
\tau_n \in \tr(k_n)$ then $(\tau_1, \ldots, \tau_n) \in \tr(k_1 + \cdots +
k_n)$.
\end{itemize}
For the purposes of this text we will need no further description of what a
tree is.  But it is also possible, as you might expect, to describe a tree
as a graph%
\index{tree!graph@as graph}%
\index{graph!tree@of tree}
of a certain kind, and this alternative, `concrete', description
can be comforting.  To say exactly \emph{what} kind of graph is more
delicate than meets the eye, and in papers on operads is often done only
vaguely.  I have therefore put a graph-theoretic definition of tree, and a
proof of its equivalence to our usual one, in Appendix~\ref{app:trees}.

\begin{figure}
\centering
\setlength{\unitlength}{1mm}
\begin{picture}(70,24)
\cell{0}{5.5}{l}{\textrm{(b)}}
\cell{15.6}{0}{bl}{\epsfig{file=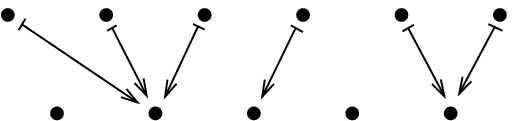}}
\cell{0}{22}{l}{\textrm{(a)}}
\cell{10}{20}{bl}{\epsfig{file=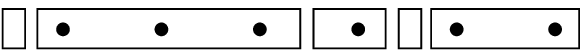}}
\end{picture}
% \hand{35}{30}
\caption{Two pictures of a map $6 \protect\go 5$ in $\scat{D}$}
\label{fig:map-in-D}
\end{figure}

\index{tree!map of|(}
In preparation for looking at maps in $\Pd{2}$---maps between trees---let
us look again at maps in $\Pd{1} = \scat{D}$.%
\index{augmented simplex category $\scat{D}$}
 An object of $\scat{D}$ is a
natural number.  A map is, as observed in the previous section, a finite
sequence $(m_1, \ldots, m_n)$ of natural numbers; the domain of such a map
is $m_1 + \cdots + m_n$ and the codomain is $n$.  If we view natural
numbers as finite sequences of $\bullet$'s then what a map does is to take
a finite sequence of $\bullet$'s (the domain), partition it into a finite
number of (possibly empty) segments, and replace each segment by a single
$\bullet$ (giving the codomain).  For example, Fig.~\ref{fig:map-in-D}
illustrates the map
\[
(0,3,1,0,2): 6 \go 5
\]
in two different ways: in~(a) as a partition, and in~(b) as a function.

\begin{figure}
\centering
\setlength{\unitlength}{1mm}
\begin{picture}(101,46)(0,3)
% (b)
\cell{0}{12}{l}{\textrm{(b)}}
\cell{10}{6}{bl}{\epsfig{file=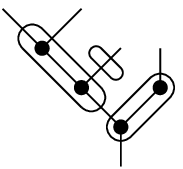}}
\cell{22}{5}{t}{\sigma}
% (a)
\cell{0}{38}{l}{\textrm{(a)}}
\cell{10}{32}{bl}{\epsfig{file=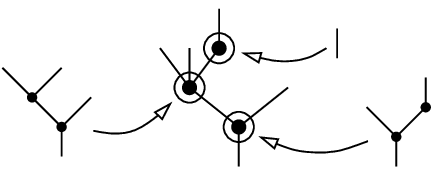}}
\cell{34}{31}{t}{\tau}
% (c)
\cell{68}{24}{l}{\textrm{(c)}}
\cell{78}{6}{bl}{\epsfig{file=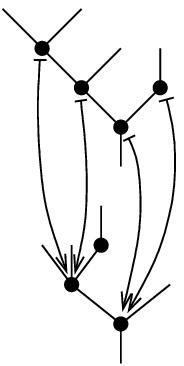}}
\cell{98}{35}{l}{\sigma}
\cell{98}{10}{l}{\tau}
\end{picture}
% \hand{100}{31}
\caption{Three pictures of a map in $\Tr(4)$}
\label{fig:map-in-Tr}
\end{figure}

Ordinarily $\scat{D} = \Pd{1}$ is described as the category of finite
totally ordered sets, but our new description leads smoothly into a
description of the category $\Tr = \Pd{2}$ of trees.  $\Tr$ is the disjoint
union $\coprod_{n\in\nat} \Tr(n)$.  An object of $\Tr(n)$ is an $n$-leafed
tree.  The set of maps in $\Tr(n)$ is 
\[
(T_2^2 1)(n) = (T_2(\tr))(n),
\]
that is, a map is an $n$-leafed tree $\tau$ in which each $k$-ary vertex
$v$ has assigned to it a $k$-leafed tree $\sigma_v$; the domain of the map
is the tree obtained by gluing the $\sigma_v$'s together in the way
dictated by the shape of $\tau$, and the codomain is $\tau$ itself.  Put
another way, what a map does is to take a tree $\sigma$ (the domain),
partition it into a finite number of (possibly trivial) subtrees, and
replace each of these subtrees by the corolla
\[
\begin{centredpic}
\begin{picture}(3,2)(-1.5,0)
% lower layer
\put(0,0){\line(0,1){1}}
\cell{0}{1}{c}{\vx}
% upper layer
\put(0,1){\line(-3,2){1.5}}
\cell{0}{1.8}{c}{\cdots}
\put(0,1){\line(3,2){1.5}}
\end{picture}
\end{centredpic}
\]
with the same number of leaves, to give the codomain $\tau$.
Fig.~\ref{fig:map-in-Tr}
depicts a certain map $\sigma \go \tau$ in $\Tr(4)$ in three different
ways: in~(a) as a 4-leafed tree $\tau$ with a $k$-leafed tree $\sigma_v$
assigned to each $k$-ary vertex $v$, in~(b) as a 4-leafed tree $\sigma$
partitioned into subtrees $\sigma_v$, and in~(c) as something looking more
like a function.  We will return to the third point of view later; for now,
just observe that there is an induced function from the vertices
of $\sigma$ to the vertices of $\tau$, in which the inverse image of a
vertex $v$ of $\tau$ is the set of vertices of $\sigma_v$.

\begin{figure}
\centering
\setlength{\unitlength}{1mm}
\begin{picture}(114,56)
% (b)
\cell{0}{11}{l}{\textrm{(b)}}
\cell{8}{3}{bl}{\epsfig{file=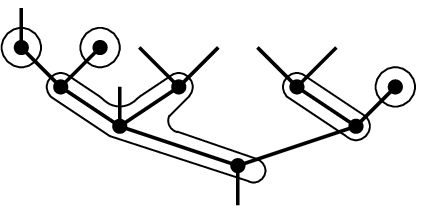}}
\cell{32}{2}{t}{\sigma}
% (a)
\cell{0}{43}{l}{\textrm{(a)}}
\cell{9}{35}{bl}{\epsfig{file=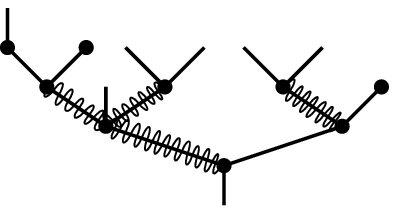}}
\cell{32}{34}{t}{\sigma}
% (c)
\cell{58}{28}{l}{\textrm{(c)}}
\cell{66}{6}{bl}{\epsfig{file=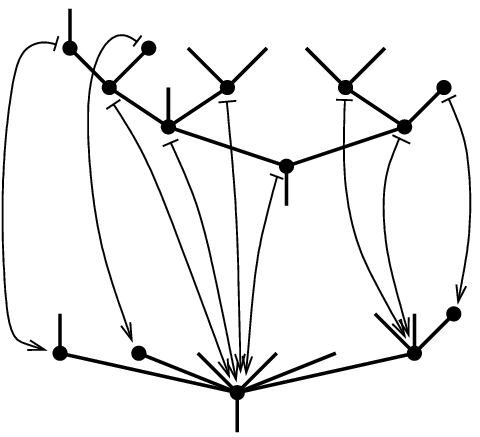}}
\cell{95.5}{27.5}{c}{\sigma}
\cell{90}{5}{t}{\tau}
\end{picture}
% \hand{90}{32}
\caption{Three pictures of an epic in $\Tr(6)$}
\label{fig:epi-in-Tr}
\end{figure}

In some texts a map of trees is described as something that `contracts some
internal edges'.  (Here an \demph{internal%
\index{edge!internal}
edge} is an edge that is not the
root or a leaf; maps of trees keep the root and leaves fixed.  To
`contract'%
\index{contraction!edge of tree@of edge of tree}%
\index{edge!contraction of}
an internal edge means to shrink it down to a vertex.)  With one
important caveat, this is what our maps of trees do: for in a map $\sigma
\go \tau$, the replacement of each partitioning subtree $\sigma_v$ by the
corolla with the same number of leaves amounts to the contraction of all
the internal edges of $\sigma_v$.  For example, Fig.~\ref{fig:epi-in-Tr}(a)
shows a tree $\sigma$ with some of its edges marked for contraction, and
Figs.~\ref{fig:epi-in-Tr}(b) and ~\ref{fig:epi-in-Tr}(c) show the
corresponding maps $\sigma \go \tau$ in two different styles (as in
Figs.~\ref{fig:map-in-Tr}(b) and~(c)); so $\tau$ is the tree obtained by
contracting the marked edges of $\sigma$.  

The caveat is that some of the $\sigma_v$'s may be the trivial tree, and
these are replaced by the 1-leafed corolla 
$
\setlength{\unitlength}{1em}
\begin{array}{c}
\begin{picture}(0,2)(0,0)
\put(0,0){\line(0,1){1}}
\cell{0}{1}{c}{\vx}
\put(0,1){\line(0,1){1}}
\end{picture}
\end{array}.
$
This does \emph{not} amount to the contraction of internal edges: it is,
rather, the addition%
\index{vertex!addition of}%
of a vertex to the middle of a (possibly external)
edge.  Any map of trees can be viewed as a combination of contractions of
internal edges and additions of vertices to existing edges.  For example,
the map illustrated in Fig.~\ref{fig:map-in-Tr} contracts two internal
edges and adds a vertex to one edge.

Analogously, any map in the category $\scat{D}$ of finite totally ordered
sets can be viewed as a combination of merging adjacent $\bullet$'s and
adding new $\bullet$'s (Fig.~\ref{fig:map-in-D}); this amounts to the
factorization%
\index{factorization of map}
of any map as a surjection followed by an injection.  So
those who define their maps between trees to be just contractions of
internal edges are doing something analogous to considering only the
surjective maps in the augmented simplex category $\scat{D}$.  Indeed, the
full subcategory of $\Tr$ consisting of just those trees in which every
vertex has exactly one edge coming up out of it is isomorphic to
$\scat{D}$; if we take only maps made out of contractions of internal edges
then we obtain the subcategory of $\scat{D}$ consisting of surjections
only.

We will come back soon to this issue of surjections and injections in
$\Tr$, with more precision.

\paragraph*{}

Some further understanding of the category of trees can be gained by
considering just those trees in which each vertex has at least two branches
coming up out of it.  I will call these `stable trees', following
Kontsevich and Manin~\cite[Definition 6.6.1]{KMGWC}.  Formally,
$\fcat{StTr}(n)$%
\glo{StTr}
is the full subcategory of $\Tr(n)$ with objects defined
by the recursive clauses
\begin{itemize}
\item $\utree\in\fcat{StTr}(1)$
\item if $n \geq 2$, $k_1, \ldots, k_n \in \nat$, and $\tau_1 \in
\fcat{StTr}(k_1), \ldots, \tau_n \in \fcat{StTr}(k_n)$ then $(\tau_1,
\ldots, \tau_n) \in \fcat{StTr}(k_1 + \cdots + k_n)$,
\end{itemize}
and an \demph{$n$-leafed stable tree}%
\index{tree!stable}
is an object of $\fcat{StTr}(n)$.
Since a stable tree can contain no subtree of the form
$
\setlength{\unitlength}{1em}
\begin{array}{c}
\begin{picture}(0,2)(0,0)
\put(0,0){\line(0,1){1}}
\cell{0}{1}{c}{\vx}
\put(0,1){\line(0,1){1}}
\end{picture}
\end{array}
$,
all maps between stable trees are `surjections', that is,
consist of just contractions of internal edges, without insertions of new
vertices.  It follows that each category $\fcat{StTr}(n)$ is finite, and so
its classifying space can be represented by a finite CW complex; this may
explain why topologists often like their trees to be stable.

\begin{figure}
\centering
\setlength{\unitlength}{1mm}
\begin{picture}(106,10)
% (a)
\cell{23}{3}{t}{\textrm{(a)}}
\cell{0}{4}{bl}{\epsfig{file=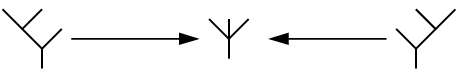}}
% (b)
\cell{86}{3}{t}{\textrm{(b)}}
\thicklines
\put(66,7){\line(1,0){40}}
\thinlines
\end{picture}
% \hand{25}{33}
\caption{(a) The category of 3-leafed stable trees, and~(b) its classifying
  space} 
\label{fig:stable-three}
\end{figure}

\begin{figure}
\centering
\setlength{\unitlength}{1mm}
\begin{picture}(97,42)
% (a)
\cell{21.5}{3}{t}{\textrm{(a)}}
\cell{0}{4}{bl}{\epsfig{file=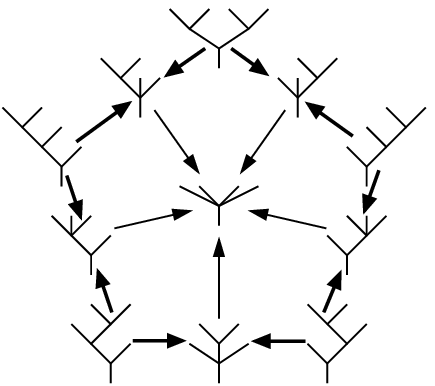}}
% (b)
\cell{80}{3}{t}{\textrm{(b)}}
\cell{63}{7.5}{bl}{\epsfig{file=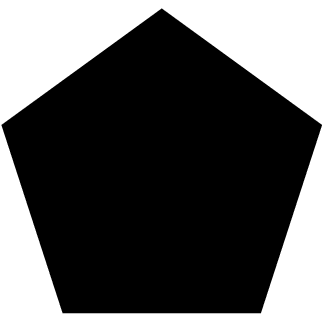}}
\end{picture}
% \hand{60}{34}
\caption{(a) The category of 4-leafed stable trees, and~(b) its classifying
  space}%
\index{pentagon}%
\index{associahedron}
\label{fig:stable-four}
\end{figure}

\begin{figure}
\centering
\begin{tabular}{c}
\epsfig{file=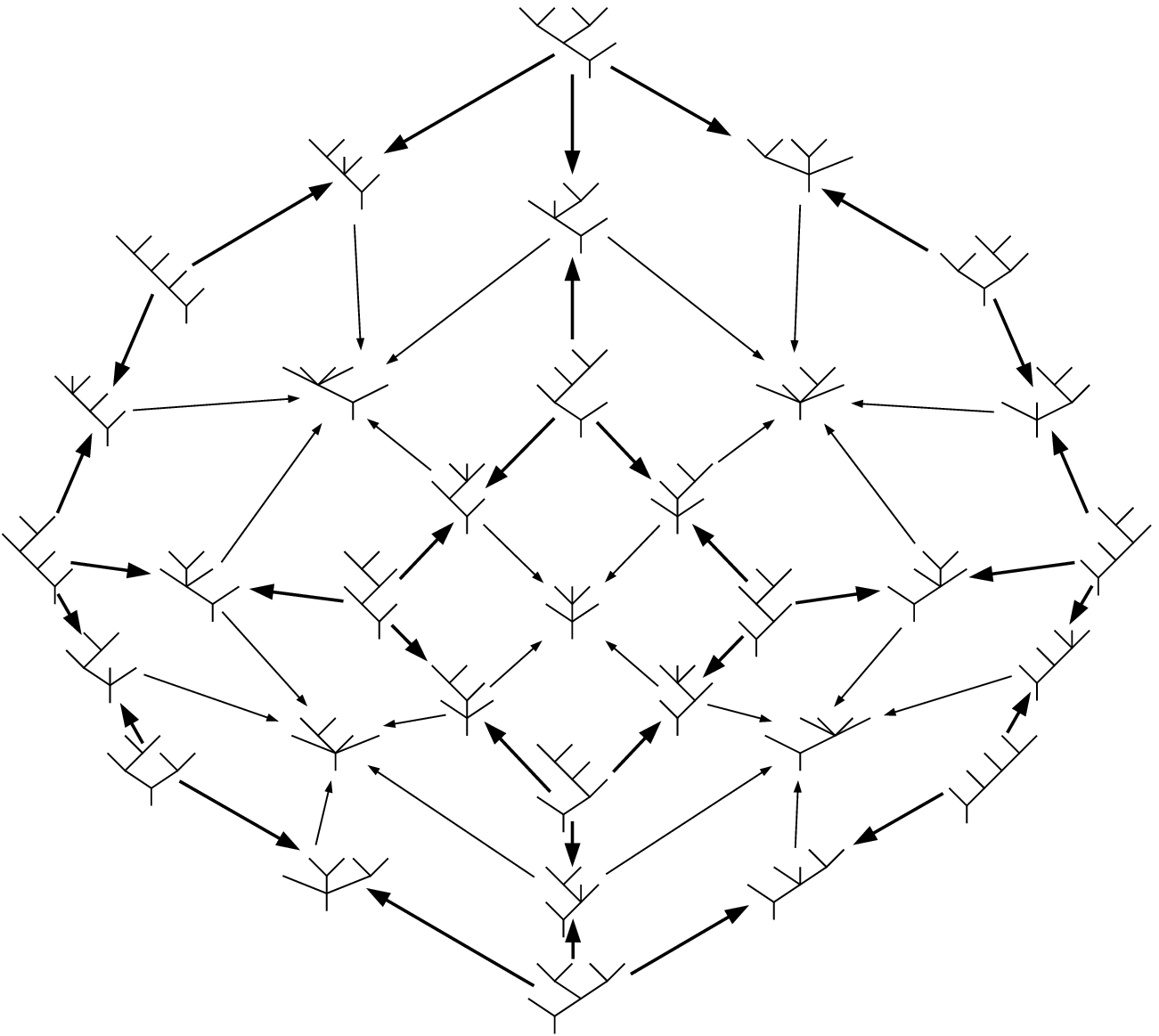,width=0.95\textwidth} \\
(a)\\
\\
\epsfig{file=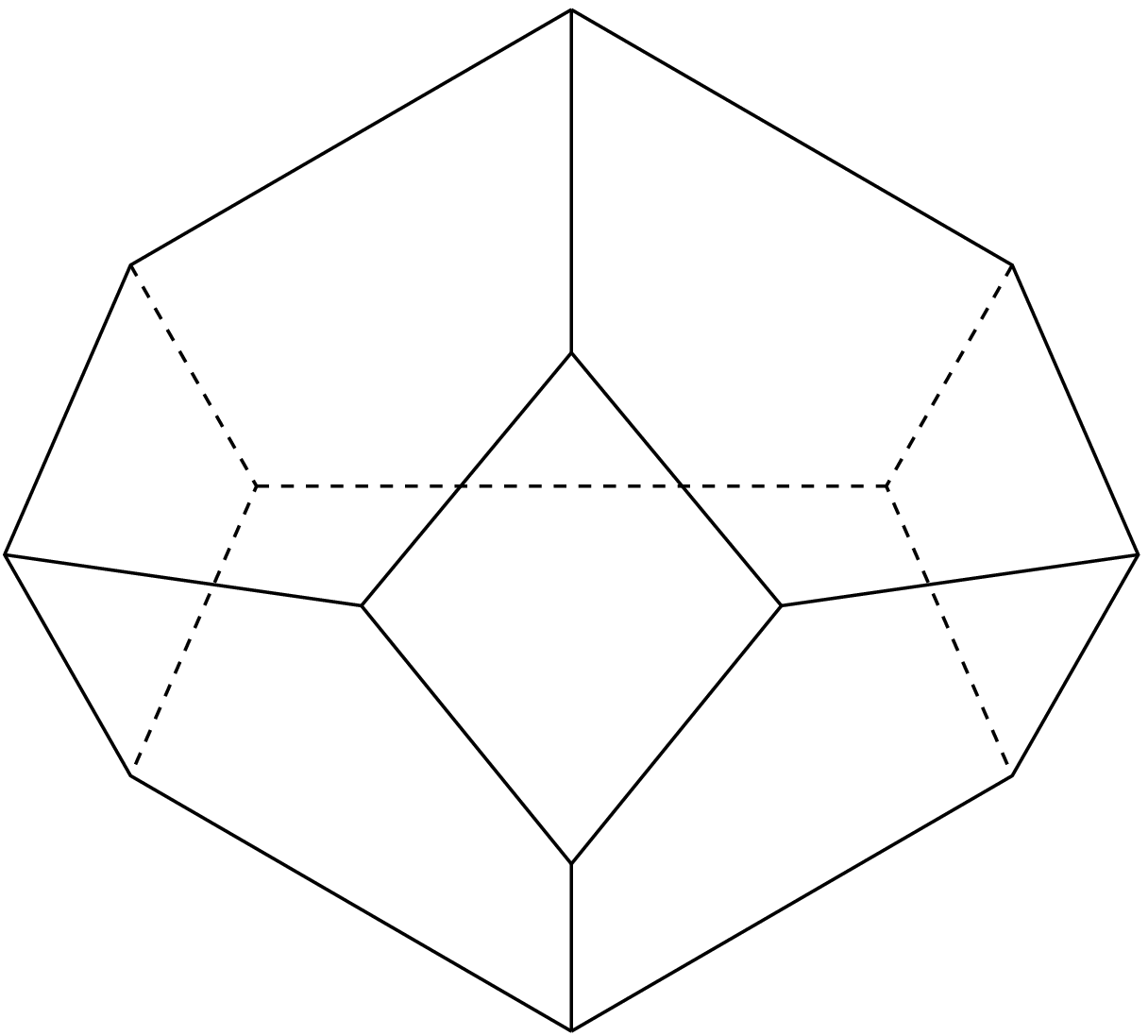,width=0.5\textwidth} \\
(b)
\end{tabular}
% \hand{170}{35}
\caption{(a) About half of the category of 5-leafed stable trees, and~(b)
the classifying space of the whole category}%
\index{associahedron}
\label{fig:stable-five}
\end{figure}

The first few categories $\fcat{StTr}(n)$ are trivial:
\begin{eqnarray*}
\fcat{StTr}(0)	&=	&\emptyset,	\\
\fcat{StTr}(1)	&=	&\{ \utree \},	\\
\fcat{StTr}(2)	&=	&
\left\{ 
% \drmk{pic of two-leafed corolla} 
\begin{centredpic}
\begin{picture}(2,2)(-1,0)
% lower layer
\put(0,0){\line(0,1){1}}
% upper layer
\cell{0}{1}{c}{\vx}
\put(0,1){\line(-1,1){1}}
\put(0,1){\line(1,1){1}}
\end{picture}
\end{centredpic}
\right\},
\end{eqnarray*}
where in each case there are no arrows except for identities.  The cases $n
= 3$, $4$, and $5$ are illustrated in Figs.~\ref{fig:stable-three}(a),
\ref{fig:stable-four}(a), and~\ref{fig:stable-five}(a).
Identity arrows are not shown, and the categories $\fcat{StTr}(n)$ are
ordered sets: all diagrams commute.  Vertices are also omitted; since the
trees are stable, this does not cause ambiguity.  Parts~(b) of the figures
show the classifying spaces of these categories, solid polytopes of
dimensions $1$, $2$ and $3$.  In the case of 5-leafed trees
(Fig.~\ref{fig:stable-five}) only about half of the category is shown,
corresponding to the front faces of the polytope; the back faces and the
terminal object of the category (the 5-leafed corolla), which sits at the
centre of the polytope, are hidden.  The whole polytope has 6 pentagonal
faces, 3 square faces, and 3-fold rotational symmetry about the central
vertical axis.

For $n\leq 5$, the classifying space $B(\fcat{StTr}(n))$ is homeomorphic to
the associahedron%
\index{associahedron}
$K_n$ (Stasheff~\cite{StaHAHI}%
\index{Stasheff, Jim}
and~\ref{eg:opd-associahedra} above), and it seems very likely that this
persists for all $n\in\nat$.  Indeed, the family of categories
$(\fcat{StTr}(n))_{n\in\nat}$ forms a sub-\Cat-operad $\fcat{STTR}$ of
$\fcat{TR} = \PD{2}$,%
\glo{TR}%
\index{tree!Cat-operad of@$\Cat$-operad of}
and the classifying%
\index{classifying space}
space functor $B: \Cat \go
\Top$ preserves finite products, so there is a (non-symmetric) topological
operad $B(\fcat{STTR})$ whose $n$th part is the classifying space of
$\fcat{StTr}(n)$.  (To make $B$ preserve finite products we must interpret
$\Top$ as the category of compactly generated or Kelley spaces: see
Segal~\cite[\S 1]{SegCSS} and Gabriel and Zisman~\cite[III.2]{GZ}.)  This
operad $B(\fcat{STTR})$ is presumably isomorphic to Stasheff's operad $K =
(K_n)_{n\in\nat}$.  A $K$-algebra is called an \demph{$A_\infty$-space},%
\index{A-@$A_\infty$-!space}
and should be thought of as an up-to-homotopy version of a topological
semigroup; the basic example is a loop space.

The categories $\fcat{StTr}(n)$ also give rise to the notion of an
$A_\infty$-algebra%
\lbl{p:A-infty-alg}
(Stasheff~\cite{StaHAHII}).%
\index{Stasheff, Jim}
 For each $n\in\nat$, there is a chain complex
$P(n)$ whose degree $k$ part is the free abelian group on the set of
$n$-leafed stable trees with $(n-k-1)$ vertices.  For instance,
\[
P(4) = 
(\cdots\rTo 0 \rTo 0 \rTo
\integers\cdot L_2
\rTo
\integers\cdot L_1 
\rTo
\integers\cdot L_0)
\]
where the sets $L_k$ are
\begin{eqnarray*}
L_2	&=	&
\left\{ 
\begin{centredpic}
\begin{picture}(3,2)(-1.5,0)
% lower layer
\put(0,0){\line(0,1){1}}
\cell{0}{1}{c}{\vx}
% upper layer
\put(0,1){\line(-3,2){1.5}}
\put(0,1){\line(-1,2){0.5}}
\put(0,1){\line(1,2){0.5}}
\put(0,1){\line(3,2){1.5}}
\end{picture}
\end{centredpic}
\right\}	\\
L_1	&=	&
\left\{ 
\begin{centredpic}
\begin{picture}(2.5,3)(-1.5,0)
% bottom layer
\put(0,0){\line(0,1){1}}
% middle layer
\cell{0}{1}{c}{\vx}
\put(0,1){\line(-1,1){1}}
\put(0,1){\line(0,1){1}}
\put(0,1){\line(1,1){1}}
% top layer
\cell{-1}{2}{c}{\vx}
\put(-1,2){\line(-1,2){0.5}}
\put(-1,2){\line(1,2){0.5}}
\end{picture}
\end{centredpic},
\begin{centredpic}
\begin{picture}(2,3)(-1,0)
% bottom layer
\put(0,0){\line(0,1){1}}
% middle layer
\cell{0}{1}{c}{\vx}
\put(0,1){\line(-1,1){1}}
\put(0,1){\line(0,1){1}}
\put(0,1){\line(1,1){1}}
% top layer
\cell{0}{2}{c}{\vx}
\put(0,2){\line(-1,2){0.5}}
\put(0,2){\line(1,2){0.5}}
\end{picture}
\end{centredpic},
\begin{centredpic}
\begin{picture}(2.5,3)(-1,0)
% bottom layer
\put(0,0){\line(0,1){1}}
% middle layer
\cell{0}{1}{c}{\vx}
\put(0,1){\line(-1,1){1}}
\put(0,1){\line(0,1){1}}
\put(0,1){\line(1,1){1}}
% top layer
\cell{1}{2}{c}{\vx}
\put(1,2){\line(-1,2){0.5}}
\put(1,2){\line(1,2){0.5}}
\end{picture}
\end{centredpic},
\begin{centredpic}
\begin{picture}(2.75,3)(-1.75,0)
% bottom layer
\put(0,0){\line(0,1){1}}
% middle layer
\cell{0}{1}{c}{\vx}
\put(0,1){\line(-1,1){1}}
\put(0,1){\line(1,1){1}}
% top layer
\cell{-1}{2}{c}{\vx}
\put(-1,2){\line(-3,4){0.75}}
\put(-1,2){\line(0,1){1}}
\put(-1,2){\line(3,4){0.75}}
\end{picture}
\end{centredpic},
\begin{centredpic}
\begin{picture}(2.75,3)(-1,0)
% bottom layer
\put(0,0){\line(0,1){1}}
% middle layer
\cell{0}{1}{c}{\vx}
\put(0,1){\line(-1,1){1}}
\put(0,1){\line(1,1){1}}
% top layer
\cell{1}{2}{c}{\vx}
\put(1,2){\line(-3,4){0.75}}
\put(1,2){\line(0,1){1}}
\put(1,2){\line(3,4){0.75}}
\end{picture}
\end{centredpic}
\right\}	\\
L_0	&=	&
\left\{ 
\begin{centredpic}
\begin{picture}(4,4)(-3,0)
% bottom layer
\put(0,0){\line(0,1){1}}
% second-bottom layer
\cell{0}{1}{c}{\vx}
\put(0,1){\line(-1,1){1}}
\put(0,1){\line(1,1){1}}
% second-top layer
\cell{-1}{2}{c}{\vx}
\put(-1,2){\line(-1,1){1}}
\put(-1,2){\line(1,1){1}}
% top layer
\cell{-2}{3}{c}{\vx}
\put(-2,3){\line(-1,1){1}}
\put(-2,3){\line(1,1){1}}
\end{picture}
\end{centredpic},
\begin{centredpic}
\begin{picture}(3,4)(-2,0)
% bottom layer
\put(0,0){\line(0,1){1}}
% second-bottom layer
\cell{0}{1}{c}{\vx}
\put(0,1){\line(-1,1){1}}
\put(0,1){\line(1,1){1}}
% second-top layer
\cell{-1}{2}{c}{\vx}
\put(-1,2){\line(-1,1){1}}
\put(-1,2){\line(1,1){1}}
% top layer
\cell{0}{3}{c}{\vx}
\put(0,3){\line(-1,1){1}}
\put(0,3){\line(1,1){1}}
\end{picture}
\end{centredpic},
\begin{centredpic}
\begin{picture}(4,3)(-2,0)
% bottom layer
\put(0,0){\line(0,1){1}}
% middle layer
\cell{0}{1}{c}{\vx}
\put(0,1){\line(-3,2){1.5}}
\put(0,1){\line(3,2){1.5}}
% top layer
\cell{-1.5}{2}{c}{\vx}
\put(-1.5,2){\line(-1,1){1}}
\put(-1.5,2){\line(1,1){1}}
\cell{1.5}{2}{c}{\vx}
\put(1.5,2){\line(-1,1){1}}
\put(1.5,2){\line(1,1){1}}
\end{picture}
\end{centredpic},
\begin{centredpic}
\begin{picture}(3,4)(-1,0)
% bottom layer
\put(0,0){\line(0,1){1}}
% second-bottom layer
\cell{0}{1}{c}{\vx}
\put(0,1){\line(-1,1){1}}
\put(0,1){\line(1,1){1}}
% second-top layer
\cell{1}{2}{c}{\vx}
\put(1,2){\line(-1,1){1}}
\put(1,2){\line(1,1){1}}
% top layer
\cell{0}{3}{c}{\vx}
\put(0,3){\line(-1,1){1}}
\put(0,3){\line(1,1){1}}
\end{picture}
\end{centredpic},
\begin{centredpic}
\begin{picture}(4,4)(-1,0)
% bottom layer
\put(0,0){\line(0,1){1}}
% second-bottom layer
\cell{0}{1}{c}{\vx}
\put(0,1){\line(-1,1){1}}
\put(0,1){\line(1,1){1}}
% second-top layer
\cell{1}{2}{c}{\vx}
\put(1,2){\line(-1,1){1}}
\put(1,2){\line(1,1){1}}
% top layer
\cell{2}{3}{c}{\vx}
\put(2,3){\line(-1,1){1}}
\put(2,3){\line(1,1){1}}
\end{picture}
\end{centredpic}
\right\}.	
\end{eqnarray*}
The differential $d$ is defined by $d(\tau) = \sum \pm \sigma$, where if
$\tau \in P(n)_k$ then the sum is over all $\sigma\in P(n)_{k-1}$ for which
there exists a map $\sigma \go \tau$.  For instance,
\[
d
% (\drmk{pic of ((ab)cd)})
\left(
\begin{centredpic}
\begin{picture}(2.5,3)(-1.5,0)
% bottom layer
\put(0,0){\line(0,1){1}}
% middle layer
\cell{0}{1}{c}{\vx}
\put(0,1){\line(-1,1){1}}
\put(0,1){\line(0,1){1}}
\put(0,1){\line(1,1){1}}
% top layer
\cell{-1}{2}{c}{\vx}
\put(-1,2){\line(-1,2){0.5}}
\put(-1,2){\line(1,2){0.5}}
\end{picture}
\end{centredpic}
\right)
=
\pm 
% \drmk{pic of (((ab)c)d)} 
\begin{centredpic}
\begin{picture}(4,4)(-3,0)
% bottom layer
\put(0,0){\line(0,1){1}}
% second-bottom layer
\cell{0}{1}{c}{\vx}
\put(0,1){\line(-1,1){1}}
\put(0,1){\line(1,1){1}}
% second-top layer
\cell{-1}{2}{c}{\vx}
\put(-1,2){\line(-1,1){1}}
\put(-1,2){\line(1,1){1}}
% top layer
\cell{-2}{3}{c}{\vx}
\put(-2,3){\line(-1,1){1}}
\put(-2,3){\line(1,1){1}}
\end{picture}
\end{centredpic}
\pm 
% \drmk{pic of ((ab)(cd))}.
\begin{centredpic}
\begin{picture}(4,3)(-2,0)
% bottom layer
\put(0,0){\line(0,1){1}}
% middle layer
\cell{0}{1}{c}{\vx}
\put(0,1){\line(-3,2){1.5}}
\put(0,1){\line(3,2){1.5}}
% top layer
\cell{-1.5}{2}{c}{\vx}
\put(-1.5,2){\line(-1,1){1}}
\put(-1.5,2){\line(1,1){1}}
\cell{1.5}{2}{c}{\vx}
\put(1.5,2){\line(-1,1){1}}
\put(1.5,2){\line(1,1){1}}
\end{picture}
\end{centredpic}.
\]
When the signs are chosen appropriately this defines an operad $P$ of chain
complexes.  A $P$-algebra is called an \demph{$A_\infty$-algebra},%
\index{A-@$A_\infty$-!algebra}
to be thought of as an up-to-homotopy differential%
\index{algebra!differential graded}
graded non-unital algebra; the
usual example is the singular chain complex of an $A_\infty$-space.  A
$P$-category is called an $A_\infty$-category%
\lbl{p:A-infty-category}%
\index{A-@$A_\infty$-!category}
(see~\ref{eg:fc-A-infty}), and consists of a collection of objects, a chain
complex $\Hom(a,b)$ for each pair $(a,b)$ of objects, maps defining binary
composition, chain homotopies witnessing that this composition is
associative up to homotopy, further homotopies witnessing that the previous
homotopies obey the pentagon law up to homotopy, and so on.
$A_\infty$-categories are cousins of weak $\omega$-categories, as we see
in~\ref{sec:non-alg-defns-n-cat}.

Finally, since the polytopes $K_n = B(\fcat{StTr}(n))$ describe higher
associativity%
\index{associativity}
conditions, they also arise in definitions of
higher-dimensional category.  For example, the pentagon%
\index{pentagon}
$K_4$ occurs in the
classical definition of bicategory~(\ref{defn:cl-bicat}), and the
polyhedron $K_5$ occurs as the `non-abelian%
\index{non-abelian 4-cocycle}
4-cocycle condition' in Gordon,%
\index{Gordon, Robert}
Power and Street's definition of tricategory~\cite{GPS}.%
\index{tricategory}

\paragraph*{}

We have already described the set $\tr(n)$ of $n$-leafed trees.  Maps
$\sigma \go \tau$ between trees are described by induction on the structure
of $\tau$:
\begin{itemize}
\item if $\tau=\utree$ then there is only one map into $\tau$; it has
domain $\utree$ and we write it as $1_\utree: \utree \go \utree$
\item if $\tau = (\tau_1, \ldots, \tau_n)$ for $\tau_1 \in \tr(k_1)$,
\ldots, $\tau_n \in \tr(k_n)$ then a map $\sigma \go \tau$ consists of
trees $\rho \in \tr(n), \rho_1 \in \tr(k_1), \ldots, \rho_n \in \tr(k_n)$
such that $\sigma = \rho \of (\rho_1, \ldots, \rho_n)$, together with maps
\[
\rho_1 \goby{\theta_1} \tau_1, 
\ \ \ldots,\ \  
\rho_n \goby{\theta_n} \tau_n,
\]
and we write this map as
\begin{equation}	\label{eq:form-of-map-of-trees}
\sigma = \rho \of (\rho_1, \ldots, \rho_n)
\goby{!_\rho * (\theta_1, \ldots, \theta_n)}
(\tau_1, \ldots, \tau_n) = \tau.
\end{equation}
\end{itemize}
It follows easily that the $n$-leafed corolla $\nu_n = (\utree, \ldots,
\utree)$ is the terminal object of $\Tr(n)$: the unique map from
$\sigma\in\tr(n)$ to $\nu_n$ is $!_\sigma * (1_\utree, \ldots, 1_\utree)$.

The rest of the structure of the $\Cat$-operad $\TR$ can be described in a
similarly explicit recursive fashion.  

To make precise the intuition that a map of trees is a function of
some sort, functors
\[
V: \Tr \go \Set,
\diagspace
E: \Tr^\op \go \Set%
\glo{vxftr}\glo{edgeftr}%
\index{vertex!functor}%
\index{edge!functor}%
\index{tree!vertex functor}%
\index{tree!edge functor}
\]
can be defined, encoding what happens on vertices and edges respectively.
Both functors turn out to be faithful, which means that a map of trees is
completely determined by its effect on either vertices or edges.  The
following account of $V$ and $E$ is just a sketch.

\begin{figure}
\centering
\setlength{\unitlength}{1mm}
\begin{picture}(88,44)(0,-1)
% (a)
\cell{11}{2}{t}{\textrm{(a)}}
\cell{0}{4}{bl}{\epsfig{file=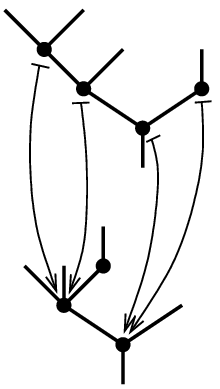}}
\cell{23}{21}{c}{\scriptstyle V(\theta)}
% middle bit
\cell{42}{37}{c}{\sigma}
\put(42,35){\vector(0,-1){23}}
\cell{42}{10}{c}{\tau}
\cell{43}{23.5}{l}{\theta}
% (b)
\cell{77}{2}{t}{\textrm{(b)}}
\cell{66}{4}{bl}{\epsfig{file=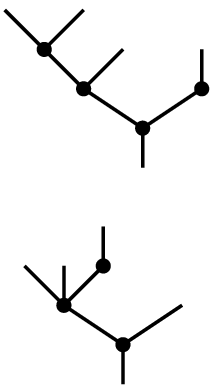}}
% (b): labels on lower tree
\cell{77.5}{6}{c}{\scriptstyle 1}
\cell{74}{10}{c}{\scriptstyle 2}
\cell{69}{14.5}{c}{\scriptstyle 3}
\cell{71.5}{15.5}{c}{\scriptstyle 4}
\cell{75.5}{14}{c}{\scriptstyle 5}
\cell{77.5}{19}{c}{\scriptstyle 6}
\cell{83}{10}{c}{\scriptstyle 7}
% (b): labels on upper tree
\cell{79.5}{28.5}{c}{\scriptstyle 1}
\cell{76}{32}{c}{\scriptstyle 2}
\cell{66.5}{41}{c}{\scriptstyle 3}
\cell{74}{41}{c}{\scriptstyle 4}
\cell{79.5}{37}{c}{\scriptstyle 5,6}
\cell{87.5}{37}{c}{\scriptstyle 7}
% (b): arrow between trees
\put(72,21){\vector(0,1){10}}
\cell{71}{26}{r}{\scriptstyle E(\theta)}
\end{picture}
% \hand{50}{37}
\caption{The effect on~(a) vertices and~(b) edges of a certain map of
4-leafed trees}
\label{fig:edge-functor-trees}
\end{figure}

The more obvious of the two is the vertex functor $V$, defined on objects
by
\begin{itemize}
\item $V(\utree) = \emptyset$
\item $V((\tau_1, \ldots, \tau_n)) = 1 + V(\tau_1) + \cdots + V(\tau_n)$.
\end{itemize}
The edge functor $E$ can be defined by first defining a functor
\[
E_n: \Tr(n)^\op \go (n+1)/\Set
\]
for each $n\in\nat$, where $(n+1)/\Set$ is the category of sets equipped
with $(n+1)$ ordered marked points.  This definition is again by induction,
the idea being that $E_n$ associates to a tree its edge-set with the $n$
input edges and the one output edge (root) distinguished.
Fig.~\ref{fig:edge-functor-trees}
illustrates a map $\theta: \sigma \go \tau$ in $\Tr(4)$; part~(a) ($=$
Fig.~\ref{fig:map-in-Tr}(c)) shows its effect $V(\theta)$ on vertices;
part~(b) shows $E(\theta)$, taking $E(\tau) = \{1, \ldots, 7\}$ and
labelling the image of $i \in \{1, \ldots, 7\}$ under $E(\theta)$ by an
$i$ on the edge $(E(\theta))(i)$ of $\sigma$.

A map of trees will be called surjective if it is built up from
contractions of internal edges (the analogues of degeneracy%
\index{degeneracy map}
maps in
$\scat{D}$).  Formally, the \demph{surjective}%
\index{surjective map of trees}%
\index{tree!map of!surjective}
maps in $\Tr$ are defined
by:
\begin{itemize}
\item $1_\utree: \utree \go \utree$ is surjective
\item with notation as in~\bref{eq:form-of-map-of-trees}, $!_\rho *
(\theta_1, \ldots, \theta_n)$ is surjective if and only if each $\theta_i$
is surjective and $\rho \neq \utree$.
\end{itemize}
The crucial part is the last: the unique map $!_\rho$ from $\rho\in\tr(n)$
to the corolla $\nu_n$ is made up of edge-contractions just as long as
$\rho$ is not the unit tree $\utree$.  

Dually, a map of trees is \demph{injective}%
\index{injective map of trees}%
\index{tree!map of!injective}
if, informally, it is built up
from adding vertices to the middle of edges (the analogues of face%
\index{face map}
maps in
$\scat{D}$).  Formally,
\begin{itemize}
\item $1_\utree: \utree \go \utree$ is injective
\item with notation as above, $!_\rho * (\theta_1, \ldots, \theta_n)$ is
injective if and only if each $\theta_i$ is injective and $\rho$ is either
$\nu_n$ or $\utree$ (the latter only being possible if $n = 1$).
\end{itemize}

The punchline is that the various possible notions of a map of trees being
`onto' (respectively, `one-to-one') all coincide:
\begin{propn}
The following conditions on a map $\theta: \sigma \go \tau$ in $\Tr$ are
equivalent: 
\begin{enumerate}
\item	\lbl{item:tree-epic}
$\theta$ is epic
\item	\lbl{item:tree-surj}
$\theta$ is surjective
\item	\lbl{item:tree-V-surj}
$V(\theta)$ is surjective
\item	\lbl{item:tree-E-inj}
$E(\theta)$ is injective (sic).
\end{enumerate}
Moreover, if each condition is replaced by its dual then the equivalence
persists. 
\end{propn}
\begin{proof}
Omitted.  `Moreover' is not just an application of formal duality, since
surjectivity and injectivity are not formal duals.  
\done
\end{proof}%
\index{tree!map of|)}

\index{pasting diagram!opetopic!category of|(}
We finish this section by re-considering briefly what we have done with
trees, but this time with 2-pasting diagrams instead.  This gives a very
good impression of the category of $n$-pasting diagrams for arbitrary $n$.

\begin{figure}
\centering
\setlength{\unitlength}{1mm}
\begin{picture}(100,30.5)(0,-8)
\cell{0}{0}{bl}{\epsfig{file=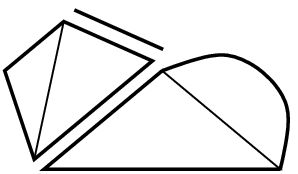}}
\cell{14.5}{-1}{t}{\sigma}
\cell{14.5}{-6}{c}{\textrm{(a)}}
\cell{49}{0}{bl}{\epsfig{file=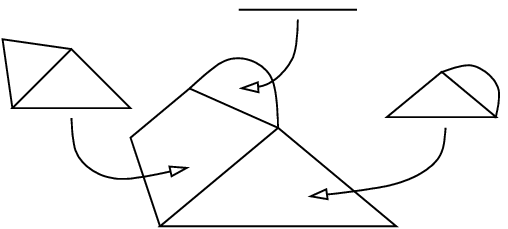}}
\cell{77}{-0.5}{t}{\tau}
\cell{74.5}{-6}{c}{\textrm{(b)}}
\end{picture}
% \hand{30}{41}
\caption{Two pictures of a map in $\Pd{2}$}
\label{fig:Pd-two-map}
\end{figure}

The objects of the category $\Pd{2} = \Tr$ are the opetopic 2-pasting
diagrams.  A map $\theta: \sigma \go \tau$ in $\Pd{2}$ takes a 2-pasting
diagram $\sigma$, partitions it into a finite number of sub-pasting
diagrams, and replaces each of these sub-pasting diagrams by the 2-opetope
\[
\topeqn{}{}{}{}{\Downarrow} 
\]
with the same number of input edges, to give the codomain $\tau$.  Another
way to put this is that each of the 2-opetopes $v$ making up the pasting
diagram $\tau$ has assigned to it a 2-pasting diagram $\sigma_v$ with the
same number of input edges, and when the $\sigma_v$'s are pasted
together according to the shape of $\tau$, the result is the pasting
diagram $\sigma$.  Fig.~\ref{fig:Pd-two-map}
shows a map in $\Pd{2}$, in~(a) as a partition of $\sigma$ and in~(b) as a
family $(\sigma_v)$ indexed over the regions $v$ of $\tau$; these
correspond precisely to the tree pictures in Figs.~\ref{fig:map-in-Tr}(b)
and~(a) respectively.  

More generally, a map of $n$-pasting diagrams consists of the replacement
of some sub-pasting diagrams by their bounding opetopes.  When the
sub-pasting diagrams are non-trivial, this amounts to the removal of some
internal faces of codimension one.  For $n=2$, faces of codimension one are
edges, and the trivial case is the replacement of the 2-pasting diagram
\[
\topebasen{}
\]
by the 2-opetope
\[
\topean{}{}{\Downarrow}
.
\]
In Fig.~\ref{fig:Pd-two-map}, there are two instances of edge-deletion and
one instance of inflating the trivial 2-pasting diagram.  This is the
distinction between epics%
\index{tree!map of!surjective}%
\index{surjective map of trees}
(degeneracy%
\index{degeneracy map}
maps) and monics%
\index{tree!map of!injective}%
\index{injective map of trees}
(face%
\index{face map}
maps) in
$\Pd{2}$, as we saw for trees.

The classical connections between trees and $A_\infty$-structures%
\index{A-@$A_\infty$-!algebra}%
\index{A-@$A_\infty$-!space}%
\index{A-@$A_\infty$-!category}
can, of
course, be phrased equally in terms of 2-pasting diagrams.  This is
probably the natural approach if we want to incorporate
$A_\infty$-structures into higher-dimensional algebra.  Stable%
\index{tree!stable}
trees
correspond to 2-pasting diagrams not containing any copies of either of the
2-opetopes
\[
\topezn{}{\Downarrow},
\diagspace
\topean{}{}{\Downarrow}
\]
---in other words, those that can be drawn using only straight lines.

\begin{figure}
\centering
\setlength{\unitlength}{1mm}
\begin{picture}(51,40)(-23,0)
\cell{0}{0}{bl}{\epsfig{file=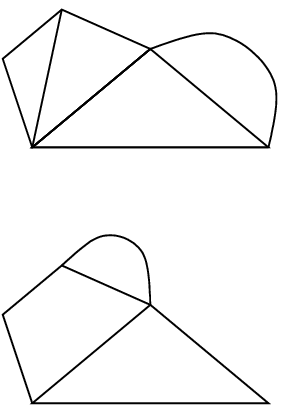}}
% lower half
%   edge labels
\cell{15}{1.5}{c}{\scriptstyle 1}
\cell{10}{5}{c}{\scriptstyle 2}
\cell{1}{4}{c}{\scriptstyle 3}
\cell{3}{13}{c}{\scriptstyle 4}
\cell{9}{12}{c}{\scriptstyle 5}
\cell{15}{16}{c}{\scriptstyle 6}
\cell{21.5}{6}{c}{\scriptstyle 7}
%   region labels
\cell{15}{5}{c}{\scriptstyle a,d}
\cell{6}{8}{c}{\scriptstyle b,c}
% connecting arrows
\put(13,24){\vector(0,-1){5}}
\cell{12}{21.5}{r}{\scriptstyle V(\theta)}
\put(17,19){\vector(0,1){5}}
\cell{18}{21.5}{l}{\scriptstyle E(\theta)}
% upper half
%   edge labels
\cell{15}{27.5}{c}{\scriptstyle 1}
\cell{10}{31}{c}{\scriptstyle 2}
\cell{1}{30}{c}{\scriptstyle 3}
\cell{3}{39}{c}{\scriptstyle 4}
\cell{13}{39}{c}{\scriptstyle 5,6}
\cell{26}{37}{c}{\scriptstyle 7}
%   region labels
\cell{15}{31}{c}{\scriptstyle a}
\cell{9}{35}{c}{\scriptstyle b}
\cell{3}{35}{c}{\scriptstyle c}
\cell{23}{33}{c}{\scriptstyle d}
% Arrow on LHS
\cell{-20}{31}{c}{\sigma}
\cell{-20}{5}{c}{\tau}
\put(-20,28){\vector(0,-1){20}}
\cell{-21}{18}{r}{\theta}
\end{picture}
% \hand{50}{42}
\caption{The effects of $V$ and $E$ on a map of 2-pasting diagrams}
\label{fig:V-and-E-for-two-Pd}
\end{figure}

Because of the inversion of dimensions, the `vertex functor' $V: \Tr =
\Pd{2} \go \Set$ assigns to a 2-pasting diagram its set of faces, or
regions, or constituent 2-opetopes.  (Compare~\ref{eg:constituents}.)  The
`edge functor' $E: \Pd{2}^\op \go \Set$ is still aptly named.  Nothing
interesting happens on the \emph{actual} vertices of 2-pasting diagrams.
Fig.~\ref{fig:V-and-E-for-two-Pd}
shows the effects of $V$ and $E$ on the map $\theta: \sigma \go \tau$ of
Fig.~\ref{fig:Pd-two-map}: the regions of $\sigma$ are labelled $a, b, c,
d$, and the image of the region $a$ under the function $V(\theta)$ is also
labelled $a$, and so on; similarly, the edges of $\tau$ are labelled $1,
\ldots, 7$ and their images under $E(\theta)$ are labelled correspondingly.
The same example was also shown in Fig.~\ref{fig:edge-functor-trees}, using
trees.

The pictures of 2-pasting diagrams can be taken seriously, that is,
geometrically realized.  This leads quickly into the geometric realization
of opetopic sets, and so to the underlying (`singular') opetopic set of a
topological space, one of the motivating examples of a weak
$\omega$-category.  We come to this in the next section.%
\index{pasting diagram!opetopic!category of|)}

\section{Opetopic sets}
\lbl{sec:ope-sets}

Opetopes were defined by Baez and Dolan in order to give a definition of
weak $n$-category.  Their definition
has been subject to various
modifications by various other people, all of the form `a weak $n$-category
is an opetopic set with certain properties'.  The next two sections are a
discussion of the general features of such definitions, not concentrating
on any version in particular.

In this section I will explain what an opetopic set is.  Again, there are
various proposed definitions,%
\index{opetopic!set!definitions of}
most of which have been proved equivalent by
Cheng%
\index{Cheng, Eugenia}
(see the Notes).  Rather than giving any particular one of them, I
will list some properties satisfied by the category $\scat{O}$%
\glo{catofopes}%
\index{opetope!category of}
of opetopes,
an opetopic set being a presheaf on $\scat{O}$.  Using this, I will show
how every topological space has an underlying opetopic set, to be thought
of as its `singular opetopic set' or `fundamental $\omega$-groupoid'.  This
will motivate the definition of weak $n$-category in the next section.

Opetopic sets should be thought of as something like simplicial sets.  A
simplicial set is a presheaf on the category $\Delta$ of simplices; an
opetopic set is a presheaf on the category $\scat{O}$ of opetopes.
Actually, it might be more apposite to compare opetopic sets to presheaves
on the category $\Delta_\mr{inj}$%
\glo{Deltainj}
of nonempty finite totally ordered sets
and order-preserving injections, rather than whole simplicial sets: we only
consider face%
\index{face map}
maps between opetopes, not degeneracies.%
\index{degeneracy map}

Here is an informal description of the category $\scat{O}$ of opetopes.
The set of objects is the set $\coprod_{n\in\nat} O_n$ of all opetopes of
all dimensions.  A map $\omega' \go \omega$ is an embedding of $\omega'$ as
a face of $\omega$.  For example, there are four maps in $\scat{O}$ from
the unique 1-opetope
\[
\topebase{}
\]
into the 2-opetope
\[
\omega = 
\topec{}{}{}{}{\Downarrow},
\]
corresponding to the three input edges and one output edge of $\omega$.
There are also four maps from the unique 0-opetope $\,\gzersu\,$ into
$\omega$, corresponding to its four vertices.  Along with the identity
$1_\omega$, that enumerates all of the maps into $\omega$ in $\scat{O}$.
Similarly, there are 21 maps whose codomain is the 3-opetope $\omega$
illustrated in Fig.~\ref{fig:3-ope-as-cod}:
\begin{figure}
\[
\begin{array}{c}
\setlength{\unitlength}{1mm}
\begin{picture}(38,21)
\cell{0}{0}{bl}{\epsfig{file=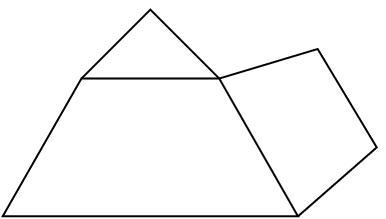}}
\da{15}{7}{0}
\da{15}{17}{0}
\da{30}{9}{-65}
\end{picture}
\end{array}
\diagspace
\begin{array}{c}
\epsfig{file=threearrow.eps}
\end{array}
\diagspace
\begin{array}{c}
\setlength{\unitlength}{1mm}
\begin{picture}(38,21)
\cell{0}{0}{bl}{\epsfig{file=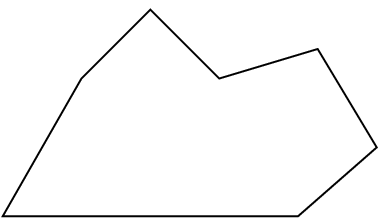}}
\da{19}{7}{0}
\end{picture}
\end{array}
\]
% \hand{50}{43}
\caption{A 3-opetope with 21 sub-opetopes}
\label{fig:3-ope-as-cod}
\end{figure}
seven maps from the unique 0-opetope, nine from the unique 1-opetope, one
from the 2-opetope with two input edges, two from the 2-opetope with three
input edges, one from the 2-opetope with six input edges, and one from
$\omega$ itself (the identity $1_\omega$).

In order to prove results about the relation between opetopic sets and
topological spaces, we will need to know some specific properties of the
category $\scat{O}$.  They are listed here, along with a few more
properties that will not be needed but add detail to the picture.  So,
$\scat{O}$ is a small category such that:
\begin{enumerate}
\item	\lbl{item:O-ax-objects}
the set of objects of $\scat{O}$ is a disjoint union of subsets
  $(O_n)_{n\in\nat}$; we write $\dim(\omega) = n$ if $\omega \in O_n$
\item	\lbl{item:O-ax-inc-dim}
if $\omega' \goby{\rho} \omega$ is a map in $\scat{O}$ then
  $\dim(\omega') \leq \dim(\omega)$ 
\item	\lbl{item:O-ax-ids}
if $\omega' \goby{\rho} \omega$ is a map in $\scat{O}$ with
  $\dim(\omega') = \dim(\omega)$ then $\rho=1_\omega$
\item every map in $\scat{O}$ is monic
\item every map in $\scat{O}$ is a composite of maps of the form $\omega'
  \go \omega$ with $\dim(\omega) = \dim(\omega') + 1$
\item if $\dim(\omega) = \dim(\omega') + 1$ then every map $\omega' \go
  \omega$ can be classified as either a `source%
\index{source!embedding}
embedding' or a `target%
\index{target!embedding}
  embedding' 
\item	\lbl{item:O-ax-finite}
if $\omega \in \scat{O}$ with $\dim(\omega) \geq 1$ then the set of pairs
  $(\omega',\rho)$ with $\dim(\omega) = \dim(\omega') + 1$ and $\rho:
  \omega' \go \omega$ is finite, and there is exactly one such pair for
  which $\rho$ is a target embedding.
\end{enumerate}
Many of these properties can be compared to properties of $\Delta_\mr{inj}$. 
They also imply that every object of $\scat{O}$ is the codomain of only
finitely many maps---in other words, has only a finite set of sub-opetopes.

An \demph{opetopic set}%
\index{opetopic!set}
is by definition a functor $X: \scat{O}^\op \go
\Set$.  If $Y: \Delta^\op \go \Set$ is a simplicial%
\index{simplicial set}
set and
$\upr{n}\in\Delta$ then an element $y$ of $Y\upr{n}$ is usually depicted as
a label attached to an $n$-simplex, and the images of $y$ under the various
face maps as labels on the various faces of the $n$-simplex, as in
\[
\setlength{\unitlength}{1mm}
\begin{picture}(30,22)(-3,-3)
\cell{0}{0}{bl}{\epsfig{file=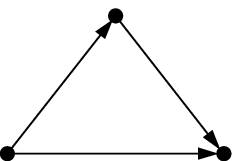}}
\cell{0}{0}{tr}{y_0}
\cell{12}{16}{b}{y_1}
\cell{24}{0}{tl}{y_2}
\cell{6}{9}{br}{y_{01}}
\cell{18}{9}{bl}{y_{12}}
\cell{12}{0}{t}{y_{02}}
\cell{12}{6}{c}{y}
\end{picture}
% \hand{25}{44}
\]
for $n=2$.  Similarly, an opetopic set can be thought of as a system of
labelled opetopes: if $\omega\in O_n$ and $x \in X(\omega)$ then $x$ is
depicted as a label on the $n$-opetope $\omega$, and the images $(X\rho)(x)
\in X(\omega')$ of $x$ under the various face maps $\rho: \omega' \go
\omega$ as labels on the faces of $\omega$.  For example, if $\omega$ is
the 2-opetope with 3 input edges then $x \in X(\omega)$ is drawn as
\[
\setlength{\unitlength}{1mm}
\begin{picture}(38,22)(-3,-3)
\cell{0}{0}{bl}{\epsfig{file=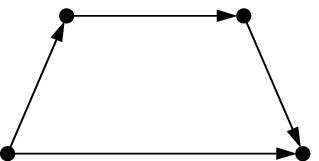}}
\cell{0}{0}{tr}{a_0}
\cell{6}{16}{br}{a_1}
\cell{26}{16}{bl}{a_2}
\cell{32}{0}{tl}{a_3}
\cell{3}{9}{r}{f_1}
\cell{16}{16}{b}{f_2}
\cell{29}{9}{l}{f_3}
\cell{16}{0}{t}{f}
\cell{16}{8}{c}{\Downarrow x}
\end{picture}
% \hand{25}{45}
\]
where $a_0, a_1, a_2, a_3$ are the elements of $X(\gzeros{})$ induced from
$x$ by the four different maps from the unique 0-opetope $\gzeros{}$ to
$\omega$ in $\scat{O}$, and similarly $f_1, f_2, f_3, f$ for the unique
1-opetope.  The elements of $\coprod_{\omega\in O_n} X(\omega)$ are called
the \demph{$n$-cells}%
\index{cell!opetopic set@of opetopic set}%
of $X$, for $n\in\nat$.

Opetopic sets can also be (informally) defined without reference to a
category of opetopes.  Thus, an opetopic set is a commutative diagram of
sets and functions
\begin{equation}	\label{diag:ope-set-comm}
\begin{diagram}[height=1.5em]
\vdots		&\vdots			&\vdots	\\
X'_2		&			&X_2	\\
\dTo<s>{\ \ \ \ \,t}	&
\rdTo(2,4) \ldTo(2,4)	&
\dTo<{s\ \ \ \ \,}>t\\
		&			&	\\
		&			&	\\
X'_1		&			&X_1	\\
\dTo<s>{\ \ \ \ \,t}	&
\rdTo(2,4) \ldTo(2,4)	&
\dTo<{s\ \ \ \ \,}>t\\
		&			&	\\
		&			&	\\
X_0 = X'_0	&			&X_0	\\
\end{diagram}%
\glo{sinopeset}\glo{tinopeset}
\end{equation}
where for each $n\geq 1$, the set $X'_n$ and the functions $s: X'_n \go
X'_{n-1}$ and $t: X'_n \go X_{n-1}$ are defined from the sets $X_n,
X'_{n-1}, X_{n-1}, \ldots, X'_1, X_1, X_0$ and the functions $s,t$ between
them in the way now explained.  As usual for directed graphs or globular
sets, elements of $X_0$ and $X_1$ are called \demph{0-cells} and
\demph{1-cells} respectively, and drawn as labelled points or intervals.
An element of $X'_1$ is called a \demph{$1$-pasting%
\index{pasting diagram!opetopic!labelled}
diagram in $X$}, and
consists of a diagram
\begin{equation}	\label{eq:ope-1-pd}
\gfsts{a_0}\gones{f_1}\gblws{a_1}\gones{f_2} 
\diagspace \cdots \diagspace 
\gones{f_k}
\glsts{a_k}
\end{equation}
of cells in $X$ ($k\in\nat$).  An element $\alpha \in X_2$ is called a
\demph{2-cell}; its source $s \alpha$ is a 1-pasting diagram in $X$, and
its target $t\alpha$ a 1-cell.  If $s\alpha$ is the 1-pasting diagram
of~\bref{eq:ope-1-pd} and $t\alpha$ is
\[
\gfsts{a_0} \gones{g} \glsts{a_k}
\]
then $\alpha$ is drawn as
\begin{equation}	\label{eq:two-ope}
\raisebox{-6.5mm}{%
\setlength{\unitlength}{1mm}
\begin{picture}(36,15)(0,-2)
% Zero-cell marks
\cell{0}{0}{c}{\zmark}
\cell{6}{8}{c}{\zmark}
\cell{36}{0}{c}{\zmark}
% One-cell arrows
\put(0,0){\vector(3,4){6}}
\put(6,8){\vector(3,1){9}}
\put(30,8){\vector(3,-4){6}}
\put(0,0){\vector(1,0){36}}
% Two-cell arrow
\cell{18}{4.5}{c}{\Downarrow}
% \put(18,7){\vector(0,-1){5}}
% Dot-dot-dot
\cell{22}{9.5}{c}{\cdots}
% Labels
\cell{-2.5}{0}{c}{a_0}
\cell{4}{9}{c}{a_1}
\cell{39}{0}{c}{a_k.}
\cell{1}{5}{c}{\scriptstyle f_1}
\cell{10}{11.5}{c}{\scriptstyle f_2}
\cell{35.5}{5}{c}{\scriptstyle f_k}
\cell{18}{-1.5}{c}{\scriptstyle g}
\cell{20.5}{4.5}{c}{\alpha}
\end{picture}}
\end{equation}
An element of $X'_2$ is a \demph{2-pasting diagram in $X$}, that is,
consists of a finite diagram of cells of the form~\bref{eq:two-ope} pasted
together, a typical example being
\begin{equation}	\label{eq:two-pd}
\raisebox{-10.5mm}{%
\setlength{\unitlength}{1mm}
\begin{picture}(50,21)(-5,-3)
% Zero-cell marks
\cell{0}{0}{c}{\zmark}
\cell{-5}{7.5}{c}{\zmark}
\cell{0}{12.5}{c}{\zmark}
\cell{7.5}{10}{c}{\zmark}
\cell{22.5}{10}{c}{\zmark}
\cell{26.25}{15}{c}{\zmark}
\cell{36.25}{12.5}{c}{\zmark}
\cell{43.75}{7.5}{c}{\zmark}
\cell{36.25}{2.5}{c}{\zmark}
\cell{30}{0}{c}{\zmark}
% 1-cell arrows
\put(0,0){\vector(-2,3){5}}
\put(-5,7.5){\vector(1,1){5}}
\put(0,12.5){\vector(3,-1){7.5}}
\put(7.5,10){\vector(1,0){15}}
\put(22.5,10){\vector(3,4){3.75}}
\put(26.25,15){\vector(4,-1){10}}
\put(36.25,12.5){\vector(3,-2){7.5}}
\put(43.75,7.5){\vector(-3,-2){7.5}}
\put(36.25,2.5){\line(-5,-2){6.25}}
\put(30,0){\vector(-3,-1){0}}
\put(36.25,12.5){\vector(0,-1){10}}
\put(0,0){\vector(3,4){7.5}}
\put(22.5,10){\vector(3,-4){7.5}}
\put(0,0){\vector(1,0){30}}
% 2-cell arrows
\cell{14}{5}{c}{\Downarrow}
\da{0}{6}{50}
\da{30}{7.5}{-50}
\da{39}{6.5}{-90}
% Labels
\cell{-1}{-1}{c}{\scriptstyle a_0}
\cell{-7}{7.5}{c}{\scriptstyle a_1}
\cell{0}{14}{c}{\scriptstyle a_2}
\cell{8.5}{11.5}{c}{\scriptstyle a_3}
\cell{21}{11}{c}{\scriptstyle a_4}
\cell{27}{16.5}{c}{\scriptstyle a_5}
\cell{37}{14}{c}{\scriptstyle a_6}
\cell{46}{7.5}{c}{\scriptstyle a_7}
\cell{38}{1.5}{c}{\scriptstyle a_8}
\cell{30}{-1.5}{c}{\scriptstyle a_9}
\cell{-3.5}{3}{c}{\scriptstyle f_1}
\cell{-3.5}{11}{c}{\scriptstyle f_2}
\cell{4}{12.5}{c}{\scriptstyle f_3}
\cell{15}{11.5}{c}{\scriptstyle f_4}
\cell{23.5}{13.5}{c}{\scriptstyle f_5}
\cell{31.5}{15}{c}{\scriptstyle f_6}
\cell{41.5}{11}{c}{\scriptstyle f_7}
\cell{41.5}{4}{c}{\scriptstyle f_8}
\cell{34}{0.25}{c}{\scriptstyle f_9}
\cell{34.5}{6.5}{c}{\scriptstyle f_{10}}
\cell{25}{4}{c}{\scriptstyle f_{11}}
\cell{5}{4}{c}{\scriptstyle f_{12}}
\cell{15}{-1.5}{c}{\scriptstyle f_{13}}
\cell{17}{5}{c}{\scriptstyle \alpha_1}
\cell{1}{8}{c}{\scriptstyle \alpha_2}
\cell{28.5}{9.5}{c}{\scriptstyle \alpha_3}
\cell{39}{8.5}{c}{\scriptstyle \alpha_4}
\end{picture}.}
\end{equation}
Note that the arrows go in compatible directions: for instance, the target or
output edge $f_{11}$ of $\alpha_3$ is a source or input edge of
$\alpha_1$.  The source of this element of $X'_2$ is 
\[
\gfsts{a_0} \gones{f_1}
\diagspace\cdots\diagspace 
\gones{f_9} \glsts{a_9} 
\ \in X'_1, 
\]
and the target is $f_{13} \in X_1$.  An element $\Gamma\in X_3$ is called a
\demph{3-cell}; if, for instance, $s\Gamma$ is the 2-pasting
diagram~\bref{eq:two-pd} then $t\Gamma$ is of the form~\bref{eq:two-ope}
with $k=9$ and $g = f_{13}$, and we picture $\Gamma$ as being a label on
the evident 3-opetope (whose output face is itself labelled $\alpha$ and
whose input faces are labelled $\alpha_1, \alpha_2, \alpha_3, \alpha_4$).
And so it continues.

\index{fundamental!omega-groupoid@$\omega$-groupoid|(}%
\index{topological space!opetopic set from|(}%
\index{opetopic!set!topological space@from topological space|(}
Perhaps the easiest example of an opetopic set is that arising from a
topological space.  We can define this rigorously using only a few of the
properties of $\scat{O}$ listed above.
% on p.~\pageref{item:O-ax-objects}.  
But first we need to `recall' two constructions, one from topology and one
from category theory.

The topology is the cone construction.  This is the canonical way of
embedding a given space into a contractible space, and amounts formally
to a functor $\fcat{Cone}$%
\index{cone}
and a (monic) natural transformation $\iota$ as
shown:
\[
\Top
\ctwo{\id}{\fcat{Cone}}{\iota}
\Top.%
\glo{iotacone}\glo{Cone}
\]
Given a space $E$, the contractible space $\fcat{Cone}(E)$ is the pushout
\[
\begin{diagram}[size=2em]
E	&\rTo^{e\goesto (e,1)}	&E \times [0,1]	\\
\dTo	&			&\dTo		\\
1	&\rTo			&\NWpbk\fcat{Cone}(E)\\
\end{diagram}
\]
in $\Top$, where $1$ is the one-point space.  If $E$ is empty then
$\fcat{Cone}(E) = 1$; otherwise $\fcat{Cone}(E)$ is $E\times [0,1]$ with
all points of the form $(e,1)$ identified.  The inclusion $\iota_E: E
\rIncl \fcat{Cone}(E)$ sends $e$ to $(e,0)$.

The category theory is the construction from any functor
\[
J: \scat{C} \go \Eee
\]
of a pair of adjoint functors
\begin{equation}	\label{eq:induced-adjn}
\begin{diagram}
\Eee	&
\pile{\rTo^{\scriptstyle U}_\top\\ \lTo_{\scriptstyle F}}	&
\ftrcat{\scat{C}^\op}{\Set}.
\end{diagram}
\end{equation}
Here $\scat{C}$ is any small category and $\Eee$ any category with small
colimits.  The right adjoint $U$ is defined by
\[
(UE)(C) = \Eee(JC, E)
\]
($E\in\Eee, C\in\scat{C}$).  The left adjoint $F$ is the left Kan%
\index{Kan, Daniel!extension}
extension
of $J$ along the Yoneda embedding:
\[
\begin{slopeydiag}
\scat{C}&	&\rTo^{\mr{Yoneda}}&	&\ftrcat{\scat{C}^\op}{\Set}	\\
	&\rdTo<J&		&\ldGet>F&				\\
	&	&\Eee.		&	&				\\
\end{slopeydiag}
\]
Explicitly, $F$ is given by the coend formula
\[
FX = \int^{C\in\scat{C}} XC \times JC
\]
($X \in \ftrcat{\scat{C}^\op}{\Set}$).  As the Yoneda embedding is full and
faithful, this is a `genuine' extension: $F(\scat{C}(\dashbk,C)) \iso JC$.
The best-known example---and probably the best remedy for readers new to
Kan extensions and coends---involves%
\index{coend}
simplicial sets.  Here
\[
(\scat{C} \goby{J} \Eee) 
=
(\Delta \goby{J} \Top)
\]
where $J$ sends $[n] \in \Delta$ to the standard $n$-simplex $\Delta^n$.
Then in the adjunction
\[
\begin{diagram}
\Top	&
\pile{\rTo^{\scriptstyle U}_\top\\ \lTo_{\scriptstyle F}}	&
\ftrcat{\Delta^\op}{\Set},
\end{diagram}
\]
$U$ sends a space to its underlying (singular) simplicial set and $F$ is
geometric%
\index{geometric realization!simplicial set@of simplicial set}
realization.  The coend formula above becomes the formula more
familiar to topologists,
\[
FX = ( \coprod_{n\in\nat} X_n \times \Delta^n ) /\sim
% F(X) = \left( \coprod_{n\in\nat} X_n \times \Delta^n \right) /\sim
\]
(Segal~\cite[\S 1]{SegCSS}, Adams~\cite[p.~58]{Ad}) The isomorphism
$F(\Delta(\dashbk,[n])) \iso \Delta^n$ asserts that the realization of the
simplicial $n$-simplex $\Delta(\dashbk,[n])$ is the topological $n$-simplex
$\Delta^n$.

\index{geometric realization!opetopic set@of opetopic set|(}%
\index{opetopic!set!geometric realization of|(}
To define both the underlying opetopic set of a topological space and,
conversely, the geometric realization of an opetopic set, it is therefore
only necessary to define a functor $J: \scat{O} \go \Top$.  We do this
under the assumption that $\scat{O}$ is a category with
properties~\bref{item:O-ax-objects}--\bref{item:O-ax-ids}
(p.~\pageref{item:O-ax-objects}).  The idea is, of course, that $J$ assigns
to each $n$-opetope $\omega$ the topological space $J(\omega)$ looking like
our usual picture of $\omega$ (and homeomorphic to the closed $n$-disk).
In the underlying opetopic set $UE$ of a space $E$, an element of
$(UE)(\omega)$ is a continuous map from $J(\omega)$ into $E$; conversely,
if $X$ is an opetopic set then the space $FX$ is formed by gluing together
copies of the spaces $J(\omega)$ according to the recipe $X$.

For $n\in\nat$, let $\scat{O}(n)$%
\glo{Oatmostn}
be the full subcategory of $\scat{O}$
with object-set $\bigcup_{k\leq n} O_k$.  We will construct for each $n$ a
functor $J_n: \scat{O}(n) \go \Top$ such that the diagram
\[
\begin{diagram}[height=2em]
\scat{O}(0)	&\rIncl		&\scat{O}(1)	&\rIncl	&\scat{O}(2)	&
\rIncl		&\cdots		\\
		&\rdTo(4,2)>{\!\!\!\!J_0}&	&\rdTo>{\!\!J_1}&\dTo>{J_2}&
\ldots		&		\\
		&		&		&	&\Top		&
		&		\\
\end{diagram}
\]
commutes.  Since $\scat{O}$ is the colimit of the top row, this will induce
a functor $J$ of the form desired.

The unique 0-opetope $\blob$ is drawn as a one-point space, so we define
$J_0$ to have constant value $1$.  Let $n\in\nat$ and suppose that we have
defined $J_n$.  By the mechanism described above, $J_n$ induces a geometric
realization functor
\[
F_n: \ftrcat{\scat{O}(n)^\op}{\Set} \go \Top.
\]
The functor $J_{n+1}$ is defined as follows.  Its value on the subcategory
$\scat{O}(n)$ of $\scat{O}(n+1)$ is the same as that of $J_n$.  Any object
$\omega$ of $\scat{O}$ induces a functor
\[
\begin{array}{rrcl}
\scat{O}(\dashbk, \omega)|_{\scat{O}(n)}:&
\scat{O}(n)^\op				&
\go					&
\Set,					\\
					&
\chi					&
\goesto					&
\scat{O}(\chi,\omega),			\\
\end{array}
\]
and if $\omega\in O_{n+1}$ then we put
\[
J_{n+1}(\omega) 
= 
\fcat{Cone}(F_n(\scat{O}(\dashbk, \omega)|_{\scat{O}(n)})).
\]
This defines $J_{n+1}$ on objects.  By assumptions~\bref{item:O-ax-inc-dim}
and~\bref{item:O-ax-ids} on $\scat{O}$, the only remaining maps on which
$J_{n+1}$ needs to be defined are those of the form $\rho: \omega' \go
\omega$ where $\omega'$ is an object of $\scat{O}(n)$ and $\omega\in
O_{n+1}$.  This is done by taking $J_{n+1}(\rho)$ to be the composite
\[
J_n(\omega') \goiso F_n(\scat{O}(\dashbk, \omega')|_{\scat{O}(n)})
\goby{F_n(\rho_*)} F_n(\scat{O}(\dashbk, \omega)|_{\scat{O}(n)})
\rIncl^{\iota} J_{n+1}(\omega)
\]
where the isomorphism is the one noted above and $\rho_*$ is composition
with $\rho$. 

To see why this is sensible, consider the first inductive steps.  The
geometric realization functor
\[
F_0: \ftrcat{\scat{O}(0)^\op}{\Set} \iso \Set \go \Top
\]
realizes sets as discrete spaces.  If $\omega$ is the unique 1-opetope then
\[
\scat{O}(\dashbk,\omega)|_{\scat{O}(0)}:
\scat{O}(0)^\op
\go
\Set
\]
has value $2$, since there are two maps from the unique 0-opetope $\omega'$
into $\omega$; hence
\[
J_1(\omega) = \fcat{Cone}(\scat{O}(\dashbk,\omega)|_{\scat{O}(0)})
\]
is the cone on the discrete 2-point space, which is homeomorphic to the
unit interval $[0,1]$.  The two maps $\omega' \parpairu \omega$ are sent by
$J_1$ to the two maps $1 \parpairu [0,1]$ picking out the endpoints.  For
the next inductive step, $n=1$, note that 
\[
F_1: \ftrcat{\scat{O}(1)^\op}{\Set} \go \Top
\]
is the usual functor geometrically%
\index{graph!directed!geometric realization of}
realizing directed graphs.  If $\omega$
is the 2-opetope with $r$ input edges then
$\scat{O}(\dashbk,\omega)|_{\scat{O}(1)}$ is the directed graph
\[
\topeq{}{}{}{}{}
\]
with $r+1$ vertices and $r+1$ edges, whose geometric realization is the
circle $S^1$; hence $J_1(\omega)$ is $\fcat{Cone}(S^1)$, homeomorphic to
the closed disk $D^2$.%
\index{geometric realization!opetopic set@of opetopic set|)}%
\index{opetopic!set!geometric realization of|)}
\index{fundamental!omega-groupoid@$\omega$-groupoid|)}%
\index{topological space!opetopic set from|)}%
\index{opetopic!set!topological space@from topological space|)}

It should certainly be true in general that if $\omega$ is an $n$-opetope
then $J(\omega)$ is homeomorphic to the $n$-disk $D^n$, but we do not have
enough information about $\scat{O}$ to prove that here.

A similar construction can be tried with strict $\omega$-categories in
place of spaces.  The idea now is that every $n$-opetope gives rise to a
strict $n$-category---hence%
\lbl{p:ope-to-n-cat}%
\index{omega-category@$\omega$-category!opetopic set from}%
\index{opetopic!set!strict omega-category@from strict $\omega$-category}
to a strict $\omega$-category in which the only cells of dimension greater
than $n$ are identities---and every map between opetopes gives rise to a
strict functor between the resulting strict $\omega$-categories.  For
example, the strict 2-category associated to the 2-opetope
\[
\topec{}{}{}{}{\Downarrow}
\]
is freely generated by this diagram, in other words, by 0-cells $a_0$,
$a_1$, $a_2$, $a_3$, 1-cells $f_i: a_{i-1} \go a_i$ ($i = 1, 2, 3$) and
$f: a_0 \go a_3$, and a 2-cell $f_3\of f_2\of f_1 \go f$.  However, the
formal construction is more difficult than that for topological spaces
because there is a question of orientation.  I do not know whether
properties~\bref{item:O-ax-objects}--\bref{item:O-ax-finite} of $\scat{O}$
suffice to do the construction, but we would certainly need to know more
about $\scat{O}$ in order to say anything really useful.

\section{Weak $n$-categories: a sketch}
\lbl{sec:ope-n}

The opetopic sets arising from topological spaces and from strict
$\omega$-categories should all be weak $\omega$-categories.  `Should' means
that if someone proposes to define%
\index{n-category@$n$-category!definitions of}
a weak $\omega$-category as an opetopic
set with certain properties then these two families of examples must surely
be included---or if not, their concept of weak $\omega$-category is very
different from mine.  So let us now look at the properties we might ask of
an opetopic set in order for it to qualify as a weak $\omega$-category.

The situation is roughly like that in~\ref{sec:non-alg-notions}, where we
characterized the (plain) multicategories that arise%
\index{monoidal category!multicategory@\vs.\ multicategory}
from monoidal
categories and so were able to re-define a monoidal category as a
multicategory%
\index{multicategory!representable}
with properties.  The main observation was that if $a_1,
\ldots, a_k$ are objects of a monoidal category $A$ then in the underlying
multicategory $C$, the canonical map
\[
\setlength{\unitlength}{1em}
\begin{picture}(4,8)(-2,0)
\cell{0}{6}{b}{\tinputsvert{a_1}{a_k}}
\cell{0}{6}{t}{\tusualvert{}}
\cell{0}{2}{t}{\toutputvert{(a_1 \otimes\cdots\otimes a_k)}}
\end{picture}
\]
is `universal'%
\index{universal!map in multicategory}
as a map with domain $a_1, \ldots, a_k$, and this determines
the tensor $(a_1 \otimes\cdots\otimes a_k)$ up to isomorphism.  It follows
that a monoidal category can be re-defined as a multicategory in which
every sequence of objects is the domain of some `universal' map.  Actually,
there is a choice of notions of `universality':
in~\ref{sec:non-alg-notions} we spoke of both universal maps and the more
general pre-universal%
\index{pre-universal!map in multicategory}
maps, and if we use pre-universals then we must add
to this re-definition the requirement that the composite of pre-universal
maps is pre-universal.

Similarly, a weak $\omega$-category can be thought of as an opetopic set in
which there are `enough universals'.  Consider, for example, the opetopic
set $X$ arising from a \emph{strict} $\omega$-category $A$.  The 0- and
1-cells of $X$ are the same as the 0- and 1-cells of $A$, and a 2-cell in
$X$ of the form~\bref{eq:two-ope} (p.~\pageref{eq:two-ope}) is a 2-cell
\[
a_0 \ctwomult{f_k \of\cdots\of f_1}{g}{} a_k
\]
in $A$.  Now, given such a string of 1-cells $f_1, \ldots, f_k$, there is a
distinguished 2-cell $\epsln$ in $X$ of this form: the one corresponding to
the identity 2-cell
\[
a_0 \ctwomult{f_k \of\cdots\of f_1}{f_k \of\cdots\of f_1}{1} a_k
\]
in $A$.  We would like to pin down some universal property of the 2-cells
of $X$ arising in this way.  The obvious approach is to say something like
`every 2-cell of $X$ of the form~\bref{eq:two-ope} factors uniquely through
$\epsln$', but to say `factors' we need to know about composition of
\emph{2-cells}, even though at present we are only trying to discuss
composition of \emph{1-cells}\ldots

This problem of `downwards induction' means that it is easier to define
weak $n$-category, for finite $n$, than weak $\omega$-category.

Here, then, is a sketch of a definition of weak $n$-category, roughly that
proposed by Baez%
\index{Baez, John}
and Dolan%
\index{Dolan, James}
in~\cite{BDHDA3}.  (A summary can also be found
as Definition \textbf{X} in my~\cite{SDN}.)  Let $n\in\nat$.  In a moment I
will say what it means for a cell of an opetopic set to be `pre-universal'
(or `universal'%
\index{universal!cell of n-category@cell of $n$-category}
in the terminology of the sources just cited), a condition
depending on $n$.  A \demph{weak $n$-category}%
\index{n-category@$n$-category!definitions of!opetopic}
is an opetopic set $X$ such
that
\begin{itemize}
\item every pasting diagram is the source of a pre-universal cell
\item the composite of pre-universal cells is pre-universal.
\end{itemize}
The first condition means that if $n\geq 1$ and $\Phi$ is an $n$-pasting
diagram in $X$ then there exists a pre-universal $(n+1)$-cell $\epsln$ with
source $\Phi$:
\[
\Phi \goby{\epsln} \phi.
\]
(In the notation of diagram~\bref{diag:ope-set-comm},
p.~\pageref{diag:ope-set-comm}, we have $\epsln\in X_{n+1}$, $\Phi =
s\epsln \in X'_n$, and $\phi = t\epsln \in X_n$.)  Think of $\phi$ as
a---or, with a pinch of salt, `the'---composite of the cells making up the
pasting diagram $\Phi$, and $\epsln$ as asserting that $\phi$ is a
composite of $\Phi$.  This suggests the correct meaning of the second
condition: that if
\[
\Psi \goby{\epsln} \psi
\]
is such that the $(n+1)$-cell $\epsln$ and all of the $n$-cells making up
the $n$-pasting diagram $\Psi$ are pre-universal, then the $n$-cell $\psi$
is also pre-universal.

What should it mean for a cell to be pre-universal?%
\index{pre-universal!cell of n-category@cell of $n$-category}
 For a start, we define
pre-universality in such a way that all cells of dimension greater than $n$
in a weak $n$-category are trivial.  This means that if $k\geq n$ then
every $k$-pasting diagram $\Phi$ is the source of precisely one
$(k+1)$-cell, whose target is to be thought of as the composite of $\Phi$.
Taking $k=n$, we obtain a composition of $n$-cells obeying strict laws.

Now, an $n$-cell
\[
\Phi \goby{\epsln} \phi
\]
is pre-universal if and only if every $n$-cell
\[
\Phi \goby{\alpha} \psi
\]
with source $\Phi$ factors uniquely through $\epsln$---in other words,
there is a unique $n$-cell $\ovln{\alpha}$ such that $\alpha$ is the
composite
\[
\Phi \goby{\epsln} \phi \goby{\ovln{\alpha}} \psi.
\]
Equivalently, $\epsln$ is pre-universal if for every $(n-1)$-cell $\psi$
parallel to $\phi$, composition with $\epsln$ induces a bijection
\begin{equation}	\label{eq:hom-sets-bjn}
X(\phi,\psi) \goiso X(\Phi,\psi)
\end{equation}
of hom-sets.  We think of $\phi$ as a composite of $\Phi$, or as `the'
composite if it is understood that it is only defined up to isomorphism.

At the next level down, we want to say that an $(n-1)$-cell $\Phi
\goby{\epsln} \phi$ is pre-universal if for every $(n-2)$-cell $\psi$
parallel to $\phi$, composition with $\epsln$ induces an
equivalence~\bref{eq:hom-sets-bjn} of hom-categories.  The difficulties are
that we must first put a category structure on the domain and codomain,
and, more significantly, that `composition' with $\epsln$ is not actually a
functor, since composition of $(n-1)$-cells is only defined up to
isomorphism.  At this point we really need to set up some more language,
and this would carry us beyond the scope of this informal account, so I
refer the curious reader to the texts cited in the Notes.

The case of weak $2$-categories with only one $0$-cell is the familiar one
of monoidal categories as representable%
\index{multicategory!representable}
multicategories~(\ref{sec:non-alg-notions}).  An opetopic set $X$ with only
one 0-cell and trivial above dimension 2 is a set $C_0$ (the 1-cells of
$X$) together with a set $C(a_1, \ldots, a_n; a) $ for each sequence
$a_1, \ldots, a_n, a$ of elements of $C_0$.  If $X$ is a weak 2-category
then there is, as explained above, a composition for 2-cells, obeying
strict associativity and unit laws; this gives a multicategory $C$ with the
indicated hom-sets.  A 2-cell of $X$ is pre-universal if and only if the
corresponding arrow of $C$ is pre-universal; the axiom that every 1-pasting
diagram of $X$ is the source of a pre-universal 2-cell is the axiom that
every sequence of objects of $C$ is the source of a pre-universal arrow;
the only other axiom is that the composite of pre-universals is
pre-universal.  So a one-object weak 2-category is nothing other than a
representable multicategory, which in turn is the same thing as a monoidal
category (Theorem~\ref{thm:rep-multi}).  See Cheng~\cite{CheOBC} for
details. 

Although the approach described above is modelled on that of Baez%
\index{Baez, John}
and
Dolan~\cite{BDHDA3},%
\index{Dolan, James}
if you look at their paper then you will see
immediately that it has features quite different from any that we have
considered.  Here is a description of what they actually do.

The most striking feature is that they work throughout with
\emph{symmetric} multicategories, and that these play essentially the same
role as our generalized multicategories.  (In fact they call symmetric
multicategories `typed%
\index{operad!typed}
operads', and symmetric operads `(untyped) operads';
this is like the `coloured operad' terminology discussed on
p.~\pageref{p:col-opd}.)  The most important thing that they do with
symmetric multicategories is this: given a symmetric multicategory $C$,
they construct a new symmetric multicategory $C^+$,%
\glo{BDslicemti}
the \demph{slice}%
\index{slice!symmetric multicategory@of symmetric multicategory}
of
$C$, such that
\begin{equation}	\label{eq:sym-slice}
\Alg(C) \iso 
(\textrm{symmetric multicategories over }
C
\textrm{ with the same object-set}).
\end{equation}
The left-hand side is the category of algebras for $C$ as a
\emph{symmetric} multicategory (p.~\pageref{p:defn-sym-alg}).  An object of
the category on the right-hand side is a map $D \goby{f} C$ where $D$ is a
symmetric multicategory with $D_0 = C_0$ and $f$ is a map of symmetric
multicategories with $f_0 = 1_{C_0}$; a map is a commutative triangle.

To describe $C^+$ explicitly it is helpful to use the following
terminology: a \demph{reduction%
\index{reduction law}
law} in a multicategory $D$ is an equation
stating what the composite of some family of arrows in $D$ is.  So
reduction laws in $D$ correspond to trees of arrows in $D$, but we think of
a reduction law as also knowing what the composite of this tree is; thus,
the collection of reduction laws in $D$ describes the composition in $D$
completely (indeed, very redundantly).  Here `tree' must be understood in
the symmetric sense: branches are allowed to cross over, but the
topologically-obvious identifications are made so that any tree is
equivalent to a `combed%
\index{tree!combed}
tree'---one where all the crossings are at the top.
Now:
\begin{itemize}
\item an object of $C^+$ is an arrow of $C$
\item an arrow of $C^+$ is a reduction law in $C$
\item a reduction law in $C^+$ is a way of assembling a family of reduction
  laws in $C$ to form a new reduction law.  
\end{itemize}
So if $\theta_1, \ldots, \theta_n, \theta$ are arrows of $C$ then an arrow
$\theta_1, \ldots, \theta_n \go \theta$ in $C^+$ is a tree%
\index{tree!vertices ordered@with vertices ordered}
with $n$
vertices, totally ordered, such that if $\theta_i$ is written at the $i$th
vertex then the evident composite can be formed in $C$ and is equal to
$\theta$.  

\begin{example}
  Consider the terminal%
\index{multicategory!terminal}
symmetric multicategory $1$, whose algebras are the
  commutative monoids~(\ref{eg:opd-terminal}).  Then according to the
  explicit description, the objects of $1^+$ are the natural numbers; an
  arrow $k_1, \ldots, k_n \go k$ in $1^+$ is a $k$-leafed tree with $n$
  vertices, totally ordered, such that the $i$th vertex has $k_i$ branches
  coming up out of it; composition is by substituting trees into vertices
  of other trees.  So the arrows $k_1, \ldots, k_r \go k$ describe the ways
  of combining one $k_1$-ary operation, one $k_2$-ary operation, \ldots,
  and one $k_n$-ary operation of a generic symmetric multicategory in order
  to obtain a $k$-ary operation.  Hence $1^+$ is the symmetric
  multicategory%
\index{multicategory!symmetric!operads@for operads}
whose algebras are symmetric operads; we met it as
  $\cat{O}'$ in~\ref{eg:sym-multi-for-opds}.
\end{example}

\begin{example}
  Let $I$%
\glo{initsymopd}
be the initial%
\index{operad!initial}
symmetric operad, which has only one object and one arrow
(the identity).  The explicit description says that $I^+$ also has only one
object, in other words, is an operad.  The reduction laws of $I$ say that
if you take $n$ copies of the identity arrow, permute them somehow, then
compose them, then you obtain the identity arrow.  So an $n$-ary operation
of the operad $I^+$ is an element of the symmetric group $S_n$, and $I^+$
is the operad%
\index{operad!symmetries@of symmetries}
$\SymOpd$ of symmetries~(\ref{eg:opd-Sym}).
\end{example}

Baez and Dolan define an \demph{$n$-dimensional opetope}%
\index{opetope}
to be an object of
$I^{+\cdots +}$, where there are $n$ $+$'s above the $I$.  For $n=0$ and
$n=1$ this looks fine: there is only one opetope in each of these
dimensions.  But a 2-dimensional opetope is an object of $I^{++}$, that is,
an arrow of $I^+$, and the set of such is the disjoint union of all the
symmetric groups---not%
\lbl{p:opetope-discrepancy}
the expected answer, $\nat$.  In their world, there is not just
one 2-dimensional opetope looking like
\[
\topec{}{}{}{}{\Downarrow},
\]
but 6,%
\index{opetope!proliferation of}
corresponding to the $3!$ permutations of the 3 input edges.
Similarly, as observed by Cheng~\cite[1.3]{CheWOM},%
\index{Cheng, Eugenia}
there are $311040$
copies of the 3-dimensional opetope shown in Fig.~\ref{fig:3-ope-as-cod}
(p.~\pageref{fig:3-ope-as-cod}).  (Verifying this is a useful exercise for
understanding the slice construction.)  We come back to this discrepancy
later.

Their next step is to define opetopic set.  We omit this; it
requires a little more multicategory theory.  Then, as sketched above, they
define what it means for a cell of an opetopic set to be universal,
hence what it means for an opetopic set to be a weak $n$-category.

Actually, they do more.  Their definition of opetope uses $I$ as its
starting point, but it is possible to start with any other symmetric
multicategory $C$ instead, to obtain a definition of `$n$-dimensional
$C$-opetope'.%
\index{opetope!generalized}
 There is a corresponding notion of `$C$-opetopic set',
and of what it means for a $C$-opetopic set to be an `$n$-coherent%
\index{coherent!algebra}%
\index{n-coherent algebra@$n$-coherent algebra}
$C$-algebra' (in the case $C=I$, a weak $n$-category).  For example, let
$C$ be the terminal operad $1$.  Then an $n$-coherent $1$-algebra is what
they call a `stable%
\index{n-category@$n$-category!stable}%
\index{stabilization}
$n$-category', or might also be called a
symmetric monoidal $n$-category, as in the periodic table
(p.~\pageref{p:periodic-table}).  A $0$-coherent $C$-algebra is just a
$C$-algebra, so a stable $0$-category is a commutative monoid, as it should
be.

\index{multicategory!symmetric vs. generalized@symmetric \vs.\ generalized|(}
We have discussed elsewhere (p.~\pageref{p:sym-discussion}) the difference
between symmetric
and generalized multicategories.  In short, the
generalized multicategory approach allows the geometry to be represented
faithfully, whereas the symmetric approach destroys the geometry and
squashes everything into one dimension; but symmetric multicategories are
perhaps easier to get one's hands on.  Cheng%
\index{Cheng, Eugenia}
offers the following
analogy~\cite[1.4]{CheWOM}: she does not like tidying her desk.%
\index{desk-tidying}
 As she
works away and produces more and more pages of notes, each page gets put
into its natural position on the desktop, with notes on related topics
side-by-side, the sheet currently being written on at the front, and so on.
If she is forced to clear up her desk then she must destroy this natural
configuration and stack the papers in some arbitrary order; but having put
them into a single stack, they are much easier to carry around.

Despite symmetric multicategories not falling readily into our scheme of
generalized multicategories, there are many points of similarity between
Baez and Dolan's constructions and ours.

Slicing is an example.  Their construction of the slice%
\index{slice!symmetric multicategory@of symmetric multicategory}
$C^+$ of a
symmetric multicategory $C$ proceeds in two stages.  First they show how to
slice%
\index{slice!symmetric multicategory by algebra@of symmetric multicategory by algebra}
a multicategory by one of its algebras: given a symmetric
multicategory $D$ and a $D$-algebra $X$, they construct a symmetric
multicategory $D/X$ such that
\begin{equation}	\label{eq:sym-slice-alg}
\Alg(D/X) \eqv \Alg(D)/X.
\end{equation}
Then they show how to construct, for any set $S$, a symmetric multicategory
$\fcat{Mti}_S$ such that
\[
\Alg(\fcat{Mti}_S) \eqv \fcat{SymMulticat}_S
\]
where the right-hand side is the category whose objects are the symmetric
multicategories with object-set $S$ and whose maps are those leaving the
object-set fixed.  (They write $D/X$ as $X^+$ and leave $\fcat{Mti}_S$
nameless.  We wrote $\fcat{Mti}_S$ as $\cat{O}'_S$
in~\ref{eg:sym-multi-for-opds}.)  The slice multicategory $C^+$ of $C$ is
then defined by
\begin{equation}	\label{eq:slice-defn}
C^+ = \fcat{Mti}_{C_0}/C,
\end{equation}
so that
\[
\Alg(C^+) \eqv 
\Alg(\fcat{Mti}_{C_0})/C \eqv 
\fcat{SymMulticat}_{C_0}/C,
\]
as required.  These constructions are mirrored in our world of generalized
multicategories.  First, we saw in Proposition~\ref{propn:slice-multicat}
how to slice%
\index{slice!generalized multicategory by algebra@of generalized multicategory by algebra}
a multicategory by an algebra: if $T$ is a cartesian monad on
a cartesian category $\Eee$, $D$ is a $T$-multicategory, and $X$ is a
$D$-algebra, then there is a $T$-multicategory $D/X$
satisfying~\bref{eq:sym-slice-alg}.  Second, assuming that $\Eee$ and $T$
are suitable, Theorem~\ref{thm:free-fixed} tells us that for any $S\in\Eee$
there is a cartesian monad $T^+_S$ on a cartesian category $\Eee^+_S$ such
that
\[
\Alg(T^+_S) \eqv (\Eee,T)\hyph\Multicat_S,
\]
and so if we take $\fcat{Mti}_S$%
\index{generalized multicategory!generalized multicategories@for generalized multicategories}
to be the terminal $T^+_S$-multicategory
then
\[
\Alg(\fcat{Mti}_S) \eqv (\Eee,T)\hyph\Multicat_S
\]
by~\ref{eg:alg-terminal}.  So we can define the $T^+_S$-multicategory $C^+$
by the same formula~(\ref{eq:slice-defn}) as in the symmetric case, to
reach the analogous conclusion
\[
\Alg(C^+) \eqv 
\Alg(\fcat{Mti}_{C_0})/C \eqv 
(\Eee,T)\hyph\Multicat_{C_0}/C.
\]

This neatly illustrates the contrast between the two approaches.  With
generalized multicategories, slicing raises the dimension: if $C$ is a
$T$-multicategory then $C^+$ is a $T^+_{C_0}$-multicategory.  When we draw
the pictures this dimension-shift seems perfectly natural.  On the other
hand, if $C$ is a symmetric multicategory then $C^+$ is again a symmetric
multicategory, and this is certainly convenient.  Note also that there is a
flexibility in the generalized approach lacking in the symmetric approach,
revealed if we attempt to drop the restriction fixing the object-set.  That
is, suppose that given a ($T$- or symmetric) multicategory $C$ we wish to
define a new multicategory%
\index{multicategory!symmetric!multicategories@for multicategories}%
\index{generalized multicategory!generalized multicategories@for generalized multicategories}
$C^+$ whose algebras are \emph{all}
multicategories over $C$.  This is easy for generalized multicategories: we
just replace $\Eee^+_S$ and $T^+_S$ by $\Eee^+$ and $T^+$ in the
construction above, using Theorem~\ref{thm:free-main}, and $C^+$ is a
$T^+$-multicategory.  But it is impossible in the symmetric theory, as
there is no symmetric multicategory whose algebras are all symmetric
multicategories.

Another point where the two approaches proceed in analogous ways is in the
construction of opetopes.  As discussed, Baez and Dolan considered
so-called $C$-opetopes for any symmetric multicategory $C$, not just the
case $C=I$ that generates the ordinary opetopes.  The analogous point for
generalized multicategories is that the opetopic categories $\Eee_n$ and
monads $T_n$ were generated from $\Eee_0 =\Set$ and $T_0=\id$
(Definition~\ref{defn:opetopic-monads}), but we could equally well have
started from any other suitable category $\Eee$ and suitable monad $T$ on
$\Eee$.  This would give a new sequence $(O_{T,n})_{n\in\nat}$ of objects
of $\Eee$---the objects of `$n$-dimensional $T$-opetopes'.%
\index{opetope!generalized}
\index{multicategory!symmetric vs. generalized@symmetric \vs.\ generalized|)}

We now come to the discrepancy between Baez and Dolan's opetopes and ours,
mentioned on p.~\pageref{p:opetope-discrepancy}.  A satisfactory resolution
has been found by Cheng%
\lbl{p:BD-tweak}%
\index{Cheng, Eugenia}
\cite{CheWOM, CheWCO}, and can be explained roughly as follows.  The
problem in Baez and Dolan's approach is that the slicing process results in
a loss of information.  For example, as we observed, the symmetric
multicategory $I^{++}$ has 6%
\index{opetope!proliferation of}
objects drawn as
\[
\topec{}{}{}{}{\Downarrow},
\]
but formally it contains nothing to say that they are in any way related.
That information about symmetries has been discarded in the slicing
process.  With each new slicing another layer of information about
symmetries is lost, so that even a simple 3-dimensional opetope is
reproduced into hundreds of thousands of apparently unrelated copies.

The remedy is to work with structures slightly more sophisticated than
symmetric multicategories, capable of holding just a little more
information.  For the short while that we discuss them, let us call them
\demph{enhanced symmetric multicategories};%
\lbl{p:enhanced}%
\index{multicategory!enhanced}%
\index{multicategory!symmetric!enhanced}
the definition can be found in Cheng~\cite[2.1]{CheWOM}%
\index{Cheng, Eugenia}
or~\cite[1.1]{CheWCO}.  An enhanced symmetric multicategory is like a
symmetric multicategory, but as well as having the usual objects and
arrows, it also has \demph{morphisms} between objects, making objects and
morphisms into a category.  The morphisms are \emph{not} the same as the
arrows, although naturally there are axioms relating them; note in
particular that whereas an arrow has a finite sequence of objects as its
domain, a morphism has only one.  An enhanced symmetric multicategory in
which the category of objects and morphisms is discrete is an ordinary
symmetric multicategory.  (This explains Cheng's terminology: for her,
enhanced symmetric multicategories are merely `symmetric multicategories',
and symmetric multicategories are `object-discrete symmetric
multicategories'.  Compare also~\ref{eg:mti-sym} above.)

The virtue of enhanced symmetric multicategories is that like symmetric
multicategories, they can be sliced, but unlike symmetric multicategories,
slicing retains all the information about symmetries.  So in the enhanced
version of $I^{++}$, all 6 objects shaped like the 2-dimensional opetope
above are isomorphic---in other words, there is essentially only one of
them.  In fact, Cheng%
\index{Cheng, Eugenia}
has proved~\cite{CheWCO} that there is a natural
one-to-one correspondence between her modified version of Baez and Dolan's
opetopes and the opetopes defined in the present text.  They also
correspond to the multitopes%
\index{multitope}
of Hermida,%
\index{Hermida, Claudio}
Makkai%
\index{Makkai, Michael}
and Power:%
\index{Power, John}
see the Notes at
the end of this chapter.

\section{Many in, many out}
\lbl{sec:many}%
\index{many in, many out|(}

Multicategories and, more specifically, opetopes, are based on the idea of
operations taking many inputs and one output.  We could, however, decide to
allow multiple outputs too, so that 2-opetopes are replaced by shapes like
\[
\begin{array}{c}
\setlength{\unitlength}{1mm}
\begin{picture}(52,24)
\cell{0}{0}{bl}{\epsfig{file=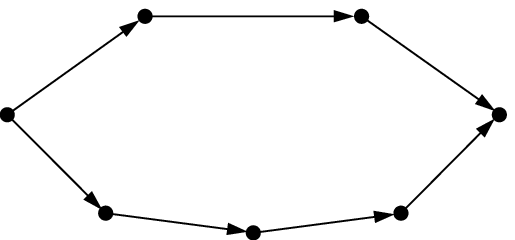}}
\cell{26}{12}{c}{\Downarrow}
\end{picture}
\end{array},
% \hand{25}{46},
\]
3-opetopes by shapes looking something like soccer balls, and so on.

The idea of doing higher category theory with cells shaped like this is
attractive, since they encompass all the other shapes in common use:
globular, simplicial, cubical and opetopic.  In other words, these shapes
are intended to represent all possible composable diagrams of cells in an
$n$- or $\omega$-category, and so we are aiming for a universal, canonical,
totally unbiased approach, and that seems very healthy.

There is, however, a fundamental obstruction.  It turns out that it does
not make sense to talk about `these shapes', at least in dimensions
higher than 2.  They are simply not well-defined.  

To make these statements precise we define an analogue of opetopic set
with multiple outputs as well as multiple inputs, then show that these
structures, unlike opetopic sets, do not form a presheaf category.  So
there is no category of many-in, many-out shapes analogous to the category
$\scat{O}$ of opetopes.  Actually, it is not quite opetopic sets of which
we define an analogue, but rather $n$-dimensional versions in which all
cells have dimension at most $n$.  We need only go up to dimension 3%
\index{three-dimensional structures@3-dimensional structures}
to
encounter the obstruction.  (It is no coincidence that the obstruction to
simple-minded coherence for $n$-categories also appears in dimension 3.)
These $n$-dimensional analogues of opetopic sets are the $n$-computads%
\index{computad|(}
of
Street~\cite{StrLIC, StrCS};%
\index{Street, Ross}
let us define them for $n\leq 3$.

The category $0\hyph\fcat{Cptd}$ of \demph{0-computads} is $\Set$.

A \demph{1-computad} is a directed graph; if $\scat{C}(1)$ is the category
$(0 \parpairu 1)$ then
\[
1\hyph\fcat{Cptd} = \ftrcat{\scat{C}(1)^\op}{\Set}.
\]
We have the familiar adjunction
\[
\begin{diagram}[width=4em]
\Cat	&
\pile{\rTo^{\scriptstyle U_1}_\top\\ \lTo_{\scriptstyle F_1}}	&
1\hyph\fcat{Cptd}
\end{diagram}
\]
(induced, if you like, by the functor $\scat{C}(1) \go \Cat$ sending 0 to
the terminal category and 1 to the category consisting of a single arrow,
using the mechanism described on p.~\pageref{eq:induced-adjn}).  There is
also a functor $P_1: 1\hyph\fcat{Cptd} \go \Set$ sending a directed graph
to the set of parallel pairs of edges in it.

A \demph{2-computad} is a 1-computad $X$ together with a set $X_2$ and a
function $\xi: X_2 \go P_1 U_1 F_1 (X)$.  A map $(X, X_2, \xi) \go (X',
X'_2, \xi')$ is a pair $(X \go X', X_2 \go X'_2)$ of maps making the
obvious square commute.  Concretely, a 2-computad consists of a collection
$X_0$ of 0-cells, a collection $X_1$ of 1-cells between 0-cells, and a
collection of $X_2$ of 2-cells $\alpha$ of the form
\[
\begin{diagram}[height=1.5em]
% 	&	&	&	&	&	&	&	&	\\
	&	&	&	&\,	&\ldots	&	&	&	\\
	&	&a_1	&\ruTo(2,1)<{f_2}&&	&\,	&	&	\\
	&\ruTo<{f_1}&	&	&	&	&	&\rdTo>{f_n}&	\\
a_0=b_0	&	&	&	&\Downarrow\alpha&&	&	&a_n=b_m\\
	&\rdTo<{g_1}&	&	&	&	&	&\ruTo>{g_m}&	\\
	&	&b_1	&	&	&	&\,	&	&	\\
	&	&	&\rdTo(2,1)<{g_2}&\,&\cdots&	&	&	\\
\end{diagram}
% \hand{35}{47}
\]
where $n, m \in \nat$, $a_i, b_j \in X_0$, and $f_i, g_j \in X_1$.
(The possibility exists that $n = m = 0$: this will prove critical.)  Hence 
\[
2\hyph\fcat{Cptd} \eqv \ftrcat{\scat{C}(2)^\op}{\Set}
\]
where $\scat{C}(2)$ is the evident category with one object for each pair
$(n,m)$ of natural numbers and then two further objects (one for each of
dimensions $1$ and $0$).  There is an adjunction
\[
\begin{diagram}[width=4em]
\strcat{2}	&
\pile{\rTo^{\scriptstyle U_2}_\top\\ \lTo_{\scriptstyle F_2}}	&
2\hyph\fcat{Cptd}
\end{diagram}
\]
between strict 2-categories and 2-computads (induced, if you like, by the
evident functor $\scat{C}(2) \go \strcat{2}$).  There is also a functor
$P_2: 2\hyph\fcat{Cptd} \go \Set$ sending a 2-computad to the set of
parallel pairs of 2-cells in it.

Similarly, a \demph{3-computad} is a 2-computad $X$ 
together with a set $X_3$ and a function $X_3 \go P_2 U_2 F_2 (X)$, and
maps are defined in the obvious way.  But there is no category
$\scat{C}(3)$ of `computopes of dimension $\leq 3$':

\begin{propn}	\lbl{propn:3-computads}
The category of 3-computads is not a presheaf category.
\end{propn}
The following proof is due to Makkai%
\index{Makkai, Michael}
and Zawadowski%
\index{Zawadowski, Marek}
(private communication)
and to Carboni%
\index{Carboni, Aurelio}
and Johnstone~\cite{CJcorr}.%
\index{Johnstone, Peter}
\begin{prooflike}{Proof of~\ref{propn:3-computads}}
The category $3\hyph\fcat{Cptd}$ is a certain comma category, the Artin%
\index{Artin gluing}%
\index{gluing}
gluing (p.~\pageref{p:Artin}) of the functor
\[
P_2 U_2 F_2: 2\hyph\fcat{Cptd} \go \Set.
\]
Proposition~\ref{propn:Artin-gluing} tells us that $3\hyph\fcat{Cptd}$ is a
presheaf category only if $P_2 U_2 F_2$ preserves wide pullbacks.
We show that it does not even preserve ordinary pullbacks (which by Carboni
and Johnstone~\cite[4.4(ii)]{CJ} implies that $3\hyph\fcat{Cptd}$, far from
being a presheaf category, is not even locally cartesian closed).  

Let $\cat{S}$ be the full subcategory of $2\hyph\fcat{Cptd}$ consisting of
those 2-computads with only one 0-cell and no 1-cells.  Plainly $\cat{S}
\eqv \Set$.  The inclusion $\cat{S} \rIncl 2\hyph\fcat{Cptd}$ preserves
pullbacks, so it suffices to show that the functor $\cat{S} \go \Set$
obtained by restricting $P_2 U_2 F_2$ to $\cat{S}$ does not preserve
pullbacks.  

If $X\in\cat{S}$ and $C \in \strcat{2}$ then a map $X \go U_2(C)$ consists
of an object $c$ of $C$ together with a 2-cell
\[
c \ctwomult{1_c}{1_c}{} c
\]
of $C$ for each 2-cell of $X$.  Example~\ref{eg:str-mon-comm}%
\index{Eckmann--Hilton argument}
tells us that a strict 2-category with only one 0-cell and one 1-cell is
just a commutative monoid, so identifying $\cat{S}$ with $\Set$,
\[
U_2 F_2: \Set \go \Set
\]
is the free commutative monoid functor $M$.  Since any two 2-cells of a
2-computad in $\cat{S}$ are parallel, we now just have to show that the
endofunctor 
\[
A \goesto M(A) \times M(A)
\]
of $\Set$ does not preserve pullbacks.  This follows from the fact that $M$
itself does not preserve pullbacks, which we proved in
Example~\ref{eg:comm-not-cart}.  
\done
\end{prooflike}%
\index{many in, many out|)}%
\index{computad|)}

\begin{notes}

Opetopes were defined by Baez and Dolan~\cite{BDHDA3}.  Hermida, Makkai and
Power defined their own version~\cite{HMP1, HMP2, HMP3, HMPHDM}, calling
them multitopes, as did I~\cite[4.1]{GOM}, re-using the name `opetope'.
After tweaking Baez and Dolan's definition (see p.~\pageref{p:BD-tweak}),
Cheng%
\index{Cheng, Eugenia}
proved all three notions equivalent \cite{CheWOM, CheWCO}.  On
ordering%
\index{order!non-canonical}%
\index{pasting diagram!opetopic!ordering of}
the opetopes making up a pasting diagram, she writes:
\begin{quote}
  The three different theories arise, essentially, from three different
  ways of tackling this issue.  Baez and Dolan propose listing them in
  every order possible, giving one description \emph{for each} ordering.
  Hermida, Makkai and Power propose picking one order at random.  Leinster
  proposes \emph{not} picking any order.
\end{quote}
(Unpublished document, but compare Cheng~\cite[1.3]{CheWOM}.)  Burroni%
\index{Burroni, Albert}
has
also considered `$\omega$-multigraphes', which are similar to opetopic sets
\cite{BurHDW, BurHDWA}.  Variants of Baez and Dolan's definition of weak
$n$-category have been proposed by Hermida, Makkai and Power (\latin{op.\
cit.}), Makkai~\cite{MakMOC},%
\index{Makkai, Michael}
and Cheng \cite{CheWOM, CheCOCOS, CheACU}.

Categories of trees have been used in a variety of situations: see, for
instance, Borcherds~\cite{Borch},%
\index{Borcherds, Richard}
Ginzburg%
\index{Ginzburg, Victor}
and Kapranov~\cite{GiKa},%
\index{Kapranov, Mikhail}
Kontsevich%
\index{Kontsevich, Maxim}
and Manin~\cite{KMGWC},%
\index{Manin, Yuri}
and Soibelman%
\index{Soibelman, Yan}
\cite{SoiMTC, SoiMBC}.
Some authors consider only stable trees, or consider all trees but only
surjective maps between them.  Our maps of trees are the identity on
leaf-sets, but some authors allow the leaves to be moved; doubtless this is
natural for certain parts of higher-dimensional algebra.

For more on $A_\infty$-spaces, $A_\infty$-algebras, and
$A_\infty$-categories see, for example, the book of Markl, Shnider and
Stasheff~\cite{MSS}.

The explanation~(\ref{sec:many}) of computads as being like opetopic sets
is a historical inversion; computads were introduced by Street
in~\cite{StrLIC}.  (His original computads are what we call 2-computads
here.)  The question of whether higher-dimensional computads form a
presheaf category has been wreathed in confusion: Schanuel%
\index{Schanuel, Stephen}
observed
correctly that 2-computads form a presheaf category (unpublished), Carboni%
\index{Carboni, Aurelio}
and Johnstone%
\index{Johnstone, Peter}
claimed incorrectly that this extended to $n$-computads for all
$n\in\nat$~\cite[4.6]{CJ}, Makkai and Zawadowski proved that 3-computads do
\emph{not} form a presheaf category (unpublished), and Carboni and
Johnstone published a variant of their proof in a
corrigendum~\cite{CJcorr}; meanwhile, Batanin%
\index{Batanin, Michael}
found a notion of generalized
computad and claimed that they form a presheaf category~\cite{BatCFM},
before retracting and refining that claim~\cite{BatCSO} in the light of
Makkai and Zawadowski's argument.  It should be added that defining
`$n$-computad' for $n>3$ is much easier if one makes Carboni and
Johnstone's error, and is otherwise not straightforward. 

There is a good deal more to be said about `many in, many out' than I have
said here.  See, for instance, the work of Szabo~\cite{Sza}%
\index{Szabo, Manfred}
and
Koslowski~\cite{KosMAP}%
\index{Koslowski, J\"urgen}
on polycategories,%
\index{polycategory}
and of Gan~\cite{Gan}%
\index{Gan, Wee Liang}
on
dioperads;%
\index{dioperad}
see also the remarks on PROs and PROPs herein.

\end{notes}

\part{$n$-Categories}
\label{part:n-categories}

\chapter{Globular Operads}
\lbl{ch:globular}

\chapterquote{%
Once Theofilos was [\ldots] painting a mural in a Mytilene baker's shop
[\ldots]  As was his habit, he had depicted the loaves of bread upright in
their trays, like heraldic emblems on an out-thrust shield---so that no one
could be in any doubt that they were loaves of bread and very fine ones,
too.  The irate baker pointed out that in real life loaves thus placed
would have fallen to the floor.  `No,' replied Theofilos---surely with that
calm, implacable self-certainty which carried him throughout what most
people would call a miserable life---`only real loaves fall down.  Painted
ones stay where you put them.'}
{\emph{The Athenian}~\cite{Athenian}}

\noindent
In the next two chapters we explore one possible definition of weak
$\omega$-category.  Its formal shape is very simple: we take
the category \Eee\ of globular sets and the free strict $\omega$-category
monad $T$ on \Eee, construct a certain $T$-operad $L$, and define a weak
$\omega$-category to be an $L$-algebra. 

In order to see that this is a \emph{reasonable} definition, and to get a
feel for the concepts involved, we proceed at a leisurely pace.  The
present chapter is devoted to contemplation of the monad
$T$~(\ref{sec:free-strict}), of $T$-operads ($=$~globular
operads,~\ref{sec:glob-opds}), and of their algebras~(\ref{sec:glob-algs}).
In Chapter~\ref{ch:a-defn} we define the particular globular operad $L$ and
look at its algebras---that is, at weak $\omega$-categories.  To keep the
explanation in this chapter uncluttered, the discussion of the
finite-dimensional case ($n$-categories) is also deferred to
Chapter~\ref{ch:a-defn}; we stick to $\omega$-categories here.

Globular operads are an absolutely typical example of generalized operads.
Pictorially, they are typical in that $T1$ is a family of shapes (globular
pasting diagrams), and an operation in a $T$-operad is naturally drawn as
an arrow with data fitting one of these shapes as its input.  (Compare
plain operads, where the `shapes' are mere finite sequences.)  Technically,
globular operads are typical in that $T$ is, in the terminology of
Appendix~\ref{app:special-cart}, a finitary familially representable monad
on a presheaf category.  Hence the explanations contained in this chapter
can easily be adapted to many other species of generalized operad.

\section{The free strict $\omega$-category monad}
\lbl{sec:free-strict}

In~\ref{sec:cl-strict} we defined a strict $\omega$-category as a globular
set equipped with extra structure, and a strict $\omega$-functor as a map
of globular sets preserving that structure.  There is consequently a
forgetful functor from the category \strcat{\omega} of strict
$\omega$-categories and strict $\omega$-functors to the category
\ftrcat{\scat{G}^\op}{\Set} of globular sets.  In
Appendix~\ref{app:free-strict} it is shown that this forgetful functor has
a left adjoint, that the adjunction is monadic, and that the induced monad
$(T, \mu, \eta)$%
\glo{strommon}%
\index{omega-category@$\omega$-category!strict!free}
on \ftrcat{\scat{G}^\op}{\Set} is cartesian.  

To understand this in pictorial terms, we start by considering $T1$, the
free strict $\omega$-category on the terminal globular set
\[
1 = (\cdots \pile{\rTo \\ \rTo} 1 \pile{\rTo \\ \rTo} \cdots
\pile{\rTo \\ \rTo} 1).	\\
\] 
The free strict $\omega$-category functor takes a globular set and creates
all possible formal composites in it.  A typical element of $(T1)(2)$ looks
like
\begin{equation}	\label{eq:typical-glob-pd}
\gfstsu
\gfoursu
\gzersu
\gonesu
\gzersu
\gtwosu
\glstsu,
\end{equation}
where each $k$-cell drawn represents the unique member of $1(k)$.  Note
that although this picture contains 4 dots representing 0-cells of $1$,
they actually all represent the same 0-cell; of course, $1$ only \emph{has}
one 0-cell.  The same goes for the 1- and 2-cells.  So we have drawn a
flattened-out version of the true, twisted, picture.  We call an element of
$(T1)(m)$ a \demph{(globular) $m$-pasting%
\index{pasting diagram!globular}
diagram} and write $T1=\pd$.%
\glo{pd}

Since the theory of strict $\omega$-categories includes identities, there
is for each $m\geq 2$ an element of $\pd(m)$ looking
like~\bref{eq:typical-glob-pd}.  Although the pictures look the same, they
are regarded as \emph{different}%
\lbl{p:degen-pds}%
\index{pasting diagram!globular!degenerate}
pasting diagrams for different values of $m$; the sets $\pd(m)$ and
$\pd(m')$ are considered disjoint when $m \neq m'$.  When it comes to
understanding the definition of weak $n$-category, this point will be
crucial.  

In the globular set $1$, all cells are endomorphisms---in other words, the
source and target maps are equal.  It follows that the same is true in
the globular set $\pd$.  We write $\bdry: \pd(m+1) \go \pd(m)$%
\glo{bdry}
instead of
$s$ or $t$, and call $\bdry$ the \demph{boundary}%
\index{boundary!pasting diagram@of pasting diagram}
operator.  For instance,
the boundary of the 2-pasting diagram~\bref{eq:typical-glob-pd} is the
1-pasting diagram
\[
\gfstsu
\gonesu
\gzersu
\gonesu
\gzersu
\gonesu
\glstsu.
\]

It is easy to describe%
\index{pasting diagram!globular!free monoid@via free monoids}
the globular set $\pd$ explicitly: writing
$\blank^*$%
\glo{starfreemon}
for the free monoid functor on $\Set$, we have $\pd(0)=1$ and
\lbl{p:pd-description}%
$\pd(m+1) = \pd(m)^*$.  That is, an $(m+1)$-pasting diagram is a
sequence of $m$-pasting diagrams.  For example, the $2$-pasting diagram
depicted in~\bref{eq:typical-glob-pd} is the sequence
\[
(\gfstsu\gonesu\gzersu\gonesu\gzersu\gonesu\glstsu, \ 
\glstsu, \ 
\gfstsu\gonesu\glstsu)
\]
of $1$-pasting diagrams, so if we write the unique element of $\pd(0)$ as
$\blob$ then~\bref{eq:typical-glob-pd} is the double sequence
\[
((\blob, \blob, \blob), (), (\blob)) 
\in \pd(2).
\]
The boundary map $\bdry: \pd(m+1) \go \pd(m)$ is defined inductively by
\[
\left(
\pd(m+1) \goby{\bdry} \pd(m)
\right) 
= 
\left(
\pd(m) \goby{\bdry} \pd(m-1)
\right)^*  
\]
($m\geq 1$).  The correctness of this description of $\pd$ follows from the
results of Appendix~\ref{app:free-strict}.

Having described \pd\ as a globular set, we turn to its strict
$\omega$-category structure: how pasting diagrams may be composed.

Typical binary compositions are
\begin{equation}	\label{eq:bin-1-comp}
\left.
\begin{array}{c}
\gfstsu\gthreesu\gzersu\gonesu\gzersu\gtwosu\glstsu	\\
\ofdim{1}	\\
\gfstsu\gonesu\gzersu\gonesu\gzersu\gthreesu\glstsu
\end{array}
\right.
\ =\ 
\gfstsu\gthreesu\gzersu\gonesu\gzersu\gfoursu\glstsu,
\end{equation}
illustrating the composition function 
$
\ofdim{1}: \pd(2) \times_{\pd(1)} \pd(2) \go \pd(2),
$
and
\begin{equation}	\label{eq:bin-0-comp}
\gfst{}\gspecialtwo\glst{}
\ofdim{0}
\gfst{}\gthree{}{}{}{}{}\glst{}
\ =\ 
\gfst{}\gspecialtwo\gfbw{}\gthree{}{}{}{}{}\glst{}\!,
\end{equation}
illustrating the composition function $ \ofdim{0}: \pd(3) \times_{\pd(0)}
\pd(3) \go \pd(3).  $ These compositions are possible because the
boundaries match: in~\bref{eq:bin-1-comp}, the 1-dimensional boundaries of
the two pasting diagrams on the left-hand side are equal, and similarly for
the 0-dimensional boundaries~\bref{eq:bin-0-comp}---indeed, this is
inevitable as there is only one 0-pasting diagram.

(The arguments on the left-hand side of~\bref{eq:bin-1-comp} are stacked
vertically rather than horizontally just to make the picture more
compelling; strictly speaking we should have written
\[
\gfstsu\gonesu\gzersu\gonesu\gzersu\gthreesu\glstsu
\ \ofdim{1} \ 
\gfstsu\gthreesu\gzersu\gonesu\gtwosu\glstsu
\]
instead.  The same applies to~\bref{eq:bin-0-comp}.)

A typical nullary%
\index{nullary!composite}
composition (identity) is
\begin{equation}	\label{eq:typical-id}
\gfstsu\gonesu\gzersu\gonesu\gzersu\gonesu\glstsu
\diagspace \goesto \diagspace 	
\gfstsu\gonesu\gzersu\gonesu\gzersu\gonesu\glstsu,
\end{equation}
illustrating the identity function $i: \pd(1) \go \pd(2)$ in the strict
$\omega$-category $\pd$.  (Recall the remarks above on degenerate%
\index{pasting diagram!globular!degenerate}
pasting
diagrams.)  So the left-hand side of~\bref{eq:typical-id} is a 1-pasting
diagram $\pi\in\pd(1)$, and the right-hand side is 2-pasting diagram $1_\pi
= i(\pi) \in \pd(2)$.  We have $\bdry(1_\pi) = \pi$.%
\lbl{p:bdry-degen-pd}

We will need to consider not just binary and nullary composition in \pd,
but composition `indexed' by arbitrary shapes, in the sense now explained.
The first binary composition~\bref{eq:bin-1-comp} above is indexed by the
2-pasting diagram $\gfstsu\gthreesu\glstsu$, in that we were composing one
2-cell with another by joining along their bounding 1-cells.  The
composition can be represented as
\begin{equation}	%\label{pic:bin-1-comp-expanded}
\setlength{\unitlength}{1cm}
\begin{array}{c}
\begin{picture}(10,3.6)%
\put(1,3){\gzersu\gthreesu\gzersu\gonesu\gzersu\gtwosu\gzersu}%
\put(1,0.4){\gzersu\gonesu\gzersu\gonesu\gzersu\gthreesu\gzersu}%
\put(5.5,1.8){\gzersu\gthreesu\gzersu}%
\put(9,1.8){\gzersu\gtwosu\gzersu}%
\cell{5.8}{2.15}{b}{\triangledown}
\cell{5.8}{1.55}{t}{\vartriangle}
\cell{9.3}{2.0}{b}{\triangledown}
\qbezier(4,3.4)(5.8,4.4)(5.8,2.35)
\qbezier(4,0)(5.8,-1)(5.8,1.35)
\qbezier(6.6,2.2)(9.3,4)(9.3,2.2)
\end{picture}
\end{array}
.
% \absentpiccy{pain35b.ps}.
\end{equation}
In general, the ways of composing pasting diagrams are indexed by pasting
diagrams themselves: for instance,
\begin{equation}	\label{pic:gen-comp}
\setlength{\unitlength}{1cm}
\begin{array}{c}
\begin{picture}(10,4.2)%
\put(4,2.1){\gzersu\gthreesu\gzersu\gtwosu\gzersu\gonesu\gzersu}%
\put(0,3){\gzersu\gonesu\gzersu\gonesu\gzersu}%
\put(2.5,3.7){\gzersu\gonesu\gzersu\gthreesu\gzersu}%
\put(5,4){\gzersu\gonesu\gzersu\gonesu\gzersu\gonesu\gzersu}%
\put(1,0.5){\gzersu\gtwosu\gzersu\gfoursu\gzersu}%
\put(9,2.1){\gzersu\gtwosu\gzersu}%
\cell{4.3}{2.45}{b}{\triangledown}
\cell{4.3}{1.85}{t}{\vartriangle}
\cell{5.2}{2.3}{b}{\triangledown}
\cell{6.2}{2.2}{b}{\triangledown}
\cell{9.3}{2.3}{b}{\triangledown}
\qbezier(1.9,3.1)(4.3,3.1)(4.3,2.65)
\qbezier(4.4,3.5)(5.2,3.1)(5.2,2.5)
\qbezier(5.8,3.8)(6.2,3.1)(6.2,2.4)
\qbezier(2.9,0.3)(4.3,0.6)(4.3,1.65)
\qbezier(6.9,2.4)(9.3,3.5)(9.3,2.5)
\end{picture}
\end{array}
% \absentpiccy{pain36.ps}
\end{equation}
represents the composition
\renewcommand{\qbeziermax}{80}%
\begin{eqnarray*}
&
\left(\begin{array}{c}
	\gfstsu\gonesu\gzersu\gonesu\glstsu	\\
	\ofdim{1}	\\
	\gfstsu\gtwosu\gzersu\gfoursu\glstsu
\end{array}\right)
\ \ofdim{0}\ 
\gfstsu\gonesu\gzersu\gthreesu\glstsu
\ \ofdim{0}\ 
\gfstsu\gonesu\gzersu\gonesu\gzersu\gonesu\glstsu
\\
= &
\gfstsu\gtwosu\gzersu\gfoursu\gzersu\gonesu\gzersu%
\gthreesu\gzersu\gonesu\gzersu\gonesu\gzersu\gonesu\glstsu
\end{eqnarray*}%
\renewcommand{\qbeziermax}{150}%
(with the same pictorial convention on positioning the arguments as
previously).

\paragraph*{}

This describes the free strict $\omega$-category $\pd$ on the terminal
globular set.  Before progressing to free strict $\omega$-categories in
general, let us pause to consider an alternative way of representing
pasting diagrams, due to Batanin~\cite{BatMGC}:%
\index{Batanin, Michael!globular operads@on globular operads}
as
trees.

First a warning:%
\index{tree!different types of}
these are not the same trees as appear elsewhere in this
text (\ref{sec:trees},~for instance).  Not only is there a formal
difference, but also the two kinds of trees play very different roles.  The
exact connection remains unclear, but for the purposes of understanding
what is written here they can be regarded as entirely different species.

The idea is that, for instance, the 2-pasting diagram
\[
\gfstsu%
\gfoursu%
\gzersu%
\gonesu%
\gzersu%
\gtwosu%
\glstsu
\]
can be portrayed as the tree 
\begin{equation}	\label{diag:Batanin-tree}
\begin{tree}
\node &\node            &\node &      &      &\node \\
      &\rt{1} \dn \lt{1}&      &      &      & \dn  \\
      &\node            &      &\node &      &\node \\
      &                 &\rt{2}&\dn   &\lt{2}&      \\
      &                 &      &\node &      &      \\
\end{tree}
\end{equation}
according to the following method.  The pasting diagram is 3 1-cells
long, so the tree begins life as
\[
\begin{tree}
\node	&	&\node	&	&\node	\\
	&\rt{2}	&\dn	&\lt{2}	&	\\
	&	&\node\makebox[0em]{\ \ .}	&	&	\\
\end{tree}
\]
Then the first column is 3 2-cells high, the second 0, and the third 1, so
it grows to~\bref{diag:Batanin-tree}.  Finally, there are no 3-cells so it
stops there.

Formally, an \demph{$m$-stage level%
\index{tree!level}
tree} ($m\in\nat$) is a
diagram
\[
\tau(m)\go\tau(m-1)\go\cdots\go\tau(1)\go\tau(0)=1
\]
in the skeletal category $\scat{D}$ of finite (possibly empty) totally
ordered sets~(\ref{eg:str-mon-D}); we write $\fcat{lt}(m)$ for the set of
all $m$-stage level trees.  The element of $\fcat{lt}(2)$
in~\bref{diag:Batanin-tree} corresponds to a certain diagram $4\go 3\go 1$
in $\scat{D}$, for example.  (Note that if $\tau$ is an $m$-stage tree with
$\tau(m)=0$ then the height of the picture of $\tau$ will be less than
$m$.)  The \demph{boundary}%
\index{boundary!level tree@of level tree}
$\bdry\tau$ of an $m$-stage tree $\tau$ is the
$(m-1)$-stage tree obtained by removing all the nodes at height $m$, or
formally, truncating
\[
\tau(m)\go\tau(m-1)\go\cdots\go\tau(1)\go\tau(0)
\]
to
\[
\mbox{\hspace{5.2em}}
\tau(m-1)\go\cdots\go\tau(1)\go\tau(0).
\]
This defines a diagram
\begin{equation}	%\label{eq:lt}
\cdots 
\goby{\bdry} 
\fcat{lt}(m) 
\goby{\bdry} 
\fcat{lt}(m-1) 
\goby{\bdry}
\cdots 
\goby{\bdry} 
\fcat{lt}(0)
\end{equation}
of sets and functions, hence a globular set $\fcat{lt}$ with $s = t =
\bdry$.    

\begin{propn}
There is an isomorphism of globular sets $\pd \iso \fcat{lt}$.
\end{propn}
\begin{proof}
There is one 0-stage tree, and informally it is clear that an $(m+1)$-stage
tree amounts to a finite sequence of $m$-stage trees (placed side by side
and with a root node adjoined).  Formally, take an $(m+1)$-stage tree
$\tau$, write $\tau(1) = \{1, \ldots, r\}$, and for $1\leq i \leq r$ and $0
\leq p \leq m$, let
\[
\tau_i(p)
=
\{ j \in \tau(p+1)
\such
\bdry^p (j) = i \}.
\]
Then we have a finite sequence $(\tau_1, \ldots, \tau_r)$ of $m$-stage
trees and so, inductively, a finite sequence of $m$-pasting diagrams, which
is an $(m+1)$-pasting diagram.  It is easy to check that this defines an
isomorphism.  
\done
\end{proof}

Composition and identities in the strict $\omega$-category $\pd$ can also
be expressed in the pictorial language of trees, in a simple way: see
Batanin~\cite{BatMGC} or Leinster~\cite[Ch.~II]{SHDCT}.

\paragraph*{}

Let us now see what $T$ does to an arbitrary globular set $X$.  An $m$-cell
of $TX$ is a formal pasting-together of cells of $X$ of dimension at most
$m$: for instance, a typical element of $(TX)(2)$ looks like%
\index{pasting diagram!globular!labelled}
\begin{equation}	\label{eq:labelled-glob-pd}
\gfst{a}%
\gfour{f}{f'}{f''}{f'''}{\alpha}{\alpha'}{\alpha''}%
\grgt{b}%
\gone{g}%
\glft{c}%
\gtwo{h}{h'}{\gamma}%
\glst{d},
% \absentpiccy{pain38.ps}
\end{equation}
where $a,b,c,d \in X(0)$, $f,f',f'',f''',g,h,h' \in X(1)$,
$\alpha,\alpha',\alpha'',\gamma \in X(2)$, and $s(\alpha)=f$, $t(\alpha)=f'$,
and so on. 

We can describe the functor $T$ explicitly; in the terminology of
Appendix~\ref{app:special-cart}, we give a `familial%
\index{familial representability!free strict omega-category functor@of free strict $\omega$-category functor}
representation'.  

First we associate to each pasting diagram $\pi$ the globular set
$\rep{\pi}$%
\glo{reppi}%
\index{pasting diagram!globular!globular set from}
that `looks like $\pi$'.  If $\pi$ is the unique $0$-pasting
diagram then
\[
\rep{\pi} 
= 
(
\ 
\cdots
\parpairu
\emptyset 
\parpairu
\emptyset 
\parpairu 
1
).
\]
Inductively, suppose that $m\geq 0$ and $\pi\in\pd(m+1)$: then $\pi =
(\pi_1, \ldots, \pi_r)$ for some $r\in\nat$ and $\pi_1, \ldots, \pi_r \in
\pd(m)$, and we put 
\begin{equation}	\label{eqn:rep-constr}
\rep{\pi} = (\ \cdots \parpairu \coprod_{i=1}^{r} \rep{\pi_i}(1) 
		\parpairu \coprod_{i=1}^{r} \rep{\pi_i}(0)
		\parpairu \{0, 1, \ldots, r\} ).
\end{equation}
The source and target maps in all but the bottom dimension are the evident
disjoint unions, and in the bottom dimension they are defined at
$x\in\rep{\pi_i}(0)$ by 
\[
s(x) = i-1, \diagspace t(x) = i.
\]
For example, if $\pi$ is the $2$-pasting diagram~\bref{eq:typical-glob-pd}
then $\rep{\pi}$ is of the form
\[
\cdots 
\parpairu \emptyset 
\parpairu 4
\parpairu 7
\parpairu 4
\]
where `4' means a 4-element set, etc.  This reflects the fact that the
picture~\bref{eq:typical-glob-pd} contains 4 0-cells, 7 1-cells, 4 2-cells,
and no higher cells.

The case of degenerate%
\index{pasting diagram!globular!degenerate}
pasting diagrams deserves attention.  If
$\pi\in\pd(m)$ then $1_\pi \in \pd(m+1)$ is represented by the same picture
as $\pi$; formally, $\rep{1_\pi} = \rep{\pi}$.%
\lbl{p:degen-rep}
In fact, if $\sigma \in \pd(m+1)$ then $\rep{\sigma}(m+1) = \emptyset$ if
and only if $\sigma = 1_{\bdry\sigma}$.

An $m$-cell of $TX$ is meant to be an
`$m$-pasting diagram labelled%
\index{pasting diagram!globular!labelled}
by cells of $X$', that is, an $m$-pasting
diagram $\pi$ together with a map $\rep{\pi} \go X$ of globular sets, which
suggests that there is an isomorphism
\begin{equation}	\label{eq:pd-rep-of-T}
(TX)(m) \iso \coprod_{\pi\in\pd(m)} 
\ftrcat{\scat{G}^\op}{\Set}(\rep{\pi}, X).
\end{equation}
This is proved as Proposition~\ref{propn:pds-formula}.  

For each $m$-pasting diagram $\pi$ ($m\geq 1$) there are source%
\index{pasting diagram!globular!source and target inclusions}
and target
inclusions $\rep{\bdry\pi} \parpairu \rep{\pi}$.%  
\lbl{p:rep-source-target}
For instance, when $\pi$ is the 2-pasting diagram
of~\bref{eq:typical-glob-pd}, these embed
$\rep{\bdry\pi}$ (a string of 3 1-cells) as the top and bottom edges of
$\rep{\pi}$.  The formal definition is straightforward and left as an
exercise.  Given a globular set $X$, these embeddings induce functions
$(TX)(m) \parpairu (TX)(m-1)$ for each $m$, so that $TX$ becomes a globular
set.  So we now have the desired explicit description of the free strict
$\omega$-category functor $T$.

Finally, $T$ is not just a functor but a monad.  The multiplication turns a
pasting diagram of pasting diagrams of cells of some globular set $X$ into
a single pasting diagram of cells of $X$ by `erasing the joins';
compare~\bref{pic:gen-comp}.  The unit realizes a single cell of $X$ as a
(trivial) pasting diagram of cells of $X$.

We could also try to describe the multiplication and unit explicitly in
terms of the family $(\rep{\pi})_{m\in\nat, \pi\in\pd(m)}$ `representing'
$T$.  This can be done, but takes appreciable effort and seems to be both
very complicated and not especially illuminating; as discussed in the
Notes to Appendix~\ref{app:special-cart}, the full theory of familially
representable monads on presheaf categories is currently beyond us.  But
for the purposes of this chapter, we have all the description of $T$ that
we need.

\section{Globular operads}
\lbl{sec:glob-opds}

A \demph{globular%
\index{globular operad}
operad} is a $T$-operad.  The purpose of this section is
to describe globular operads pictorially.  The more general
$T$-multicategories are not mentioned until the end of
Chapter~\ref{ch:other-defns}, and there only briefly.

A globular operad $P$ is a $T$-graph 
\begin{equation}	\label{diag:collection}
\begin{slopeydiag}
	&	&P	&	&	\\
	&\ldTo<d&	&\rdTo	&	\\
\pd = T1&	&	&	&1
\end{slopeydiag}
\end{equation}
equipped with composition and identity operations satisfying associativity
and identity axioms.  (In a standard abuse of language, we use $P$ to mean
either the whole operad or just the globular set at the apex of the
diagram.)  We consider each part of this description in turn.

A $T$-graph~\bref{diag:collection} whose object-of-objects is $1$ will be
called a \demph{collection}.%
\lbl{p:defn-collection}%
\index{collection}
So a collection is merely a globular set over $\pd$, and consists of a set
$P(\pi)$ for each $m\in\nat$ and $m$-pasting diagram $\pi$, together with a
pair of functions $P(\pi) \parpair{s}{t} P(\bdry\pi)$ (when $m\geq 1$)
satisfying the usual globularity equations.
(Formally,~\ref{propn:pshf-slice} tells us that a presheaf on $\scat{G}$
over $\pd$ is the same thing as a presheaf on the category of elements of
$\pd$, whose objects are pasting diagrams and whose arrows are generated by
those of the form $\bdry\pi \parpair{\sigma}{\tau} \pi$ subject to the
duals of the globularity equations.)

If we were discussing plain%
\index{globular operad!plain operad@\vs.\ plain operad}
rather than globular operads then $\pd =
T1$ would be replaced by $\nat$, and a collection would be a sequence
$(P(k))_{k\in\nat}$ of sets.  In that context we think of an element
$\theta$ of $P(k)$ as an operation of arity $k$ (even though it does not
actually act on anything before a $P$-algebra is specified) and draw it as
\[
k
\left\{
\mbox{\rule[-4ex]{0em}{8ex}}
\right. 
% \!\!\!\!\!\!
\begin{centredpic}
\begin{picture}(6,4)(-1,-2)
\cell{0}{0}{l}{\tusual{\theta}}
\cell{0}{0}{r}{\tinputsslft{}{}{}}
\cell{4}{0}{l}{\toutputrgt{}}
\end{picture}
\end{centredpic}
.
\]
Similarly, if $P$ is a (globular) collection and $\pi$ an $m$-pasting
diagram for some $m\in\nat$, we think of an element of $P(\pi)$ as an
`operation of arity $\pi$' and draw it as an arrow whose input is (a
picture of) $\pi$ and whose output is a single $m$-cell.  For instance, if
$m=2$ and
\[
\pi = \gfstsu\gthreesu\gzersu\gtwosu\glstsu
\]
then $\theta\in P(\pi)$ is drawn as
\[
\begin{diagram}[width=4em]
\gfstsu\gthreesu\gzersu\gtwosu\glstsu	&
\rGlobopd^\theta			&
\gfstsu\gtwosu\glstsu .			% \makebox[0em]{\ .}\\
\end{diagram}
\]
So $\theta$ is thought of as an operation capable of taking data shaped
like the pasting diagram $\pi$ as input and producing a single 2-cell as
output.  This is a figurative description but, as we shall see, becomes
literal when an algebra for $P$ is present.

Composition%
\index{globular operad!composition in}
in a globular operad $P$ is a map $\comp: P\of P \go P$ of
collections.  The collection $P\of P$ is the composite down the left-hand
diagonal of the diagram
\[
\begin{slopeydiag}
   &       &   &       &   &       &P\of P\Spbk&  &   \\
   &       &   &       &   &\ldTo  &      &\rdTo  &   \\
   &       &   &       &TP &       &      &       &P, \\
   &       &   &\ldTo<{Td}&&\rdTo<{T!}&   &\ldTo>d&   \\
   &       &T\pd&      &   &       &\pd   &       &   \\
   &\ldTo<{\mu_1}&&    &   &       &      &       &   \\
\pd&       &   &       &   &       &      &       &   \\
\end{slopeydiag}
\]
and a typical element of $(P\of P)(2)$ is depicted as
\begin{equation}	\label{pic:compn-in-operad}
\setlength{\unitlength}{1cm}
\begin{array}{c}
\begin{picture}(9.8,5)(-0.8,0)%
\put(3.5,2.5){\gzersu\gthreesu\gzersu\gtwosu\gzersu}%
\put(0,4){\gzersu\gfoursu\gzersu\gonesu\gzersu\gtwosu\gzersu}%
\put(0,1.2){\gzersu\gonesu\gzersu\gonesu\gzersu\gthreesu\gzersu}%
\put(3.5,0.4){\gzersu\gtwosu\gzersu\gtwosu\gzersu}%
\put(8,2.5){\gzersu\gtwosu\glstsu .}%
\cell{3.8}{2.85}{b}{\triangledown}
\cell{3.8}{2.25}{t}{\vartriangle}
\cell{4.6}{2.45}{t}{\vartriangle}
\cell{8.3}{2.7}{b}{\triangledown}
\cell{-0.1}{4}{r}{\pi_{1}=}
\cell{-0.1}{1.2}{r}{\pi_{2}=}
\cell{3.4}{0.4}{r}{\pi_{3}=}
\cell{3.4}{2.5}{r}{\pi=}
\qbezier(2.8,4.1)(3.8,3.8)(3.8,3.05)
\qbezier(2.8,1.2)(3.8,1.4)(3.8,2.05)
\qbezier(4.2,0.8)(4.6,1.6)(4.6,2.25)
\qbezier(5.4,2.7)(8.3,3.5)(8.3,2.9)
\cell{3.9}{3.7}{c}{\theta_1}
\cell{3.7}{1.3}{c}{\theta_2}
\cell{4.7}{1.4}{c}{\theta_3}
\cell{7.3}{3.3}{c}{\theta}
\end{picture}
\end{array}
% 
% \absentpiccy{pain42.ps}.
\end{equation}
Here $\theta_1 \in P(\pi_1)$, $\theta_2 \in P(\pi_2)$, $\theta_3 \in
P(\pi_3)$, $\theta \in P(\pi)$, and it is meant to be implicit that
$\theta_1$, $\theta_2$, and $\theta_3$ match on their sources and targets:
$t\theta_1 = s\theta_2$ and $tt\theta_1 = ss\theta_3$.  The
left-hand half of the diagram (containing the $\theta_i$'s) is an element
of the fibre over $\pi$ of the map $T!: (TP)(2) \go \pd(2)$, and the
right-hand half ($\theta$) is an element of the fibre over $\pi$ of the map
$d: P(2) \go \pd(2)$ (that is, an element of $P(\pi)$), so the whole
diagram is a 2-cell of $P\of P$.  More precisely, it is an element of $(P\of
P)(\pi\of(\pi_1, \pi_2, \pi_3))$, where 
\[
\pi\of(\pi_{1},\pi_{2},\pi_{3}) = 
\gzersu\gfoursu\gzersu\gonesu\gzersu%
\gfoursu\gzersu\gtwosu\gzersu\gtwosu\gzersu
\]
is the composite of $\pi$ with $\pi_1$, $\pi_2$, $\pi_3$ in the
$\omega$-category \pd.  So, the composition function $\comp$ of the
globular operad $P$ sends the data assembled in~\bref{pic:compn-in-operad}
to an element $\theta \of (\theta_1, \theta_2, \theta_3) \in
P(\pi\of(\pi_{1},\pi_{2},\pi_{3}))$, which may be drawn as
\[
\setlength{\unitlength}{1cm}
\begin{picture}(10,3)%
\put(2,1.5){\gzersu\gfoursu\gzersu\gonesu\gzersu%
\gfoursu\gzersu\gtwosu\gzersu\gtwosu\gzersu}%
\put(9,1.5){\gzersu\gtwosu\glstsu.}%
\cell{9.3}{1.7}{b}{\triangledown}%
\cell{1.8}{1.5}{r}{\pi\of(\pi_{1},\pi_{2},\pi_{3})=}
\cell{8.3}{2.1}{br}{\theta\of (\theta_{1},\theta_{2},\theta_{3})}
\qbezier(6.6,1.7)(9.3,2.4)(9.3,1.9)
\end{picture}	
% 
% \absentpiccy{pain43.ps}.
\]

(The `linear' notation $\pi\of(\pi_{1},\pi_{2},\pi_{3})$ and
$\theta\of(\theta_{1},\theta_{2},\theta_{3})$ should not be taken too
seriously: there is evidently no canonical order in which to put the
$\pi_i$'s.)

That $\comp: P\of P \go P$ is a map of globular sets says that composition
is compatible with source and target.  In the example above,
\begin{eqnarray*}
s(\theta\of(\theta_{1},\theta_{2},\theta_{3}))	&
=	& 
s\theta \of (s\theta_1, s\theta_3),	\\
% \diagspace
t(\theta\of(\theta_{1},\theta_{2},\theta_{3}))	&
=	& 
t\theta \of (t\theta_2, t\theta_3),
\end{eqnarray*}
where the composite 
\[
s\theta \of (s\theta_1, s\theta_3) \in 
P(\gfstsu\gonesu\gzersu\gonesu
\gzersu\gonesu\gzersu\gonesu
\gzersu\gonesu\glstsu)
\]
is as shown:
\[
\setlength{\unitlength}{1cm}
\begin{picture}(9.8,3.8)(-1.1,0.4)%
\put(3.5,2.5){\gzersu\gonesu\gzersu\gonesu\gzersu}%
\put(0,3.95){\gzersu\gonesu\gzersu\gonesu\gzersu\gonesu\gzersu}%
\put(2.9,0.5){\gzersu\gonesu\gzersu\gonesu\gzersu}%
\put(7.9,2.5){\gzersu\gonesu\gzersu}%
\cell{3.9}{2.75}{b}{\triangledown}
\cell{4.8}{2.35}{t}{\vartriangle}
\cell{8.3}{2.7}{b}{\triangledown}
\cell{3.35}{2.55}{r}{\bdry\pi=}
\cell{-0.2}{4}{r}{\bdry\pi_1=}
\cell{2.7}{0.55}{r}{\bdry\pi_3=}
\qbezier(2.9,4)(3.8,3.8)(3.9,2.95)
\qbezier(3.9,0.8)(4.8,1.5)(4.8,2.15)
\qbezier(5.4,2.7)(8.3,3.5)(8.3,2.9)
\cell{4.0}{3.7}{c}{s\theta_1}
\cell{4.85}{1.4}{c}{s\theta_3}
\cell{7.2}{3.3}{c}{s\theta}
\end{picture}
\]
and $t\theta \of (t\theta_2, t\theta_3)$ similarly.

To see what identities%
\index{globular operad!identities in}
in a globular operad $P$ are, let $\iota_{m} \in
\pd(m)$%
\glo{iotasingle}%
\index{pasting diagram!globular!unit}
be the $m$-pasting diagram
\lbl{p:looking-single}%
looking like a single $m$-cell.  (Formally, $\iota_0$ is the unique element
of $\pd(0)$ and
\[
\iota_{m+1} = (\iota_m) \in (\pd(m))^* = \pd(m+1),
\]
using the description of \pd\ on p.~\pageref{p:pd-description}).  Then
the map $1 \goby{\eta_1} \pd$ sends the unique $m$-cell of $1$ to
$\iota_m$, so the identities function $\ids: 1 \go P$ consists of an
element $1_m \in P(\iota_m)$%
\glo{idglobopd}
for each $m\in\nat$.  For instance, the
2-dimensional identity operation $1_2$ of $P$ is drawn as 
\[
\begin{diagram}[size=4em]
\gfstsu\gtwosu\glstsu			&
\rGlobopd^{1_2}				&
\gfstsu\gtwosu\glstsu .			\\
\end{diagram}
\]
That $\ids$ is a map of globular sets says that $s(1_m) = 1_{m-1} = t(1_m)$
for all $m\geq 1$.

Finally, the composition and identities in $P$ are required to obey
associativity and identity laws.  Together these say that there is only
one way of composing any `tree' of operations of the operad: for instance,
if
\[
\begin{diagram}[height=1.5em]
\star	&	&	&	&	&	&	\\
	&\rdTo>{\theta_{11}}&&	&	&	&	\\
\star	&\rTo^{\theta_{12}}&\star&&	&	&	\\
	&\ruTo>{\theta_{13}}&&\rdTo(2,3)>{\theta_1}&&&	\\
\star	&	&	&	&	&	&	\\
	&	&	&	&\star	&\rTo^{\theta}&\gfstsu\gtwosu\glstsu\\
	&	&	&\ruTo>{\theta_2}&&	&	\\
	&	&\star	&	&	&	&	\\
\end{diagram}
\]
is a diagram of the same general kind as~(\ref{pic:compn-in-operad}), with
each $\star$ representing a 2-pasting diagram, then
\[
\theta \of 
(\theta_1 \of (\theta_{11},\theta_{12},\theta_{13}),\theta_2) 
=
(\theta\of(\theta_1,\theta_2)) 
\of 
(\theta_{11},\theta_{12},\theta_{13},1_2).
\]

We have now unwound all of the data and axioms for a globular operad.
Although it may seem complicated on first reading, it is summed up simply:
a globular operad is a collection of operations together with a unique
composite for any family of operations that might plausibly be composed.

\section{Algebras for globular operads}
\lbl{sec:glob-algs}%
\index{globular operad!algebra for|(}%
\index{algebra!globular operad@for globular operad|(}

We now confirm what was suggested implicitly in the previous section: that
if $P$ is a globular operad then a $P$-algebra structure on a globular set
$X$ consists of a function
\[
\ovln{\theta}:
\{ \textrm{labellings of } \pi \textrm{ by cells of } X \}
\go
X(m)
\]%
\index{pasting diagram!globular!labelled|(}%
for each number $m$, $m$-pasting diagram $\pi$, and operation $\theta\in
P(\pi)$, satisfying sensible axioms.

So, fix a globular operad $P$.  According to the general definition, an
algebra for $P$ is an algebra for the monad $T_P$ on the category of
globular sets, which is defined on objects $X \in
\ftrcat{\scat{G}^\op}{\Set}$ by
\[
\begin{slopeydiag}
	&		&T_P X\Spbk	&		&	\\
	&\ldTo		&		&\rdTo		&	\\
TX	&		&		&		&P	\\
	&\rdTo<{T!}	&		&\ldTo>{d}	&	\\
	&		&\pd,		&		&	\\
\end{slopeydiag}
\]
that is, by
\[
(T_P X)(m) \iso
\coprod_{\pi\in\pd(m)}
P(\pi) \times \ftrcat{\scat{G}^\op}{\Set}(\rep{\pi},X).
\]
So a $P$-algebra is a globular set $X$ together with a function
\[
h_\pi: P(\pi) \times \ftrcat{\scat{G}^\op}{\Set}(\rep{\pi},X)
\go
X(m)
\]
for $m\in\nat$ and $\pi\in P(m)$, satisfying axioms.  Writing
$h_\pi(\theta, \dashbk)$ as $\ovln{\theta}$%
\glo{actionglobopd}
and recalling that a map
$\rep{\pi} \go X$ of globular sets is a `labelling of $\pi$ by cells of
$X$', we see that this is exactly the description above.
For example, suppose that 
\begin{equation}	\label{eq:lute}
\pi = 
\gfst{}\gfour{}{}{}{}{}{}{}\grgt{}\gone{}\glst{}
\in \pd(2),
\end{equation}
that $\theta\in P(\pi)$, and that
\begin{equation}	\label{eq:labelled-lute}
\mathbf{a} =
\gfst{a}%
\gfour{f}{f'}{f''}{f'''}{\alpha}{\alpha'}{\alpha''}%
\grgt{b}%
\gone{g}%
\glst{c}
\end{equation}
is a diagram of cells in $X$: then $\ovln{\theta}$ assigns to this diagram a
2-cell $\ovln{\theta}(\mathbf{a})$ of $X$.%
\index{pasting diagram!globular!labelled|)}

What are the axioms?  First, $h: T_P X \go X$ must be a map of globular
sets, which says that $\ovln{s(\theta)} = s \of \ovln\theta$ and
$\ovln{t(\theta)} = t \of \ovln\theta$.  So in our example,
$\ovln{\theta}(\mathbf{a})$ is a 2-cell of the form
\[
\gfst{d}\gtwo{k}{k'}{}\glst{e}
\]
where 
\begin{eqnarray*}
k	&= 	
	&\ovln{s(\theta)}(\gfsts{a}\gones{f}\gblws{b}\gones{g}\glsts{c}),\\
k'	&= 	
	&\ovln{t(\theta)}(\gfsts{a}\gones{f''}\gblws{b}\gones{g}\glsts{c}),\\
d	&=	&\ovln{ss(\theta)}(\gzeros{a}),			\\
e	&=	&\ovln{tt(\theta)}(\gzeros{c}).
\end{eqnarray*}

Second, $h: T_P X \go X$ must obey the usual axioms for an algebra for a
monad.  These say that composition in the operad is interpreted in the
model (algebra) as ordinary composition of functions, and identities
similarly.

An example for composition:%
\index{globular operad!composition in}
take 2-pasting diagrams
\[
\begin{array}{c}
\pi_1 = \gfstsu\gfoursu\gzersu\gonesu\glstsu,
\diagspace
\pi_2 = \gfstsu\gthreesu\glstsu,		\\
\pi = \gfstsu\gtwosu\gzersu\gtwosu\glstsu
\end{array}
\]
and write
\[
\pi \of (\pi_1, \pi_2) = 
\gfstsu\gfoursu\gzersu\gonesu\gzersu\gthreesu\glstsu.
\]
Let
\[
\begin{array}{c}
\theta_1 \in P(\pi_1), 
\diagspace
\theta_2 \in P(\pi_2),				\\
\theta \in P(\pi)
\end{array}
\]
be operations of $P$ satisfying $tt(\theta_1) = ss(\theta_2) \in
P(\gzeros{})$.  Let
\[
\begin{array}{c}
\mathbf{a}_1 = 
\gfsts{a}
\gfours{}{}{}{}{\alpha}{\alpha'}{\alpha''}
% \gfours{f}{f'}{f''}{f'''}{\alpha}{\alpha'}{\alpha''}
\grgts{b}
\gones{}
% \gones{g}
\glsts{c},
\diagspace
\mathbf{a}_2 = 
\gfsts{c}
\gthrees{}{}{}{\gamma}{\gamma'}
% \gthrees{h}{h'}{h''}{\gamma}{\gamma'}
\glsts{d},					\\
\mathbf{a} = 
\gfsts{a}
\gfours{}{}{}{}{\alpha}{\alpha'}{\alpha''}
% \gfours{f}{f'}{f''}{f'''}{\alpha}{\alpha'}{\alpha''}
\grgts{b}
\gones{}
% \gones{g}
\glfts{c}
\gthrees{}{}{}{\gamma}{\gamma'}
% \gthrees{h}{h'}{h''}{\gamma}{\gamma'}
\glsts{d}
\end{array}
\]
be diagrams of cells in $X$ (from which the 1-cell labels have been
omitted).  Then there is a composite operation $\theta \of (\theta_1,
\theta_2) \in P(\pi \of (\pi_1, \pi_2))$, and the composition-compatibility
axiom on the algebra $X$ says that
\[
\ovln{\theta \of (\theta_1, \theta_2)} (\mathbf{a}) 
=
\ovln{\theta} 
\left(
\gfst{}
\gtwo{}{}{\!\!\!\!\!\!\!\ovln{\theta_1}(\mathbf{a_1})}
\gblw{}
\gtwo{}{}{\!\!\!\!\!\!\!\ovln{\theta_2}(\mathbf{a_2})}
\glst{}
\right).
\]

An example for identities:%
\index{globular operad!identities in}
if $\alpha$ is a 2-cell of $X$ then
$\ovln{1_2}(\alpha) = \alpha$.  In general, the $m$-pasting diagram
$\iota_m$ (defined on p.~\pageref{p:looking-single}) satisfies
$\rep{\iota_m} \iso \scat{G}(\dashbk, m)$, and the identity axiom says that
\[
\ovln{1_m}: \ftrcat{\scat{G}^\op}{\Set}(\rep{\iota_m},X) \go X(m)
\]
is the canonical (Yoneda) isomorphism.  

We meet a non-trivial example of a globular operad in the next chapter.
Its algebras are, by definition, the weak $\omega$-categories.  As a
trivial example for now, the terminal globular operad $P$ is characterized
by $P(\pi)$ having exactly one element for each pasting diagram $\pi$, and
for the general reasons given in~\ref{eg:alg-terminal}, a $P$-algebra is
exactly a $T$-algebra, that is, a \emph{strict}%
\index{omega-category@$\omega$-category!strict!operad for}
$\omega$-category.%
\index{globular operad!algebra for|)}%
\index{algebra!globular operad@for globular operad|)}

\begin{notes}

Globular operads were introduced by Batanin~\cite{BatMGC}.%
\index{Batanin, Michael!globular operads@on globular operads}
 He studied them
in the wider context of `monoidal%
\index{monoidal globular category}
globular categories', and considered `operads in \cat{C}' for any monoidal
globular category \cat{C}.  There is a particular monoidal globular
category $\mathit{Span}$ such that operads in $\mathit{Span}$ are exactly
the globular operads of this chapter, which are the only kind of operads we
need in order to define weak $\omega$-categories.

The realization that Batanin's operads in $\mathit{Span}$ are just
$T$-operads, for $T$ the free strict $\omega$-category monad on globular
sets, was first recorded in my paper of~\cite{GOM}, and subsequently
explained in more detail in my~\cite{SHDCT} and~\cite{OHDCT}.

\end{notes}

\chapter{A Definition of Weak $n$-Category}
\lbl{ch:a-defn}

\chapterquote{%
Vico took it for granted that the first language of humanity was in the
form of hieroglyphics; that is, of metaphors and animated figures [\ldots]
He had intimated war with just `five real words': a frog, a mouse, a bird,
a ploughshare, and a bow}{%
Eco~\cite{Eco}}

\noindent
Algebraic structures are often defined in a way that suggests conflict:
generators \vs.\ relations, operations \vs.\ equations, composition%
\index{composition!coherence@\vs.\ coherence}%
\index{coherence!composition@\vs.\ composition}
\vs.\
coherence.  For example, in the definition of bicategory one equips a
2-globular set first with various composition operations, then with
coherence isomorphisms to ensure that some of the derived compositions are,
in fact, essentially the same.  One imagines the two sides pulling against
each other: more operations make the structure bigger and wilder, more
equations or coherence cells make it smaller and more tame.

With this picture in mind, the most obvious way to go about defining weak
$n$-category is to set up a family of higher-dimensional composition
operations subject to a family of higher-dimensional coherence constraints.
This is the strategy in Batanin's and Penon's proposed definitions, both of
which we discuss in the Chapter~\ref{ch:other-defns}.  But it is not our
strategy in this chapter.

In the definition proposed here, no distinction is made between composition
and coherence.  They are seen as two aspects of a single idea,
`contraction', not as opposing forces.  This unified approach is in many
ways more simple and graceful: one idea instead of two.

Contractions%
\index{contraction!map of globular sets@on map of globular sets}
are explained in~\ref{sec:contr}.  A map of globular sets may have the
property of being contractible, which viewed topologically means something
like being injective on homotopy; if so, it admits at least one
contraction, which is something like a homotopy%
\index{homotopy!lifting}
lifting.  These definitions
lead to definitions of contractibility of, and contraction on, a globular
operad.  By considering some low-dimensional situations, we see how
contraction alone generates a natural theory of weak $\omega$-categories. 

In~\ref{sec:defn-L} weak $\omega$-categories are defined formally as
algebras for the initial globular operad equipped with a contraction.  We
look at some examples, including the fundamental $\omega$-groupoid
of a topological space.

The finite-dimensional case, weak $n$-categories, is a shade less easy than
the infinite-dimensional case because we have to take care in the top
dimension.  In~\ref{sec:wk-n} we define weak $n$-categories and look at
various ways of constructing weak $n_1$-categories from weak
$n_2$-categories, for different (possibly infinite) values of $n_1$ and
$n_2$.  For instance, if $a$ and $b$ are 0-cells of an $n$-category $X$
then there is a `hom-$(n-1)$-category' $X(a,b)$.

The first test for a proposed definition of $n$-category is that it does
something sensible when $n\leq 2$.  We show in~\ref{sec:wk-2} that ours
passes: weak $0$-categories are sets, weak $1$-categories are categories,
and weak $2$-categories are unbiased bicategories.

\section{Contractions}
\lbl{sec:contr}

Given the language of globular operads, we need one more concept in order
to express our definition of weak $\omega$-category: contractions.  First
we define contraction on a map of globular sets, then we define contraction
on a globular operad, then we see what this has to do with the theory of
weak $\omega$-categories.

\begin{defn}
Let $X$ be a globular set and let $m\in\nat$.  Two cells $\alpha^-,
\alpha^+ \in X(m)$ are \demph{parallel}%
\index{parallel}
if $m=0$ or if $m\geq 1$,
$s(\alpha^-) = s(\alpha^+)$, and $t(\alpha^-) = t(\alpha^+)$.
\end{defn}

\begin{figure}
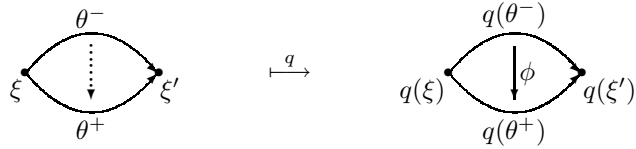

\[
\gfst{\xi}\gtwodotty{\theta^-}{\theta^+}{}\glst{\xi'}
\mbox{\hspace{2.5em}}
\stackrel{q}{\goesto}
\mbox{\hspace{2.5em}}
\gfst{q(\xi)}\gtwo{q(\theta^-)}{q(\theta^+)}{\phi}\glst{q(\xi')}
\]
\caption{Effect of a contraction $\kappa$, shown for $m=2$.  The dotted
  arrow is $\kappa_\phi(\theta^-, \theta^+)$}
\label{fig:contr-on-map}
\end{figure}

\begin{defn}	\lbl{defn:omega-map-contraction}
Let $q: X \go Y$ be a map of globular sets.  For $m\geq 1$ and $\phi\in
Y(m)$, write (Fig.~\ref{fig:contr-on-map})
\begin{eqnarray*}
\mr{Par}_q(\phi)	&
=	&
\{ (\theta^-, \theta^+) \in X(m-1) \times X(m-1) \such 
\theta^- \textrm{ and } \theta^+ \textrm{ are parallel,} 
\\ & &
q(\theta^-) = s(\phi), \textrm{ } q(\theta^+)= t(\phi)
\}.%
\glo{Par}
\end{eqnarray*}
A \demph{contraction}%
\index{contraction!map of globular sets@on map of globular sets}
$\kappa$ on $q$ is a family of functions
\[
\left(
\mr{Par}_q(\phi) \goby{\kappa_\phi} X(m) 
\right)_{m\geq 1, \phi\in Y(m)}
\]
such that for all $m\geq 1$, $\phi\in Y(m)$, and $(\theta^-, \theta^+) \in
\mr{Par}_q(\phi)$, 
\[
s(\kappa_\phi(\theta^-, \theta^+)) = \theta^-,
\ \ 
t(\kappa_\phi(\theta^-, \theta^+)) = \theta^+,
\ \ 
q(\kappa_\phi(\theta^-, \theta^+)) = \phi.
\]
\end{defn}
So a map admits a contraction when it is `injective%
\index{homotopy!injective on}
on homotopy' in some
oriented sense.

\begin{defn}	\lbl{defn:omega-opd-contraction}
A \demph{contraction}%
\index{contraction!collection@on collection}
on a collection $(P \goby{d} T1)$ is a contraction on
the map $d$.  A \demph{contraction}%
\index{contraction!globular operad@on globular operad}
on a globular operad is a contraction
on its underlying collection.  A map, collection or operad is
\demph{contractible}%
\index{contractible}
if it admits a contraction.
\end{defn}
Explicitly, if $\pi$ is a pasting diagram then $\mr{Par}_d(\pi)$ is the set
of parallel pairs $(\theta^-, \theta^+)$ of elements of $P(\bdry\pi)$.  A
contraction assigns to each pasting diagram $\pi$ and parallel $\theta^-,
\theta^+ \in P(\bdry\pi)$ an element of $P(\pi)$ with source $\theta^-$ and
target $\theta^+$.  We usually write $\mr{Par}_P(\pi)$%
\glo{Parcoll}
instead of $\mr{Par}_d(\pi)$; then a contraction $\kappa$
on $P$ consists of a function 
\[
\kappa_\pi: \mr{Par}_P(\pi) \go P(\pi)
\]
for each pasting diagram $\pi$, satisfying the source and target axioms.

To understand what contractions have to do with weak $\omega$-categories,
let us forget these definitions for a while and ask: what should the operad
for weak $\omega$-categories be?%
\index{omega-category@$\omega$-category!weak!theories of}
 In other words, how should we pick a
globular operad $L$ in such a way that $L$-algebras might reasonably be
called weak $\omega$-categories?  The elements of $L(\pi)$ are all the
possible ways of composing a labelled diagram of shape $\pi$ in an
arbitrary weak $\omega$-category, so the question is: what such ways should
there be?

There are many sensible answers.  Even for weak 2-categories there is an
infinite family of definitions, all equivalent (Chapter~\ref{ch:monoidal}).
We saw, for instance, that we could choose to start with $100$ different
specified ways of composing diagrams of shape
\[
\gfstsu\gonesu\gzersu%\gone{}\gblws{}
\gonesu\glstsu,
\]
just as long as we made them all coherently isomorphic.  This particular
choice seems bizarre, and in choosing our theory $L$ we try to make it in
some sense canonical.

So let us decide what operations to put into the theory of weak
$\omega$-categories, starting at the bottom dimension and working our way
up.  

There should not be any operations in dimension $0$ except for the
identity: the only way of obtaining a $0$-cell in a weak $\omega$-category
is to start with that same $0$-cell and do nothing.  So if $\blob$ denotes
the unique $0$-pasting diagram then we want $L(\blob) = 1$.

In dimension $1$, let us be unbiased%
\index{omega-category@$\omega$-category!unbiased|(}
and specify for each $k\in\nat$ a
single way of composing a string
\[
\gfstsu\gonesu
\diagspace \cdots \diagspace 
\gonesu\glstsu
\]
of $k$ 1-cells.  These specified compositions can be built up to
make more complex ones: if we denote the specified $k$-fold composition of
1-cells in a weak $\omega$-category by
\[
\gfsts{a_0}\gones{f_1}
\diagspace
\cdots
\diagspace
\gones{f_k}\glsts{a_k}
\diagspace
\goesto
\diagspace
\gfsts{a_0}
\gones{(f_k \of\cdots\of f_1)}
\glsts{a_k}
\]
(as we did for unbiased bicategories) then an example of a built-up
operation is
\[
\gfsts{a_0}\gones{f_1}
\gblws{a_1}\gones{f_2}
\gblws{a_2}\gones{f_3}
\gblws{a_3}\gones{f_4}
\glsts{a_4}
\diagspace\goesto\diagspace
\gfsts{a_0}\gones{((f_4\of f_3\of f_2)\of f_1)}\glsts{a_4}.
\]
There should be no other ways of obtaining a 1-cell in a weak
$\omega$-category.  So if $\chi_k$%
\glo{chi}
is the 1-pasting diagram made up of $k$
1-cells then we want $L(\chi_k)$ to be the set $\tr(k)$ of $k$-leafed
trees~(\ref{eg:opd-of-trees}).

What operations in a weak $\omega$-category result in a 2-cell? 
First,
ordinary composition makes a 2-cell from diagrams of shapes such as
\[
\rho = \gfstsu\gfoursu\glstsu,
\diagspace
\sigma = \gfstsu\gtwosu\gzersu\gtwosu\glstsu.
\]
Choosing to be unbiased again, we specify an operation $\theta\in L(\rho)$
such that $s(\theta) = t(\theta) = \id$.  When acting in a weak
$\omega$-category, $\theta$ takes as input a diagram of 2-cells $\alpha_i$
and produces as output a single 2-cell as shown:
% $\phi$ as shown:
%
\begin{equation}	\label{diag:2-dim-high}
\gfsts{a}
\gfours{f}{g}{h}{k}{\alpha}{\beta}{\gamma}
\glsts{b}
\diagspace \goesto \diagspace
\gfsts{a}
\gtwos{f}{k}{} %{\phi}
\glsts{b}.
\end{equation}
Similarly, we specify one operation of arity $\sigma$ whose source and
target are both the specified element of $L(\chi_2)$ (which is ordinary
binary composition of 1-cells):
\begin{equation}	\label{diag:2-dim-wide}
\gfsts{a}
\gtwos{f}{g}{\alpha}
\gblws{b}
\gtwos{f'}{g'}{\alpha'}
\glsts{c}
\diagspace \goesto \diagspace
\gfsts{a}
\gtwos{(f' \of f)}{(g' \of g)}{} %{\phi}
\glsts{c}.
\end{equation}
Second, there are coherence 2-cells.  For instance, a string of three
1-cells gives rise to an associativity 2-cell:
\begin{equation}	\label{diag:2-dim-ass}
\gfstsu{}\gones{f_1}\gzersu\gones{f_2}\gzersu\gones{f_3}\gzersu
\diagspace
\goesto
\diagspace
\gfstsu
\gtwos{((f_3 \of f_2) \of f_1)}{(f_3 \of (f_2 \of f_1))}{}
\glstsu
\end{equation}
and similarly in more complex cases: 
% a string $(f_1, f_2, f_3, f_4)$ of 1-cells gives a 2-cell
%
\begin{equation}	\label{diag:2-dim-coh}
\gfstsu{}\gones{f_1}\gzersu\gones{f_2}\gzersu%
\gones{f_3}\gzersu\gones{f_4}\gzersu
\diagspace
\goesto
\diagspace
\gfstsu
\gtwos{(1\of (f_4 \of f_3 \of f_2) \of 1 \of f_1)}%
{(f_4 \of (f_3 \of (f_2 \of 1 \of f_1)))}%
{}
\glstsu.
\end{equation}
Being unbiased once more, we \emph{specify} operations giving 2-cells in
each of these two ways; we do not, for instance, insist that the coherence
cell on the right-hand side of~\bref{diag:2-dim-coh} should be equal to
some composite of associativity and identity coherence cells.

Here comes the crucial point.  It looks as if these two kinds of operations
resulting in 2-cells, composition%
\index{composition!coherence@\vs.\ coherence}%
\index{coherence!composition@\vs.\ composition}
and coherence, are quite
different---complementary, even.  But they can actually be regarded as two
instances of the same thought, and that makes matters much simpler.

First recall from~\ref{sec:free-strict} that any 1-dimensional picture of a
pasting diagram can also be taken to represent a (degenerate) element of
$\pd(2)$.  In particular, the left-hand sides of~\bref{diag:2-dim-ass}
and~\bref{diag:2-dim-coh} can be regarded as degenerate elements of
$\pd(2)$: they are $1_{\chi_3}$ and $1_{\chi_4}$ respectively
(see~\bref{eq:typical-id}).  Thus, each
of~\bref{diag:2-dim-high}--\bref{diag:2-dim-coh} portrays an element of
$L(\pi)$ for some $\pi\in\pd(2)$.

Now note that four times over, we have taken a 2-pasting diagram $\pi$ and
elements $\theta^-, \theta^+ \in L(\bdry\pi)$ and decreed that $L(\pi)$
should contain a specified element $\theta$ with source $\theta^-$ and
target $\theta^+$.  In the first two cases this is obvious; in the
third~(\ref{diag:2-dim-ass}) we took $\pi = 1_{\chi_3}$ (hence $\bdry\pi =
\chi_3$) and the evident two elements $\theta^-, \theta^+ \in L(\chi_3)$; the
last is similar.

Which pairs $(\theta^-, \theta^+)$ of 1-dimensional operations should we
use to generate the 2-dimensional operations?  The simplest possible
answer, and the properly unbiased%
\index{omega-category@$\omega$-category!unbiased|)}
one, is `all of them'.  So the principle
is:
\begin{quote}
  Let $\pi \in \pd(2)$ and $\theta^-, \theta^+ \in L(\bdry\pi)$.  Then
  there is a specified element $\theta\in L(\pi)$ satisfying $s(\theta) =
  \theta^-$ and $t(\theta) = \theta^+$.
\end{quote}
`Specified' means that we take these operations $\theta$ as primitive and
generate the new 2-dimensional operations in $L$ freely using its operadic
structure, just as one dimension down we took the $k$-fold composition
operations ($k\in\nat$) as primitive and generated from them derived
1-dimensional operations, indexed by trees.  In dimension 2 the
combinatorial situation is more difficult, and I will not attempt an
explicit description of $L(\pi)$ for $\pi\in\pd(2)$.  A categorical
description takes its place; in the next section, $L$ is defined as the
universal operad containing elements specified in this way.

Because our principle takes in both composition%
\index{composition!coherence@\vs.\ coherence}%
\index{coherence!composition@\vs.\ composition}
and coherence at once, it
produces operations traditionally regarded as a hybrid of the two.  For
instance, there is a specified operation of the form
\[
\gfsts{a}
\gtwos{f}{g}{\alpha}
\gblws{b}
\gtwos{f'}{g'}{\alpha'}
\glsts{c}
\diagspace \goesto \diagspace
\gfsts{a}
\gtwos{((f' \of 1) \of f)}{(g' \of g)}{} %{\beta}
\glsts{c},
\]
where traditionally operations such as this would be built up from
horizontal composition~(\ref{diag:2-dim-wide}) and coherence cells.  The
spirit of our definition is that composition and coherence are \emph{not}
separate entities: they are two sides of the same coin.

What we have said for dimension 2 applies equally in all dimensions, with
just one small refinement: there can only be a $\theta$ satisfying
$s(\theta) = \theta^-$ and $t(\theta) = \theta^+$ if $\theta^-$ and
$\theta^+$ are parallel.  (This is trivial in dimension 2 because $L(0) =
1$.)  So the general principle is:
\begin{quote}
  Let $n\in\nat$, let $\pi \in \pd(n)$, and let $\theta^-, \theta^+$ be
  parallel elements of $L(\bdry\pi)$.  Then there is a specified element
  $\theta\in L(\pi)$ satisfying $s(\theta) = \theta^-$ and $t(\theta) =
  \theta^+$.
\end{quote}
In other words: 
\begin{quote}
  $L$ is equipped with a contraction.%
\index{contraction!globular operad@on globular operad}
\end{quote}

Everything that we did in choosing $L$ is encapsulated in this statement.
We started from the identity element of $L(0)$, part of the operad
structure of $L$.  Then we chose to specify $k$-fold compositions of
1-cells, which amounted to applying the contraction with $\pi = \chi_k$.
Then the operad structure of $L$ gave derived 1-dimensional operations such
as $(f_1, f_2, f_3) \goesto ((f_3 \of f_2) \of f_1)$.  Then we chose to
specify one operation of arity $\pi$ with any given source and target, for
each 2-pasting diagram $\pi$.  Then the operad structure of $L$ gave
derived 2-dimensional operations.  All in all, we chose $L$ to be the
universal operad equipped with a contraction.

\section{Weak $\omega$-categories}
\lbl{sec:defn-L}

We have decided that the operad for weak $\omega$-categories ought to come
with a specified contraction, and that it ought to be `universal',
`minimal', or `freely generated' as such.  Precisely, it ought to be
initial in the category of operads equipped with a contraction.  
\begin{defn}	\lbl{defn:OC}
The category $\fcat{OC}$%
\glo{OC}
of (globular) \demph{operads-with-contraction}%
\index{globular operad!contraction@with contraction}
has
as objects all pairs $(P,\kappa)$ where $P$ is a globular operad and
$\kappa$ a contraction on $P$, and as maps $(P, \kappa) \go (P', \kappa')$
all operad maps $f: P \go P'$ preserving contractions:
\[
f(\kappa_\pi(\theta^-, \theta^+)) 
=
\kappa'_\pi (f\theta^-, f\theta^+) 
\]
whenever $m\geq 1$, $\pi\in\pd(m)$, and $(\theta^-, \theta^+) \in
\mr{Par}_P(\pi)$.  
\end{defn}

\begin{propn}	\lbl{propn:OC-initial}
The category $\fcat{OC}$ has an initial object.
\end{propn}
\begin{proof}
Appendix~\ref{app:initial}. 
\done
\end{proof}
We write $(L, \lambda)$%
\glo{initL}
for the initial object of $\fcat{OC}$.  In the
previous section we described $L$ explicitly in low dimensions and saw
informally how to construct it in higher dimensions: if $L$ is known up to
and including dimension $n-1$, then $L(n)$ is obtained by first closing
under contraction then closing under $n$-dimensional operadic composition
and identities.

\begin{defn}	\lbl{defn:weak-omega-cat}%
\index{omega-category@$\omega$-category!weak}
A \demph{weak $\omega$-category} is an $L$-algebra.
\end{defn}
We write $\wkcat{\omega}$%
\glo{wkomegaCat}
for $\Alg(L)$.  Since the algebras construction
is functorial (p.~\pageref{p:Alg-functorial}), this category is determined
uniquely up to isomorphism.  

Observe that the maps in $\wkcat{\omega}$ preserve the operations from
$L$---that is, the weak $\omega$-category structure---\emph{strictly}.%
\index{omega-functor@$\omega$-functor, strict}%
\index{omega-category@$\omega$-category!weak!map of}
 In
this text we do not go as far as a definition of \emph{weak}
$\omega$-functor, nor do we reach a definition of (weak) equivalence of
weak $\omega$-categories.  These are serious omissions, and the reader may
feel cheated that we are stopping when we have barely begun, but that is
the state of the art.

Definition~\ref{defn:weak-omega-cat} is just one of many proposed
definitions of weak $\omega$-category.  In my~\cite{SDN}, it is
`Definition~\textbf{L1}'.  Its place among other definitions is discussed
in Chapter~\ref{ch:other-defns}.

\begin{example}	\lbl{eg:wk-omega-cat-contr-opd}
If $P$ is a contractible operad then any contraction $\kappa$ on $P$ gives
rise to a unique map $f: (L, \lambda) \go (P, \kappa)$ in $\fcat{OC}$,
whose underlying map $f: L \go P$ of operads induces a functor
\begin{equation}	\label{eq:P-omega-cat}
\Alg(P) \go \wkcat{\omega}.
\end{equation}
In other words, any algebra for an operad-with-contraction is canonically a
weak $\omega$-category.  All the examples below of weak $\omega$-categories
are constructed in this way.  

Functor~\bref{eq:P-omega-cat} is always faithful (being the identity on
underlying globular sets) and is full if for for each pasting diagram
$\pi$, the function $f_\pi: L(\pi) \go P(\pi)$ is surjective (exercise).
\end{example}

\begin{example}	\lbl{eg:wk-omega-cat-str}
The algebras for the terminal%
\index{globular operad!terminal}
globular operad $1$ are the strict%
\index{omega-category@$\omega$-category!strict!operad for}
$\omega$-categories (by~\ref{eg:alg-terminal}), so the unique map $L \go 1$
induces a functor
\begin{equation}	\label{eq:str-to-wk}
\strcat{\omega} \go \wkcat{\omega}
\end{equation}
---a strict $\omega$-category is a special weak $\omega$-category.

The operad $1$ admits a unique contraction, and thus becomes the
\emph{terminal} object of $\fcat{OC}$.  Contractibility of $L$ implies that
$L(\pi)$ is nonempty for each pasting diagram $\pi$, so by the observations
of the previous example, functor~\bref{eq:str-to-wk} is full and faithful.
We would expect this: if $C$ and $D$ are strict $\omega$-categories then
there ought to be a single notion of `strict map $C \go D$', independent of
whether $C$ and $D$ are regarded as strict or as weak $\omega$-categories.
\end{example}

\begin{example}	\lbl{eg:wk-omega-cat-indisc}
A directed graph is \demph{indiscrete}%
\index{graph!directed!indiscrete}
if for all objects $x$ and $y$,
there is exactly one edge from $x$ to $y$.  Such a graph has a unique
category%
\index{category!indiscrete}
structure.  The $\omega$-categorical analogue of this observation
is that any \demph{contractible}%
\index{contractible!globular set}%
\index{globular set!contractible}
globular set $X$---one for which the
unique map $X \go 1$ is contractible---has a weak $\omega$-category
structure.  To see the analogy, note that contractibility of a globular set
says that any two parallel $n$-cells $x$ and $y$ have at least one
$(n+1)$-cell $f: x \go y$ between them, and indiscreteness of a graph says
that any two objects $x$ and $y$ have a map $f: x \go y$ between them and
any two parallel arrows have an equality between them.

A contractible globular set $X$ acquires the structure of a weak
$\omega$-category as follows.  Recall from~\ref{sec:endos} that there is an
endomorphism%
\index{globular operad!endomorphism}%
\index{endomorphism!globular operad}%
\index{globular operad!algebra for}%
\index{algebra!globular operad@for globular operad}
operad $\END(X)$, and that if $P$ is any operad then a
$P$-algebra structure on $X$ is just an operad map $P \go \END(X)$.  In
particular, $X$ is canonically an $\END(X)$-algebra.  So
by~\ref{eg:wk-omega-cat-contr-opd}, it is enough to prove that
contractibility of the globular set $X$ implies contractibility of the
operad $\END(X)$.  

It follows from the definition of $\END$ that for any globular set $X$ and
pasting diagram $\pi\in\pd(n)$, an element of $(\END(X))(\pi)$ is a
sequence of functions
\[
(f_n, f_{n-1}^-, f_{n-1}^+, f_{n-2}^-, f_{n-2}^+, \ldots, f_0^-, f_0^+)
\]
making the diagram
\begin{equation}	\label{eq:endo-element}
% \[
\begin{diagram}[scriptlabels,height=2.5em]
(TX)(\pi)	&\pile{\rTo^s\\ \rTo_t}	&
(TX)(\bdry\pi)	&\pile{\rTo^s\\ \rTo_t}	&
(TX)(\bdry^2\pi)&\pile{\rTo^s\\ \rTo_t} &
\ \ \cdots \ \ 	&\pile{\rTo^s\\ \rTo_t} &
(TX)(\bdry^n\pi)	\\
\dTo~{f_n}				&	&
\dTo<{f_{n-1}^-} \dTo>{f_{n-1}^+}	&	&
\dTo<{f_{n-2}^-} \dTo>{f_{n-2}^+}	&	&
					&	&
\dTo<{f_0^-} \dTo>{f_0^+}		\\
X(n)		&\pile{\rTo^s\\ \rTo_t}	&
X(n-1)		&\pile{\rTo^s\\ \rTo_t}	&
X(n-2)		&\pile{\rTo^s\\ \rTo_t} &
\ \ \cdots \ \ 	&\pile{\rTo^s\\ \rTo_t} &
X(0)		\\
\end{diagram}
% \]
\end{equation}
commute serially---that%
\index{serially commutative}
is,
\[
s \of f_i^- = f_{i-1}^- \of s, 
\diagspace
t \of f_i^+ = f_{i-1}^+ \of t
\]
for all $i \in \{ 1, \ldots, n \}$, interpreting both $f_n^-$ and $f_n^+$
as $f_n$.  (Compare Batanin~\cite[Prop.~7.2]{BatMGC}.)%
\index{Batanin, Michael!globular operads@on globular operads}
 Contractibility of
$\END(X)$ says that given all of such a serially commutative diagram except
for $f_n$, there exists a function $f_n$ completing it.  This holds if $X$
is contractible.

Different contractions on $X$ induce different contractions on $\END(X)$,
hence different weak $\omega$-category structures on $X$.  They should all
be `equivalent', and if $X$ is nonempty then the weak $\omega$-category $X$
should be equivalent to $1$, but we do not attempt to make this precise.
\end{example}

Our final example of a weak $\omega$-category is one of the principal
motivations for the subject. 
\begin{example}	\lbl{eg:wk-omega-cat-Pi}
Any topological space $S$ gives rise to a weak $\omega$-category
$\Pi_\omega S$,%
\glo{Piomega}
its \demph{fundamental $\omega$-groupoid}.%
\index{fundamental!omega-groupoid@$\omega$-groupoid}
  Indeed, there
is a product-preserving functor
\[
\Pi_\omega: \Top \go \wkcat{\omega}.
\]

To show this, we first establish the relationship between spaces and
globular sets.  As in~\ref{eg:n-glob-set-Pi}, the Euclidean disks $D^n$
define a functor $\scat{G} \go \Top$, and by the mechanism described on
p.~\pageref{eq:induced-adjn}, this induces a pair of adjoint functors
\[
\begin{diagram}
\Top	&
\pile{\rTo^{\scriptstyle \Pi_\omega}_\top\\ 
	\lTo_{\scriptstyle | \cdot |}}	&
\ftrcat{\scat{G}^\op}{\Set}.
\end{diagram}
\]
The right adjoint is given on a space $S$ by
\[
(\Pi_\omega S)(n) = \Top(D^n, S),
\]
so $\Pi_\omega S$ is just like the singular simplicial set of $S$, but with
disks in place of simplices.  The left adjoint is \demph{geometric%
\index{geometric realization!globular set@of globular set}
realization}, given on a globular set $X$ by the coend%
\index{coend}
formula
\[
|X| = \int^{n\in \scat{G}} X(n) \times D^n.
\]
For example, if $\pi$ is an $n$-pasting diagram then $|\rep{\pi}|$ is the
space resembling the usual picture of $\pi$; it is a finite CW-complex made
up of cells of dimension at most $n$, and is contractible, as can be proved
by induction.

To give the globular set $\Pi_\omega S$ the structure of a weak
$\omega$-category we define a contractible operad $P$ (independently of
$S$) and show that $\Pi_\omega S$ is naturally a $P$-algebra.
By~\ref{eg:wk-omega-cat-contr-opd}, this suffices.

Our operad $P$ is analogous to the universal%
\index{operad!universal for loop spaces}
operad for iterated loop
spaces~(\ref{eg:opd-univ-loop}).  If $\pi$ is an $n$-pasting diagram then
an element of $P(\pi)$ is a map from $D^n$ to $|\rep{\pi}|$ respecting the
boundaries; for instance, if $\pi$ is the $1$-pasting diagram $\chi_k$
consisting of $k$ arrows in a row then an element of $P(\pi)$ is an
endpoint-preserving map $[0,1] \go [0,k]$.  In general, an element of
$P(\pi)$ is a sequence of maps
\[
\theta =
(\theta_n, 
\theta_{n-1}^-, \theta_{n-1}^+, 
\theta_{n-2}^-, \theta_{n-2}^+, 
\ldots, 
\theta_0^-, \theta_0^+)
\]
making the diagram
\[
\begin{diagram}[scriptlabels,height=2.5em]
|\rep{\pi}|		&\pile{\lTo\\ \lTo}	&
|\rep{\bdry\pi}|	&\pile{\lTo\\ \lTo}	&
|\rep{\bdry^2\pi}|	&\pile{\lTo\\ \lTo} &
\ \ \cdots \ \ 	&\pile{\lTo\\ \lTo} &
|\rep{\bdry^n\pi}|=1	\\
\uTo~{\theta_n}					&	&
\uTo<{\theta_{n-1}^-} \uTo>{\theta_{n-1}^+}	&	&
\uTo<{\theta_{n-2}^-} \uTo>{\theta_{n-2}^+}	&	&
						&	&
\uTo<{\theta_0^-} \uTo>{\theta_0^+}		\\
D^n		&\pile{\lTo\\ \lTo}	&
D^{n-1}		&\pile{\lTo\\ \lTo}	&
D^{n-2}		&\pile{\lTo\\ \lTo} &
\ \ \cdots \ \ 	&\pile{\lTo\\ \lTo} &
D^0 = 1		\\
\end{diagram}
\]
commute serially (as in~\bref{eq:endo-element}).  The maps along the
top are induced by the source and target inclusions $\rep{\bdry^{i+1}\pi}
\parpairu \rep{\bdry^i\pi}$ described on p.~\pageref{p:rep-source-target};
they are all injective, so $\theta_n$ determines the whole of $\theta$.
The obvious restriction maps $P(\pi) \parpairu P(\bdry\pi)$ make $P$ into a
collection.

We have to show that the collection $P$ is contractible, has the structure
of an operad, and acts on the globular set $\Pi_\omega S$ for any space
$S$.  Contractibility amounts to the condition that if $\pi$ is an
$n$-pasting diagram then any continuous map from the unit $(n-1)$-sphere
into $|\rep{\pi}|$ extends to the whole unit $n$-ball $D^n$, which holds
because $|\rep{\pi}|$ is contractible.  An action of $P$ on $\Pi_\omega S$
consists of a function
\[
% P(\pi) \times 
\ovln{\theta}:
\ftrcat{\scat{G}^\op}{\Set}(\rep{\pi}, \Pi_\omega S)
\go
(\Pi_\omega S)(n)
\]
for each $n\in\nat$, $\pi\in\pd(n)$, and $\theta \in P(\pi)$, satisfying
axioms involving the operad structure of $P$ (yet to be defined).  By the
adjunction above, such a function can equally be written as
\[
\ovln{\theta}:
\Top(|\rep{\pi}|, S)
\go 
\Top(D^n, S),
\]
and we take $\ovln{\theta}$ to be composition with $\theta_n$.

All that remains is to endow the collection $P$ with the structure of an
operad and check some axioms.  As suggested by the similarity between the
diagrams in this example and the last, $P$ can be constructed as an
endomorphism operad.  By the Yoneda Lemma, an element of $P(\pi)$ is a
sequence of natural transformations
\[
(\alpha_n, 
\alpha_{n-1}^-, \alpha_{n-1}^+, 
\ldots, 
\alpha_0^-, \alpha_0^+)
\]
making the diagram
\[
\begin{diagram}[scriptlabels,height=2.5em]
\Top(|\rep{\pi}|, \dashbk)		&\pile{\rTo\\ \rTo}	&
\Top(|\rep{\bdry\pi}|, \dashbk)	&\pile{\rTo\\ \rTo}	&
\ \ \cdots \ \ 	&\pile{\rTo\\ \rTo} &
\Top(|\rep{\bdry^n\pi}|, \dashbk)	\\
\dTo~{\alpha_n}					&	&
\dTo<{\alpha_{n-1}^-} \dTo>{\alpha_{n-1}^+}	&	&
						&	&
\dTo<{\alpha_0^-} \dTo>{\alpha_0^+}			\\
\Top(D^n, \dashbk)			&\pile{\rTo\\ \rTo}	&
\Top(D^{n-1}, \dashbk)		&\pile{\rTo\\ \rTo}	&
\ \ \cdots \ \ 	&\pile{\rTo\\ \rTo} &
\Top(D^0, \dashbk)			\\
\end{diagram}
\]
commute serially.  But for any $m\in\nat$ and $\rho\in\pd(m)$ there are
isomorphisms
\[
\Top(|\rep{\rho}|, S) 
\iso
\ftrcat{\scat{G}^\op}{\Set}(\rep{\rho}, \Pi_\omega S)
\iso 
(T(\Pi_\omega S))(\rho)
\]
and
\[
\Top(D^m, S)
\iso
(\Pi_\omega S)(m),
\] 
both natural in $S\in\Top$, so an element of $P(\pi)$ is a sequence of
families of functions
\begin{equation}	\label{eq:Pi-data}
\left(
(f_{S,n})_{S\in\Top}, 
(f_{S,n-1}^-)_{S\in\Top}, (f_{S,n-1}^+)_{S\in\Top}, 
\ldots, 
(f_{S,0}^-)_{S\in\Top}, (f_{S,0}^+)_{S\in\Top}
\right)
\end{equation}
natural in $S$ and making the diagram 
\[
\begin{diagram}[scriptlabels,height=2.5em]
(T(\Pi_\omega S))(\pi)	&\pile{\rTo^s\\ \rTo_t}	&
(T(\Pi_\omega S))(\bdry\pi)	&\pile{\rTo^s\\ \rTo_t}	&
\ \ \cdots \ \ 	&\pile{\rTo^s\\ \rTo_t} &
(T(\Pi_\omega S))(\bdry^n\pi)		\\
\dTo~{f_{S,n}}				&	&
\dTo<{f_{S,n-1}^-} \dTo>{f_{S,n-1}^+}	&	&
					&	&
\dTo<{f_{S,0}^-} \dTo>{f_{S,0}^+}	\\
(\Pi_\omega S)(n)		&\pile{\rTo^s\\ \rTo_t}	&
(\Pi_\omega S)(n-1)		&\pile{\rTo^s\\ \rTo_t}	&
\ \ \cdots \ \ 	&\pile{\rTo^s\\ \rTo_t} &
(\Pi_\omega S)(0)		\\
\end{diagram}
\]
commute serially for each $S\in\Top$.  This now looks
like~\bref{eq:endo-element}, an operation in an endomorphism operad.
Indeed, composition with $T$ induces a cartesian monad $T_*$ on the
category
\[
\ftrcat{\Top}{\ftrcat{\scat{G}^\op}{\Set}},
\]
and if $R\in\Top$ then evaluation at $R$ induces a strict map of monads
\[
\mr{ev}_R: 
(\ftrcat{\Top}{\ftrcat{\scat{G}^\op}{\Set}}, T_*)
\go
(\ftrcat{\scat{G}^\op}{\Set}, T)
\]
hence a functor
\[
(\mr{ev}_R)_*: 
T_* \hyph\Operad
\go
T \hyph \Operad.
\]
We therefore have a $T_*$-operad $\END(\Pi_\omega)$ and a globular operad
\[
P' = (\mr{ev}_\emptyset)_* (\END(\Pi_\omega)), 
\]
and a few calculations reveal that an element of $P'(\pi)$ consists of data
as in~\bref{eq:Pi-data}.  So $P'(\pi) \iso P(\pi)$, giving $P$ the
structure of an operad.  Further checks show that the operad structure is
compatible with the action on $\Pi_\omega S$ described above.
\end{example}

\section{Weak $n$-categories}
\lbl{sec:wk-n}

We now imitate what we did for $\omega$-categories to obtain a
definition of weak $n$-category.  This is straightforward except for one
subtlety.  We then explore the relationship between weak $n$-categories for
different values of $n$, and in particular we see how a weak $n$-category
can be regarded as a weak $\omega$-category trivial above dimension $n$.

First recall from~\ref{sec:cl-strict} that an $n$-globular set is a
presheaf on the category $\scat{G}_n$ and that there is a forgetful functor
from $\strcat{n}$ to $\ftrcat{\scat{G}_n^\op}{\Set}$.
Theorem~\ref{thm:n-forgetful-properties} implies that this induces a
cartesian monad $\gm{n}$%
\glo{gmn}%
\index{n-category@$n$-category!strict!free}
on $\ftrcat{\scat{G}_n^\op}{\Set}$.  A
$\gm{n}$-graph whose object-of-objects is $1$ will be called an
\demph{$n$-collection},%
\index{collection!n-@$n$-}
a $\gm{n}$-operad will be called an
\demph{$n$-globular operad}%
\index{n-globular operad@$n$-globular operad}
or simply an \demph{$n$-operad}, and the
category of $n$-operads will be written $n\hyph\Operad$.%
\glo{nOperad}

The subtlety concerns operations in the top%
\index{top dimension}
dimension.  In a
bicategory, for example, there are various ways of composing 2-cell
diagrams of shape
\[
\gfstsu\gfoursu\gzersu\gtwosu\gzersu\gthreesu\glstsu,
\]
but each such way is uniquely determined once one has chosen how to compose
the string of three 1-cells along the top and the string of three 1-cells
along the bottom.  This is essentially the `all diagrams commute'
coherence%
\index{coherence!bicategories@for bicategories}
theorem.  In creating the theory of weak $\omega$-categories we never
declared any two operations to be equal: the contraction on the operad $L$
deferred the relation to the next dimension.  In the finite-dimensional
case we cannot always defer to the next dimension, and we need a notion of
contraction sensitive to that.

\begin{defn}
Let $n\in\nat$.  A map $q: X \go Y$ of $n$-globular sets is \demph{tame}%
\index{tame!map of n-globular sets@map of $n$-globular sets}
if
any two parallel $n$-cells $\alpha^-, \alpha^+$ of $X$ satisfying
$q(\alpha^-) = q(\alpha^+)$ are equal.
\end{defn}

\begin{defn}
Let $n\in\nat$.  A \demph{precontraction}%
\index{precontraction|(}
on a map $q: X \go Y$ of
$n$-globular sets is a family of functions
\[
\left(
\mr{Par}_q(\phi) \goby{\kappa_\phi} X(m) 
\right)_{1 \leq m \leq n, \phi\in Y(m)}
\]
satisfying the same conditions as those in
Definition~\ref{defn:omega-map-contraction}.  A \demph{contraction} is a
precontraction on a tame map.
\end{defn}
Tameness means injectivity in the top dimension (relative to sources and
targets), so contractibility continues to mean something like `injectivity%
\index{homotopy!injective on}
on homotopy', as in the infinite-dimensional case.  

\begin{defn}
An $n$-collection $(P \goby{d} \gm{n}1)$ is \demph{tame}%
\index{tame!n-collection or n-operad@$n$-collection or $n$-operad}
if the map $d$ is
tame, and a \demph{(pre)contraction}%
\index{contraction!n-collection or n-operad@on $n$-collection or $n$-operad}%
\index{precontraction|)}
on the $n$-collection is a
(pre)contraction on the map $d$.
\end{defn}
Tameness of an $n$-collection $P$ means that if for all $\pi\in\pd(n)$,
parallel elements of $P(\pi)$ are equal, or equivalently the function
\[
(s,t): P(\pi) \go \mr{Par}_P(\pi)
\]
is injective.  If $P$ is precontractible then $(s,t)$ is surjective, with
one-sided inverse $\kappa_\pi$.  Hence a precontraction is a contraction if
and only if $(s,t)$ is bijective for all $\pi\in\pd(n)$.  This implies that
a contractible $n$-operad is entirely%
\index{top dimension}
determined by its lower
$(n-1)$-dimensional part; if $\pi\in\pd(n)$ then $P(\pi)$ must be the set of
parallel pairs of elements of $P(\bdry\pi)$.  

\begin{example}	\lbl{eg:sesqui}
A \demph{sesquicategory}%
\index{sesquicategory}
is a category $X$ together with a functor
$\mr{HOM}: X^\op \times X \go \Cat$ such that
\[
\begin{diagram}[height=2em]
X^\op \times X	&\rTo^{\mr{HOM}}	&\Cat		\\
		&\rdTo<{\Hom}		&\dTo>{\mr{ob}}	\\
		&			&\Set		\\
\end{diagram}
\]
commutes.  The objects of $X$ are called the \demph{0-cells} of the
sesquicategory, and if $a$ and $b$ are 0-cells then the objects and arrows
of the category $\mr{HOM}(a,b)$ are called \demph{1-cells} and
\demph{2-cells} respectively.  Any strict 2-category has an underlying
sesquicategory.

Concretely, a sesquicategory consists of a 2-globular set together with
structure and axioms saying that any diagram of shape
\[
\gfstsu\gonesu
\diagspace 
\cdots
\diagspace 
\gonesu\glstsu
\]
has a unique composite 1-cell and any diagram of shape%
\index{whiskering}
\[
\gfstsu\gonesu\ \cdots\ \gonesu\gzersu
\gdotssu
\gzersu\gonesu\ \cdots\ \gonesu\glstsu
\]
has a unique composite 2-cell.  (The sequences of cells shown can have any
length, including zero.)  See Leinster~\cite[III.1]{SHDCT} or
Street~\cite[\S 2]{StrCS}%
\index{Street, Ross}
for a more detailed presentation.  What a sesquicategory does not have is a
canonical horizontal%
\index{composition!horizontal}
composition of
2-cells: given a diagram
\begin{equation}	\label{diag:ssq-horiz}
\gfsts{a}
\gtwos{f}{g}{\alpha}
\gfbws{a'}
\gtwos{f'}{g'}{\alpha'}
\glsts{a''},
\end{equation}
the only resulting 2-cells of the form
\[
\gfsts{a}
\gtwos{f'\sof f}{g'\of g}{}
\glsts{a''}
\]
are the derived composites $(\alpha' g) \of (f' \alpha)$ and $(g' \alpha)
\of (\alpha' f)$, which in general are not equal.

There is a 2-operad $\fcat{Ssq}$ whose algebras are sesquicategories.  If
$\pi$ is a 0- or 1-pasting diagram then $\fcat{Ssq}(\pi) = 1$.  A 2-pasting
diagram is a finite sequence $(k_1, \ldots, k_n)$ of natural numbers, as
observed in~\ref{sec:free-strict}; pictorially, $n$ is the width of the
diagram and $k_i$ the height of the $i$th column.  Writing $\mb{m} = \{ 1,
\ldots, m\}$, an element of $\fcat{Ssq}(k_1, \ldots, k_n)$ is a total order%
\index{order!non-canonical}
on the disjoint union $\mb{k_1} + \cdots + \mb{k_n}$ that restricts to the
standard order on each $\mb{k_i}$; hence
\[
|\fcat{Ssq}(k_1, \ldots, k_n)|
=
\frac{(k_1 + \cdots + k_n)!}{k_1! k_2! \cdots k_n!}.
\]

A precontraction on $\fcat{Ssq}$ amounts to a choice of element of
$\fcat{Ssq}(k_1, \ldots, k_n)$ for each $n, k_1, \ldots, k_n \in \nat$, and
since this set always has at least one element, $\fcat{Ssq}$ is
precontractible.  But it is not contractible, or equivalently not tame: for
since all elements of $\fcat{Ssq}(k_1, \ldots, k_n)$ are parallel, this
would say that $|\fcat{Ssq}(k_1, \ldots, k_n)| = 1$ for all $n, k_1,
\ldots, k_n$, which is false.
\end{example}

The categories $n\hyph\fcat{OP}$%
\glo{nOP}
of $n$-operads-with-precontraction and
$n\hyph\fcat{OC}$%
\glo{nOC}
of $n$-operads-with-contraction are defined analogously
to $\fcat{OC}$~(\ref{defn:OC}).  

\begin{propn}	\lbl{propn:n-OC-initial}%
\index{n-globular operad@$n$-globular operad!initial with contraction}
For each $n\in\nat$, the category $n\hyph\fcat{OC}$ has an initial object.
\end{propn}
We write $(L_n, \lambda_n)$%
\glo{initLn}
for the initial object.  The proof comes
later~(\ref{cor:Ln-from-L}).

\begin{defn}%
\index{n-category@$n$-category!weak}
Let $n\in\nat$.  A \demph{weak $n$-category} is an $L_n$-algebra.
\end{defn}
We write $\wkcat{n}$%
\glo{wkncat}
for $\Alg(L_n)$, the category of weak $n$-categories
and strict $n$-functors.  This contains $\strcat{n}$ as a full subcategory,
just as in the infinite-dimensional case~(\ref{eg:wk-omega-cat-str}).

\index{n-category@$n$-category!weak!varying n@for varying $n$|(}
We now embark on a comparison of the theories of $n$-categories for varying
values of $n$, including $n=\omega$.  Among other things we prove
Proposition~\ref{propn:n-OC-initial}, the strategy for which is as follows.
Imagine $L_n$ being built from the bottom dimension upwards, just as we did
for $L$: then it is clear that $L_n$ and $L$ should be the same up to and
including dimension $n-1$.  They are different in dimension%
\index{top dimension}
$n$, as
discussed above in the case $n=2$: if $\pi$ is an $n$-pasting diagram then
an element of $L_n(\pi)$ is a parallel pair of elements of $L_n(\bdry\pi)$.
So we define an $n$-operad from $L$ by discarding all the operations of
dimension $n$ and higher, then adding in parallel pairs as the operations
in dimension $n$, and we prove that this is the initial
$n$-operad-with-contraction.

Let $m, M \in \nat \cup \{\omega\}$ with $m\leq M$.  We compare $m$- and
$M$-dimensional structures, using the following natural conventions for the
$\omega$-dimensional case:
\[
\scat{G}_\omega = \scat{G},
\ \ 
\gm{\omega} = T,
\ \ 
\omega\hyph\Operad = T\hyph\Operad,
\ \ 
\omega\hyph\fcat{OP} = \omega\hyph\fcat{OC} = \fcat{OC}.%
\glo{Gomega}\glo{omegaOperad}\glo{omegaOC}
\]

We start with globular sets, and the adjunction 
\[
\begin{diagram}[height=2em]
\ftrcat{\scat{G}_M^\op}{\Set}	\\
\dTo<{R^M_m}	\ladj	\uTo>{I^M_m}	\\
\ftrcat{\scat{G}_m^\op}{\Set}.
\end{diagram}
\]
Here $R_M^m$%
\glo{RMm}
is \demph{restriction}:%
\index{restriction}%
\index{globular set!restriction of}
if $X$ is an $M$-globular set then
$(R^M_m X)(k) = X(k)$ for all $k\leq m$.  Its right adjoint $I^M_m$%
\glo{IMm}
forms
the \demph{indiscrete}%
\index{globular set!indiscrete}
$M$-globular set on an $m$-globular set $Y$, 
\[
\cdots 
\ 
\pile{\rTo^1\\ \rTo_1}
S
\pile{\rTo^1\\ \rTo_1}
S
\pile{\rTo^{\mr{pr}_1}\\ \rTo_{\mr{pr}_2}}
Y(m)
\pile{\rTo^s\\ \rTo_t}
Y(m-1)
\pile{\rTo^s\\ \rTo_t}
\ 
\cdots 
\ 
\pile{\rTo^s\\ \rTo_t}
Y(0)
\]
where $S$ is the set of parallel pairs of elements of $Y(m)$.  This
adjunction is familiar in the case $M=1$, $m=0$: $R^1_0$ assigns to a
directed graph its set of objects (vertices), and $I^1_0$ forms the
indiscrete graph on a set~(\ref{eg:wk-omega-cat-indisc}).  Formally,
$R^M_m$ is composition with, and $I^M_m$ right Kan extension along, the
obvious inclusion $\scat{G}_m \rIncl \scat{G}_M$.

Next we move to the level of operads.  The functor $R^M_m$ is naturally a
weak map of monads
\[
(\ftrcat{\scat{G}_M^\op}{\Set}, \gm{M}) 
\go 
(\ftrcat{\scat{G}_m^\op}{\Set}, \gm{m})
\]
(as shown in Appendix~\ref{app:free-strict}), and preserves limits.  Hence
$R^M_m$ is naturally a cartesian colax map of monads, and by
Example~\ref{eg:cart-adjts}, there is an induced adjunction
\[
\begin{diagram}[height=2em]
M\hyph\Operad	\\
\dTo<{(R^M_m)_*} \ladj	\uTo>{(I^M_m)_*}	\\
m\hyph\Operad.
\end{diagram}%
\index{globular operad!restriction of}%
\index{globular operad!indiscrete}
\]
The functor $(R^M_m)_*$ forgets the top $(M-m)$ dimensions of an
$M$-operad, and $(I^M_m)_*$ is defined on an $m$-operad $Q$ at a
$k$-pasting diagram $\pi$ by
\[
\left((I^M_m)_* Q \right)
(\pi)
=
\left\{
\begin{array}{ll}
\mr{Par}_Q(\bdry^{k-m-1}\pi)	&
			\textrm{if } k > m	\\
Q(\pi)			&\textrm{if } k \leq m.
\end{array}
\right.
\]

Now we bring in precontractions.  A precontraction on an $M$-operad $P$
restricts to a precontraction on the $m$-operad $(R^M_m)_* P$, and
similarly for $(I^M_m)_*$, so the previous adjunction lifts to an
adjunction
\[
\begin{diagram}[height=2em]
M\hyph\fcat{OP}	\\
\dTo<{(R^M_m)_*} \ladj	\uTo>{(I^M_m)_*}	\\
m\hyph\fcat{OP}.
\end{diagram}
\]

So far, we have proved:
\begin{propn}	\lbl{propn:change-of-dimension}
For any $m, M \in \nat \cup \{\omega\}$ with $m\leq M$, there are
adjunctions
\[
\begin{diagram}[height=2em]
\ftrcat{\scat{G}_M^\op}{\Set}	\\
\dTo<{R^M_m}	\ladj	\uTo>{I^M_m}	\\
\ftrcat{\scat{G}_m^\op}{\Set}
\end{diagram}
\diagspace
\begin{diagram}[height=2em]
M\hyph\Operad	\\
\dTo<{(R^M_m)_*} \ladj	\uTo>{(I^M_m)_*}	\\
m\hyph\Operad
\end{diagram}
\diagspace
\begin{diagram}[height=2em]
M\hyph\fcat{OP}	\\
\dTo<{(R^M_m)_*} \ladj	\uTo>{(I^M_m)_*}	\\
m\hyph\fcat{OP}
\end{diagram}
\]
defined by restriction $R^M_m$ and indiscrete extension $I^M_m$.  
\done
\end{propn}

Finally, we consider contractions themselves.  We have seen that an
$n$-operad-with-contraction is entirely determined by its underlying
$(n-1)$-operad-with-precontraction.  In fact, the two types of structure
are the same: 

\begin{propn}	\lbl{propn:OP-OC}
For any positive integer $n$, the third adjunction of
Proposition~\ref{propn:change-of-dimension} in the case $M=n$,
$m=n-1$ restricts to an equivalence
\[
\begin{diagram}[height=2em]
n\hyph\fcat{OC}	\\
\dTo<{(R^n_{n-1})_*} \eqv	\uTo>{(I^n_{n-1})_*}	\\
(n-1)\hyph\fcat{OP}.
\end{diagram}
\]
\end{propn}
\begin{proof}
Any adjunction $\cat{C} \pile{\rTo_\bot \\ \lTo} \cat{D}$ restricts to an
equivalence between the full subcategory of $\cat{C}$ consisting of those
objects at which the unit of the adjunction is an isomorphism and the full
subcategory of $\cat{D}$ defined dually.  In the case at hand we have
$(R^n_{n-1})_* \of (I^n_{n-1})_* = 1$, and the counit is the identity
transformation.  Given an $n$-operad $P$ with precontraction $\kappa$, the
unit map
\[
P \go (I^n_{n-1})_* (R^n_{n-1})_* P
\]
is the identity in dimensions less than $n$, and in dimension $n$ is made
up of the functions
\[
(s,t): P(\pi) \go \mr{Par}_P(\pi)
\]
($\pi\in\pd(n)$).  So $P$ belongs to the relevant subcategory of
$n\hyph\fcat{OP}$ if and only if this map is a bijection for all
$\pi\in\pd(n)$, that is, $\kappa$ is a contraction.
\done
\end{proof}

\begin{example}	\lbl{eg:Gray}
Let $\fcat{Ssq}$ be the 2-operad for sesquicategories~(\ref{eg:sesqui}).%
\index{sesquicategory}
The resulting 3-operad $\fcat{Gy} = (I^3_2)_* (\fcat{Ssq})$ is given on
$k$-pasting diagrams $\pi$ by
\[
\fcat{Gy}(\pi)
=
\left\{
\begin{array}{ll}
\fcat{Ssq}(\bdry\pi) \times \fcat{Ssq}(\bdry\pi)	&
\textrm{if } k=3,	\\
\fcat{Ssq}(\pi)		&
\textrm{if } k\leq 2.
\end{array}
\right.
\]
A $\fcat{Gy}$-algebra is what we will call a \demph{Gray-category}.%
\index{Gray-category}
 It
consists of a 3-globular set $X$ with a sesquicategory structure on its 0-,
1- and 2-cells and further structure in the top dimension: given a
3-pasting diagram labelled by cells of $X$, the ways of composing it to
yield a single 3-cell correspond one-to-one with the ways of composing both
the 2-dimensional diagram at its source and the 2-dimensional diagram at
its target.  

Let us explore the 3-dimensional operations.  By considering 3-pasting
diagrams $\pi$ for which $\bdry^2\pi$ is a single 1-cell, we see that for
any two 0-cells $a$ and $b$, the 2-globular set $X(a,b)$ has the structure
of a strict 2-category.  More subtly, take the 2-pasting diagram
\[
\rho 
=
\gfstsu\gtwosu\gzersu\gtwosu\glstsu,%
\index{composition!horizontal}
\]
and recall that $\bdry(1_\rho) = \rho$ (p.~\pageref{p:bdry-degen-pd}).  We
saw in~\ref{eg:sesqui} that $\fcat{Gy}(\rho) = \fcat{Ssq}(\rho)$
has exactly two elements $\theta_1$, $\theta_2$, sending the data
in~\bref{diag:ssq-horiz} to the derived composites
\[
\gamma_1 = (\alpha' g) \of (f' \alpha),
\diagspace
\gamma_2 = (g' \alpha) \of (\alpha' f) 
\]
respectively.  So $\fcat{Gy}(1_\rho) = \{ \theta_1, \theta_2 \}^2$.  The
element $(\theta_1, \theta_2)$ of $\fcat{Gy}(1_\rho)$ sends a cell
diagram~\bref{diag:ssq-horiz} in a Gray-category to a 3-cell of the form
\[
\gfst{a}
\gthreecell{f'\of f}{g'\of g}{\gamma_1}{\gamma_2}{}
\glst{a''}.
\]
Similarly, the data of~\bref{diag:ssq-horiz} is sent by $(\theta_2,
\theta_1)$ to a 3-cell $\gamma_2 \go \gamma_1$, and by $(\theta_i,
\theta_i)$ to the identity 3-cell on $\gamma_i$.  Since these are the only
four elements of $\fcat{Gy}(1_\rho)$, the two 2-cells $\gamma_2 \oppairu
\gamma_1$ are mutually inverse.  In other words, in a Gray-category there
is a specified isomorphism between the two derived horizontal
compositions
of 2-cells.

Here Gray-categories are treated as structures in their own right, but they
were introduced by Gordon,%
\index{Gordon, Robert}
Power and Street~\cite{GPS} as a special kind of tricategory%
\index{tricategory}
(their notion of weak 3-category).  The embedding of
$\fcat{Gy}$-algebras into the class of tricategories is non-canonical,
amounting to the choice of one of the two derived horizontal compositions
of 2-cells.  `Left-handed%
\index{handedness}
Gray-categories' and `right-handed
Gray-categories' therefore form (different) subclasses of the class of all
tricategories.  Gordon, Power and Street made the arbitrary decision to
consider the left-handed version (let us say), and proved the important
result that every tricategory is equivalent%
\index{coherence!tricategories@for tricategories}
to some left-handed
Gray-category.  By duality, the same result also holds for right-handed
Gray-categories.  

We do not set up a notion of weak equivalence of our weak $n$-categories,
so cannot attempt an analogous coherence theorem, but we can show how to
realize Gray-categories as weak 3-categories.  Choose a precontraction on
$\fcat{Ssq}$.  By the last proposition, this induces a contraction on
$(I^3_2)_*(\fcat{Ssq}) = \fcat{Gy}$, hence a map $L_3 \go \fcat{Gy}$ of
3-operads, hence a functor from Gray-categories to weak 3-categories (full
and faithful, in fact).  The final functor depends on the precontraction
chosen.  A precontraction is~(\ref{eg:sesqui}) a choice for each $n, k_1,
\ldots, k_n$ of a total order on the set $\mb{k_1} + \cdots + \mb{k_n}$
restricting to the standard order on each $\mb{k_i}$, and there are two
particularly obvious ones: order the set $\{1, \ldots, n\}$ either
backwards or forwards, then order $\mb{k_1} + \cdots + \mb{k_n}$
lexicographically.  The two embeddings induced are the analogues in our
unbiased world of the left- and right-handed embeddings in the biased world
of tricategories.
\end{example}

We can now read off results relating the theories of $n$- and
$\omega$-categories.
\begin{cor}	\lbl{cor:Ln-from-L}
$(I^n_{n-1})_* (R^\omega_{n-1})_* (L, \lambda)$ is an initial object of
$n\hyph\fcat{OC}$, for any positive integer $n$. 
\end{cor}
\begin{proof}
By~\ref{propn:change-of-dimension} and~\ref{propn:OP-OC}, we have a diagram
\[
\begin{diagram}[height=2.5em] %[height=2em,width=6em]
\omega\hyph\fcat{OC}	&	&	\\
\dTo<{(R^\omega_{n-1})_*} \ladj \uTo>{(I^\omega_{n-1})_*}	&
	&	\\
(n-1)\hyph\fcat{OP}	&
\pile{\lTo^{(R^n_{n-1})_*}_{\eqv}\\ \\ \rTo_{(I^n_{n-1})_*}}	&
n \hyph \fcat{OC},	\\
\end{diagram}
\]
and left adjoints and equivalences preserve initial objects.
\done
\end{proof}
This proves Proposition~\ref{propn:n-OC-initial}, the existence of an
initial $n$-operad-with-contraction, for $n\geq 1$, and tells us that
\[
(L_n, \lambda_n) 
\iso
(I^n_{n-1})_* (R^\omega_{n-1})_* (L, \lambda).
\]
(The case $n=0$ is done explicitly in~\ref{sec:wk-2}.)  So $L_n$ is
constructed from $L$ by first forgetting all of $L$ above dimension $n-1$,
then adjoining the only possible family of $n$-dimensional operations that
will make the resulting $n$-operad contractible.  In particular, $L_n$ and
$L$ agree up to and including dimension $n-1$, as do their associated
contractions: 
\begin{cor}	\lbl{cor:lower-same}
For any positive integer $n$, there is an isomorphism 
\[
(R^n_{n-1})_* (L_n, \lambda_n) 
\iso
(R^\omega_{n-1})_* (L, \lambda)
\]
of $(n-1)$-operads-with-precontraction.
\done
\end{cor}

The ideas we have discussed suggest two alternative definitions of weak
$n$-category, which we now formulate and prove equivalent to the main one.

The first starts from the thought that the cells of dimension at most $n$
in a weak $\omega$-category do not usually form a weak $n$-category, but
they should do if the composition of $n$-cells is strict enough.  
\begin{defn}%
\index{tame!algebra for globular operad}%
\index{algebra!globular operad@for globular operad!tame}%
\index{globular operad!algebra for!tame}
Let $n\in\nat$ and let $P$ be an $n$-operad.  A $P$-algebra $X$ is
\demph{tame} if 
\[
\ovln{\theta^-} = \ovln{\theta^+}: (TX)(\pi) \go X(n)
\]
for any $\pi\in\pd(n)$ and parallel $\theta^-, \theta^+ \in P(\pi)$.
\end{defn}
Any algebra for a tame $n$-operad, and in particular any algebra for a
contractible $n$-operad, is tame.  A tame $((R^\omega_n)_* L)$-algebra
ought to be the same thing as a weak $n$-category.  

To state this precisely, note that the unit of the adjunction
$(R^n_{n-1})_* \ladj (I^n_{n-1})_*$ gives a map
% \[
% 
\begin{equation}	\label{eq:tame-unit}
\alpha:
(R^\omega_n)_* L 
\go 
(I^n_{n-1})_* (R^n_{n-1})_* (R^\omega_n)_* L
\iso
L_n
\end{equation}
of operads, which induces a functor
\begin{equation}	\label{eq:tame-eqv}
\wkcat{n} \go \Alg((R^\omega_n)_* L).
\end{equation}
\begin{propn}	\lbl{propn:tame-n-cat}
The functor~\bref{eq:tame-eqv} restricts to an equivalence between
$\wkcat{n}$ and the full subcategory of $\Alg((R^\omega_n)_* L)$ consisting
of the tame algebras.
\end{propn}
The proof uses a rather technical lemma, whose own proof is
straightforward: 
\begin{lemma}
Let $S$ and $S'$ be monads on a category $\cat{C}$ and let $\psi: S \go
S'$ be a natural transformation commuting with the monad structures.  Write
$\cat{D}$ for the full subcategory of $\cat{C}^S$ consisting of the
$S$-algebras $(SX \goby{h} X)$ for which $h$ factors through $\psi_X$.  If
each component of $\psi$ is split epi then the induced functor
$\cat{C}^{S'} \go \cat{C}^S$ restricts to an equivalence $\cat{C}^{S'} \go
\cat{D}$.  
\done
\end{lemma}

\begin{prooflike}{Proof of Proposition~\ref{propn:tame-n-cat}}
The $n$-operads $(R^\omega_n)_* L$ and $L_n$ induce respective monads
$(\gm{n})_{(R^\omega_n)_* L}$ and $(\gm{n})_{L_n}$ on
$\ftrcat{\scat{G}_n^\op}{\Set}$, and the map $\alpha$
of~\bref{eq:tame-unit} induces a transformation $\psi$ from the first
monad to the second, commuting with the monad structures.
If $X$ is any $n$-globular set then the map $\psi_{X}$ of $n$-globular
sets is split epi: in dimension $k$ it is the function
\[
\coprod_{\pi\in\pd(k)} L(\pi) \times (\gm{n} X)(\pi)
\goby{\coprod \alpha_\pi \times 1} 
\coprod_{\pi\in\pd(k)} L_n(\pi) \times (\gm{n} X)(\pi),
\]
and $\alpha_\pi$ is bijective when $k<n$ (Corollary~\ref{cor:lower-same}),
so it is enough to show that $\alpha_\pi$ is surjective when $k=n$; this in
turn is true because
\[
\alpha_\pi = (s,t): L(\pi) \go \mr{Par}_L(\pi)
\]
and $L$ is precontractible.  

The lemma now applies, and we have only to check that $\cat{D}$ is the
subcategory consisting of the tame $((R^\omega_n)_* L)$-algebras.  Since 
$\psi$ is the identity in dimensions less than $n$, an
$((R^\omega_n)_* L)$-algebra $X$ is in $\cat{D}$ if and only if there is a
factorization
\[
\begin{diagram}[height=2em]
L(\pi) \times (\gm{n} X)(\pi)	&
\rTo^{\alpha_\pi \times 1}	&
L_n(\pi) \times (\gm{n} X)(\pi)	\\
				&
\rdTo<{\mr{action}}		&
\dGet				\\
				&
				&
X(n)
\end{diagram}
\]
in the category of sets for each $\pi\in\pd(n)$.  We have already seen that
$\alpha_\pi$ identifies two elements of $L(\pi)$ just when they are
parallel, so this is indeed tameness.
\done 
\end{prooflike}

\begin{cor}	\lbl{cor:restriction}
Let $n\in\nat$ and $N\in \nat\cup \{\omega\}$, with $n\leq N$.  If $X$
is a weak $N$-category with the property that for all $\pi\in\pd(n)$ and
parallel $\theta^-, \theta^+ \in L_N(\pi)$, 
\[
\ovln{\theta^-} = \ovln{\theta^+}: (\gm{N} X)(\pi) \go X(n),
\]
then its $n$-dimensional restriction $R^N_n X$ inherits the structure of a
weak $n$-category.
\end{cor}
\begin{proof}
The $L_N$-algebra structure on $X$ induces an $((R^N_n)_* L_N)$-algebra
structure on $R^N_n X$ (see~\ref{sec:change}).  If $n=N$ the result is
trivial; otherwise $(R^N_n)_* L_N \iso (R^\omega_n)_* L$
by~\ref{cor:lower-same}, and the result follows
from~\ref{propn:tame-n-cat}.  \done
\end{proof}

The second alternative definition says that a weak $n$-category is a weak
$\omega$-category with only identity cells in dimensions higher than $n$.
We show that this is equivalent to the main definition of weak
$n$-category.  More generally, we show that if $n \leq N$ then a weak
$N$-category trivial above dimension $n$ is the same thing as a weak
$n$-category.  The most simple case is that a discrete category is the same
thing as a set.

\begin{defn}
Let $n\in\nat$ and $N\in \nat\cup \{\omega\}$, with $n\leq N$.  An
$N$-globular set $X$ is \demph{$n$-dimensional}%
\index{globular set!n-dimensional@$n$-dimensional}
if for all $n\leq k < N$,
the maps
\[
X(k+1) \parpair{s}{t} X(k)
\]
are equal and bijective.
\end{defn}
Any $n$-dimensional $N$-globular set is isomorphic to a `strictly
$n$-dimensional' $N$-globular set, that is, one of the form
\begin{equation}	\label{diag:disc-glob-set}
\cdots
\parpair{1}{1}
X(n)
\parpair{1}{1}
X(n)
\parpair{s}{t}
X(n-1)
\parpair{s}{t}
\ 
\cdots
\ 
\parpair{s}{t}
X(0).
\end{equation}

\begin{thm}	\lbl{thm:n-dim}
Let $n\in\nat$ and $N\in \nat\cup \{\omega\}$, with $n\leq N$.  There
is an equivalence of categories
\[
\wkcat{n}
\eqv
(n \textup{-dimensional weak } N \textup{-categories})
\]
where the right-hand side is a full subcategory of $\wkcat{N}$.
\end{thm}
To prove this we construct the discrete weak $N$-category on a weak
$n$-category, then show that the $N$-categories so arising are exactly the
$n$-dimensional ones.

We start with the discrete construction in the setting of \emph{strict}%
\index{n-category@$n$-category!weak vs. strict@weak \vs.\ strict}
higher categories, and derive from it the weak version.  Almost all of the
steps involved are thought-free applications of previously-established
theory.  Let
\[
D^N_n: 
\ftrcat{\scat{G}_n^\op}{\Set}
\go 
\ftrcat{\scat{G}_N^\op}{\Set}%
\glo{DNn}
\]
be the functor sending an $n$-globular set $X$ to the $N$-globular
set of~\bref{diag:disc-glob-set}.  This lifts to a functor
\[
\strcat{n} \rIncl \strcat{N},%
\index{n-category@$n$-category!strict!discrete}
\]
and by Lemma~\ref{lemma:lax-map-mnds-is-ftr}, there is a corresponding lax
map of monads
\[
(D^N_n, \psi): 
(\ftrcat{\scat{G}_n^\op}{\Set}, \gm{n})
\go 
(\ftrcat{\scat{G}_N^\op}{\Set}, \gm{N}), 
\]
which, since $D^N_n$ preserves finite limits, induces in turn a functor
\[
(D^N_n)_*:
n\hyph\Operad
\go
N\hyph\Operad.
\]
Explicitly, if $Q$ is an $n$-operad and $\pi$ a $k$-pasting diagram then
\[
((D^N_n)_* Q)(\pi)
=
\left\{
\begin{array}{ll}
Q(\bdry^{k-n} \pi)	&\textrm{if } k \geq n	\\
Q(\pi)			&\textrm{if } k\leq n,
\end{array}
\right.
\]
from which it follows that $(D^N_n)_*$ lifts naturally to a functor
\[
(D^N_n)_*: 
n\hyph\fcat{OC}
\go
N\hyph\fcat{OC}.
\]
Since $L_N$ is initial in $N\hyph\fcat{OC}$, there is a canonical map of
operads
\[
L_N \go (D^N_n)_* L_n,
\]
inducing a functor
\[
\Alg((D^N_n)_* L_n)
\go
\Alg(L_N) = \wkcat{N}.
\]
But we also have from~\ref{sec:change} a functor
\[
\wkcat{n} = \Alg(L_n)
\go
\Alg((D^N_n)_* L_n),
\]
and so obtain a composite functor
\begin{equation}	\label{eq:wk-disc-functor}
D^N_n: \wkcat{n} \go \wkcat{N},%
\glo{DNnweak}
\end{equation}
as required.  A weak $N$-category isomorphic to $D^N_n Y$ for some weak
$n$-category $Y$ will be called a \demph{discrete}%
\index{n-category@$n$-category!weak!discrete}
weak $N$-category on a weak $n$-category.

On underlying globular sets, the discrete functor~\bref{eq:wk-disc-functor}
is merely the original functor $D^N_n$, so~\bref{eq:wk-disc-functor}
restricts to a functor
\begin{equation}	\label{eq:n-dim-eqv}
D^N_n:
\wkcat{n}
\go
(n \textup{-dimensional weak } N \textup{-categories}).
\end{equation}
A sharper statement of Theorem~\ref{thm:n-dim} is that this is an
equivalence of categories.

\begin{prooflike}{Proof of~\ref{thm:n-dim}}
We show that restriction%
\index{restriction}
$R^N_n$ is inverse to the functor $D^N_n$
of~\bref{eq:n-dim-eqv}.  

First, if $X$ is an $n$-dimensional weak $N$-category then
Corollary~\ref{cor:restriction} applies to give $R^N_n X$ the structure of
a weak $n$-category.  For let $\pi\in\pd(n)$ and let $\theta^-$, $\theta^+$
be parallel elements of $L_N(\pi)$.  By contractibility of $L_N$, there
exists $\theta \in L_N(1_\pi)$ satisfying $s(\theta) = \theta^-$ and
$t(\theta) = \theta^+$.  Since any $(n+1)$-cell of $X$ has the same source
and target, we have
\[
\ovln{\theta^-} 
=
\ovln{s(\theta)}
=
s \of \ovln{\theta}
=
t \of \ovln{\theta}
=
\ovln{t(\theta)}
=
\ovln{\theta^+}.
\]
So by~\ref{cor:restriction}, $R^N_n$ induces a functor
\[
R^N_n:
(n \textup{-dimensional weak } N \textup{-categories})
\go
\wkcat{n}.
\]

The composite functor $R^N_n \of D^N_n$ on $\wkcat{n}$ is the identity,
ultimately because the same is true for strict $n$-categories.  Conversely,
let $X$ be an $n$-dimensional weak $N$-category.  We may assume that $X$ is
strictly $n$-dimensional~\bref{diag:disc-glob-set}, which means that there
is an equality $D^N_n R^N_n X = X$ of globular sets.  It remains only to
check that the two $L_N$-algebra structures on $X$ agree; but certainly
they agree in dimensions $n$ and lower, and $n$-dimensionality
of $X$ guarantees that they agree in dimensions higher than $n$ too.
\done
\end{prooflike}

We finish with a method for turning higher-dimensional categories into
lower-dimensional ones.  It is an analogue of the path space construction
in topology (or with slightly different analogies, the loop space
construction or desuspension): given an $n$-category, we forget the
$0$-cells and decrease all the dimensions by $1$ to produce an
$(n-1)$-category.

Again, we start by doing it in the strict%
\index{n-category@$n$-category!weak vs. strict@weak \vs.\ strict}
setting.  A strict $n$-category
is a category enriched in $\strcat{n}$, and a finite-limit-preserving
functor $\cat{V} \go \cat{W}$ between categories with finite limits induces
a finite-limit-preserving functor $\cat{V}\hyph\Cat \go \cat{W}\hyph\Cat$,
so the functor $\Cat \go \Set$ sending a category to its set of arrows
induces a finite-limit-preserving functor
\[
\strcat{n} \go \strcat{(n-1)}
\]
for each positive integer $n$.  This is a lift of the functor
\[
\begin{array}{rrcl}
J_n:	&
\ftrcat{\scat{G}_n^\op}{\Set}		&
\go	&
\ftrcat{\scat{G}_{n-1}^\op}{\Set},	\\%
\glo{Jn}
	&
X	&
\goesto	&
\left(
X(n)
\parpair{s}{t}
\ 
\cdots
\ 
\parpair{s}{t}
X(1)
\right).
\end{array}
\]
By~\ref{lemma:lax-map-mnds-is-ftr}, $J_n$ has the structure of a lax map of
monads
\[
\left(
\ftrcat{\scat{G}_n^\op}{\Set}, \gm{n}
\right)
\go 
\left(
\ftrcat{\scat{G}_{n-1}^\op}{\Set}, \gm{n-1}
\right),
\]
which, since $J_n$ preserves finite limits, induces a functor
\[
(J_n)_*: 
n\hyph\Operad
\go
(n-1)\hyph\Operad.
\]
To describe $(J_n)_*$ explicitly we use the \demph{suspension}%
\index{suspension!globular pasting diagram@of globular pasting diagram}
operator
$\Sigma: \pd(k) \go \pd(k+1)$,%
\glo{suspglobpd}
defined by $\Sigma\pi = (\pi)$.  Here
$\pd(k+1)$ is regarded as the free monoid on $\pd(k)$ and $(\pi)$ is a
sequence of length one.  An example explains the name:
\[
\Sigma
\left(
\gfst{}\gthree{}{}{}{}{}\gfbw{}\gone{}\glst{}
\right)
\ 
=
\ 
\gfst{}\gspecialone{}{}{}{}{}{}{}{}{}\glst{}.
\]
Now, if $P$ is an $n$-operad, $0\leq k\leq n-1$, and $\pi\in\pd(k)$, 
we have
\[
((J_n)_* P)(\pi) = P(\Sigma \pi),
\]
and using the equation $\bdry\Sigma\pi = \Sigma\bdry\pi$, we find that
$(J_n)_*$ lifts naturally to a functor
\[
(J_n)_*: 
n\hyph\fcat{OC} 
\go 
(n-1)\hyph\fcat{OC}.
\]
Just as for the discrete construction, this induces a functor
\[
J_n: \wkcat{n} \go \wkcat{(n-1)}
\]
whose effect on the underlying globular sets is the original `shift'
functor $J_n$.

The $(n-1)$-category $J_n X$ arising from an $n$-category $X$ is called the
\demph{localization}%
\index{localization}
of $X$.  Its underlying $(n-1)$-globular set is the
disjoint union over all $a, b \in X(0)$ of the $(n-1)$-globular sets
$X(a,b)$ defined in the proof of~\ref{propn:str-n-cats-comparison}.  Since
the functor $\gm{n-1}$ preserves
coproducts~(\ref{thm:n-forgetful-properties}),%
\index{coproduct!preserved by monad}%
\index{monad!coproduct-preserving}
it follows from the lemma
below that the $(n-1)$-category structure on $J_n X$ amounts to an
$(n-1)$-category structure on each $X(a,b)$, in both the strict and the
weak settings.  So localization defines functors
\begin{eqnarray*}
\strcat{n}      &\go        &(\strcat{(n-1)})\hyph\Gph,        \\
\wkcat{n}       &\go        &(\wkcat{(n-1)})\hyph\Gph.         
\end{eqnarray*}
In the strict case, an $n$-category is a graph of $(n-1)$-categories
together with composition functors obeying simple laws; this is just
ordinary enrichment.%
\index{enrichment!define n-category@to define $n$-category}
 In the weak case it is much more difficult to say
what extra structure is needed.

\begin{lemma}  \label{lemma:coprod-algs}
Let $S$ be a cartesian monad on a presheaf category $\cat{E}$, such that
the functor part of $S$ preserves coproducts.  Let $P$ be an $S$-operad
and let $(X_i)_{i\in I}$ be a family of objects of $\cat{E}$.  Then a
$P$-algebra structure on $\coprod X_i$ amounts to a $P$-algebra structure
on each $X_i$.  
\end{lemma}
\begin{proof}
It is enough to show that the functor $S_P$ preserves coproducts.  Since
$S_P$ is defined using $S$ and pullback, and pullbacks interact well with
coproducts in presheaf categories, this follows from the same property of
$S$.
\done
\end{proof}

Finally, localization%
\index{localization}
works just as well for $\omega$-categories.  The
localization functors $\strcat{n} \go \strcat{(n-1)}$ induce in the limit
an endofunctor of $\strcat{\omega}$, which on underlying globular sets is
the endofunctor $J$ of $\ftrcat{\scat{G}}{\Set}$ forgetting $0$-cells.  So
exactly as in the finite-dimensional case, a weak $\omega$-category $X$
gives rise to a family $(X(a,b))_{a, b \in X(0)}$ of weak
$\omega$-categories.%
\index{n-category@$n$-category!weak!varying n@for varying $n$|)}

\section{Weak $2$-categories}
\lbl{sec:wk-2}%
\index{bicategory!unbiased|(}

A polite person proposing a definition of weak $n$-category should explain
what happens when $n=2$.  With our definition, $\wkcat{2}$ turns out to be
equivalent to $\UBistr$, the category of small unbiased bicategories and
unbiased strict functors.

Observe that since the maps in $\wkcat{2}$ are \emph{strict} functors, we
obtain an equivalence with $\UBistr$, not $\UBiwk$ or $\UBilax$; and unlike
its weak and lax siblings, $\UBistr$ is not equivalent to the analogous
category of classical%
\index{bicategory!unbiased vs. classical@unbiased \vs.\ classical}
bicategories (or at least, the obvious functor is not
an equivalence).  So we cannot conclude that $\wkcat{2}$ is equivalent to
$\Bistr$.  Nevertheless, the results of~\ref{sec:notions-bicat} mean that
it is fair to regard classical bicategories as essentially the same as
unbiased bicategories, and therefore, by the results below, essentially the
same as weak 2-categories.  One would expect that if the definition of weak
functor between $n$-categories were in place, $\Biwk$ would be equivalent
to the category of weak 2-categories and weak functors.

\begin{thm}	\lbl{thm:wk-2}
There are equivalences of categories
\begin{eqnarray*}
\wkcat{0}	&\eqv	&\Set,	\\
\wkcat{1}	&\eqv	&\Cat,	\\
\wkcat{2}	&\eqv	&\UBistr.
\end{eqnarray*}
\end{thm}

So far we have ignored weak $0$-categories;%
\index{zero-category@0-category}
indeed, we have not even proved
that there is an initial 0-operad-with-contraction.  A $0$-globular set is
a set and the monad $\gm{0}$ is the identity, so a $0$-operad is a monoid
and an algebra for a $0$-operad is a set acted on by the corresponding
monoid.  There is a unique precontraction on every $0$-operad, which is a
contraction just when the corresponding monoid has cardinality $1$.  So
$0\hyph\fcat{OC}$ is the category of one-element monoids, any object $L_0$
of which is initial, and
\[
\wkcat{0} \eqv \Set.
\]

A weak $1$-category%
\index{one-category@1-category}
is a $1$-dimensional weak $2$-category
(Theorem~\ref{thm:n-dim}), so the middle equivalence of
Theorem~\ref{thm:wk-2} will follow from the last.  It is, however, easy
enough to prove directly.  We have just seen that a
$0$-operad-with-precontraction is a monoid, so the initial such is also the
terminal such.  The equivalence $1\hyph\fcat{OC} \eqv 0\hyph\fcat{OP}$
of~\ref{propn:OP-OC} then tells us that the initial
$1$-operad-with-contraction is the terminal $1$-operad.  But algebras for
the terminal $\gm{1}$-operad are just $\gm{1}$-algebras, so
\[
\wkcat{1} \eqv \Cat.
\]

It is not prohibitively difficult to prove the 2-dimensional equivalence
result explicitly, as in Leinster~\cite[4.8]{OHDCT}; 2-operads are just
about manageable.  Here we use an abstract method instead, taking advantage
of some earlier calculations.

Notation: we write 
\begin{itemize}
\item $W$%
\glo{Wfreemon}
for both the free monoid monad on $\Set$ and the free strict
monoidal category monad on $\Cat$
\item $\cat{V}\hyph\Gph$ for the category of graphs%
\index{graph!enriched}
enriched in a given
finite product category $\cat{V}$~(\ref{defn:V-gph}), and $\Gph$%
\glo{abbrevGph}
for
$\Set\hyph\Gph$
\item $\Sigma: \cat{V} \go \cat{V}\hyph\Gph$%
\glo{suspobjgraph}%
\index{suspension|(}
for the functor sending an
object $V$ of $\cat{V}$ to the one-object $\cat{V}$-graph whose single
hom-set is $V$
\item $\fc_\cat{V}$%
\glo{fcV}
for the free%
\index{category!free enriched}
$\cat{V}$-enriched category monad on
$\cat{V}\hyph\Gph$ (when it exists).
\end{itemize}

A 1-globular set is a directed graph and $\gm{1}$ is the free category monad
$\fc$ of Chapter~\ref{ch:fcm}, so a $1$-operad is an $\fc$-operad,%
\index{fc-operad@$\fc$-operad}
which is
an $\fc$-multicategory with only one 0-cell and one horizontal 1-cell.  A
precontraction $\kappa$ on an $\fc$-operad assigns to each $r\in\nat$ and
pair $(f,f')$ of vertical 1-cells a 2-cell
\[
\begin{fcdiagram}
\bullet	&\rTo		&\bullet&\rTo		&\ 	&\cdots	
&\ 	&\rTo		&\bullet	\\
\dTo<{f}&		&	&		&\Downarrow\kappa_r(f,f')&
&	&		&\dTo>{f'}	\\
\bullet	&		&	&		&\rTo	&	
&	&		&\bullet	\\
\end{fcdiagram}
\]
with $r$ 1-cells along the top.  Recall from~\ref{eg:fcm-cl-opd} that there
is an embedding
\[
\Operad \rIncl \fc\hyph\Operad
\]
identifying plain operads with $\fc$-operads having only one 0-cell, one
vertical 1-cell and one horizontal 1-cell.  There we called the embedding
$\Sigma$; here we call it $\Sigma_*$, because it is induced by the weak map
of monads $\Sigma: (\Set, W) \go (\Gph, \fc)$.  If $P$ is a plain operad
then a precontraction on $\Sigma_* P$ consists of an element $\kappa_r \in
P(r)$ for each $r\in\nat$.  Take the plain operad $\tr$ of trees%
\index{tree!operad of|(}
and the
$r$-leafed corolla $\nu_r\in\tr(r)$ for each
$r\in\nat$~(\ref{eg:opd-of-trees}): then using the fact that $\tr$ is the
free operad containing an operation of each arity, it is easy to show that
the corresponding $1$-operad-with-precontraction $\Sigma_* \tr$ is initial.
So $L_2 = (I^2_1)_* \Sigma_* \tr$, and we have
\begin{equation}	\label{eq:wk-2-L2}
\wkcat{2} \iso \Alg((I^2_1)_* \Sigma_* \tr).
\end{equation}

On the other hand, we saw earlier that the theory of unbiased bicategories
is also generated by the operad of trees.  Specifically, we showed
in~\ref{sec:notions-bicat} that
\[
\UBistr 
\iso 
1\hyph\Bistr
=
\fcat{CatAlg}_\mr{str} I_* \tr.%
\index{Sigma-bicategory@$\Sigma$-bicategory}%
\index{bicategory!Sigma-@$\Sigma$-}
\]
The functors $I_*$ and $\fcat{CatAlg}_\mr{str}$ can be described as
follows.  We have maps of monads
\[
(\Set, W)
\goby{I}
(\Cat, W)
\goby{\Sigma}
(\Cat\hyph\Gph, \fc_\Cat)
\]
where $I$ is the indiscrete%
\index{category!indiscrete}
category functor (p.~\pageref{p:indiscrete});
$I$ is lax and $\Sigma$ is weak.  Recalling from~\ref{eg:mti-Cat} that a
$(\Cat, W)$-operad is a $\Cat$-operad, we obtain induced functors
\[
\Operad 
\goby{I_*} 
\Cat\hyph\Operad
\goby{\Sigma_*}
\fc_\Cat \hyph\Operad
\goby{\Alg}
\CAT^\op.
\]
The functor $I_*$ is the same as the one used in Chapter~\ref{ch:monoidal},
and the composite of the last two functors is $\fcat{CatAlg}_\mr{str}$, so
\[
\UBistr \iso \Alg(\Sigma_* I_* \tr).%
\index{tree!operad of|)}
\]
Comparing with~\bref{eq:wk-2-L2}, it is enough to prove
\begin{lemma}
For any plain operad $P$, there is an isomorphism of categories
\[
\Alg((I^2_1)_* \Sigma_* P)
\iso
\Alg(\Sigma_* I_* P).
\]
\end{lemma}
\begin{proof}
We have three monads on the category $\Gph\hyph\Gph$ of 2-globular sets:
first, $\fc\hyph\Gph$, the result of applying the 2-functor
$\blank\hyph\Gph: \CAT \go \CAT$ to the monad $\fc$ on $\Gph$; second,
$\fc_\Gph$, the free $\Gph$-enriched category monad; third, $\gm{2}$.  We
show in the proof of~\ref{propn:free-enr} that $\gm{2}$ is the result of
gluing $\fc\hyph\Gph$ to $\fc_\Gph$ by a distributive%
\index{distributive law}
law
\[
\lambda: 
(\fc\hyph\Gph) \of \fc_\Gph
\go
\fc_\Gph \of (\fc\hyph\Gph).
\]
We also saw in Lemma~\ref{lemma:distrib-corr} that a distributive law gives
rise to a monad `$\twid{S}$', which in this case is the monad $\fc_\Cat$ on
$\Cat\hyph\Gph$, and in Lemma~\ref{lemma:distrib-iso-algs} that it gives
rise to a lax map of monads, 
% `$(U, \psi)$', 
which in this case is of the form
\[
(U, \psi): 
(\Cat\hyph\Gph, \fc_\Cat) 
\go 
(\Gph\hyph\Gph, \gm{2})
\]
where $U$ is the forgetful functor.  

It is straightforward to check that there is an equality of natural
transformations
\[
\begin{diagram}[size=2em]
\Set		&\rTo^W		&\Set		\\
\dTo<\Sigma	&\neeq		&\dTo>\Sigma	\\
\Gph		&\rTo^{\gm{1}}	&\Gph		\\
\dTo<{I^2_1}	&\nent		&\dTo>{I^2_1}	\\
\Gph\hyph\Gph	&\rTo_{\gm{2}}	&\Gph\hyph\Gph	\\
\end{diagram}
\diagspace
=
\diagspace
\begin{diagram}[size=2em]
\Set		&\rTo^W		&\Set		\\
\dTo<I		&\nent		&\dTo>I		\\
\Cat		&\rTo^W		&\Cat		\\
\dTo<\Sigma	&\neeq		&\dTo>\Sigma	\\
\Cat\hyph\Gph	&\rTo^{\fc_\Cat}&\Cat\hyph\Gph	\\
\dTo<U		&\nent\psi	&\dTo>U		\\
\Gph\hyph\Gph	&\rTo_{\gm{2}}	&\Gph\hyph\Gph	\\
\end{diagram}
\]
where the unmarked transformations are the ones referred to above.  So
\[
(I^2_1)_* \Sigma_* = U_* \Sigma_* I_*:
\Operad \go \gm{2}\hyph\Operad,
\]%
\index{suspension|)}%
and it is enough to prove that for any $\fc_\Cat$-operad $Q$,
\[
\Alg(Q) \iso \Alg(U_* Q).
\]
This follows immediately from Proposition~\ref{propn:shape-distrib}.  
\done
\end{proof}
That completes the proof of Theorem~\ref{thm:wk-2}.%
\index{bicategory!unbiased|)}

In~\ref{sec:contr} we said what we wanted the operad $L$ to look like
in low dimensions: if $\blob$ is the unique $0$-pasting diagram then
$L(\blob)$%
\lbl{p:L-blob}
should be a one-element set, and if $\chi_k$ is the $1$-pasting diagram made
up of $k$ arrows then $L(\chi_k)$ should be $\tr(k)$.  We now know that our
wishes were met: for by~\ref{propn:change-of-dimension}, the 2-dimensional
restriction $(R^\omega_2)_* L$ is the initial 2-operad-with-precontraction,
and we have shown this to be $\Sigma_* \tr$.

\index{three-category@3-category!definitions of|(}%
\index{tricategory!non-algebraic nature of|(}
What about $3$-categories?  It should be possible to write down an explicit
definition of unbiased tricategory (similar to that of unbiased
bicategory,~\ref{defn:lax-bicat}) and to prove that the category of
unbiased tricategories and strict maps is equivalent to $\wkcat{3}$.  This
would be a lot of work, and it is not clear that the result would have any
advantage over the abstract definition of weak $3$-category.  

We could also try to compare our weak $3$-categories with the tricategories 
of Gordon,%
\index{Gordon, Robert}
Power and Street~\cite{GPS}.  The analogous comparison one dimension down,
between weak $2$-categories and classical bicategories, is already
difficult because we do not have a notion of weak functor between
$n$-categories (see the beginning of this section).  A further difficulty
is that Gordon, Power and Street's definition is not quite algebraic;%
\index{algebraic theory!tricategories are not}
put
another way, the forgetful functor
\[
(\textrm{tricategories } + \textrm{ strict maps})
\go
\ftrcat{\scat{G}^\op}{\Set}
\]
seems highly unlikely to be monadic.  If true, this means that there can be
no 3-operad whose algebras are precisely tricategories (in apparent
contradiction to Batanin~\cite[p.~94]{BatMGC}).%
\index{Batanin, Michael!globular operads@on globular operads}

The reason why the theory of tricategories is not quite algebraic is as
follows.  Most of the definition of tricategory consists of some data
subject to some equations, but a small part does not: in items~(TD5)
and~(TD6), it is stipulated that certain transformations of bicategories
are equivalences.  This is not an algebraic axiom, as there are many
different choices of weak inverses and none has been specified.  Compare
the fact that the forgetful functor from non-empty sets to $\Set$ is not
monadic%
\index{monadic adjunction}
(indeed, has no left adjoint), in contrast to the forgetful functor
from pointed sets to $\Set$.  To make the definition algebraic we would
have to add in as data a weak inverse for each of these equivalences,
together with two invertible modifications witnessing that it is a weak
inverse, and then add more coherence axioms (saying, among other things,
that this data forms an \emph{adjoint} equivalence).  The result would be
an even more complicated, but conceptually pure, notion of tricategory.%
\index{three-category@3-category!definitions of|)}%
\index{tricategory!non-algebraic nature of|)}

\begin{notes}

Contractions were introduced in my~\cite{SHDCT}.  When I wrote and made
public the first version of that paper I believed that I was explaining
Batanin's notion of contraction,%
\index{contraction!notions of}
but in fact I was inventing a new one;
see~\ref{sec:alg-defns-n-cat} for an explanation of the difference.  

Some of the results here on $n$- and 2-categories appeared in
my~\cite{OHDCT}.  

See Crans~\cite[2.3]{CraTPG}%
\index{Crans, Sjoerd}
for a completely elementary definition of
Gray-category~(\ref{eg:Gray}).%
\index{Gray-category}
 I learned that the operad for
Gray-categories can be defined from the operad for sesquicategories from
Batanin~\cite[p.~94]{BatMGC}.%
\index{Batanin, Michael!globular operads@on globular operads}

\end{notes}

\chapter{Other Definitions of Weak $n$-Category}
\lbl{ch:other-defns}

\chapterquote{%
Zounds!  I was never so bethump'd with words!}{%
Shakespeare~\cite{Sha}}

\index{n-category@$n$-category!definitions of|(}
\noindent
The definition of weak $n$-category studied in the previous chapter is, of
course, just one of a host of proposed definitions.  Ten of them were
described in my~\cite{SDN} survey, all except one in formal, precise terms.
However, the format of that paper did not allow for serious discussion of
the interrelationships, and one might get the impression from it that the
ten definitions embodied eight or so completely different approaches to the
subject.

I hope to correct that impression here.  Fundamentally, there seem to be
only two approaches.

In the first, a weak $n$-category is regarded as a presheaf%
\index{presheaf!structure@with structure}
with
\emph{structure}.%
\index{structure vs. properties@structure \vs.\ properties}%
\index{properties vs. structure@properties \vs.\ structure}
Usually `presheaf' means $n$-globular set,
and `structure' means $S$-algebra structure for some monad $S$, often
coming from a globular operad.  The definition studied in the previous
chapter is of this type.  

In the second approach, a weak $n$-category is regarded as a presheaf with
\emph{properties}.%
\index{presheaf!properties@with properties}
 There is no hope that a weak $n$-category could be
defined as an $n$-globular set with properties, so the category
$\scat{C}_n$ on which we are taking presheaves must be larger than
$\scat{G}_n$; presheaves on $\scat{C}_n$ must somehow have composition
built in.  The case $n=1$ makes this clear.  In the first approach, a
category is defined as a directed graph (presheaf on $\scat{G}_1$) with
structure.  In the second, a category cannot be defined a category as a
presheaf-with-properties on $\scat{G}_1$, but it can be defined as a
presheaf-with-properties on the larger category $\Delta$:%
\index{simplex category $\Delta$}
this is the
standard characterization of a category by its
nerve~(p.~\pageref{p:defn-nerve}).%
\index{nerve!category@of category}

There are other descriptions of the difference between the two approaches.
In the sense of the introduction to Chapter~\ref{ch:monoidal}, the first is
algebraic%
\index{algebraic theory!n-categories@of $n$-categories}
(the various types of composition in a weak $n$-category are
\emph{bona fide} operations) and the second is non-algebraic (composites
are not determined uniquely, only up to equivalence).  Or, the first
approach can be summarized as `take the theory of strict $n$-categories,
weaken%
\index{weakening!theory of n-categories@for theory of $n$-categories}
it, then take models for the weakened theory', and the second as
`take weak models for the theory of strict $n$-categories'.%
\index{n-category@$n$-category!definitions of|)}

The definitions following the two approaches are discussed
in~\ref{sec:alg-defns-n-cat} and~\ref{sec:non-alg-defns-n-cat}
respectively.  What I have chosen to say about each definition (and which
definitions I have chosen to say anything about at all) is dictated by how
much I feel capable of saying in a simple and not too technical way; the
emphasis is therefore rather uneven.  In particular, there is more on the
definitions close to that of the previous chapter than on those further
away.

\section{Algebraic definitions}
\lbl{sec:alg-defns-n-cat}%
\index{algebraic theory!n-categories@of $n$-categories|(}

Here we discuss the definitions proposed by Batanin, Penon, Trimble, and
May.  We continue to use the notation of the previous chapter: $T$ is the
free strict $\omega$-category monad, $\pd = T1$, and so on.

\minihead{Batanin's definition}%
\index{Batanin, Michael!definition of n-category@definition of $n$-category|(}%
\index{n-category@$n$-category!definitions of!Batanin's|(}

The definition of weak $\omega$-category studied in the previous chapter is
a simplification of Batanin's~\cite{BatMGC} definition.  There are two main
differences.  The less significant is bias: where our definition treats
composition of all shapes equally, Batanin's gives special status to binary
and nullary compositions.  For instance, our weak $2$-categories are
unbiased bicategories, but his are classical, biased, bicategories.
The more significant difference is conceptual.  We integrated composition%
\index{composition!coherence@\vs.\ coherence}%
\index{coherence!composition@\vs.\ composition}
and coherence into the single notion of contraction; Batanin keeps the two
separate.  This makes his definition more complicated to state, but more
obvious from the traditional point of view.

Composition is handled as follows.  The unit map $\eta_1: 1 \go T1 = \pd$
picks out, for each $m\in\nat$, the $m$-pasting diagram $\iota_m$ looking
like a single $m$-cell.  Define 
a collection $\fcat{binpd} \rIncl \pd$%
\index{pasting diagram!globular!binary}
by
\[
\fcat{binpd}(m)
=
\{ \iota_m \of_0 \iota_m,\,
\iota_m \of_1 \iota_m,\,
\ldots,\,
\iota_m \of_m \iota_m \}
\subseteq \pd(m)
\]
where $\of_p$ is composition in the strict $\omega$-category $\pd$ and
$\iota_m \of_m \iota_m$ means $\iota_m$; then $\fcat{binpd}$
consists of binary and unary diagrams such as
\[
\begin{array}{rclcrcl}
\iota_1 \of_0 \iota_1 &= &
\gfstsu\gonesu\gzersu\gonesu\glstsu,	&&
\iota_1 &=&
\gfstsu\gonesu\glstsu,	\\
\iota_2 \of_0 \iota_2 &= & 
\gfstsu\gtwosu\gzersu\gtwosu\glstsu,	&&	
\iota_2 \of_1 \iota_2 &= &
\gfstsu\gthreesu\glstsu.
\end{array}
\]
\begin{defn}
Let $P$ be a globular operad.  A \demph{system%
\index{system of compositions}%
\index{compositions, system of}
of (binary) compositions}
in $P$ is a map of collections $\fcat{binpd} \go P$, written
\[
\left(\iota_m \of_p \iota_m	\in \fcat{binpd}(m) \right)
\ \goesto\ 
\left(\delta^m_p \in P(\iota_m \of_p \iota_m)\right),%
\glo{deltasys}
\]
such that $\delta^m_m$ is the identity operation $1_m \in P(\iota_m)$ for
each $m\in\nat$.  
\end{defn}
So, for instance, $s(\delta^2_0) = t(\delta^2_0) = \delta^1_0$ and
$s(\delta^2_1) = t(\delta^2_1) = 1_1$.  (If we wanted to do an unbiased
version of Batanin's definition, we could replace $\fcat{binpd}$ by
$\pd$.)  

\begin{example}	\lbl{eg:contr-gives-sys}
A contraction $\kappa$ on an operad $P$ canonically determines a system
$\delta^\blob_\blob$ of compositions, defined inductively by
\[
\delta^m_p =
\left\{
\begin{array}{ll}
\kappa_{\iota_m \sof_p \iota_m} (\delta^{m-1}_p, \delta^{m-1}_p)	&
\textrm{if }
p < m,	\\
1_m	&
\textrm{if }
p = m.
\end{array}
\right.
\]
\end{example}

\begin{figure}
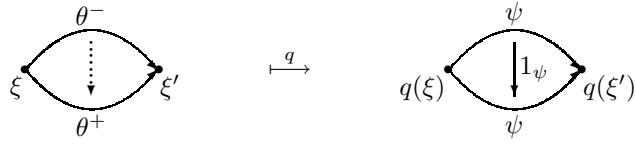

\[
\gfst{\xi}\gtwodotty{\theta^-}{\theta^+}{}%
\glst{\xi'}
\mbox{\hspace{2.5em}}
\stackrel{q}{\goesto}
\mbox{\hspace{2.5em}}
\gfst{q(\xi)}\gtwo{\psi}{\psi}{1_\psi}\glst{q(\xi')}
\]
\caption{Effect of a coherence $\zeta$, shown for $k=1$.  The dotted arrow
  is $\zeta(\theta^-, \theta^+)$, and $\psi = q(\theta^-) = q(\theta^+)$}
\label{fig:coh-on-map}
\end{figure}

To describe coherence we use the notion of a \demph{reflexive%
\index{reflexive globular set}%
\index{globular set!reflexive}
globular
set}, that is, a globular set $Y$ together with functions
\[
\cdots \ 
\lTo^i	Y(k+1) 
\lTo^i	Y(k)
\lTo^i
\ \cdots \ 
\lTo^i	Y(0),
\]%
\glo{irefl}%
written $i(\psi) = 1_\psi$, such that $s(1_\psi) = t(1_\psi) = \psi$ for
each $k\geq 0$ and $\psi\in Y(k)$.  Reflexive globular sets form a presheaf
category $\ftrcat{\scat{R}^\op}{\Set}$.%
\glo{Rrefl}
 The underlying globular set of a
strict $\omega$-category is canonically reflexive, taking the identity
cells $1_\psi$.  This applies in particular to $T1 = \pd$; if $\pi \in
\pd(k)$ then $1_\pi \in \pd(k+1)$ is the degenerate $(k+1)$-pasting
diagram represented by the same picture as $\pi$.
\begin{defn}
Let $X$ be a globular set, $Y$ a reflexive globular set, and $q: X \go Y$ a
map of globular sets.  For $k\geq 0$ and $\psi\in Y(k)$, write
(Fig.~\ref{fig:coh-on-map}) 
\begin{eqnarray*}
\mr{Par}'_q(k)	&
=	&
\{ (\theta^-, \theta^+) \in X(k) \times X(k) \such	
\theta^- \textrm{ and } \theta^+ \textrm{ are parallel, }
\\  &	&
q(\theta^-) = q(\theta^+) \}.
\end{eqnarray*}
A \demph{coherence}%
\index{coherence!map of globular sets@on map of globular sets}
$\zeta$ on $q$ is a family of functions
\[
\left(
\mr{Par}'_q(k) \goby{\zeta_k} X(k+1)
\right)_{k\geq 0}
\]
such that for all $k\geq 0$ and $(\theta^-, \theta^+) \in \mr{Par}'_q(k)$,
writing $\zeta_k$ as $\zeta$,
\[
s(\zeta(\theta^-, \theta^+)) = \theta^-,
\ \ 
t(\zeta(\theta^-, \theta^+)) = \theta^+,
\ \ 
q(\zeta(\theta^-, \theta^+)) = 1_{q(\theta^-)} (= 1_{q(\theta^+)}).
\]
\end{defn}

\begin{example}
Any contraction on a map canonically determines a coherence, as is clear
from a comparison of Figs.~\ref{fig:contr-on-map}
(p.~\pageref{fig:contr-on-map}) and~\ref{fig:coh-on-map}.  Formally,
$\mr{Par}'_q(k) = \coprod_{\psi\in Y(k)} \mr{Par}_q(1_\psi)$ and a
contraction $\kappa$ determines the coherence $\zeta$ given by
$\zeta(\theta^-, \theta^+) = \kappa_{1_\psi}(\theta^-, \theta^+)$ where
$\psi = q(\theta^-) = q(\theta^+)$.  The class of contractible maps is
closed under composition, but the class of maps admitting a coherence is
not.
\end{example}

\begin{defn}	\lbl{defn:coh-coll}
A \demph{coherence}%
\index{coherence!collection@on collection}
on a collection $(P \goby{d} T1)$ is a coherence on the
map $d$, and a \demph{coherence}%
\index{coherence!globular operad@on globular operad}
on an operad is a coherence on its
underlying collection.  A map, collection, or operad is \demph{coherent}%
\index{coherent!map}%
\index{coherent!globular operad}
if
it admits a coherence.  
\end{defn}
Explicitly, a coherence on an operad $P$ assigns to each $\pi\in\pd(k)$ and
parallel pair $\theta^-, \theta^+ \in P(\pi)$ an element $\zeta(\theta^-,
\theta^+)$ of $P(1_\pi)$ with source $\theta^-$ and target $\theta^+$.  Let
$X$ be a $P$-algebra; then since $\rep{\pi} = \rep{1_\pi}$
(p.~\pageref{p:degen-rep}), we have $(TX)(\pi) = (TX)(1_\pi)$, so if
$\mathbf{x} \in (TX)(\pi)$ then there is a $(k+1)$-cell
\[
\ovln{\zeta(\theta^-, \theta^+)} (\mathbf{x}):
\ovln{\theta^-}(\mathbf{x})
\go 
\ovln{\theta^+}(\mathbf{x})
\]
in $X$ connecting the two `composites' $\ovln{\theta^-}(\mathbf{x})$ and
$\ovln{\theta^+}(\mathbf{x})$ of $\mathbf{x}$.  Taking $\pi = \iota_k$ and
$\theta^- = \theta^+ = 1_k$, this provides in particular a reflexive
structure on the underlying globular set of $X$.
\begin{example}	\lbl{eg:coh-vs-contr}
Any contractible operad is coherent (by the previous example), but not
conversely.  For instance, there is an operad $R$ whose algebras are
reflexive globular sets; it is uniquely determined by
\[
R(\pi) =
\left\{
\begin{array}{ll}
1		&
\textrm{if } \pi \in 
\{ \iota_k, 1_{\iota_{k-1}}, 1_{1_{\iota_{k-2}}}, \ldots \}
\\
\emptyset	&\textrm{otherwise}
\end{array}
\right.
\]
($\pi\in\pd(k)$).  This operad is coherent (trivially) but not contractible
(since some of the sets $R(\pi)$ are empty).  So a given globular set $X$
is an algebra for some coherent operad if and only if it admits a reflexive
structure; on the other hand, by~\ref{eg:wk-omega-cat-contr-opd}, it is an
algebra for some contractible operad if and only if it admits a weak
$\omega$-category structure.
\end{example}

\index{contraction!notions of|(}%
A coherence is what Batanin%
\index{contraction!Batanin's sense@in Batanin's sense}
calls a contraction.  As we have seen, our
contractions are more powerful, providing both a coherence and a system of
compositions.  

Rather confusingly, at several points in Batanin~\cite{BatMGC} the word
`contractible' is used as an abbreviation for `contractible [coherent] and
admitting a system of compositions'.  In particular, weak $\omega$- and
$n$-categories are often referred to as algebras for a `universal
contractible operad'.  This is not meant literally: `universal' means
weakly%
\index{weakly initial}
initial (in other words, there is at least one map from it to any
other contractible operad), and the operad $R$ of~\ref{eg:coh-vs-contr} is
in Batanin's terminology contractible, so any genuine `universal
contractible operad' $P$ satisfies $P(\pi) = \emptyset$ for almost all
pasting diagrams $\pi$.  A $P$-algebra is then nothing like an
$\omega$-category.  Indeed, $R$ is in Batanin's terminology the initial
operad equipped with a contraction, and its algebras are mere reflexive
globular sets.  The system of compositions is a vital ingredient; left out,
the situation degenerates almost entirely.

I believe it is the case that given a system of compositions and a
coherence on an operad, a contraction can be built.  This is the
non-canonical converse to the canonical process in the other direction; the
situation is like that of biased \vs.\ unbiased bicategories.  So it
appears that despite an abuse of terminology and two different definitions
of contractibility, the term `contractible operad' means exactly the same
in Batanin's work as here.%
\index{contraction!notions of|)}

The two ingredients---composition and coherence---can be combined to make a
definition of weak $\omega$-category in several possible ways:
\begin{itemize}
\item
Imitate the definition of the previous chapter.  In other words, take the
category of globular operads equipped with both a system of compositions
and a coherence, prove that it has an initial object $(B,%
\glo{BBatopd}
\delta^\blob_\blob, \kappa)$, and define a weak $\omega$-category as a
$B$-algebra.  This is Definition \textbf{B1} in my~\cite{SDN} survey.
\item
Use Batanin's operad $K$, constructed in his Theorem~8.1.  He proves that
$K$ is weakly initial in the full subcategory of $T\hyph\Operad$ consisting
of the coherent operads admitting a system of compositions.  Weak
initiality does not characterize $K$ up to isomorphism, so one needs some
further information about $K$ in order to use this definition.  It seems to
be claimed that $K$ is initial in the category of operads equipped with a
system of compositions and a coherence (Remark~2 after the proof of
Theorem~8.1), but it does not seem obvious that this claim is true,
essentially because of the set-theoretic complement taken in the proof of
Lemma~8.1.
\item 
Define a weak $\omega$-category as a pair $(P, X)$ where $P$ is a coherent
operad admitting a system of compositions and $X$ is a $P$-algebra.  Given
such a pair $(P, X)$, we can choose a system of compositions and a
coherence on $P$, and this turns $X$ into a $B$-algebra---that is, a weak
$\omega$-category in the sense of the first method.  The present method has
some variants: we might insist that $P(\blob) = 1$, where $\blob$ is the
unique $0$-pasting diagram, or we might drop the condition that $P$ admits
a system of compositions and replace it with the more relaxed requirement
that $P(\pi) \neq \emptyset$ for all pasting diagrams $\pi$ (giving
Batanin's Definition~8.6 of `weak $\omega$-categorical object').
\end{itemize}

Batanin's weak $\omega$-categories can be compared with the weak
$\omega$-categories of the previous chapter.  We have already shown that a
contraction on a globular operad gives rise canonically to a system of
compositions and a coherence
(\ref{eg:contr-gives-sys},~\ref{eg:coh-vs-contr}).  This is true in
particular of the operad $L$, so there is a canonical map $B \go L$ of
operads, inducing in turn a canonical functor from $L$-algebras to
$B$-algebras.  Conversely, $B$ is non-canonically contractible, so there is
a non-canonical functor in the other direction.  In the case $n=2$, this
is the comparison of biased and unbiased bicategories.%
\index{Batanin, Michael!definition of n-category@definition of $n$-category|)}%
\index{n-category@$n$-category!definitions of!Batanin's|)}

\minihead{Penon's definition}%
\index{Penon, Jacques!definition of n-category@definition of $n$-category|(}%
\index{n-category@$n$-category!definitions of!Penon's|(}

The definition of weak $\omega$-category proposed by Penon~\cite{Pen} does
not use the language of operads, but is nevertheless close in spirit to the
definition of Batanin.  It can be stated very quickly.

\begin{defn}
Let $q: X \go Y$ be a map of reflexive globular sets.  A coherence $\zeta$
on $q$ is \demph{normal}%
\index{normal coherence}%
\index{coherence!normal}
if $\zeta(\theta, \theta) = 1_\theta$ for all
$k\geq 0$ and $\theta\in X(k)$.
\end{defn}
Penon calls a normal coherence an \'etirement,%
\index{etirement@\'etirement}
or stretching.%
\index{stretching}
This might
seem to conflict with the contraction terminology, but it is only a matter
of viewpoint: $X$ is being shrunk, $Y$ stretched.

\begin{defn}
An \demph{$\omega$-magma}%
\index{magma}%
\index{omega-magma@$\omega$-magma}
is a reflexive globular set $X$ equipped with a
binary composition function $\of_p: X(m) \times_p X(m) \go X(m)$ for each
$m > p \geq 0$, satisfying the source and target axioms
of~\ref{defn:strict-n-cat-glob}\bref{item:s-t-comp}.
\end{defn}
An $\omega$-magma is a very wild structure, and a strict $\omega$-category
very tame; weak $\omega$-categories are somewhere in between.  

Let $\cat{Q}$ be the category whose objects are quadruples $(X, Y, q,
\zeta)$ with $X$ an $\omega$-magma, $Y$ a strict $\omega$-category, $q: X
\go Y$ a map of $\omega$-magmas, and $\zeta$ a normal coherence on $q$.
Maps
\[
(X, Y, q, \zeta) \go (X', Y', q', \zeta')
\]
in $\cat{Q}$ are pairs $(X \go X', Y \go Y')$ of maps commuting with all
the structure present.  There is a forgetful functor $U$ from $\cat{Q}$ to
the category $\ftrcat{\scat{R}^\op}{\Set}$ of reflexive globular sets,
sending $(X, Y, q, \zeta)$ to the underlying reflexive globular set of $X$.
This has a left adjoint $F$, and a weak $\omega$-category is
defined as an algebra for the induced monad $U\of F$ on
$\ftrcat{\scat{R}^\op}{\Set}$.

Observe that $U$ takes the underlying reflexive globular set of $X$, not of
$Y$.  The object $(X, Y, q, \zeta)$ of $\cat{Q}$ should therefore be
regarded as $X$ (not $Y$) equipped with extra structure, making it perhaps
more apt to think of an object of $\cat{Q}$ as a shrinking rather than a
stretching.

Recent work of Cheng%
\index{Cheng, Eugenia}
relates Penon's definition to the definition of the
previous chapter through a series of intermediate definitions; older work
of Batanin~\cite{BatPMW}%
\index{Batanin, Michael}%
\index{Batanin, Michael!definition of n-category@definition of $n$-category}
relates Penon's definition to his own.  One point
can be explained immediately.  Take the operad $B$ equipped with its system
of compositions $\delta^\blob_\blob$ and its coherence $\kappa$.  There is
a unique reflexive structure on the underlying globular set of $B$ for
which $\kappa$ is normal, namely $1_\theta = \zeta(\theta, \theta)$.  Also,
the system of compositions in the operad $B$ makes its underlying globular
set into an $\omega$-magma: given $0\leq p < m$ and $(\theta_1, \theta_2)
\in B(m) \times_{B(p)} B(m)$, put
\[
\theta_1 \ofdim{p} \theta_2 
=
\delta^m_p \of (\theta_1, \theta_2)
\]
where the $\of$ on the right-hand side is operadic composition.  This gives
$B \go T1$ the structure of an object of $\cat{Q}$.  Now, if $X$ is any
$B$-algebra, we have a pullback square
\[
\begin{diagram}[size=1.7em]
	&	&T_B X\Spbk&	&	\\
	&\ldTo  &	&\rdTo	&	\\
TX	&	&	&	&B	\\
	&\rdTo<{T!}&	&\ldTo  &	\\
	&	&T1,	&	&	\\
\end{diagram}
\]
and the $\cat{Q}$-object structure on $B \go T1$ induces a $\cat{Q}$-object
structure on $T_B X \go TX$.  So any $B$-algebra gives rise to an object
of $\cat{Q}$ and hence, applying the comparison functor for the monad $U\of
F$, a Penon weak $\omega$-category.%
\index{Penon, Jacques!definition of n-category@definition of $n$-category|)}%
\index{n-category@$n$-category!definitions of!Penon's|)}

\minihead{Trimble's and May's definitions}%
\index{Trimble, Todd|(}%
\index{n-category@$n$-category!definitions of!Trimble's|(}%
\index{n-category@$n$-category!flabby|(}%
\index{flabby n-category@flabby $n$-category|(}
\index{fundamental!n-groupoid@$n$-groupoid|(}

Trimble has also proposed a simple definition of $n$-category, unpublished
but written up as Definition \textbf{Tr} in my~\cite{SDN} survey.  He was
not quite so ambitious as to attempt a fully weak notion of $n$-category;
rather, he sought a notion just general enough to capture the fundamental
$n$-groupoids of topological spaces.  He called his structures `flabby
$n$-categories'.

Trimble's definition uses simple operad language.  Let $\cat{A}$ be a
category with finite products and let $P$ be a (non-symmetric) operad in
$\cat{A}$.  Extending slightly the terminology of
p.~\pageref{p:P-category}, a \demph{categorical $P$-algebra}%
\index{categorical algebra for operad}
is an
$\cat{A}$-graph $X$ together with a map
\begin{equation}	\label{eq:cat-alg-action}
P(k) \times X(x_0, x_1) \times \cdots \times X(x_{k-1}, x_k)
\go 
X(x_0, x_k)
\end{equation}
for each $k\in\nat$ and $x_0, \ldots, x_k \in X_0$, satisfying the evident
axioms.  The category of categorical $P$-algebras is written
$\fcat{CatAlg}(P)$,%
\glo{TrimCatAlg}
and itself has finite products.  If $F: \cat{A} \go
\cat{A'}$ is a finite-product-preserving functor then there is an induced
operad $F_* P$ in $\cat{A'}$ and a finite-product-preserving functor
\[
\widehat{F}: 
\fcat{CatAlg}(P) \go \fcat{CatAlg}(F_* P).
\]
  
The topological content consists of two observations: first
(Example~\ref{eg:opd-Trimble}), that if $E$ is the operad of path%
\index{operad!path reparametrizations@of path reparametrizations}
reparametrizations then there is a canonical functor
\[
\Xi: \Top \go \fcat{CatAlg}(E),
\]
and second, that taking path-components defines a finite-product-preserving
functor $\Pi_0: \Top \go \Set$.

Applying the operadic constructions recursively to the topological data, we
define for each $n\in\nat$ a category $\fcat{Flabby}\hyph n\hyph\Cat$ with
finite products and a functor
\[
\Pi_n: \Top \go \fcat{Flabby}\hyph n\hyph\Cat 
\]%
\glo{Pin}%
preserving finite products: $\fcat{Flabby}\hyph 0\hyph\Cat = \Set$, and
$\Pi_{n+1}$ is the composite functor
\[
% \Pi_{n+1}
% =
% \left(
\Top 
\goby{\Xi} 
\fcat{CatAlg}(E) 
\goby{\widehat{\Pi_n}} 
\fcat{CatAlg}((\Pi_n)_* E)
=
\fcat{Flabby}\hyph (n+1)\hyph\Cat.
% \right).
\]
That completes the definition.

Operads in $\Top$ can be regarded as $T$-operads, where $T$ is the free
topological%
\index{monoid!topological!free}
monoid monad (\ref{eg:mon-free-topo-monoid},~\ref{eg:mti-Top}).
The operadic techniques used in Trimble's definition are then expressible
in the language of generalized operads, and I make the following
\begin{claim}
  For each $n\in\nat$, there is a contractible globular $n$-operad whose
  category of algebras is equivalent to $\fcat{Flabby}\hyph n\hyph\Cat$.
\end{claim}
I hope to prove this elsewhere.  It implies that every flabby $n$-category
is a weak $n$-category in the sense of the previous chapter.  The
$n$-operad concerned is something like an $n$-dimensional version of the
operad $P$ of Example~\ref{eg:wk-omega-cat-Pi} (fundamental weak
$\omega$-groupoids), but a little smaller.

\index{May, Peter!definition of n-category@definition of $n$-category|(}%
\index{n-category@$n$-category!definitions of!May's|(}
Alternatively, we can generalize Trimble's definition by considering
operads in an arbitrary symmetric monoidal category.  This leads us to the
definition of enriched $n$-category proposed by May~\cite{MayOCA} (although
that is not what led May there).%
\index{Trimble, Todd|)}%
\index{n-category@$n$-category!definitions of!Trimble's|)}%
\index{n-category@$n$-category!flabby|)}%
\index{flabby n-category@flabby $n$-category|)}
\index{fundamental!n-groupoid@$n$-groupoid|)}

Let $\cat{B}$ and $\cat{G}$ be symmetric monoidal categories and suppose
that $\cat{B}$ acts on $\cat{G}$ in a way compatible with both monoidal
structures.  (Trivial example: $\cat{B} = \cat{G}$ with tensor as action.)
Let $P$ be an operad in $\cat{B}$.  Then we can define a
\demph{categorical%
\index{categorical algebra for operad}
$P$-algebra in $\cat{G}$} as a $\cat{G}$-graph $X$ together with maps as
in~\bref{eq:cat-alg-action}, with the first $\times$ replaced by the action
$\odot$ and the others $\times$'s by $\otimes$'s, subject to the inevitable
axioms.  For example, if $\cat{B} = \Top$ and $\cat{G} = \fcat{Flabby}\hyph
n\hyph\Cat$ with monoidal structures given by products then the functor
$\Pi_n: \cat{B} \go \cat{G}$ induces an action $\odot$ of $\cat{B}$ on
$\cat{G}$ by $B \odot G = \Pi_n(B) \times G$, and a categorical $P$-algebra
in $\cat{G}$ is what we previously called a categorical $((\Pi_n)_*
P)$-algebra.

The idea now is that if we can find some substitute for the functor $\Xi:
\Top \go \fcat{CatAlg}(E)$ then we can imitate Trimble's recursive
definition in this more general setting.  So, May starts with a `base'
symmetric monoidal category $\cat{B}$ and an operad $P$ in $\cat{B}$, each
carrying certain extra structure and satisfying certain extra properties,
the details of which need not concern us here.  He then considers symmetric
monoidal categories $\cat{G}$, restricting his attention to just those that
are `good' in the sense that they too have certain extra structure and
properties, including that they are acted on by $\cat{B}$.  Then the point
is that
\begin{itemize}
\item the trivial example $\cat{G} = \cat{B}$ is good
\item if $\cat{G}$ is good then so is the category of categorical
  $P$-algebras in $\cat{G}$.
\end{itemize}
So given a category $\cat{B}$ and operad $P$ in $\cat{B}$ as above, we can
define for each $n\in\nat$ the category $\cat{B}(n; P)$ of
\demph{$n$-$P$-categories enriched in $\cat{B}$}%
\index{enrichment!of n-categories@of $n$-categories}%
\index{n-category@$n$-category!enriched}
as follows: 
\begin{itemize}
\item $\cat{B}(0; P) = \cat{B}$
\item $\cat{B}(n+1; P) = (\textrm{categorical } P \textrm{-algebras in } 
\cat{B}(n; P) )$.
\end{itemize}
An operad $P$ in $\cat{B}$ is called an \demph{$A_\infty$-operad}%
\index{A-@$A_\infty$-!operad}
if for each $k\in\nat$, the object $P(k)$ of $\cat{B}$ is weakly equivalent
to the unit object; here `weakly equivalent' refers to a Quillen model%
\index{model category}
category structure on $\cat{B}$, which is part of the assumed structure.
May proposes that when $P$ is an $A_\infty$-operad, $n$-$P$-categories
enriched in $\cat{B}$ should be called weak $n$-categories enriched in
$\cat{B}$.  His definition, like Trimble's, aims for a slightly different
target from most of the other definitions.  It has enrichment built in; he
writes
\begin{quote}
  In all of the earlier approaches, $0$-categories are understood to be
  sets, whereas we prefer a context in which $0$-categories come with their
  own homotopy%
\index{homotopy-algebraic structure}%
\index{homotopy!zero-categories@of 0-categories}
theory.
\end{quote}
So, for instance, we might take the category $\cat{B}$ of $0$-categories to
be a convenient category of topological spaces, or the category of
simplicial sets, or a category of chain complexes.%
\index{algebraic theory!n-categories@of $n$-categories|)}
\index{May, Peter!definition of n-category@definition of $n$-category|)}%
\index{n-category@$n$-category!definitions of!May's|)}

\section{Non-algebraic definitions}
\lbl{sec:non-alg-defns-n-cat}

Most of the proposed non-algebraic definitions of weak $n$-category can be
expressed neatly using a generalization of the standard nerve construction,
which describes a category as a simplicial set with properties.  We discuss
nerves in general, then the definitions of Joyal, Tamsamani, Simpson, Baez
and Dolan (and others), Street, and Leinster (Definition $\mathbf{L'}$
of~\cite{SDN}).  We finish by looking at some structures approximating to
the idea of a weak $\omega$-category in which all cells of dimension $2$
and higher are invertible: $A_\infty$-categories, Segal categories, and
quasi-categories.  These have been found especially useful in geometry.

\minihead{Nerves}%
\index{nerve|(}

The nerve idea allows us to define species of mathematical structures by
saying on the one hand what the structures look like locally, and on the
other how the local pieces are allowed to be fitted together.  We consider
it in some generality.

The setting is a category $\cat{C}$ of `mathematical structures' with a
small subcategory $\scat{C}$ of `local pieces' or `building blocks'.  The
condition that every object of $\cat{C}$ is put together from objects of
$\scat{C}$ is called density, defined in a moment.

\begin{propn}
Let $\scat{D}$ be a small category and $\cat{D}$ a category with small
colimits.  The following conditions on a functor $I: \scat{D} \go \cat{D}$
are equivalent:
\begin{enumerate}
\item	\lbl{item:dense-coend}
for each $Y\in\cat{D}$, the canonical map
\[
\int^{D\in\scat{D}} \cat{D}(ID, Y) \times ID
\go
Y
\]
is an isomorphism
\item	\lbl{item:dense-ff}
the functor
\[
\begin{array}{rcl}
\cat{D}	&\go		&\ftrcat{\scat{D}^\op}{\Set},	\\
Y	&\goesto	&\cat{D}(I\dashbk, Y)
\end{array}
\]
is full and faithful.
\end{enumerate}
\end{propn}
\begin{proof}
See Mac Lane~\cite[X.6]{MacCWM}.
\done
\end{proof}
A functor $I: \scat{D} \go \cat{D}$ is \demph{dense}%
\index{dense}
if it satisfies the
equivalent conditions of the Proposition, and a subcategory $\scat{C}$ of a
category $\cat{C}$ is \demph{dense} if the inclusion functor $\scat{C}
\rIncl \cat{C}$ is dense.  Condition~\bref{item:dense-coend} formalizes the
idea that objects of $\cat{C}$ are pasted-together objects of $\scat{C}$.
Condition~\bref{item:dense-ff} is what we use in examples to prove density.

\begin{example}
Let $n\in\nat$, let $n\hyph\fcat{Mfd}$ be the category of smooth
$n$-manifolds%
\index{manifold}
and smooth maps, and let $\scat{U}_n$ be the subcategory
whose objects are all open subsets of $\mathbb{R}^n$ and whose maps $f: U
\go U'$ are diffeomorphisms from $U$ to an open subset of $U'$.  Then
$\scat{U}_n$ is dense in $n\hyph\fcat{Mfd}$: every manifold is a pasting of
Euclidean open sets.
\end{example}

\begin{example}
Let $k$ be a field, let $\fcat{Vect}_k$ be the category of all vector%
\index{vector space}
spaces over $k$, and let $\fcat{Mat}_k$%
\glo{Matk}
be the category whose objects are
the natural numbers and whose maps $m\go n$ are $n\times m$ matrices%
\index{matrix}
over
$k$.  There is a natural inclusion $\fcat{Mat}_k \rIncl \fcat{Vect}_k$
(sending $n$ to $k^n$), and $\fcat{Mat}_k$ is then dense in
$\fcat{Vect}_k$.  This reflects the fact that any vector space is the
colimit (pasting-together) of its finite-dimensional subspaces, which in
turn is true because the theory of vector spaces is finitary (its
operations take only a finite number of arguments).  
\end{example}

\begin{example}
The inclusion $\Delta \rIncl \Cat$%
\index{simplex category $\Delta$}
(p.~\pageref{p:defn-nerve}) is also
dense.  This says informally that a category is built out of objects,
arrows, commutative triangles, commutative tetrahedra, \ldots, and source,
target and identity functions between them.  We would not expect the
commutative tetrahedra and higher-dimensional simplices to be necessary,
and indeed, if $\Delta_2$ denotes the full subcategory of $\Delta$
consisting of the objects $\upr{0}$, $\upr{1}$ and $\upr{2}$ then the
inclusion $\Delta_2 \rIncl \Cat$ is also dense.
\end{example}
Generalizing the terminology of this example, if $\scat{C}$ is a dense
subcategory of $\cat{C}$ and $X$ is an object of $\cat{C}$ then the
\demph{nerve}%
\index{nerve!general notion}
(over $\scat{C}$) of $X$ is the presheaf $\cat{C}(\dashbk,
X)$ on $\scat{C}$.

For us, the crucial point about density is that it allows the mathematical
structures (objects of $\cat{C}$) to be viewed as
presheaves-with-properties%
\index{presheaf!properties@with properties}
on the category $\scat{C}$ of local pieces.
That is, if $\scat{C}$ is dense in $\cat{C}$ then $\cat{C}$ is equivalent
to the full subcategory of $\ftrcat{\scat{C}^\op}{\Set}$ consisting of
those presheaves isomorphic to the nerve of some object of $\cat{C}$.
Sometimes the presheaves arising as nerves can be characterized
intrinsically, yielding an alternative definition of $\cat{C}$ by
presheaves on $\scat{C}$.

\begin{example}%
\index{vector space}
A vector space over a field $k$ can be \emph{defined} as a presheaf $V:
\fcat{Mat}_k^\op \go \Set$ preserving finite limits.
\end{example}

\begin{example}%
\index{manifold!alternative definition of}
A smooth $n$-manifold can be \emph{defined} as a presheaf on $\scat{U}_n$
with certain properties.
\end{example}

\begin{example}	\lbl{eg:nerve-cat-chars}
A category can be \emph{defined} as a presheaf on $\Delta$%
\index{simplex category $\Delta$}
(or indeed on $\Delta_2$) with certain properties.  There are various ways
to express those properties: for instance, categories are functors
$\Delta^\op \go \Set$ preserving finite limits, or preserving certain
pullbacks, or they are simplicial sets in which every inner horn%
\index{horn!filler for}
has a
unique filler.  (A horn is \demph{inner}%
\index{horn!inner}
if the missing face is not the first or the last one.)
\end{example}

\index{nerve!n-category@of $n$-category|(}
The structures we want to define are weak $n$-categories, and the strategy
is:
\begin{itemize}
\item find a small dense subcategory $\scat{C}$ of the category of
  \emph{strict} $n$-categories
\item find conditions on a presheaf on $\scat{C}$ equivalent to it being a
  nerve of a strict $n$-category
\item relax%
\index{weakening!theory of n-categories@for theory of $n$-categories}
those conditions to obtain a definition of \emph{weak}
  $n$-category.
\end{itemize}
Most of the definitions of weak $n$-category described below can be
regarded as implementations of this strategy.  A different way to put it is
that we seek an intrinsic characterization of presheaves on $\scat{C}$ of
the form
\[
C
\ \goesto\ 
\{ \textrm{weak functors } C \go Y \}
\]
for some weak $n$-category $Y$.  Of course, we start from a position of not
knowing what a weak $n$-category or functor is, but we choose the
conditions on presheaves to fit the usual intuitions.%
\index{nerve|)}

\minihead{Joyal's definition}%
\index{Joyal, Andr\'e!definition of n-category@definition of $n$-category|(}%
\index{n-category@$n$-category!definitions of!Joyal's|(}%

Perhaps the most obvious implementation is to take the local pieces to be
all globular%
\index{pasting diagram!globular}
pasting diagrams.  So, let $\Delta_\omega$%
\glo{Deltaomega}%
\index{simplex category $\Delta$!omega-dimensional analogue of@$\omega$-dimensional analogue of}
be the category
with object-set $\coprod_{m\in\nat} \pd(m)$ and hom-sets
\[
\Delta_\omega (\sigma, \pi)
=
\strcat{\omega}(F\rep{\sigma}, F\rep{\pi})
\iso
\ftrcat{\scat{G}^\op}{\Set} (\rep{\sigma}, T\rep{\pi})
\]
where $F: \ftrcat{\scat{G}^\op}{\Set} \go \strcat{\omega}$ is the free
strict $\omega$-category functor and, as in Chapter~\ref{ch:a-defn}, $T$ is
the corresponding monad.  There is an inclusion functor $\Delta_\omega
\rIncl \strcat{\omega}$, and a typical object of the corresponding
subcategory of $\strcat{\omega}$ is the $\omega$-category naturally
depicted as
\[
\gfstsu\gfoursu\gzersu\gonesu\gzersu\gthreesu\glstsu
\]
---that is, freely generated by the 0-, 1- and 2-cells shown, and with only
identity cells in dimensions 3 and above.  

The subcategory $\Delta_\omega$ of $\strcat{\omega}$ is dense.  That the
induced functor $\strcat{\omega} \go \ftrcat{\Delta_\omega^\op}{\Set}$ is
faithful follows from $\Delta_\omega$ containing the trivial pasting
diagrams
\[
\gzersu\,,
\diagspace
\gfstsu\gonesu\glstsu,
\diagspace
\gfstsu\gtwosu\glstsu,
\diagspace
\ldots
\]
representing single $m$-cells.  That it is full follows from
$\Delta_\omega$ containing maps such as 
\begin{equation}	\label{eq:comp-inducers}
\gfstsu\gtwosu\glstsu
\ \goby{f}\ 
\gfstsu\gonesu\glstsu
\diagspace
\textrm{and}
\diagspace
\gfstsu\gtwosu\glstsu
\ \goby{g}\ 
\gfstsu\gthreesu\glstsu
\end{equation}
which, for appropriately chosen $f$ and $g$, induce 2-cell identities and
vertical 2-cell composition respectively.

A strict $\omega$-category can therefore be defined as a presheaf on
$\Delta_\omega$ with properties.  Joyal has proposed~\cite{JoyDDT} a way of
describing and then relaxing those properties: a strict $\omega$-category
is a presheaf on $\Delta_\omega$ for which `every inner horn%
\index{horn!filler for}
has a unique
filler',
and a weak $\omega$-category is defined by simply dropping the
uniqueness.

We can do the same with $\omega$ replaced by any finite $n$ (taking care in
the top%
\index{top dimension}
dimension).  Recall from p.~\pageref{p:degen-rep} that if $1_\pi$
denotes the $(m+1)$-pasting diagram resembling an $m$-pasting diagram $\pi$
then $\rep{\pi} \iso \rep{1_\pi}$: so $\Delta_n$ is equivalent to its full
subcategory consisting of just the $n$-pasting diagrams.  For instance,
$\Delta_1$ is equivalent to the usual category $\Delta$ of $1$-pasting
diagrams
\[
\gzersu\,,
\diagspace
\gfstsu\gonesu\glstsu,
\diagspace
\gfstsu\gonesu\gzersu\gonesu\glstsu, 
\diagspace
\ldots,
\]
and we recover the standard nerve construction for categories.  

Joyal also noted a duality.%
\index{duality!intervals vs. ordered sets@intervals \vs.\ ordered sets}
 The category $\Delta$ is equivalent to
the opposite of the category $\scat{I}$ of \demph{finite strict intervals},
that is, finite totally ordered sets with distinct least and greatest
elements (to be preserved by the maps).  Generalizing this, he defined a
category $\scat{I}_\omega$ of `finite disks',%
\index{disk!Joyal's sense@in Joyal's sense}
equivalent to the opposite of
$\Delta_\omega$.  
% (See the Notes at the end of the chapter for further
% references.)  
So his weak $\omega$-categories are functors $\scat{I}_\omega \go \Set$
satisfying horn-filling conditions.%
\index{Joyal, Andr\'e!definition of n-category@definition of $n$-category|)}%
\index{n-category@$n$-category!definitions of!Joyal's|)}%

\minihead{Tamsamani's and Simpson's definitions}%
\index{Simpson, Carlos!definition of n-category@definition of $n$-category|(}%
\index{n-category@$n$-category!definitions of!Simpson's|(}%
\index{Tamsamani, Zouhair!definition of n-category@definition of $n$-category|(}%
\index{n-category@$n$-category!definitions of!Tamsamani's|(}%

A very similar story can be told for the definitions proposed by
Tamsamani~\cite{TamSNN} and Simpson~\cite{SimCMS}.  Observe that in the
proof of the density of $\Delta_\omega$ in $\strcat{\omega}$, we did not
use many of the pasting diagrams, so we can replace $\Delta_\omega$ by a
smaller category.  Tamsamani and Simpson consider just `cuboidal'%
\index{pasting diagram!cubical}
pasting
diagrams such as
\[
\gfstsu\gfoursu\gzersu\gfoursu\gzersu\gfoursu\gzersu\gfoursu\glstsu
\]
and similarly `cuboidal' maps between the strict $\omega$-categories that
they generate.

Let us restrict ourselves to the $n$-dimensional case, since that is a
little easier.  There is a functor
\[
I: \Delta^n \go \strcat{n}
\]%
\index{simplex category $\Delta$}%
which, for instance, when $n=2$, sends $(\upr{4}, \upr{3})$ to the free
strict 2-category on the diagram above.  In general, each $(r_1, \ldots,
r_n) \in \nat^n$ determines an $n$-pasting diagram $\pi_{r_1, \ldots,
r_n}$, given inductively by
\[
\pi_{r_1, \ldots, r_n}
=
\left(\pi_{r_2, \ldots, r_n}, \ldots, \pi_{r_2, \ldots, r_n}\right)
\]
with $r_1$ terms on the right-hand side, and then
\[
I(\upr{r_1}, \ldots, \upr{r_n})
=
F\rep{\pi_{r_1, \ldots, r_n}}.
\]
To describe $I$ on maps, take, for instance, $n=2$ and the map $(\id,
\delta): (\upr{1}, \upr{1}) \go (\upr{1}, \upr{2})$ in which $\delta$ is
the injection omitting $1 \in \upr{2}$ from its image; then $I(\id,
\delta)$ is the map $g$ of~\bref{eq:comp-inducers}.

By exactly the same argument as for Joyal's definition, the functor
$\Delta^n \go \strcat{n}$ is dense.  A strict $n$-category is therefore the
same thing as a presheaf on $\Delta^n$ (a `multisimplicial%
\index{multisimplicial set}%
\index{simplicial set!multi-}
set') with
properties, and relaxing those properties gives a definition of weak
$n$-category.

Nerves of strict $n$-categories are characterized among functors
$(\Delta^n)^\op \go \Set$ by the properties that the functor is degenerate
in certain ways (to give us $n$-categories rather%
\index{n-tuple category@$n$-tuple category!degenerate}
than $n$-tuple
categories) and, more significantly, that certain pullbacks are preserved.
Tamsamani sets up a notion of equivalence,%
\index{equivalence!n-categories@of $n$-categories}
and defines weak $n$-category by
asking only that the pullbacks are preserved up to equivalence.  Simpson
does the same, but with a more stringent notion of equivalence that he
calls `easy%
\index{equivalence!easy}
equivalence'.  It is indeed easier, and is nearly the same as
the notion of contractibility%
\index{contractible!map of globular sets}
of a map of globular sets: see my
survey~\cite{SDN} for details.  In the special case of one-object
2-categories, Tamsamani's definition gives the homotopy monoidal categories
of Section~\ref{sec:non-alg-notions}, and Simpson's gives the same but with
the extra condition that the functors $\xi^{(k)}$ of
Proposition~\ref{propn:simp-eqs}, which for homotopy%
\index{monoidal category!homotopy}
monoidal categories
are required to be equivalences, are \emph{genuinely} surjective on
objects.%
\index{Simpson, Carlos!definition of n-category@definition of $n$-category|)}%
\index{n-category@$n$-category!definitions of!Simpson's|)}%
\index{Tamsamani, Zouhair!definition of n-category@definition of $n$-category|)}%
\index{n-category@$n$-category!definitions of!Tamsamani's|)}%

\minihead{Opetopic definitions}%
\index{Baez, John!definition of n-category@definition of $n$-category|(}%
\index{Dolan, James!definition of n-category@definition of $n$-category|(}%
\index{n-category@$n$-category!definitions of!opetopic|(}%

We have already~(\ref{sec:ope-n}) looked at the opetopic definitions of
weak $n$-category: that of Baez and Dolan and subsequent variants.  They
are all of the form `a weak $n$-category is an opetopic set with
properties', for varying meanings of `opetopic set' and varying lists of
properties.  Here we see how this fits in with the nerve idea.  The
situation is, as we shall see, slightly different from that in the
definitions of Joyal, Tamsamani, and Simpson.

We start with the category $\scat{O}$ of opetopes.  As mentioned
on p.~\pageref{p:ope-to-n-cat}, there is an embedding $\scat{O} \rIncl
\strcat{\omega}$, and this induces a functor 
\[
U: \strcat{\omega} \go \ftrcat{\scat{O}^\op}{\Set}
\]
sending a strict $\omega$-category to its underlying opetopic set.

This functor is faithful, because there is for each $m\in\nat$ an
$m$-opetope resembling a single globular $m$-cell.  It is not, however,
full.  To see this, note that if $F: A \go B$ is a strict map of strict
$\omega$-categories then the induced map $U(F): U(A) \go U(B)$ preserves
universality%
\index{universal!preservation}
of cells: for instance, if $f$ and $g$ are abutting 1-cells in
$A$ then $U(F)$ sends the canonical 2-cell
\[
\topeb{f}{g}{g\of f}{\Downarrow}
\]
in $U(A)$ to the canonical 2-cell 
\[
\topeb{Ff}{Fg}{(Fg)\of (Ff)}{\Downarrow}
\]
in $U(B)$.  But not every map $U(A) \go U(B)$ of opetopic sets preserves
universality; indeed, any \emph{lax} map $A \go B$ of strict
$\omega$-categories ought (in principle, at least) to induce a map $U(A)
\go U(B)$, and this will preserve universality if and only if the lax map
is weak.  Compare the relationship between monoidal categories and plain
multicategories (\ref{eg:map-mti-mon},~\ref{sec:non-alg-notions}).

So $\scat{O}$ is not dense in $\strcat{\omega}$, and correspondingly $U$
does not define an equivalence between $\strcat{\omega}$ and a full
subcategory of $\ftrcat{\scat{O}^\op}{\Set}$.  But with a slight
modification, the nerve idea can still be made to work.  For the above
arguments suggest that $U$ defines an equivalence between
\[
(\textrm{strict } \omega \textrm{-categories } + \textrm{ weak maps})
\]
and a full subcategory of 
\[
(\textrm{opetopic sets } + \textrm{universality-preserving maps}),
\]
and it is then, as usual, a matter of identifying the characteristic
properties of those opetopic sets arising from strict $\omega$-categories,
then relaxing the properties to obtain a definition of weak
$\omega$-category.  So a weak $\omega$-category is defined as an opetopic
set with properties, and a weak $\omega$-functor as a map of opetopic sets
preserving universality.%
\index{Baez, John!definition of n-category@definition of $n$-category|)}%
\index{Dolan, James!definition of n-category@definition of $n$-category|)}%
\index{n-category@$n$-category!definitions of!opetopic|)}%

\minihead{Street's definition}%
\index{Street, Ross!definition of n-category@definition of $n$-category|(}%
\index{n-category@$n$-category!definitions of!Street's|(}%
\index{simplicial set!n-category as@$n$-category as|(}

The definition of weak $\omega$-category proposed by Street has the
distinctions of being the first and probably the most tentatively phrased;
it hides in the last paragraph of his paper of~\cite{StrAOS}.  It was part
of the inspiration for Baez and Dolan's definition, and has much in common
with it, but uses simplicial rather than opetopic sets.

Street follows the nerve idea explicitly.  He first constructs an embedding
\[
I: \Delta \rIncl \strcat{\omega},
\]
where $I(m)$ is the $m$th `oriental',%
\index{oriental}
the free strict%
\index{omega-category@$\omega$-category!strict!free}
$\omega$-category on
an $m$-simplex.  For example, $I(3)$ is the strict $\omega$-category
generated freely by 0-, 1-, 2- and 3-cells
\[
\setlength{\unitlength}{1em}
\begin{picture}(21,4.6)
\cell{0}{0}{bl}{%
\begin{picture}(7,4.6)(-0.5,-0.8)
% 0-dim labels
\cell{0}{0}{c}{a_0}
\cell{1}{3}{c}{a_1}
\cell{5}{3}{c}{a_2}
\cell{6}{0}{c}{a_3}
% arrows
\put(0.05,0.4){\vector(1,3){0.75}}
\put(1.4,3){\vector(1,0){3}}
\put(5.2,2.65){\vector(1,-3){0.75}}
\put(0.5,0){\vector(1,0){4.9}}
\qbezier(0.4,0.3)(2.5,1.5)(4.6,2.7)
\put(4.6,2.7){\vector(3,2){0}}
% 1-dim labels
\cell{0.3}{1.5}{r}{\scriptstyle f_{01}}
\cell{3}{3.2}{b}{\scriptstyle f_{12}}
\cell{5.7}{1.5}{l}{\scriptstyle f_{23}}
\cell{3}{-0.2}{t}{\scriptstyle f_{03}}
\cell{2.2}{0.7}{c}{\scriptstyle f_{02}}
% 2-dim arrows
\cell{1.5}{2}{c}{\rotatebox{35}{$\Uparrow$}}
\cell{3.7}{1}{c}{\Uparrow}
% 2-dim labels
\cell{1.8}{2.2}{l}{\scriptstyle \alpha_{012}}
\cell{4}{0.9}{l}{\scriptstyle \alpha_{023}}
\end{picture}}
\cell{10.5}{2.7}{t}{\Rrightarrow}
\cell{10.5}{3}{b}{\Gamma}
\cell{14}{0}{bl}{%
\begin{picture}(7,4.6)(-0.5,-0.8)
% 0-dim labels
\cell{0}{0}{c}{a_0}
\cell{1}{3}{c}{a_1}
\cell{5}{3}{c}{a_2}
\cell{6}{0}{c}{a_3\makebox[0em][l]{.}}
% arrows
\put(0.05,0.4){\vector(1,3){0.75}}
\put(1.4,3){\vector(1,0){3}}
\put(5.2,2.65){\vector(1,-3){0.75}}
\put(0.5,0){\vector(1,0){4.9}}
\qbezier(5.6,0.3)(3.5,1.5)(1.4,2.7)
\put(5.6,0.3){\vector(3,-2){0}}
% 1-dim labels
\cell{0.3}{1.5}{r}{\scriptstyle f_{01}}
\cell{3}{3.2}{b}{\scriptstyle f_{12}}
\cell{5.7}{1.5}{l}{\scriptstyle f_{23}}
\cell{3}{-0.2}{t}{\scriptstyle f_{03}}
\cell{3.8}{0.7}{c}{\scriptstyle f_{13}}
% 2-dim arrows
\cell{4.5}{2}{c}{\rotatebox{-35}{$\Uparrow$}}
\cell{2.3}{1}{c}{\Uparrow}
% 2-dim labels
\cell{4.2}{2.2}{r}{\scriptstyle \alpha_{123}}
\cell{2}{0.9}{r}{\scriptstyle \alpha_{013}}
\end{picture}}
\end{picture}
% \hand{30}{66}.
\]
(Orientation needs care.)  This induces a functor 
\[
U: \strcat{\omega} \go \ftrcat{\Delta^\op}{\Set}
\]
and, roughly speaking, a weak $\omega$-category is defined as a simplicial
set with horn-filling%
\index{horn!filler for}
properties.  

This is, however, a slightly inaccurate account.  For similar reasons to
those in the opetopic case, $\Delta$ is not dense in $\strcat{\omega}$; the
functor $U$ is again faithful but not full.  Street's original solution was
to replace $\ftrcat{\Delta^\op}{\Set}$ by the category $\fcat{Sss}$ of
\demph{stratified%
\index{simplicial set!stratified}
simplicial sets}, that is, simplicial sets equipped with
a class of distinguished cells in each dimension (to be thought of as
`universal',%
\index{universal!cell of n-category@cell of $n$-category}
`hollow',%
\index{hollow}
or `thin').%
\index{thin}
 The underlying simplicial set of a
strict $\omega$-category has a canonical stratification, so $U$ lifts to a
functor
\[
U': \strcat{\omega} \go \fcat{Sss},
\]
and $U'$ \emph{is} full and faithful.  Detailed work by Street
\cite{StrAOS,StrFN} and Verity%
\index{Verity, Dominic}
(unpublished) gives precise conditions for
an object of $\fcat{Sss}$ to be in the image of $U'$.  So a strict
$\omega$-category is the same thing as a simplicial set equipped with a
class of distinguished cells satisfying some conditions.  One of the
conditions is that certain horns have unique fillers, and dropping the
uniqueness gives Street's proposed definition of weak $\omega$-category.

The most vexing aspect of this proposal is that extra structure is required
on the simplicial set.  It would seem more satisfactory if, as in the
opetopic approach, the universal cells could be recognized intrinsically.
A recent paper of Street~\cite{StrWOC} aims to repair this apparent defect,
proposing a similar definition in which a weak $\omega$-category is
genuinely a simplicial set with properties.%
\index{Street, Ross!definition of n-category@definition of $n$-category|)}%
\index{n-category@$n$-category!definitions of!Street's|)}%
\index{simplicial set!n-category as@$n$-category as|)}
\index{nerve!n-category@of $n$-category|)}

\minihead{Contractible multicategories}%
\index{multicategory!contractible|(}%
\index{contractible!multicategory|(}%
\index{globular multicategory|(}
\index{n-category@$n$-category!definitions of!contractible multicategory|(}

The last definition of weak $n$-category that we consider was introduced as
Definition $\mathbf{L'}$ in my~\cite{SDN} survey.  Like the other
definitions in this section, it is non-algebraic and can be described in
terms of nerves.  The nerve description seems, however, to be rather
complicated (the shapes involved being a combination of globular and
opetopic) and not especially helpful, so we approach it from another
angle instead.  

The idea is that a weak $\omega$-category is meant to be a `weak%
\index{algebra!monad@for monad!weak}%
\index{weakening!theory of n-categories@for theory of $n$-categories}
algebra'
for the free strict $\omega$-category monad on globular sets.  We saw
in~\ref{eg:multi-alg} that for any cartesian monad $T$, a \emph{strict}
$T$-algebra is the same thing as a $T$-multicategory whose domain map is
the identity---in other words, with underlying graph of the form
\[
\begin{diagram}[size=1.7em]
	&	&TX	&	&	\\
	&\ldTo<1&	&\rdTo>h&	\\
TX	&	&	&	&X.	\\  
\end{diagram}
\]
To define `weak $T$-algebra' we relax the condition that the domain map is
the identity, asking only that it be an equivalence in some sense.
Contractibility together with surjectivity on 0-cells is a reasonable
notion of equivalence: contractible means something like `injective%
\index{homotopy!injective on}
on
homotopy', and any contractible map surjective on $0$-cells is surjective
on $m$-cells for all $m\in\nat$.  For $1$-dimensional structures, it means
full, faithful and surjective on objects.  We also ask that the domain map
is injective on $0$-cells, expressing the thought that $0$-cells in an
$\omega$-category should not be composable.

Definition $\mathbf{L'}$ says, then, that a weak $\omega$-category is a
globular multicategory $C$ whose domain map $C_1 \go TC_0$ is bijective on
$0$-cells and contractible. 

Any weak $\omega$-category in the sense of the previous chapter gives rise
canonically to one in the sense of $\mathbf{L'}$.  This is the
`multicategory%
\index{generalized multicategory!elements@of elements}
of elements' construction of~\ref{sec:alg-fibs}: if $(T_L X
\goby{h} X)$ is an $L$-algebra then there is a commutative diagram
\[
% \begin{diagram}[size=1.7em]
\begin{slopeydiag}
	&		&T_L X		&		&	\\
	&\ldTo		&		&\rdTo>h	&	\\
TX	&		&\dTo		&		&X	\\
\dTo<{T!}&		&L		&		&\dTo>!	\\
	&\ldTo		&		&\rdTo		&	\\
T1	&		&		&		&1	\\
\end{slopeydiag}
% \end{diagram}
\]
the left-hand half of which is a pullback square and the top part of which
forms a multicategory $C^X$ with $C^X_0 = X$ and $C^X_1 = T_L X$.  The map
$L \go T1$ is contractible (by definition of $L$) and bijective on
$0$-cells (because, as we saw on p.~\pageref{p:L-blob}, if $\blob$ denotes
the $0$-pasting diagram then $L(\blob) = 1$).  Contractibility and
bijectivity on 0-cells are stable under pullback, so the multicategory
$C^X$ is a weak $\omega$-category in the sense of $\mathbf{L'}$.

Unpicking this construction explains further the idea behind $\mathbf{L'}$.
An $m$-cell of $T_L X$ is a pair 
\[
(\theta, \mathbf{x}) 
\in
\coprod_{\pi\in\pd(m)}
L(\pi) 
\times
\ftrcat{\scat{G}^\op}{\Set}(\rep{\pi}, X),
\]
which lies over $\mathbf{x} \in (TX)(m)$ and $\ovln{\theta}(\mathbf{x}) \in
X(m)$.  It is usefully regarded as a `way of composing' the labelled
pasting diagram $\mathbf{x}$.  (Contractibility guarantees that there are
plenty of ways of composing.)  Among all weak $\omega$-categories in the
sense of $\mathbf{L'}$, those of the form $C^X$ have the special feature
that the set of ways of composing a labelled pasting diagram $\mathbf{x}
\in (TX)(\pi)$ depends only on the pasting diagram $\pi$, not on the
labels: it is just $L(\pi)$.

So in the definition of the previous chapter, the ways of composition
available in a weak $\omega$-category are prescribed once and for all; in
the present definition, they are allowed to vary from $\omega$-category to
$\omega$-category.  This is precisely analogous to the difference between
the loop%
\index{loop space!machine}
space machinery of Boardman--Vogt%
\index{Boardman, Michael}%
\index{Vogt, Rainer}
and May%
\index{May, Peter}
(with fixed parameter
spaces forming an operad) and that of Segal%
\index{Segal, Graeme}
(with a variable, flabby,%
\index{flab}
structure): see Adams~\cite[p.~60]{Ad}.%
\index{multicategory!contractible|)}%
\index{contractible!multicategory|)}%
\index{globular multicategory|)}
\index{n-category@$n$-category!definitions of!contractible multicategory|)}

\minihead{Locally groupoidal structures}%
\index{n-category@$n$-category!locally groupoidal|(}%
\index{locally groupoidal structures|(}

Many of the weak $\omega$-categories of interest in geometry have the
property that all cells of dimension $2$ and higher are equivalences
(weakly invertible).%
\index{invertibility}
 Some examples were given in `Motivation for
Topologists'.  There are several structures in use that, roughly speaking,
aim to formalize the idea of such a weak $\omega$-category.  I will
describe three of them here.

A weak $\omega$-category in which all cells (of dimension $1$ and higher)
are equivalences is called a \demph{weak $\omega$-groupoid}.%
\index{omega-groupoid@$\omega$-groupoid}
  This is, of
course, subject to precise definitions of weak $\omega$-category and
equivalence.  From the topological viewpoint one of the main purposes of
$\omega$-groupoids is to model homotopy%
\index{homotopy!type}%
\index{fundamental!omega-groupoid@$\omega$-groupoid}
types of spaces (see Grothendieck's
letter of~\cite{GroPS}, for instance), so it is reasonable to replace
$\omega$-groupoids by spaces, or perhaps simplicial%
\index{simplicial set!omega-category from@$\omega$-category from}
sets or chain
complexes.%
\index{chain complex!omega-category from@$\omega$-category from}
  The structures we seek are, therefore, graphs $(X(x,x'))_{x, x'
\in X_0}$ of spaces, simplicial sets, or chain complexes, together with
extra data determining some kind of weak composition.

We have already seen one version of this: an $A_\infty$-category%
\index{A-@$A_\infty$-!category}
(p.~\pageref{p:A-infty-category}) is a graph $X$ of chain complexes together
with various composition maps
\[
X(x_{k-1}, x_k) \otimes \cdots \otimes X(x_0, x_1)
\go
X(x_0, x_k)
\]
parametrized by elements of the operad $A_\infty$.  There is a similar
notion with spaces in place of complexes.  These are algebraic definitions
(so properly belong in the previous section).

A similar but non-algebraic notion is that of a Segal
category (sometimes
called by other names: see the Notes below).  Take a bisimplicial%
\index{bisimplicial set}%
\index{simplicial set!bi-}
set,
expressed as a functor
\[
X: \Delta^\op \go \ftrcat{\Delta^\op}{\Set},
\]
and suppose that the simplicial set $X\upr{0}$ is discrete (that is, the
functor $X\upr{0}$ is constant).  Write $X_0$ for the set of points
(constant value) of $X\upr{0}$ .  Then for each $k\in\nat$, the simplicial
set $X\upr{k}$ decomposes naturally as a coproduct
\[
 X\upr{k} 
\iso 
\coprod_{x_0, \ldots, x_k \in X_0} 
X(x_0, \ldots, x_k),
\]
and for each $k\in\nat$ and $x_0, \ldots, x_k \in X_0$, there is a natural
map
\[
X(x_0, \ldots, x_k)
\go
X(x_{k-1}, x_k) \times \cdots \times X(x_0, x_1).
\]
A bisimplicial set $X$ is called a \demph{Segal%
\index{Segal, Graeme!category}
 category} if $X\upr{0}$ is
discrete and each of these canonical maps is a weak equivalence of
simplicial sets, in the homotopy-theoretic sense.  This definition is very
closely related to the definitions of $n$-category proposed by Simpson%
\index{Simpson, Carlos!definition of n-category@definition of $n$-category}%
\index{n-category@$n$-category!definitions of!Simpson's}%
and
Tamsamani,%
\index{Tamsamani, Zouhair!definition of n-category@definition of $n$-category}%
\index{n-category@$n$-category!definitions of!Tamsamani's}%
 and in particular to the definition of homotopy%
\index{bicategory!homotopy}
bicategory in
Chapter~\ref{ch:monoidal}.

Finally, there are the quasi-categories of Joyal, Boardman, and Vogt.  We
have a diagram
\[
\begin{diagram}[size=1.7em]
	&	&\{\textrm{quasi-categories}\}	&	&	\\
	&\ruIncl&			&\luIncl&	\\
\{\textrm{categories}\}&&		&	&
\{\textrm{Kan complexes}\}\\
	&\luIncl&			&\ruIncl&	\\
	&	&\{\textrm{groupoids}\}&	&	\\
\end{diagram}%
\index{groupoid}
\]
of classes of simplicial sets, in which
\begin{itemize}
\item categories are identified with simplicial sets in which every inner
  horn%
\index{horn!filler for}
has a unique filler~(\ref{eg:nerve-cat-chars})
\item \demph{Kan%
\index{Kan, Daniel!complex}
complexes} are simplicial sets in which every horn
  has at least one filler (the principal example being the underlying
  simplicial set of a space)
\item at the intersection, groupoids are simplicial sets in which every
  horn has a unique filler
\item at the union, \demph{quasi-categories}%
\index{quasi-category}%
\index{category!quasi-}
are simplicial sets in which
  every inner horn has at least one filler.
\end{itemize}
Large amounts of the theory of ordinary categories can be reproduced for
quasi-categories, although requiring much longer proofs.  Given a
simplicial set $X$ and elements $x, x' \in X\upr{0}$, there is a simplicial
set $X(x, x')$, the analogue of the space of paths from $x$ to $x'$ in
topology, and it can be shown that if $X$ is a quasi-category then each
$X(x, x')$ is a Kan complex.  So quasi-categories do indeed approximate the
idea of a weak $\omega$-category in which all cells of dimension at least
$2$ are invertible.%
\index{n-category@$n$-category!locally groupoidal|)}%
\index{locally groupoidal structures|)}

\begin{notes}

An extensive bibliography and historical discussion of proposed
definitions of $n$-category is in my survey paper~\cite{SDN}. 

Various people have confused the contractions%
\index{contraction!notions of}
of Chapter~\ref{ch:a-defn}
with the contractions of Batanin, principally me (see the Notes to
Chapter~\ref{ch:a-defn}) but also Berger~\cite[1.20]{Ber}.%
\index{Berger, Clemens}
 I apologize to
Batanin for stealing his word, using it for something else, then
renaming his original concept~(\ref{defn:coh-coll}).

Berger~\cite{Ber}%
\index{Berger, Clemens}
has investigated nerves of $\omega$-categories in detail,
making connections to various proposed definitions of weak
$\omega$-categories, especially Joyal's.  Other work on higher-dimensional
nerves has been done by Street%
\index{Street, Ross}
(\cite[\S 10]{StrCS} and references therein)
and Duskin%
\index{Duskin, John}
\cite{DusSM1, DusSM2}.  Joyal's definition has also been
illuminated by Makkai%
\index{Makkai, Michael}
and Zawadowski~\cite{MZ}%
\index{Zawadowski, Marek}
and Batanin%
\index{Batanin, Michael}
and
Street~\cite{BSUPM}.%
\index{Street, Ross}

The observation that most geometrically interesting $\omega$-categories are
locally groupoidal was made to me by Bertrand To\"en.%
\index{To\"en, Bertrand}
 His~\cite{TV} paper
with Vezzosi%
\index{Vezzosi, Gabriele}
gives an introduction to Segal categories, as well as further
references.  In particular, they point to Dwyer,%
\index{Dwyer, William}
Kan%
\index{Kan, Daniel}
and Smith~\cite{DKS},%
\index{Smith, Jeffrey}
where Segal categories were called special%
\index{bisimplicial set!special}%
\index{simplicial set!bi-!special}%
\index{special}
bisimplicial sets, and
Schw\"anzl%
\index{Schwanzl@Schw\"anzl, Roland}
and Vogt~\cite{SV},%
\index{Vogt, Rainer}
where they were called $\Delta$-categories.%
\index{Delta-category@$\Delta$-category}%
\index{category!Delta-@$\Delta$-}

Joyal's quasi-categories are a renaming of Boardman%
\index{Boardman, Michael}
and Vogt's%
\index{Vogt, Rainer}
`restricted%
\index{Kan, Daniel!complex}
Kan complexes'~\cite{BV}.  His work on quasi-categories
remains unpublished, but has been presented in seminars since 1997 or
earlier.

I thank Sjoerd Crans for useful conversations on the details of Batanin's
definition, Jacques Penon for the observation that coherent maps do not
compose, and Michael Batanin and Ross Street for useful comments on their
respective definitions.

\end{notes}

\part*{Appendices}

\appendix

% \addtocontents{toc}{\contentsline {chapter}{\numberline {}}{}}
\ucontents{part}{Appendices}

\chapter{Symmetric Structures}
\lbl{app:sym}

\chapterquote{%
For my birthday I got a humidifier and a de-humidifier \ldots\  I put them
in the same room and let them fight it out}{%
Steven Wright}

\noindent
Here we meet an alternative definition of symmetric multicategory and
prove it equivalent to the usual one.  This has two purposes.  First,
the alternative definition is in some ways nicer and more natural than the
usual one, avoiding as it does the delicate matter of formulating the
symmetry axioms (\ref{eg:opd-Sym},~\ref{defn:sym-mti}).  Second, it will be
used in Appendix~\ref{app:special-cart} to show that every symmetric
multicategory gives rise to a $T$-multicategory for each $T$ belonging to a
certain large class of cartesian monads.

Actually, we start~(\ref{sec:comm-mons}) with an alternative definition of
commutative monoid.  Although this could hardly be shorter than the
standard definition, it acts as a warm-up to the alternative definition of
symmetric multicategory~(\ref{sec:sym-multis}).

\section{Commutative monoids}
\lbl{sec:comm-mons}%
\index{monoid!commutative|(}

Here we define `fat commutative monoids' and prove that they are
essentially the same as the ordinary kind.  The idea is to have some direct
way of summing arbitrary finite families $(a_x)_{x\in X}$ of elements of a
commutative monoid $A$, not just ordered sequences $a_1, \ldots, a_n$.
\begin{defn}	\lbl{defn:fat-cm}
A \demph{fat commutative monoid}%
\index{monoid!commutative!fat}
is a set $A$ equipped with a function
$\sum_X: A^X \go A$%
\glo{fatsum}
for each finite set $X$, satisfying the axioms below.
We write elements of $A^X$ as families $(a_x)_{x\in X}$ (where $a_x \in
A$) and $\sum_{X} (a_x)_{x\in X}$ as $\sum_{x\in X} a_x$.  The axioms are:
\begin{itemize}
\item for any map $s: X \go Y$ of finite sets and any family $(a_x)_{x\in
X}$ of elements of $A$, 
\[
\sum_{y\in Y} \sum_{x \in s^{-1} \{y\}} a_x
=
\sum_{x\in X} a_x
\]
\item for any one-element set $X$ and any $a\in A$, 
\[
a = \sum_{x\in X} a.
\]
\end{itemize}
A \demph{map} $A \go A'$ of fat commutative monoids is a function $f: A \go
A'$ such that for all finite sets $X$, the square
\[
\begin{diagram}[size=2em]
A^X		&\rTo^{f^X}	&A'^X		\\
\dTo<{{\textstyle\Sigma}_X}	&		&\dTo>{{\textstyle\Sigma}_X}\\
A		&\rTo_{f}	&A'		\\
\end{diagram}
\]
commutes.  
\end{defn}
Observe the following crucial property immediately:
\begin{lemma}	\lbl{lemma:fat-cm}
Let $A$ be a fat commutative monoid, $s: X \go Y$ a bijection between
finite sets, and $(b_y)_{y\in Y}$ an indexed family of elements of $A$.
Then 
\[
\sum_{x\in X} b_{s(x)} = \sum_{y\in Y} b_y.
\]
\end{lemma}
\begin{proof}
Define a family $(a_x)_{x\in X}$ by $a_x = b_{s(x)}$.  Then
\[
\sum_{x\in X} a_x 
=
\sum_{y\in Y} \sum_{x\in s^{-1}\{y\}} a_x 
=
\sum_{y\in Y} b_y,
\]
by the first and second axioms respectively.
\done
\end{proof}

This allows us to take the expected liberties with notation.  If, for
example, we have a family of elements $a_{v,w,x}$ indexed over finite sets
$V$, $W$ and $X$, then we may write $\sum_{v\in V, w\in W, x\in X}
a_{v,w,x}$ without ambiguity; the sum could `officially' be interpreted as
either of
\[
\sum_{((v,w),x) \in (V\times W)\times X} a_{v,w,x}
\diagspace
\textrm{or}
\diagspace
\sum_{(v,(w,x)) \in V\times (W\times X)} a_{v,w,x}
\]
(or some further possibility), but these expressions are equal.

Those concerned with foundations might feel uneasy about the idea of
specifying a function $\sum_X: A^X \go A$ `for each finite%
\lbl{p:quantification-qualms}
set $X$'.  The remedy is to choose a small full subcategory $\scat{F}$ of
the category of finite sets and functions, such that $\scat{F}$ contains at
least one object of each finite cardinality, and interpret `finite set' as
`object of $\scat{F}\,$'; everything works just as well.  In particular,
you might choose to replace the category of finite sets with its skeleton
whose objects are the natural numbers $\mb{n} = \{1, \ldots, n\}$, and this
might seem like a simplifying move, but it can actually make fat
commutative monoids confusing to work with: for instance, whereas
bijections
\[
X \goby{s} Y \goby{t} Z
\]
can be composed in only one possible order, permutations $s, t \in S_n$ can
be composed in two.  I will stick with `all finite sets'.

Write $\fcat{FatCommMon}$%
\glo{FatCommMon}
and $\fcat{CommMon}$ for the categories of fat
and ordinary commutative monoids, respectively.

\begin{thm}	\lbl{thm:fat-cm-eqv}
There is an isomorphism of categories 
\[
\fcat{FatCommMon} \iso
\fcat{CommMon}.
\]
\end{thm}
\begin{proof}
We show that both sides are isomorphic to the category $\fcat{UCommMon}$ of
unbiased commutative monoids.  By definition, an \demph{unbiased
commutative monoid}%
\index{monoid!commutative!unbiased}
is a set $A$ equipped with an $n$-ary addition
operation
\[
\begin{array}{rcl}
A^n			&\go		&A			\\
(a_1, \ldots, a_n)	&\goesto	&(a_1 + \cdots + a_n)
\end{array}
\]
for each $n\in\nat$, satisfying the three axioms displayed in
Example~\ref{eg:opd-terminal} (written there with $\cdot$ instead of $+$
and $x$'s instead of $a$'s).  We have $\fcat{UCommMon} \iso
\fcat{CommMon}$, easily.

Given a fat commutative monoid $(A, \sum)$, define an unbiased
commutative monoid structure $+$ on $A$ by
\[
(a_1 + \cdots + a_n) 
= 
\sum_{x\in \{1, \ldots, n\} } a_x.
\]
The first two axioms for an unbiased commutative monoid follow from the two
axioms for a fat commutative monoid, and the third follows from
Lemma~\ref{lemma:fat-cm}.

Conversely, take an unbiased commutative monoid $(A,+)$ and define a fat
commutative monoid structure $\sum$ on $A$ as follows.  For any finite set
$X$, let $n_X\in\nat$ be the cardinality of $X$ and choose a bijection $t_X:
\{1, \ldots, n_X\} \goiso X$; then define $\sum_X: A^X \go A$ by
\[
\sum_{x\in X} a_x = (a_{t_X(1)} + \cdots + a_{t_X(n_X)}).
\]
By commutativity, this definition is independent of the choice of $t_X$.
Clearly the axioms for a fat commutative monoid are satisfied.

It is straightforward to check that these two processes are mutually
inverse and extend to an isomorphism of categories.
% Lemma~\ref{lemma:fat-cm} is again useful here.  
\done
\end{proof}%
\index{monoid!commutative|)}

\section{Symmetric multicategories}
\lbl{sec:sym-multis}%
\index{multicategory!symmetric|(}

As in the previous section, we reformulate a notion of symmetric structure
by moving from finite sequences to finite families: so in a `fat
symmetric multicategory', maps look like
\[
(a_x)_{x\in X} \goby{\theta} b
\]
rather than
\[
a_1, \ldots, a_n \goby{\theta} b.
\]

\begin{defn}	\lbl{defn:fat-sm}
A \demph{fat symmetric multicategory}%
\index{multicategory!symmetric!fat}
$A$ consists of
\begin{itemize}
\item a set $A_0$,%
\glo{A0fat}
whose elements are called the \demph{objects} of $A$
\item for each finite set $X$, family $(a_x)_{x\in X}$ of objects, and
object $b$, a set $C((a_x)_{x\in X}; b)$,%
\glo{fathomset}
whose elements $\theta$ are called
\demph{maps} in $A$ and written
\[
(a_x)_{x\in X} \goby{\theta} b
\]
\item for each function $s: X \go Y$ between finite sets, family
$(a_x)_{x\in X}$ of objects, family $(b_y)_{y\in Y}$ of objects, and object
$c$, a function
\[
A((b_y)_{y\in Y}; c)
\times
\prod_{y\in Y} A((a_x)_{x\in s^{-1}\{y\}}; b_y)
\go
A((a_x)_{x\in X}; c),
\]
called \demph{composition} and written
\[
(\phi, (\theta_y)_{y\in Y}) \goesto \phi \of (\theta_y)_{y\in Y}%
\glo{fatcomp}
\]
\item for each one-element set $X$ and object $a$, an \demph{identity} map
\[
1_a^X \in A((a)_{x\in X}; a),%
\glo{fatids}
\]
\end{itemize}
satisfying
\begin{itemize}
\item associativity: if $X \goby{s} Y \goby{t} Z$ are functions between
finite sets and 
\[
(a_x)_{x\in s^{-1}\{y\}} \goby{\theta_y} b_y,
\diagspace
(b_y)_{y\in t^{-1}\{z\}} \goby{\phi_z} c_z,
\diagspace
(c_z)_{z\in Z} \goby{\psi} d
\]
are maps in $A$ (for $y\in Y$, $z\in Z$), then
\[
(\psi \of (\phi_z)_{z\in Z}) \of (\theta_y)_{y\in Y}
=
\psi \of (\phi_z \of (\theta_y)_{y\in t^{-1}\{z\}})_{z\in Z}
\]
\item left identity axiom: if $X$ is a finite set, $Y$ a one-element set,
and $\theta: (a_x)_{x\in X} \go b$ a map in $A$, then
\[
1_b^Y \of (\theta)_{y\in Y} = \theta
\]
\item right identity axiom: if $X$ is a finite set and $\theta: (a_x)_{x\in
X} \go b$ a map in $A$, then 
\[
\theta \of (1_{a_x}^{\{x\}})_{x\in X} = \theta.
\]
\end{itemize}
A \demph{map} $f: A \go A'$ of fat symmetric multicategories consists of
\begin{itemize}
\item a function $f: A_0 \go A'_0$ 
\item for each finite set $X$, family $(a_x)_{x\in X}$ of objects of $A$,
and object $b$ of $A$, a function
\[
A( (a_x)_{x\in X}; b )
\go
A'( (fa_x)_{x\in X}; fb ),
\]
also written as $f$,
\end{itemize}
such that
\begin{itemize}
\item $f(\phi \of (\theta_y)_{y\in Y}) = f(\phi) \of (f(\theta_y))_{y\in Y}$
whenever these composites make sense
\item $f(1_a^X) = 1_{fa}^X$ whenever $a$ is an object of $A$ and $X$ is a
one-element set.   
\end{itemize}
This defines a category $\fcat{FatSymMulticat}$.%
\glo{FatSymMulticat}
\end{defn}

As promised, this definition avoids the delicate symmetry axioms present in
the traditional version.  The following lemma, analogous to
Lemma~\ref{lemma:fat-cm}, shows that the symmetric group actions really are
hiding in there.
\begin{lemma}	\lbl{lemma:fat-sm}
Let $A$ be a fat symmetric multicategory.  Then any bijection $s: X \go Y$
between finite sets and map $\phi: (b_y)_{y\in Y} \go c$ in $A$ give rise
to a map $\phi\cdot s:%
\glo{fatcdot}
(b_{s(x)})_{x\in X} \go c$ in $A$.  This
construction satisfies
\[
\phi\cdot (s\of r) = (\phi\cdot s)\cdot r,
\diagspace
\theta \cdot 1_Y = \theta,
\]
where $W \goby{r} X \goby{s} Y$ in the first equation.  Moreover, if $f$ is a
map of fat symmetric multicategories then $f(\phi\cdot s) = f(\phi)\cdot s$
whenever these expressions make sense.  
\end{lemma}
\begin{proof}
Take $s$ and $\phi$ as in the statement.  Define a family $(a_x)_{x\in X}$
by $a_x = b_{s(x)}$.  For each $y\in Y$ we have the map
\[
1_{b_y}^{s^{-1}\{y\}}: 
(b_y)_{x\in s^{-1}\{y\}} \go b_y;
\]
but $b_y = a_{s^{-1}(y)}$, so the domain of this map is $(a_x)_{x\in
s^{-1}\{y\}}$.  We may therefore define
\[
\phi\cdot s = 
\phi\of \left( 1_{b_y}^{s^{-1}\{y\}} \right)_{y\in Y}:
(a_x)_{x\in X} \go c.
\]
The two equations follow from the associativity and identity axioms
for a fat symmetric multicategory, respectively.  That maps of fat
symmetric multicategories preserve $\cdot$ is immediate from the
definitions. 
\done
\end{proof}
This lemma is very useful in the proof below that fat and ordinary
symmetric multicategories are essentially the same.  Like
Lemma~\ref{lemma:fat-cm}, it also allows us to take notational liberties.
If, for example, we have a family of objects $a_{v,w,x}$ indexed over
finite sets $V$, $W$ and $X$, then we may safely speak of `maps
\[
(a_{v,w,x})_{v\in V, w\in W, x\in X} \go b
\]
in $A$'; it does not matter whether the indexing set in the domain is
meant to be $(V\times W)\times X$ or $V\times (W\times X)$ (or some other
3-fold product), as the canonical isomorphism between them induces a
canonical isomorphism of hom-sets,
\[
A((a_{v,w,x})_{((v,w),x) \in (V\times W) \times X} ; b)
\goiso
A((a_{v,w,x})_{(v,(w,x)) \in V\times (W \times X)} ; b).
\]

\begin{example}
A \demph{fat symmetric operad}%
\index{operad!symmetric!fat}
$P$ is, of course, a fat symmetric
multicategory with only one object; it consists of a set $P(X)$ for each
finite set $X$, together with composition and identity operations.  

Many well-known examples of symmetric operads are naturally regarded as fat
symmetric operads.  For instance, there is a fat symmetric operad $\ldisks$
where an element of $\ldisks(X)$ is an $X$-indexed family $(\alpha_x)_{x\in
X}$ of disjoint little%
\index{operad!little disks}
disks inside the unit disk
(compare~\ref{eg:opd-little-disks}).  Or, for any set $S$ there is a fat
symmetric operad $\END(S)$%
\index{operad!endomorphism}%
\index{endomorphism!symmetric operad}
defined by $ (\END(S))(X) = \Set(S^X, S)$
(compare~\ref{eg:opd-End}).  Or, there is a fat symmetric operad $O$ in
which $O(X)$ is the set of total orders%
\index{order!operad of orders}\index{operad!orders@of orders}
 on $X$, with composition done
lexicographically; under the equivalence we are about to establish, it
corresponds to the ordinary operad $\SymOpd$ of
symmetries.%
\index{operad!symmetries@of symmetries}
(This last example appeared in Beilinson%
\index{Beilinson, Alexander}
and Drinfeld%
\index{Drinfeld, Vladimir}
\cite[1.1.4]{BeDr}
and~\ref{eg:opd-Sym} above.)
\end{example}

\begin{thm}	\lbl{thm:fat-sm-eqv}
There is a canonical equivalence of categories
\[
\ovln{\blank}: \fcat{FatSymMulticat} \goby{\eqv} \fcat{SymMulticat}.
\]
\end{thm}

\begin{proof}
We define the functor $\ovln{\blank}$ and show that it is full, faithful
and (genuinely) surjective on objects.  Details are omitted.

To define $\ovln{\blank}$ on objects, take a fat symmetric multicategory
$A$.  The symmetric multicategory $\ovln{A}$ has the same objects as $A$
and hom-sets 
\[
\ovln{A} (a_1, \ldots, a_n; b) 
=
A((a_x)_{x\in [1,n]}; b),
\]
where for $m,n\in \nat$ we write 
\[
[m,n] = \{ l \in \nat \such m \leq l \leq n \}.
\]
For composition, take maps
\[
\begin{array}{c}
a_1^1, \ldots, a_1^{k_1} \goby{\theta_1} b_1,
\diagspace \ldots, \diagspace 
a_n^1, \ldots, a_n^{k_n} \goby{\theta_n} b_n,	\\
b_1, \ldots, b_n \goby{\phi} c
\end{array}
\]
in $\ovln{A}$.  Define objects $a_1, \ldots, a_{k_1 + \cdots + k_n}$ by the
equation of formal sequences
\[
(a_1, \ldots, a_{k_1 + \cdots + k_n}) 
=
(a_1^1, \ldots, a_1^{k_1}, \ldots, a_n^1, \ldots, a_n^{k_n}).
\]
For each $x\in [1,n]$ there is an obvious bijection
\[
t_x: 
[k_1 + \cdots + k_{x-1} + 1, k_1 + \cdots + k_{x-1} + k_x]
\goiso
[1, k_x]
\]
defined by subtraction: so the map 
\[
(a_{k_1 + \cdots + k_{x-1} + y})_{y\in [1,k_x]}
=
(a_x^y)_{y\in [1,k_n]}
\goby{\theta_x}
b_x
\]
in $A$ gives rise to a map
\[
(a_z)_{z\in [k_1 + \cdots + k_{x-1} + 1, k_1 + \cdots + k_{x-1} + k_x]}
\goby{\theta_x \cdot t_x}
b_x
\]
in $A$, by Lemma~\ref{lemma:fat-sm}.  It now makes sense to define
composition in $\ovln{A}$ by
\[
\phi \of (\theta_1, \ldots, \theta_n)
=
\phi \of (\theta_x \cdot t_x)_{x\in [1,n]},
\]
since the domain of this map is $(a_z)_{z\in [1, k_1 + \cdots + k_n]}$.
Identities in $\ovln{A}$ are easier: for $a\in A$, put $1_a =
1_a^{[1,1]}$.  The structure $\ovln{A}$ just defined really is a symmetric
multicategory, as is straightforward to prove with the aid of
Lemma~\ref{lemma:fat-sm}.  

The definition of the functor $\ovln{\blank}$ on morphisms and the proof
of functoriality are also straightforward.

Now we show that $\ovln{\blank}$ is surjective on objects.  For each finite
set $X$, choose a bijection $s_X: [1,n_X] \go X$, where $n_X =
\mr{card}(X)$.  In the case that $X=[m+1,m+n]$ for some $m,n\in\nat$,
choose $s_X$ to be the obvious bijection (add $m$).  Let $C$ be a symmetric
multicategory; our task is to define a fat symmetric multicategory $A$ such
that $\ovln{A} = C$.  (This is the uphill direction and is bound to require
more work.)  We define the objects of $A$ to be the objects of $C$.  If
$(a_x)_{x\in X}$ is a finite family of objects then we define
\[
A((a_x)_{x\in X}; b)
=
C(a_{s_X(1)}, \ldots, a_{s_X(n_X)}; b).
\]
If $X$ is a one-element set and $a$ an object then we put
\[
1_a^X = 1_a \in C(a;a) = A((a)_{x\in X}; a).
\]
The definition of composition is clear in principle but fiddly in practice,
hence omitted.  The idea is that a composite in $A$ can almost be defined
as a composite in $C$, but because the bijections $s_X$ were chosen at
random, we have to apply a symmetry after composing in $C$---the unique
symmetry that makes the domain come out right.  The axioms for the fat
symmetric multicategory $A$ then follow, with some effort, from the
ordinary symmetric multicategory axioms on $C$.

We also have to show that $\ovln{A} = C$.  Certainly their object-sets are
equal, and their hom-sets are equal because
\[
\ovln{A}(a_1, \ldots, a_n; b)
=
A((a_x)_{x\in [1,n]}; b)
=
C(a_{s_{[1,n]}(1)}, \ldots, a_{s_{[1,n]}(n)}; b)
\]
and we chose $s_{[1,n]}$ to be the identity.  Composition in $\ovln{A}$ was
defined using the obvious bijections
\[
[m+1, m+k] \goiso [1,k]
\]
for certain values of $m$ and $k$, and to show that it coincides with
composition in $A$ we use the fact that $s_{[m+1,m+k]}$ was also chosen to
be the obvious bijection.

Next, $\ovln{\blank}$ is full.  Let $A$ and $A'$ be fat symmetric
multicategories and $h: \ovln{A} \go \ovln{A'}$ a map of ordinary symmetric
multicategories; we define a map $f: A \go A'$ such that $\ovln{f} = h$.
On objects, $f(a) = h(a)$.  Given a map 
$
\theta: (a_x)_{x\in X} \go b
$
in $A$, choose a bijection $s: [1,n] \go X$, where $n=\mr{card}(X)$.  Then
we have the map
\[
\theta\cdot s: (a_{s(y)})_{y\in [1,n]} \go b
\]
in $A$, that is, we have 
\[
\theta\cdot s: a_{s(1)}, \ldots, a_{s(n)} \go b
\]
in $\ovln{A}$; so we obtain the map
\[
h(\theta\cdot s): fa_{s(1)}, \ldots, fa_{s(n)} \go fb
\]
in $\ovln{A'}$, that is, 
\[
h(\theta\cdot s): (fa_{s(y)})_{y\in [1,n]} \go fb
\]
in $A'$.  It therefore makes sense to define
\[
f(\theta) = h(\theta\cdot s) \cdot s^{-1}: (fa_x)_{x\in X} \go fb.
\]
This definition is independent of the choice of $s$: note that any other
choice is of the form $s\of \sigma$ for some $\sigma\in S_n$, then use
Lemma~\ref{lemma:fat-sm} and the fact that $h$ preserves symmetric group
actions.  It is straightforward to check that $f$ really is a map of fat
symmetric multicategories and that $\ovln{f} = h$.  

Finally, $\ovln{\blank}$ is faithful.  Let $A \parpair{f}{g} A'$ be a pair
of maps of fat symmetric multicategories satisfying $\ovln{f} = \ovln{g}$.
Certainly $f$ and $g$ agree on objects.  Given $\theta: (a_x)_{x\in X} \go
b$ in $A$, choose a bijection $s: [1,n] \go X$, where $n=\mr{card}(X)$;
then we have a map 
\[
\theta\cdot s: a_{s(1)}, \ldots, a_{s(n)} \go b
\]
in $\ovln{A}$.  Using Lemma~\ref{lemma:fat-sm},
\[
f(\theta)
=
f(\theta\cdot s \cdot s^{-1})
=
f(\theta\cdot s) \cdot s^{-1}
=
\ovln{f}(\theta\cdot s) \cdot s^{-1},
\]
and similarly $g(\theta) = \ovln{g}(\theta\cdot s) \cdot s^{-1}$, so
$f(\theta) = g(\theta)$. 
\done
\end{proof}%
\index{multicategory!symmetric|)}

\begin{notes}

Something very close to the notion of fat symmetric multicategories
appeared in Beilinson%
\index{Beilinson, Alexander}
and Drinfeld%
\index{Drinfeld, Vladimir}
\cite[\S 1.1]{BeDr}, under the name of
`pseudo-tensor%
\index{pseudo-tensor category}%
\index{category!pseudo-tensor}
categories'.  The idea has probably appeared elsewhere too.
Beilinson and Drinfeld insisted that the domain $(a_x)_{x\in X}$ of a map
should be a \emph{non-empty} finite family of objects, and correspondingly
that the functions called $s: X\go Y$ in Definition~\ref{defn:fat-sm}
should be surjective.  This amounts to excluding the possibility of
nullary%
\index{nullary!arrow}
maps, which we have no reason to do.  They also handled one-member families
slightly differently.

\end{notes}

\chapter{Coherence for Monoidal Categories}
\lbl{app:unbiased}%
\index{Sigma-monoidal category@$\Sigma$-monoidal category|(}%
\index{monoidal category!Sigma-@$\Sigma$-|(}%
\index{coherence!monoidal categories@for monoidal categories|(}

\noindent
Here we prove the `descriptive' coherence theorems for unbiased and
classical monoidal categories and functors
(\ref{thm:diag-coh-umc},~\ref{thm:diag-coh-mc}).

Unbiased and classical monoidal categories were defined
concretely~(\ref{defn:mon-cat},~\ref{defn:lax-mon-cat}), whereas
$\Sigma$-monoidal categories were defined
abstractly~(\ref{defn:Sigma-mon-cat}).  To carry out the comparisons, we
put the definition of $\Sigma$-monoidal category into concrete terms.  Let
$\Sigma \in \Set^\nat$.  Unwinding the abstract definition, we find that a
$\Sigma$-monoidal category is a triple $A = (A, \otimes, \delta)$%
\glo{cohdelta}
consisting of
\begin{itemize}
\item a small category $A$
\item a functor $\otimes_\tau: A^n \go A$%
\glo{otimestree}
for each $n\in\nat$,
$\tau\in (F\Sigma)(n)$
\item an isomorphism
\[
(\delta_{\tau, \tau'})_{a_1, \ldots, a_n}:
\otimes_\tau(a_1, \ldots, a_n)
\goiso
\otimes_{\tau'}(a_1, \ldots, a_n)
\]
in $A$ for each $n\in\nat, \tau, \tau' \in (F\Sigma)(n), a_i \in A$
(usually just written $\delta_{\tau, \tau'}$)
\end{itemize}
satisfying 
\begin{description}
\item[\astyle{MC1}]	% \lbl{ax:obj:nat}
$(\delta_{\tau, \tau'})_{a_1, \ldots, a_n}$ is natural in $a_1,
\ldots, a_n \in A$, for each $n\in\nat, \tau, \tau' \in (F\Sigma)(n)$
\item[\astyle{MC2}] 	% \lbl{ax:obj:func}
$\delta_{\tau', \tau''} \of \delta_{\tau, \tau'} = \delta_{\tau, \tau''}$
and $1 = \delta_{\tau, \tau}$, for all $n\in\nat, \tau, \tau', \tau'' \in
(F\Sigma)(n)$
\item[\astyle{MC3}]	% \lbl{ax:obj:alg-obj}
$\otimes_{\tau \sof (\tau_1, \ldots, \tau_n)} = (A^{k_1 + \cdots + k_n}
\goby{\otimes_{\tau_1} \times\cdots\times \otimes_{\tau_n}} A^n
\goby{\otimes_\tau} A)$ for all $n, k_i \in\nat, \tau\in (F\Sigma)(n),
\tau_i \in (F\Sigma)(k_i)$; and $1_A = \otimes_\utree$
\item[\astyle{MC4}] 	% \lbl{ax:obj:alg-map}
the diagram
\begin{diagram}[width=2em]
\begin{array}[t]{l}
\otimes_\tau ( \otimes_{\tau_1}(a_1^1, \ldots, a_1^{k_1}), \ldots, \\
\otimes_{\tau_n}(a_n^1, \ldots, a_n^{k_n}))			
\end{array}
&
\rEquals							&
\otimes_{\tau\sof(\tau_1, \ldots, \tau_n)} ( a_1^1, \ldots, a_n^{k_n})	\\
\dTo<{\otimes_\tau(\delta_{\tau_1,\tau'_1}, \ldots, 
\delta_{\tau_n,\tau'_n})}					&
								&
									\\
\begin{array}{l}
\otimes_\tau ( \otimes_{\tau'_1}(a_1^1, \ldots, a_1^{k_1}), \ldots, \\
\otimes_{\tau'_n}(a_n^1, \ldots, a_n^{k_n}))			
\end{array}
&
								&
\dTo>{\delta_{\tau\sof(\tau_1, \ldots, \tau_n), 
\tau'\sof(\tau'_1, \ldots, \tau'_n)}}					\\
\dTo<{\delta_{\tau,\tau'}}					&
								&
									\\
\begin{array}[b]{l}
\otimes_{\tau'} ( \otimes_{\tau'_1}(a_1^1, \ldots, a_1^{k_1}), \ldots, \\
\otimes_{\tau'_n}(a_n^1, \ldots, a_n^{k_n}))			
\end{array}
&
\rEquals							&
\otimes_{\tau'\sof(\tau'_1, \ldots, \tau'_n)} ( a_1^1, \ldots, a_n^{k_n})\\
\end{diagram}
commutes, for all $n, k_i \in\nat, \tau,\tau'\in (F\Sigma)(n),
\tau_i,\tau'_i\in (F\Sigma)(k_i), a_i^j\in A$.
\end{description}
A lax monoidal functor $(P, \pi): (A,\otimes,\delta) \go
(A',\otimes,\delta)$ (where, in an abuse of notation, we write
$\otimes$ for the tensor and $\delta$ for the coherence maps in both $A$
and $A'$) consists of
\begin{itemize}
\item a functor $P: A \go A'$
\item a map
\[
(\pi_\tau)_{a_1, \ldots, a_n}:
\otimes_\tau(Pa_1, \ldots, Pa_n) \go P\otimes_\tau (a_1, \ldots, a_n)
\]
(usually just written $\pi_\tau$) for each $n\in\nat, \tau\in (F\Sigma)(n),
a_i\in A$ 
\end{itemize}
satisfying
\begin{description}
\item[\astyle{MF1}] 	% \lbl{ax:map:nat}
$(\pi_\tau)_{a_1, \ldots, a_n}$ is natural in $a_1, \ldots, a_n \in A$
\item[\astyle{MF2}]	% \lbl{ax:map:nat-tau}
the diagram
\begin{diagram}[size=2em]
\otimes_\tau (Pa_1, \ldots, Pa_n)	&
\rTo^{\delta_{\tau,\tau'}}		&
\otimes_{\tau'}(Pa_1, \ldots, Pa_n)	\\
\dTo<{\pi_\tau}				&
					&
\dTo>{\pi_{\tau'}}			\\
P\otimes_\tau (a_1, \ldots, a_n)	&
\rTo_{P\delta_{\tau,\tau'}}		&
P\otimes_{\tau'} (a_1, \ldots, a_n)	\\
\end{diagram}
commutes, for all $n\in\nat, \tau,\tau'\in (F\Sigma)(n), a_i\in A$
\item[\astyle{MF3}] 	% \lbl{ax:map:coh}
the diagram
\begin{diagram}[width=2em]
\begin{array}[t]{l}
\otimes_\tau (\otimes_{\tau_1} (Pa_1^1, \ldots, Pa_1^{k_1}), \ldots,\\
\otimes_{\tau_n} (Pa_n^1, \ldots, Pa_n^{k_n}))		
\end{array}
&
\rEquals						&
\otimes_{\tau\sof(\tau_1, \ldots, \tau_n)} (Pa_1^1, \ldots, Pa_n^{k_n})	\\
\dTo<{\otimes_\tau (\pi_{\tau_1}, \ldots, \pi_{\tau_n})}&
							&
							\\
\begin{array}{l}
\otimes_\tau ( P\otimes_{\tau_1}(a_1^1, \ldots, a_1^{k_1}), \ldots,\\
P\otimes_{\tau_n}(a_n^1, \ldots, a_n^{k_n}) )		
\end{array}
&
							&
\dTo>{\pi_{\tau\sof(\tau_1, \ldots, \tau_n)}}		\\
\dTo<{\pi_\tau}						&
							&
							\\
\begin{array}[b]{l}
P\otimes_\tau (\otimes_{\tau_1}(a_1^1, \ldots, a_1^{k_1}), \ldots,\\
\otimes_{\tau_n}(a_n^1, \ldots, a_n^{k_n}))		
\end{array}
&
\rEquals						&
P\otimes_{\tau\sof(\tau_1, \ldots, \tau_n)} (a_1^1, \ldots, a_n^{k_n})\\
\end{diagram}
commutes for all $n, k_i\in\nat, \tau\in (F\Sigma)(n), \tau_i\in
(F\Sigma)(k_i), a_i^j\in A$, and the diagram
\begin{diagram}[size=2em]
Pa	&\rEquals	&\otimes_\utree Pa	\\
\dTo<1	&		&\dTo>{\pi_1}		\\
Pa	&\rEquals	&P\otimes_\utree a	\\
\end{diagram}
commutes for all $a\in A$ (where in the expression `$\pi_1$', the $1$ is
the unit of the operad $F\Sigma$).
\end{description}
A weak (respectively, strict) monoidal functor is a lax monoidal functor
$(P, \pi)$ in which all the $\pi_\tau$'s are isomorphisms (respectively,
identities).  

In the cases at hand, the object $\Sigma$ of $\Set^\nat$ is either the
terminal object $1$ (for unbiased monoidal categories) or the object
$\Sigma_\mr{c}$ defined in Theorem~\ref{thm:diag-coh-mc} (for classical
monoidal categories).  The \Set-operad $F\Sigma$ is then, respectively,
either the operad $\tr$ of all trees~(\ref{eg:opd-of-trees}) or the operad
$\ctr$ of classical trees~(\ref{eg:opd-of-cl-trees}).

Sections~\ref{sec:app-UMC} and~\ref{sec:app-MC} prove coherence for
unbiased and classical monoidal categories, respectively.  The unbiased
case contains some scary-looking expressions but is completely
straightforward: one only needs to be awake, not clever.  By contrast, the
classical case requires guile, cunning and trickery---attributes not
displayed here, but vital to the proofs of the coherence theorems upon
which we rely.

All of the results that we prove continue to hold in the more general
setting of $\Sigma$-bicategories~(\ref{sec:notions-bicat}).%
\index{bicategory!Sigma-@$\Sigma$-}%
\index{Sigma-bicategory@$\Sigma$-bicategory}
 The proofs
need only superficial changes of the kind described
in~\ref{sec:notions-bicat}; for instance, the category $A$ becomes a
\Cat-graph $B$ and tensor becomes composition.

\section{Unbiased monoidal categories}
\lbl{sec:app-UMC}%
\index{monoidal category!unbiased|(}

The plan is to define a functor
\[
J: 1\hyph\MClax \go \UMClax,
\]
to prove that it is an isomorphism (by showing in turn that it is injective
on objects, surjective on objects, faithful, and full), and then to prove
that it restricts to isomorphisms
\[
1\hyph\MCwk \goiso \UMCwk,
\diagspace
1\hyph\MCstr \goiso \UMCstr
\]
at the weak and strict levels.

\index{tree|(}
Recall from~\ref{eg:opd-of-trees} that $\utree\in\tr(1)$ denotes
the `null' or identity tree, and that $\nu_n \in \tr(n)$ denotes the
simplest $n$-leafed tree 
% \drmk{pic of corolla}.  
$\begin{centredpic}
\begin{picture}(3,2)(-1.5,0)
% lower layer
\put(0,0){\line(0,1){1}}
\cell{0}{1}{c}{\vx}
% upper layer
\put(0,1){\line(-3,2){1.5}}
\cell{0}{1.8}{c}{\cdots}
\put(0,1){\line(3,2){1.5}}
\end{picture}
\end{centredpic}$.
(Recall also that $\nu_1 =
\setlength{\unitlength}{1em}
\begin{array}{c}
\begin{picture}(0,2)(0,0)
\put(0,0){\line(0,1){1}}
\cell{0}{1}{c}{\vx}
\put(0,1){\line(0,1){1}}
\end{picture}
\end{array}
\neq
\setlength{\unitlength}{1em}
\begin{array}{c}
\begin{picture}(0,1)(0,0)
\put(0,0){\line(0,1){1}}
\end{picture}
\end{array}$.) 
We will be composing trees---that is, composing in the operad $\tr$---and
in particular we will often use the composite tree
\[
\nu_n \of (\nu_{k_1}, \ldots, \nu_{k_n}) = 
% \drmk{picture}.
\begin{centredpic}
\begin{picture}(6,3)(-3,0)
% bottom layer
\put(0,0){\line(0,1){1}}
% middle layer
\cell{0}{1}{c}{\vx}
\put(0,1){\line(-2,1){2}}
\put(0,1){\line(2,1){2}}
\cell{0}{1.8}{c}{\cdots}
% top layer
\cell{-2}{2}{c}{\vx}
\put(-2,2){\line(-1,1){1}}
\put(-2,2){\line(1,1){1}}
\cell{-2}{2.8}{c}{\cdots}
\cell{2}{2}{c}{\vx}
\put(2,2){\line(-1,1){1}}
\put(2,2){\line(1,1){1}}
\cell{2}{2.8}{c}{\cdots}
\end{picture}
\end{centredpic}.
% \dr{36}{18}.
\]

Our first task is to define a functor $J: 1\hyph\MClax \go \UMClax$.
\begin{description}
\item[On objects] Let $(A, \otimes, \delta)$ be an object of
$1\hyph\MClax$.  The unbiased monoidal category $J(A, \otimes, \delta)$ is
given by taking the underlying category to be $A$, the $n$-fold tensor
$\otimes_n: A^n \go A$ to be $\otimes_{\nu_n}$, and the coherence maps to
be
\begin{eqnarray*}
\lefteqn{\gamma_{((a_1^1, \ldots, a_1^{k_1}), \ldots, (a_n^1, \ldots, a_n^{k_n}))}
=
(\delta_{\nu_n \sof (\nu_{k_1}, \ldots, \nu_{k_n}), 
\nu_{k_1 + \cdots + k_n}})_{a_1^1, \ldots, a_n^{k_n}}:}	\\
&((a_1^1 \otimes\cdots\otimes a_1^{k_1}) \otimes\cdots\otimes
(a_n^1 \otimes\cdots\otimes a_n^{k_n}))
\go
(a_1^1 \otimes\cdots\otimes a_n^{k_n})
\end{eqnarray*}
and
\[
\iota_a = (\delta_{\utree, \nu_1})_a: a \go (a).
\]
\item[On maps] Let $(P, \pi): A \go A'$ be a map in $1\hyph\MClax$.  The
lax monoidal functor $J(P, \pi): J(A) \go J(A')$ is given by taking the
same underlying functor $P$, and by taking the coherence map
\[
\pi_{a_1, \ldots, a_n}:
(Pa_1 \otimes \cdots \otimes Pa_n) 
\go
P(a_1 \otimes\cdots\otimes a_n)
\]
(re-using the letter $\pi$, in another slight abuse) to be
$(\pi_{\nu_n})_{a_1, \ldots, a_n}$. 
\end{description}

\begin{lemma}
This defines a functor $J: 1\hyph\MClax \go \UMClax$.
\end{lemma}

\begin{proof} 
We have to check three things:
\begin{itemize}
\item $J(A,\otimes,\delta)$ as defined above really is an unbiased
monoidal category---in other words, the axioms in
Definition~\ref{defn:lax-mon-cat} hold.  Naturality and invertibility of
$\gamma$ and $\iota$  follow from the same properties for $\delta$. The
associativity axiom holds because both routes around the square are
\[
(\delta_{ \nu_n\sof (\nu_{m_1}, \ldots, \nu_{m_n}) \sof (\nu_{k_1^1},
\ldots, \nu_{k_n^{m_n}}), \nu_{k_1^1 + \cdots + k_n^{m_n}}})_%
{a_{1,1,1}, \ldots, a_{n,m_n,k_n^{m_n}}},
\]
as can be shown from axioms~\astyle{MC2} and~\astyle{MC4}.  The identity
axioms hold by similar reasoning.
\item $J(P, \pi)$ as defined above really is a lax monoidal functor between
unbiased monoidal categories (Definition~\ref{defn:u-lax-mon-ftr}).
Naturality of $\pi_{a_1, \ldots, a_n}$ in the $a_i$'s follows from the same
naturality for $\pi_{\nu_n}$.  The coherence axioms can be deduced from
axioms \astyle{MF1}--\astyle{MF3}.
\item $J$ preserves composition and identities.  This is trivial.
\done
\end{itemize}
\end{proof}

\begin{lemma}	\lbl{lemma:u-inj}
The functor $J: 1\hyph\MClax \go \UMClax$ is injective on objects.
\end{lemma}

\begin{proof}
Suppose that $(A, \otimes, \delta)$ and $(A', \otimes', \delta')$ are
$1$-monoidal categories with $J(A, \otimes, \delta) = J(A', \otimes',
\delta')$.  (Just for once, the tensor of $A'$ is written $\otimes'$ rather
than $\otimes$, and similarly $\delta'$ rather than $\delta$.)  Then:
\begin{itemize}
\item $A = A'$ immediately.
\item $\otimes_\tau = \otimes'_\tau$ for all $n\in\nat$ and
$\tau\in\tr(n)$.  This is proved by induction on $\tau$, using the
description of $\tr$ in~\ref{eg:opd-of-trees}.  Since there will be many
inductions on trees in this appendix, I will write this one out in full and
leave the others to the virtuous reader.  First, $\otimes_\utree = 1_A =
\otimes'_\utree$ by~\astyle{MC3}.  Second, take $\tau_1\in\tr(k_1), \ldots,
\tau_n\in\tr(k_n)$. Then
\begin{eqnarray*}
\otimes_{(\tau_1, \ldots, \tau_n)}	&=	&
\otimes_{\nu_n \sof (\tau_1, \ldots, \tau_n)}	\\
&=&\otimes_{\nu_n} \of (\otimes_{\tau_1} \times\cdots\times
\otimes_{\tau_n}),   
\end{eqnarray*}
again by~\astyle{MC3}.  But $\otimes_{\nu_n} = \otimes'_{\nu_n}$ since
$J(A, \otimes, \delta) = J(A', \otimes', \delta')$, and $\otimes_{\tau_i} =
\otimes'_{\tau_i}$ by inductive hypothesis, so $\otimes_{(\tau_1, \ldots,
\tau_n)} = \otimes'_{(\tau_1, \ldots, \tau_n)}$, completing the induction.
\item $\delta_{\tau,\tau'} = \delta'_{\tau,\tau'}$ for all $\tau, \tau' \in
\tr(n)$.  Since $\delta_{\tau,\tau'} = \delta_{\tau',\nu_n}^{-1} \of
\delta_{\tau,\nu_n}$, it is enough to prove this in the case
$\tau'=\nu_n$; and that is done by another short induction on $\tau$.
\done
\end{itemize}
\end{proof}

\begin{lemma}
The functor $J: 1\hyph\MClax \go \UMClax$ is surjective on objects.
\end{lemma}

\begin{proof}
Take an unbiased monoidal category $(A, \otimes, \gamma, \iota)$.  Attempt
to define a $1$-monoidal category $(A, \otimes, \delta)$ as follows:
\begin{itemize}
\item The underlying category $A$ is the same.
\item The tensor $\otimes_\tau: A^n \go A$, for $\tau\in\tr(n)$, is defined
inductively on $\tau$ by $\otimes_\utree = 1_A$ and 
\[
\otimes_{(\tau_1, \ldots, \tau_n)} = 
(A^{k_1 + \cdots + k_n} 
\goby{\otimes_{\tau_1} \times\cdots\times \otimes_{\tau_n}}
A^n \goby{\otimes_n} A).
\]
\item The coherence isomorphisms are defined by $\delta_{\tau, \tau'} =
\delta_{\tau'}^{-1} \of \delta_\tau$, where in turn
\[
\delta_\tau: \otimes_\tau(a_1, \ldots, a_n) \goiso 
(a_1 \otimes\cdots\otimes a_n)
\]
is defined by taking $\delta_\utree = \iota$ and taking $\delta_{(\tau_1,
\ldots, \tau_n)}$ to be the composite
\begin{eqnarray*}
&&
(\otimes_{\tau_1}(a_1^1, \ldots, a_1^{k_1}) 
\ \otimes\ \cdots\ \otimes\ 
\otimes_{\tau_n}(a_n^1, \ldots, a_n^{k_n}) )	\\ &
\goby{(\delta_{\tau_1} \otimes\cdots\otimes \delta_{\tau_n})} &	
((a_1^1 \otimes\cdots\otimes a_1^{k_1}) \otimes\cdots\otimes
(a_n^1 \otimes\cdots\otimes a_n^{k_n}))		\\
&
\goby{\gamma}	&	
(a_1^1 \otimes\cdots\otimes a_n^{k_n}).
\end{eqnarray*}
\end{itemize}
This does indeed satisfy the axioms for a $1$-monoidal category:
\begin{description}
\item[\astyle{MC1}, \astyle{MC2}] Immediate.
\item[\astyle{MC3}] That $\otimes_\utree = 1_A$ is immediate.  The other
equation can be proved by induction on $\tau$, or by using the fact that
$\tr$ is the free operad on $1 \in\Set^\nat$ and that
$(\Cat(A^n,A))_{n\in\nat}$ forms an operad.
\item[\astyle{MC4}] It is enough to show that this axiom is satisfied when
$\tau'=\nu_n, \tau'_1=\nu_{k_1}, \ldots, \tau'_n=\nu_{k_n}$: in other
words, that
\begin{diagram}[width=2em]
\begin{array}[t]{l}
\otimes_\tau ( \otimes_{\tau_1}(a_1^1, \ldots, a_1^{k_1}), \ldots, \\
\otimes_{\tau_n}(a_n^1, \ldots, a_n^{k_n}))			
\end{array}
&
\rEquals							&
\otimes_{\tau\sof(\tau_1, \ldots, \tau_n)} ( a_1^1, \ldots, a_n^{k_n})	\\
\dTo<{\otimes_\tau(\delta_{\tau_1}, \ldots, 
\delta_{\tau_n})}						&
								&
\dTo>{\delta_{\tau\sof(\tau_1,\ldots,\tau_n)}}				\\
\begin{array}{l}
\otimes_\tau ( (a_1^1 \otimes\cdots\otimes a_1^{k_1})
\otimes\cdots	\\
\otimes (a_n^1 \otimes\cdots\otimes a_n^{k_n}) )
\end{array}
&
								&
(a_1^1 \otimes\cdots\otimes a_n^{k_n})					\\
\dTo<{\delta_\tau}						&
								&
\uTo>{\delta_{\nu_n \sof (\nu_{k_1},\ldots,\nu_{k_n})} = \gamma}	\\
\begin{array}[b]{l}
((a_1^1 \otimes\cdots\otimes a_1^{k_1}) \otimes\cdots \\
\otimes (a_n^1 \otimes\cdots\otimes a_n^{k_n}))			
\end{array}
&
\rEquals							&
\begin{array}[b]{r}
((a_1^1 \otimes\cdots\otimes a_1^{k_1}) \otimes\cdots	\\
\otimes (a_n^1 \otimes\cdots\otimes a_n^{k_n}))				
\end{array}
\\
\end{diagram}
commutes.  This is done by induction on $\tau$, using the associativity
axiom for unbiased monoidal categories.
\end{description}
Finally, $J(A,\otimes,\delta) = (A,\otimes,\gamma,\iota)$:
\begin{itemize}
\item Clearly the underlying categories agree, both being $A$.
\item We have 
\[
\otimes_{\nu_n} 
= \otimes_{(\utree, \ldots, \utree)} 
= \otimes_n \of (\otimes_\utree \times\cdots\times \otimes_\utree)
= \otimes_n,
\]
so the tensor products agree.
\item Certainly $\delta_\utree = \iota$.  Also
\[
\delta_{\nu_k} 
= \gamma \of (\delta_\utree \otimes\cdots\otimes \delta_\utree)
= \gamma \of (\iota \otimes\cdots\otimes \iota)
= 1
\]
for any $k$, so 
\[
\delta_{\nu_n \sof (\nu_{k_1}, \ldots, \nu_{k_n})} 
= \gamma \of (\delta_{\nu_{k_1}} \otimes\cdots\otimes \delta_{\nu_{k_n}})
= \gamma.
\]
Hence the coherence maps $\iota$ and $\gamma$ also agree.
\done
\end{itemize}
\end{proof}

\begin{lemma}
The functor $J: 1\hyph\MClax \go \UMClax$ is faithful.
\end{lemma}

\begin{proof}
Suppose that $A \parpair{(P,\pi)}{(Q,\chi)} A'$ in $1\hyph\MClax$ with
$J(P,\pi) = J(Q,\chi)$.  Then:
\begin{itemize}
\item $P=Q$ immediately.
\item $\pi_\tau = \chi_\tau$ for all $\tau\in\tr(n)$, by an induction on
$\tau$ similar to the one written out in the proof of~\ref{lemma:u-inj}.
\done
\end{itemize}
\end{proof}

\begin{lemma}	\lbl{lemma:u-full}
The functor $J: 1\hyph\MClax \go \UMClax$ is full.
\end{lemma}

\begin{proof}
Let $A, A' \in 1\hyph\MClax$ and let $J(A) \goby{(P,\pi)} J(A')$ be a map
in $\UMClax$.  Attempt to define a map $A \goby{(P,\pi)} A'$ in
$1\hyph\MClax$ (where as usual we abuse notation by recycling the name
$(P,\pi)$) as follows:
\begin{itemize}
\item The underlying functor $P$ is the same.
\item The coherence maps $\pi_\tau$ are defined by induction on $\tau$.  We
take $\pi_\utree = 1$ and take
\[
\otimes_{(\tau_1, \ldots, \tau_n)} (Pa_1^1, \ldots, Pa_n^{k_n}) 
\goby{\pi_{(\tau_1, \ldots, \tau_n)}}
P\otimes_{(\tau_1, \ldots, \tau_n)} (a_1^1, \ldots, a_n^{k_n}) 
\]
to be the composite
\begin{eqnarray*}
\lefteqn{
(\otimes_{\tau_1}(Pa_1^1, \ldots, Pa_1^{k_1}) 
\ \otimes\ \cdots\ \otimes\ 
\otimes_{\tau_n}(Pa_n^1, \ldots, Pa_n^{k_n}) )}	\\ 
&
\goby{(\pi_{\tau_1} \otimes\cdots\otimes \pi_{\tau_n})}	&
\begin{array}[t]{l}
(P\otimes_{\tau_1}(a_1^1, \ldots, a_1^{k_1}) \ \otimes\ \cdots\\
\otimes\ P\otimes_{\tau_n}(a_n^1, \ldots, a_n^{k_n}))	
\end{array}
\\
&
\goby{\pi_{\otimes_{\tau_1}(a_1^1, \ldots, a_1^{k_1}), \ldots,
\otimes_{\tau_n}(a_n^1, \ldots, a_n^{k_n})}}	&
P(a_1^1 \otimes\cdots\otimes a_n^{k_n}).
\end{eqnarray*}
\end{itemize}

This does satisfy the axioms above for a lax monoidal functor between
$1$-monoidal categories:
\begin{description}
\item[\astyle{MF1}] Immediate.
\item[\astyle{MF2}] It is enough to prove this when $\tau'=\nu_n$; in
other words, that
\begin{diagram}[size=2em]
\otimes_\tau (Pa_1, \ldots, Pa_n)	&\rTo^{\delta_\tau}	&
(Pa_1 \otimes\cdots\otimes Pa_n)	\\
\dTo<{\pi_\tau}				&			&
\dTo>{\pi_{a_1,\ldots,a_n}}		\\
P\otimes_\tau(a_1, \ldots, a_n)		&
\rTo_{P\delta_\tau}			&
P(a_1 \otimes\cdots\otimes a_n)		\\
\end{diagram}
commutes.  This is done by induction on $\tau$, using the coherence
axioms for lax monoidal functors between unbiased monoidal categories.
\item[\astyle{MF3}] This is, inevitably, another induction on $\tau$.
\end{description}

Finally, $J(P,\pi) = (P,\pi)$:
\begin{itemize}
\item The underlying functors agree, both being $P$.
\item We have
\begin{eqnarray*}
(\pi_{\nu_n})_{a_1, \ldots, a_n} &=	&
\pi_{\otimes_\utree(a_1), \ldots, \otimes_\utree(a_n)}
\of
((\pi_\utree)_{a_1} \otimes\cdots\otimes (\pi_\utree)_{a_n})	\\
	&=	&
\pi_{a_1, \ldots, a_n} \of (1 \otimes\cdots\otimes 1)	\\
	&=	&
\pi_{a_1, \ldots, a_n},
\end{eqnarray*}
so the coherence maps also agree.
\done
\end{itemize}
\end{proof}%
\index{tree|)}

Theorem~\ref{thm:diag-coh-umc} now follows in the lax case: $\UMClax \iso
1\hyph\MClax$.  The weak and strict cases follow from:

\begin{lemma}	\lbl{lemma:u-restrict}
The isomorphism $J: 1\hyph\MClax \goiso \UMClax$ restricts to 
isomorphisms
\[
1\hyph\MCwk \goiso \UMCwk,
\diagspace
1\hyph\MCstr \goiso \UMCstr.
\]
\end{lemma}

\begin{proof}
Trivial.
\done
\end{proof}%
\index{monoidal category!unbiased|)}

\section{Classical monoidal categories}
\lbl{sec:app-MC}%
\index{monoidal category!classical|(}

The strategy is the same as in the unbiased setting: we define
a functor
\[
J: \Sigma_\mr{c}\hyph\MClax \go \MClax,
\]
prove it is an isomorphism, then prove that it restricts to isomorphisms at
the levels of weak and strict maps.

\index{tree!classical|(}
We will use the operad \ctr\ of classical or unitrivalent
trees~(\ref{eg:opd-of-cl-trees}).  In particular, we have trees
\[
\nu_0 = \nuzeropic \in \ctr(0), 
\diagspace
\nu_2 = \nutwopic \in \ctr(2), 
\]
the identity tree $\utree\in\ctr(1)$, and composite trees
\[
\begin{array}{c}
\nu_2 \of (\nu_2, \utree) = \assleftpic \in \ctr(3),
\diagspace
\nu_2 \of (\utree, \nu_2) = \assrightpic \in \ctr(3),\\
\nu_2 \of (\nu_0, \utree) = \lambdapic \in \ctr(1),
\diagspace
\nu_2 \of (\utree, \nu_0) = \rhopic \in \ctr(1).
\end{array}
\]

To prove our result, we will certainly want to use the fact that `all
diagrams commute' in a classical monoidal category, and some similar
statement for monoidal functors.  In other words, the proof of our unified
coherence theorem for classical monoidal categories and functors will
depend on pre-established coherence theorems.  A rather vague description
of those theorems was given in~\ref{sec:mon-cats}.  Using the language of
trees we can now be more precise.

So, let $A = (A, \otimes, I, \alpha, \lambda, \rho)$ be a classical
monoidal category.  We define for each $n\in\nat$ and $\tau\in\ctr(n)$ a
functor $\otimes_\tau: A^n \go A$.  This definition is by induction on
$\tau$ (using the description of $\ctr$ in~\ref{eg:opd-of-cl-trees}), as
follows:
\begin{itemize}
\item $\otimes_\utree: A^1 \go A$ is the canonical isomorphism
\item $\otimes_\nuzeropic: A^0 \go A$ `is' the unit object $I$ of $A$
\item if $\tau_1\in\ctr(k_1)$ and $\tau_2\in\ctr(k_2)$ then 
$\otimes_{(\tau_1, \tau_2)}$ is the composite
\[
A^{k_1 + k_2} \goby{\otimes_{\tau_1} \times \otimes_{\tau_2}}
A^2 \goby{\otimes} A.
\]
\end{itemize}
The coherence results from~\ref{sec:mon-cats} can now be stated as:
\begin{itemize}
\item Let $A$ be a classical monoidal category.  Then for each $n\in\nat$
and $\tau, \tau' \in \ctr(n)$, there is a unique natural transformation
$A^n \ctwomult{\otimes_\tau}{\otimes_{\tau'}}{} A$ built out of $\alpha$,
$\lambda$ and $\rho$.
\item Let $(P, \pi): A \go A'$ be a lax monoidal functor.  Then for each
$n\in\nat$ and $\tau, \tau' \in \ctr(n)$, there is a unique natural
transformation
\begin{diagram}
A^n		&\rTo^{\otimes_\tau}	&A		\\
\dTo<{P^n}	&\nent			&\dTo>P		\\
A'^n		&\rTo_{\otimes_{\tau'}}	&A'
\end{diagram}
built out of the coherence isomorphisms $\alpha, \lambda$ and $\rho$ of $A$ and
$A'$ and the coherence maps $\pi_{\dashbk,\dashbk}$ and $\pi_\cdot$.
\end{itemize}
I will refer to these as \demph{informal coherence}%
\index{coherence!monoidal categories@for monoidal categories!informal}
for classical monoidal
categories and functors.  They could be made precise by defining `built out
of', but I hope the reader will be content to use them in their present
imprecise form.

Our first task is to define a functor $J: \Sigma_\mr{c}\hyph\MClax \go
\MClax$.

\begin{description}
\item[On objects] Let $(A, \otimes, \delta)$ be an object of
$\Sigma_\mr{c}\hyph\MClax$, as described in the introduction to this
appendix.  Then the classical monoidal category $J(A, \otimes, \delta)$ is
given by taking the underlying category to be $A$, the tensor to be
$\otimes_\nutwopic$, the unit object to be $\otimes_\nuzeropic$, and
the coherence isomorphisms to be
\[
\alpha = \delta_{\tau_1, \tau'_1},
\diagspace
\lambda = \delta_{\tau_2, \utree},
\diagspace
\rho = \delta_{\tau_3, \utree}
\]
where
\[
\tau_1 = \assleftpic,
\diagspace
\tau'_1 = \assrightpic,
\diagspace
\tau_2 = \lambdapic,
\diagspace
\tau_3 = \rhopic.
\]
\item[On maps] Let $(P, \pi): A \go A'$ be a map in
$\Sigma_\mr{c}\hyph\MClax$, also as described in the introductory section.
Then the lax monoidal functor $J(P,\pi): J(A) \go J(A')$ is given by taking
the same underlying category $P$ and coherence maps
\[
\pi_{\dashbk, \dashbk} = \pi_{\nutwopic},
\diagspace
\pi_\cdot = \pi_{\nuzeropic}.
\] 
% %
\end{description}

\begin{lemma}
This defines a functor $J: \Sigma_\mr{c}\hyph\MClax \go \MClax$.
\end{lemma}

\begin{proof}
We have to check three things:
\begin{itemize}
\item $J(A, \otimes, \delta)$ as defined above really is a classical
monoidal category.  Naturality and invertibility of $\alpha$, $\lambda$ and
$\rho$ follow from the same properties of $\delta$.  The
pentagon~(\ref{defn:mon-cat}) commutes because each route around it is
$\delta_{\tau,\tau'}$, where
\[
\tau = 
\begin{array}{c}\epsfig{file=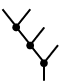}\end{array},
\diagspace
\tau' = 
\begin{array}{c}\epsfig{file=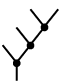}\end{array}.
\]
Similarly, the triangle commutes because both ways round it are
$\delta_{\sigma,\sigma'}$, where
\[
\sigma = 
\begin{array}{c}\epsfig{file=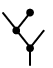}\end{array},
\diagspace
\sigma' = \nutwopic.
\]
\item $(J,\pi)$ as defined above really is a lax monoidal functor between
classical monoidal categories.  The axioms can be deduced from
\astyle{MF1}--\astyle{MF3}. 
\item $J$ preserves composition and identities.  This is trivial.%
\done
\end{itemize}
\end{proof}

\begin{lemma}	\lbl{lemma:c-inj}
The functor $J: \Sigma_\mr{c}\hyph\MClax \go \MClax$ is injective on
objects. 
\end{lemma}

\begin{proof}
Suppose that $(A, \otimes, \delta)$ and $(A', \otimes', \delta')$ are
$\Sigma_\mr{c}$-monoidal categories with $J(A, \otimes, \delta) = J(A',
\otimes', \delta')$.  Then:
\begin{itemize}
\item $A = A'$ immediately.
\item $\otimes_\tau = \otimes'_\tau$ for all $n\in\nat$ and
$\tau\in\tr(n)$.  This is proved by induction on $\tau$, using the
definition of $\ctr$ in~\ref{eg:opd-of-cl-trees}.  As in the previous
section, there will be many proofs by induction on the structure of a tree;
and as a sample of the technique in the classical case, I will write this
one out in full.  First, $\otimes_\utree = 1_A = \otimes'_\utree$
by~\astyle{MC3}.  Second, $\otimes_\nuzeropic = I = \otimes'_\nuzeropic$
since $J(A, \otimes, \delta) = J(A', \otimes', \delta')$.  Third, take
$\tau_1 \in \ctr(k_1)$ and $\tau_2 \in \ctr(k_2)$.  Then
\[
\otimes_{(\tau_1, \tau_2)} 
= \otimes_{\nutwopic \sof (\tau_1, \tau_2)}
= \otimes_\nutwopic \of (\otimes_{\tau_1} \times \otimes_{\tau_2}),
\]
again by~\astyle{MC3}.  But $\otimes_\nutwopic = \otimes'_\nutwopic$ since
$J(A, \otimes, \delta) = J(A', \otimes', \delta')$, and $\otimes_{\tau_i} =
\otimes'_{\tau_i}$ by inductive hypothesis, so $\otimes_{(\tau_1, \tau_2)}
= \otimes'_{(\tau_1, \tau_2)}$, completing the induction.
\item $\delta_{\tau, \tau'} = \delta'_{\tau, \tau'}$ for all $\tau, \tau'
\in \ctr(n)$.  This follows from informal coherence for classical
monoidal categories---in particular, from the fact that there is \emph{at
least} one natural transformation $\otimes_\tau \go \otimes_{\tau'}$ built
up from $\alpha$, $\lambda$ and $\rho$.  
\done
\end{itemize}
\end{proof}

\begin{lemma}
The functor $J: \Sigma_\mr{c}\hyph\MClax \go \MClax$ is surjective on
objects. 
\end{lemma}

\begin{proof}
Take a classical monoidal category $(A,\otimes,I,\alpha,\lambda,\rho)$.
Attempt to define a $\Sigma_\mr{c}$-monoidal category $(A,\otimes,\delta)$
as follows:
\begin{itemize}
\item The underlying category $A$ is the same.
\item The tensor $\otimes_\tau: A^n \go A$, for $\tau\in\ctr(n)$, is
defined inductively on $\tau$ by $\otimes_\utree = 1_A$,
$\otimes_\nuzeropic = I$, and
\[
\otimes_{(\tau_1, \tau_2)} = 
(A^{k_1 + k_2} \goby{\otimes_{\tau_1} \times \otimes_{\tau_2}}
A^2 \goby{\otimes} A).
\]
\item The coherence isomorphism $\delta_{\tau,\tau'}: \otimes_\tau \go
\otimes_{\tau'}$ is the canonical natural isomorphism in the statement of
informal coherence. 
\end{itemize}
This does indeed satisfy the axioms for a $\Sigma_\mr{c}$-monoidal
category, as listed in the introductory section:
\begin{description}
\item[\astyle{MC1}, \astyle{MC2}] Immediate.
\item[\astyle{MC3}] That $\otimes_\utree = 1_A$ is immediate.  The other
equation can be proved by induction on $\tau$, or by using the fact that
$\ctr$ is the free operad on $\Sigma_\mr{c} \in\Set^\nat$ and that
$(\Cat(A^n,A))_{n\in\nat}$ forms an operad.
\item[\astyle{MC4}] Follows immediately from informal coherence.  
\end{description}
Finally, $J(A,\otimes,\delta) = (A, \otimes, I, \alpha, \lambda, \rho)$:
\begin{itemize}
\item Clearly the underlying categories agree, both being $A$.
\item We have 
\[
\otimes_\nutwopic = \otimes_{(\utree, \utree)}
= \otimes \of (\otimes_\utree \times \otimes_\utree)
= \otimes
\]
and $\otimes_\nuzeropic = I$, so the tensor products and unit objects
agree.
\item To prove that the coherence isomorphisms agree, we have to check that
\[
\delta_{\assleftpic, \assrightpic} = \alpha,
\diagspace
\delta_{\lambdapic, \utree} = \lambda,
\diagspace
\delta_{\rhopic, \utree} = \rho.
\]
This can be done either by a short calculation or by applying informal
coherence. 
\done
\end{itemize}
\end{proof}

\begin{lemma}
The functor $J: \Sigma_\mr{c}\hyph\MClax \go \MClax$ is faithful.
\end{lemma}

\begin{proof}
Suppose that $A \parpair{(P,\pi)}{(Q,\chi)} A'$ in
$\Sigma_\mr{c}\hyph\MClax$ with $J(P,\pi) = J(Q,\chi)$.  Then:
\begin{itemize}
\item $P=Q$ immediately.
\item $\pi_\tau = \chi_\tau$ for all $\tau\in\ctr(n)$.  This can be proved
either by an induction on $\tau$ similar to the one written out in the
proof of Lemma~\ref{lemma:c-inj}, or by applying informal coherence for lax
monoidal functors.
\done
\end{itemize}
\end{proof}

\begin{lemma}
The functor $J: \Sigma_\mr{c}\hyph\MClax \go \MClax$ is full.
\end{lemma}

\begin{proof}
Let $A, A' \in \Sigma_\mr{c}\hyph\MClax$ and let $J(A) \goby{(P,\pi)}
J(A')$ be a map in $\MClax$.  Attempt to define a map $A \goby{(P,\pi)} A'$
in $\Sigma_\mr{c}\hyph\MClax$ (where again we abuse notation by
re-using the name $(P,\pi)$) as follows:
\begin{itemize}
\item The underlying functor $P$ is the same.
\item The coherence maps $\pi_\tau$ are defined by induction on $\tau$.  We
take $\pi_\utree = 1$; we take
\[
(\pi_\nuzeropic: \otimes_\nuzeropic \go P\otimes_\nuzeropic)
=
(\pi_\cdot: I \go PI);
\]
and we take the component of $\pi_{(\tau_1, \tau_2)}$ at $a_1^1,
\ldots, a_1^{k_1}, a_2^1, \ldots, a_2^{k_2}$
to be the composite
\begin{eqnarray*}
\lefteqn{
(\otimes_{\tau_1}(Pa_1^1, \ldots, Pa_1^{k_1}) 
\ \otimes\ 
\otimes_{\tau_2}(Pa_2^1, \ldots, Pa_2^{k_2}) )}	\\ &
\goby{(\pi_{\tau_1} \otimes \pi_{\tau_2})}	&
(P\otimes_{\tau_1}(a_1^1, \ldots, a_1^{k_1}) \ \otimes\ 
P\otimes_{\tau_2}(a_2^1, \ldots, a_2^{k_2}))	\\
&
\goby{\pi_{\otimes_{\tau_1}(a_1^1, \ldots, a_1^{k_1}),
\otimes_{\tau_2}(a_2^1, \ldots, a_2^{k_2})}} 	&
P(\otimes_{\tau_1}(a_1^1, \ldots, a_1^{k_1}) \ \otimes\ 
\otimes_{\tau_2}(a_2^1, \ldots, a_2^{k_2})).
\end{eqnarray*}
\end{itemize}
This satisfies axioms \astyle{MF1}--\astyle{MF3} for a lax monoidal
functor between $\Sigma_\mr{c}$-monoidal categories, by informal coherence.
We then have $J(P,\pi) = (P, \pi)$:
\begin{itemize}
\item The underlying functors agree, both being $P$.
\item We have
\[
(\pi_\nutwopic)_{a_1,a_2}
= \pi_{\otimes_\utree(a_1), \otimes_\utree(a_2)} \of
((\pi_\utree)_{a_1} \otimes (\pi_\utree)_{a_2})
= \pi_{a_1, a_2} \of (1 \otimes 1)
= \pi_{a_1, a_2}
\]
and $\pi_\nuzeropic = \pi_\cdot$, so the coherence maps also agree.
\done
\end{itemize}
\end{proof}%
\index{tree!classical|)}

Theorem~\ref{thm:diag-coh-mc} now follows in the lax case: $\MClax \iso
\Sigma_\mr{c}\hyph\MClax$.  The weak and strict cases follow from:

\begin{lemma}
The isomorphism $J: \Sigma_\mr{c}\hyph\MClax \goiso \MClax$ restricts to
isomorphisms
\[
\Sigma_\mr{c}\hyph\MCwk \goiso \MCwk,
\diagspace
\Sigma_\mr{c}\hyph\MCstr \goiso \MCstr.
\]
\end{lemma}

\begin{proof}
Trivial.
\done
\end{proof}%
\index{Sigma-monoidal category@$\Sigma$-monoidal category|)}%
\index{monoidal category!Sigma-@$\Sigma$-|)}%
\index{coherence!monoidal categories@for monoidal categories|)}%
\index{monoidal category!classical|)}

\begin{notes}

References can be found in the Notes to Chapter~\ref{ch:monoidal}.

\end{notes}

\chapter{Special Cartesian Monads}
\lbl{app:special-cart}

\chapterquote{%
Pictures can't say `ain't'}{%
Worth~\cite{Wor}}%
\index{monad!cartesian|(}

\noindent
We have met many monads, most of them cartesian.  Some had special
properties beyond being cartesian---for instance, some were the monads
arising from operads, and, as will be explained, some admitted a certain
explicit representation.  Here we look at these special kinds of cartesian
monad and prove some results supporting the theory in the main text.

First~(\ref{sec:opds-alg-thys}) we look at the monads arising from
plain operads.  Monads are algebraic theories, so we can ask which
algebraic theories come from operads.  The answer turns out to be the
strongly regular theories~(\ref{eg:opd-sr}).

In~\ref{sec:alt-app} we saw that the monad arising from an
operad is cartesian.
% ~(\ref{propn:ind-monad-cart}).  
It now follows that the monad corresponding to a strongly regular theory is
cartesian, a fact we used in many of the examples in
Chapter~\ref{ch:gom-basics}.  More precisely, we saw
in~\ref{cor:T-opd-colax} that a monad on $\Set$ arises from a plain operad
if and only if it is cartesian and `augmented over the free monoid monad',
meaning that there exists a cartesian natural transformation from it into
the free monoid functor, commuting with the monad structures.  On the other
hand, not every cartesian monad on $\Set$ possesses such an augmentation,
as we shall see.

One possible drawback of the generalized operad approach to
higher-dimensional category theory is that it can involve monads that are
rather hard to describe explicitly.  For instance, the sequence $T_n$ of
`opetopic' monads (Chapter~\ref{ch:opetopic}) was generated recursively
using nothing more than the existence of free operads, and to describe
$T_n$ explicitly beyond low values of $n$ is difficult.  The second and
third sections of this appendix ease this difficulty.

The basic result is that if a $\Set$-valued functor preserves infinitary or
`wide' pullbacks---which is not asking much more than it be
cartesian---then it is the coproduct of a family of representable functors.
This goes some way towards providing the explicit form desired in the
previous paragraph.  In~\ref{sec:fam-rep-Set} we look at monads on $\Set$
whose functor parts are `familially representable' in this sense.
(Finitary familially representable functors can also be viewed as a
non-symmetric version of the analytic%
\index{functor!analytic}
functors of Joyal~\cite{JoyFAE}.)%
\index{Joyal, Andr\'e}
In~\ref{sec:fam-rep-pshf} we examine monads $(T,\mu,\eta)$ on presheaf
categories whose functor parts $T$ satisfy an analogous condition.  All the
opetopic%
\index{opetopic!monad}
monads $T_n$ are of this form, and this enables us to prove
in~\ref{sec:cart-sym} that every symmetric multicategory gives rise
naturally to a $T_n$-multicategory for each $n\in\nat$.

\begin{figure}
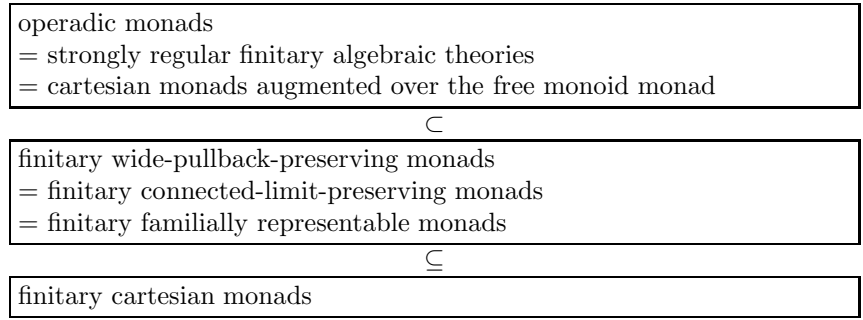

\begin{tabular}{c}
\fbox{\parbox{.9\textwidth}{operadic monads\\
$=$ strongly regular finitary algebraic theories\\
$=$ cartesian monads augmented over the free monoid monad}}\\
$\subset$\\
\fbox{\parbox{.9\textwidth}{finitary wide-pullback-preserving monads\\
$=$ finitary connected-limit-preserving monads\\
$=$ finitary familially representable monads}}\\
$\subseteq$\\
\fbox{\parbox{.9\textwidth}{finitary cartesian monads}}
\end{tabular}
\caption{Classes of finitary cartesian monads on $\Set$}
\label{fig:cart-monads}
\end{figure}

The situation for the various types of finitary cartesian monad on $\Set$ is
depicted in Fig.~\ref{fig:cart-monads}.
The terminology is defined and the proofs are given below.
Example~\ref{eg:fam-rep-not-opdc} shows that the first inclusion is proper;
I do not know whether the second is too.

\section{Operads and algebraic theories}
\lbl{sec:opds-alg-thys}%
\index{algebraic theory!plain operad as|(}
\index{operad!algebraic theory@as algebraic theory|(}

Some algebraic theories can be described by plain operads, and some cannot.
In this section we will prove that the theories that can are precisely
those that admit strongly regular presentations, in the sense of
Example~\ref{eg:opd-sr}.

What exactly this means is as follows.  Any plain operad $P$ gives rise
to a monad $(T_P, \mu^P, \eta^P)$ on $\Set$, as we saw in~\ref{sec:algs}.
Say that a monad on $\Set$ is \demph{operadic}%
\index{monad!operadic}
if it is isomorphic to the
monad $(T_P, \mu^P, \eta^P)$ arising from some plain operad $P$.  On
the other hand, any algebraic theory gives rise to a monad on $\Set$, as we
will discuss in a moment.  Say that a monad on $\Set$ is \demph{strongly
regular}%
\index{monad!strongly regular}%
\index{strongly regular theory}
if it is isomorphic to the monad arising from some strongly
regular finitary theory.  Our result is:
\begin{thm}	\lbl{thm:sr-operadic}
A monad on $\Set$ is operadic if and only if it is strongly regular. 
\end{thm}
This means that if we have before us an algebraic theory presented by
operations and strongly regular equations then we may deduce immediately
that it can be described by an operad.  Moreover, since any operadic monad
is cartesian (Proposition~\ref{propn:ind-monad-cart}), we obtain the useful
corollary:
\begin{cor}
Any strongly regular monad on $\Set$ is cartesian.  
\done
\end{cor}
This is used throughout Section~\ref{sec:cart-monads} to generate examples
of cartesian monads on $\Set$.  It is also proved in Carboni and
Johnstone~\cite{CJ}, Proposition~3.2 (part of which is wrong, as discussed
on p.~\pageref{p:CJ-error}, but the part we are using here is right).  The
other half of Theorem~\ref{thm:sr-operadic} says that any operad can be
presented by a system of operations and strongly regular equations.  This
may not be so useful, but gives the story a tidy ending.

The definition of an algebraic theory and the way that one gives rise to a
monad on $\Set$ are well known, and a detailed account can be found in, for
instance, Manes~\cite{Manes} or Borceux~\cite{Borx2}.  Things are easier
when the theory is strongly regular, and this case is all we will need, so
I will describe it in full.

Let $\Sigma \in \Set^\nat$.  Write $F\Sigma$ for the free plain operad on
$\Sigma$ (or for the underlying object of $\Set^\nat$), as constructed
explicitly in~\ref{sec:om-further}.  We think of $\Sigma$ as a
`signature':%
\index{signature}
$\Sigma(n)$ is the set of primitive $n$-ary operations, and $(F\Sigma)(n)$
is the set of $n$-ary operations derived from the primitive ones.  A
\demph{strongly regular presentation of an algebraic theory}%
\index{strongly regular theory}
is a pair
$(\Sigma, E)$ where $\Sigma \in \Set^\nat$ and $E$ is a family
$(E_n)_{n\in\nat}$ with $E_n \sub (F\Sigma)(n) \times (F\Sigma)(n)$.  We
think of $E$ as the system of equations.  For example, the usual
presentation of the theory of semigroups is
\[
\Sigma(n) = 
\left\{
\begin{array}{ll}
\{\sigma \}	&\textrm{if } n=2	\\
\emptyset	&\textrm{otherwise,}
\end{array}
\right.
\diagspace
E_n = 
\left\{
\begin{array}{ll}
\{ (\sigma\of(\sigma, 1), \sigma\of(1, \sigma)) \}&\textrm{if } n = 3	\\
\emptyset					&\textrm{otherwise.}	\\
\end{array}
\right.
\]
(Implicitly, we are using `algebraic theory'%
\index{algebraic theory}
to mean `finitary,
single-sorted algebraic theory'.)

The operad $F\Sigma$ induces a monad $(T_{F\Sigma}, \mu^{F\Sigma},
\eta^{F\Sigma})$ on $\Set$.  We have 
\[
T_{F\Sigma} A 
= 
\coprod_{n\in\nat} (F\Sigma)(n) \times A^n
\]
for any set $A$; this is usually called the set of `$\Sigma$-terms in $A$'.
Bringing in equations, if $(\Sigma, E)$ is a strongly regular
presentation of an algebraic theory then the induced monad
$(T_{(\Sigma,E)}, \mu^{(\Sigma,E)}, \eta^{(\Sigma,E)})$ on $\Set$ is given
as follows.  For any set $A$,
\[
T_{(\Sigma,E)} A = (T_{F\Sigma} A)/\sim_A%
\glo{simA}
\]
where $\sim_A$ is the smallest equivalence relation on $T_{F\Sigma}A$
such that
\begin{itemize}
\item $\sim_A$ is a congruence:%
\index{congruence}
if $\sigma\in\Sigma(n)$ and we have $\tau_i
\in (F\Sigma)(k_i)$, $\hat{\tau}_i \in (F\Sigma)(\hat{k}_i)$, $a_i^j,
\hat{a}_i^j \in A$ satisfying
\[
(\tau_i, a_i^1, \ldots, a_i^{k_i}) \sim_A 
(\hat{\tau}_i, \hat{a}_i^1, \ldots, \hat{a}_i^{\hat{k}_i})
\]
for each $i = 1, \ldots, n$, then 
\[
(\sigma\of(\tau_1, \ldots, \tau_n), a_1^1, \ldots, a_n^{k_n})
\sim_A
(\sigma\of(\hat{\tau}_1, \ldots, \hat{\tau}_n), 
\hat{a}_1^1, \ldots, \hat{a}_n^{\hat{k}_n})
\]
\item the equations are satisfied: if $(\tau, \hat{\tau}) \in E_n$ and
$\tau_1 \in (F\Sigma)(k_1), \ldots, \tau_n \in (F\Sigma)(k_n)$, $a_i^j \in
A$, then
\[
(\tau\of(\tau_1, \ldots, \tau_n), a_1^1, \ldots, a_n^{k_n})
\sim_A
(\hat{\tau}\of(\tau_1, \ldots, \tau_n), a_1^1, \ldots, a_n^{k_n}).
\]
\end{itemize}
This defines an endofunctor $T_{(\Sigma,E)}$ of $\Set$ and a natural
transformation $\epsln: T_{F\Sigma} \go T_{(\Sigma,E)}$ (the quotient map).
The multiplication $\mu^{(\Sigma,E)}_A$ and unit $\eta^{(\Sigma,E)}_A$ of
the monad are the unique maps making
\begin{equation}	\label{eq:qt-mult-unit}
\begin{diagram}[size=2em]
T_{F\Sigma}^2 A		&\rTo^{\mu^{F\Sigma}_A}	&T_{F\Sigma} A	&
			\lTo^{\eta^{F\Sigma}_A}	&A		\\
\dTo<{(\epsln*\epsln)_A}&			&\dTo~{\epsln_A}&
						&\dEquals	\\
T_{(\Sigma,E)}^2 A	&\rGet_{\mu^{(\Sigma,E)}_A}&T_{(\Sigma,E)} A&
			\lGet_{\eta^{(\Sigma,E)}_A}&A,		\\
\end{diagram}
\end{equation}
commute.  By definition, a monad on $\Set$ is strongly regular if and only
if it is isomorphic to a monad of the form $(T_{(\Sigma,E)},
\mu^{(\Sigma,E)}, \eta^{(\Sigma,E)})$.

We have now given sense to the terms in the theorem and can start to prove
it.  Fix a strongly regular presentation $(\Sigma, E)$ of an algebraic
theory.  The first task is to simplify the clauses above defining $\sim_A$,
which we do by reducing to the case $A=1$.  Note that $\sim_1$ is an
equivalence relation on $\coprod_{n\in\nat} (F\Sigma)(n)$.

\begin{lemma}
Let $A$ be a set, $n, \hat{n} \in \nat$, $\tau \in (F\Sigma)(n)$,
$\hat{\tau} \in (F\Sigma)(\hat{n})$.
\begin{enumerate}
\item 	\lbl{item:one-to-A}
If $\tau \sim_1 \hat{\tau}$ then $n = \hat{n}$ and $(\tau, a_1, \ldots,
a_n) \sim_A (\hat{\tau}, a_1, \ldots, a_n)$ for all $a_1, \ldots, a_n \in
A$. 
\item 	\lbl{item:A-to-one}
If $a_1, \ldots, a_n, \hat{a}_1, \ldots, \hat{a}_{\hat{n}} \in A$ and $(\tau,
a_1, \ldots, a_n) \sim_A (\hat{\tau}, \hat{a}_1, \ldots, \hat{a}_{\hat{n}})$
then $n=\hat{n}$, $a_1=\hat{a}_1, \ldots, a_n=\hat{a}_{\hat{n}}$, and $\tau
\sim_1 \hat{\tau}$.
\end{enumerate}
\end{lemma}

\begin{proof}
To prove~\bref{item:one-to-A}, define a relation $\approx$ on
$\coprod_{m\in\nat} (F\Sigma)(m)$ as follows: for $\phi\in(F\Sigma)(m)$ and
$\hat{\phi} \in (F\Sigma)(\hat{m})$,
\[
\phi \approx \hat{\phi}
\iff
m = \hat{m} 
\textrm{ and }
(\phi, a_1, \ldots, a_m) \sim_A (\hat{\phi}, a_1, \ldots, a_m)
\textrm{ for all }
a_1, \ldots, a_m \in A.
\]
Then $\approx$ is an equivalence relation since $\sim_A$ is.  If we can
prove that $\approx$ also satisfies the two conditions for which $\sim_1$
is minimal then we will have $\sim_1 \sub \approx$,
proving~\bref{item:one-to-A}.  So, for the first condition ($\approx$ is a
congruence): suppose that $\sigma\in\Sigma(m)$ and that for each $i=1,
\ldots, m$ we have $\phi_i, \hat{\phi}_i \in (F\Sigma)(k_i)$ with $\phi_i
\approx \hat{\phi}_i$.  Then for all $a_1^1, \ldots, a_m^{k_m} \in A$, the
fact that $\phi_i \approx \hat{\phi}_i$ implies that
\[
(\phi_i, a_i^1, \ldots, a_i^{k_i}) 
\sim_A
(\hat{\phi}_i, a_i^1, \ldots, a_i^{k_i}),
\]
and so the fact that $\sim_A$ is a congruence implies that
\[
(\sigma\of(\phi_1, \ldots, \phi_m), a_1^1, \ldots, a_m^{k_m})
\sim_A
(\sigma\of(\hat{\phi}_1, \ldots, \hat{\phi}_m), a_1^1, \ldots, a_m^{k_m}).
\]
Hence $\sigma\of(\phi_1, \ldots, \phi_m) \approx \sigma\of(\hat{\phi}_1,
\ldots, \hat{\phi}_m)$, as required.  For the second condition (the
equations are satisfied): if $(\phi, \hat{\phi}) \in E_m$ then
$\phi\approx\hat{\phi}$, by taking $\tau_1, \ldots, \tau_n$ all to be the
identity $1 \in (F\Sigma)(1)$ in the second condition for $\sim_A$.  

Several of the proofs here are of this type: we have an equivalence
relation $\sim$ (in this case $\sim_1$) defined to be minimal such that
certain conditions are satisfied, and we want to prove that if $x \sim y$
then some conclusion involving $x$ and $y$ holds.  To do this we define a
new equivalence relation $\approx$ by `$x \approx y$ if and only if the
conclusion holds'; we then show that $\approx$ satisfies the conditions for
which $\sim$ was minimal, and the result follows.  These proofs, like
diagram chases in homological algebra, are more illuminating to write than
to read, so the remaining similar ones (including~\bref{item:A-to-one}) are
omitted.  None of them is significantly more complicated than the one
just done.  
\done
\end{proof}

\begin{cor}	\lbl{cor:functors-same}
For any set $A$, the quotient map $\epsln_A: T_{F\Sigma}A \go T_{(\Sigma,
E)}A$ induces a bijection 
\[
\coprod_{n\in\nat} ((F\Sigma)(n)/\sim_1) \times A^n
\goiso
T_{(\Sigma,E)}A.
\]
\ \done
\end{cor}

Having expressed the equivalence relation $\sim_A$ for an arbitrary set $A$
in terms of the single equivalence relation $\sim_1$, we now show that
$\sim_1$ has a congruence property of an operadic kind.
\begin{lemma}	\lbl{lemma:operadic-congruence}
If $\tau, \hat{\tau} \in (F\Sigma)(n)$ with $\tau \sim_1 \hat{\tau}$, and
$\tau_i, \hat{\tau}_i \in (F\Sigma)(k_i)$ with $\tau_i \sim_1 \hat{\tau}_i$
for each $i = 1, \ldots, n$, then
\[
\tau \of (\tau_1, \ldots, \tau_n) \sim_1 
\hat{\tau} \of (\hat{\tau}_1, \ldots, \hat{\tau}_n).
\]
\end{lemma}

\begin{proof}
Define a relation $\approx$ on $\coprod_{m\in\nat}(F\Sigma)(m)$ as follows:
if $\phi\in (F\Sigma)(m)$ and $\hat{\phi}\in (F\Sigma)(\hat{m})$ then $\phi
\approx \hat{\phi}$ if and only if $m=\hat{m}$ and
\[
\begin{array}{ll}
\textrm{for all }  	&
\phi_i, \hat{\phi}_i \in (F\Sigma)(k_i) 
\textrm{ with }
\phi_i \sim_1 \hat{\phi}_i
\ (i = 1, \ldots, n),	\\
\textrm{we have }	&
\phi \of (\phi_1, \ldots, \phi_m)
\sim_1
\hat{\phi} \of (\hat{\phi}_1, \ldots, \hat{\phi}_m).
\end{array}
\]
Then show that $\sim_1 \sub \approx$ in the usual way. 
\done
\end{proof}
Define $P_{(\Sigma,E)}(n) = (F\Sigma)(n)/\sim_1$ for each $n$: then
Lemma~\ref{lemma:operadic-congruence} tells us that there is a unique
operad structure on $P_{(\Sigma,E)}$ such that the quotient map $F\Sigma
\go P_{(\Sigma,E)}$ is a map of operads.  So starting from a strongly
regular presentation of a theory we have constructed an operad, and this
operad induces the same monad as the theory:
\begin{cor}	\lbl{cor:sr-implies-operadic}
There is an isomorphism of monads $T_{P_{(\Sigma,E)}} \iso
T_{(\Sigma,E)}$.  Hence any strongly regular monad is operadic.
\end{cor}
\begin{proof}
We have only to prove the first sentence; the second follows immediately.
By Corollary~\ref{cor:functors-same}, there is an isomorphism of functors
$T_{P_{(\Sigma,E)}} \iso T_{(\Sigma,E)}$ making the diagram
\[
\begin{slopeydiag}
		&	&T_{F\Sigma}	&		&	\\
		&\ldTo	&		&\rdTo>{\epsln}	&	\\
T_{P_{(\Sigma,E)}}&	&\rTo^{\diso}	&		&T_{(\Sigma,E)}\\
\end{slopeydiag}
\]
commute, where the left-hand arrow is the quotient map.  The multiplication
$\mu^{P_{(\Sigma,E)}}$ of the monad $T_{P_{(\Sigma,E)}}$ comes from
composition in the operad $P_{(\Sigma,E)}$, and this in turn comes via the
quotient map from composition in the operad $F\Sigma$, so it follows from
diagram~\bref{eq:qt-mult-unit} (p.~\pageref{eq:qt-mult-unit}) that
$\mu^{P_{(\Sigma,E)}}$ corresponds under the isomorphism to
$\mu^{(\Sigma,E)}$.  The same goes for units; hence the result.  \done
\end{proof}

This concludes the proof of the more `useful' half of
Theorem~\ref{thm:sr-operadic}.  To prove the converse, we first establish
one more fact about strongly regular presentations in general.
\begin{lemma}	\lbl{lemma:univ-sr}
Let $(\Sigma,E)$ be a strongly regular presentation of an algebraic
theory.  Then the quotient map $\epsln: F\Sigma \go P_{(\Sigma,E)}$ is the
universal operad map out of $F\Sigma$ with the property that $\epsln(\tau)
= \epsln(\hat{\tau})$ for all $(\tau,\hat{\tau}) \in E_n$.
\end{lemma}
`Universal' means that if $\zeta: F\Sigma \go Q$ is a map of operads
satisfying $\zeta(\tau) = \zeta(\hat{\tau})$ for all $(\tau,\hat{\tau}) \in
E_n$ then there is a unique operad map $\ovln{\zeta}: P_{(\Sigma,E)} \go Q$
such that $\ovln{\zeta} \of \epsln = \zeta$.
\begin{proof}
Universality can be re-formulated as the following condition: $\sim_1$ is
the smallest equivalence relation on $\coprod_{n\in\nat}$ that is an
operadic congruence (in other words, such that
Lemma~\ref{lemma:operadic-congruence} holds) and satisfies $\tau \sim_1
\hat{\tau}$ for all $(\tau,\hat{\tau}) \in E_n$.  The proof is by the usual
method.  \done
\end{proof}

Now fix an operad $P$.  Define $\Sigma_P \in \Set^\nat$ by
\[
\Sigma(n) = \{\sigma_\theta \such \theta\in P(n) \}
\]
where $\sigma_\theta$ is a formal symbol.  Let $E_P$ be the system of
equations with elements
\[
(\sigma_\theta \of (\sigma_{\theta_1}, \ldots, \sigma_{\theta_n}),
\sigma_{\theta\sof (\theta_1, \ldots, \theta_n)})
\in 
(E_P)_{k_1 + \cdots + k_n}
\]
for each $\theta\in P(n)$, $\theta_i \in P(k_i)$, and 
\[
(1, \sigma_1) \in (E_P)_1.
\]
Then $(\Sigma_P, E_P)$ is a strongly regular presentation of an algebraic
theory. 

There is a map $\epsln: F\Sigma_P \go P$ of operads, a component of the
counit of the adjunction between $\Operad$ and $\Set^\nat$.
Concretely~(\ref{sec:om-further}), $\epsln$ is described by the inductive
clauses
\begin{eqnarray*}
\epsln(1)	&=	&1,	\\
\epsln(\sigma_\theta \of (\tau_1, \ldots, \tau_n))	&
=	&
\theta \of (\epsln(\tau_1), \ldots, \epsln(\tau_n))	
\end{eqnarray*}
for $\theta\in P(n), \tau_i\in (F\Sigma_P)(k_i)$.
\begin{lemma}	\lbl{lemma:univ-operadic}
$\epsln: F\Sigma_P \go P$ is the universal operad map out of $F\Sigma_P$
with the property that $\epsln(\tau) = \epsln(\hat{\tau})$ for all
$(\tau,\hat{\tau}) \in (E_P)_n$.
\end{lemma}
\begin{proof}
That $\epsln$ does have the property is easily verified.  For universality,
take a map $\zeta: F\Sigma \go Q$ of operads such that $\zeta(\tau) =
\zeta(\hat{\tau})$ for all $(\tau,\hat{\tau}) \in (E_P)_n$; we want there
to be a unique operad map $\ovln{\zeta}: P \go Q$ such that $\ovln{\zeta}
\of \epsln = \zeta$.  Since $\epsln(\sigma_\theta) = \theta$ for all
$\theta \in P(n)$, the only possibility is $\ovln{\zeta}(\theta) =
\zeta(\sigma_\theta)$, and it is easy to check that this map $\ovln{\zeta}$
does satisfy the conditions required.  \done
\end{proof}

\begin{cor}	\lbl{cor:operadic-implies-sr}
There is an isomorphism of monads $T_{(\Sigma_P,E_P)} \iso T_P$.  Hence any
operadic monad is strongly regular.
\end{cor}
\begin{proof}
Lemmas~\ref{lemma:univ-sr} and~\ref{lemma:univ-operadic} together imply
that there is an isomorphism of operads $P_{(\Sigma_P,E_P)} \iso P$; so
there is an isomorphism of monads $T_{P_{(\Sigma_P,E_P)}} \iso T_P$.
Corollary~\ref{cor:sr-implies-operadic} implies that
$T_{P_{(\Sigma_P,E_P)}} \iso T_{(\Sigma_P,E_P)}$.  The result follows.
\done
\end{proof}

The proof of Theorem~\ref{thm:sr-operadic} is now complete.

Similar results can be envisaged for other kinds of operad.  For example,
any symmetric%
\index{operad!symmetric}
operad induces a monad on $\Set$ whose algebras are exactly
the algebras for the operad; an obvious conjecture is that the monads
arising in this way are those that can be presented by finitary operations
and equations in which the same variables appear, without repetition, on
each side (but not necessarily in the same order).  Another possibility is
to produce results of this kind for operads in a symmetric%
\index{operad!symmetric monoidal category@in symmetric monoidal category}
monoidal
category; compare Example~\ref{eg:opd-sr-enr}.  We would like, for
instance, a general principle telling us that the theory of Lie%
\index{Lie algebra}
algebras
can be described by a symmetric operad of vector spaces on the grounds that
its governing equations
\begin{eqnarray*}
{[}x,y] + [y,x] 				&=	&0,	\\
{[}x,[y,z]] + [y,[z,x]] + [z,[x,y]] 		&=	&0
\end{eqnarray*}
are both `good': the same variables are involved, without repetition, in
each summand.  This principle is well-known informally, but as far as I
know has not been proved.

A more difficult generalization would be to algebraic%
\index{algebraic theory!presheaf category@on presheaf category}
theories on presheaf
categories.  For example, take the theory of strict $\omega$-categories%
\index{omega-category@$\omega$-category!strict!theory of}
presented in the algebraic way~(\ref{defn:strict-n-cat-glob}).  This is not
a `strongly regular' presentation in the simple-minded sense, because the
interchange%
\index{interchange}
law
\[
(\beta' \ofdim{p} \beta) \ofdim{q} (\alpha' \ofdim{p} \alpha)
=
(\beta' \ofdim{q} \alpha') \ofdim{p} (\beta \ofdim{q} \alpha)
\]
does permute the order of the variables.  Yet somehow this equation should
be regarded as `good', since when the appropriate picture is drawn---e.g.
\lbl{p:sr-presheaf}%
\[
\left(
\begin{array}{c}
\gfstsu \gtwos{}{}{\alpha} \glstsu \\
\of \\
\gfstsu \gtwos{}{}{\alpha'} \glstsu
\end{array}
\right)
\ \of\ 
\left(
\begin{array}{c}
\gfstsu \gtwos{}{}{\beta} \glstsu \\
\of \\
\gfstsu \gtwos{}{}{\beta'} \glstsu
\end{array}
\right)
=
\begin{array}{c}
\left(
\gfstsu \gtwos{}{}{\alpha} \glstsu
\ \of\ 
\gfstsu \gtwos{}{}{\beta} \glstsu
\right)					\\
\of					\\
\left(
\gfstsu \gtwos{}{}{\alpha'} \glstsu
\ \of\ 
\gfstsu \gtwos{}{}{\beta'} \glstsu
\right)					
\end{array}
\]
for $p=1$ and $q=0$---there is no movement in the positions of the cells.
So we hope for some notion of a strongly regular theory on a presheaf
category, and a result saying that the monad induced by such a theory is at
least cartesian, and perhaps operadic in some sense.  This would render
unnecessary the \latin{ad hoc} calculations of
Appendix~\ref{app:free-strict} showing that the free strict
$\omega$-category monad is cartesian.  However, such a general result has
yet to be found.%
\index{algebraic theory!plain operad as|)}
\index{operad!algebraic theory@as algebraic theory|)}

\section{Familially representable monads on $\Set$}
\lbl{sec:fam-rep-Set}

\index{representable functor|(}
Not every set-valued functor is representable, but every set-valued functor
is a colimit of representables.  It turns out that an intermediate
condition is relevant to the theory of operads: that of being a
\emph{coproduct} of representables.  Such a functor is said to be
`familially representable'.  

Familial representability has been studied in, among other places, an
important paper of Carboni%
\index{Carboni, Aurelio}
and Johnstone~\cite{CJ}.%
\index{Johnstone, Peter}
 This section is not
much more than an account of some of their results.  Our goal is to
describe monads on $\Set$ whose functor parts are familially representable
and whose natural transformation parts are cartesian; we lead up to this by
considering, more generally, familially representable functors into $\Set$.

But first, as a warm-up, consider ordinary representability.  Here follow
some very well-known facts, presented in a way that foreshadows the
material on familial representability.

Let us say that a category $\cat{A}$ has the \demph{adjoint functor
property} if any limit-preserving functor from $\cat{A}$ to a locally small
category has a left adjoint.  (`Limit-preserving' really means
`small-limit-preserving'.)  For example, the Special Adjoint Functor
Theorem states that any complete, well-powered, locally small category with
a cogenerating set has the adjoint functor property.  In particular, $\Set$
has the adjoint functor property.

\begin{propn}	\lbl{propn:eqv-condns-rep}
Let $\cat{A}$ be a category with the adjoint functor property.  The
following conditions on a functor $T: \cat{A} \go \Set$ are equivalent:
\begin{enumerate}
\item	\lbl{item:rep-lim}
$T$ preserves limits
\item 	\lbl{item:rep-adj}
$T$ has a left adjoint
\item 	\lbl{item:rep-rep}
$T$ is representable.
\end{enumerate}
\end{propn}
\begin{proof}
\begin{description}
\item[\bref{item:rep-lim}\implies\bref{item:rep-adj}] Adjoint functor
property.
\item[\bref{item:rep-adj}\implies\bref{item:rep-rep}] Take a left
adjoint $S$ to $T$: then
$T$ is represented by $S1$.
\item[\bref{item:rep-adj}\implies\bref{item:rep-rep}] Standard.
\done
\end{description}
\end{proof}

Suppose we want to define the subcategory of $\ftrcat{\cat{A}}{\Set}$
consisting of only those things that are representable.  We know what it
means for a functor $\cat{A} \go \Set$ to be representable.  By rights, a
natural transformation
\begin{equation}	\label{eq:rep-transf}
\cat{A}(X, \dashbk) \go \cat{A}(X', \dashbk)
\end{equation}
should qualify as a map in our subcategory just when it is of the form
$\cat{A}(f, \dashbk)$ for some map $f:X' \go X$ in $\cat{A}$.  But the
Yoneda Lemma implies that \emph{all} natural
transformations~\bref{eq:rep-transf} are of this form, so in fact the
appropriate subcategory is just the full subcategory of representable
functors.

Representable functors are certainly cartesian, so we might reasonably ask
when a natural transformation between representables is cartesian.%
\index{transformation!cartesian}
 The
answer is `seldom': the natural transformation~\bref{eq:rep-transf}
induced by a map $f: X' \go X$ is cartesian just when $f$ is an
isomorphism.%
\index{representable functor|)}

The basic result on familial representability is like
Proposition~\ref{propn:eqv-condns-rep}, but involves a weaker set of
equivalent conditions.  To express them we need some terminology.

Let $\cat{A}$ be a category, $I$ a set, and $(X_i)_{i\in I}$
a family of objects of $\cat{A}$: then there is a functor
\[
\coprod_{i\in I} \cat{A}(X_i, \dashbk): \cat{A} \go \Set.
\]
Such a functor, or one isomorphic to it, is said to be \demph{familially
representable};%
\index{familial representability!functor into Set@of functor into $\Set$}
$(X_i)_{i\in I}$ is the \demph{representing family}.  Note
that we can recover $I$ as the value of the functor at the terminal object
of $\cat{A}$, if it has one.

A category is \demph{connected}%
\index{connected}
if it is nonempty and any functor from it
into a discrete category is constant.  A \demph{connected limit}%
\index{limit!connected}
is a limit
of a functor whose domain is a connected category.  The crucial fact about
connected limits is:
\begin{lemma}	\lbl{lemma:conn-lims}
Connected limits commute with small coproducts in $\Set$.
\done
\end{lemma}

\begin{cor}	\lbl{cor:conn-lims}
For any set $I$, the forgetful functor $\Set/I \go \Set$ preserves and
reflects connected limits.
\end{cor}
\begin{proof}
$\Set/I$ is equivalent to $\Set^I$, and under this equivalence the
forgetful functor becomes the functor $\Sigma: \Set^I \go \Set$ defined by
$\Sigma((X_i)_{i\in I}) = \coprod_{i\in I} X_i$.
Lemma~\ref{lemma:conn-lims} tells us that $\Sigma$ preserves connected
limits.  Moreover, $\Set^I$ has all connected limits and $\Sigma$ reflects
isomorphisms, so preservation implies reflection.  \done
\end{proof}

One type of connected limit is of particular importance to us.  A
\demph{wide%
\lbl{p:defn-wide-pb}
pullback}%
\index{wide pullback}%
\index{pullback}
is a limit of shape
\begin{diagram}[height=2em]
\gzersu	&\gzersu&\gzersu	&\ldots	&	&\gzersu	&\ldots\\
	&\rdTo(4,2)&\rdTo(3,2)	&\rdTo(2,2)&\ldTo(1,2)&		&	\\
	&	&		&	&\gzeros{}&		&	\\
\end{diagram}
where the top row indicates a set of any cardinality.  Formally, for any
set $K$ let $\scat{P}_K$%
\lbl{p:defn-wide-pb-shape}\glo{widepbshape}
be the (connected) category whose objects are the elements of
$K$ together with a further object $\infty$ and whose only non-identity
maps are a single map $k \go \infty$ for each $k\in K$; a wide pullback is
a limit of a functor with domain $\scat{P}_K$ for some set $K$.  If $K$ has
cardinality two then this is just an ordinary pullback.

If $T: \cat{A} \go \cat{B}$ is a functor and $\cat{A}$ has a terminal
object then we write $\slind{T}: \cat{A} \go \cat{B}/T1$ for the functor
\[
X \goesto \bktdvslob{TX}{T!}{T1}.
\]
We recover $T$ from $\slind{T}$ as the composite
\begin{equation}	\label{eq:slind-comp}
T = (\cat{A} \goby{\slind{T}} \cat{B}/T1 \goby{\mr{forgetful}} \cat{B}).
\end{equation}

\begin{propn}	\lbl{propn:eqv-condns-fam-rep}
Let \cat{A} be a category with the adjoint functor property and a terminal
object.  The following conditions on a functor $T: \cat{A} \go \Set$ are
equivalent:
\begin{enumerate}
\item	\lbl{item:fam-rep-conn}
$T$ preserves connected limits
\item	\lbl{item:fam-rep-wide}
$T$ preserves wide pullbacks
\item	\lbl{item:fam-rep-lim}
$\slind{T}$ preserves limits
\item	\lbl{item:fam-rep-adj}
$\slind{T}$ has a left adjoint
\item	\lbl{item:fam-rep-fam-rep}
$T$ is familially representable.
\end{enumerate}
\end{propn}
\begin{proof}
\begin{description}
\item[\bref{item:fam-rep-conn}\implies\bref{item:fam-rep-wide}]
A wide pullback is a connected limit.
\item[\bref{item:fam-rep-wide}\implies\bref{item:fam-rep-lim}] The
forgetful functor $\Set/T1 \go \Set$ reflects wide pullbacks
(Corollary~\ref{cor:conn-lims}), and $T$ preserves them, so by
equation~\bref{eq:slind-comp}, $\slind{T}$ preserves them too.  $\slind{T}$
also preserves terminal objects, trivially.  So $\slind{T}$ preserves all
limits.
\item[\bref{item:fam-rep-lim}\implies\bref{item:fam-rep-adj}]
Adjoint functor property.
\item[\bref{item:fam-rep-adj}\implies\bref{item:fam-rep-fam-rep}]
Take a left adjoint $S$ to $\slind{T}$.  For each $i\in T1$
we have adjunctions
\[
\begin{diagram}[scriptlabels]
\cat{A}	&
\pile{\rTo^{\slind{T}}_\top\\  \\ \lTo_{S}}	&
\Set/T1	&
\pile{\rTo^{\blank_i}_\top \\ \\ \lTo_{i_!}}	&
\Set,	\\
\end{diagram}
\]
where $\blank_i$ takes the fibre over $i$ and $i_! K$ is the function $K
\go T1$ constant at $i$.  So
\[
TX	\iso	\coprod_{i\in T1} (\slind{T}X)_i	
	\iso	\coprod_{i\in T1} \Set(1, (\slind{T}X)_i)	
	\iso	\coprod_{i\in T1} \cat{A}(S i_! 1, X)
\]
naturally in $X \in \cat{A}$.

\item[\bref{item:fam-rep-fam-rep}\implies\bref{item:fam-rep-conn}]
Representables preserve limits~(\ref{propn:eqv-condns-rep}) and coproducts
commute with connected limits in $\Set$~(\ref{lemma:conn-lims}), so any
coproduct of representables preserves connected limits.
\done
\end{description}
\end{proof}

Let $(X_i)_{i\in I}$ and $(X'_{i'})_{i'\in I'}$ be families of objects in a
category $\cat{A}$.  A function $\phi: I \go I'$ together with a family of
maps $(X'_{\phi(i)} \goby{f_i} X_i)_{i\in I}$ induces%
\lbl{p:transf-fam-rep}
a natural transformation%
\index{transformation!familially representable functors@of familially representable functors}
\begin{equation}	\label{eq:fam-rep-transf}
\coprod_{i \in I} \cat{A}(X_i, \dashbk)
\go
\coprod_{i' \in I'} \cat{A}(X'_{i'}, \dashbk),
\end{equation}
and the Yoneda Lemma implies that, in fact, all natural
transformations~\bref{eq:fam-rep-transf} arise uniquely in this way.  Let
$\Fam(\cat{A})$%
\glo{Fam}
be the category whose objects are families $(X_i)_{i\in I}$
of objects of $\cat{A}$ and whose maps $(X_i)_{i\in I} \go (X'_{i'})_{i'\in
I'}$ are pairs $(\phi, f)$ as above.  (This is the free coproduct
cocompletion of $\cat{A}^\op$.)  There is a functor
\begin{equation}	\label{eq:fam-rep-Yoneda}
\begin{array}{rrcl}
\fcat{y}:	&\Fam(\cat{A})	&\go	&\ftrcat{\cat{A}}{\Set}	\\
		&(X_i)_{i\in I}	&\goesto&
\coprod_{i\in I} \cat{A} (X_i, \dashbk),%
\glo{FamYon}
\end{array}
\end{equation}
and we have just seen that it is full and faithful.  It therefore defines
an equivalence between $\Fam(\cat{A})$ and the full subcategory of
$\ftrcat{\cat{A}}{\Set}$ formed by the familially representable functors.

A transformation $(\phi, f)$ as above is cartesian%
\lbl{p:cart-transf-fam-rep}%
\index{transformation!cartesian}
just when each of its pieces
\[
\cat{A}(f_i, \dashbk):
\cat{A}(X_i, \dashbk)
\go
\cat{A}(X'_{\phi(i)}, \dashbk)
\]
($i\in I$) is cartesian, which, from our earlier result, happens just when
each $f_i: X'_{\phi(i)} \go X_i$ is an isomorphism.

We have already encountered many familially representable endofunctors of
$\Set$:
\begin{example}
For any plain operad $P$, the functor part of the induced monad $T_P$
on $\Set$ is familially representable.  Indeed,
\[
T_P(X) 
\iso \coprod_{n\in\nat} P(n) \times X^n
\iso \coprod_{i\in I} [X_i, X]
\]
where $I=\coprod_{n\in\nat} P(n)$, $X_i$ is an $n$-element set for $i\in
P(n)$, and $[X_i,X]$%
\glo{homsetset}
denotes the hom-set $\Set(X_i,X)$.
\end{example}
In this example all of the sets $X_i$ are finite and, as we saw
in~\ref{sec:opds-alg-thys}, the monad $T_P$ corresponds to a finitary
algebraic theory.  Here is the general result on finiteness.%
\index{finitary}%
\index{algebraic theory!finitary}
\begin{propn}	\lbl{propn:fam-rep-Set-finitary}
The following are equivalent conditions on a family $(X_i)_{i\in I}$ of
sets:
\begin{enumerate}
\item	\lbl{item:fam-rep-finitary}
the functor $\coprod_{i\in I} [X_i, \dashbk]: \Set \go \Set$ is finitary
\item	\lbl{item:fam-rep-finite}
the set $X_i$ is finite for each $i\in I$
\item	\lbl{item:fam-rep-fin-operadic}
there is a sequence $(P(n))_{n\in\nat}$ of sets such that 
\[
\coprod_{i\in I} [X_i, \dashbk] 
\iso 
\coprod_{n\in\nat} P(n) \times (\dashbk)^n.
\]
\end{enumerate}
\end{propn}
\begin{proof}
\begin{description}
\item[\bref{item:fam-rep-finitary}\implies\bref{item:fam-rep-finite}] 
For each $i\in I$ we have a pullback square
\[
\begin{diagram}[height=2.5em,width=4em]
[X_i, \dashbk]	&\rIncl^{\mr{copr}_i}	&\coprod_{i\in I} [X_i, \dashbk]\\
\dTo		&			&\dTo>{\mr{pr}}			\\
\Delta 1	&\rIncl_{\mr{copr}_i}	&\Delta I			\\
\end{diagram}
\]
in $\ftrcat{\Set}{\Set}$, where $\Delta K: \Set \go \Set$ denotes the
constant functor with value $K$.  We know that all four functors except
perhaps $[X_i, \dashbk]$ are finitary; since pullbacks commute with
filtered colimits in $\Set$ it follows that $[X_i, \dashbk]$ is finitary
too.  This in turn implies that $X_i$ is finite, as is well-known (Ad\'amek
and Rosick\'y, \cite[p.~9]{AR}).
\item[\bref{item:fam-rep-finite}\implies\bref{item:fam-rep-fin-operadic}]
For each $n\in\nat$, put
\[
P(n) = \{i \in I \such \textrm{cardinality}(X_i) = n \}.
\]
By choosing for each $i\in I$ a bijection $X_i \goiso \{1, \ldots, n\}$, we
obtain isomorphisms
\[
\coprod_{i\in I} [X_i, X] 
\iso
\coprod_{n\in\nat} \coprod_{i \in P(n)} [X_i, X]
\iso 
\coprod_{n\in\nat} \coprod_{i \in P(n)} X^n
\iso
\coprod_{n\in\nat} P(n) \times X^n
\]
natural in $X \in \Set$.
\item[\bref{item:fam-rep-fin-operadic}\implies\bref{item:fam-rep-finitary}]
If~\bref{item:fam-rep-fin-operadic} holds then $\coprod_{i\in I} [X_i,
\dashbk]$ is a colimit of endofunctors of $\Set$ of the form $(\dashbk)^n$
($n\in\nat$); each of these is finitary, so the whole functor is finitary.
\done
\end{description}
\end{proof}
So an endofunctor of $\Set$ is finitary and familially representable if and
only if it is, in the sense of~\bref{item:fam-rep-fin-operadic}, operadic.
The story for \emph{monads} is quite different, as we see
below~(\ref{eg:fam-rep-not-opdc}).  A hint of the subtlety is that the
isomorphism constructed in the proof of
\bref{item:fam-rep-finite}\implies\bref{item:fam-rep-fin-operadic} is not
canonical: it involves an arbitrary choice of a total ordering of $X_i$,
for each $i\in I$.

Let us turn, then, to monads on $\Set$.  A \demph{familially representable
monad}%
\index{familial representability!monad on Set@of monad on $\Set$}%
\index{monad!familially representable}
$(T, \mu, \eta)$ on $\Set$ is one whose functor part $T$ is
familially representable and whose natural transformation parts $\mu$ and
$\eta$ are cartesian.  (For example, any operadic monad on $\Set$ is
familially representable.)  To see what $\mu$ and $\eta$ look like in terms
of the representing family of $T$, we need to consider composition of
familially representable functors.

So, let 
\[
S \iso \coprod_{h\in H} [W_h, \dashbk],
\diagspace
T \iso \coprod_{i\in I} [X_i, \dashbk]
\]
be familially representable endofunctors of $\Set$.
Conditions~\bref{item:fam-rep-conn} and~\bref{item:fam-rep-wide} of
Proposition~\ref{propn:eqv-condns-fam-rep} make it clear that $T\of S$ must
be familially representable, but what is the representing family?  For any
set $X$ we have
\begin{eqnarray*}
TSX	&\iso	&\coprod_{i\in I} 
		\left[X_i, \coprod_{h\in H} [W_h, X] \right]	\\
	&\iso	&\coprod_{i\in I} \coprod_{g\in [X_i, H]}
		\left[\coprod_{x\in X_i} W_{g(x)}, X \right].	
\end{eqnarray*}
Hence
\[
T\of S \iso \coprod_{j\in J} [Y_j, \dashbk]
\]
where
\[
J = \coprod_{i\in I} [X_i, H] = TH,
\diagspace
Y_{i,g} = \coprod_{x\in X_i} W_{g(x)} 
\]
for $i\in I$ and $g\in [X_i, H]$ (that is, $(i,g) \in J$).  In particular, the
case $S=T$ gives
\[
T^2 \iso \coprod_{j\in J} [Y_j, \dashbk]
\]
where
\[
J = \coprod_{i\in I} [X_i, I] = TI,
\diagspace
Y_{i,g} = \coprod_{x\in X_i} X_{g(x)}. 
\]
The familial representation of the identity functor is
\[
\id \iso \coprod_{i\in 1} [1, \dashbk].
\]

A familially representable monad on $\Set$ therefore consists of
\begin{itemize}
\item a family $(X_i)_{i\in I}$ of sets, inducing the functor $T =
\coprod_{i\in I} [X_i, \dashbk]$
\item a function
\[
m: \{ (i, g) \such i \in I, g \in [X_i, I] \}  \go  I
\]
and for each $i\in I$ and $g\in [X_i, I]$, a bijection
\[
X_{m(i,g)} \goiso \coprod_{x\in X_i} X_{g(x)},
\]
inducing the natural transformation $\mu: T^2 \go T$ whose component
$\mu_X$ at a set $X$ is the composite
\begin{eqnarray*}
\coprod_{i\in I, g\in [X_i, I]} 
\left[ \coprod_{x\in X_i} X_{g(x)}, X \right] 
&
\goiso
&
\coprod_{i\in I, g\in [X_i, I]} [X_{m(i,g)}, X]
\\
&
\goby{\mr{canonical}}
&
\coprod_{k\in I} [X_k, X]
\end{eqnarray*}
\item an element $e\in I$
such that $X_e$ is a one-element set, inducing the
natural transformation $\eta: 1 \go T$ whose component $\eta_X$ at a set
$X$ is the composite
\[
X \goiso [X_e,X] \rIncl \coprod_{i\in I} [X_i, X],
\]
\end{itemize}
such that $\mu$ and $\eta$ obey associativity and unit laws.  The monad
is finitary just when all the $X_i$'s are finite.

\begin{example}
Let $(T, \mu, \eta)$ be the free monoid monad on $\Set$.  Then
\[
TX = \coprod_{i\in\nat} [i,X],
\]
so in this case $I=\nat$ and $X_i = i$.  (We use $i$ to denote an
$i$-element set, say $\{1, \ldots, i\}$.)  The multiplication function
\[
m: \{ (i, g) \such i \in \nat, g\in [i, \nat] \}  \go  \nat
\]
is given by
\[
m(i, g) = g(1) + \cdots + g(i),
\]
and there is an obvious choice of bijection
\[
X_{g(1) + \cdots + g(i)} \goiso X_{g(1)} + \cdots + X_{g(i)}.
\]
The unit element $e$ is $1\in\nat$.
\end{example}

This explicit form for familially representable monads on $\Set$ will
enable us to prove, for instance, that a commutative monoid is naturally an
algebra for any finitary familially representable monad on $\Set$
(Example~\ref{eg:cm-ftr-Set}).

We have now proved almost all of the equalities and inclusions in
Fig.~\ref{fig:cart-monads}.  All that remains is to prove that the first
inclusion is proper: not every finitary familially representable monad on
$\Set$ is operadic (strongly%
\index{strongly regular theory!familially representable monad@\vs.\ familially representable monad}%
\index{familial representability!strong regularity@\vs.\ strong regularity}
regular).  At this point we come to an error in Carboni and
Johnstone~\cite{CJ}.  In their Proposition~3.2, they prove that every
strongly regular monad on $\Set$ is familially representable; but they also
claim the converse, which is false.%
\lbl{p:CJ-error}  
The error is in the final sentence of the proof.  Having started with a
familially representable monad $(T, \mu, \eta)$, they construct a strongly
regular monad---let us call it $(T',\mu',\eta')$---and claim that the
monads $(T, \mu, \eta)$ and $(T', \mu', \eta')$ are isomorphic; and while
there is indeed a (non-canonical) isomorphism between the functors $T$ and
$T'$, there is not in general an isomorphism that commutes with the monad
structures.  The following counterexample is from Carboni and Johnstone's
corrigenda~\cite{CJcorr}.

\begin{example}		\lbl{eg:fam-rep-not-opdc}
The monad $(T, \mu, \eta)$ on $\Set$ corresponding to the theory of
monoids%
\index{monoid!anti-involution@with anti-involution}
with an anti-involution
(\ref{eg:mon-mon-with-anti-inv},~\ref{eg:mti-anti-inv}) is familially
representable and finitary, but not operadic.

The free monoid with anti-involution on a set $X$ is the set of expressions
$x_1^{\sigma_1} x_2^{\sigma_2} \cdots x_n^{\sigma_n}$ where $n\geq 0$, $x_i
\in X$, and $\sigma_i \in \{+1, -1\}$: so 
\[
TX \iso \coprod_{n\in\nat} (\{+1, -1\} \times X)^n
\iso \coprod_{i\in I} [X_i, X]
\]
where $I = \coprod_{n\in\nat} \{+1, -1\}^n$ and $X_{(\sigma_1, \ldots,
\sigma_n)} = n$.  The unit map at $X$ is 
\[
\begin{array}{rrcl}
\eta_X:	&X	&\go		&TX	\\
	&x	&\goesto	&(+1, x) \in (\{+1,-1\} \times X)^1.
\end{array}
\]
So far this is the same as for monoids with \emph{involution}
(\ref{eg:mon-mon-with-inv},~\ref{eg:mti-inv}), but the multiplication is
different: its component $\mu_X$ at $X$ is the map
\[
\coprod_{n\in\nat} \left(
\{+1,-1\} \times 
\coprod_{k\in\nat} (\{+1,-1\} \times X)^k 
\right)^n
\go
\coprod_{m\in\nat} (\{+1,-1\} \times X)^m
\]
given by
\[
\begin{array}{r}
\mu_X\left(
(\sigma_1, (\sigma_1^1, x_1^1, \ldots, \sigma_1^{k_1}, x_1^{k_1})),
\ldots,
(\sigma_n, (\sigma_n^1, x_n^1, \ldots, \sigma_n^{k_n}, x_n^{k_n}))
\right)	\\
= 
(\hat{\sigma}_1^1, \hat{x}_1^1, \ldots, 
	\hat{\sigma}_1^{k_1}, \hat{x}_1^{k_1},
\ldots,
\hat{\sigma}_n^1, \hat{x}_n^1, \ldots, 
	\hat{\sigma}_n^{k_n}, \hat{x}_n^{k_n})
\end{array}
\]
where for $i\in \{1, \ldots, n\}$,
\[
(\hat{\sigma}_i^1, \hat{x}_i^1, \ldots, 
	\hat{\sigma}_i^{k_i}, \hat{x}_i^{k_i})
=
\left\{
\begin{array}{ll}
(\sigma_i^1, x_i^1, \ldots, \sigma_i^{k_i}, x_i^{k_i})	&
\textrm{if } \sigma_i = +1	\\
(-\sigma_i^{k_i}, x_i^{k_i}, \ldots, -\sigma_i^1, x_i^1)&
\textrm{if } \sigma_i = -1.	\\
\end{array}
\right.
\]
The functor $T$ is plainly familially representable, and finitary
by~\ref{propn:fam-rep-Set-finitary}.  The unit and multiplication are
cartesian by the observations on p.~\pageref{p:cart-transf-fam-rep}.  So
the monad is familially representable and finitary, as claimed.

On the other hand, the following argument---reminiscent of a standard proof
of Brouwer's%
\index{Brouwer's fixed point theorem}
fixed point theorem---shows that it is not operadic.  For
suppose it were.  Then by Proposition~\ref{propn:ind-monad-cart}, there is
a cartesian natural transformation $\psi$ from $T$ to the free monoid monad
$S$, commuting with the monad structures.  There is also a cartesian
natural transformation $\theta: S \go T$ commuting with the monad
structures: in the notation of p.~\pageref{p:transf-fam-rep}, it is given
by the function
\[
\begin{array}{rrcl}
\phi:	&\nat	&\go		&\{+1, -1\}^n		\\
	&n	&\goesto	&(+1, \ldots, +1)
\end{array}
\]
together with the identity map $f_n: n \go n$ for each $n\in\nat$; on
categories of algebras, it induces the functor $\Set^T \go \Set^S$
forgetting anti-involutions.  So we have a commutative triangle
\[
\begin{slopeydiag}
	&		&T		&		&	\\
	&\ruTo<\theta	&		&\rdTo>\psi	&	\\
S	&		&\rTo_{\psi\theta}&		&S	\\
\end{slopeydiag}
\]
of cartesian natural transformations, each commuting with the monad
structures.  The members of the representing family of $S$ all have
different cardinalities, so any cartesian natural endomorphism of $S$ is
actually an automorphism.  Hence in the induced triangle of functors on
categories of algebras,
\[
\begin{diagram}[width=2em,height=2em,scriptlabels,noPS]
	&	&\textrm{(monoids with anti-involution)}&	&	\\
	&\ldTo<{\mr{forgetful}}&			&\luTo	&	\\
\fcat{Monoid}&	&\lTo					&	&
\fcat{Monoid},								\\
\end{diagram}
\]
the functor along the bottom is an isomorphism.  This implies that the
forgetful functor is surjective on objects.  So every monoid admits an
anti-involution, and in particular, every monoid is isomorphic to its
opposite (the same monoid with the order of multiplication reversed); but,
for instance, the monoid of endomorphisms of a two-element set does not
have this property.  This is the desired contradiction.
\end{example}

\section{Familially representable monads on presheaf categories}
\lbl{sec:fam-rep-pshf}

The cartesian monads arising in this book are typically monads on presheaf
categories.  Very often they are familially representable (in a sense
defined shortly), and this section provides some of the theory of such
monads.  After some general preliminaries, we concentrate on the special
case of presheaf categories $\ftrcat{\scat{B}^\op}{\Set}$ where $\scat{B}$
is discrete, arriving eventually at an explicit description of finitary
familially representable monads on such categories.  This will be used
in the next section to provide a link between symmetric and generalized
multicategories. 

Let $\cat{A}$ be a category and $\scat{B}$ a small category.  What should
it mean for a functor $T: \cat{A} \go \ftrcat{\scat{B}^\op}{\Set}$ to be
familially representable?  If $\scat{B}$ is discrete then the answer is
clear: for each $b\in\scat{B}$ the $b$-component 
\[
T_b = T(\dashbk)(b): \cat{A} \go \Set
\]
should be familially representable.  If $\scat{B}$ is not discrete then the
answer is not quite so obvious: we certainly do want each $T_b$ to be
familially representable, but we also want the representing family to vary
functorially as $b$ varies in $\scat{B}$.  Thus:
\begin{defn}	\lbl{defn:fam-rep-pshf}
Let $\cat{A}$ be a category and $\scat{B}$ a small category.  A functor $T:
\cat{A} \go \ftrcat{\scat{B}^\op}{\Set}$ is \demph{familially
representable}%
\index{familial representability!functor into presheaf category@of functor into presheaf category}
if there exists a functor $R$ making the diagram
\[
\begin{diagram}[height=2em]
\scat{B}^\op	&\rGet^R		&\Fam(\cat{A})		\\
		&\rdTo<{\ovln{T}}	&\dTo>{\fcat{y}}	\\
		&			&\ftrcat{\cat{A}}{\Set}	\\
\end{diagram}
\]
commute up to natural isomorphism.  Here $\ovln{T}$ is the transpose of $T$
and $\fcat{y}$ is the `Yoneda' functor of~\bref{eq:fam-rep-Yoneda}
(p.~\pageref{eq:fam-rep-Yoneda}).
\end{defn}
Note that when $\scat{B}$ is the terminal category, this is compatible with
the definition of familial representability for set-valued functors.  Note
also that since $\fcat{y}$ is full and faithful, the functor $R$ is
determined uniquely up to isomorphism, if it exists.

This definition is morally the correct one but in practical terms
needlessly elaborate: if each $T_b$ is familially representable then the
representing family \emph{automatically} varies functorially in $b$, as the
equivalence~\bref{item:fam-rep-pshf-fr-pw}
$\Leftrightarrow$~\bref{item:fam-rep-pshf-fr} in the following result
shows.
\begin{propn}	\lbl{propn:eqv-condns-fam-rep-pshf}
Let $\cat{A}$ be a category with the adjoint functor property and a
terminal object.  Let $\scat{B}$ be a small category.  The following
conditions on a functor $T: \cat{A} \go \ftrcat{\scat{B}^\op}{\Set}$ are
equivalent: 
\begin{enumerate}
\item	\lbl{item:fam-rep-pshf-conn}
$T$ preserves connected limits
\item	\lbl{item:fam-rep-pshf-wide}
$T$ preserves wide pullbacks
\item	\lbl{item:fam-rep-pshf-fr-pw}
for each $b\in \scat{B}$, the functor $T_b: \cat{A} \go \Set$ is familially
representable 
\item 	\lbl{item:fam-rep-pshf-fr}
$T$ is familially representable.
\end{enumerate}
\end{propn}
In our applications of this result, $\cat{A}$ will be a presheaf category.
Such an $\cat{A}$ does have the adjoint functor property: for by
Borceux~\cite[4.7.2]{Borx1}, it satisfies the hypotheses of the Special
Adjoint Functor Theorem.
\begin{proof}
\begin{description}
\item[\bref{item:fam-rep-pshf-conn} $\Leftrightarrow$
\bref{item:fam-rep-pshf-wide} $\Leftrightarrow$
\bref{item:fam-rep-pshf-fr-pw}] 
Since limits in a presheaf category are computed pointwise, $T$ preserves
limits of a given shape if and only if $T_b$ preserves them for each
$b\in\scat{B}$.  So these implications follow from
Proposition~\ref{propn:eqv-condns-fam-rep}.
\item[\bref{item:fam-rep-pshf-fr-pw} $\Rightarrow$ \bref{item:fam-rep-pshf-fr}]
Choose for each $b\in \scat{B}$ a representing family $(X_{b,i})_{i\in
I(b)}$ and an isomorphism
\[
\alpha_b: T_b 
\goiso 
\fcat{y}((X_{b,i})_{i\in I(b)}) = 
\coprod_{i\in I(b)} \cat{A}(X_{b,i}, \dashbk).
\]
Since $\fcat{y}$ is full and faithful, the assignment $b \goesto R(b) = 
(X_{b,i})_{i\in I(b)}$ extends uniquely to a functor $R$ such that
$\alpha$ is a natural isomorphism $\ovln{T} \goiso \fcat{y} \of R$.  
\item[\bref{item:fam-rep-pshf-fr} $\Rightarrow$
\bref{item:fam-rep-pshf-fr-pw}] 
Trivial.
\done
\end{description}
\end{proof}

A familially representable functor $\cat{A} \go
\ftrcat{\scat{B}^\op}{\Set}$ is determined by a functor $R: \scat{B}^\op
\go \Fam(\cat{A})$.  Explicitly, such a functor $R$ consists of
\begin{itemize}
\item a functor $I: \scat{B}^\op \go \Set$ 
\item for each $b \in \scat{B}$, a family $(X_{b,i})_{i\in I(b)}$ of objects
of $\cat{A}$
\item for each $b' \goby{g} b$ in $\scat{B}$ and each $i\in I(b)$, a map
\[
X_{g,i}: X_{b', (Ig)(i)} \go X_{b,i}
\]
in $\cat{A}$
\end{itemize}
such that
\begin{itemize}
\item if $b'' \goby{g'} b' \goby{g} b$ in $\scat{B}$ and $i\in I(b)$ then
$X_{g'\sof g, i} = X_{g, i} \of X_{g', (Ig)(i)}$
\item if $b\in \scat{B}$ then $X_{1_b,i} = 1_{X_{b,i}}$.
\end{itemize}
The resulting functor $T: \cat{A} \go \ftrcat{\scat{B}^\op}{\Set}$ is given
at $X \in \cat{A}$ and $b \in \scat{B}$ by
\[
(TX)(b) = \coprod_{i\in I(b)} \cat{A} (X_{b,i}, X).
\]

A good example of a familially representable functor into a presheaf
category is the free strict $\omega$-category%
\index{familial representability!free strict omega-category functor@of free strict $\omega$-category functor}
functor $T$
discussed in Chapter~\ref{ch:globular} and Appendix~\ref{app:free-strict}.
Here is the 1-dimensional version.

\begin{example}
Let $\scat{B} = (0 \parpair{\sigma}{\tau} 1)$, so that
$\ftrcat{\scat{B}^\op}{\Set}$ is the category \fcat{DGph} of directed
graphs; let $\cat{A} = \fcat{DGph}$ too.  Then the free%
\index{category!free (fc)@free ($\fc$)}
category functor
$T: \fcat{DGph} \go \fcat{DGph}$ is familially representable.  The directed
graph $I$ is defined by $I(0)=1$ and $I(1)=\nat$.  The families
$(X_{b,i})_{i\in I(b)}$ ($b\in \scat{B}$) are defined as follows.  For
$b=1$ and $i\in\nat$, let $X_{1,i}$ be the graph
\[
\gfsts{0} \gones{} \gfbws{1} \gones{} 
\diagspace \cdots \diagspace 
\gones{} \glsts{i}
\]
with $(i+1)$ vertices and $i$ edges.  For $b=0$, writing $1=\{j\}$, let
\[
X_{0,j} = \gzeros{} = X_{1,0}.
\]
The maps $0 \parpair{\sigma}{\tau} 1$ in $\scat{B}$ induce, for each
$i\in\nat$, the graph maps
\[
X_{0,j} \parpairu X_{1,i}
\]
sending $X_{0,j}$ to the first and last vertex of $X_{1,i}$, respectively.
This data does represent $T$: if $X$ is a directed graph then
\[
\coprod_{i\in I(1)} \fcat{DGph}(X_{1,i}, X)	
=
\coprod_{i\in \nat} 
\{ \textrm{strings of } i \textrm{ arrows in } X \}
= 
(TX)_1
\]
and
\[
\coprod_{j\in I(0)} \fcat{DGph}(X_{0,j}, X)	
=
X_0 
=
(TX)_0,
\]
as required.
\end{example}

A short digression: the \demph{Artin%
\lbl{p:Artin}
gluing}%
\index{Artin gluing}%
\index{gluing}
of a functor $T: \cat{A} \go \cat{B}$ is the category
$\cat{B}\gluing T$ in which an object is a triple $(X, Y, \pi)$ with
$X\in\cat{A}$, $Y\in\cat{B}$, and $\pi: Y \go TX$ in $\cat{B}$, and a map
$(X, Y, \pi) \go (X', Y', \pi')$ is a pair of maps $(X\go X', Y\go Y')$
making the evident square commute.
\begin{propn}[Carboni--Johnstone]	\lbl{propn:Artin-gluing}%
\index{Carboni, Aurelio}%
\index{Johnstone, Peter}
Let $\cat{A}$ and $\cat{B}$ be presheaf categories and $T: \cat{A} \go
\cat{B}$ a functor.  Then $\cat{B}\gluing T$ is a presheaf category if
and only if $T$ preserves wide pullbacks.  
\end{propn}
\begin{proof}
See Carboni and Johnstone~\cite[4.4(v)]{CJ}.  Their proof of `if' is a
little roundabout, so the following sketch of a direct version may be of
interest.

By Proposition~\ref{propn:eqv-condns-fam-rep-pshf}, $T$ is familially
representable.  Suppose that $\cat{A} = \ftrcat{\scat{A}^\op}{\Set}$ and
$\cat{B} = \ftrcat{\scat{B}^\op}{\Set}$, and represent $T$ by the families
$(X_{b,i})_{i\in I(b)}$, as above.  An object of $\cat{B}\gluing T$
consists of functors $X: \scat{A}^\op \go \Set$ and $Y: \scat{B}^\op \go
\Set$ and a map
\[
\pi_b: Yb \go \coprod_{i\in I(b)} 
\ftrcat{\scat{A}^\op}{\Set} (X_{b,i}, X)
\]
for each $b\in \scat{B}$, satisfying naturality axioms.  Equivalently, it
consists of
\begin{enumerate}
\item a functor $X: \scat{A}^\op \go \Set$,
\item a family $(Y(b,i))_{b\in\scat{B}, i\in I(b)}$ of sets, functorial in
$b$, and
\item	\lbl{item:gluing-object-three}
for each $b\in\scat{B}$, $i\in I(b)$ and $y\in Y(b,i)$, a natural
transformation $X_{b,i} \go X$, 
\end{enumerate}
satisfying axioms; indeed,~\bref{item:gluing-object-three} can equivalently
be replaced by
\begin{enumerate}
\newcommand{\temptheenumi}{\theenumi}
\renewcommand{\theenumi}{\ref{item:gluing-object-three}$'$}
\item
\renewcommand{\theenumi}{\temptheenumi} 
for each $b\in\scat{B}$, $i\in I(b)$, $a\in\scat{A}$ and $x\in X_{b,i}(a)$,
a function $Y(b,i) \go X(a)$.
\end{enumerate}
Equivalently, it is a functor $\scat{C}^\op \go \Set$, where $\scat{C}$ is
the category whose object-set $\scat{C}_0$ is the disjoint union
\[
\scat{C}_0 = \scat{A}_0 + \{(b, i) \such b\in\scat{B}_0, i\in I(b)\},
\]
whose maps are given by
\begin{eqnarray*}
\scat{C}(a', a)		&=	&\scat{A}(a', a),			\\
\scat{C}((b,i), (b',i'))&=	&\{g\in\scat{B}(b,b') \such (Ig)(i')=i\},\\
\scat{C}(a, (b,i))	&=	&X_{b,i}(a),				\\
\scat{C}((b,i), a)	&=	&\emptyset
\end{eqnarray*}
($a, a' \in \scat{A}_0$, $b, b'\in \scat{B}_0$, $i\in I(b)$, $i'\in
I(b')$), and whose composition and identities are the evident ones.
\done  
\end{proof}

\index{power of Set@power of $\Set$|(}%
\index{Set, power of@$\Set$, power of|(}
Returning to the main story of familially representable functors into
presheaf categories, let us restrict our attention to functors of the form
\[
T: \Set^C \go \Set^B
\]
where $C$ and $B$ are sets.  This makes the calculations much easier but
still provides a wide enough context for our applications. 

A familially representable functor $T: \Set^C \go \Set^B$ consists of 
\begin{itemize}
\item a family $(I(b))_{b\in B}$ of sets
\item for each $b\in B$, a family $(X_{b,i})_{i\in I(b)}$ of objects of
$\Set^C$, 
\end{itemize}
and the actual functor $T$ is then given by 
\[
(TX)(b) 
\iso
\coprod_{i\in I(b)} \Set^C(X_{b,i}, X) 
\iso
\coprod_{i\in I(b)} \prod_{c \in C} [X_{b,i}(c), X(c)]
\]
for each $X\in \Set^C$ and $b\in B$.

Let $T': \Set^C \go \Set^B$ be another such functor, with representing
families $(X_{b,i'})_{i'\in I'(b)}$ ($b\in B$).  A natural transformation%
\index{transformation!familially representable functors@of familially representable functors}
\[
\Set^C \ctwomult{T'}{T}{} \Set^B
\]
consists merely of a natural transformation
\[
\Set^C \ctwomult{T'_b}{T_b}{} \Set
\]
for each $b\in B$.  So by the results of the previous section, a
transformation $T' \go T$ is described by
\begin{itemize}
\item for each $b\in B$, a function $\phi_b: I'(b) \go I(b)$
\item for each $b\in B$ and $i' \in I'(b)$, a map 
$
f_{b,i'}: X_{b, \phi_b(i')} \go X'_{b,i'}
$
in $\Set^C$,
\end{itemize}
and the induced map $(T'X)(b) \go (TX)(b)$ is the composite
\[
\begin{array}{rl}
						&
\coprod_{i'\in I'(b)} \Set^C (X'_{b,i'}, X)	\\
\goby{\coprod_{i'\in I'(b)} f_{b,i'}^*}		&
\coprod_{i'\in I'(b)} \Set^C (X_{b, \phi_b(i')}, X)	\\
\goby{\mr{canonical}}				&
\coprod_{i\in I(b)} \Set^C (X_{b,i}, X)
\end{array}
\]
($X\in \Set^C$, $b\in B$).  The transformation is cartesian%
\index{transformation!cartesian}
if and only if
each map $f_{b,i'}$ is an isomorphism.

Just as for functors $\Set \go \Set$, the theory tells us that the class of
familially representable functors between presheaf categories is closed
under composition.  Working out the representing family of a composite
seems extremely complicated for presheaf categories in general, but is
manageable in our restricted context of direct powers of $\Set$.  So, take
sets $D$, $C$ and $B$ and familially representable functors
\[
\Set^D \goby{S} \Set^C \goby{T} \Set^B
\]
given by
\[
(SW)(c) \iso \coprod_{h\in H(c)} \Set^D (W_{c,h}, W),
\diagspace
(TX)(b) \iso \coprod_{i\in I(b)} \Set^C (X_{b,i}, X)
\]
($W\in \Set^D$, $X\in \Set^C$, $c\in C$, $b\in B$).  Then for $W\in \Set^D$
and $b\in B$,
\begin{eqnarray*}
(TSW)(b)	&\iso	&
\coprod_{i\in I(b)} \Set^C (X_{b,i}, SW)			\\
	&\iso	&
\coprod_{i\in I(b)} \prod_{c\in C} 
[X_{b,i}(c), \coprod_{h\in H(c)} \Set^D (W_{c,h}, W)]		\\
	&\iso	&
\coprod_{i\in I(b)} \prod_{c\in C} 
\coprod_{g\in [X_{b,i}(c), H(c)]} \prod_{x\in X_{b,i}(c)}
\Set^D (W_{c, g(x)}, W)						\\
	&\iso	&
\coprod_{i\in I(b)} \coprod_{\gamma\in \Set^C(X_{b,i}, H)} 
\prod_{c\in C} \prod_{x\in X_{b,i}(c)}
\Set^D (W_{c, \gamma_c(x)}, W)						\\
	&\iso	&
\coprod_{(i,\gamma)\in J(b)} \Set^D (Y_{b,i,\gamma}, W),
\end{eqnarray*}
where 
\begin{eqnarray*}
J(b)		&=	&
\coprod_{i\in I(b)} \Set^C (X_{b,i}, H) \iso (TH)(b),	\\
Y_{b,i,\gamma}	&=	&
\coprod_{c\in C, x\in X_{b,i}(c)} W_{c,\gamma_c(x)}.
\end{eqnarray*}
The identity functor on $\Set^B$ is also familially representable: for
$X\in \Set^B$, 
\[
X(b) \iso \coprod_{i\in 1} \Set^B (\delta_b, X)
\]
where, treating $B$ as a discrete category, $\delta_b = B(\dashbk,b) \in
\Set^B$.

Assembling these descriptions gives a description of monads on $\Set^B$
whose functor parts are familially representable and whose natural
transformation parts are cartesian---\demph{familially representable
monads}%
\index{familial representability!monad on power of Set@of monad on power of $\Set$}%
\index{monad!familially representable}
on $\Set^B$, as we call them.  Such a monad consists of:
\begin{itemize}
\item for each $b\in B$, a set $I(b)$ and a family $(X_{b,i})_{i\in I(b)}$
of objects of $\Set^B$, inducing the functor
\[
\begin{array}{rl}
T: 		&\Set^B \go \Set^B,		\\
(TX)(b) = 	&\coprod_{i\in I(b)} \Set^B (X_{b,i}, X)
\end{array}
\]
\item for each $b\in B$, a function
\[
m_b: 
\{(i,\gamma) \such i\in I(b), \gamma\in \Set^B (X_{b,i}, X) \}
\go
I(b),
\]
and for each $b\in B$, $i\in I(b)$ and $\gamma\in\Set^B (X_{b,i}, X)$, an
isomorphism 
\[
X_{b, m_b(i,\gamma)} 
\goiso
\coprod_{c\in B, x\in X_{b,i}(c)} X_{c,\gamma_c(x)},
\]
inducing the natural transformation $\mu: T^2 \go T$ whose component
\[
\mu_{X,b}: (T^2 X)(b) \go (TX)(b)
\]
is the composite
\[
\begin{array}{rl}
	&
\coprod_{i\in I(b), \gamma\in \Set^B (X_{b,i}, I)}
\Set^B (\coprod_{c\in B, x\in X_{b,i}(c)} X_{c,\gamma_c(x)}, X)	\\
\goiso	&
\coprod_{i\in I(b), \gamma\in \Set^B (X_{b,i}, I)}
\Set^B (X_{b, m_b(i,\gamma)}, X)				\\
\goby{\mr{canonical}}	&
\coprod_{k\in I(b)} \Set^B (X_{b,k}, X)
\end{array}
\]
\item for each $b\in B$, an element $e_b \in I(b)$ such that $X_{b, e_b}
\iso \delta_b$, inducing the natural transformation $\eta: 1 \go T$ whose
component
\[
\eta_{X,b}: X(b) \go (TX)(b)
\]
is the composite
\[
X(b) 
\goiso 
\Set^B (X_{b,e_b}, X)
\rIncl
\coprod_{i\in I(b)} \Set^B (X_{b,i}, X),
\]
\end{itemize}
such that $\mu$ and $\eta$ obey associativity and unit laws.  This is
complicated, but rest assured that it gets no worse.

The aim of this section was to describe finitary familially representable
monads on categories of the form $\Set^B$, and we are nearly there.  All
that remains is `finitary'.%
\index{finitary}
\begin{propn}	\lbl{propn:fr-pshf-finitary}
Let $C$ and $B$ be sets.  Let $T: \Set^C \go \Set^B$ be a familially
representable functor with representing families $(X_{b,i})_{i\in I(b)}$
($b\in B$).  The following
% conditions on $T$ 
are equivalent:
\begin{enumerate}
\item	\lbl{item:fin-fam-rep-pshf-finitary}
$T$ is a finitary functor
\item	\lbl{item:fin-fam-rep-pshf-finite-power}
$\coprod_{c\in C} X_{b,i} (c)$ is a finite set for each $b\in B$ and $i\in
I(b)$. 
\end{enumerate}
\end{propn}
\begin{proof}
$T$ is finitary if and only if the functor
\[
T_b = \coprod_{i\in I(b)} \Set^C (X_{b,i}, \dashbk): 	
\Set^C	\go \Set
\]
is finitary for each $b\in B$.  By the arguments in the proof of
Proposition~\ref{propn:fam-rep-Set-finitary}, $T_b$ is finitary if and only
if the functor
\[
\Set^C (X_{b,i}, \dashbk): \Set^C \go \Set
\]
is finitary for each $i\in I(b)$.  But $\Set^C (X_{b,i}, \dashbk)$ is
finitary if and only if $\coprod_{c\in C} X_{b,i}(c)$ is a finite set
(Ad\'amek%
\index{Adamek, Jiri@Ad\'amek, Ji\v{r}\'\i}
and Rosick\'y~\cite[p.~9]{AR}).%
\index{Rosicky, Jiri@Rosick\'y, Ji\v{r}\'\i}
 \done
\end{proof}

\section{Cartesian structures from symmetric structures}
\lbl{sec:cart-sym}%
\index{monoid!commutative|(}

A commutative monoid is a structure in which every finite family of
elements has a well-defined sum.  So if $T = \coprod_{i\in I} [X_i,
\dashbk]$ is any finitary familially representable monad on $\Set$ then any
commutative monoid $A$ is naturally a $T$-algebra via the map
\[
TA = \coprod_{i\in I} [X_i, A] \go A
\]
whose $i$-component sends a family $(a_x)_{x\in X_i}$ to $\sum_{x\in X_i}
a_x$.  

This is the simplest case of the theme of this section: how symmetric
structures give rise to cartesian structures.  We prove two main results.
The first is that if $B$ is any set, $T$ any finitary familially
representable monad on $\Set^B$, and $A$ any commutative monoid, then the
constant family $(A)_{b\in B}$ is naturally a $T$-algebra.  (Just now we
looked at the case $B=1$.)  It follows that any commutative monoid is
naturally a $T_n$-operad, where $T_n$ is the $n$th opetopic
monad~(\ref{sec:opetopes}).  We can also extract a rigorous definition of
the set-with-multiplicities of $n$-opetopes underlying a given
$n$-dimensional opetopic pasting diagram.

The second main result is that any symmetric multicategory $A$ gives rise
to a $T$-multicategory.  Here $T$ is, again, any finitary familially
representable monad on $\Set^B$ for some set $B$, and the object-of-objects
of the induced $T$-multicategory is $(A_0)_{b\in B}$.  So this is like the
result for commutative monoids but one level up; it states that symmetric
multicategories play some kind of universal role for generalized
multicategories, despite (apparently) not being generalized multicategories
themselves.  A corollary is that any symmetric multicategory is naturally a
$T_n$-multicategory for each $n$, previously stated as
Theorem~\ref{thm:sm-opetopic} and pictorially very plausible---see
p.~\pageref{eq:Tn-sym-arrow}.

We start with commutative monoids.  The first main result mentioned above
is reasonably clear informally.  I shall, however, prove it with some care
in preparation for the proof of the second main result, which involves many
of the same thoughts in a more complex setting.  

Notation: for any set $B$ there are adjoint functors
\[
\begin{diagram}[scriptlabels]
\Set^B	&\pile{\rTo^\Sigma_\bot\\ \lTo_\Delta}	&\Set,\\
\end{diagram}%
\glo{Sigmacoproduct}
\]
where $\Sigma X = \coprod_{b\in B} Xb$ and $(\Delta A)(b) = A$.  

\begin{thm}	\lbl{thm:cm-main}
Let $B$ be a set and $T$ a finitary familially representable monad on
$\Set^B$.  Then there is a canonical functor
\[
\fcat{CommMon} \go (\Set^B)^T
\]
making the diagram
\[
\begin{diagram}[size=2em]
\fcat{CommMon}		&\rTo		&(\Set^B)^T		\\
\dTo<{\mr{forgetful}}	&
					&
\dTo>{\mr{forgetful}}	\\
\Set			&\rTo_\Delta	&\Set^B			\\
\end{diagram}
\]
commute.
\end{thm}
(The standard notation is awkward here: $(\Set^B)^T$ is the category of
algebras for the monad $T$ on the category $\Set^B$.)
\begin{proof}
To make the link between commutative monoids and familial representability
we use the fat%
\index{monoid!commutative!fat}
commutative monoids of~\ref{sec:comm-mons}, which come
equipped with explicit operations for summing arbitrary finite families
(not just finite sequences) of elements.  By Theorem~\ref{thm:fat-cm-eqv},
it suffices to prove the present theorem with `\fcat{CommMon}' replaced by
`\fcat{FatCommMon}'.

Represent the functor $T$ by families $(X_{b,i})_{i\in I(b)}$ ($b\in B$),
the unit $\eta$ of the monad by $e$, and the multiplication $\mu$ by $m$
and a nameless isomorphism, as in the description at the end
of~\ref{sec:fam-rep-pshf}.

For any set $A$ and any $b\in B$ we have
\[
(T\Delta A)(b) 
\iso 
\coprod_{i\in I(b)} (X_{b,i}, \Delta A) 
\iso
\coprod_{i\in I(b)} [\Sigma X_{b,i}, A],
\]
and by Proposition~\ref{propn:fr-pshf-finitary}, the set $\Sigma X_{b,i}$ is
finite.  So given a fat commutative monoid $A$, we have a map
\[
\theta^A_B: (T\Delta A)(b) \go (\Delta A)(b) = A
\]
whose component at $i\in I(b)$ is summation
\[
\sum_{\Sigma X_{b,i}}: [\Sigma X_{b,i}, A] \go A. 
\]
If we can show that $\theta^A$ is a $T$-algebra structure on $\Delta A$
then we are done: the functoriality is trivial.

For the multiplication axiom we have to show that the square
\[
\begin{diagram}[height=2em]
(T^2 \Delta A)(b)	&\rTo^{\mu_{\Delta A, b}}	&(T\Delta A)(b)	\\
\dTo<{(T\theta^A)_b}	&			&\dTo>{\theta^A_b}	\\
(T\Delta A)(b)		&\rTo_{\theta^A_b}		&A		\\
\end{diagram}
\]
commutes, for each $b\in B$.  We have
\[
(T^2 \Delta A)(b) 
\iso 
\coprod_{(i,\gamma)\in J(b)} 
[\Sigma Y_{b,i,\gamma}, A]
\]
where
\begin{eqnarray*}
J(b)	&=	&
\{ (i,\gamma) \such i\in I(b), \gamma\in \Set^B(X_{b,i}, I) \},	\\
Y_{b,i,\gamma}	&=	&
\coprod_{c\in B, x\in X_{b,i}(c)} X_{c,\gamma_c(x)}.
\end{eqnarray*}
Let $i\in I(b)$ and $\gamma\in \Set^B(X_{b,i}, I)$.  Write $k = m_b(i,\gamma)
\in I$; then part of the data for the monad is an isomorphism $X_{b,k} \goiso
Y_{b,i,\gamma}$.  The clockwise route around the square has
$(i,\gamma)$-component
\[
\begin{diagram}[height=2em]
[\Sigma Y_{b,i,\gamma}, A]	&\rTo^\diso	&[\Sigma X_{b,k}, A]	\\
				&		&
\dTo>{\sum_{\Sigma X_{b,k}}}						\\
				&		&A,			\\
\end{diagram}
\]
and by Lemma~\ref{lemma:fat-cm}, this is just $\sum_{\Sigma
Y_{b,i,\gamma}}$.  For the anticlockwise route, an element of $[\Sigma
Y_{b,i,\gamma}, A]$ is a family $(a_{c,x,d,w})_{c,x,d,w}$ indexed over
\[
c\in B, x\in X_{b,i}(c), d\in B, w\in X_{c, \gamma_c(x)} (d).
\]
For each $c\in B$ and $x\in X_{b,i}(c)$ we therefore have a family
\[
(a_{c,x,d,w})_{d,w} \in [\Sigma X_{x, \gamma_c(x)}, A],
\]
and 
\[
\theta^A_c ((a_{c,x,d,w})_{d,w}) = 
\sum_{d,w} a_{c,x,d,w} \in A,
\]
so
\[
(T\theta^A)_b (a_{c,x,d,w})_{c,x,d,w} =
\left( \sum_{d,w} a_{c,x,d,w} \right)_{c,x}
\in
[\Sigma X_{b,i}, A].
\]
So the anticlockwise route sends $(a_{c,x,d,w})_{c,x,d,w}$ to $\sum_{c,x}
\sum_{d,w} a_{c,x,d,w}$, and the clockwise route sends it to
$\sum_{c,x,d,w} a_{c,x,d,w}$; the two are equal.

For the unit axiom we have to show that
\[
\begin{diagram}[height=2em]
(\Delta A)(b)	&\rTo^{\eta_{\Delta A, b}}	&(T\Delta A)(b)		\\
		&\rdTo<1			&\dTo>{\theta^A_b}	\\
		&				&A			\\
\end{diagram}
\]
commutes, for each $b\in B$.  Write $k = e_b\in I$; we know that $X_{b,k}
\iso \delta_b$, and so $\Sigma X_{b,k}$ is a one-element set.  So if $a\in A$
then 
\[
\theta^A_b (\eta_{\Delta A, b} (a))
=
\theta^A_b ((a)_{x\in\Sigma X_{b,k}})
=
\sum_{x\in\Sigma X_{b,k}} a
=
a,
\]
as required.
\done
\end{proof}

\begin{example}	\lbl{eg:cm-ftr-opetopic}
Let $n\in\nat$; take $B$ to be the set $O_{n+1}$ of $(n+1)$-opetopes and
$T$ to be the $(n+1)$th opetopic%
\index{opetopic!monad}
monad $T_{n+1}$ on $\Set/O_{n+1} \eqv
\Set^{O_{n+1}}$.  By Proposition~\ref{propn:free-refined} and induction,
$T_{n+1}$ is finitary and preserves wide pullbacks, so is familially
representable by Proposition~\ref{propn:eqv-condns-fam-rep-pshf}.
Theorem~\ref{thm:cm-main} then gives a canonical functor
\[
\fcat{CommMon} 
\go 
(\Set/O_{n+1})^{T_{n+1}}
\eqv
T_n \hyph\Operad,
\]
proving Theorem~\ref{thm:cm-ftr-opetopic}.
\end{example}

\begin{example}	\lbl{eg:cm-ftr-Set}
For any finitary familially representable monad $T$ on $\Set$, there is a
canonical functor $\fcat{CommMon} \go \Set^T$ preserving underlying sets.
Most such $T$'s that we have met have been operadic,%
\index{monad!operadic}
in which case
something stronger is true: there is a canonical functor from the category
of not-necessarily-commutative monoids to $\Set^T$, as noted after
Proposition~\ref{propn:ind-monad-cart}.  But, for instance, if $T$ is the
monad corresponding to the theory of monoids%
\index{monoid!anti-involution@with anti-involution}
with an
anti-involution (Example~\ref{eg:fam-rep-not-opdc}) then the commutativity is
necessary.  Concretely, Theorem~\ref{thm:cm-main} produces the functor
\[
\begin{array}{rcl}
\fcat{CommMon} 	&\go 	&(\textrm{monoids with an anti-involution}) 	\\
(A,\cdot,1)	&\goesto&(A, \cdot, 1, \blank^\circ)			\\
\end{array}
\]
defined by taking the anti-involution $\blank^\circ$ to be the identity;
that this does define an anti-involution is exactly commutativity.
\end{example}

Write $M$%
\glo{Mfreecommmon}
for the free commutative monoid monad on $\Set$.
Lemma~\ref{lemma:lax-map-mnds-is-ftr} tells us that
Theorem~\ref{thm:cm-main} is equivalent to:
\begin{cor}	\lbl{cor:cm-lax-map}
Let $B$ be a set and $T$ a finitary familially representable monad on
$\Set^B$.  Then there is a canonical natural transformation
\[
\begin{diagram}
\Set^B 		&\rTo^T		&\Set^B		\\
\uTo<\Delta	&\sent \psi^T	&\uTo>\Delta	\\
\Set		&\rTo_M		&\Set		\\
\end{diagram}
\]
with the property that $(\Delta, \psi^T)$ is a lax map of monads $(\Set,M)
\go (\Set^B,T)$.  
\done
\end{cor}

When $B=1$, this corollary provides a natural transformation $\psi^T: T \go
M$ commuting with the monad structures.  It is straightforward to check
that if $\alpha: T' \go T$ is a cartesian natural transformation commuting
with the monad structures then $\psi^T \of \alpha = \psi^{T'}$.  So $M$ is
the vertex (codomain) of a cone on the inclusion functor
\[
\begin{array}{rl}
	&(\textrm{finitary familially representable monads on } \Set	\\
	&+ \textrm{ cartesian transformations commuting with the monad
	structures}) 							\\
\rIncl	&(\textrm{monads on } \Set					\\
	&+ \textrm{ transformations commuting with the monad structures}), 
\end{array}
\]
in which the coprojections are the transformations $\psi^T$.  In fact, it
is a colimit cone: the theory of commutative monoids plays a universal role
for finitary familially representable monads on $\Set$, despite not being
familially representable itself.  Since this fact will not be used, I will
not prove it; the main tactic is to consider the free algebraic theory on a
single $n$-ary operation.

As explained on p.~\pageref{p:colax-lax-mate}, we can translate between lax
and colax maps of monads using mates.  Applying this to
Corollary~\ref{cor:cm-lax-map} gives:
\begin{cor}	\lbl{cor:cm-colax-map}
Let $B$ be a set and $T$ a finitary familially representable monad on
$\Set^B$.  Then there is a canonical natural transformation
\[
\begin{diagram}
\Set^B 		&\rTo^T		&\Set^B		\\
\dTo<\Sigma	&\swnt \phi^T	&\dTo>\Sigma	\\
\Set		&\rTo_M		&\Set		\\
\end{diagram}
\]
with the property that $(\Sigma, \phi^T)$ is a colax map of monads
$(\Set^B,T) \go (\Set,M)$.
\end{cor}
\begin{proof}
Take $\phi^T$ to be the mate of $\psi^T$ under the adjunction $\Sigma \ladj
\Delta$. 
\done
\end{proof}

\begin{example}	\lbl{eg:constituents}
Let $n\in\nat$ and take $B=O_n$ and $T=T_n$,%
\index{opetopic!monad}
as in
Example~\ref{eg:cm-ftr-opetopic} but with the indexing shifted.  An object
$X$ of $\Set^{O_n} \eqv \Set/O_n$ is thought of as a set of labelled
$n$-opetopes.  An element of $T_n X$ (or rather, of its underlying set
$\Sigma T_n X$) is then an $X$-labelled $n$-pasting diagram;%
\index{pasting diagram!opetopic}
on the other
hand, an element of $M \Sigma X$ is a finite set-with-multiplicities of
labels (disregarding shapes completely).  So there ought to be a forgetful
function $\Sigma T_n X \go M \Sigma X$, and there is: $\phi^{T_n}_X$.

When $X$ is the terminal object of $\Set^{O_n}$ we have $\Sigma X = O_n$
and $\Sigma T_n X = O_{n+1}$, so $\phi^{T_n}_X$ is a map $O_{n+1} \go M
O_n$.  This sends an $n$-pasting diagram to the set-with-multiplicities%
\index{pasting diagram!opetopic!ordering of}
of
its constituent $n$-opetopes: for example, the $2$-pasting
diagram~\bref{diag:bigger-two-pd} on p.~\pageref{diag:bigger-two-pd}
(stripped of its labels) is sent to the set-with-multiplicities
\[
\left[ \ 
\topez{}{\Downarrow},\ 
\topea{}{}{\Downarrow},\ 
\topec{}{}{}{}{\Downarrow},\ 
\topec{}{}{}{}{\Downarrow},\ 
\topeeu{\Downarrow}\ 
\right]
\]
of 2-opetopes.
\end{example}%
\index{monoid!commutative|)}

\index{multicategory!symmetric vs. generalized@symmetric \vs.\ generalized|(}
We now consider symmetric multicategories.  Much of what follows is similar
to what we did for commutative monoids but at a higher level of complexity.

\begin{thm}	\lbl{thm:sm-main}
Let $B$ be a set and $T$ a finitary familially representable monad on
$\Set^B$.  Then there is a canonical functor 
\[
\fcat{FatSymMulticat} \go T\hyph\Multicat%
\index{multicategory!symmetric!fat}
\]
making the diagram
\[
\begin{diagram}[size=2em]
\fcat{FatSymMulticat}	&\rTo		&T\hyph\Multicat	\\
\dTo<{\blank_0}		&		&\dTo>{\blank_0}	\\
\Set			&\rTo_\Delta	&\Set^B			\\
\end{diagram}
\]
commute, where in both cases $\blank_0$ is the functor assigning to a
multicategory its object of objects.
\end{thm}
Before we prove the Theorem let us gather a corollary and some examples.
\begin{cor}
Theorem~\ref{thm:sm-main} holds with `\fcat{FatSymMulticat}' replaced by
`\fcat{SymMulticat}' and `canonical' by `canonical up to isomorphism'.  
\end{cor}
\begin{proof}
Follows from Theorem~\ref{thm:fat-sm-eqv}.
\done
\end{proof}

\begin{example}	\lbl{eg:sm-opetopic}
For any $n\in\nat$, let $B$ be the set $O_n$ of $n$-opetopes and let $T_n$
be the $n$th opetopic%
\index{opetopic!monad}
monad.  Then, as observed
in~\ref{eg:cm-ftr-opetopic}, $T_n$ is finitary and familially
representable.  So the Corollary produces a functor
\[
\fcat{SymMulticat} \go T_n\hyph\Multicat,
\]
proving Theorem~\ref{thm:sm-opetopic}. 
\end{example}

\begin{example}
Taking $B=1$, any symmetric multicategory is naturally a $T$-multicategory
for any finitary familially representable monad $T$ on $\Set$.  If $T$ is
operadic%
\index{monad!operadic}
then we do not need the symmetries: the canonical natural
transformation from $T$ to the free monoid monad induces a functor
\[
\Multicat \go T\hyph\Multicat.
\]
This is analogous to the situation for commutative monoids described
in~\ref{eg:cm-ftr-Set}.
\end{example}

We finish with a proof of Theorem~\ref{thm:sm-main}.  Undeniably it is
complicated, but there is no great conceptual difficulty; the main struggle
is against drowning in notation.  It may help to keep
Example~\ref{eg:sm-opetopic} in mind.

\begin{prooflike}{Proof of Theorem~\ref{thm:sm-main}}
Let $T$ be a finitary familially representable monad on $\Set^B$.  As in
the description just before Proposition~\ref{propn:fr-pshf-finitary},
represent $T$ by families $(X_{b,i})_{i\in I(b)}$ ($b\in B$), the unit
$\eta$ by $e$, and the multiplication $\mu$ by $m$ and a bijection
\[
s_{b,i,\gamma}: 
X_{b, m_b(i, \gamma)}
\goiso
\coprod_{c\in B, x\in X_{b,i}(c)} X_{c, \gamma_c(x)}
\]
for each $b\in B$, $i\in I$, and $\gamma\in \Set^B(X_{b,i}, I)$.  Given
such $b$, $i$, and $\gamma$, and given $d\in B$, write
\[
s_{b,i,\gamma,d}(v) = (c(v), x(v), w(v))
\] 
for any $v\in X_{b, m_b(i, \gamma)}(d)$; so here $c(v) \in B$, $x(v) \in
X_{b,i}(c)$, and $w \in X_{c, \gamma_c(x)}(d)$.

Since $T$ is finitary, Proposition~\ref{propn:fr-pshf-finitary} tells us
that the set $\Sigma X_{b,i}$ is finite for each $b\in B$ and $i\in I(b)$.

Let $P$ be a fat symmetric multicategory.  Our main task is to define from
$P$ a $T$-multicategory $C$.  

Define $C_0 = \Delta P_0$ (as we must).  Then for each $b\in B$ we have
\begin{itemize}
\item $C_0(b) = P_0$
\item $(TC_0)(b) = \coprod_{i\in I(b)} [\Sigma X_{b,i}, P_0]$, so an
element of $(TC_0)(b)$ consists of an element $i\in I(b)$ together with a
family $(q_{c,x})_{c,x}$ of objects of $P$ indexed over $c\in B$ and $x\in
X_{b,i}(c)$.
\end{itemize}

Define the $T$-graph
\[
\begin{slopeydiag}
	&		&C_1	&		&	\\
	&\ldTo<\dom	&	&\rdTo>\cod	&	\\
TC_0	&		&	&		&C_0	\\
\end{slopeydiag}
\]
by declaring an element of $C_1(b)$ with domain $(i, (q_{c,x})_{c,x}) \in
(TC_0)(b)$ and codomain $r\in (C_0)(b)$ to be a map
\[
(q_{c,x})_{c,x} \goby{\theta} r
\]
in $P$.

To define $\comp: C_1\of C_1 \go C_1$ we first have to compute $C_1 \of
C_1$.  With a little effort we find that an element of $(C_1 \of C_1)(b)$
consists of the following data (Fig.~\ref{fig:sm-comp-data}):
\begin{figure}
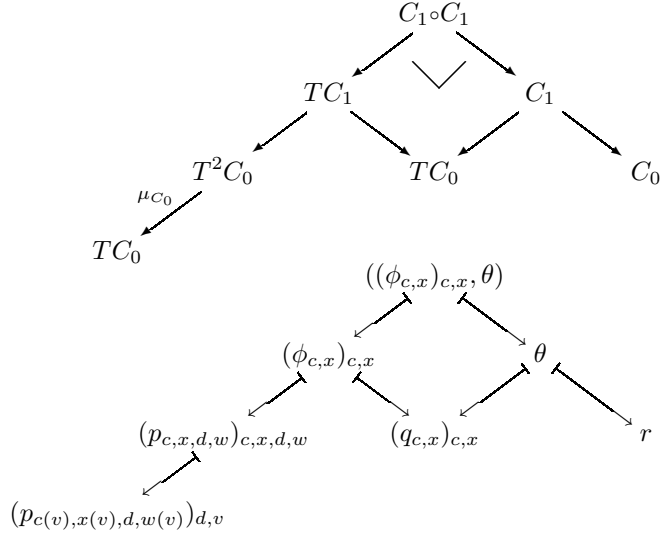

\[
\begin{array}{r}
\begin{diagram}[width=2em,height=1.5em,scriptlabels,tight]
   &       &   &       &   &       &C_1\of C_1\Spbk&& &       &   \\
   &       &   &       &   &\ldTo  &      &\rdTo  &   &       &   \\
   &       &   &       &TC_1&      &      &       &C_1&       &   \\
   &       &   &\ldTo  &   &\rdTo  &      &\ldTo  &   &\rdTo  &   \\
   &       &T^2 C_0&   &   &       &TC_0  &       &   &       &C_0\\
   &\ldTo<{\mu_{C_0}}&&&   &       &      &       &   &       &   \\
TC_0&      &   &       &   &       &      &       &   &       &   \\
\end{diagram}
\\
\begin{diagram}[width=2em,height=1.5em,scriptlabels,tight]
   &       &   &       &   &       &((\phi_{c,x})_{c,x},\theta)&&&&\\
   &       &   &       &   &\ldGoesto&    &\rdGoesto& &       &   \\
   &       &   &       &(\phi_{c,x})_{c,x}&&&     &\theta&    &   \\
   &       &   &\ldGoesto& &\rdGoesto&    &\ldGoesto& &\rdGoesto& \\
   &       &(p_{c,x,d,w})_{c,x,d,w}&&&&(q_{c,x})_{c,x}&&&     &r  \\
   &\ldGoesto& &       &   &       &      &       &   &       &   \\
(p_{c(v),x(v),d,w(v)})_{d,v}&&&&&  &      &       &   &       &   \\
\end{diagram}
\end{array}
\]
\caption{Data for composition in $C$.  Each entry in the lower diagram is
an element of the $b$-component of the corresponding object in the upper
diagram.}
\label{fig:sm-comp-data}
\end{figure}
\begin{itemize}
\item an element $i\in I(b)$ and a function $\gamma: X_{b,i} \go I$
\item a map $(q_{c,x})_{c,x} \goby{\theta} r$ in $P$
\item a family of maps $((p_{c,x,d,w})_{d,w} \goby{\phi_{c,x}}
q_{c,x})_{c,x}$ in $P$, where the inner indexing is over $d\in B$ and $w\in
X_{c, \gamma_c(x)}(d)$.
\end{itemize}
Going down the left-hand slope of the diagram in
Fig.~\ref{fig:sm-comp-data}, the image of this element of $(C_1\of C_1)(b)$
in $(TC_0)(b)$ is the element $m_b(i,\gamma)$ of $I(b)$ together with the
family $(p_{c(v),x(v),d,w(v)})_{d,v}$ indexed over $d\in B$ and $v\in
X_{b,m_b(i,\gamma)}(d)$.  So to define composition, we must derive from our
data a map
\begin{equation}	\label{eq:comp-form}
(p_{c(v),x(v),d,w(v)})_{d,v} \go r
\end{equation}
in $P$.  But we have the composite
\[
\theta\of (\phi_{c,x})_{c,x}: (p_{c,x,d,w})_{c,x,d,w} \go r
\]
in $P$ and the bijection
\begin{eqnarray*}
\Sigma s_{b,i,\gamma}: %&
\Sigma X_{b, m_b(i, \gamma)}	&
\goiso	&
\Sigma \left(
\coprod_{c\in B, x\in X_{b,i}(c)} X_{c, \gamma_c(x)}	
\right) \\
	%&
(d,v)	&
\goesto	&
(c(v), x(v), d, w(v)),
\end{eqnarray*}
so $(\theta\of (\phi_{c,x})_{c,x}) \cdot (\Sigma s_{b,i,\gamma})$ is a map
of the form~\bref{eq:comp-form}; and this is what we define the composite
to be.  

To define $\ids: C_0 \go C_1$, first note that if $p\in C_0(b) = P_0$ then
$\eta_{C_0,b}(p) \in (TC_0)(b)$ consists of the element $e_b\in I(b)$
together with the family $(p)_{c,u}$ indexed over $c\in B$ and $u\in
U = X_{b,e_b}(c)$, as in the diagrams
\[
\begin{diagram}[width=2em,height=1.5em,scriptlabels,tight]
	&		&C_0	&	&	\\
	&\ldTo<{\eta_{C_0}}&	&\rdTo>1&	\\
TC_0	&		&	&	&C_0	\\
\end{diagram}
\diagspace
\begin{diagram}[width=2em,height=1.5em,scriptlabels,tight]
	&		&p	&	&	\\
	&\ldGoesto	&	&\rdGoesto&	\\
(p)_{c,u}&		&	&	&p.	\\
\end{diagram}
\]
Since $U$ is a one-element set, we may define $\ids(p)\in C_1(b)$ to be the
identity map
\[
1_p^U: (p)_{c,u} \go p
\]
in $P$, and this has the correct domain and codomain.

We have now defined all the data for the $T$-multicategory $C$.  It only
remains to check the associativity and identity axioms, and these follow
from the associativity and identity axioms on the fat symmetric
multicategory $P$.  So we have defined the desired functor
\[
\fcat{FatSymMulticat} \go T\hyph\Multicat
\]
on objects.

The rest is trivial.  If $f: P \go P'$ is a map of fat symmetric
multicategories and $C$ and $C'$ are the corresponding $T$-multicategories
then we have to define a map $h: C \go C'$.  We take $h_0 = \Delta f_0$ (as
we must) and define $h$ to act on maps as $f$ does.  Since $f$ preserves
all the structure, so does $h$.  Functoriality is immediate.  
\done
\end{prooflike}%
\index{multicategory!symmetric vs. generalized@symmetric \vs.\ generalized|)}
\index{power of Set@power of $\Set$|)}%
\index{Set, power of@$\Set$, power of|)}
\index{monad!cartesian|)}

\begin{notes}

The main result of Section~\ref{sec:opds-alg-thys}---that the theories
described by plain operads are exactly the strongly regular ones---must
exist as a subconscious principle, at least, in the mind of anyone who has
worked with operads.  Nevertheless, this is as far as I know the first
proof.  A brief sketch proof was given in my~\cite[4.6]{GOM}.

The theory of familially representable monads on presheaf categories
presented here is clearly unsatisfactory.  The proof of
Theorem~\ref{thm:sm-main} is so complicated that it is at the limits of
tolerability, for this author at least.  Moreover, we have not even
attempted to describe explicitly a cartesian monad structure on a
familially representable endofunctor on a general presheaf category
$\ftrcat{\scat{B}^\op}{\Set}$ (where $\scat{B}$ need not be discrete), nor
to describe explicitly algebras for such monads.  So the theory works, but
only just.

Other ideas on the relation between cartesian and symmetric structures can
be found in Weber~\cite{Web}%
\index{Weber, Mark}
and Batanin~\cite{BatCSO}.%
\index{Batanin, Michael}

\end{notes}

\chapter{Free Multicategories}
\lbl{app:free-mti}%
\index{generalized multicategory!free|(}

\noindent
In this appendix we define what it means for a monad or a category to be
`suitable' and prove in outline the results on free multicategories
stated in~\ref{sec:free-mti}.

First we need some terminology.  Let $\Eee$ be a cartesian category,
$\scat{I}$ a small category, $D: \scat{I} \go \Eee$ a functor for which a
colimit exists, and $(D(I) \goby{p_I} L)_{I\in\scat{I}}$ a colimit cone.
Any map $L' \go L$ gives rise to a new functor $D': \scat{I} \go \Eee$ and
a new cone $(D'(I) \goby{p'_I} L')_{I\in\scat{I}}$ by pullback: if
$\Delta L$ denotes the functor $\scat{I} \go \Eee$ constant at $L$ then 
\begin{diagram}[size=2em]
D'		&\rTo	&D		\\
\dTo<{p'}	&	&\dTo>p		\\
\Delta L'	&\rTo	&\Delta L	\\
\end{diagram}
is a pullback square in the functor category $\ftrcat{\scat{I}}{\Eee}$.  We
say that the colimit $(D(I) \goby{p_I} L)_{I\in\scat{I}}$ is \demph{stable
under pullback}%
\index{stable under pullback}
if for any map $L' \go L$ in $\Eee$, the resulting cone
$(D'(I) \goby{p'_I} L')_{I\in\scat{I}}$ is also a colimit.  

The maps $p_I$ in a colimit cone $(D(I) \goby{p_I} L)_{I\in\scat{I}}$ are
called the \demph{coprojections}%
\index{coprojection}
of the colimit, so we say that the
colimit of $D$ `has monic coprojections' if each $p_I$ is monic.

A category is said to have \demph{disjoint finite coproducts}%
\index{disjoint coproduct}%
\index{coproduct!disjoint}
if it has
finite coproducts, these coproducts have monic coprojections, and for any
pair $A, B$ of objects, the square
\begin{diagram}[size=2em]
0	&\rTo		&B		\\
\dTo	&		&\dTo		\\
A	&\rTo		&A+B		\\
\end{diagram}
is a pullback.

Let $\omega$ be the natural numbers with their usual ordering. A \demph{nested
sequence}%
\index{nested sequence}
in a category \Eee\ is a functor $\omega\go\Eee$ in which the image
of every map in $\omega$ is monic; in other words, it is a diagram
\[
A_0 \monic A_1 \monic \cdots
\]
in \Eee, where $\monic$ indicates a monic.  A functor that
preserves pullbacks also preserves monics, so it makes sense for such a
functor to `preserve colimits of nested sequences'. 

Let $\scat{I}$ and $\scat{J}$ be small categories and $\Eee$ a category
with all limits of shape $\scat{I}$ and colimits of shape $\scat{J}$.  We
say that limits of shape $\scat{I}$ and colimits of shape $\scat{J}$
\demph{commute}%
\index{commuting of limit with colimit}
in $\Eee$ if for each functor $P: \scat{I} \times \scat{J}
\go \Eee$, the canonical map
\begin{equation}	\label{eq:lim-colim-map}
\lim_{\rightarrow\scat{J}} \lim_{\leftarrow\scat{I}} P
\go
\lim_{\leftarrow\scat{I}} \lim_{\rightarrow\scat{J}} P
\end{equation}
is an isomorphism.  In particular, let $\scat{I}$ be the 3-object category
such that limits over $\scat{I}$ are pullbacks, and let $\scat{J} =
\omega$; we say that pullbacks and colimits of nested sequences commute in
$\Eee$ if this canonical map is an isomorphism for all functors $P$ such
that $P(I,\dashbk): \omega \go \Eee$ is a nested sequence for each $I \in
\scat{I}$.

A category $\Eee$ is \demph{suitable}%
\index{suitable}
if
\begin{itemize}
\item $\Eee$ is cartesian
\item $\Eee$ has disjoint finite coproducts, and these are stable under
pullback
\item $\Eee$ has colimits of nested sequences; these commute with
pullbacks and have monic coprojections.
\end{itemize}
A monad $(T, \mu, \eta)$ is \demph{suitable}%
\index{suitable}
if
\begin{itemize}
\item $(T, \mu, \eta)$ is cartesian
\item $T$ preserves colimits of nested sequences.
\end{itemize}

\section{Proofs}
\lbl{sec:free-mti-proofs}

We sketch proofs of each of the results stated in~\ref{sec:free-mti}.

\begin{quotedthm}{Theorem~\ref{thm:free-gen}}
Any presheaf category is suitable.  Any finitary cartesian monad on a
cartesian category is suitable.
\end{quotedthm}
\begin{proof}
First note that the category $\omega$ is filtered.  The suitability of
$\Set$ then reduces to a collection of standard facts.  Presheaf categories
are also suitable, as limits and colimits in them are computed pointwise.  
The second sentence is trivial.  
\done
\end{proof}

Before we embark on the proofs of the main theorems, here is the main idea.
A $T$-multicategory with object-of-objects $E$ is a monoid in the monoidal
category $\Eee/(TE \times E)$ (p.~\pageref{p:slice-monoidal}), so a free
$T$-multicategory is a free%
\index{monoid!free}
monoid of sorts.  The usual formula for the
free monoid on an object $X$ of a monoidal category is $\coprod_{n\in\nat}
X^{\otimes n}$.  But this only works if the tensor product preserves
countable coproducts on each side, and this is only true in our context if
$T$ preserves countable coproducts, which is often not the case---consider
plain multicategories, for instance.  So we need a more subtle
construction.  What we actually do corresponds to taking the colimit of the
sequence
\[
\coprod_{k=0}^{0} X^{\otimes k}
\rMonic
\coprod_{k=0}^{1} X^{\otimes k}
\rMonic
\coprod_{k=0}^{2} X^{\otimes k}
\rMonic
\cdots
\]
in the case that the monoidal category \emph{does} have coproducts
preserved by the tensor; in the general case we replace $\coprod_{k=0}^{n}
X^{\otimes k}$ by $X^{(n)}$, defined recursively by
\[
X^{(0)} = I,
\diagspace 
X^{(n+1)} = I + (X \otimes X^{(n)}).
\]

\begin{quotedthm}{Theorem~\ref{thm:free-main}}
Let $T$ be a suitable monad on a suitable category $\Eee$.  Then the
forgetful functor
\[
\Cartpr\hyph\Multicat \go \Eee^+ = \Cartpr\hyph\Graph
\]
has a left adjoint, the adjunction is monadic, and if $T^+$ is the induced
monad on $\Eee^+$ then both $T^+$ and $\Eee^+$ are suitable.
\end{quotedthm}

\begin{proof}
We proceed in four steps:
\begin{enumerate}
\item 	\lbl{item:ftr}
construct a functor $F: \Eee^+ \go \Cartpr\hyph\Multicat$
\item 	%\lbl{item:adjn}
construct an adjunction between $F$ and the forgetful functor $U$
\item 	%\lbl{item:newsuit}
check that $\Eee^+$ and $T^+$ are suitable
\item 	%\lbl{item:monadic}
check that the adjunction is monadic.
\end{enumerate}
Each step goes as follows.
\begin{enumerate}
\item \emph{Construct a functor $F: \Eee^+ \go \Cartpr\hyph\Multicat$}.
Let $X$ be a $T$-graph.  For each $n\in\nat$, define a $T$-graph
\[
\begin{slopeydiag}
	&	&X_1^{(n)}&	&	\\
	&\ldTo<{d_n}&	&\rdTo>{c_n}&	\\
TX_0	&	&	&	&X_0	\\
\end{slopeydiag}
\]
by
\begin{itemize}
\item $X_1^{(0)}=X_0$, $d_0=\eta_{X_0}$, and $c_0=1$
\item $X_1^{(n+1)} = X_0 + X_1\of X_1^{(n)}$ (where $\of$ is 1-cell
composition in the bicategory $\Sp{\Eee}{T}$), with the obvious choices of
$d_{n+1}$ and $c_{n+1}$.
\end{itemize}
For each $n\in\nat$, define a map $i_n: X_1^{(n)} \go X_1^{(n+1)}$ by taking
\begin{itemize}
\item $i_0: X_0 \go X_0 + X_1 \of X_0$ to be first coprojection
\item $i_{n+1} = 1_{X_0} + (1_{X_1} * i_n)$.
\end{itemize}
Then the $i_n$'s are monic, and by taking $X_1^*$ to be the colimit of 
\[
X_1^{(0)} \rMonic^{i_0} X_1^{(1)} \rMonic^{i_1} \cdots
\]
we obtain a $T$-graph
\[
FX =
\left(
\begin{diagram}[width=1.7em,height=1.7em,scriptlabels,noPS]
	&	&X_1^*	&	&	\\
	&\ldTo<d&	&\rdTo>c&	\\
TX_0	&	&	&	&X_0	\\
\end{diagram}
\right).
\]
This $T$-graph $FX$ naturally has the structure of a $T$-multicategory.
The identities map $X_0 \go X_1^*$ is just the coprojection $X_1^{(0)}
\monic X_1^*$.  Composition comes from canonical maps $X_1^{(m)} \of
X_1^{(n)} \go X_1^{(m+n)}$ (defined by induction on $m$ for each fixed
$n$), which piece together to give a map $X_1^* \of X_1^* \go X_1^*$.  It
is the definition of composition that needs most of the suitability
axioms.

We have now described what $F$ does to objects, and extension to morphisms
is straightforward.

(The colimit of the nested sequence of $X_1^{(n)}$'s appears, in light
disguise, as the recursive description of the free plain%
\index{multicategory!free}
multicategory
monad in~\ref{sec:om-further}: $X_1^{(n)}$ is the set of formal expressions
that can be obtained from the first clause on
p.~\pageref{p:free-plain-clauses} and up to $n$ applications of the second
clause.)

\item \emph{Construct an adjunction between $F$ and $U$}.  We do this by
constructing unit and counit transformations and verifying the triangle
identities.  Both transformations are the identity on the object of objects
(`$X_0$'), so we only need define them on the object of arrows.  For the
unit, $\eta^+$, if $X\in\Eee^+$ is a $T$-graph then $\eta^+_X: X_1 \go
X_1^*$ is the composite
\[
X_1 \goiso X_1 \of X_0 \rMonic^{\mr{copr}_2} X_0 + X_1 \of X_0 = X_1^{(1)}
\rMonic X_1^*.
\]
For the counit, $\epsln^+$, let $C$ be a $T$-multicategory; write $X =
UC$ and use the notation $X_1^{(n)}$ and $X_1^*$ as in part~\bref{item:ftr}.
We need to define a map $\epsln^+_C: X_1^* \go C$.  To do this,
define for each $n\in\nat$ a map $\epsln^+_{C,n}: X_1^{(n)} \go C_1$ by 
\begin{itemize}
\item $\epsln^+_{C,0} = (C_0 \goby{\ids} C_1)$
\item $\epsln^+_{C,n+1} = 
(C_0 + C_1 \of X_1^{(n)} 
\goby{1+(1*\epsln^+_{C,n})} 
C_0 + C_1 \of C_1 
\goby{q} 
C_1)$, where $q$ is $\ids$ on the first summand and $\comp$ on the second,
\end{itemize}
and then take $\epsln^+_C$ to be the induced map on the colimit.
Verification of the triangle identities is straightforward.

\item \emph{Check that $\Eee^+$ and $T^+$ are suitable}.  The forgetful
functor
\[
\begin{array}{rcl}
\Eee^+ 	&\go		&\Eee \times \Eee	\\
X	&\goesto	&(X_0, X_1)		\\
\end{array}
\]
creates pullbacks and colimits.  This implies that $\Eee^+$ possesses
pullbacks, finite coproducts and colimits of nested sequences, and that
they behave as well as they do in $\Eee$.  So $\Eee^+$ is suitable, and it
is now straightforward to check that $T^+$ is suitable too.

\item \emph{Check that the adjunction is monadic}.  We apply the Monadicity
Theorem (Mac Lane~\cite[VI.7]{MacCWM}) by checking that $U$ creates
coequalizers for $U$-absolute-coequalizer pairs.  This is a completely
routine procedure and works for any cartesian (not necessarily suitable)
$\Eee$ and $T$.
\done
\end{enumerate}
\end{proof}

The fixed-object version of the theorem is now easy to deduce:
\begin{quotedthm}{Theorem~\ref{thm:free-fixed}}
Let $T$ be a suitable monad on a suitable category $\Eee$, and let $E \in
\Eee$.  Then the forgetful functor
\[
\Cartpr\hyph\Multicat_E \go \Eee^+_E = \Eee/(TE \times E)
\]
has a left adjoint, the adjunction is monadic, and if $T^+_E$ is the induced
monad on $\Eee^+_E$ then both $T^+_E$ and $\Eee^+_E$ are suitable.
\end{quotedthm}
\begin{proof}
In the adjunction $(F,U,\eta^+,\epsln^+)$ constructed in the proof
of~\ref{thm:free-main}, each of $F$, $U$, $\eta^+$ and $\epsln^+$ leaves
the object of objects unchanged.  The adjunction therefore restricts to an
adjunction between the subcategories $\Eee^+_E$ and
$\Cartpr\hyph\Multicat_E$, and the restricted adjunction is also monadic.

All we need to check, then, is that $\Eee^+_E$ and $T^+_E$ are suitable.
For $\Eee^+_E$, it is enough to know that the slice of a suitable category
is suitable, and to prove this we need only note that for any $E' \in
\Eee$, the forgetful functor $\Eee/E' \go \Eee$ creates both pullbacks and
colimits.  Suitability of $T^+_E$ follows from suitability of $T^+$ since
the inclusion $\Eee^+_E \rIncl \Eee^+$ preserves and reflects both
pullbacks and colimits of nested sequences.  \done
\end{proof}

Finally, we prove the result stating that if $\Eee$ and $T$ have
certain special properties beyond suitability then those properties are
inherited by $\Eee^+$ and $T^+$, or, in the fixed-object case, by
$\Eee^+_E$ and $T^+_E$.
\begin{quotedthm}{Proposition~\ref{propn:free-refined}}
If $\Eee$ is a presheaf category and the functor $T$ preserves wide
pullbacks then the same is true of $\Eee^+$ and $T^+$ in
Theorem~\ref{thm:free-main}, and of $\Eee^+_E$ and $T^+_E$ in
Theorem~\ref{thm:free-fixed}.  Moreover, if $T$ is finitary then so are
$T^+$ and $T^+_E$.
\end{quotedthm}
(Wide pullbacks were defined on p.~\pageref{p:defn-wide-pb}.)
\begin{proof}
First we show that $\Eee^+$ and $\Eee^+_E$ are presheaf categories.  Let
$G: \Eee \go \Eee$ be the functor defined by $G(E) = T(E) \times E$: then
in the terminology of p.~\pageref{p:Artin}, $\Eee^+$ is the Artin gluing
$\Eee\gluing G$.  Since $T$ preserves wide pullbacks, so too does $G$, and
Proposition~\ref{propn:Artin-gluing} then implies that $\Eee^+$ is a
presheaf category.  The fixed-object case is easier: the slice of a
presheaf category is a presheaf category~(\ref{propn:pshf-slice}).

Next we show that $T^+$ and $T^+_E$ preserve wide pullbacks.  Recall that
in the proof of~\ref{thm:free-main}, free $T$-multicategories were
constructed using coproducts, colimits of nested sequences, pullbacks, and
the functor $T$.  (The last two of these were hidden in the notation `$X
\of X_1^{(n)}$'.)  It is therefore enough to show that all four of these
entities commute with wide pullbacks.  The last two are immediate, and
Lemma~\ref{lemma:conn-lims} implies that coproducts commute with wide
pullbacks in $\Set$ (and so in any presheaf category); all that remains is
to prove that colimits of nested sequences commute with wide pullbacks in
$\Set$.  This is actually not true in general, but a slightly weaker
statement is true and suffices.  Specifically, let $\scat{I}$ be a category
of the form $\scat{P}_K$ (p.~\pageref{p:defn-wide-pb-shape}), so that a
limit over $\scat{I}$ is a wide pullback; let $\scat{J} = \omega$; and let
$P: \scat{I} \times \scat{J} \go \Set$ be a functor such that
\[
\begin{diagram}[height=2em]
P(I,J)	&\rTo	&P(I,J')	\\
\dTo	&	&\dTo		\\
P(I',J)	&\rTo	&P(I',J')	\\
\end{diagram}
\]
is a pullback square for each pair of maps $(I \go I', J \go J')$.  Then
the canonical map~\bref{eq:lim-colim-map} (p.~\pageref{eq:lim-colim-map})
is an isomorphism.  In the case at hand these squares are of the form
\[
\begin{diagram}[height=2em]
X_1^{(n)}		&\rTo^{i_n}	&X_1^{(n+1)}		\\
\dTo<{f_1^{(n)}}	&		&\dTo>{f_1^{(n+1)}}	\\
Y_1^{(n)}		&\rTo_{i_n}	&Y_1^{(n+1)}		\\
\end{diagram}
\]
where $f: X \go Y$ is some map of $T$-graphs, and it is easily checked that
such squares are pullbacks.

For `moreover' we have to show that $T^+$ and $T^+_E$ preserve filtered
colimits if $T$ does.  Just as in the previous paragraph, this reduces to
the statement that filtered colimits commute with pullbacks in $\Set$,
whose truth is well-known (Mac Lane~\cite[IX.2]{MacCWM}).  \done
\end{proof}

\begin{notes}

These proofs first appeared in my~\cite{GECM}.  They have much in common
with the free monoid construction of Baues,%
\index{Baues, Hans-Joachim}
Jibladze%
\index{Jibladze, Mamuka}
and Tonks~\cite{BJT}.%
\index{Tonks, Andy}
\index{generalized multicategory!free|)}

\end{notes}

\chapter{Definitions of Tree}
\lbl{app:trees}%
\index{tree!graph@as graph|(}%
\index{graph!tree@of tree|(}

\chapterquote{%
I met this guy\\
\makebox[0em][l]{and he looked like he might have been a hat-check clerk}\\
at an ice rink\\
which in fact\\
he turned out to be}{%
Laurie Anderson~\cite{Laurie}}

\noindent
Trees appear everywhere in higher-dimensional algebra.  In this text they
were defined in a purely abstract way~(\ref{eg:opd-of-trees}): $\tr$ is the
free plain operad on the terminal object of $\Set^\nat$, and an $n$-leafed
tree is an element of $\tr(n)$.  But for the reasons laid out at the
beginning of~\ref{sec:trees}, I give here a `concrete', graph-theoretic,
definition of (finite, rooted, planar) tree and sketch a proof that it is
equivalent to the abstract definition.

\section{The equivalence}

The main subtlety is that the trees we use are not quite finite graphs in
the usual sense: some of the edges have a vertex at only one of their ends.
(Recall from~\ref{eg:opd-of-trees} that in a tree, an edge with a free end
is not the same thing as an edge ending in a vertex.)  This suggests the
following definitions.

\begin{figure}
\centering
\setlength{\unitlength}{1mm}
\begin{picture}(105,50)(0,-5)
\cell{0}{0}{bl}{%
\begin{picture}(50,45)
\cell{0}{0}{bl}{\epsfig{file=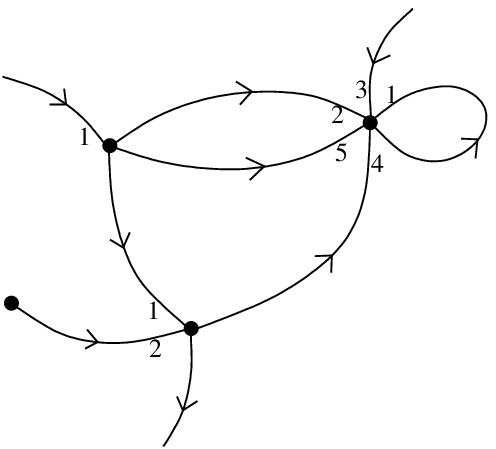}}
\cell{20}{3}{c}{o}
\cell{5}{38}{c}{i_1}
\cell{38}{45}{t}{i_2}
\end{picture}}
\cell{25}{-5}{b}{\textrm{(a)}}
\cell{62}{0}{bl}{%
\begin{picture}(43,45)(-1,0)
\cell{0}{0}{bl}{\epsfig{file=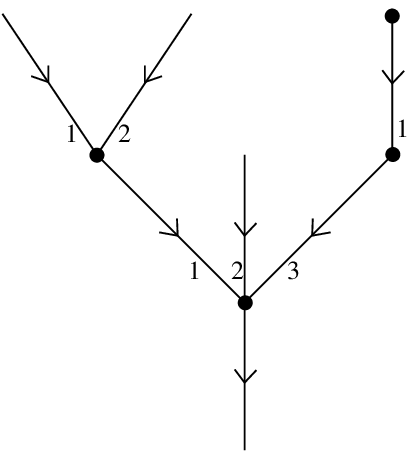}}
\cell{26}{4}{l}{o}
\cell{1}{41}{c}{i_1}
\cell{20}{41}{c}{i_2}
\cell{26}{27}{l}{i_3}
\end{picture}}
\cell{83.5}{-5}{b}{\textrm{(b)}}
\end{picture}
% \hand{50}{26}
\caption{(a) Input-output graph with 4 vertices and 2 input edges $i_1,
i_2$, (b) combinatorial tree with 4 vertices and 3 input edges $i_1, i_2,
i_3$.  In both, the numbers indicate the order on the edges arriving at
each vertex.}
\label{fig:comb-graphs}
\end{figure}

\begin{defn}
\item A (planar) \demph{input-output graph}%
\index{graph!input-output}
(Fig.~\ref{fig:comb-graphs}(a))
consists of
\begin{itemize}
\item a finite set $V$ (the \demph{vertices})
\item a finite set $E$ (the \demph{edges}), a subset $I \sub E$ (the
\demph{input edges}), and an element $o \in E$ (the \demph{output edge})
\item a function $s: E\without I \go V$ (\demph{source}) and a function $t:
E\without \{o\} \go V$ (\demph{target})
\item for each $v \in V$, a total order $\leq$ on $t^{-1}\{v\}$.
\end{itemize}
\end{defn}
We write $v \goby{e}$ to mean that $e$ is a non-input edge with $s(e) = v$,
and similarly $\goby{e} v'$ to mean that $e$ is a non-output edge with
$t(e) = v'$, and of course \mbox{$v \goby{e} v'$} to mean that $e$ is a
non-input, non-output edge with $s(e) = v$ and $t(e) = v'$.

A tree is roughly speaking a connected, simply connected graph, and the
following notion of path allows us to express this.  

\begin{defn}
A \demph{path}%
\index{path!graph@in graph}
from a vertex $v$ to an edge $e$ in an input-output graph is
a diagram
\[
v= v_1 \goby{e_1} v_2 \goby{e_2} 
\ \cdots \ 
\goby{e_{l-1}} v_l \goby{e_l = e}
\]
in the graph.  That is, a path from $v$ to $e$ consists of
\begin{itemize}
\item an integer $l\geq 1$
\item a sequence $(v_1, v_2, \ldots, v_l)$ of vertices with $v_1=v$
\item a sequence $(e_1, \ldots, e_{l-1}, e_l)$ of edges with $e_l=e$
\end{itemize}
such that
\[
v_1 = s(e_1), \ 
t(e_1) = v_2 = s(e_2), \ 
\ldots, \ 
t(e_{l-1}) = v_l = s(e_l)
\]
and all of these sources and targets are defined.
\end{defn}

\begin{defn}
A \demph{combinatorial tree}%
\index{tree!combinatorial}
is an input-output graph such that for every
vertex $v$, there is precisely one path from $v$ to the output edge.
\end{defn}
Fig.~\ref{fig:comb-graphs}(b) shows a combinatorial tree.  The ordering of
the edges arriving at each vertex encodes the planar embedding.  `Tree' is
an abbreviation for `finite, rooted,%
\index{root of tree}
planar tree'.  If we were doing
symmetric operads then we would use non-planar trees, if we were doing
cyclic operads then we would use non-rooted trees, and so on.

A combinatorial tree is essentially the same thing as a tree in our sense,
where `essentially the same thing' refers to the obvious notion of
isomorphism between combinatorial trees.  We write $\fcat{combtr}(n)$ for
the set of isomorphism classes of combinatorial trees with $n$ input edges.

\begin{propn}
For each $n\in\nat$, there is a canonical bijection $\tr(n) \iso
\fcat{combtr}(n)$. 
\end{propn}

With a little more work we could define an operad structure on
$(\fcat{combtr}(n))_{n\in\nat}$ and turn the proposition into the stronger
statement that the operads $\tr$ and $\fcat{combtr}$ are isomorphic.  With
more work still we could define maps between combinatorial trees and
so define a \Cat-operad $\fcat{COMBTR}$ in which $\fcat{COMBTR}(n)$ is the
category of $n$-leafed combinatorial trees; then we could prove that this
\Cat-operad is equivalent to the \Cat-operad $\TR$.

\begin{prooflike}{Sketch proof}
The strategy is to define functions 
\[
\tr(n) \oppair{\Phi}{\Psi} 
\{ \textrm{combinatorial trees with } n \textrm{ input edges} \}
\]
for each $n\in\nat$, such that if $G \iso G'$ then $\Psi(G) =
\Psi(G')$, and such that $\Psi(\Phi(\tau)) = \tau$ and $\Phi(\Psi(G))
\iso G$ for each tree $\tau$ and combinatorial tree $G$.  The proposition
follows.  The definition of $\Phi(\tau)$ is by induction on the structure
of the tree $\tau$, and the definition of $\Psi(G)$ is by induction on the
number of vertices of the combinatorial tree $G$.  All the details are
straightforward. 
\done
\end{prooflike}

\begin{notes}

Some pointers to the literature on trees can be found in the Notes to
Chapter~\ref{ch:opetopic}.%
\index{tree!graph@as graph|)}%
\index{graph!tree@of tree|)}

\end{notes}

\chapter{Free Strict $n$-Categories}
\lbl{app:free-strict}%
\index{n-category@$n$-category!strict!free|(}%
\index{omega-category@$\omega$-category!strict!free|(}

\noindent
Here we prove that the forgetful functor
\[
(\textrm{strict } n \textrm{-categories} )
\go
(n\textrm{-globular sets})
\]
is monadic and that the induced monad is cartesian and finitary.  We also
prove the analogous results for $\omega$-categories.  We used monadicity
and that the monad is cartesian in Part~\ref{part:n-categories}, in order
to be able to define and understand globular operads and weak $n$- and
$\omega$-categories.  We use the fact that the monad is finitary for
technical purposes in Appendix~\ref{app:initial}.

It is frustrating to have to prove this theorem, for two separate reasons.
The first is that the proof can \emph{almost} be made trivial: an adjoint
functor theorem tells us that the forgetful functor has a left adjoint, a
monadicity theorem tells us that it is monadic, and a routine calculation
tells us that it is finitary.  However, we have no result of the form
`given an adjunction, its induced monad is cartesian if the right adjoint
satisfies certain conditions', and so in order to prove that the monad is
cartesian we are forced to actually construct the whole adjunction
explicitly.  The second is that there ought to be some way of simply
looking at the theory of strict $\omega$-categories, presented as an
algebraic theory on the presheaf category of globular sets (as
in~\ref{defn:strict-n-cat-glob}), and applying some general principle to
deduce the theorem immediately; see p.~\pageref{p:sr-presheaf}.  But
again, we currently have no way of doing this.

Despite these frustrations, the proof is quite easy.  By exploiting the
definition of strict $n$- and $\omega$-categories by iterated
enrichment~(\ref{defn:strict-n-cat-enr}), we can reduce it to the proof of
some straightforward statements about enriched categories.

Here is the idea.  Suppose we know how to construct the free strict
2-category on a 2-globular set.  Then we can construct the free strict
3-category on a 3-globular set $X$ in two steps:
\begin{itemize}
  \item For each $x, x' \in X(0)$, take the 2-globular set $X(x,x')$ of
  cells whose 0-dimensional source is $x$ and whose 0-dimensional target is
  $x'$, replace it by the free strict 2-category on $X(x,x')$, and
  reassemble to obtain a new 3-globular set $Y$.  The $0$- and $1$-cells of
  $Y$ are the same as those of $X$, and a typical 3-cell of $Y$ looks
  like
  \[
\setlength{\unitlength}{1mm}
\begin{picture}(70,56)(-35,-28)
\cell{0}{0.5}{c}{\epsfig{file=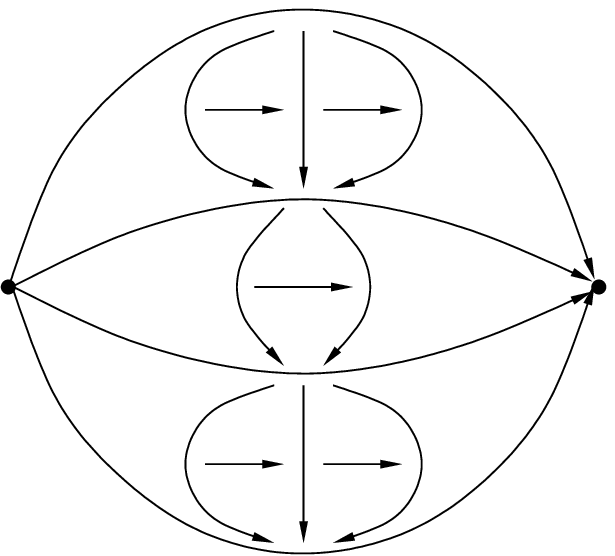}}
\cell{-32}{0}{r}{x}
\cell{32}{0}{l}{x'}
\cell{15}{26}{c}{\scriptstyle f_0}
\cell{15}{9}{c}{\scriptstyle f_1}
\cell{15}{-9}{c}{\scriptstyle f_2}
\cell{15}{-25}{c}{\scriptstyle f_3}
\cell{-14}{17}{c}{\scriptstyle \alpha_0}
\cell{2}{15}{c}{\scriptstyle \alpha_1}
\cell{14}{17}{c}{\scriptstyle \alpha_2}
\cell{-9}{0}{c}{\scriptstyle \beta_0}
\cell{9}{0}{c}{\scriptstyle \beta_1}
\cell{-14}{-18}{c}{\scriptstyle \gamma_0}
\cell{2}{-20}{c}{\scriptstyle \gamma_1}
\cell{14}{-18}{c}{\scriptstyle \gamma_2}
\cell{-6}{19}{b}{\scriptstyle \Gamma_1}
\cell{6}{19}{b}{\scriptstyle \Gamma_2}
\cell{0}{1}{b}{\scriptstyle \Delta_1}
\cell{-6}{-17}{b}{\scriptstyle \Theta_1}
\cell{6}{-17}{b}{\scriptstyle \Theta_2}
\end{picture}
%   \hand{70}{48}
  \]
  where $x$ and $x'$ are 0-cells of $X$, the $f_i$'s are 1-cells, the
  $\alpha_i$'s, $\beta_i$'s and $\gamma_i$'s are 2-cells, and the
  $\Gamma_i$'s, $\Delta_i$'s and $\Theta_i$'s are 3-cells.
  
  \item Write $\cat{A}$ for the category of 2-globular sets, think of $Y$
  as a family $(Y(x,x'))_{x, x'\in X(0)}$ of objects of $\cat{A}$, and let
  $Z$ be the free $\cat{A}$-enriched category on $Y$.  The $0$-cells of $Z$
  are the same as those of $X$ and $Y$, and a typical 3-cell of $Z$ looks
  like
  \[
  \epsfig{file=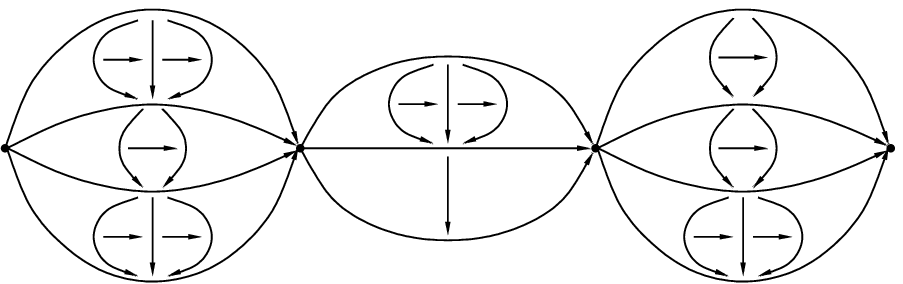}
%   \hand{45}{49}
  \]
  where the cells making up this diagram are cells of $X$.  So $Z$ is the
  free strict 3-category on $X$.
\end{itemize}
We therefore construct the free strict $n$-category functor inductively on
$n$, and pass to the limit to reach the infinite-dimensional version. 

In the first section,~\ref{sec:free-enr}, we prove the necessary results on
free enriched categories.  In the second,~\ref{sec:free-n}, we apply them
to establish the existence and properties of the free strict $n$- and
$\omega$-category functors.  Since the construction is explicit, we are
then able to verify the intuitively plausible formula~\bref{eq:pd-rep-of-T}
(p.~\pageref{eq:pd-rep-of-T}) for the free functor in terms of pasting
diagrams.

\section{Free enriched categories}
\lbl{sec:free-enr}%
\index{category!free enriched|(}%
\index{enrichment!category@of category!finite product category@in finite product category|(}

Recall from~\ref{sec:cl-enr} that for every category $\cat{A}$ there is a
category $\cat{A}\hyph\Gph$ of $\cat{A}$-graphs, and that if $\cat{A}$ has
finite products then there is also a category $\cat{A}\hyph\Cat$ of
$\cat{A}$-enriched categories and a forgetful functor $\cat{A}\hyph\Cat \go
\cat{A}\hyph\Gph$.  Any finite-product-preserving functor $Q: \cat{B} \go
\cat{A}$ between categories with finite products induces a (strictly)
commutative square
\[
\begin{diagram}[size=2em]
\cat{B}\hyph\Cat	&\rTo	&\cat{A}\hyph\Cat	\\
\dTo			&	&\dTo			\\
\cat{B}\hyph\Gph	&\rTo	&\cat{A}\hyph\Gph,	\\
\end{diagram}
\]
and in particular induces an unambiguous functor $\cat{B}\hyph\Cat \go
\cat{A}\hyph\Gph$.  So if $T$ is a monad on a category $\cat{A}$ with
finite products then the forgetful functor $\cat{A}^T \go \cat{A}$ induces
a forgetful functor $\cat{A}^T\hyph\Cat \go \cat{A}\hyph\Gph$.  This
functor will play an important role.

A monad will be called \demph{coproduct-preserving}%
\index{coproduct!preserved by monad}%
\index{monad!coproduct-preserving}
if its functor part
preserves all (small) coproducts.

\begin{propn}	\lbl{propn:free-enr}
Let $\cat{A}$ be a presheaf category.
\begin{enumerate}
\item	\lbl{item:free-enr-simple}
The forgetful functor 
\[
\cat{A}\hyph\Cat \go \cat{A}\hyph\Gph
\]
is monadic, and the induced monad is cartesian, finitary and
coproduct-preserving.
\item	\lbl{item:free-enr-graphs}
For any monad $T$ on $\cat{A}$, the forgetful functor 
\[
\cat{A}^T\hyph\Gph \go \cat{A}\hyph\Gph
\]
is monadic, and the induced monad is cartesian (respectively, finitary or
coproduct-preserving) if $T$ is.
\item	\lbl{item:free-enr-both}
For any coproduct-preserving monad $T$ on $\cat{A}$, the forgetful functor 
\[
\cat{A}^T\hyph\Cat \go \cat{A}\hyph\Gph
\]
is monadic; the induced monad $\newmnd{T}$ is coproduct-preserving, and is
cartesian (respectively, finitary) if $T$ is.
\end{enumerate}
Moreover, $\newmnd{T}$ 
is given on $\cat{A}$-graphs $X$ by $(\newmnd{T}X)_0 = X_0$ and
\[
(\newmnd{T}X)(x,x') 
=
\coprod_{x=x_0, x_1, \ldots, x_{r-1}, x_r=x'}
T(X(x_0, x_1)) \times\cdots\times T(X(x_{r-1},x_r))
\]
($x,x' \in X_0$), where the coproduct is over all $r\in\nat$ and sequences
$x_0, \ldots, x_r$ of elements of $X_0$ satisfying $x_0 = x$ and $x_r =
x'$.
\end{propn}
The hypothesis that $\cat{A}$ is a presheaf category is excessive, but
makes the proof easier and serves our purpose.
\begin{proof}
Part~\bref{item:free-enr-simple} is simple: free $\cat{A}$-enriched
categories are constructed just as free ordinary categories are.
The induced monad $\fc_\cat{A}$ on $\cat{A}\hyph\Gph$ is given on
$\cat{A}$-graphs $X$ by $(\fc_\cat{A} X)_0 = X_0$ and
\[
(\fc_\cat{A}X)(x,x')
=
\coprod_{x=x_0, x_1, \ldots, x_{r-1}, x_r=x'}
X(x_0, x_1) \times\cdots\times X(x_{r-1},x_r)
\]
($x,x' \in X_0$).  Everything works as in the familiar case $\cat{A}=\Set$
because $\cat{A}$ is a $\Set$-valued functor category.

Part~\bref{item:free-enr-graphs} is also straightforward: $\blank\hyph\Gph$
defines a strict map $\CAT \go \CAT$ of strict 2-categories (ignoring
set-theoretic worries), and so turns the adjunction between $\cat{A}^T$ and
$\cat{A}$ into an adjunction between the corresponding categories of
graphs.  Explicitly, the induced monad $T_*$ on $\cat{A}\hyph\Gph$ is given
on $\cat{A}$-graphs $X$ by $(T_* X)_0 = X_0$ and $(T_* X)(x,x') =
T(X(x,x'))$.  That $T_*$ inherits the properties of $T$ is easily checked.

To prove~\bref{item:free-enr-both} and `moreover' it suffices to put a
monad structure on the composite functor $\fc_\cat{A} \of T_*$, to
construct an isomorphism of categories
\[
(\cat{A}\hyph\Gph)^{\fc_\cat{A} \of T_*} \iso \cat{A}^T\hyph\Cat
\]
commuting with the forgetful functors to $\cat{A}\hyph\Gph$, and to check
that the monad $\fc_\cat{A} \of T_*$ is cartesian if $T$ is.  (Recall the
two-step strategy of the introduction.)  The monad structure on
$\fc_\cat{A} \of T_*$ consists of the monad structures on $T_*$ and
$\fc_\cat{A}$ glued together by a distributive%
\index{distributive law}
law~(\ref{defn:distrib-law})
\[
\lambda: T_* \of \fc_\cat{A} \go \fc_\cat{A} \of T_*.
\]
The functor $\fc_\cat{A} \of T_*$ is given by the formulas for $\newmnd{T}$
in `moreover'; the composite the other way round is given by $((T_* \of
\fc_\cat{A})(X))_0 = X_0$ and
\[
((T_* \of \fc_\cat{A})X)(x,x') 
=
\coprod_{x=x_0, x_1, \ldots, x_{r-1}, x_r=x'}
T(X(x_0, x_1) \times\cdots\times X(x_{r-1},x_r)).
\]
So we can define $\lambda_X$ by the evident natural maps.  It is
straightforward to check that $\lambda$ is indeed a distributive law and
that $\lambda$ is cartesian if $T$ is, so the monad $\fc_\cat{A} \of T_*$
is cartesian if $T$ is~(\ref{lemma:distrib-gives-monad}).  All that remains
is the isomorphism.  That follows from~\ref{lemma:distrib-iso-algs} once we
know that the monad on $(\cat{A}\hyph\Gph)^{T_*} \iso \cat{A}^T\hyph\Gph$
corresponding to the distributive law $\lambda$ is $\fc_{\cat{A}^T}$, and
that too is easily checked.  \done
\end{proof}%
\index{category!free enriched|)}%
\index{enrichment!category@of category!finite product category@in finite product category|)}

\section{Free $n$- and $\omega$-categories}
\lbl{sec:free-n}%
\index{enrichment!define n-category@to define $n$-category|(}

Since we are using the definition of strict $n$- and $\omega$-categories by
iterated enrichment, it is convenient to replace the category of
$n$-globular sets by the equivalent category $n\hyph\Gph$ of $n$-graphs%
\index{n-graph@$n$-graph|(}
introduced in~\ref{propn:str-n-cats-comparison}.  There is a forgetful
functor $U_n: \strc{n} \go n\hyph\Gph$ defined by taking $U_0$ to be the
identity and $U_{n+1}$ to be the diagonal of the commutative square
\begin{equation}	\label{diag:cat-gph-square}
\begin{diagram}[size=2em]
(\strc{n})\hyph\Cat	&\rTo	&(n\hyph\Gph)\hyph\Cat	\\
\dTo			&	&\dTo			\\
(\strc{n})\hyph\Gph	&\rTo	&(n\hyph\Gph)\hyph\Gph.	\\
\end{diagram}
\end{equation}
(We abbreviate our previous notation for the category of strict
$n$-categories, \strcat{n}, to \strc{n}.)  There is also a restriction
functor $R_n: (n+1)\hyph\Gph \go n\hyph\Gph$%
\glo{Rngraphs}
for each $n\in\nat$, defined
by $R_0(X) = X_0$ and $R_{n+1} = R_n\hyph\Gph$.  These functors fit
together into a strictly commutative diagram
\[
\begin{diagram}[height=2em,width=3em] %[size=3em,tight]
\cdots		&\strc{(n+1)}	&\rTo^{S_n}	
&\strc{n}	&\rTo^{S_{n-1}}	&\ &\cdots&\ &\rTo^{S_0}&\strc{0}	\\
		&\dTo>{U_{n+1}}	&		
&\dTo>{U_n}	&		&&	&&		&\dTo>{U_0}	\\
\cdots		&(n+1)\hyph\Gph	&\rTo_{R_n}	
&n\hyph\Gph	&\rTo_{R_{n-1}}	&\ &\cdots&\ &\rTo_{R_0}&0\hyph\Gph,	\\
\end{diagram}
\]
where the functors $S_n$ are the usual ones
(p.~\pageref{p:forgetful-strict-n}).  Passing to the limit gives a category
$\omega\hyph\Gph$,%
\glo{omegaGph}%
\index{omega-graph@$\omega$-graph}
the \demph{$\omega$-graphs}, and a forgetful functor $U:
\strc{\omega} \go \omega\hyph\Gph$.

Analogous functors can, of course, be defined with
$\ftrcat{\scat{G}_n^\op}{\Set}$ in place of $n\hyph\Gph$ and
$\ftrcat{\scat{G}^\op}{\Set}$ in place of $\omega\hyph\Gph$.  It is
straightforward to check that there are equivalences of categories
\[
n\hyph\Gph \eqv \ftrcat{\scat{G}_n^\op}{\Set},
\diagspace
\omega\hyph\Gph \eqv \ftrcat{\scat{G}^\op}{\Set}
\]
commuting with all these functors.  We are therefore at liberty to use
graphs in place of globular sets.

\begin{thm}	\lbl{thm:n-forgetful-properties}
For each $n\in\nat$, the forgetful functor $\strc{n} \go
\ftrcat{\scat{G}_n^\op}{\Set}$ is monadic and the induced monad is
cartesian, finitary, and coproduct-preserving.
\end{thm}
\begin{proof}
We replace this forgetful functor by $U_n$ and use induction.  $U_0$ is
the identity.  Given $n\in\nat$, write $T_n$ for the monad induced by $U_n$
and its left adjoint (which exists by inductive hypothesis).  Then
$U_{n+1}$ is the diagonal of the square~\bref{diag:cat-gph-square}, and
under the equivalence $\strc{n} \eqv (n\hyph\Gph)^{T_n}$ becomes the
forgetful functor
\[
(n\hyph\Gph)^{T_n}\hyph\Cat \go (n\hyph\Gph)\hyph\Gph.
\]
The result now follows from
Proposition~\ref{propn:free-enr}\bref{item:free-enr-both}.
\done
\end{proof}

We want to deduce the same result for $\omega$-dimensional structures, and
morally this should be immediate from their definition by limits.  The only
problem is that $\strc{\omega}$ and $\omega\hyph\Gph$ are defined as
\emph{strict}, or 1-categorical, limits in the 2-category $\CAT$, and
properties such as adjointness are 2-categorical.  So, for instance, the
fact that each $U_n$ has a left adjoint $F_n$ does not \latin{a priori}
guarantee that $U$ has a left adjoint, since the squares
\begin{equation}	\label{diag:F-square}
\begin{diagram}[size=2em]
\strc{(n+1)}	&\rTo^{S_n}	&\strc{n}	\\
\uTo<{F_{n+1}}	&		&\uTo>{F_n}	\\
(n+1)\hyph\Gph	&\rTo_{R_n}	&n\hyph\Gph	\\
\end{diagram}
\end{equation}
are only known to commute up to (canonical) isomorphism.  

A satisfactory resolution would involve the theory of weak%
\index{limit!weak}
limits in a
2-category.  Here, however, we use a short and nasty method, exploiting
some special features of the situation.  

The key is that each of the functors $S_n$ has the following (easily
proved) isomorphism-lifting property: if $C\in \strc{(n+1)}$ and $j: S_n(C)
\goiso D$ is an isomorphism in $\strc{n}$, then there exists an isomorphism
$i: C \goiso C'$ in $\strc{(n+1)}$ such that $S_n C' = D$ and $S_n i = j$.
This allows us to choose left adjoints $F_0$, $F_1$, \ldots\ successively
so that the squares~\bref{diag:F-square} are strictly commutative.  

Observe also that the categories $\strc{n}$ have all (small) limits and
colimits and the functors $S_n$ preserve them, as follows by induction
using standard facts about enriched categories.  Together with the
isomorphism-lifting property, this implies that $\strc{\omega}$ has all
limits and colimits and that a (co)cone in the category $\strc{\omega}$ is
a (co)limit if and only if its image in each of the categories $\strc{n}$
is a (co)limit.  The same is true with $\Gph$ in place of $\Cat$, easily.

\begin{thm}	\lbl{thm:omega-forgetful-properties}
The forgetful functor $\strc{\omega} \go \ftrcat{\scat{G}^\op}{\Set}$ is
monadic and the induced monad is cartesian, finitary, and
coproduct-preserving.
\end{thm}
\begin{proof}
It is equivalent to prove the same properties of $U: \strc{\omega} \go
\omega\hyph\Gph$.  Choose left adjoints $F_n$ to the $U_n$'s so that the
squares~\bref{diag:F-square} commute, and let $F$ be the induced functor
$\omega\hyph\Gph \go \strc{\omega}$; this is left adjoint to $U$.  With the
aid of the Monadicity Theorem we see that all the properties of $U$
remaining to be proved concern limits and colimits in $\strc{\omega}$ and
$\omega\hyph\Gph$, and by the observations above, they are implied by the
corresponding properties of $U_n$
(Theorem~\ref{thm:n-forgetful-properties}).  \done
\end{proof}%
\index{enrichment!define n-category@to define $n$-category|)}

We can now read off an explicit formula for the free strict
$\omega$-category monad $T$.  Define the globular set $\pd$%
\index{pasting diagram!globular}
by taking
$\pd(0)$ to be a one-element set and $\pd(n+1)$ to be the free monoid on
$\pd(n)$, and define for each pasting diagram $\pi$ a globular set
$\rep{\pi}$, as in~\ref{sec:free-strict}.  Write $\gm{n}$ for the free
strict $n$-category monad on the category $\ftrcat{\scat{G}_n^\op}{\Set}$
of $n$-globular sets, and $X_{\trunc{n}}$%
\glo{lowertrunc}%
\index{restriction}
for the $n$-globular set obtained
by forgetting all the cells of a globular set $X$ above dimension $n$.

\begin{propn}	\label{propn:pds-formula}%
\index{familial representability!free strict omega-category functor@of free strict $\omega$-category functor}
For globular sets $X$ and $n\in\nat$, there is an isomorphism
\[
(TX)(n) 
\iso 
\coprod_{\pi\in\pd(n)} 
\ftrcat{\scat{G}^\op}{\Set} (\rep{\pi}, X)
\]
natural in $X$.
\end{propn}
\begin{proof}
If $\pi\in\pd(n)$ then the globular set $\rep{\pi}$ is empty above
dimension $n$.  Also
\[
(TX)(n) = (TX)_{\trunc{n}}(n) \iso (\gm{n} X_{\trunc{n}})(n)
\]
by construction of $\gm{n}$ and $T$.  So the claimed isomorphism is equivalent
to
\[
(\gm{n} X_{\trunc{n}})(n) 
\iso 
\coprod_{\pi\in\pd(n)} 
\ftrcat{\scat{G}_n^\op}{\Set} (\rep{\pi}_{\trunc{n}}, X_{\trunc{n}}).
\]
This is a statement about $n$-globular sets; let us translate it into one
about $n$-graphs.

First, if $n\in\nat$ and $Z$ is an $n$-graph then the corresponding
$n$-globular set has a set of $n$-cells, which we write as $Z(n)$.  If
$n=0$ then this is given by $Z(0) = Z$, and then inductively,
\[
Z(n+1) \iso \coprod_{z,z'\in Z_0} (Z(z,z'))(n).
\]
Second, if $n\in\nat$ and $\pi\in\pd(n)$ then there is an $n$-graph
$\twid{\pi}$ corresponding to the $n$-globular set $\rep{\pi}_{\trunc{n}}$.
If $n=0$ and $\pi$ is the unique element of $\pd(0)$ then $\twid{\pi} = 1$.
If $\pi\in\pd(n+1)$ then $\pi = (\pi_1, \ldots, \pi_r)$ for some $r\in\nat$
and $\pi_i \in \pd(n)$, and in this case $\twid{\pi}$ is given by
\[
(\twid{\pi})_0 = \{ 0, \ldots, r \}
\]
and
\[
\twid{\pi}(i,j)
=
\left\{
\begin{array}{ll}
\twid{\pi_j}	&\textrm{if } j = i+1	\\
0		&\textrm{otherwise.}
\end{array}
\right.
\]
The claimed isomorphism is therefore equivalent to 
\[
(\gm{n} Y)(n) 
\iso 
\coprod_{\pi\in\pd(n)}
n\hyph\Gph (\twid{\pi}, Y)
\]
for $n\in\nat$ and $n$-graphs $Y$.  

For $n=0$ this is trivial.  Then inductively, using `moreover' of
Proposition~\ref{propn:free-enr} for the second step,
\begin{eqnarray*}
\lefteqn{(\gm{n+1} Y)(n+1)
=
\coprod_{y,y' \in Y_0} ((\gm{n+1} Y)(y,y'))(n)}				\\
% 	&	&(\gm{n+1} Y)(n+1)					\\
% 	&=	&\coprod_{y,y' \in Y_0} 
% 			((\gm{n+1} Y)(y,y'))(n)				\\
	&\iso	&\coprod_{r\in\nat, y_0, \ldots, y_r \in Y_0}
\left(\bkthack
	\gm{n}(Y(y_0,y_1)) \times\cdots\times \gm{n}(Y(y_{r-1},y_r))
\right) (n)								\\
	&\iso	&\coprod_{r\in\nat, y_0, \ldots, y_r \in Y_0}
(\gm{n}(Y(y_0,y_1)))(n) \times\cdots\times (\gm{n}(Y(y_{r-1},y_r)))(n)	\\
	&\iso	&\coprod_{\begin{array}{c}\scriptstyle
			r\in\nat, y_0, \ldots, y_r \in Y_0,\\ \scriptstyle
			\pi_1, \ldots, \pi_r \in \pd(n)
			\end{array}}
		n\hyph\Gph(\twid{\pi_1}, Y(y_0,y_1))
		\times\cdots\times
		n\hyph\Gph(\twid{\pi_r}, Y(y_{r-1},y_r))		\\
	&\iso	&\coprod_{r\in\nat, \pi_1, \ldots, \pi_r \in \pd(n)}
		(n+1)\hyph\Gph(\twid{(\pi_1, \ldots, \pi_r)}, Y)	\\
	&\iso	&\coprod_{\pi\in\pd(n+1)}
		(n+1)\hyph\Gph(\twid{\pi},Y),
\end{eqnarray*}
as required.
\done
\end{proof}%
\index{n-graph@$n$-graph|)}

\begin{notes}

The material here first appeared in my thesis~\cite[App.~C]{OHDCT}.%
\index{n-category@$n$-category!strict!free|)}%
\index{omega-category@$\omega$-category!strict!free|)}

\end{notes}

\chapter{Initial Operad-with-Contraction}
\lbl{app:initial}%
\index{globular operad!contraction@with contraction|(}

\chapterquote{%
There existed another ending to the story of O.}{%
R\'eage~\cite{Rea}}

\noindent
We prove Proposition~\ref{propn:OC-initial}: the category $\fcat{OC}$ of
operads-with-contraction has an initial object.  This was needed in
Chapter~\ref{ch:a-defn} for the definition of weak $\omega$-category.

The explanation in~\ref{sec:contr} suggests an explicit construction of the
initial operad-with-contraction: ascend through the dimensions, at each
stage freely adding in elements obtained by contraction and then freely
adding in elements obtained by operadic composition.  Here we take a
different approach, exploiting a known existence theorem.

\section{The proof}
\lbl{sec:initial-proof}

The following result appears to be due to Kelly~\cite[27.1]{KelUTT}.%
\index{Kelly, Max}
\begin{thm}	\lbl{thm:comb-mon}
Let
\begin{diagram}[size=2em]
\cat{D}	&\rTo	&\cat{C}	\\
\dTo	&	&\dTo>V		\\
\cat{B}&\rTo_U	&\cat{A}	\\
\end{diagram}
be a (strict) pullback diagram in \fcat{CAT}. If $\cat{A}$ is locally
finitely presentable and each of $U$ and $V$ is finitary and monadic then
the functor $\cat{D} \go \cat{A}$ is also monadic.
\done
\end{thm}
Actually, all we need is:
\begin{cor}
In the situation of Theorem~\ref{thm:comb-mon}, \cat{D} has an initial object.
\end{cor}
\begin{proof}
A locally finitely presentable category is by definition cocomplete, so
\cat{A} has an initial object. The functor $\cat{D} \go \cat{A}$ has a left
adjoint (being monadic), which applied to the initial object of \cat{A} gives
an initial object of \cat{D}. 
\done
\end{proof}

Let $T$ be the free strict $\omega$-category monad on the category
$\ftrcat{\scat{G}^\op}{\Set}$ of globular sets, as in
Chapters~\ref{ch:globular} and~\ref{ch:a-defn}.  Write $\fcat{Coll}$ for
the category $\ftrcat{\scat{G}^\op}{\Set}/\pd$ of collections
(p.~\pageref{p:defn-collection}), $\fcat{CC}$ for the category of
collections-with-contraction%
\index{collection!contraction@with contraction}
(defined by replacing `operad' by `collection'
throughout Definition~\ref{defn:OC}), and $\fcat{Operad}$ for the category
of globular operads.  Then there is a (strict) pullback diagram
\begin{diagram}[size=2em]
\fcat{OC}	&\rTo	&\fcat{Operad}	\\
\dTo		&	&\dTo>V		\\
\fcat{CC}	&\rTo_U	&\fcat{Coll}	\\
\end{diagram}
in \fcat{CAT}, made up of forgetful functors.

To prove that $\fcat{OC}$ has an initial object, we verify the hypotheses
of Theorem~\ref{thm:comb-mon}.  The only non-routine part is showing that
we can freely add a contraction%
\index{contraction!free}
to any collection.

\minihead{\fcat{Coll} is locally finitely presentable} 

Since $\fcat{Coll}$ is a slice of a presheaf category, it is itself a
presheaf category~(\ref{propn:pshf-slice}) and so locally finitely
presentable (Borceux~\cite[Example 5.2.2(b)]{Borx2}).

\minihead{$U$ is finitary and monadic} 

It is straightforward to calculate that $U$ creates filtered colimits; and
since \fcat{Coll} possesses all filtered colimits, $U$ preserves them too.
It is also easy to calculate that $U$ creates coequalizers for $U$-split
coequalizer pairs.  So we have only to show that $U$ has a left adjoint.

Let $P$ be a collection.  We construct a new collection $FP$, a contraction
$\kappa^P$ on $FP$, and a map $\alpha_P: P \go FP$, together having the
appropriate universal property; so the functor $P \goesto (FP, \kappa^P)$
is left adjoint to $U$, with $\alpha$ as unit.  The definitions of $FP$ and
$\alpha_P$ are by induction on dimension:
\begin{itemize}
\item if $\pi$ is the unique element of $\pd(0)$ then $(FP)(\pi) = P(\pi)$
\item if $n\geq 1$ and $\pi\in\pd(n)$ then
$
% \begin{equation}	\label{eq:FP}
(FP)(\pi)
=
P(\pi) + 
\mr{Par}_{FP}(\pi)
$
\item $\alpha_{P,\pi}: P(\pi) \rIncl (FP)(\pi)$ is inclusion as the first
summand, for all $\pi$
\item if $n\geq 1$ and $\pi\in\pd(n)$ then the source map $s: (FP)(\pi) \go
(FP)(\bdry\pi)$ is defined on the first summand of $(FP)(\pi)$ as the
composite
\[
P(\pi) \goby{s} P(\bdry\pi) \goby{\alpha_{P,\bdry\pi}} (FP)(\bdry\pi)
\]
and on the second summand as first projection; the target map is defined
dually.
\end{itemize}
The globularity equations hold, so $FP$ forms a collection.  Clearly
$\alpha_P: P \go FP$ is a map of collections.  The contraction $\kappa^P$
on $FP$ is defined by taking 
\[
\kappa^P_\pi: \mr{Par}_{FP}(\pi) \go (FP)(\pi)
\]
($n\geq 1$, $\pi\in\pd(n)$) to be inclusion as the second summand.  It is
easy to check that $FP$, $\kappa^P$ and $\alpha^P$ have the requisite
universal property: so $U$ has a left adjoint.

\minihead{$V$ is finitary and monadic}

The functor $T$ is finitary, by~\ref{thm:omega-forgetful-properties}.  This
implies by~\ref{thm:free-gen} that the monad $T$ is suitable, and so
by~\ref{thm:free-fixed} that $V$ is monadic.  It also implies by the
`moreover' of~\ref{propn:free-refined} that the monad induced by $V$ and
its left adjoint is finitary.  If a category has colimits of a certain
shape and a monad on it preserves colimits of that shape, then so too does
the forgetful algebra functor; hence $V$ is finitary, as required.  

\begin{notes}

I thank Steve Lack and John Power for telling me that the result I
needed,~\ref{thm:comb-mon}, was in Kelly~\cite{KelUTT}, and Sjoerd Crans
for telling me exactly where.%
\index{globular operad!contraction@with contraction|)}

\end{notes}

\addtocontents{toc}{\contentsline {chapter}{ }{ }}

\backmatter

\normalsize

\chapter{Glossary of Notation}

\chapterquote{%
`It's a revealing thing, an author's index of his own work,' she informed
me.  `It's a shameless exhibition---to the \emph{trained} eye.'}{%
Vonnegut~\cite{Von}}

\noindent
Page numbers indicate where the term is defined; multiple page numbers
mean that the term is used in multiple related senses.  Boxes ($\Box$)
stand for some or all of `lax', `colax', `wk', and `str'.

\subsection*{Latin letters}

\begin{glosslist}
\guse{\gr{Ab}}{\Ab}{category of abelian groups}
\guse{\gr{Algplainmulti}, \gr{Alggenmti}}%
{\Alg\blank}{category of algebras for a multicategory}
\guse{\gr{Algwk}--\gr{Algstr}}{\Alg_\Box\blank}%
{category of algebras for a $\Cat$-operad}
\guse{\gr{BBatopd}}{B}{initial operad with coherence and system of compositions}
\guse{\gr{Bilax}}{\fcat{Bicat}_\Box}{category of bicategories}
\guse{\gr{blankBi}}{\blank\hyph\fcat{Bicat}_\Box}%
{category of `bicategories' according to some theory}
\guse{\gr{Braid}}{\fcat{Braid}}{category of braids}
\guse{\gr{CartMndlax}, \gr{CartMndwk}}{\fcat{CartMnd}_\Box}%
{2-category of cartesian monads and lax maps}
\guse{\gr{Cat}, \gr{Cat2cat}}{\Cat}{(2-)category of small categories}
\guse{\gr{wkdblCat2}}{\Cat_2}{weak double category of categories}
\guse{\gr{Catinternal}}{\Cat\blank}{category of internal categories}
\guse{\gr{VCat}}{\blank\hyph\Cat}{category of enriched categories}
\guse{\gr{CatAlgwk}, \gr{TrimCatAlg}}{\fcat{CatAlg}\blank}%
{category of categorical algebras}
\guse{\gr{CatOperad}, \gr{CatOperad2cat}}{\Cat\hyph\Operad}%
{(2-)category of $\Cat$-operads}
\guse{\gr{bicatco}}{\blank^\mr{co}}{dual of bicategory, reversing 2-cells}
\guse{\gr{codcat}, \gr{codgenmti}}{\cod}{codomain map of (multi)category}
\guse{\gr{CommMon}}{\fcat{CommMon}}{category of commutative monoids}
\guse{\gr{compcat}, \gr{compenr}, \gr{compgenmti}}{\comp}%
{composition map in (multi)category}
\guse{\gr{Cone}}{\fcat{Cone}\blank}{cone on space}
\guse{\gr{ctr}}{\ctr}{operad of classical trees}
\guse{\gr{disccat}, \gr{discintcat}}{D}{discrete category}
\guse{\gr{augsimplexcat}}{\scat{D}}{augmented simplex category}
\guse{\gr{disk}}{D^m}{$m$-dimensional closed disk/ball}
\guse{\gr{littleddisks}}{\ldisks_d}{little $d$-disks operad}
\guse{\gr{DNn}, \gr{DNnweak}}{D^N_n}%
{discrete $N$-dimensional structure on $n$-dim'l structure}
\guse{\gr{DFibcat}}{\fcat{DFib}\blank}{category of discrete fibrations}
\guse{\gr{domcat}, \gr{domgenmti}}{\dom}{domain map of (multi)category}
\guse{\gr{DOpfibcat}, \gr{DOpfibgen}}%
{\fcat{DOpfib}\blank}{category of discrete opfibrations}
\guse{\gr{edgeftr}}{E}{edge functor}
\guse{\gr{nthopecat}}{\Eee_n}{$n$th opetopic category}
\guse{\gr{Endplainmulti}, \gr{Endplainoperad}, \gr{Endgenmti}}%
{\End\blank}%
{endomorphism operad or multicategory}
\guse{\gr{Fam}}{\Fam\blank}{category of families of objects}
\guse{\gr{FatCommMon}}{\fcat{FatCommMon}}{category of fat commutative monoids}
\guse{\gr{FatSymMulticat}}{\fcat{FatSymMulticat}}%
{category of fat symmetric multicategories}
\guse{\gr{fc}, \gr{fcprecise}}{\fc}{free category monad}
\guse{\gr{fcV}}{\fc_{\cat{V}}}{free $\cat{V}$-enriched category monad}
\guse{\gr{G}, \gr{Gomega}}{\scat{G}, \scat{G}_\omega}%
{category on which globular sets are presheaves}
\guse{\gr{Gn}}{\scat{G}_n}{category on which $n$-globular sets are presheaves}
\guse{\gr{abbrevGph}}{\Gph}{category of directed graphs}
\guse{\gr{VGph}}{\blank\hyph\Gph}{category of enriched graphs}
\guse{\gr{nGph}}{n\hyph\Gph}{category of $n$-graphs}
\guse{\gr{TGraph}}{\blank\hyph\Graph}{category of generalized graphs}
\guse{\gr{HMonCat}}{\fcat{HMonCat}}{category of homotopy monoidal categories}
\guse{\gr{Homgenmti}}{\HOM}{`hom' generalized graph}
\guse{\gr{indisccat}, \gr{indiscgenmti}}{I}{indiscrete (multi)category}
\guse{\gr{unitobj}}{I}{unit in monoidal category}
\guse{\gr{initsymopd}}{I}{initial symmetric operad}
\guse{\gr{istrncat}, \gr{irefl}}{i}{identities map}
\guse{\gr{idcatmodule}}{I_A}{identity module over category $A$}
\guse{\gr{IMm}}{I^M_m}{indiscrete $M$-category on $m$-category}
\guse{\gr{idscat}, \gr{idsenr}, \gr{idsgenmti}}{\ids}%
{identities map in (multi)category}
\guse{\gr{Jn}}{J_n}{localization of $n$-dimensional structure}
\guse{\gr{assK}}{K}{operad of associahedra}
\guse{\gr{initL}}{(L,\lambda)}{initial globular operad with contraction}
\guse{\gr{initLn}}{(L_n, \lambda_n)}{initial globular $n$-operad with contraction}
\guse{\gr{ninebicat}}{\fcat{LaxBicat}_\Box}{category of lax bicategories}
\guse{\gr{ninemoncat}}{\fcat{LaxMonCat}_\Box}{category of lax monoidal categories}
\guse{\gr{Mfreecommmon}}{M}{free commutative monoid monad}
\guse{\gr{Matk}}{\fcat{Mat}_k}{category of natural numbers and matrices over $k$}
\guse{\gr{Mndlax}}{\fcat{Mnd}_\Box}{2-category of monads}
\guse{\gr{RMod}}{\blank\hyph\fcat{Mod}}{category of modules}
\guse{\gr{Monofcat}}{\fcat{Mon}\blank}{category of monoids in monoidal category}
\guse{\gr{fcmtiMon}}{\Mon\blank}{\fc-multicategory of monads in \fc-multicategory}
\guse{\gr{MClax}}{\fcat{MonCat}_\Box}{category of monoidal categories}
\guse{\gr{blankMonCatwk},\gr{blankMonCatlax}}{\blank\hyph\MCwk}%
{category of `monoidal categories' according to some theory}
\guse{\gr{Monoid}}{\fcat{Monoid}}{category of monoids}
\guse{\gr{Multicat}, \gr{Multicat2cat}}{\Multicat}%
{(2-)category of plain multicategories}
\guse{\gr{TMulticat}, \gr{TMulticat2cat}}{\blank\hyph\Multicat}%
{(2-)category of generalized multicategories}
\guse{\gr{Multicatfixedobjs}}{\blank\hyph\Multicat_E}{category of generalized
  multicategories with objects $E$}
\guse{\gr{nat}}{\nat}{natural numbers (including $0$)}
\guse{\gr{catofopes}}{\scat{O}}{category of opetopes}
\guse{\gr{Oatmostn}}{\scat{O}(n)}{category of opetopes of dimension at most $n$}
\guse{\gr{nopes}}{O_n}{set of $n$-opetopes}
\guse{\gr{obcat}}{\ob}{object-set of category}
\guse{\gr{OC}}{\fcat{OC}}{category of globular operads with contraction}
\guse{\gr{nOC}}{n\hyph\fcat{OC}}{category of globular $n$-operads with contraction}
\guse{\gr{nOP}}{n\hyph\fcat{OP}}{category of globular $n$-operads with precontraction}
\guse{\gr{catop}}{\blank^\op}{opposite of category}
\guse{\gr{bicatop}}{\blank^\op}{dual of bicategory, reversing 1-cells}
\guse{\gr{Operad}}{\Operad}{category of plain operads}
\guse{\gr{VOperad}}{\blank\hyph\Operad}{category of operads in symmetric multicategory}
\guse{\gr{TOperad}}{\blank\hyph\Operad}{category of generalized operads}
\guse{\gr{nOperad}}{n\hyph\Operad}{category of globular $n$-operads}
\guse{\gr{widepbshape}}{\scat{P}_K}{shape of $K$-ary wide pullback}
\guse{\gr{Par}}{\mr{Par}_q\blank}{parallel pairs of cells for map $q$}
\guse{\gr{Parcoll}}{\mr{Par}_P\blank}{parallel pairs of cells in collection $P$}
\guse{\gr{PDn}}{\PD{n}}{structured category of opetopic $n$-pasting diagrams}
\guse{\gr{Pdn}}{\Pd{n}}{category of opetopic $n$-pasting diagrams}
\guse{\gr{pd}}{\pd}{strict $\omega$-category of globular pasting diagrams}
\guse{\gr{pri}}{\mr{pr}_i}{$i$th projection}
\guse{\gr{Rrefl}}{\scat{R}}{category on which reflexive globular sets are presheaves}
\guse{\gr{Rngraphs}}{R_n}{restriction of $(n+1)$-graph to $n$-graph}
\guse{\gr{RMm}}{R^M_m}{restriction of $M$-dimensional structure to
  $m$ dimensions}
\guse{\gr{fcRing}, \gr{fcRing2}}{\fcat{Ring}, \fcat{Ring}_2}%
{\fc-multicategory of rings}
\guse{\gr{globsce}, \gr{sinopeset}}{s}{source map}
\guse{\gr{dsphere}}{S^d}{$d$-sphere}
\guse{\gr{Sstrncat}}{S_n}{functor forgetting $(n+1)$-cells of $(n+1)$-category}
\guse{\gr{symgp}}{S_n}{symmetric group on $n$ letters}
\guse{\gr{Set}}{\Set}{category of sets}
\guse{\gr{fcSet2}}{\fcat{Set}_2}{\fc-multicategory of sets}
\guse{\gr{strictcover}, \gr{strictcoverbicat}}{\st}%
{strict cover}
\guse{\gr{ninebicat}}{\strcat{2}_\Box}{category of strict 2-categories}
\guse{\gr{StrMonCat}, \gr{ninemoncat}}{\fcat{StrMonCat}_\Box}%
{category of strict monoidal categories}
\guse{\gr{strcatn}}{\strcat{n}}{category of strict $n$-categories}
\guse{\gr{strcatomega}}{\strcat{\omega}}{category of strict $\omega$-categories}
\guse{\gr{Struc}}{\blank\hyph\Struc}{category of structured categories}
\guse{\gr{StTr}}{\fcat{StTr}}{category of stable trees}
\guse{\gr{SymMulticat}}{\fcat{SymMulticat}}{category of symmetric multicategories}
\guse{\gr{SymOpd}, \gr{SymOpdsym}}{\SymOpd}{operad of symmetries}
\guse{\gr{Tmonad}}{T = (T, \mu, \eta)}{usual name for cartesian monad}
\guse{\gr{strommon}}{(T,\mu,\eta)}%
{free strict $\omega$-category monad 
(in Ch.\ \ref{ch:globular}, \ref{ch:a-defn}, \ref{ch:other-defns})}
\guse{\gr{nthopemonad}}{T_n}{$n$th opetopic monad}
\guse{\gr{gmn}}{\gm{n}}{free strict $n$-category monad}
\guse{\gr{globtgt}, \gr{opetgt}, \gr{tinopeset}}{t}{target map}
\guse{\gr{Top}, \gr{Topcat}}{\Top}{($\omega$-)category of topological spaces}
\guse{\gr{Topstar}}{\Top_*}{category of based spaces}
\guse{\gr{TR}}{\TR}{$\Cat$-operad of trees}
\guse{\gr{Tr}}{\Tr}{category of trees}
\guse{\gr{tr}}{\tr}{operad of trees}
\guse{\gr{fgtslice}}{U_E}{forgetful functor $\Eee/E \go \Eee$}
\guse{\gr{ninebicat}}{\fcat{UBicat}_\Box}{category of unbiased bicategories}
\guse{\gr{ninemoncat}}{\fcat{UMonCat}_\Box}{category of unbiased monoidal categories}
\guse{\gr{Vmonmti}}{V}{forgetful functor from monoidal categories to multicategories}
\guse{\gr{vxftr}}{V}{vertex functor}
\guse{\gr{Wfreemon}}{W}{free monoid monad}
\guse{\gr{wkncat}}{\wkcat{n}}{category of weak $n$-categories}
\guse{\gr{wkomegaCat}}{\wkcat{\omega}}{category of weak $\omega$-categories}
\guse{\gr{FamYon}}{\fcat{y}}{`Yoneda' embedding for familially representable functors}
\end{glosslist}

\subsection*{Greek letters}

\begin{glosslist}

\guse{\gr{monalpha}, \gr{bicatass}}{\alpha}{associativity isomorphism}
\guse{\gr{gammalaxmoncat}, \gr{gammabicat}}{\gamma}%
{coherence map in lax bicategory}
\guse{\gr{Delta}}{\Delta}{category of nonempty finite totally ordered sets}
\guse{\gr{Deltainj}}{\Delta_\mr{inj}}{subcategory of $\Delta$ consisting
  only of injections}
\guse{\gr{Deltadiag}}{\Delta}{diagonal functor}
\guse{\gr{Deltaomega}}{\Delta_\omega}{category of globular pasting diagrams}
\guse{\gr{cohdelta}}{\delta}{coherence isomorphisms in a $\Sigma$-monoidal category}
\guse{\gr{deltasys}}{\delta^\blob_\blob}{system of binary compositions}
\guse{\gr{bdry}}{\bdry}{boundary of globular pasting diagram}
\guse{\gr{iotalaxmoncat}, \gr{iotabicat}}{\iota}%
{coherence map in lax bicategory}
\guse{\gr{iotacone}}{\iota}{embedding of space in its cone}
\guse{\gr{iotasingle}}{\iota_m}{globular pasting diagram resembling single $m$-cell}
\guse{\gr{monlambda}, \gr{bicatlambda}}{\lambda}%
{left unit isomorphism}
\guse{\gr{initL}}{(L,\lambda)}{initial globular operad with contraction}
\guse{\gr{initLn}}{(L_n, \lambda_n)}{initial globular $n$-operad with contraction}
\guse{\gr{corollan}}{\nu_n}{$n$-leafed corolla}
\guse{\gr{arbprod}}{\prod}{product}
\guse{\gr{fundn}, \gr{Pi1top}, \gr{Pitwobicat}, \gr{fundubicat}, \gr{Pin}}%
{\Pi_n}%
{fundamental $n$-groupoid of space}
\guse{\gr{fundomega}, \gr{Piomega}}{\Pi_\omega}%
{fundamental $\omega$-groupoid of space}
\guse{\gr{monrho}, \gr{bicatrho}}{\rho}{right unit isomorphism}
\guse{\gr{fatsum}}{\sum}{summation in fat commutative monoid}
\guse{\gr{SigmaMCbi}, \gr{Sigmamtifc}, \gr{suspglobpd}, \gr{suspobjgraph}}%
{\Sigma}{suspension}
\guse{\gr{Sigmacoproduct}}{\Sigma}{coproduct of family of sets}
\guse{\gr{Sigmac}}{\Sigma_\mr{c}}{signature for classical, biased structures}
\guse{\gr{chi}}{\chi_k}{1-pasting diagram made up of $k$ abutting arrows}
\guse{\gr{Omegaloop}}{\Omega}{loop space functor}
\guse{\gr{omegaGph}}{\omega\hyph\Gph}{category of $\omega$-graphs}
\guse{\gr{omegaOC}}{\omega\hyph\fcat{OC}, \omega\hyph\fcat{OP}}%
{category of globular operads with contraction}
\guse{\gr{omegaOperad}}{\omega\hyph\Operad}{category of globular operads}
\guse{\gr{wkomegaCat}}{\wkcat{\omega}}{category of weak $\omega$-categories}
\end{glosslist}

\subsection*{Other symbols}

\begin{glosslist}

\guse{\gr{C0cat}, \gr{C0multicat}, \gr{C0genmti}, \gr{A0fat}}%
{\blank_0}{objects of (multi)category}
\guse{\gr{C1cat}, \gr{C1genmti}}{\blank_1}{arrows of (multi)category}

\guse{\gr{thetabar}, \gr{thetabaropd}, \gr{actionglobopd}}{\ovln{\blank}}%
{action of multicategory on algebra}
\guse{\gr{mate}}{\ovln{\blank}}{mate}
\guse{\gr{reppi}}{\rep{\blank}}{globular set representing pasting diagram}

\guse{\gr{lowertrunc}}{\blank_{\trunc{n}}}{restriction of globular set to $n$ dimensions}
\guse{\gr{starfreemon}}{\blank^*}{free monoid monad}
\guse{\gr{starenr}, \gr{mndmapindftr}, \gr{laxmtiindftr},
  \gr{colaxmtiindftr}}%
{\blank_*}%
{induced functor (various contexts)}

\guse{\gr{plusalgtomulti}}{\blank^+}{functor turning algebras into multicategories}
\guse{\gr{BDslicemti}}{\blank^+}{slice of symmetric multicategory}
\guse{\gr{plusofcaty}}{\Eee^+}{category of generalized graphs in $\Eee$}
\guse{\gr{plusofmonad}}{T^+}{free $T$-multicategory monad}
\guse{\gr{pluscatyfo}}{\Eee^+_E}{category of generalized graphs with
  object-of-objects $E$}
\guse{\gr{plusmonadfo}}{T^+_E}{free ($T$-multicategory with object-of-objects
  $E$) monad}

\guse{\gr{ftrcat}, \gr{ftrbicat}}{\ftrcat{\ }{\ }}{functor (bi)category}
\guse{\gr{homsetset}}{[\ ,\ ]}{set of functions between sets}
\guse{\gr{homset}}{\cat{A}(A,B)}{hom-set in category}
\guse{\gr{multicathomset}}{C(a_1, \ldots, a_n; a)}{hom-set in plain multicategory}
\guse{\gr{fathomset}}{C((a_x)_{x\in X}; b)}{hom-set in fat symmetric multicategory}

\guse{\gr{vertcompnat}}{\of}{vertical composition of natural transformations}
\guse{\gr{nfoldcompbicat}}{\of_n}{$n$-fold composition}
\guse{\gr{nfoldcompbicatinfix}}{(f_n \of\cdots\of f_1)}{$n$-fold composite}
\guse{\gr{multicomposite}, \gr{opdcomposite}}%
{\theta \of (\theta_1, \ldots, \theta_n)}%
{composite in plain multicategory or operad}
\guse{\gr{fatcomp}}{\phi \of (\theta_y)_{y\in Y}}{composite in fat symmetric multicategory}
\guse{\gr{circlei}}{\of_i}{composition in $i$th place (`circle-$i$')}
\guse{\gr{ofdim}}{\of_p}{composition in $n$-category by gluing along $p$-cells}
\guse{\gr{horizcompnat}, \gr{bicatstar}}{*}{horizontal composition}
\guse{\gr{nfoldstar}}{(\alpha_n * \cdots * \alpha_1)}{$n$-fold horizontal
  composition} 

\guse{\gr{terminal}}{1}{terminal object}
\guse{\gr{onestrncat}, \gr{multiidentity}, \gr{opdidentity},
  \gr{idglobopd}, \gr{fatids}}%
{1}%
{identity map}

\guse{\gr{otimes}}{\otimes}{tensor product in monoidal category}
\guse{\gr{otimesninfix}}{(a_1 \otimes \cdots \otimes a_n)}{$n$-fold tensor of
  objects $a_i$}
\guse{\gr{otimesn}}{\otimes_n}{$n$-fold tensor product}
\guse{\gr{otimestree}}{\otimes_\tau}{`tensor' product of type $\tau$}

\guse{\gr{dotsigma}, \gr{fatcdot}}{\cdot}%
{action of permutation or bijection on map}

\guse{\gr{slicecat}}{\cat{E}/E}{slice category }
\guse{\gr{sliceofmonadbyalg}}{T/h}{slice of monad $T$ by algebra $h$}
\guse{\gr{catelts}}{\scat{A}/X}%
{category of elements of presheaf $X$ on $\scat{A}$}
\guse{\gr{genmtielts}}{C/h}%
{multicategory of elements of algebra $h$ for $C$}

\guse{\gr{utree}}{\utree}{trivial (unit) tree}
\guse{\gr{tupleoftrees}}{(\tau_1, \ldots, \tau_n)}{tree formed by joining
  $\tau_1, \ldots, \tau_n$}
\guse{\gr{inducedmonad}}{(T^C, \mu^C, \eta^C)}{monad induced by
  $T$-multicategory $C$}
\guse{\gr{catpower}}{\cat{E}^S}{power of category}
\guse{\gr{Sp}, \gr{Spwkdbl}}{\Sp{\Eee}{T}}%
{bicategory or weak double category of $T$-spans}

\guse{\gr{lwrn}}{\lwr{n}}{$n$-element totally ordered set}
\guse{\gr{upr}}{\upr{n}}{$(n+1)$-element totally ordered set}
\guse{\gr{bang}}{!}{unique map to terminal object}
\guse{\gr{compbang}}{e_!}{composition with map $e$}
\guse{\gr{rightangle}}{\mbox{\huge$\lrcorner$}}{indicates that square is a pullback}

\guse{\gr{iso}}{\iso}{isomorphism}
\guse{\gr{eqvcat}, \gr{eqvinbicat}}{\eqv}{equivalence}
\guse{\gr{simA}}{\sim}{equivalence under a system of equations}
\guse{\gr{ladj}}{\ladj}{is left adjoint to}

\guse{\gr{coend}}{\int}{coend}

\guse{\gr{wedge}}{\wej}{wedge of spaces}
\guse{\gr{binprod}}{\times}{binary product}
\guse{\gr{pbtimes}}{\times_C}{pullback over object $C$}
\guse{\gr{bincoprod}}{+}{binary coproduct}
\guse{\gr{arbprod}}{\prod}{product}
\guse{\gr{arbcoprod}}{\coprod}{coproduct}

\end{glosslist}

% \newcommand{\seealso}[1]{see{\textit{also }#1}}

% I put the following into .emacs:
% 
% (global-set-key [f1] [?% return ?% return return right backspace up ?\\ ?i ?n ?d ?e ?x ?{ ?} return ?% left left left])

\index{truncation|see{restriction}}
\index{unitrivalent|see{tree, classical}}
\index{multiple outputs|see{many in, many out}}
\index{T-structured@$T$-structured category|see{structured category}}
\index{opfibration|see{fibration}}
\index{T-multicategory@$T$-multicategory|see{generalized multicategory}}
\index{multicategory!generalized|see{generalized multicategory}}
\index{T-operad@$T$-operad|see{generalized operad}}
\index{operad!generalized|see{generalized operad}}
\index{natural transformation|see{transformation}}
\index{transformation!of monoidal categories|see{monoidal transformation}}
\index{PROP|see{PRO}}
\index{bimodule|see{module}}
\index{cubical!n-category@$n$-category|see{$n$-tuple category}}
\index{lax something|see{something, lax}}
\index{colax something|see{something}}
\index{weak something|see{something, weak}}
\index{strict something|see{something, strict}}
\index{infinity-category@$\infty$-category|see{$\omega$-category}}
\index{left-handed|see{handedness}}
\index{rectification|see{coherence}}
\index{triple|see{monad}}
\index{poset|see{order}}
\index{morphism|see{arrow; cell}}
\index{category!enriched|see{enrichment}}
\index{multicategory!enriched|see{enrichment}}
\index{transformation!enriched|see{enrichment}}
\index{Godement law|see{interchange}}
\index{multicategory!enriched|see{enrichment}}
\index{identities, exclusion of|see{pseudo-operad}}
\index{multicategory!T-@$T$-|see{generalized multicategory}} 
\index{operad!T-@$T$-|see{generalized operad}} 
\index{symmetry|see{handedness}}
\index{functor!monoidal categories@of monoidal categories|see{monoidal
  functor}} 
\index{roof|see{span}}
% \index{category!structured|see{structured category}}
\index{T-n@$T_n$|see{opetopic monad}}
\index{generalized multicategory!symmetric multicategory@\vs.\
  symmetric multicategory|see{multicategory, symmetric \vs.\ generalized}}
% \index{opetopic!pasting diagram|see{pasting diagram, opetopic}}
\index{unbiased something|see{something, unbiased}}
\index{operad!contraction@with contraction|see{globular operad, with
    contraction}} 
\index{operad!globular|see{globular operad}}
\index{n-operad@$n$-operad|see{$n$-globular operad}}
\index{globular operad!n-@$n$-|see{$n$-globular operad}}
\index{monad!gluing of|see{distributive law}}
\index{filler|see{horn}}
\index{topological space!n-category from@$n$-category from|see{fundamental}}
\index{Gamma-space@$\Gamma$-space|see{$\Delta$-space}}
\index{map|see{arrow; cell}}

% \index{|see{}}

\index{associativity|seealso{$n$-category, weak \vs.\ strict}}
\index{two-category@2-category|seealso{bicategory}}
\index{bicategory|seealso{fundamental 2-groupoid; 2-category}}
\index{three-category@3-category|seealso{tricategory}}
\index{tricategory|seealso{3-category}}
\index{simplex category $\Delta$|seealso{augmented simplex category}}
\index{categorification|seealso{weakening}}
\index{representable functor|seealso{familial representability}}
\index{order|seealso{simplex category; augmented simplex
    category}}
\index{plus construction $\blank^+$|seealso{slice}}
\index{multicategory|seealso{generalized multicategory; monoidal category; operad }}
\index{operad|seealso{generalized operad; monad, operadic; multicategory}}
\index{monoidal category|seealso{bicategory}}
\index{arrow|seealso{cell}}
\index{groupoid|seealso{fundamental; $n$-groupoid; $\omega$-groupoid}}
\index{graph|seealso{$n$-graph}}
\index{globular set|seealso{$n$-globular set}}
\index{object|seealso{cell}}
\index{disk|seealso{operad, little disks}}
\index{multicategory!underlying|seealso{monoidal category}}
\index{homotopy-algebraic structure|seealso{$A_\infty$; loop space; monoidal category, homotopy}}
\index{generalized operad|seealso{generalized multicategory}}
\index{path|seealso{operad of path reparametrizations}}
\index{opalgebra|seealso{coalgebra}}
\index{coalgebra|seealso{opalgebra}}
\index{universal|seealso{pre-universal}}
\index{omega-category@$\omega$-category|seealso{$n$-category}}
\index{Street, Ross|seealso{tricategory}}
\index{Power, John|seealso{tricategory}}
\index{order!non-canonical|seealso{handedness}}
\index{top dimension|seealso{tame}}
\index{braid|seealso{monoidal category, braided}}
\index{fc-multicategory@$\fc$-multicategory|seealso{double category, weak}}
\index{equivalence|seealso{biequivalence}}
\index{n-groupoid@$n$-groupoid|seealso{fundamental}}
\index{omega-groupoid@$\omega$-groupoid|seealso{fundamental}}
\index{monoid!commutative|seealso{Eckmann--Hilton}}
\index{monoidal category!multicategory@\vs.\ multicategory|seealso{multicategory, underlying}}
\index{omega-category@$\omega$-category!strict!free|seealso{globular operad}}
\index{Johnstone, Peter|seealso{strongly regular theory}}
\index{Carboni, Aurelio|seealso{strongly regular theory}}

% \index{|seealso{}}

\small
\addcontentsline{toc}{chapter}{Index}
\printindex

\end{document}